%% file: main.tex
\documentclass[11pt,twoside]{book}

\usepackage{amsmath, amsthm,amsfonts,amssymb,mathtools,mathdots,scalerel,mathrsfs,stmaryrd,cmll}
\usepackage[utf8]{inputenc}
\usepackage[T1]{fontenc}
\usepackage{mlmodern, tikz-cd, quiver}
\usepackage[many]{tcolorbox}
\usepackage{emptypage}

\usepackage{enumitem, imakeidx, tocloft, csquotes, lipsum, tabularx}
\usepackage[hidelinks, backref=page]{hyperref}
\usepackage{tikz, fancyhdr, xparse, xcolor}
\usepackage[british]{babel}
\usepackage[a4paper, top=4.5cm, bottom=4.5cm, left=4cm, right=4cm, asymmetric]{geometry}
\usepackage[textwidth=3cm, textsize=small, colorinlistoftodos]{todonotes}

\usepackage{ahtitle, titlesec}
\usepackage{etoolbox}
\makeatletter
\patchcmd{\ttlh@hang}{\parindent\z@}{\parindent\z@\leavevmode}{}{}
\patchcmd{\ttlh@hang}{\noindent}{}{}{}
\makeatother

\input{macros}

\input{style}

\title{Combinatorics of higher-categorical diagrams}

\author{Amar Hadzihasanovic}

\institution{Tallinn University of Technology} 

\makeindex[name=counterex, columns=1, title={Index of counterexamples}]
\makeindex

\begin{document}

{$\quad$}

\vspace{15pt}

\maketitle 

\noindent\makebox[\textwidth][r]{%
	\begin{minipage}[t]{.75\textwidth}\begin{guide}
\small \cemph{Abstract.}
This is a book on higher-categorical diagrams, including pasting diagrams.
It aims to provide a thorough and modern reference on the subject, collecting, revisiting and expanding results scattered across the literature, informed by recent advances and practical experience with higher-dimensional diagram rewriting.

We approach the subject as a kind of directed combinatorial topology: a diagram is a map from a ``directed cell complex'', encoded combinatorially as a face poset together with orientation data.
Unlike previous expositions, we adopt from the beginning a functorial viewpoint, focussing on morphisms and categorical constructions.
We do not tie ourselves to a specific model of higher categories, and instead treat diagrams as independent combinatorial structures that admit functorial interpretations in various contexts.

Topics covered include the theory of layerings of diagrams; acyclicity properties and their consequences; constructions including Gray products, suspensions, and joins; special shapes such as globes, oriented simplices, cubes, and positive opetopes; the interpretation of diagrams in strict $\omega$\nbd categories and their geometric realisation as simplicial and CW complexes; and Steiner's theory of directed chain complexes.
\end{guide}\end{minipage}}

\vspace{20pt}

\makeaftertitle

\normalsize \pagestyle{empty}

\cleardoublepage 
\setcounter{tocdepth}{1}

\newgeometry{left=5.5cm, right=5.5cm}
\tableofcontents
\restoregeometry

\pagestyle{firstpage}
\include{introduction}
\cleardoublepage

\pagestyle{mainpage}
\include{order}
\include{ogposets}
\include{molecules}
\include{layerings}
\include{omegacats}
\include{maps}
\include{constructions}
\include{acyclic}
\include{special}
\include{geometric}
\include{steiner}

\cleardoublepage
\pagestyle{plain}
\include{backmatter}

\end{document}

%% file: macros.tex
\newcommand\eqdef{\coloneqq}
\newcommand\nbd{\nobreakdash-\hspace{0pt}}
\newcommand\idd[1]{\mathrm{id}_{#1}}
\newcommand\bigid[1]{\mathrm{Id}_{#1}}
\newcommand\invrs[1]{#1^{-1}}

\newcommand\after{\circ}

\newcommand\incl{\hookrightarrow}
\newcommand\incliso{\stackrel{\sim}{\hookrightarrow}}
\newcommand\iso{\stackrel{\sim}{\rightarrow}}

\newcommand\restr[2]{{#1}{\raisebox{0pt}{$|_{#2}$}}}

\newcommand\set[1]{\left\{ {#1} \right\}}
\newcommand\size[1]{\left|{#1}\right|}
\newcommand\powerset[1]{\mathscr{P}{#1}}
\newcommand\isocl[1]{[#1]}
\newcommand\height[1]{\mathrm{h}({#1})}

\newcommand\clas[1]{\mathscr{#1}}

\newcommand\pfw[1]{{#1}_*}
\newcommand\pb[1]{{#1}^*}

\newcommand{\overbar}[1]{{\mkern 1.5mu\overline{\mkern-1.5mu#1\mkern-1.5mu}\mkern 1.5mu}}

\newcommand\posnat{\mathbb{N} \setminus \set{0}}

\newcommand\slice[2]{{#1}/{\raisebox{-2pt}{$#2$}}}
\newcommand\join{\,{\star}\,}

\newcommand\opp[1]{#1^\mathrm{op}}
\newcommand\coo[1]{#1^\mathrm{co}}
\newcommand\optot[1]{#1^\circ}
\newcommand\dual[2]{\fun{D}_{#1}{#2}}

\newcommand\cat[1]{\mathbf{#1}}
\newcommand\smcat[1]{\mathscr{#1}}

\newcommand\fun[1]{\mathsf{#1}}

\newcommand\property[1]{\mathrm{#1}}

\newcommand\order[2]{#2^{(#1)}}

\newcommand\ogpos{\cat{ogPos}}
\newcommand\ogposloc{\ogpos_\mathit{le}}
\newcommand\ogposbot{\ogpos^+}
\newcommand\otgpos{\cat{otgPos}}

\newcommand\poscat{\cat{Pos}}

\newcommand\gph{\cat{Gph}}

\newcommand\posclos{\poscat_\mathit{cl}}
\newcommand\posbot{\poscat^+}
\newcommand\posclosbot{\posclos^+}
\newcommand\posloc{\poscat_\mathit{le}}

\newcommand\rngph[1]{{#1}\cat{Gph}_\mathit{ref}}
\newcommand\romegagph{\rngph{\omega}}
\newcommand\ncat[1]{{#1}\cat{Cat}}
\newcommand\omegacat{\ncat{\omega}}
\newcommand\polmap{\cat{Pol}_{\downarrow}}

\newcommand\thetacat{\Theta}
\newcommand\simplexcat{\Delta}

\newcommand\rdcpx{\cat{RDCpx}_{=}}
\newcommand\rdcpxmap{\cat{RDCpx}_{\downarrow}}
\newcommand\rdcpxcomap{\cat{RDCpx}_{\uparrow}}
\newcommand\rdcpxiso{\cat{RDCpx}_\mathit{iso}}
\newcommand\rdcpxmapfa{\rdcpxmap^\mathit{fa}}
\newcommand\rdcpxmapac{\rdcpxmap^\mathit{ac}}
\newcommand\rdcpxcomapac{\rdcpxcomap^\mathit{ac}}

\newcommand\sset{\cat{sSet}}
\newcommand\ktop{k\cat{Top}}
\newcommand\rcpx{\cat{RCpx}_{\downarrow}}
\newcommand\cwpos{\cat{cwPos}}

\newcommand\chaug{\cat{Ch}^+}
\newcommand\dchaug{\cat{DCh}^+}
\newcommand\stcpx{\cat{DCh}^+_{\mathit{St}}}
\newcommand\sstcpx{\cat{DCh}^+_{\mathit{sSt}}}

\newcommand\hasse[1]{\mathscr{H}{#1}}
\newcommand\hasseo[1]{\vec{\mathscr{H}}{#1}}

\DeclareMathOperator{\clos}{cl}
\newcommand\clset[1]{\mathrm{cl}\set{#1}}

\newcommand\codim[2]{\mathrm{codim}_{#2}(#1)}
\newcommand\maxel[1]{\mathscr{M}\!\mathit{ax}\,#1}
\newcommand\minel[1]{\mathscr{M}\!\mathit{in}\,#1}
\newcommand\grade[2]{#2_{#1}}
\DeclareMathOperator{\lydim}{lydim}
\DeclareMathOperator{\frdim}{frdim}
\newcommand{\cdim}[1]{\mathrm{dim}_{#1}}

\newcommand\bound[2]{\partial_{#1}^{#2}}
\newcommand\faces[2]{\Delta_{#1}^{#2}}
\newcommand\cofaces[2]{\nabla_{#1}^{#2}}
\newcommand\inter[1]{\mathrm{int}\,#1}

\newcommand\cp[1]{\,{\scriptstyle\#}_{#1}\,}

\newcommand\gencp[1]{\,{\scriptstyle\widehat{\#}}_{#1}\,}
\newcommand\cpsub[2]{\triangleright_{#1}^{#2}}

\newcommand\cpiso[2]{\,{\scriptstyle\#}^{#2}_{#1}\,}
\newcommand\celto{\Rightarrow}
\newcommand\compos[1]{\langle#1\rangle}

\newcommand\submol{\sqsubseteq}

\newcommand\subs[3]{#1[#2/#3]}

\newcommand\graph[1]{\mathscr{G}#1}
\newcommand\flow[2]{\mathscr{F}_{#1}{#2}}
\newcommand\extflow[2]{\overbar{\mathscr{F}}_{#1}{#2}}
\newcommand\maxflow[2]{\mathscr{M}_{#1}{#2}}

\newcommand\layerings[2]{\mathscr{L}\!\mathit{ay}_{#1}{#2}}
\newcommand\orderings[2]{\mathscr{O}\!\mathit{rd}_{#1}{#2}}
\newcommand\lto[2]{\fun{o}_{#1, #2}}
\newcommand\step[2]{#2^{\mathsf{st}(#1)}}

\newcommand\prech{\preceq_\mathit{h}}
\newcommand\precv{\preceq_\mathit{v}}
\newcommand\precflow{\preceq}

\newcommand\gray{\otimes}
\newcommand\sus[1]{\fun{S}{#1}}
\newcommand\inj[1]{{#1} \join }
\newcommand\inr[1]{\join {#1} }

\newcommand\augm[1]{{{#1}_\bot}}
\newcommand\dimin[1]{{#1}_{\not\bot}}

\newcommand\ordcpx[1]{{#1}^\Delta}
\newcommand\realis[1]{\|{#1}\|}

\newcommand\geosim[1]{\Delta^{#1}}

\newcommand\nondeg[1]{\mathscr{N}\!\mathit{d}\,{#1}}
\newcommand\lnk[2]{\mathit{Lk}_{#1}\,{#2}}

\newcommand\gener[1]{\mathscr{#1}}
\newcommand\cwcom[2]{({#1},\gener{#2})}
\newcommand\fpos[1]{\mathfrak{F}{#1}}

\newcommand\chain[2]{{#1}_{#2}}
\newcommand\der{\mathrm{d}}
\newcommand\eau{\mathrm{e}}
\newcommand\freeab[1]{\mathbb{Z}{#1}}

\newcommand\dir[1]{{#1}^{\rightarrow}}
\newcommand\dfreeab[1]{\vec{\mathbb{Z}}{#1}}
\newcommand\freemon[1]{\mathbb{N}{#1}}

\newcommand\linea[1]{\fun{\lambda}{#1}}
\newcommand\nufun[1]{\fun{\nu}{#1}}
\newcommand\spanset[1]{\langle {#1} \rangle}
\DeclareMathOperator{\Ima}{Im}
\DeclareMathOperator{\supp}{supp}

\newcommand\gltab[3]{{#3}_{#1}^{#2}}
\newcommand\batom[1]{\langle #1 \rangle}

\newcommand\thearrow[1]{{#1}\vec{I}}
\newcommand\globe[1]{O^{#1}}
\newcommand\thetac[1]{\vartheta_{#1}}
\newcommand\disk[2]{D_{#1, #2}}

\newcommand\molec{\mathit{Mol}}
\newcommand\molecin[1]{\slice{\molec}{#1}}
\newcommand\atom{\mathit{Atom}}
\newcommand\atomin[1]{\slice{\atom}{#1}}

\newcommand\treeroot{\bullet}
\newcommand\branch[1]{\land({#1})}

\newcommand\simplex[1]{\vec{\Delta}^{#1}}
\newcommand\simface[1]{d^{#1}}
\newcommand\simdeg[1]{s^{#1}}
\newcommand\cube[1]{\vec{\square}^{#1}}
\newcommand\cubeface[2]{\delta^{#1}_{#2}}
\newcommand\cubedeg[1]{\sigma^{#1}}
\newcommand\cubeconn[2]{\gamma^{#1}_{#2}}
\newcommand\cubecp[1]{\, {\circ}_{#1} \,}
\newcommand\cubecomp[1]{\mu_{#1}}

\newcommand\omegatit{\texorpdfstring{$\omega$}{omega}}

\newcommand\skel[2]{\sigma_{\leq {#1}}#2}
\newcommand\cpable[1]{\times_{#1}}
\DeclareMathOperator{\cspan}{span}

\newcommand\ldiag[1]{\ell({#1})}

\newcommand\cemph[1]{{\color{\mycolor}\emph{#1}}}

%% file: style.tex
\newtheoremstyle{ittheorem}
  {\topsep}   
  {\topsep}   
  {\itshape}  
  {0pt}       
  {\sffamily \itshape \bfseries} 
  { ---}         
  {5pt plus 1pt minus 1pt} 
  {}          

\newtheoremstyle{itdfn}
  {\topsep}   
  {\topsep}   
  {}  
  {0pt}       
  {\sffamily \itshape \bfseries} 
  {}         
  {5pt plus 1pt minus 1pt} 
  {\thmnumber{#2}{\thmnote{\normalfont\ \ %
  {\sffamily(#3)}.}}}          

\newtheoremstyle{itrmk}
  {0.5\topsep}   
  {0.5\topsep}   
  {\normalfont}  
  {0pt}       
  {\sffamily \itshape} 
  { ---}         
  {5pt plus 1pt minus 1pt} 
  {}          
  
\newtheoremstyle{itexm}
  {0.5\topsep}      
  {0.5\topsep}      
  {\normalfont}     
  {0pt}             
  {\sffamily \itshape \bfseries \color{\mycolor}}        
  {\\}               
  {5pt plus 1pt minus 1pt} 
  {\thmname{#1} \thmnumber{#2}{\thmnote{\normalfont\ \ %
  {\sffamily(#3)}.}}}           
  
\makeatletter
  \renewcommand\@upn{\textit}
\makeatother

\newcommand\mycolor{red!75!black}

\tcbset{
    boxstyle/.style={%
        colframe=\mycolor, colback=gray!10,%
        left=2mm, right=2mm,%
        sharp corners,%
        breakable=true,%
        borderline west={1pt}{0pt}{\mycolor},%
        enhanced, boxrule=0pt, frame hidden}
    }
    
\theoremstyle{ittheorem}
\newtheorem{thm}{Theorem}[section]
\newtheorem{prop}[thm]{Proposition}
\newtheorem{cor}[thm]{Corollary}
\newtheorem{lem}[thm]{Lemma}

\theoremstyle{itdfn}
\newtheorem{dfn}[thm]{}
\theoremstyle{itrmk}
\newtheorem{rmk}[thm]{Remark}
\newtheorem{comm}[thm]{Comment}
\theoremstyle{itexm}
\newtheorem{exm}[thm]{Example}

\tcolorboxenvironment{exm}{boxstyle, frame hidden}

\newenvironment{guide}
{
    \begin{tcolorbox}[boxstyle, frame hidden]
}{
    \end{tcolorbox}
}

\setlist{leftmargin=20pt,itemsep=0pt,topsep=1ex}

\renewcommand*{\backref}[1]{}
\renewcommand*{\backrefalt}[4]{%
	\ifcase #1 Not cited.%
	\or        Cited on p.~#2.%
	\else      Cited on pp.~#2.%
    \fi}

\fancypagestyle{mainpage}
{%
    \fancyhf{}

    \fancyhead[CO] {\small \itshape \rightmark}
    \fancyhead[RO,LE] {\oldstylenums{\thepage}}
    \fancyhead[CE] {\textsc{\leftmark}}
    \setlength{\headheight}{14pt}
}

\fancypagestyle{firstpage}
{%
    \fancyhf{}
    
    \fancyfoot[C] {\oldstylenums{\thepage}}
}

\fancypagestyle{plain}
{%
    \fancyhf{}
    
    \fancyfoot[C] {\oldstylenums{\thepage}}
}

\titleformat{\chapter}
 {\large \scshape \filcenter}{\thechapter.}{1em}{}
\titlespacing*{\chapter}
{0pt}{40pt}{20pt}

\titleformat{\section}
 {\normalsize \itshape \filcenter}{\thesection.}{1em}{}
\titlespacing*{\section}
{0pt}{4.5ex plus 1ex minus .2ex}{2.5ex plus .2ex}

\renewcommand{\cftsecpagefont}{\mdseries}

\makeatletter \renewcommand{\cftsecfillnum}[1]{%
  {\cftsecleader}\nobreak
  \makebox[\@pnumwidth][\cftpnumalign]{\cftsecpagefont \oldstylenums{#1}}\cftsecafterpnum\par
} \makeatother
\makeatletter \renewcommand{\cftchapfillnum}[1]{%
  {\cftchapleader}\nobreak
  \makebox[\@pnumwidth][\cftpnumalign]{}\cftchapafterpnum\par
} \makeatother

%% file: introduction.tex
\chapter*{Introduction}
\thispagestyle{firstpage}
\addcontentsline{toc}{section}{Introduction}

As higher categories become more pervasive in mathematics, physics, and computer science, so does the use of \emph{diagrammatic reasoning} with its unique blend of the topological and the algebraic.
However, beyond the reasonably well-understood case of dimension 2, this proceeds with uncertain steps.
The interplay of soundness of diagrammatic proofs, coherence of higher-algebraic structures, and strictification results is often misunderstood, with the consequence that diagrammatic reasoning in higher dimensions, where strictification results as strong as in dimension 2 are not available, is avoided or approached tentatively.
Furthermore, the classical literature on the topic, while groundbreaking at its time, was notoriously error-prone and is now somewhat outdated, yet no comprehensive modern reference has emerged to supplant it. 
The consequence is that results that could be derived from general statements are commonly re-proven in special cases, and those proofs remain scattered across the technical sections of many different papers.

This book is an attempt to provide a systematic treatise and build a modern reference on higher-categorical diagrams, including pasting diagrams.
It stems from my attempt, over several years, to develop a ``convenient'' framework for higher-dimensional diagram rewriting, where ``convenience'' should be evaluated on multiple levels:
expressiveness for real use cases of higher-dimensional rewriting, such as presentations of higher-dimensional algebraic theories and their rewriting-theoretic coherence proofs;
connection to practically used models of higher categories, monoidal categories, and other higher structures;
soundness for homotopical algebra;
and feasibility of implementation in proof assistants or other computational aids.

All of these aspects, which are prominent in my current and past research activities, have influenced the perspective and the choices made in the writing of this book, yet are not, for the most part, explicitly featured in it.
This is not, indeed, a research monograph, although it does contain original results: it is a reference book, revisiting all the classical topics from a better vantage point.
Its aim is not to break new grounds, but to set a trustworthy yet flexible foundation for future developments.

As was realised in the early days of higher category theory, the theory of higher-categorical diagrams is a form of \emph{directed combinatorial topology}: the ``shape'' of an $n$\nbd dimensional diagram is, in its essence, the data of an $n$\nbd dimensional cell complex, together with an orientation on its cells of every dimension which subdivides their boundary into two ``halves'' --- an \emph{input}, or \emph{source} half, and an \emph{output}, or \emph{target} half --- in such a way that both halves are \emph{composable} arrangements of cells for the higher-categorical structure in which the diagram is meant to be interpreted.
The characterisation of such composable arrangements is traditionally called the \emph{pasting problem}, and the result that a certain class of composable arrangements admits a composite --- usually, a \emph{unique} composite, in a suitable sense --- is called a \emph{pasting theorem}.

The early days of higher-categorical diagrams were focussed on the pasting problem for strict $\omega$\nbd categories.
While this model remains actively studied, and has had a recent resurgence driven by the \emph{polygraph} approach to higher-dimensional rewriting, it is fair to say that it is less prominent than weaker models which are in general sound for homotopical algebra.\footnote{
Something that is known to be false of the strict model \cite[Theorem 4.4.2]{simpson2009homotopy}.}
The good news is that the study of higher-categorical diagrams can safely be decoupled from any particular model of higher categories, simply by adopting the ``categorical point of view'': instead of \emph{identifying} diagrams with special instances of some other structure, we can study them independently, then consider functorial interpretations in various contexts.
In this way, the exact same class of diagrams can in principle be used to reason in and about strict higher-dimensional categories, weak higher-dimensional categories, or topological spaces.

The systematic adoption of the ``categorical point of view'', hence the study not only of diagram shapes but of their \emph{morphisms} as well, is the single main technical innovation in this book, as this point of view is conspicuously absent in all the classical references.\footnote{
	An exception that will be extensively discussed is \cite{steiner2004omega}, which does consider morphisms.
However, since the purpose of this work, for technical reasons, is to describe a \emph{full} subcategory of the category of strict $\omega$\nbd categories, those morphisms are affected by the same ``topological unsoundness'' problem that affects strict $\omega$\nbd categories.}
This simple shift in perspective naturally overcomes some of the issues and restrictions of earlier approaches.

For example, most of these approaches impose some form of \emph{acyclicity} condition on the shapes of diagrams.
This has the effect of forbidding very simple non-composable shapes already in dimension 1, and commonly occurring composable shapes starting from dimension 3, such as
\[\begin{tikzcd}
	\bullet && \bullet
	\arrow[curve={height=-18pt}, from=1-1, to=1-3]
	\arrow[curve={height=-18pt}, from=1-3, to=1-1]
\end{tikzcd}\quad \text{and} \quad
\begin{tikzcd}[sep=tiny]
	&& {\bullet} &&&& {\bullet} \\
	{\bullet} &&& {\bullet} && {\bullet} &&& {\bullet} \\
	& {\bullet} &&&&&& {\bullet}
	\arrow[Rightarrow, from=2-4, to=2-6]
	\arrow[curve={height=6pt}, from=2-1, to=3-2]
	\arrow[""{name=0, anchor=center, inner sep=0}, curve={height=12pt}, from=3-2, to=2-4]
	\arrow[""{name=1, anchor=center, inner sep=0}, curve={height=12pt}, from=2-6, to=3-8]
	\arrow[curve={height=6pt}, from=3-8, to=2-9]
	\arrow[""{name=2, anchor=center, inner sep=0}, curve={height=-12pt}, from=1-7, to=2-9]
	\arrow[curve={height=-6pt}, from=2-6, to=1-7]
	\arrow[""{name=3, anchor=center, inner sep=0}, curve={height=-12pt}, from=2-1, to=1-3]
	\arrow[curve={height=-6pt}, from=1-3, to=2-4]
	\arrow[from=3-2, to=1-3]
	\arrow[from=1-7, to=3-8]
	\arrow[curve={height=-6pt}, shorten >=7pt, Rightarrow, from=3-2, to=3]
	\arrow[curve={height=6pt}, shorten <=7pt, Rightarrow, from=0, to=1-3]
	\arrow[curve={height=-6pt}, shorten <=7pt, Rightarrow, from=1, to=1-7]
	\arrow[curve={height=6pt}, shorten >=7pt, Rightarrow, from=3-8, to=2]
\end{tikzcd} \]
respectively.
Moreover, acyclic shapes tend to not be stable under various useful constructions: typically, stronger acyclicity properties are not stable under arbitrary duals, and weaker acyclicity properties are not stable under pasting or under Gray products.

The restriction to acyclic shapes turns out to be an artefact: since these approaches lack good notions of morphisms of diagram shapes, they only consider ``subshapes'', and need these to form a strict $\omega$\nbd category, which is only guaranteed by acyclicity.\footnote{
Indeed, the only exception in \emph{not} requiring acyclicity in general, which is \cite{steiner1993algebra}, obtains only a ``partial'' $\omega$\nbd category of subshapes, which is then algebraically freely extended; but this defeats in part the purpose of a combinatorial description.}
By considering more general morphisms, we can dispose of the acyclicity condition: for comparison, the \emph{linear subgraphs} of a directed graph only form a category when the directed graph is acyclic, but general \emph{paths} in a directed graph always form a category.
The morphism-based approach also restores the status of \emph{pasting} as a universal construction, that is, as a certain pushout of inclusions of diagram shapes, in such a way that ``pasting together diagram shapes'' satisfies the equations of strict $\omega$\nbd categories up to unique isomorphism, simply as a consequence of the fact that different sequences of pastings compute colimit cones under the same functor.

It is common knowledge among categorically-minded mathematicians that there is often a trade-off between \emph{nice objects} and \emph{nice categories}: the first form less nice categories, the latter include some less nice objects.
Having established that acyclic shapes are at the \emph{nice objects, bad category} end of the spectrum, the question is what class of diagram shapes strikes the best balance.
Perhaps unsurprisingly in hindsight, I ended up where many combinatorial topologists end up, looking towards \emph{regular cell complexes}.

A regular cell complex is a cell complex whose every closed cell is embedded homeomorphically into the complex, so, in particular, the image of each closed cell remains a closed topological ball.
A classical result establishes that regular cell complexes are ``combinatorial'' in the sense that they can be reconstructed up to homeomorphism from their \emph{face poset}, that is, the data of what cells lie in the boundary of each cell.
Most significantly, regular cell complexes are closed under a number of constructions --- products, joins, suspensions, gluings --- and the results can be readily computed at the level of face posets.

The theory of \emph{regular directed complexes} that I present in the book can be seen as the directed version of the theory of regular cell complexes, where face posets are supplemented with appropriate orientation data.
Since the boundary of a ``directed cell'' is subdivided into two halves, we give ``regularity'' a stricter interpretation in the directed world, requiring not only that the entire boundary of an $(n+1)$\nbd cell be an $n$\nbd sphere, but that each of the two halves be a closed $n$\nbd ball; all of these properties can be achieved by purely combinatorial means.
Remarkably, the good stability properties of regular cell complexes translate to good stability properties of regular directed complexes under the $n$\nbd categorical counterparts of topological constructions, such as \emph{Gray products} as a counterpart to cartesian products.

If we try to identify the most useful notions of morphism of regular directed complexes, it turns out that there are two natural choices, dual to each other, characterised by the fact that they induce strict functors of strict $\omega$\nbd categories \emph{covariantly} or \emph{contravariantly}.
We call these \emph{maps} and \emph{comaps}, respectively.
Interestingly, both have underlying order-preserving maps of posets, so both induce continuous maps in the geometric realisation of a regular directed complex via the usual order complex construction.
Maps admit a ternary factorisation system, separating them further into three classes: \emph{final maps} which can only ``collapse'' cells; \emph{surjective local embeddings} which can only rigidly identify sets of cells; and \emph{inclusions} which embed a complex into a larger complex.
Comaps, on the other hand, are dual to \emph{subdivisions}: they can merge groups of cells into a single cell, preserving the overall boundary.

This refined study of morphisms between diagram shapes is entirely original to this book, and I am confident that it will play a role in the quest for good algebraic models of non-strict higher categories. 
Indeed, when restricted to shapes of pasting diagrams, these classes of morphisms separately embody different operations that we expect in such models: inclusions are dual to faces; final maps are dual to units or degeneracies; and comaps are a restricted form of composition.
A careful study may reveal how these can interact while avoiding the pitfall of ``strict Eckmann--Hilton''.

The definition of regular directed complex, which determines the shapes of non-necessarily composable diagrams, depends on the definition of \emph{molecule}, which determines the shapes of composable diagrams.
Another point of departure with most earlier approaches\footnote{
	But not \cite{steiner1993algebra}, the main influence of this book, and its follow-ups.}
	is that molecules are defined \emph{inductively}, as produced by a number of ``constructors'', rather than by \emph{characterisation}, that is, by listing axioms which, when satisfied, guarantee that a certain face-poset-like structure determines a well-formed pasting diagram.
In other words, we take a \emph{synthetic} rather than an \emph{analytic} approach.
This choice is the obvious one if one is interested in higher-dimensional diagram rewriting as a computational tool, since the definition of well-formed cell shapes in terms of constructors translates smoothly into a definition of constructors for a higher-dimensional rewrite system.
In fact, the development of the theory of molecules has proceeded alongside their implementation as data structures in \texttt{rewalt}, a Python library for computational higher-dimensional diagram rewriting.\footnote{
See \url{http://rewalt.readthedocs.io}. All Hasse diagrams and string diagrams in the book were generated with \texttt{rewalt}, which supports TikZ output.}

Another advantage of the inductive approach is that many restricted classes of shapes that have been usefully studied in their own right --- globes, oriented simplices, oriented cubes, positive opetopes, and the objects of the $\thetacat$ category which I will refer to as \emph{thetas} --- and whose specific combinatorics, as needed for the purpose of characterisation, can vary wildly between each other, all turn out to admit very simple definitions as inductive subclasses of molecules.
This allows us to take a uniform approach to these classes of shapes, repurposing results proven for the class of all molecules, and obtain concise proofs of facts that would otherwise need a copious amount of \emph{ad hoc} combinatorics.

Indeed, what I hope will be one of the main contributions of this book is simply distilling certain recurring ideas into a small number of powerful technical lemmas and proof methods, stated at the right level of generality, and applicable to a wide assortment of problems related to higher-categorical diagrams.
I will single out three in particular, whose recurrence should become apparent to the reader as they go through the book.
\begin{enumerate}
	\item The first is the notion of \emph{layering} --- a pasting decomposition into ``layers'', each containing a single cell whose dimension is larger than the pasting dimension --- and its associated proof method, \emph{induction on layering dimension}.
		The idea of layering is certainly present, and used to great effect, in much of the earlier literature; but its in-depth study, in relation to the notion of \emph{ordering} of a diagram shape, is original.
	\item The second is the notion of \emph{oriented thinness}, which is a combinatorial property enjoyed by regular directed complexes, that on its own guarantees that ``cellular chains'' canonically form a chain complex.
		This is a very simple property, which I singled out in some of my earlier work, and turns out to be surprisingly powerful in proving statements of a local character. 

	\item The third is the notion of \emph{generalised pasting}, which gives conditions under which molecules can be glued together along \emph{portions} of their input or output boundaries, producing another molecule.
		To my knowledge, this is entirely original.
		Its main use is producing what, in all fairness, is the first \emph{readable} proof that pasting diagram shapes are closed under Gray products, but it turns out to be very useful beyond that, especially in its restricted form of \emph{pasting at a submolecule}.
\end{enumerate}


\section*{What is not in the book}

Before entering into more detail about the content and structure of the book, I will spend a few paragraphs on what is \emph{not} in it, to avoid later disappointing the reader who is searching for their favourite topic.
The following is a non-exhaustive list of topics that one may, within reason, be hoping to find, but will not find in this book.

\begin{itemize}
	\item \emph{A detailed comparison with other formalisms for pasting diagrams.}

		This is meant as a \emph{reference book}: a source of useful definitions, methods, and results for the researcher who is using higher-categorical diagrams in their work.
		As a consequence, I am only interested in using the most convenient framework for achieving these results.
		This framework is the result of many adjustments made over several years and is, in my opinion and to the best of my knowledge, the best available.
		Unlike \cite{forest2022unifying}, this is not meant as a \emph{survey}, so except for a few pointers I make no attempt to formally compare regular directed complexes to parity complexes, pasting schemes, or any other class of structures used to describe diagram shapes.
		
		A partial exception is that I go into detail about the connection to Steiner's theory of \emph{augmented directed chain complexes}.
		This is for several reasons.
		First of all, Steiner theory is not strictly about \emph{diagrams}, and has a sufficiently unique character to be considered of independent interest.
		Secondly, the most popular definition of Gray products of strict $\omega$\nbd categories, and the only available definition of joins of strict $\omega$\nbd categories, both go through Steiner theory, so addressing it is necessary in order to compare these operations on regular directed complexes with those on strict $\omega$\nbd categories.
		Finally, the use of Steiner theory has been particularly popular recently in the theory of higher categories, so it seems indeed useful to detail to what extent the theory presented in this book subsumes it.

	\item \emph{Proofs in higher category theory, combinatorial topology, commutative algebra, etc that are not strictly related to the combinatorics of diagrams.}

		I follow a double standard for proofs.
		For anything that concerns regular directed complexes, or any of the surrounding combinatorial structures, I attempted to be detailed and careful \emph{in excess} of the expected standard for research-level mathematics.
		In general, I have tried to skip a step or claim it is ``straightforward'' only when it \emph{honestly} is straightforward, at the level of ``spelling out two sides of an equation and verifying that they are equal''.
		This is because the field of higher-categorical combinatorics is notoriously treacherous, tied historically to some famous retractions and corrections, and it is my wish to make this book as trustworthy as possible.

		On the other hand, on anything which is \emph{not} strictly related to the combinatorial framework presented in the book, I have relied systematically on citations to other sources: the theory touches on several other fields, and it would not be reasonable to redevelop all the background from scratch.
		In many instances, the cited sources will differ to some extent in their setup and notation: for instance, the book uses the ``single-set'' definition of strict $\omega$\nbd category, unlike the sources cited for some related results.
		In all these cases, I trust the reader to make the necessary adjustments.

	\item \emph{Any specific development of diagrams in multiple categories, higher operads, or other higher structures.}

		The theory developed in the book is only about diagrams in higher-dimensional categories.
		On the other hand, it is my experience that the diagram shapes that appear in other higher structures, as well as the ``composable'' subclasses, are typically \emph{restrictions} of those that appear in higher-dimensional categories.
		For example, the shapes of diagrams in multiple (double, triple, $\ldots$) categories are typically oriented cubical shapes, which are also valid in higher-dimensional categories; the difference appears at the level of the \emph{labelling}, where different faces of cubical in multiple categories are labelled in different sets of cells (horizontal, vertical, $\ldots$) rather than the same set.
		Similarly, the shapes of diagrams in (non-symmetric) operads are ``many-to-one'', \emph{opetopic} shapes, which are restrictions of the ``many-to-many'' shapes appearing in higher-dimensional categories.
		In a somewhat different case, the diagram shapes that can be interpreted in planar \emph{polycategories} are the same as those that can be interpreted in monoidal categories, yet the class of \emph{composable} shapes is restricted to those whose ``graph of connections'' is acyclic.
		Thus there is no reason why the theory presented here cannot be repurposed for one's favourite higher structure.

	\item \emph{Any theory specifically about string diagrams, or their higher-dimensional generalisations, such as ``manifold diagrams''.}

		The main perspective adopted on diagrams in this book is the ``pasting diagram'' perspective, where a $k$\nbd dimensional cell in an $n$\nbd category is, indeed, pictured as $k$\nbd dimensional topological cell.
		To a certain extent, the ``string diagram'' perspective is simply the Poincar\'e-dual one, where a $k$\nbd dimensional cell in an $n$\nbd category is pictured as an $(n-k)$\nbd dimensional topological cell.
		As far as this goes, one can definitely read the book with string diagrams in their mind, and in fact we do \emph{make ample use} of string diagrams: the ability to exchange dimension and codimension is precious for understanding the top-dimensional structure of diagrams in dimension higher than 3, where spatial intuition falters.

		On the other hand, in the same way as the connection to cell complexes or the homotopy hypothesis are somewhat more natural to the ``pasting diagram'' perspective, there are many aspects that are specific to the ``string diagram'' perspective, especially when it comes to higher categories with dualisable cells, and their connection to the \emph{tangle} and \emph{cobordism} hypotheses.
		None of these are addressed here.
		Fortunately, the recent book \cite{dorn2021framed} develops a form of ``directed combinatorial topology'' that is tailored to this perspective, and makes a perfect companion to this book.
		How exactly to make the two combinatorial frameworks converge is an interesting question for future work.

	\item \emph{Any specific discussion of diagrams in weak $(\infty, n)$\nbd categories.}

		This is a very active area of research, and one in which I am particularly interested \cite{hadzihasanovic2020diagrammatic}.
		In particular, a number of recent works have been devoted to the question of whether certain strict $n$\nbd categories that are presented by combinatorial diagram shapes, and are ``freely generated'' in the sense of polygraphs, are still ``free'' in the suitable homotopically-coherent sense when seen as weak $(\infty, n)$\nbd categories \cite{maehara2023orientals} \cite{gagna2023nerves} \cite{campion2023infty}.

		In this case, the reason why no such material is included in the book is not because it would not be a good thematic fit.
		Rather, I believe that it is too recent to be included in a reference book.
		Perhaps it will be part of a future edition.

	\item \emph{A discussion of computational aspects of diagram rewriting.}

		This is an increasingly popular topic \cite{corbyn2024homotopy} on which I am personally active \cite{hadzihasanovic2023higher}, and, again, it is not covered because the research is too recent and too active.
\end{itemize}


\section*{Structure of the book}

The book consists of 11 chapters, which are further subdivided into a total of 39 sections.
Each chapter opens with an extensive introduction, which outlines its content and contextualises it within the book and in conversation with the earlier literature.
Each section, in turn, has a shorter introduction providing more specific guidance.
The main body of the book is followed by a bibliography, an index of terms and notations, and an index of counterexamples.
The 11 chapters can be roughly subdivided into three groups:
\begin{itemize} 
	\item \emph{foundational chapters} (Chapter \ref{chap:order}, Chapter \ref{chap:ogposets}, Chapter \ref{chap:molecules}, Chapter \ref{chap:omegacats}, and Chapter \ref{chap:maps}),

	\item \emph{technical chapters} (Chapter \ref{chap:layerings}, Chapter \ref{chap:constructions}, and Chapter \ref{chap:acyclic}),

	\item \emph{special interest chapters} (Chapter \ref{chap:special}, Chapter \ref{chap:geometric}, and Chapter \ref{chap:steiner}).
\end{itemize}
The foundational chapters are where we set up the framework.
They are, in a certain sense, \emph{definition-centric}: while several useful results appear here, they tend to be aimed towards supporting or justifying a definition.
For example, Theorem \ref{thm:morphisms_of_atoms_are_injective} is a highly non-trivial result whose proof relies on machinery from all previous chapters, but its main \emph{purpose} is justifying the definitions of \emph{combinatorial diagram} and of \emph{map} of regular directed complexes.
Readers may want to skim through these chapters on a first reading, focussing on the definitions, the examples, and the narrative, in order to get an intuitive sense of how the theory works, and how it can be used for their own purposes.

Specifically, Chapter \ref{chap:order} presents the order-theoretic notions that will be used in the rest of a book (except for a few complements that are postponed to Section \ref{sec:pushforwards}).
Chapter \ref{chap:ogposets} introduces \emph{oriented graded posets}, the underlying combinatorial structure that we will use to model diagram shapes, as well as their relation to chain complexes.
Chapter \ref{chap:molecules} introduces \emph{molecules}, the subclass of oriented graded posets which model pasting diagram shapes, that is, those shapes that are composable in the algebra of strict $\omega$\nbd categories, and that will form the building blocks of all other diagram shapes.
Chapter \ref{chap:omegacats} constructs strict $\omega$\nbd categories out of oriented graded posets.
Then, it defines \emph{regular directed complexes}, which model general diagram shapes and are arguably the main character of the book.
Chapter \ref{chap:maps} defines and studies the two natural notions of morphisms of regular directed complexes, \emph{maps} and \emph{comaps}.

By contrast, the technical chapters are \emph{theorem-centric}: their purpose is to prove non-trivial results about higher-categorical diagrams and devise methods applicable both to further theory-building, or to direct practice.
These chapters are probably most interesting to readers who are actively working on higher category theory, higher-dimensional rewriting, or higher algebra, and are looking for specific technical facts that will help them in their research. 

Chapter \ref{chap:layerings} is centred on the problem of recognising what portions of a diagram are \emph{rewritable}, in the sense that they can be substituted with another diagram with the same boundary, producing another well-formed diagram.
This is of course a foundational problem for higher-dimensional rewriting, but it also shows up in a number of other technical questions.
Chapter \ref{chap:constructions} focusses on constructions and operations on oriented graded posets, and whether they preserve the classes of molecules and regular directed complexes.
These include \emph{generalised pastings}, \emph{Gray products}, \emph{suspensions}, \emph{joins}, and direction-reversing \emph{duals}.
Chapter \ref{chap:acyclic} is centred on the problem of what oriented graded posets present strict $\omega$\nbd categories that are \emph{freely generated} in the sense of polygraphs.
It provides a most general, technical answer, then proceeds to consider some simpler sufficient conditions, studying their further consequences and stability under various operations.
Some of the results of this chapter were produced in collaboration with Diana Kessler.

Foundational and technical chapters form the core of the book.
Special interest chapters are a sort of coda: they could be characterised as \emph{application-centric}.
Each of them explores a particular facet of the theory of higher-categorical diagrams, technically and conceptually independent of the others.

Chapter \ref{chap:special} focusses on some special classes of diagram shapes that have been considered as the underlying combinatorics of models of $(\infty, n)$\nbd categories: \emph{globes}, \emph{thetas}, \emph{oriented simplices}, \emph{oriented cubes}, and \emph{positive opetopes}.
Chapter \ref{chap:geometric} looks at connections with combinatorial topology, and contains the proof that regular directed complexes are combinatorial presentations of regular CW complexes.
Finally, Chapter \ref{chap:steiner} examines the connection between the theory of regular directed complexes and Steiner's theory of augmented directed chain complexes, on which the now-standard definition of the Gray product of strict $\omega$\nbd categories is based.

This completes a very brief overview of the content; I invite you to read the chapter-specific introductions for more detail.


\section*{Main sources of the book}

While there is much original content in the book, an equal amount is inspired, borrowed, or adapted from other sources.
Since even the proofs that are more closely based on these sources have, at the very least, been rethought, I did not in general give precise attributions in the text, except for those results whose proof is \emph{only} referenced.
Here I would like to give an overview of the main sources that I have used, and how their content relates to the content of the book.
The order is chronological.

\begin{itemize}
	\item \emph{The algebra of oriented simplexes} \cite{street1987algebra}.
		I feel compelled to include this because, even though I may not have used its technical content directly, this is the paper that opened up the entire field of higher-categorical combinatorics, and its influence is all-encompassing.

	\item \emph{An $n$\nbd categorical pasting theorem} \cite{power1991pasting}.
		The original pasting theorem used a definition of diagram shapes that is less combinatorial, more concretely topological, involving embeddings into $n$\nbd dimensional space.
		It has an important influence on Section \ref{sec:in_low_dim}: its \emph{domain replacement condition} corresponds, in our framework, to the admissibility of \emph{multiple substitution} of submolecules, and Power's proof that this holds automatically in low dimension is a key insight for our proof of Theorem \ref{thm:dim3_frame_acyclic}.

	\item \emph{$\infty$\nbd groupoids and homotopy types} \cite{kapranov1991infty}.
		This is the notorious paper that mistakenly claimed a proof of the homotopy hypothesis for strict $\omega$\nbd categories with weak inverses.
		However, it is also the first, and for a long time the \emph{only} article to have thought of using one of the combinatorial frameworks for pasting diagrams to define a ``shape category'' for higher structures, hence also one of the only ones that considered interesting morphisms of diagram shapes, analogous to our ``maps''.
		The authors' cavalier attitude towards proofs meant that they did not really have much control of either the shapes (they used acyclic shapes, namely, Johnson's composable pasting schemes \cite{johnson1989combinatorics}, without realising that they are not closed under all the operations they needed, as observed in \cite[Discussion A.2]{henry2019non}) or their morphisms (they assumed that the ``co-degeneracy'' maps between shapes could reproduce every identity in their realisations as strict $n$\nbd categories, which they crucially cannot \cite[\S 8.55]{hadzihasanovic2020diagrammatic}).
		Nevertheless, with all its flaws, the article remains ahead of its time and a treasure trove of ideas.

	\item \emph{The algebra of directed complexes} \cite{steiner1993algebra}.
		This is, by quite a distance, the single most important influence of this book. 
		In fact, the entire book can be fairly characterised as an expansion and commentary on the results and ideas planted in this article.
		From the definition of \emph{molecule}, to the notion of \emph{frame dimension}, to the idea of \emph{splitness} which here has become \emph{frame-acyclicity}, and the recognition that it, and not acyclicity, is key to algebraic freeness of the presented strict $\omega$\nbd categories, the influence of this article is everywhere in the book.
		It is, in my opinion, the most far-sighted article ever written on higher-categorical diagrams, and I hope that its insights, which have been somewhat overshadowed by Steiner's later work on directed chain complexes, will be recognised over time.

	\item \emph{Higher-dimensional word problems with applications to equational logic} \cite{burroni1993higher}.
		While the connection between higher-categorical diagrams and cell complexes was already understood, this is the article the brought \emph{rewriting} into the mix as an additional facet of a topological-categorical-computational triality.
		Presumably, without this article and the school that stemmed from it, rewriting-theoretic ideas such as \emph{substitution} and \emph{rewritable submolecules} --- which, beyond their immediate applicability, turn out to be essential technical tools for studying higher-dimensional diagrams --- would not have been so obvious or prominent.

	\item \emph{Pasting in multiple categories} \cite{steiner1998pasting}.
		While not quite as influential as Steiner's earlier article, this is where the notion of \emph{round} molecule comes from, as well as the realisation that roundness leads to good topological properties.

	\item \emph{Omega-categories and chain complexes} \cite{steiner2004omega}.
		This is the original ``Steiner theory'' article and is still the main source on directed chain complexes.

	\item \emph{On positive opetopes, positive opetopic cardinals and positive opetopic sets} \cite{zawadowski2007positive}.
		In the world of opetopes, \emph{regularity} becomes \emph{positivity}, that is, the property of each cell having at least one input face.
		By recognising the good combinatorial-topological properties of this restriction, and then also by considering collapsing maps between positive opetopes as a means to obtain degeneracies on opetopic sets without ``nullary representability'', Zawadowski's work has been an important influence on this book's framework.

	\item {\setlength{\spaceskip}{0.2em plus 0.05em minus 0.0667em}
	\emph{Joint et tranches pour les $\infty$\nbd cat{\'e}gories strictes} \cite{ara2020joint}}.
		This work, which has revitalised research on and around Steiner theory by mending some gaps and expanding its technical reach, has also been my main reference on $\omega$\nbd categorical constructions such as joins and suspensions.

	\item \emph{Type theoretical approaches to opetopes} \cite{thanh2022type}.
		This is the main source of the approach to opetopes used in Section \ref{sec:opetopes}.
\end{itemize}
I have previously presented parts of the framework and results of the book in \cite{hadzihasanovic2020combinatorial}, \cite{hadzihasanovic2020diagrammatic}, \cite{hadzihasanovic2021smash}, 
\cite{hadzihasanovic2023data}, \cite{hadzihasanovic2023higher}.


\section*{How to cite this book}

It is very likely that the book will be updated more than once with corrections and new material.
Please refer to the exact version when citing.

\clearpage
\thispagestyle{empty}

\chapter*{Acknowledgements}
\thispagestyle{firstpage}
\addcontentsline{toc}{section}{Acknowledgements}

First of all, I would like to thank Cl\'emence Chanavat and Diana Kessler for their feedback and support at various stages during writing.

My understanding of higher-categorical combinatorics, as well as other categorical, combinatorial, topological, and algebraic aspects that played a role in it, has been influenced, over the years, by conversations with many people.
Only naming ones that I can recount, I would like to thank
Samson Abramsky,
Dimitri Ara,
Cameron Calk,
J.~Scott Carter,
Bob Coecke,
Pierre-Louis Curien,
Antonin Delpeuch,
Christoph Dorn,
Eric Finster,
Simon Forest,
Harry Gindi,
Stefano Gogioso,
Eric Goubault,
Yves Guiraud,
Philip Hackney,
Masahito Hasegawa,
Simon Henry,
Alex Kavvos,
Aleks Kissinger,
Joachim Kock,
Guillaume Laplante-Anfossi,
Louise Leclerc,
Chaitanya Leena Subramaniam,
Fosco Loregian,
F\'elix Loubaton,
Georges Maltsiniotis,
Ioannis Markakis,
Dan Marsden,
Paul-Andr\'e Melli\`es,
Fran\c{c}ois M\'etayer,
Samuel Mimram,
Viktoriya Ozornova,
Robin Piedeleu,
David Reutter,
Morgan Rogers,
Mario Rom\'an,
Martina Rovelli,
Chiara Sarti,
Pawe\l{} Soboci\'nski,
Georg Struth,
Jamie Vicary, and
Noam Zeilberger
for these conversations.

During the writing of this book, I was supported by the ESF funded Estonian IT Academy research measure (project 2014-2020.4.05.19-0001) and by the Estonian Research Council grant PSG764.

This book and the entire field of higher-categorical combinatorics owe very much to the work of Marek Zawadowski, who passed away suddenly shortly before my work was completed.
I only met him once, at a CIRM meeting in 2017, right after finishing my doctorate, and he left a lovely impression on me; those who knew him better speak highly of him both on a personal and on an intellectual level.
In one of our conversations, he expressed skepticism that ``many-to-many'' shapes could be a suitable foundation for higher category theory, given how complicated ``many-to-one'' shapes already were.
I would be curious to hear his opinion on this book.

\clearpage

%% file: order.tex
\chapter{Elements of order theory} \label{chap:order}
\thispagestyle{firstpage}

\begin{guide}
Consider the following pasting diagram, which may be interpreted in a 2\nbd category such as the 2\nbd category of small categories, functors, and natural transformations.
\begin{equation} \label{eq:example_diagram}
\begin{tikzcd}[sep=small]
	{{\scriptstyle x}\;\bullet} && {{\scriptstyle y}\;\bullet} && {{\scriptstyle z}\;\bullet} \\
	& {{\scriptstyle y}\;\bullet}
	\arrow[""{name=0, anchor=center, inner sep=0}, "f", curve={height=-18pt}, from=1-1, to=1-3]
	\arrow["g", from=1-3, to=1-5]
	\arrow["f"', curve={height=6pt}, from=1-1, to=2-2]
	\arrow["t"', curve={height=6pt}, from=2-2, to=1-3]
	\arrow["\alpha", shorten <=3pt, shorten >=6pt, Rightarrow, from=2-2, to=0]
\end{tikzcd}
\end{equation}
We are interested in encoding the information contained in this pasting diagram into a combinatorial structure.

To start with, we can separate the \cemph{shape} of the pasting diagrams, consisting of all the different elements in the picture --- 0\nbd cells (points), 1\nbd cells (arrows), 2\nbd cells (oriented disks) --- from its interpretation, as expressed by their labelling.
For this purpose, we assign distinct numbers to all cells of the same dimension.
\begin{equation} \label{eq:example_shape}
\begin{tikzcd}[sep=small]
	{{\scriptstyle 0}\;\bullet} && {{\scriptstyle 2}\;\bullet} && {{\scriptstyle 3}\;\bullet} \\
	& {{\scriptstyle 1}\;\bullet}
	\arrow[""{name=0, anchor=center, inner sep=0}, "3", curve={height=-18pt}, from=1-1, to=1-3]
	\arrow["2", from=1-3, to=1-5]
	\arrow["0"', curve={height=6pt}, from=1-1, to=2-2]
	\arrow["1"', curve={height=6pt}, from=2-2, to=1-3]
	\arrow["0", shorten <=3pt, shorten >=6pt, Rightarrow, from=2-2, to=0]
\end{tikzcd}
\end{equation}
We may then use the pair $(n, k)$ to refer to the $n$\nbd cell to which we assigned number $k$.
We recover (\ref{eq:example_diagram}) by supplementing (\ref{eq:example_shape}) with the labelling function
\begin{align*}
	(0, 0) & \mapsto x, \quad 
	       & (0, 1), (0, 2) & \mapsto y, \quad 
	       & (0, 3) & \mapsto z, \\
	(1, 0), (1, 3) & \mapsto f, \quad 
		       & (1, 1) & \mapsto t, \quad 
		       & (1, 2) & \mapsto g, \\
	(2, 0) & \mapsto \alpha. && &&
\end{align*}
We have now reduced the problem of encoding a pasting diagram into the problem of encoding its shape.

If we forget about the direction of cells, what (\ref{eq:example_shape}) is describing is a finite 2\nbd dimensional cell complex --- in fact, a simplicial complex --- homeomorphic to the wedge of a disk (closed 2\nbd ball) and an interval (closed 1\nbd ball).
This cell complex is, moreover, of a special kind: it is \cemph{regular}, that is, the inclusion of each generating closed cell is a homeomorphism with its image.
A classical result of combinatorial topology (Proposition 
\ref{prop:face_poset_has_an_inverse_up_to_isomorphism}) states that a regular cell complex can be uniquely reconstructed, up to homeomorphism, from its \cemph{face poset}: the poset whose
\begin{itemize}
	\item elements are the generating cells,
	\item given two distinct generating cells $x, y$, we have $x < y$ if and only if $x$ lies in the boundary of $y$.
\end{itemize}
The face poset of (\ref{eq:example_shape}) is
\[
	\input{img/face_poset.tex}
\]
where we have edges of the Hasse diagram go in the direction of the \emph{covering} relation, opposite to the order.
Notice that this poset is \emph{graded} by the dimension of the cells, that is, all cells of dimension $n$ have the same ``height'' $n$ in the Hasse diagram.
This is, in fact, the case for all face posets of regular cell complexes.

From this poset, we are able to reconstruct the shape of our pasting diagram, minus the direction of cells.
Indeed, in a higher-categorical pasting diagram, the boundary of an $n$\nbd cell is divided into an \cemph{input}, or \emph{source} half and an \cemph{output}, or \emph{target} half.
But (\ref{eq:example_shape}) satisfies a further regularity conditions: both the input and the output half of each $n$\nbd cell are, themselves, regular cell complexes homeomorphic to closed $(n-1)$\nbd balls.
In particular, they are the closures of their top-dimensional cells.

To specify the separation of the boundary of an $n$\nbd cell, thus, we just need to specify which $(n-1)$\nbd cells are \emph{input faces} and which are \emph{output faces}.
We can achieve that by labelling the edges of the face poset as either input or output edges.
Any set of two labels will work; we will use the set $\set{+, -}$, with $+$ for output and $-$ for input, but in Hasse diagrams we will portray the labelling by marking input edges in a different colour.

The \cemph{oriented face poset} of (\ref{eq:example_shape}) is then
\begin{equation} \label{eq:example_ogposet}
	\input{img/oriented_face_poset.tex}
\end{equation}
and in conjunction with the labelling function it is a complete specification of (\ref{eq:example_diagram}). \index{face poset!oriented}
We call this kind of combinatorial structure --- a finite graded poset, together with a 2\nbd valued edge-labelling of its Hasse diagram --- an \cemph{oriented graded poset}.

Since oriented graded posets are the fundamental structure with which we will work in the book, it is necessary to spend some time setting up their terminology and basic category theory.
In this chapter, we address the order-theoretic part.
Much of the content is standard or elementary, but may not be easily traced in the literature, especially for what concerns particular classes of morphisms.
\end{guide}


\section{Closed maps of posets} \label{sec:closed_maps}

\begin{guide}
In the face poset of a cell complex, an element $x$ by itself corresponds to an \emph{open} cell.
In order to obtain a \emph{closed} cell, one has to take the union of all the cells that lie in its boundary, that is, all $y$ such that $y \leq x$.
In order-theoretic terms, this is the \emph{lower set} of $x$.
Consequently, \cemph{closed subsets} which contain the lower set of each of their elements play an important role in our theory, as do \cemph{closed order-preserving maps}, with the property that both direct and inverse image preserve closed subsets.

In this section, we look at properties of the category $\posclos$ of posets and closed order-preserving maps, which is less studied than the category of all order-preserving maps.
Fortunately, \emph{colimits} are computed in the same way in the two categories, and that is what we will need the most.
\end{guide}

\begin{dfn}[Order induced on a subset] \index{poset!subsets!induced order}
	Let $U \subseteq P$ be a subset of a poset.
	The \emph{induced order} on $U$ is the restriction of the partial order on $P$ to $U$.
\end{dfn}

\begin{comm}
	We will always implicitly assume that a subset of a poset comes equipped with the induced order.
\end{comm}

\begin{dfn}[Closure of a subset] \index{poset!subsets!closure} \index{$\clos{U}$}
Let $P$ be a poset and $U \subseteq P$.
The \emph{closure of $U$} is the subset of $P$ 
\[
	\clos{U} \eqdef \set{ x \in P \mid \text{there exists $y \in U$ such that $x \leq y$} }.
\]
\end{dfn}

\begin{dfn}[Closed subset] \index{poset!subsets!closed}
Let $U$ be a subset of a poset.
We say that $U$ is \emph{closed} if $U = \clos{U}$.
\end{dfn}

\begin{lem} \label{lem:closure_is_monotonic}
Let $U, V$ be subsets of a poset.
If $U \subseteq V$ then $\clos{U} \subseteq \clos{V}$.
In particular, if $U \subseteq V$ and $V$ is closed then $\clos{U} \subseteq V$.
\end{lem}
\begin{proof}
Let $x \in \clos{U}$. 
Then there exists $y \in U$ such that $x \leq y$.
Since $U \subseteq V$, $y \in V$.
It follows that $x \in \clos{V}$.
\end{proof}

\begin{lem} \label{lem:closure_of_union}
Let $(U_i)_{i \in I}$ be a family of subsets of a poset.
Then 
\[\clos{\bigcup_{i \in I} U_i} = \bigcup_{i \in I} \clos{U_i}.\]
In particular, if all the $U_i$ are closed, so is their union.
\end{lem}
\begin{proof}
Straightforward.
\end{proof}

\begin{dfn}[The category $\poscat$] \index{$\poscat$}
We let $\poscat$ denote the category of posets and order-preserving maps.
\end{dfn}

\begin{lem} \label{lem:inverse_images_of_closed}
Let $P, Q$ be posets and $f\colon P \to Q$ a function of their underlying sets.
The following are equivalent:
	\begin{enumerate}[label=(\alph*)]
		\item $f$ is order-preserving;
		\item for all closed subsets $V \subseteq Q$, the inverse image $\invrs{f}(V) \subseteq Q$ is closed.
	\end{enumerate}
\end{lem}
\begin{proof}
	Suppose that $f$ is order-preserving, let $V \subseteq Q$ be closed, $x \in \invrs{f}(V)$, and $y \leq x$.
	Then $f(x) \in V$ and $f(y) \leq f(x)$, so $y \in \invrs{f}(V)$. 
	It follows that $\invrs{f}(V)$ is closed.

	Conversely, suppose $x \leq y$ in $P$.
	Then $y \in \invrs{f}(\clset{f(y)})$, which by assumption is closed, so $x \in \invrs{f}(\clset{f(y)})$, that is, $f(x) \in \clset{f(y)}$, so $f(x) \leq f(y)$.
\end{proof}

\begin{dfn}[Closed map] \index{map!closed}
	Let $f\colon P \to Q$ be an order-preserving map of posets.
	We say that $f$ is \emph{closed} if, for all closed subsets $U \subseteq P$, the direct image $f(U) \subseteq Q$ is closed.
\end{dfn}

\begin{lem} \label{lem:map_closed_iff_commutes_with_closure}
	Let $f\colon P \to Q$ be an order-preserving map of posets.
	The following are equivalent:
	\begin{enumerate}[label=(\alph*)]
		\item $f$ is closed;
		\item for all $U \subseteq P$, we have $\clos{f(U)} = f(\clos{U})$.
	\end{enumerate}
\end{lem}
\begin{proof}
	First of all, observe that $f(\clos{U}) \subseteq \clos{f(U)}$ always holds: if $x \in \clos{U}$, there exists $y \in U$ such that $x \leq y$.
	Then $f(x) \leq f(y) \in f(U)$.
	It follows that $f(x) \in \clos{f(U)}$.

	Suppose that $f$ is closed.
	Then $f(U) \subseteq f(\clos{U})$ and the latter is closed, so by Lemma
	\ref{lem:closure_is_monotonic}, we have $\clos{f(U)} \subseteq f(\clos{U})$.
	Conversely, assume $\clos{f(U)} = f(\clos{U})$ for all $U \subseteq P$.
	Then if $U$ is closed, $f(U) = f(\clos{U}) = \clos{f(U)}$, so $f(U)$ is closed.
\end{proof}

\begin{lem} \label{lem:map_closed_iff_maps_lowersets_to_lowersets}
	Let $f\colon P \to Q$ be an order-preserving map of posets.
	The following are equivalent:
	\begin{enumerate}[label=(\alph*)]
		\item $f$ is closed;
		\item for all $x \in P$, we have $\clset{f(x)} = f(\clset{x})$.
	\end{enumerate}
\end{lem}
\begin{proof}
	One implication is a special case of Lemma \ref{lem:map_closed_iff_commutes_with_closure}.
	For the other, let $U \subseteq P$.
	Then
	\[
		\clos{f(U)} = \clos{f\left(\bigcup_{x \in U} \set{x} \right)} =
		\bigcup_{x \in U} \clset{f(x)}
	\]
	using Lemma \ref{lem:closure_of_union} and the fact that direct images preserve all unions.
	By assumption, this is equal to
	\[
		\bigcup_{x \in U} f(\clset{x}) = f\left( \bigcup_{x \in U} \clset{x} \right) = f(\clset{U}),
	\]
	and we conclude by Lemma \ref{lem:map_closed_iff_commutes_with_closure}.
\end{proof}

\begin{dfn}[The category $\posclos$] \index{$\posclos$}
We let $\posclos$ denote the category whose objects are posets and morphisms are closed order-preserving maps.
\end{dfn}

\begin{rmk}
Evidently, $\posclos$ may be identified with a wide subcategory of $\poscat$.
\end{rmk}

\begin{dfn}[Closed embedding] \index{map!closed embedding}
	Let $f\colon P \to Q$ be an order-preserving map of posets.
	We say that $f$ is a \emph{closed embedding} if it is closed and injective.
\end{dfn}

\begin{rmk}
	Any closed subset $U \subseteq P$ determines a closed embedding $\imath\colon U \incl P$, where $U$ has the induced order.
	We will often not distinguish between a closed subset and the closed embedding that it determines.
\end{rmk}

\begin{lem} \label{lem:closed_embedding_is_order_reflecting}
	Let $\imath\colon P \incl Q$ be a closed embedding.
	Then $\imath$ is order-reflecting.
\end{lem}
\begin{proof}
	Let $x, y \in P$ and suppose $\imath(x) \leq \imath(y)$.
	Then $\imath(x) \in \clset{\imath(y)} = \imath(\clset{y})$, that is, there exists $x' \leq y$ such that $\imath(x') = \imath(x)$.
	Because $\imath$ is injective, $x = x'$, so $x \leq y$.
\end{proof}

\begin{lem} \label{lem:colimits_in_posclos}
	The category $\posclos$ has all small colimits, and they are preserved and reflected by the subcategory inclusion $\posclos \incl \poscat$.
\end{lem}
\begin{proof}
The category $\poscat$ is locally finitely presentable, so it has all small colimits.
It then suffices to show that coproducts in $\poscat$ are also coproducts in $\posclos$, and that the coequaliser in $\poscat$ of two closed order-preserving maps is a coequaliser in $\posclos$.

Given an indexed family of posets $(P_i)_{i \in I}$, the coproduct injections 
\[ \imath_i\colon P_i \incl \sum_{i \in I} P_i \] 
in $\poscat$ are all closed embeddings.
Given a diagram of closed order-preserving maps $(f_i\colon P_i \to Q)_{i \in I}$, let $f\colon \sum_{i \in I} P_i \to Q$ be the map induced by the universal property of the coproduct in $\poscat$.
For all $x \in \sum_{i \in I} P_i$, there exists a unique pair of $i \in I$ and $x' \in P_i$ such that $\imath_i(x') = x$.
Then 
\begin{align*}
	f(\clset{x}) & = f(\clset{\imath_i(x')}) = f(\imath_i(\clset{x'})) = \\
		     & = f_i(\clset{x'}) = \clset{f_i(x')} = \clset{f(x)},
\end{align*}
which by Lemma \ref{lem:map_closed_iff_maps_lowersets_to_lowersets} implies that $f$ is closed.

Next, let $f_0, f_1\colon P \to Q$ be a parallel pair of closed order-preserving maps, and let $\isocl{-}\colon Q \to C$ exhibit their coequaliser in $\poscat$.
This is a surjective order-preserving map, so every element in $C$ can be represented as $\isocl{y}$ for some $y \in Q$.
By construction of coequalisers in $\poscat$, for all $y, y' \in Q$, we have $\isocl{y} \leq \isocl{y'}$ if and only if $y \leq y'$, or there exist $m > 0$ and, for all $i \in \set{1, \ldots, m}$, elements $x_i \in P$ and $b_i \in \set{0, 1}$ such that
\[
	y \leq f_{b_1}(x_1), \quad \quad f_{1 - b_i}(x_i) \leq f_{b_{i+1}}(x_{i+1}), \quad \quad f_{1 - b_m}(x_m) \leq y'.
\]
To prove that $\isocl{-}$ is closed, we need to show that, if $\isocl{y} \leq \isocl{y'}$, then there exists $\tilde{y} \leq y'$ such that $\isocl{\tilde{y}} = \isocl{y}$.
For all $z \in Q$ such that $\isocl{z} = \isocl{y}$, let
\[
	m(z) \eqdef \begin{cases}
		0 & \text{if $z \leq y'$}, \\
		\min \set{m \mid \text{there exist $(x_i, b_i)_{i=1}^m$ exhibiting $\isocl{z} \leq \isocl{y'}$} } & \text{otherwise.}
	\end{cases}
\]
Let $\tilde{y} \in Q$ be such that $\isocl{\tilde{y}} = \isocl{y}$ and $m(\tilde{y})$ is minimal.
If $m(\tilde{y}) = 0$, then we are done.
Suppose $m \eqdef m(\tilde{y}) > 0$, and let $(x_i, b_i)_{i=1}^m$ be a sequence of minimal length exhibiting $\isocl{\tilde{y}} \leq \isocl{y}$.
Because $\tilde{y} \leq f_{b_1}(x_1)$ and $f_{b_1}$ is closed, there exists $x'_1 \leq x_1$ such that $f_{b_1}(x'_1) = \tilde{y}$.
Then 
\[
	z \eqdef f_{{1-b_1}}(x'_1) \leq f_{{1-b_1}}(x_1) \leq \begin{cases}
		y' & \text{if $m = 1$}, \\
		f_{b_2}(x_2) & \text{otherwise,}
	\end{cases}
\]
so 
\[
	\isocl{z} = \isocl{f_{{1-b_1}}(x'_1)} = \isocl{f_{b_1}(x'_1)} = \isocl{\tilde{y}} = \isocl{y}, 
\]
but $m(z) < m(\tilde{y})$, contradicting the minimality assumption for $\tilde{y}$.
This proves that $\isocl{-}$ is closed.

Given a closed map $g\colon Q \to R$ such that $g \after f_0 = g \after f_1$, let $h\colon C \to R$ be the map induced by the universal property of the coequaliser in $\poscat$.
Then, for all $\isocl{y} \in Q$, we have 
\begin{align*}
	h(\clset{\isocl{y}}) & = h(\isocl{\clset{y}}) = g(\clset{y}) = \\
			     & = \clset{g(y)} = \clset{h(\isocl{y})},
\end{align*}
proving that $h$ is closed.
\end{proof}

\begin{lem} \label{lem:pullbacks_of_closed_embeddings}
	The category $\posclos$ has pullbacks of closed embeddings, and they are preserved and reflected by the subcategory inclusion $\posclos \incl \poscat$.
	Moreover, closed embeddings are stable under pullbacks in both $\posclos$ and $\poscat$.
\end{lem}
\begin{proof}
	Let $\imath\colon V \incl Q$ be a closed embedding and $f\colon P \to Q$ an order-preserving map.
	The pullback of $\imath$ along $f$ in $\poscat$ can be constructed as the inclusion of the inverse image $\invrs{f}(\imath(V))$ into $P$.
	Since $\imath$ is a closed embedding, $\imath(V)$ is closed in $Q$, so by Lemma \ref{lem:inverse_images_of_closed} $\invrs{f}(\imath(V))$ is closed in $P$, and its inclusion into $P$ is a closed embedding.
	Moreover, the order-preserving map $\invrs{f}(\imath(V)) \to V$ is, essentially, the restriction of $f$ to a closed subset of $P$, so it is closed whenever $f$ is.
	
	It is then straightforward to prove that, when $f$ is closed, the pullback square is still a pullback square in $\posclos$, by observing that if $h = j \after g$ for order-preserving maps $h, j, g$ such that $h$ is closed and $j$ is a closed embedding, then $g$ is closed.
\end{proof}

\begin{lem} \label{lem:pushouts_of_spans_of_closed_embeddings}
	The pushout in $\posclos$ of a closed embedding along a closed embedding is a closed embedding.
	Moreover, any pushout square of closed embeddings is also a pullback square.
\end{lem}
\begin{proof}
	Consider a pushout square
\[\begin{tikzcd}
	U && Q \\
	P && R
	\arrow["\imath", hook', from=1-1, to=2-1]
	\arrow["f", from=2-1, to=2-3]
	\arrow["j", hook, from=1-1, to=1-3]
	\arrow["g", from=1-3, to=2-3]
	\arrow["\lrcorner"{anchor=center, pos=0.125, rotate=180}, draw=none, from=2-3, to=1-1]
\end{tikzcd}\]
	in $\posclos$ where $\imath$ and $j$ are closed embeddings, and let $x, y \in P$ such that $f(x) \leq f(y)$.
	By construction of pushouts in $\posclos$, either $x \leq y$, or there exist $m > 0$ and, for all $i \in \set{1, \ldots, m}$, elements $x_i \in U$ such that
\[
	x \leq \imath(x_1), \quad \quad
	\begin{cases}
		j(x_i) \leq j(x_{i+1}) & \text{if $i$ is odd}, \\
		\imath(x_i) \leq \imath(x_{i+1}) & \text{if $i$ is even},
	\end{cases} \quad \quad
	\imath(x_m) \leq y.
\]
But since both $\imath$ and $j$ are closed embeddings, by Lemma \ref{lem:closed_embedding_is_order_reflecting} they are order-reflecting, so in fact
\[
	x_1 \leq \ldots \leq x_m
\]
in $U$, hence $x \leq \imath(x_1) \leq \ldots \leq \imath(x_m) \leq y$ in $P$.
We conclude that $f$ is order-reflecting, and since $f$ is closed, it is a closed embedding.
By symmetry, $g$ is also a closed embedding.

To show that the square is a pullback square, it suffices to prove that $\imath(U) = \invrs{f}(g(Q))$.
One inclusion is immediate from commutativity of the square.
For the other, let $x \in P$ and suppose $f(x) \leq g(y)$ for some $y \in Q$.
By construction of pushouts in $\posclos$, this implies that there exists $x' \in U$ such that $x \leq \imath(x')$, and since $\imath$ is a closed embedding we conclude that $x \in \imath(U)$.
\end{proof}

\begin{prop} \label{prop:em_factorisation_in_posclos}
Let $f\colon P \to Q$ be a closed order-preserving map of posets.
Then $f$ factors as
\begin{enumerate}
	\item a surjective closed order-preserving map $\widehat{f}\colon P \to f(P)$,
	\item followed by a closed embedding $\imath\colon f(P) \incl Q$.
\end{enumerate}
This factorisation is unique up to unique isomorphism.
\end{prop}
\begin{proof}
	The set-theoretic image $f(P)$ of $P$ through a closed map is a closed subset of $Q$.
	Giving $f(P)$ the partial order induced from $Q$ makes the inclusion of $f(P)$ a closed embedding; in fact, this is the only partial order on $f(P)$ such that the inclusion is a closed embedding.
	Uniqueness of this factorisation up to unique isomorphism is then a consequence of the same property of the factorisation of functions of the underlying sets.
\end{proof}

\begin{cor} \label{cor:em_ofs_on_posclos}
The classes of
\begin{enumerate}
	\item surjective closed order-preserving maps,
	\item closed embeddings
\end{enumerate}
form an orthogonal factorisation system on $\posclos$.
\end{cor}
\begin{proof}
Both classes are evidently closed under composition and contain all isomorphisms.
The statement then follows from Proposition \ref{prop:em_factorisation_in_posclos}.
\end{proof}

\begin{dfn}[Local embedding] \index{map!local embedding} \index{local embedding!of posets}
Let $f\colon P \to Q$ be an order-preserving map of posets.
We say that $f$ is a \emph{local embedding} if, for all $x \in P$, the restriction $\restr{f}{\clset{x}}\colon \clset{x} \to Q$ is a closed embedding.
\end{dfn}

\begin{rmk}
	A local embedding is always closed.
\end{rmk}

\begin{dfn}[The category $\posloc$] \index{$\posloc$}
	We let $\posloc$ denote the category whose objects are posets and morphisms are local embeddings.
\end{dfn}


\section{Graded posets} \label{sec:graded}

\begin{guide}
	In this section, we focus on \cemph{graded posets}, whose Hasse diagram has elements stacked neatly into rows, the ``height'' of each element being determined by the length of a maximal outgoing downward path.
	As mentioned at the beginning of the chapter, face posets of regular cell complexes are graded by the topological dimension of each cell.
	We will over-extend this metaphor by speaking of the \cemph{dimension} of elements of a graded poset.

	From a technical standpoint, the main use of gradedness is that it allows us to prove things by induction on the dimension of elements.
	For this purpose, a weaker condition called \cemph{locally finite height} --- where outgoing downward paths are allowed to have different lengths, as long as their length has a uniform finite bound --- is often sufficient, and for reasons of tidiness we will consider it first.

	This section also sets up some important notation and terminology, which will recur all throughout the book, such as \cemph{faces}, \cemph{cofaces}, \cemph{intervals}, \cemph{codimension}, and the grading of subsets of a graded poset.
\end{guide}

\begin{dfn}[Size of a set] \index{size!of a set} \index{$\size{X}$}
Let $X$ be a set.
We let $|X|$ denote its \emph{size} or cardinality.
\end{dfn}

\begin{dfn}[Chain in a poset] \index{poset!chain} \index{size!of a chain} \index{$\size{c}$}
Let $P$ be a poset.
A \emph{chain in $P$} is an order-preserving map $c\colon L \to P$ where $L$ is a linear order.
We let $\size{c} \eqdef \size{c(L)}$.
\end{dfn}

\begin{dfn}[Finite height] \index{poset!of finite height} \index{$\height{P}$}
Let $P$ be a poset.
We say that $P$ has \emph{finite height} if there exists $n \in \mathbb{N}$ such that, for all chains $c$ in $P$, we have $\size{c} \leq n$.
In this case, the \emph{height} of $P$ is the integer
\[
	\height{P} \eqdef \min \set{ n \geq -1 \mid \text{for all chains $c$ in $P$, $\size{c} \leq n + 1$} }.
\]
\end{dfn}

\begin{dfn}[Maximal element] \index{poset!maximal element} \index{$\maxel{P}$}
Let $P$ be a poset, $x \in P$.
We say that $x$ is \emph{maximal in $P$} if, for all $y \in P$, if $x \leq y$ then $x = y$.
We write $\maxel{P}$ for the set of maximal elements in $P$.
\end{dfn}

\begin{dfn}[Greatest element] \index{poset!greatest element}
Let $P$ be a poset, $x \in P$.
We say that $x$ is the \emph{greatest element of $P$} if, for all $y \in P$, $y \leq x$.
\end{dfn}

\begin{lem} \label{lem:closure_of_maximal}
Let $P$ be a poset of finite height.
Then $P$ is the closure of its set of maximal elements, that is, $P = \clos{\maxel{P}}$.
\end{lem}
\begin{proof}
Since $\maxel{P} \subseteq P$, by Lemma \ref{lem:closure_is_monotonic} we have $\clos{\maxel{P}} \subseteq \clos{P} = P$.

Conversely, let $x \in P$.
We construct a chain $c\colon \mathbb{N} \to P$ as follows.
Let $c(0) \eqdef x$.
For each $i \geq 0$, if $c(i)$ is maximal in $U$, then let $c(j) \eqdef c(i)$ for all $j > i$.
Otherwise, there exists $x' > c(i)$ in $P$, and we let $c(i+1) \eqdef x'$.
Since $\size{c}$ is finite, there exist $i, j \in \mathbb{N}$ such that $i < j$ and $c(i) = c(j)$, hence $x \leq c(i)$ and $c(i)$ is maximal, so $x \in \clos{\maxel{P}}$.
\end{proof}

\begin{dfn}[Locally finite height] \index{poset!of locally finite height}
Let $P$ be a poset.
We say that $P$ has \emph{locally finite height} if, for all $x \in P$, the poset $\clset{x}$ has finite height.
\end{dfn}

\begin{rmk}
	Having locally finite height implies, but is stronger than well-foundedness: the latter stipulates that each chain under an element has finite size, but this size may be unbounded across all chains.
\end{rmk}
	
\begin{dfn}[Interval] \index{poset!interval} \index{$[x, y]$}
Let $P$ be a poset and let $x, y \in P$ such that $x \leq y$.
The \emph{interval from $x$ to $y$} is the subset
\[
	[x, y] \eqdef \set{ z \mid x \leq z \leq y }.
\]
\end{dfn}

\begin{lem} \label{lem:intervals_locally_finite_height}
Let $P$ be a poset of locally finite height.
Then every interval in $P$ has finite height.
\end{lem}
\begin{proof}
	Let $x, y \in P$ such that $x \leq y$.
	Then $[x, y] \subseteq \clset{y}$, so every chain in $[x, y]$ is a chain in $\clset{y}$.
	It follows that $\height{[x,y]}$ is bounded by $\height{\clset{y}}$.
\end{proof}

\begin{dfn}[Covering relation] \index{poset!covering relation}
Let $P$ be a poset.
Given elements $x, y \in P$, we say that $y$ \emph{covers} $x$ if $x < y$ and, for all $y' \in P$, if $x < y' \leq y$ then $y' = y$. 
\end{dfn}

\begin{dfn}[Faces and cofaces] \index{poset!faces} \index{poset!cofaces} \index{$\faces{}{} x$} \index{$\cofaces{}{} x$}
Let $P$ be a poset and $x \in P$.
The set of \emph{faces} of $x$ is
\begin{equation*}
    \faces{}{} x \eqdef \set{ y \in P \mid \text{$x$ covers $y$} }
\end{equation*}
and the set of \emph{cofaces} of $x$ is
\begin{equation*}
    \cofaces{}{} x \eqdef \set{ y \in P \mid \text{$y$ covers $x$} }.
\end{equation*}
\end{dfn}

\begin{lem} \label{lem:if_locally_finite_height_then_covers}
Let $P$ be a poset of locally finite height and $x, y \in P$ such that $x < y$.
Then there exists $x' \in [x, y] \cap \faces{}{}y$.
\end{lem}
\begin{proof}
	We construct a chain $c\colon \mathbb{N} \to [x, y]$ as follows.
	We let $c(0) \eqdef x$.
	For $i \geq 0$, if $y$ covers $c(i)$, we let $c(j) \eqdef c(i)$ for all $j > i$.
	Otherwise, there exists $x'$ such that $c(i) < x' < y$, and we let $c(i+1) \eqdef x'$.
	Since $\size{c}$ is finite, there exists $i, j \in \mathbb{N}$ such that $i < j$ and $c(i) = c(j)$, so $y$ covers $c(i)$.
\end{proof}

\begin{dfn}[Minimal element] \index{poset!minimal element} \index{$\minel{P}$}
Let $P$ be a poset.
We say that $x$ is \emph{minimal in $P$} if, for all $y \in P$, if $y \leq x$ then $y = x$.
We write $\minel{P}$ for the set of minimal elements in $P$.
\end{dfn}

\begin{dfn}[Least element] \index{poset!least element}
Let $P$ be a poset, $x \in P$.
We say that $x$ is the \emph{least element of $P$} if, for all $y \in P$, $x \leq y$.
\end{dfn}

\begin{lem} \label{lem:minimal_iff_no_faces}
Let $P$ be a poset of locally finite height, $x \in P$.
The following are equivalent:
\begin{enumerate}[label=(\alph*)]
	\item $x$ is minimal in $P$;
	\item $\faces{}{}x = \varnothing$.
\end{enumerate}
\end{lem}
\begin{proof}
	If $x$ is minimal, then there exist no $y < x$, so $\faces{}{}x = \varnothing$.
	Conversely, if $x$ is not minimal, then there exists $y < x$, so by Lemma \ref{lem:if_locally_finite_height_then_covers} there exists $x' \in \faces{}{}x$, and $\faces{}{}x \neq \varnothing$.
\end{proof}

\begin{dfn}[Directed graph] \index{directed graph}
A \emph{directed graph} is a diagram
\begin{equation*}
\mathscr{G} \eqdef
    \begin{tikzcd}
	{E_\mathscr{G}} & {V_\mathscr{G}}
	\arrow["s", shift left=1.5, from=1-1, to=1-2]
	\arrow["t"', shift right=1.5, from=1-1, to=1-2]
    \end{tikzcd}
\end{equation*}
of sets and functions.
The elements of $E_\mathscr{G}$ are called \emph{edges} and the elements of $V_\mathscr{G}$ are called \emph{vertices}.
Given an edge $e \in E_\mathscr{G}$, the \emph{source} of $e$ is $s(e)$ and the \emph{target} of $e$ is $t(e)$.
\end{dfn}

\begin{dfn}[Homomorphism of directed graphs] \index{directed graph!homomorphism}
	Let $\mathscr{G}$ and $\mathscr{H}$ be directed graphs.
	A \emph{homomorphism} $f\colon \mathscr{G} \to \mathscr{H}$ is a pair of functions
	\[
		f_V\colon V_\mathscr{G} \to V_\mathscr{H}, \quad \quad
		f_E\colon E_\mathscr{G} \to E_\mathscr{H}
	\]
	that preserve sources and targets, that is, satisfy
	\[ 
		s\after f_E = f_V \after s, \quad \quad t\after f_E = f_V \after t.
	\]
\end{dfn}

\begin{dfn}[The category $\gph$] \index{$\gph$}
There is a category $\gph$ whose objects are directed graphs and morphisms are homomorphisms of directed graphs.
\end{dfn}

\begin{dfn}[Hasse diagram] \index{poset!Hasse diagram} \index{$\hasse{P}$}
Let $P$ be a poset.
The \emph{Hasse diagram} of $P$ is the directed acyclic graph $\hasse{P}$ whose
\begin{itemize}
    \item set of vertices is the underlying set of $P$, and
    \item set of edges is 
	    \[ \set{(y, x) \mid x \in \faces{}{}y}, \]
	   with $s\colon (y, x) \mapsto y$ and $t\colon (y, x) \mapsto x$.
\end{itemize}
\end{dfn}

\begin{lem} \label{lem:local_embedding_preserves_faces}
Let $f\colon P \to Q$ be a local embedding of posets.
Then, for all $x \in P$, $f$ induces a bijection between $\faces{}{}x$ and $\faces{}{}f(x)$. 
\end{lem}
\begin{proof}
	By definition, $f$ restricts to an isomorphism between $\clset{x}$ and its image $\clset{f(x)}$.
	Such an isomorphism clearly preserves and reflects the covering relation.
\end{proof}

\begin{cor} \label{cor:local_embedding_hasse}
Let $f\colon P \to Q$ be a local embedding of posets.
Then $f$ induces a homomorphism $\hasse{f}\colon \hasse{P} \to \hasse{Q}$.
This assignment determines a functor $\hasse{}\colon \posloc \to \gph$.
\end{cor}

\begin{lem} \label{lem:hasse_diagram_paths}
Let $P$ be a poset of locally finite height and $x, y \in P$.
The following are equivalent:
\begin{enumerate}[label=(\alph*)]
	\item $x \leq y$;
	\item there exists a finite path from $y$ to $x$ in $\hasse{P}$.
\end{enumerate}
\end{lem}
\begin{proof}
	One implication is obvious.
	For the other one, suppose that $x \leq y$; we proceed by induction on $\height{[x, y]}$, which is well-defined by Lemma \ref{lem:intervals_locally_finite_height}.
	
	Notice that there is always a chain $c\colon \set{0 < 1} \to [x, y]$ with $c(0) = x$ and $c(1) = y$.
	If $\height{[x, y]} = 0$, then $\size{c} = 1$, hence $x = y$, and we have a trivial path in $\hasse{P}$.
	Suppose $\height{[x, y]} > 0$.
	Then $x < y$, and by Lemma \ref{lem:if_locally_finite_height_then_covers} there exists $x' \in [x, y]$ such that $y$ covers $x'$.
	We claim that $\height{[x, x']} < \height{[x, y]}$.
	Indeed, any chain $c$ in $[x, x']$ can be extended to a chain $c'$ in $[x, y]$ by adding a new greatest element to its domain and mapping it to $y$; then $\size{c'} = \size{c} + 1$.
	By the inductive hypothesis, there exists a finite path from $x'$ to $x$ in $\hasse{P}$.
	Concatenating it with the path that traverses the edge from $y$ to $x'$, we obtain a finite path from $y$ to $x$.
\end{proof}

\begin{lem} \label{lem:height_bounds_length_of_chains}
Let $P$ be a poset of locally finite height.
Then every path in $\hasse{P}$ has length at most $\height{P}$.
\end{lem}
\begin{proof}
Let $x_n \to x_{n-1} \to \ldots \to x_0$ be a path of length $n$ in $\hasse{P}$.
Then 
\[
	c\colon \set{0 < \ldots < n} \to \clset{x}, \quad \quad i \mapsto x_i
\]
is a chain in $P$ with $\size{c} = n+1$.
It follows that $n \leq \height{P}$.
\end{proof}

\begin{dfn}[Graded poset] \index{poset!graded|see {graded poset}} \index{graded poset}
Let $P$ be a poset of locally finite height.
We say that $P$ is \emph{graded} if, for all $x \in P$, all maximal paths starting from $x$ in $\hasse{P}$ have the same length.
\end{dfn}

\begin{dfn}[Dimension of an element] \index{graded poset!dimension!of an element} \index{$\dim{x}$}
Let $P$ be a graded poset and $x \in P$.
The \emph{dimension} of $x$ is the length $\dim{x}$ of a maximal path starting from $x$ in $\hasse{P}$.
\end{dfn}

\begin{comm}
The dimension is more commonly known in order theory as the \emph{rank} or \emph{degree} of an element.
We use this terminology to enhance the elements-as-cells metaphor.
\end{comm}

\begin{lem} \label{lem:minimal_iff_dim_0}
Let $P$ be a graded poset, $x \in P$.
The following are equivalent:
\begin{enumerate}[label=(\alph*)]
	\item $x$ is minimal in $P$;
	\item $\dim{x} = 0$.
\end{enumerate}
\end{lem}
\begin{proof}
	Suppose that $\faces{}{}x = \varnothing$.
	It follows that there are no edges in $\hasse{P}$ whose source is $x$, hence no non-constant paths starting from $x$, so $\dim{x} = 0$.
	Conversely, if there are no non-constant paths starting from $x$, then $\faces{}{}x = \varnothing$.
	We conclude by Lemma \ref{lem:minimal_iff_no_faces}.
\end{proof}

\begin{dfn}[Grading of a subset] \index{graded poset!grading}
Let $U$ be a subset of a graded poset.
For each $n \in \mathbb{N}$, we write $\grade{n}{U} \eqdef \set{ x \in U \mid \dim x = n }$.
We have $U = \sum_{n \in \mathbb{N}} \grade{n}{U}$.
\end{dfn}

\begin{dfn}[Skeleta of a graded poset] \index{skeleton!of a graded poset} \index{$\skel{n}{P}$} \index{graded poset!skeleton}
	Let $P$ be a graded poset, $n \in \mathbb{N}$.
	The \emph{$n$\nbd skeleton of $P$} is the closed subset $\skel{n}{P} \eqdef \sum_{k \leq n} \grade{k}{P}$.
\end{dfn}

\begin{comm}
	It is useful to extend the indexing of $\grade{n}{U}$ and $\skel{n}{P}$ to negative integers, by stipulating that they are both equal to $\varnothing$ when $n < 0$.
\end{comm}

\begin{lem} \label{lem:dimension_is_height}
Let $P$ be a graded poset and $x \in P$.
Then 
\[	\dim{x} = \height{\clset{x}}.\]
\end{lem}
\begin{proof}
By Lemma \ref{lem:height_bounds_length_of_chains}, since every path in $\hasse{P}$ starting from $x$ is contained in $\clset{x}$, we have $\dim{x} \leq \height{\clset{x}}$.

Let $c\colon L \to \clset{x}$ be a chain with $L$ non-empty, $n \eqdef \size{c} - 1$.
Then there exist $i_0 < \ldots < i_n$ in $L$ such that 
\[
	c(i_0) < \ldots < c(i_n) \leq x.
\]
By Lemma \ref{lem:hasse_diagram_paths}, there exist finite paths in $\hasse{P}$ from $x$ to $c(i_n)$ and non-trivial finite paths from $c(i_{j+1})$ to $c(i_j)$ for all $j \in \set{0, \ldots, n-1}$.
Concatenating them all, we obtain a path of length $> n$ starting from $x$.
It follows that $n \leq \dim{x}$, and since $c$ was arbitrary, $\height{\clset{x}} \leq \dim{x}$.
\end{proof}

\begin{lem} \label{lem:dimension_is_monotonic}
Let $P$ be a graded poset and $x, y \in P$ such that $x \leq y$.
Then
\begin{enumerate}
	\item $\dim{x} \leq \dim{y}$,
	\item $\dim{x} = \dim{y}$ if and only if $x = y$.
\end{enumerate}
\end{lem}
\begin{proof}
Take a maximal path starting from $x$ in $\hasse{P}$.
Because $P$ is graded, this path has length $\dim{x}$.
Since $x \leq y$, by Lemma \ref{lem:hasse_diagram_paths} there is a finite path from $y$ to $x$ in $\hasse{P}$; let $k$ be its length.
Concatenating the two paths necessarily produces a maximal path from $y$, whose length $k + \dim{x}$ is equal to $\dim{y}$.
It follows that $\dim{y} - \dim{x} = k \geq 0$, which is only equal to $0$ if the path from $y$ to $x$ has length $0$, that is, if $y = x$.
\end{proof}

\begin{lem} \label{lem:closure_of_set_and_closure_of_faces}
Let $P$ be a graded poset, $x \in P$.
Then 
\[ \clset{x} = \set{x} + \clos{\faces{}{} x}. \]
\end{lem}
\begin{proof}
Let $y \in \clset{x}$, that is, $y \leq x$ in $P$.
Then either $x = y$, so $y \in \set{x}$, or there exists a non-trivial path from $x$ to $y$ in $\hasse{P}$, which necessarily passes through some $y' \in \faces{}{} x$, so $y \in \clset{y'} \subseteq \clos{\faces{}{} x}$.
This proves that $\clset{x} \subseteq \set{x} \cup \clos{\faces{}{} x}$.
The converse inclusion is straightforward, as is the fact that $\set{x}$ and $\clos{\faces{}{} x}$ are disjoint.
\end{proof}

\begin{lem} \label{lem:local_embedding_into_graded}
	Let $f\colon P \to Q$ be a local embedding of posets.
	Then
	\begin{enumerate}
		\item if $Q$ has locally finite height, so does $P$,
		\item if $Q$ is graded, so is $P$.
	\end{enumerate}
\end{lem}
\begin{proof}
Let $x \in P$.
By definition, $f$ restricts to an isomorphism between $\clset{x}$ and $\clset{f(x)}$.
Then, if $\clset{f(x)}$ has finite height, so does $\clset{x}$.
Furthermore, by Corollary \ref{cor:local_embedding_hasse}, $\hasse{\restr{f}{\clset{x}}}$ is a graph isomorphism between $\hasse{\clset{x}}$ and $\hasse{\clset{f(x)}}$.
This establishes a bijection between (maximal) paths starting from $x$ in $\hasse{P}$, which are entirely contained in $\hasse{\clset{x}}$, and (maximal) paths starting from $f(x)$ in $\hasse{Q}$, which are entirely contained in $\hasse{\clset{f(x)}}$. 
\end{proof}

\begin{lem} \label{lem:closed_map_of_graded_posets_dim_non_increasing}
	Let $f\colon P \to Q$ be a closed order-preserving map of graded posets.
	Then $f$ is dimension-non-increasing.
\end{lem}
\begin{proof}
	Let $x \in P$.
	Because $\clset{f(x)} = f(\clset{x})$, every chain in $\clset{f(x)}$ factors as a chain in $\clset{x}$ followed by $f$.
	We conclude by Lemma \ref{lem:dimension_is_height}.
\end{proof}

\begin{dfn}[Dimension of a graded poset] \index{graded poset!finite-dimensional} \index{$\dim{P}$} \index{dimension!of a graded poset}
Let $P$ be a graded poset.
The \emph{dimension} of $P$ is 
\[    
	\dim P \eqdef
	\begin{cases}
		\max \left(\set{-1} \cup \set{ \dim x \mid x \in P }\right) & \text{if defined,} \\
		\infty & \text{otherwise.}
	\end{cases}
\]
We say that $P$ is \emph{finite-dimensional} if $\dim{P} < \infty$.
\end{dfn}

\begin{dfn}[Codimension of an element] \index{graded poset!codimension} \index{$\codim{x}{U}$}
Let $U$ be a finite-dimensional closed subset of a graded poset, $x \in U$.
The \emph{codimension} of $x$ in $U$ is the integer 
\[
	\codim{x}{U} \eqdef \dim{U} - \dim{x}.
\]
When $U = \clset{y}$, we write $\codim{x}{y} \eqdef \codim{x}{\clset{y}} = \dim{y} - \dim{x}$.
\end{dfn}

\begin{lem}
Let $P$ be a graded poset.
The following are equivalent:
\begin{enumerate}[label=(\alph*)]
	\item $P$ is finite-dimensional;
	\item $P$ has finite height.
\end{enumerate}
Moreover, in either case, $\dim{P} = \height{P}$.
\end{lem}
\begin{proof}
	Suppose $P$ is empty.
	Then $\dim{P} = -1$ and every chain in $P$ has size 0, so $\height{P} = -1$.

	Suppose that $P$ is non-empty and finite-dimensional, and let $c\colon L \to P$ be a chain in $P$ of size $> 0$.
	For all $i \in L$, $c$ restricts to a chain
	\[
		\restr{c}{\clset{i}}\colon \clset{i} \to \clset{c(i)}
	\]
	and by Lemma \ref{lem:dimension_is_height} $\height{\clset{c(i)}} = \dim{c(i)}$.
	It follows that
	\[
		\size{c(L)} = \size{ \bigcup_{i \in L} c(\clset{i}) } \leq \max \set{ \dim{c(i)} \mid i \in L } + 1 \leq \dim{P} + 1,
	\]
	so $P$ has finite height and $\height{P} \leq \dim{P}$.

	Conversely, suppose that $P$ is non-empty and has finite height.
	Then for all $x \in P$, we have $\dim{x} = \height{\clset{x}} \leq \height{P}$.
	It follows that $P$ is finite-dimensional and $\dim{P} \leq \height{P}$.
	Then in either case $\dim{P} = \height{P}$.
\end{proof}

\begin{dfn}[Pure graded poset] \index{graded poset!pure}
Let $P$ be a finite-dimensional graded poset.
We say that $P$ is \emph{pure} if all the maximal elements of $P$ have codimension $0$, that is, for all $x \in \maxel{P}$, $\dim{x} = \dim{P}$.
\end{dfn}

\begin{exm}[The poset of divisors of 12]
Let $P$ be the poset of divisors of 12, ordered by divisibility.
The following is a depiction of its Hasse diagram.
\[
	\input{img/hasse_divisors.tex}
\]
Then $P$ is a graded poset, with 
\[
	\grade{0}{P} = \set{1}, \quad \grade{1}{P} = \set{2, 3}, \quad \grade{2}{P} = \set{4, 6}, \quad \grade{3}{P} = \set{12};
\]
in fact, the poset of divisors of any natural number is graded, with $\dim{n}$ equal to the number of prime factors of $n$ counted with their multiplicity.
The dimension of each element can be reconstructed as their ``height'' in the Hasse diagram.

Consider the subset $\set{1, 3, 4}$.
This is not closed, since it does not contain $2 \leq 4$.
Let $U \eqdef \clset{3, 4} = \set{1, 2, 3, 4}$.
Then
\[
	\dim{U} = 2, \quad \maxel{U} = \set{3, 4}.
\]
It follows that $U$ is closed but not pure, since $3 \in \maxel{U}$ but $\dim{3} = 1$.
\end{exm}


\section{Some operations on posets} \label{sec:poset_operations}

\begin{guide}
	In this section, we consider some common constructions and operations on posets, and how they interact with the notions that we considered in previous sections: closed maps and grading.
	We observe, for instance, that the usual \cemph{product} of posets is \emph{not} a categorical product in $\posclos$, but does determine a symmetric monoidal structure.

	We also take the opportunity to clarify a potential source of confusion.
	In poset topology, face posets are conventionally assumed to have a least element, whose dimension is set to $-1$.
	This is best thought of as a ``formal'' least element, to be removed and added back at one's convenience.
	We model this as an equivalence of categories, whose two sides, \cemph{augmentation} and \cemph{diminution}, add and remove the formal least element, respectively.
	Most usefully, the product of posets with a least element determines a symmetric monoidal structure which is \emph{not} equivalent, via diminution, to the product of posets, but instead to a different monoidal product, that we call \cemph{cellular join}, as it corresponds to the topological join on face posets of regular cell complexes.
\end{guide}

\begin{dfn}[Product of posets] \index{poset!product} \index{$P \times Q$} \index{product!of posets}
	Let $I$ be a set and $(P_i)_{i\in I}$ be an $I$\nbd indexed family of posets.
	The \emph{product} $\prod_{i \in I} P_i$ is its categorical product in $\poscat$.
	Explicitly, this can be constructed as 
\begin{itemize}
	\item the set of $I$\nbd indexed sequences $(x_i)_{i \in I}$ with $x_i \in P_i$ for all $i \in I$, with
	\item the order relation $(x_i)_{i \in I} \leq (y_i)_{i \in I}$ if and only if $x_i \leq y_i$ for all $i \in I$, and 
	\item the projections defined by $\pi_j\colon (x_i)_{i \in I} \mapsto x_j$ for all $j \in I$.
\end{itemize}
	The \emph{binary product} $(P, Q) \mapsto P \times Q$, together with the \emph{terminal} poset $1$ with one element, induce a cartesian monoidal structure $(\poscat, \times, 1)$ on $\poscat$.
\end{dfn}

\begin{prop} \label{prop:posclos_monoidal}
	The cartesian monoidal structure $(\poscat, \times, 1)$ restricts to a monoidal structure on $\posclos$.
\end{prop}
\begin{proof}
	It suffices to show that, if $f\colon P \to P'$ and $g\colon Q \to Q'$ are closed, then $f \times g\colon P \times Q \to P' \times Q'$ is closed.
	This is defined by
	\[
		f \times g\colon (x, y) \mapsto (f(x), g(y)).
	\]
	Suppose that $(x', y') \leq (f(x), g(y))$.
	Then $x' \leq f(x)$ and $y' \leq g(y)$.
	Since $f$ and $g$ are closed, there exist $x'' \leq x$ and $y'' \leq y$ such that $f(x'') = x'$ and $g(y'') = y'$.
	Then $(x'', y'') \leq (x, y)$ and $(f(x''), g(y'')) = (x', y')$, so by Lemma 
	\ref{lem:map_closed_iff_maps_lowersets_to_lowersets} $f \times g$ is closed.
\end{proof}

\begin{rmk}
	Note that $P \times Q$ is \emph{not}, in general, a categorical product in $\posclos$, so the monoidal structure of Proposition \ref{prop:posclos_monoidal} is not cartesian.
	In particular, the diagonal map $P \to P \times P$ is not a closed map as soon as $P$ has at least two related elements.
\end{rmk}

\begin{lem} \label{lem:product_preserves_colimits}
	Let $P$ be a poset.
	The functors $P \times -$ and $- \times P$ preserve all small colimits in $\posclos$.
\end{lem}
\begin{proof}
	The category $\poscat$ is cartesian closed, so $P \times -$ and $- \times P$ are left adjoint functors and preserve all small colimits in $\poscat$.
	The statement then follows from Lemma \ref{lem:colimits_in_posclos}.
\end{proof}

\begin{lem} \label{lem:product_of_posets_faces}
	Let $P$, $Q$ be posets and let $(x, y) \in P \times Q$.
	Then
	\[
		\faces{}{}(x, y) = \faces{}{}x \times \set{y} + \set{x} \times \faces{}{}y.
	\]
\end{lem}
\begin{proof}
	Suppose that $(x', y') \in \faces{}{}(x, y)$.
	If $x' \neq x$, then $x' < x$, so we have $(x', y') < (x, y') \leq (x, y)$, and by definition of the covering relation $(x, y') = (x, y)$, so $y = y'$.
	If $x''$ is such that $x' < x'' \leq x$, then we have $(x', y) < (x'', y) \leq (x, y)$, so $(x'', y) = (x, y)$ and $x'' = x$.
	This proves that $x' \in \faces{}{}x$.
	Similarly, if $y' \neq y$, we deduce that $x' = x$ and $y' \in \faces{}{}y$.
	Then
	\[
		\faces{}{}(x, y) \subseteq \faces{}{}x \times \set{y} + \set{x} \times \faces{}{}y.
	\]
	Conversely, if $x' \in \faces{}{}x$ and $(x', y) < (x'', y'') \leq (x, y)$, then necessarily $y'' = y$ and $x' < x'' \leq x$, so $x'' = x$, and we conclude that $(x', y) \in \faces{}{}(x, y)$.
	Similarly, if $y' \in \faces{}{}y$, we deduce that $(x, y') \in \faces{}{}(x, y)$.
\end{proof}

\begin{prop} \label{prop:product_of_graded_is_graded}
	Let $P$, $Q$ be graded posets.
	Then
	\begin{enumerate}
		\item the product $P \times Q$ is graded,
		\item for all $(x, y) \in P \times Q$, $\dim{(x, y)} = \dim{x} + \dim{y}$.
	\end{enumerate}
\end{prop}
\begin{proof}
	Let $x \in P$ and $y \in Q$.
	We will show that $\clset{(x, y)}$ is graded and has dimension $\dim{x} + \dim{y}$ by double induction on $(\dim{x}, \dim{y})$.
	When $\dim{x} = \dim{y} = 0$, then $(x, y)$ is minimal, so $\clset{(x, y)} = \set{(x, y)}$ and $\dim{(x, y)} = 0$.
	
	Otherwise, by Lemma \ref{lem:product_of_posets_faces}, every path from $(x, y)$ in $\hasse{(P \times Q)}$ begins with an edge to $(x', y)$ with $x' \in \faces{}{}x$, or an edge to $(x, y')$ with $y' \in \faces{}{}y$.
	By the inductive hypothesis, $\clset{(x', y)}$ and $\clset{(x, y')}$ are graded with dimension $\dim{x} + \dim{y} - 1$, and we conclude.
\end{proof}

\begin{dfn}[The category $\posbot$] \index{$\posbot$}
We let $\posbot$ denote the subcategory of $\poscat$ whose
\begin{itemize}
	\item objects are posets with a least element,
	\item morphisms preserve and reflect the least element.
\end{itemize}
\end{dfn}

\begin{dfn}[The category $\posclosbot$] \index{$\posclosbot$}
We let $\posclosbot$ denote the subcategory of $\posclos$ whose
\begin{itemize}
	\item objects are posets with a least element,
	\item morphisms reflect the least element.
\end{itemize}
\end{dfn}

\begin{rmk}
A closed order-preserving map of posets with a least element always preserves the least element.
Thus we have a commutative diagram
\[\begin{tikzcd}
	\posclosbot && \posclos \\
	\posbot && \poscat
	\arrow[hook, from=1-1, to=1-3]
	\arrow[hook', from=1-1, to=2-1]
	\arrow[hook, from=2-1, to=2-3]
	\arrow[hook', from=1-3, to=2-3]
\end{tikzcd}\]
of inclusions of subcategories.
\end{rmk}

\begin{lem} \label{lem:products_restrict_to_posbot}
	The monoidal structure $(\poscat, \times, 1)$ restricts to a monoidal structure on $\posbot$ and $\posclosbot$.
\end{lem}
\begin{proof}
	If $P$, $Q$ are posets with least elements $\bot_P, \bot_Q$, then $(\bot_P, \bot_Q)$ is the least element of $P \times Q$.
	Moreover, the unit $1$ evidently has a least element.
	Finally, if $f$ and $g$ preserve and reflect the least element, then so does $f \times g$.
\end{proof}

\begin{dfn}[Augmentation of a poset] \index{poset!augmentation} \index{$\augm{P}$} \index{augmentation!of a poset}
	Let $P$ be a poset.
	The \emph{augmentation of $P$} is the poset $\augm{P}$ whose 
	\begin{itemize}
		\item underlying set is $\set{\augm{x} \mid x \in P} + \set{\bot}$,
		\item partial order is defined by
			\[
				\clset{x} \eqdef \begin{cases}
					\set{\bot} & 
					\text{if $x = \bot$}, \\
					\set{\bot} + \set{\augm{y} \mid y \leq x'} &
					\text{if $x = \augm{x'}$, $x' \in P$}.
				\end{cases}
			\]
	\end{itemize}
	By construction, $\bot$ is the least element of $\augm{P}$.
\end{dfn}

\begin{lem} \label{lem:augmentation_is_a_functor}
	Let $f\colon P \to Q$ be an order-preserving map of posets, and let
	\begin{align*}
		\augm{f}\colon \augm{P} & \to \augm{Q}, \\
			x & \mapsto \begin{cases}
				\bot & \text{if $x = \bot$,} \\
				\augm{f(x')} & \text{if $x = \augm{x'}$, $x' \in P$}.
			\end{cases}
	\end{align*}
	Then
	\begin{enumerate}
		\item $\augm{f}$ is an order-preserving map of posets,
		\item if $f$ is closed, then so is $\augm{f}$.
	\end{enumerate}
	This assignment determines a functor $\augm{(-)}\colon \poscat \to \posbot$, restricting to a functor $\augm{(-)}\colon \posclos \to \posclosbot$.
\end{lem}
\begin{proof}
	Let $x \in \augm{P}$.
	If $x = \bot$, then
	\[
		\augm{f}(\clset{x}) = \augm{f}(\set{\bot}) = \set{\bot} = \clset{\augm{f}(x)}.
	\]
	If $x = \augm{x'}$ for some $x' \in P$, then
	\begin{align*}
		\augm{f}(\clset{x}) & = \set{\bot} \cup \set{\augm{f}(\augm{y}) \mid y \leq x'} =
		\set{\bot} \cup \set{\augm{f(y)} \mid y \leq x'} \subseteq \\
				    & \subseteq \set{\bot} \cup \set{\augm{y'} \mid y' \leq f(x')} = \clset{\augm{f(x')}} = \clset{\augm{f}(x)},
	\end{align*}
	with equality if $f$ is also closed.
	It follows that $\augm{f}$ is always order-preserving, and it is closed whenever $f$ is closed. 
	By construction, $\augm{f}$ preserves and reflects the least element.
	Functoriality is straightforward.
\end{proof}

\begin{lem} \label{lem:cofaces_in_augmentation}
	Let $P$ be a poset, $x \in \augm{P}$.
	Then 
	\[
		\cofaces{}{}x = \begin{cases}
			\set{\augm{y} \mid y \in \minel{P}} & 
				\text{if $x = \bot$,} \\
			\set{\augm{y} \mid y \in \cofaces{}{}x'} &
				\text{if $x = \augm{x'}$, $x' \in P$}.
			\end{cases}
	\]
\end{lem}
\begin{proof}
	Suppose $x = \bot$ and let $y \in \cofaces{}{}\bot$.
	Then $y = \augm{y'}$ for some $y' \in P$.
	Suppose that $x' \leq y'$ in $P$.
	Then $\bot < \augm{x'} \leq \augm{y'}$ in $\augm{P}$, so $\augm{x'} = \augm{y'}$, hence $x' = y'$.
	It follows that $y' \in \minel{P}$.
	Conversely, if $y' \in \minel{P}$, we have $\bot < \augm{y'}$.
	If $\bot < z \leq \augm{y'}$, then $z = \augm{z'}$ for some $z' \leq y'$ in $P$, so $z' = y'$, hence $\augm{y'} \in \cofaces{}{}\bot$.

	The case $x = \augm{x'}$ follows immediately from the fact that $\augm{x'} \leq y$ in $\augm{P}$ if and only if $y = \augm{y'}$ and $x' \leq y'$ for some $y' \in P$.
\end{proof}

\begin{dfn}[Diminution of a poset with least element] \index{$\dimin{P}$} \index{poset!diminution} \index{diminution!of a poset}
	Let $P$ be a poset with a least element $\bot$.
	The \emph{diminution of $P$} is the subset $\dimin{P} \eqdef P \setminus \set{\bot}$ of $P$ with the induced order.
\end{dfn}

\begin{lem} \label{lem:augmentation_has_an_inverse}
	Diminution extends to a functor $\dimin{(-)}\colon \posbot \to \poscat$ which is inverse to $\augm{(-)}\colon \poscat \to \posbot$ up to natural isomorphism.
	The equivalence restricts to an equivalence between $\posclos$ and $\posclosbot$.
\end{lem}
\begin{proof}
	It suffices to show that $\augm{(-)}$ is an equivalence and that, for all posets $P$ with a least element, $\augm{(\dimin{P})}$ is isomorphic to $P$.
	Let $P$ be a poset with a least element $\bot$.
	Then
	\begin{align*}
		\varphi\colon P &\to \augm{(\dimin{P})}, \\
		x & \mapsto \begin{cases}
			\bot & \text{if $x = \bot$}, \\
			\augm{x} & \text{if $x \neq \bot$}
		\end{cases}
	\end{align*}
	is a bijection at the level of the underlying sets, and it is straightforward to check that it is order-preserving and order-reflecting, so it is an isomorphism of posets.
	This proves that $\augm{(-)}$ is essentially surjective.

	Let $f\colon \augm{P} \to \augm{Q}$ be a morphism in $\posbot$.
	We define a function $f'\colon P \to Q$ as follows.
	Let $x \in P$.
	Because $f$ reflects the least element, there exists a necessarily unique $y \in Q$ such that $f(\augm{x}) = \augm{y}$, and we let $f'(x) \eqdef y$.
	Since $f$ is order-preserving, so is $f'$, and by construction $f = \augm{f'}$.
	Moreover, if $f$ is closed, then so is $f'$.
	This proves that $\augm{(-)}$ is full onto $\posbot$, and its restriction to $\posclos$ is full onto $\posclosbot$.

	Finally, let $f, g\colon P \to Q$ be order-preserving maps of posets, and suppose $\augm{f} = \augm{g}$.
	Then for all $x \in P$ we have $\augm{f(x)} = \augm{f}(\augm{x}) = \augm{g}(\augm{x}) = \augm{g(x)}$, so $f(x) = g(x)$, hence $f = g$.
	This proves that $\augm{(-)}$ is faithful.
\end{proof}

\begin{lem} \label{lem:posclosbot_connected_colimits}
The subcategory inclusion $\posclosbot \incl \posclos$
\begin{enumerate}
	\item reflects all colimits,
	\item preserves colimits of connected diagrams.
\end{enumerate}
\end{lem}
\begin{proof}
	Let $\fun{F}\colon \smcat{C} \to \posclosbot$ be a diagram, let $\gamma$ be a cone under $\fun{F}$ with tip $P$ in $\posclosbot$, and let $\bot_P$ be the least element of $P$.

	First, suppose $\gamma$ is a colimit cone in $\posclos$.
	Let $\eta$ be another cone under $\fun{F}$ with tip $Q$ in $\posclosbot$; to prove that $\gamma$ is a colimit cone in $\posclosbot$, it suffices to show that the universal closed order-preserving map $f\colon P \to Q$ induced by the universal property of $\gamma$ in $\posclos$ reflects the least element.
	Let $x \in P$.
	Then there exists an object $c$ in $\smcat{C}$ and $y \in \fun{F}c$ such that $x = \gamma_c(y)$.
	Since $\eta_c$ reflects the least element, if $f(x) = f(\gamma_c(y)) = \eta_c(y)$ is the least element of $Q$, then $y$ is the least element of $\fun{F}c$, so $x = \gamma_c(y) = \bot_P$.

	Next, suppose that $\smcat{C}$ is connected and that $\gamma$ is a colimit cone in $\posclosbot$, and let $\eta$ be a cone under $\fun{F}$ with tip $Q$ in $\posclos$.
	First of all, we will show that $Q$ has a least element.
	Since $\smcat{C}$ is connected, there exists an object $c$ in $\smcat{C}$.
	Let $\bot_c$ be the least element of $\fun{F}c$, and let
	\[
		\bot_Q \eqdef \eta_c(\bot_c).
	\]
	We claim that $\bot_Q$ is the least element of $Q$.
	Because $\eta_c$ is closed, 
	\[
		\eta_c(\clset{\bot_c}) = \eta_c(\set{\bot_c}) = \set{\bot_Q} = \clset{\eta_c(\bot_c)} = \clset{\bot_Q}
	\]
	so $\bot_Q$ is minimal.
	Suppose that $x$ is another minimal element in $Q$.
	Then there exist $c'$ in $\smcat{C}$ and $y \in \fun{F}c'$ such that $x = \eta_{c'}(y)$, and we can take $y$ to be the least element $\bot_{c'}$ of $\fun{F}c'$.
	Because $\smcat{C}$ is connected, there is a zig-zag of morphisms
\[\begin{tikzcd}
	c \equiv c_0 & {c_1} & {c_2} & \ldots & {c_{m-2}} & {c_{m-1}} & {c_m \equiv c'} 
	\arrow["{f_0}", from=1-2, to=1-1]
	\arrow["{f_1}"', from=1-2, to=1-3]
	\arrow[from=1-4, to=1-3]
	\arrow[from=1-4, to=1-5]
	\arrow["{f_{m-2}}", from=1-6, to=1-5]
	\arrow["{f_{m-1}}"', from=1-6, to=1-7]
\end{tikzcd}\]
	connecting $c$ and $c'$ in $\smcat{C}$.
	Moreover, $\fun{F}f_i$ preserves and reflects the least element for all $i \in \set{0, \ldots, m}$, so from the fact that
	\[
		\eta_{c_{2k-2}} \after \fun{F}f_{2k-2} = \eta_{c_{2k-1}} = \eta_{c_{2k}} \after \fun{F}f_{2k-1}
	\]
	for all $k \in \set{1, \ldots, \frac{m}{2}}$, we deduce that $x = \bot$.

	Now, by Lemma \ref{lem:augmentation_has_an_inverse}, $\dimin{\gamma}$ is a colimit cone under $\dimin{\fun{F}}$ in $\posclos$, hence a colimit cone in $\poscat$ by Lemma \ref{lem:colimits_in_posclos}.
	Moreover, $\eta$ restricts to a cone $\eta'$ under $\dimin{\fun{F}}$ in $\poscat$, defined by
	\[
		\eta'_c \eqdef \restr{(\eta_c)}{\dimin{(\fun{F}c)}}
	\]
	for each object $c$ in $\smcat{C}$; this is well-defined because all maps in the image of $\fun{F}$ preserve and reflect the least element.
	Let $f'\colon \dimin{P} \to Q$ be the unique order-preserving map induced by the universal property of $\dimin{\gamma}$ in $\poscat$.
	We let
	\begin{align*}
		f\colon P & \to Q, \\
			x & \mapsto \begin{cases}
				f'(x) & \text{if $x \in \dimin{P}$}, \\
				\bot_Q & \text{if $x = \bot_P$}.
			\end{cases}
	\end{align*}
	Since all components of $\eta$ and $\gamma$ preserve the least element and $f' \after \eta'_c = \dimin{(\gamma_c)}$, we have
	\[
		f \after \eta_c = \gamma_c
	\]
	for all objects $c$ in $\smcat{C}$, and uniqueness of $f'$ implies uniqueness of $f$ with this property.
	This proves that $\gamma$ is a colimit cone in $\posclos$.
\end{proof}

\begin{rmk}
	Note that $\posclosbot \incl \posclos$ does not preserve colimits of non-connected diagrams: for example, $1$ is an initial object in $\posclosbot$, but it is not an initial object in $\posclos$.
\end{rmk}

\begin{prop} \label{prop:augmentation_preserves_graded}
	Let $P$ be a poset of locally finite height.
	Then
	\begin{enumerate}
		\item $\augm{P}$ is graded if and only if $P$ is graded,
		\item if $P$ is graded, then for all $x \in \augm{P}$
			\[
				\dim{x} = \begin{cases}
					0 &
					\text{if $x = \bot$,} \\
					\dim{x'} + 1 &
					\text{if $x = \augm{x'}$, $x' \in P$}.
				\end{cases}
			\]
	\end{enumerate}
\end{prop}
\begin{proof}
	Let $x \in \augm{P}$.
	Given a path
	\[
		x \equiv x_0 \to x_1 \to \ldots \to x_{m-1} \to x_m
	\]
	in $\hasse{\augm{P}}$, by Lemma \ref{lem:minimal_iff_no_faces} the path is maximal if and only if $x_m = \bot$.
	Moreover, by Lemma \ref{lem:cofaces_in_augmentation}, for each $i \in \set{1, \ldots, m}$,
	\begin{itemize}
		\item if $x_i = \augm{(x'_i)}$ for some $x'_i \in P$, then $x_{i-1} = \augm{(x'_{i-1})}$ for some $x'_{i-1} \in \cofaces{}{}x'_i$,
		\item if $x_m = \bot$, then $x_{m-1} = \augm{(x'_{m-1})}$ for some $x'_{m-1} \in \minel{P}$.
	\end{itemize}
	It follows by backward recursion that if $x = \bot$, then $m = 0$, so $\dim{\bot}$ is always well-defined and equal to $0$, while if $x = \augm{x'}$, then there is a maximal path
	\[
		x' \equiv x'_0 \to x'_1 \to \ldots \to x'_{m-1}
	\]
	in $\hasse{P}$.
	Conversely, if
	\[
		x' \equiv x'_0 \to x'_1 \to \ldots \to x'_{m-1}
	\]
	is a maximal path in $P$, then $x'_{m-1} \in \minel{P}$, so
	\[
		x \equiv \augm{(x'_0)} \to \augm{(x'_1)} \to \ldots \to \augm{(x'_{m-1})} \to \bot
	\]
	is a maximal path in $\hasse{\augm{P}}$.
	We conclude that $\dim{x}$ is well-defined and equal to $m$ if and only if $\dim{x'}$ is well-defined and equal to $m-1$.
\end{proof}

\begin{cor} \label{cor:diminution_preserves_graded}
	Let $P$ be a graded poset with a least element.
	Then $\dimin{P}$ is graded.
\end{cor}
\begin{proof}
	Follows from Proposition \ref{prop:augmentation_preserves_graded} and the isomorphism between $P$ and $\augm{(\dimin{P})}$.
\end{proof}

\begin{dfn}[Cellular join of posets] \index{poset!cellular join} \index{$P \join Q$} \index{join!cellular}
	Let $P$, $Q$ be posets.
	The \emph{cellular join of $P$ and $Q$} is the poset $P \join Q \eqdef \dimin{(\augm{P} \times \augm{Q})}$.
\end{dfn}

\begin{comm}
	We call this the \emph{cellular join} to distinguish it from the usual join of two posets $P$, $Q$, which is the poset whose underlying set is $P + Q$ and partial order is defined by $x \leq y$ if either
	\begin{itemize}
		\item $x, y \in P$ and $x \leq y$ in $P$,
		\item $x, y \in Q$ and $x \leq y$ in $Q$, or
		\item $x \in P$ and $y \in Q$.
	\end{itemize}
\end{comm}

\begin{prop} \label{prop:join_monoidal_structure_posclos}
	The cellular join of posets extends to an essentially unique monoidal structure $(\poscat, \join, \varnothing)$ on $\poscat$ such that
	\[
		\dimin{(-)}\colon (\posbot, \times, 1) \to (\poscat, \join, \varnothing)
	\]
	is a strong monoidal functor.
	This monoidal structure restricts to a monoidal structure on $\posclos$.
\end{prop}
\begin{proof}
	Immediate from the fact that $\dimin{(-)}$ is an equivalence, and we can transport monoidal structures along equivalences.
	Moreover, $\dimin{1} = \varnothing$.
\end{proof}

\begin{lem} \label{lem:join_preserves_connected_colimits}
	Let $P$ be a poset.
	The functors $P \join -$ and $- \join P$ preserve all colimits of connected diagrams in $\posclos$.
\end{lem}
\begin{proof}
	Let $\fun{F}$ be a connected diagram in $\posclos$ and let $\gamma$ be a colimit cone under $\fun{F}$.
	Then $\augm{\gamma}$ is a colimit cone under a connected diagram in $\posclosbot$, so by Lemma \ref{lem:posclosbot_connected_colimits} it is also a colimit cone in $\posclos$.
	By Lemma \ref{lem:product_preserves_colimits}, $\augm{P} \times \augm{\gamma}$ and $\augm{\gamma} \times \augm{P}$ are colimit cones in $\posclos$, whose image factors through the inclusion $\posclosbot \incl \posclos$ by Lemma \ref{lem:products_restrict_to_posbot}.
	It follows from Lemma \ref{lem:posclosbot_connected_colimits} that they are also colimit cones in $\posclosbot$.
	Since $\dimin{(-)}$ preserves all colimits,
	\[
		P \join \gamma = \dimin{(\augm{P} \times \augm{\gamma})}, \quad \quad
		\gamma \join P = \dimin{(\augm{\gamma} \times \augm{P})}
	\]
	are colimit cones in $\posclos$.
\end{proof}

\begin{dfn}[Elements of the cellular join] \index{$\inj{x}, \inr{y}, x \join y$}
	Let $x \in P$ and $y \in Q$.
	We introduce the notation
	\[
		\inj{x} \eqdef (\augm{x}, \bot) \quad \quad 
		\inr{y} \eqdef (\bot, \augm{y}), \quad \quad
		x \join y \eqdef (\augm{x}, \augm{y})
	\]
	for elements of $P \join Q$.
\end{dfn}

\begin{lem} \label{lem:faces_of_poset_join}
	Let $P$, $Q$ be posets and let $z \in P \join Q$.
	Then
	\[
		\cofaces{}{}z = \begin{cases}
			\set{\inj{x'} \mid x' \in \cofaces{}{}x} +
			\set{x \join y' \mid y' \in \minel{Q}} &
			\text{if $z = \inj{x}$, $x \in P$}, \\
			\set{x' \join y \mid x' \in \minel{P}} + 
			\set{\inr{y'} \mid y' \in \cofaces{}{}y} &
			\text{if $z = \inr{y}$, $y \in Q$}, \\
			\set{x' \join y \mid x' \in \cofaces{}{}x} +
			\set{x \join y' \mid y' \in \cofaces{}{}y} &
			\text{if $z = x \join y$, $x \in P$, $y \in Q$.}
		\end{cases}
	\]
\end{lem}
\begin{proof}
	Follows from the definition, Lemma \ref{lem:product_of_posets_faces}, and Lemma 
	\ref{lem:cofaces_in_augmentation} by a simple case distinction.
\end{proof}

\begin{prop} \label{lem:join_preserves_graded}
	Let $P$, $Q$ be graded posets.
	Then
	\begin{enumerate}
		\item the cellular join $P \join Q$ is graded,
		\item for all $z \in P \join Q$,
		\[
			\dim{z} = 
			\begin{cases}
				\dim{x} & \text{if $z = \inj{x}$, $x \in P$,} \\
				\dim{y} & \text{if $z = \inr{y}$, $y \in Q$,} \\
				\dim{x} + \dim{y} + 1 & \text{if $z = x \join y$, $x \in P$, $y \in Q$}.
			\end{cases}
		\]
	\end{enumerate}
\end{prop}
\begin{proof}
	By Proposition \ref{prop:product_of_graded_is_graded}, Proposition \ref{prop:augmentation_preserves_graded}, and the isomorphism between $\augm{(P \join Q)}$ and $\augm{P} \times \augm{Q}$, we have that 
	\begin{enumerate}
		\item $\augm{(P \join Q)}$ is graded, so $P \join Q$ is also graded, 
		\item for all $x \in P$ and $y \in Q$,
	\[
	\begin{cases}
		\dim{(\inj{x})} + 1 = \dim{(\augm{x}, \bot)} = (\dim{x} + 1) + 0, \\
		\dim{(\inr{y})} + 1 = \dim{(\bot, \augm{y})} = 0 + (\dim{y} + 1), \\
		\dim{(x \join y)} + 1 = \dim{(\augm{x}, \augm{y})} = (\dim{x} + 1) + (\dim{y} + 1),
	\end{cases}
	\]
	\end{enumerate}
	and the statement follows.
\end{proof}

\begin{lem} \label{lem:preservation_of_closed_embeddings}
	Let $\imath\colon U \incl P$, $j\colon V \incl Q$ be closed embeddings of posets.
	Then
	\begin{enumerate}
		\item $\imath \times j\colon U \times V \to P \times Q$ is a closed embedding,
		\item $\augm{\imath}\colon \augm{U} \to \augm{P}$ is a closed embedding,
		\item if $P$ has a least element, $\dimin{\imath}\colon \dimin{U} \to \dimin{P}$ is a closed embedding,
		\item $\imath \join j\colon U \join V \to P \join Q$ is a closed embedding.
	\end{enumerate}
\end{lem}
\begin{proof}
	By inspection of the definitions.
\end{proof}

\begin{comm}
	By virtue of Lemma \ref{lem:preservation_of_closed_embeddings}, we can let all these operations act on closed subsets of a poset, and write 
	\[
		U \times V \subseteq P \times Q, \quad \quad \augm{U} \subseteq \augm{P}, \quad \quad U \join V \subseteq P \join Q
	\]
	for a given pair of closed subsets $U \subseteq P$, $V \subseteq Q$.
\end{comm}

\clearpage
\thispagestyle{empty}

%% file: ogposets.tex
\chapter{Oriented graded posets} \label{chap:ogposets}
\thispagestyle{firstpage}

\begin{guide}
	This chapter properly introduces oriented graded posets.
	As in the introduction to Chapter \ref{chap:order}, we define an \cemph{orientation} as an edge-labelling of the Hasse diagram of a graded poset with values in the set $\set{ +, - }$, but there are other equivalent descriptions, such as a bipartition 
	\[ \faces{}{}x = \faces{}{+}x + \faces{}{-}x\] 
	of the set of faces of each element $x$.

	An orientation on a face poset defines a ``flow'' between cells, \emph{from} an input face \emph{to} a cell, and \emph{from} a cell to an \emph{output} face, something that is made precise by the notion of \cemph{oriented Hasse diagram}, where the edges labelled with $-$ are reversed.

	Once an orientation is fixed, it singles out some cells in each dimension as being ``sources'' or ``sinks'' of the higher-dimensional flow, that is, having no output or no input cofaces.
	Such cells in dimension $n$ form the top-dimensional cells of the \cemph{input} and \cemph{output $n$\nbd boundary} of an oriented graded poset.
	At this stage, there is no relation \emph{a priori} between boundaries in different dimensions; this will no longer be the case once we focus on molecules and regular directed complexes.

	After we set up the basic definitions, we move on to defining a category of oriented graded posets.
	While at a later stage we will want to consider different kinds of morphisms, for now the natural definition is a ``rigid'' one, inducing bijections between the sets of input and output faces of an element and its image.
	As hoped, these morphisms have underlying closed order-preserving maps of posets, which are in addition dimension-preserving.

	Having a graded set together with boundary operators, the reader versed in homological algebra will no doubt expect to see a \cemph{chain complex} appear at some point.
	While not every oriented graded poset determines a chain complex, there is a simple yet powerful condition, called \cemph{oriented thinness}, which ensures that this is the case.
	It will be a theorem about regular directed complexes that they satisfy this condition, but it is quite natural to explore its immediate consequences already at this stage, and this what we will do in the last section.
\end{guide}


\section{Orientation and boundaries} \label{sec:orientation}

\begin{guide}
	In this section, we introduce the fundamental definitions and terminology relative to oriented graded posets, culminating in the definition of input and output $n$\nbd boundaries.
	We then prove a number of useful lemmas on boundaries and their interaction with set-theoretic operations.
\end{guide}

\begin{dfn}[Orientation on a graded poset] \index{graded poset!orientation}
Let $P$ be a graded poset.
An \emph{orientation} on $P$ is an edge-labelling of $\hasse{P}$ with values in $\set{ +, - }$.
\end{dfn}

\begin{dfn}[Sign algebra] \index{$\alpha, \beta, \ldots$}
We will use $\alpha, \beta, \ldots$ for variables ranging over $\set{ +, - }$.
We let $-\alpha$ be $-$ if $\alpha = +$ and $+$ if $\alpha = -$.
\end{dfn}

\begin{dfn}[Oriented graded poset] \index{graded poset!oriented|see {oriented graded poset}} \index{oriented graded poset}
An \emph{oriented graded poset} is a graded poset $P$ together with an orientation on $P$.
\end{dfn}

\begin{dfn}[Orientation induced on a closed subset] \index{oriented graded poset!induced orientation}
Let $P$ be an oriented graded poset and $U \subseteq P$ a closed subset.
The \emph{induced orientation on $U$} is the restriction of the orientation on $P$ to $\hasse{U}$.
This makes $U$ with the induced order an oriented graded poset.
\end{dfn}

\begin{comm}
We will implicitly assume that any closed subset of an oriented graded poset comes with the induced order and the induced orientation.
\end{comm}

\begin{dfn}[Input and output faces and cofaces] \index{$\faces{}{+}x, \faces{}{-}x$} \index{$\cofaces{}{+}x, \cofaces{}{-}x$}
Let $P$ be an oriented graded poset and $x \in P$.
The set of \emph{input faces} of $x$ is
\begin{equation*}
    \faces{}{-} x \eqdef \set{ y \in P \mid \text{$x$ covers $y$ with orientation $-$}  }
\end{equation*}
and the set of \emph{output faces} of $x$ is
\begin{equation*}
    \faces{}{+} x \eqdef \set{ y \in P \mid \text{$x$ covers $y$ with orientation $+$}  }.
\end{equation*}
Dually, the set of \emph{input cofaces} of $x$ is
\begin{equation*}
    \cofaces{}{-} x \eqdef \set{ y \in P \mid \text{$y$ covers $x$ with orientation $-$}  }
\end{equation*}
and the set of \emph{output cofaces} of $x$ is
\begin{equation*}
    \cofaces{}{+} x \eqdef \set{ y \in P \mid \text{$y$ covers $x$ with orientation $+$}  }.
\end{equation*}
We have $\faces{}{} x = \faces{}{+}x \cup \faces{}{-} x$ and $\cofaces{}{} x = \cofaces{}{+} x \cup \cofaces{}{-} x$.
\end{dfn}

\begin{dfn}[Oriented Hasse diagram] \index{oriented graded poset!oriented Hasse diagram} \index{$\hasseo{P}$} \index{oriented Hasse diagram!of an oriented graded poset}
Let $P$ be an oriented graded poset.
The \emph{oriented Hasse diagram of $P$} is the directed graph $\hasseo{P}$ whose
\begin{itemize}
    \item set of vertices is the underlying set of $P$, and
    \item set of edges is 
	    \[ \set{ (y, x) \mid \text{$y \in \faces{}{-}x$ or $x \in \faces{}{+}y$} }, \]
	with $s\colon (y, x) \mapsto y$ and $t\colon (y, x) \mapsto x$.
\end{itemize}
\end{dfn}

\begin{exm}[An oriented Hasse diagram]
	The oriented Hasse diagram for our example (\ref{eq:example_ogposet}) is
\[
	\input{img/oriented_hasse.tex}
\]
	with the ``reversed'' edges marked in a different colour for highlight.

	This redundant representation, with input edges both coloured differently and reversed, will be our preferred depiction of oriented graded posets in the rest of the book.
	We will also usually omit dimensions of vertices in such diagrams, that is, write just $k$ for $(n, k)$, letting the height of the vertex convey the same information.
\end{exm}

\begin{dfn}[Power set] \index{$\powerset{X}$}
	Let $X$ be a set.
	We let $\powerset{X}$ denote its power set.
\end{dfn}

\begin{prop} \label{prop:data_for_ogposets}
Let $P$ be an oriented graded poset.
Then $P$ can be uniquely reconstructed from any of the following data:
\begin{enumerate}[label=(\alph*)]
	\item the functions $\faces{}{-}, \faces{}{+}\colon P \to \powerset{P}$,
	\item the functions $\cofaces{}{-}, \cofaces{}{+}\colon P \to \powerset{P}$,
	\item the oriented Hasse diagram $\hasseo{P}$ together with the function $\dim{\!}\colon P \to \mathbb{N}$.
\end{enumerate}
\end{prop}
\begin{proof}
	Given the $\faces{}{\alpha}$ functions, we reconstruct $\hasse{P}$ together with its edge-labelling as the directed graph whose vertices are the elements of $P$, with an edge $y \to x$ if and only if
	\begin{itemize}
		\item $x \in \faces{}{+}y$, in which case we label the edge $+$, or
		\item $x \in \faces{}{-}y$, in which case we label the edge $-$.
	\end{itemize}
	The case of the $\cofaces{}{\alpha}$ functions is dual.

	Given $\hasseo{P}$ and the $\dim{}$ function, we reconstruct $\hasse{P}$ together with its edge-labelling as the directed graph whose vertices are the same as those of $\hasseo{P}$, with an edge $y \to x$ if and only if $\dim{y} = \dim{x} + 1$ and
	\begin{itemize}
		\item there is an edge $y \to x$ in $\hasseo{P}$, in which case we label the edge $+$, or
		\item there is an edge $x \to y$ in $\hasseo{P}$, in which case we label the edge $-$.\qedhere
	\end{itemize}
\end{proof}

\begin{prop} \label{prop:ogposet_from_faces}
	Let $P$ be a set, $\faces{}{-}, \faces{}{+}\colon P \to \powerset{P}$ be functions, and define inductively
	\begin{align*}
		\grade{0}{P} & \eqdef \set{x \in P \mid \faces{}{-}x \cup \faces{}{+}x = \varnothing}, \\
		\grade{n}{P} & \eqdef \set{x \in P \mid \text{$\faces{}{+}x \cup \faces{}{-}x \in \powerset{\grade{n-1}{P}} \setminus \set{\varnothing}$, $\faces{}{+}x \cap \faces{}{-}x = \varnothing$}}, \quad n > 0.
	\end{align*}
	The following are equivalent:
\begin{enumerate}[label=(\alph*)]
	\item $P = \sum_{n \in \mathbb{N}} \grade{n}{P}$,
	\item $P$ admits a structure of oriented graded poset whose functions of input and output faces are $\faces{}{-}$ and $\faces{}{+}$.
\end{enumerate}
\end{prop}
\begin{proof}
	Suppose that $P$ is an oriented graded poset whose input and output faces are given by $\faces{}{-}$ and $\faces{}{+}$.
	We show that $\grade{n}{P}$, as defined here, is equal to $\set{x \in P \mid \dim{x} = n}$ by induction on $n$.
	By Lemma \ref{lem:minimal_iff_dim_0}, we have $\dim{x} = 0$ if and only if $\faces{}{}x = \varnothing$.
	For $n > 0$, we have $\dim{x} = n$ if and only if every maximal path starting from $x$ in $\hasse{P}$ has length $n$, if and only if $\faces{}{}x$ is non-empty and every maximal path starting from $y \in \faces{}{}x$ has length $n-1$.
	Moreover, $\faces{}{+}x$ and $\faces{}{-}x$ are always disjoint, since every edge in $\hasse{P}$ is given a unique label by the orientation.

	Conversely, define $\hasse{P}$ to be the directed graph whose set of vertices is $P$, and whose set of edges is 
	\[
		\set{(y, x) \mid \text{$y \in P$, $x \in \faces{}{+}x \cup \faces{}{-}x$}}
	\]
	with $s\colon (y, x) \mapsto y$ and $t\colon (y, x) \mapsto x$.
	By assumption, for all $y \in P$, there is a unique $n \in \mathbb{N}$ such that $y \in \grade{n}{P}$.
	Then if there is an edge $y \to x$ in $\hasse{P}$, necessarily $x \in \grade{n-1}{P}$.
	It follows that, for all $x \neq y \in P$, if there is a path from $y$ to $x$ in $\hasse{P}$ and $x \in \grade{m}{P}$, then $m < n$, so in particular $\hasse{P}$ is acyclic, and every directed acyclic is the Hasse diagram of a poset.
	Moreover, if $y \in \grade{n}{P}$ and $n > 0$, by assumption there exists at least one edge $y \to x$, and if $y \in \grade{0}{P}$ then there is no edge with source $y$.
	It follows that every maximal path starting from $y$ in $\hasse{P}$ has length $n$, so $P$ is graded.
	Finally, for every edge $(y, x)$ in $\hasse{P}$, since $\faces{}{+}y \cap \faces{}{-}y = \varnothing$, there is a unique $\alpha \in \set{+, -}$ such that $x \in \faces{}{\alpha}y$.
	This determines a unique orientation on $P$ such that $\faces{}{+}$ and $\faces{}{-}$ are the functions of input and output faces.
\end{proof}

\begin{comm}
	Proposition \ref{prop:ogposet_from_faces} can be used to establish the relation between oriented graded posets, the \emph{directed precomplexes} of \cite{steiner1993algebra}, and the \emph{$\omega$\nbd hypergraphs} of \cite{forest2022unifying}.
	The last two are defined, up to differences in notation, as sets $P$ together with a grading $P = \sum_{n \in \mathbb{N}} \grade{n}{P}$ and functions $\faces{}{-}, \faces{}{+}\colon P \to \powerset{P}$ that are compatible with the grading.
	However, given $n > 0$ and $x \in \grade{n}{P}$,
	\begin{itemize}
		\item neither the definition of directed precomplexes, nor the definition of $\omega$\nbd hypergraphs, require that $\faces{}{+}x \cup \faces{}{-}x \neq \varnothing$ or that $\faces{}{+}x \cap \faces{}{-}x = \varnothing$;
		\item the definition of $\omega$\nbd hypergraphs requires that $\faces{}{+}x$ and $\faces{}{-}x$ be finite.
	\end{itemize}
	Thus an oriented graded poset can be always identified with a directed precomplex satisfying extra conditions.
	If $\faces{}{\alpha}x$ is finite for all $x \in P$ and $\alpha \in \set{+, -}$, then it can also be identified with an $\omega$\nbd hypergraph satisfying extra conditions.
\end{comm}

\begin{dfn}[Input and output $n$-boundaries] \index{oriented graded poset!boundaries} \index{$\faces{n}{+}U, \faces{n}{-}U$} \index{$\bound{n}{+}U, \bound{n}{-}U$}
Let $U$ be a closed subset of an oriented graded poset.
For all $\alpha \in \set{ +, - }$ and $n \in \mathbb{N}$, let
\begin{equation*}
    \faces{n}{\alpha} U \eqdef 
    \set{ x \in U_n \mid \cofaces{}{-\alpha} x \cap U = \varnothing  }.
\end{equation*}
For each $n \in \mathbb{N}$, the \emph{input $n$\nbd boundary} of $U$ is the closed subset
\begin{equation*}
    \bound{n}{-} U \eqdef \clos{(\faces{n}{-} U)} \cup \bigcup_{k < n} \clos{\grade{k}{(\maxel{U})}}
\end{equation*}
and the \emph{output $n$\nbd boundary} of $U$ is the closed subset
\begin{equation*}
    \bound{n}{+} U \eqdef \clos{(\faces{n}{+} U)} \cup \bigcup_{k < n} \clos{\grade{k}{(\maxel{U})}}.
\end{equation*}
For $n < 0$, we let $\faces{n}{\alpha}U = \bound{n}{\alpha}U \eqdef \varnothing$.
\end{dfn}

\begin{dfn}[Notation for boundaries] \index{$\bound{n}{+}x, \bound{n}{-}x$}
We will use the following notations, for $x$ an element in an oriented graded poset, $U$ a closed subset, $n \in \mathbb{N}$, and $\alpha \in \set{ +, - }$:
\begin{equation*}
    \bound{n}{\alpha} x \eqdef \bound{n}{\alpha} \clset{{ x }}, \quad \quad
    \bound{n}{} U \eqdef \bound{n}{-} U \cup \bound{n}{+} U, \quad \quad
    \faces{n}{} U \eqdef \faces{n}{-} U \cup \faces{n}{+} U.
\end{equation*}
\end{dfn}

\begin{dfn}[Boundary and interior] \index{oriented graded poset!boundary} \index{oriented graded poset!interior} \index{$\bound{}{} U$} \index{$\inter{U}$} \index{boundary!of an oriented graded poset}
Let $U$ be a closed subset of an oriented graded poset.
The \emph{boundary} of $U$ is the subset
\[
	\bound{}{} U \eqdef \bigcup_{n < \dim{U}} \bound{n}{}U
\]
and the \emph{interior} of $U$ is its complement
\[
	\inter{U} \eqdef U \setminus \bound{}{} U.
\]
\end{dfn}

\begin{exm}[Input and output boundaries]
The definition of input and output $n$\nbd boundaries is intended to combinatorially capture \emph{globular} boundary operators when $U$ is the oriented face poset of a pasting diagram.
Let $U$ be our running example (\ref{eq:example_ogposet}).
Then
\begin{align*}
	\faces{2}{-} U & = \faces{2}{+} U = \set { (2, 0) }, \\
	\faces{1}{-} U & = \set{ (1, 0), (1, 1), (1, 2) }, 
		       & \faces{1}{+} U & = \set{ (1, 3), (1, 2) }, \\
	\faces{0}{-} U & = \set{ (0, 0) },
		       & \faces{0}{+} U & = \set{ (0, 3) }, \\
	\maxel{U} & = \set{ (1, 2), (2, 0) },
\end{align*}
while $\faces{n}{\alpha} U = \varnothing$ if $n > 2$.
Consequently, $\bound{n}{\alpha} U = U$ for all $n \geq 2$, while
\begin{align*}
	\bound{1}{-} U & = \set{ (0, 0), (0, 1), (0, 2), (0, 3), (1, 0), (1, 1), (1, 2) }, \\
	\bound{1}{+} U & = \set{ (0, 0), (0, 2), (0, 3), (1, 3), (1, 2) },
\end{align*}
corresponding to the regions 
\[
\begin{tikzcd}[sep=small]
	{\color{\mycolor}{\scriptstyle 0}\;\bullet} && {\color{\mycolor} {\scriptstyle 2}\;\bullet} && {\color{\mycolor} {\scriptstyle 3}\;\bullet} \\
	& {\color{\mycolor} {\scriptstyle 1}\;\bullet}
	\arrow[""{name=0, anchor=center, inner sep=0}, "3", curve={height=-18pt}, from=1-1, to=1-3]
	\arrow["2", color=\mycolor, from=1-3, to=1-5]
	\arrow["0"', color=\mycolor, curve={height=6pt}, from=1-1, to=2-2]
	\arrow["1"', color=\mycolor, curve={height=6pt}, from=2-2, to=1-3]
	\arrow["0", shorten <=3pt, shorten >=6pt, Rightarrow, from=2-2, to=0]
\end{tikzcd} \quad \text{ and } \quad
\begin{tikzcd}[sep=small]
	{\color{\mycolor}{\scriptstyle 0}\;\bullet} && {\color{\mycolor} {\scriptstyle 2}\;\bullet} && {\color{\mycolor} {\scriptstyle 3}\;\bullet} \\
	& { {\scriptstyle 1}\;\bullet}
	\arrow[""{name=0, anchor=center, inner sep=0}, "3", color=\mycolor, curve={height=-18pt}, from=1-1, to=1-3]
	\arrow["2", color=\mycolor, from=1-3, to=1-5]
	\arrow["0"', curve={height=6pt}, from=1-1, to=2-2]
	\arrow["1"', curve={height=6pt}, from=2-2, to=1-3]
	\arrow["0", shorten <=3pt, shorten >=6pt, Rightarrow, from=2-2, to=0]
\end{tikzcd}
\]
of (\ref{eq:example_shape}), respectively.
Finally,
\[
	\bound{0}{-} U = \set{ (0, 0) }, \quad \bound{0}{+} U = \set{ (0, 3) }.
\]
\end{exm}

\begin{lem} \label{lem:dimension_of_boundary}
Let $U$ be a closed subset of an oriented graded poset, $n \in \mathbb{N}$, and $\alpha \in \set{ +, - }$.
Then $\dim{\bound{n}{\alpha} U} \leq n$.
\end{lem}
\begin{proof}
Let $x \in \bound{n}{\alpha} U$.
By definition there exists $y$ such that $x \leq y$ and either $y \in \faces{n}{\alpha}U$, so $\dim{y} = n$, or $y \in \grade{k}{(\maxel{U})}$, and $\dim{y} = k < n$.
In either case, by Lemma \ref{lem:dimension_is_monotonic}, $\dim{x} \leq \dim{y} \leq n$.
\end{proof}

\begin{lem} \label{lem:maximal_in_boundary}
Let $U$ be a closed subset of an oriented graded poset, $n \in \mathbb{N}$, and $\alpha \in \set{ +, - }$.
Then
\begin{enumerate}
    \item $\grade{n}{(\bound{n}{\alpha}U)} = \faces{n}{\alpha}U$,
    \item $\grade{k}{( \maxel{(\bound{n}{\alpha}U)} )} = \grade{k}{(\maxel{U})}$ for all $k < n$.
\end{enumerate}
\end{lem}
\begin{proof}
Let $x \in \bound{n}{\alpha}U$.
Then by definition there exists $y$ such that $x \leq y$ and either $y \in \faces{n}{\alpha}U$ or $y \in \grade{k}{(\maxel{U})}$ for some $k < n$.
If $x$ is maximal, necessarily $x = y$, and we obtain one inclusion.
The converse inclusions are evident.
\end{proof}

\begin{lem} \label{lem:maximal_vs_faces}
Let $U$ be a closed subset of an oriented graded poset, $n \in \mathbb{N}$, and $\alpha \in \set{ +, - }$.
Then 
\begin{enumerate}
    \item $\grade{n}{(\maxel{U})} = \faces{n}{+} U \cap \faces{n}{-} U$,
    \item if $n = \dim{U}$, then $\grade{n}{(\maxel{U})} = \faces{n}{\alpha} U = \grade{n}{U}$.
\end{enumerate}
\end{lem}
\begin{proof}
Let $x \in U$, $\dim{x} = n$.
Then $x$ is maximal if and only if it has no cofaces in $U$, if and only if $\cofaces{}{-\alpha}x \cap U = \cofaces{}{\alpha}x \cap U = \varnothing$, if and only if $x \in \faces{n}{+} U \cap \faces{n}{-} U$.
If $n = \dim{U}$, then every element of $\grade{n}{U}$ is maximal in $U$, so 
\[
	\grade{n}{U} = \grade{n}{(\maxel{U})} \subseteq \faces{n}{\alpha} U \subseteq \grade{n}{U}
\]
using the first part of the proof, and we conclude that they are all equal.
\end{proof}

\begin{lem} \label{lem:boundary_inclusion}
Let $U$ be a closed subset of an oriented graded poset, $n \in \mathbb{N}$, and $\alpha \in \set{ +, - }$.
Then
\begin{enumerate}
    \item $\bound{n}{\alpha} U \subseteq U$,
    \item $\bound{n}{\alpha}U = U$ if and only if $n \geq \dim{U}$.
\end{enumerate}
\end{lem}
\begin{proof}
By definition, $\faces{n}{\alpha}U \subseteq U$ and $\grade{k}{(\maxel{U})} \subseteq U$ for all $k < n$.
Because $U$ is closed, by Lemma \ref{lem:closure_is_monotonic} it also contains their closures.
This proves that $\bound{n}{\alpha} U \subseteq U$.

Suppose $n < \dim{U}$.
By Lemma \ref{lem:dimension_of_boundary}, $\dim{\bound{n}{\alpha} U} \leq n < \dim{U}$, so $U \neq \bound{n}{\alpha} U$.
Conversely, suppose $n \geq \dim{U}$.
By Lemma \ref{lem:closure_of_maximal} and Lemma \ref{lem:closure_of_union},
\begin{equation*}
    U = \clos{\maxel{U}} = \bigcup_{k \leq \dim{U}} \clos{\grade{k}{(\maxel{U})}}.
\end{equation*}
If $n > \dim{U}$, this is included in (hence equal to) $\bound{n}{\alpha} U$.
If $n = \dim{U}$, we use Lemma \ref{lem:maximal_vs_faces} to rewrite this as
\[    
\clos{(\faces{n}{\alpha} U)} \cup \bigcup_{k < n} \clos{\grade{k}{(\maxel{U})}},
\]
which is equal to $\bound{n}{\alpha} U$.
\end{proof}

\begin{cor} \label{cor:dimension_from_boundary}
Let $U$ be a closed subset of an oriented graded poset. Then 
\begin{equation*}
	\dim{U} = 
	\begin{cases}
		\min \set{ n \geq -1 \mid \bound{n}{+} U = \bound{n}{-} U = U } & \text{if defined}, \\
		\infty & \text{otherwise}.
	\end{cases}
\end{equation*}
\end{cor}

\begin{lem} \label{lem:faces_of_union}
Let $U, V$ be closed subsets of an oriented graded poset, $n \in \mathbb{N}$, and $\alpha \in \set{ +, - }$.
Then
\begin{enumerate}
    \item $\maxel{(U \cup V)} = (\maxel{U} \cap \maxel{V}) + (\maxel{U} \setminus V) + (\maxel{V} \setminus U)$,
    \item $\faces{n}{\alpha}(U \cup V) = (\faces{n}{\alpha}U \cap \faces{n}{\alpha}V) + (\faces{n}{\alpha}U \setminus V) + (\faces{n}{\alpha}V \setminus U)$.
\end{enumerate}
\end{lem}
\begin{proof}
Follows straightforwardly from the definitions using the decomposition $U \cup V = (U \cap V) + (U \setminus V) + (V \setminus U)$.
\end{proof}

\begin{cor} \label{cor:boundary_of_union}
Let $U, V$ be closed subsets of an oriented graded poset, $n \in \mathbb{N}$, and $\alpha \in \set{ +, - }$.
Then $\bound{n}{\alpha}(U \cup V) \subseteq \bound{n}{\alpha}U \cup \bound{n}{\alpha}V$.
\end{cor}

\begin{lem} \label{lem:faces_intersection}
Let $V \subseteq U$ be closed subsets of an oriented graded poset, $n \in \mathbb{N}$, and $\alpha \in \set{ +, - }$.
Then
\begin{enumerate}
    \item $V \cap \faces{n}{\alpha}U \subseteq \faces{n}{\alpha}V$,
    \item $V \cap \grade{n}{(\maxel{U})} \subseteq \grade{n}{(\maxel{V})}$.
\end{enumerate}
\end{lem}
\begin{proof}
Let $x \in V \cap \faces{n}{\alpha}U$.
Then $\cofaces{}{-\alpha}x \cap V \subseteq \cofaces{}{-\alpha}x \cap U = \varnothing$, so $x \in \faces{n}{\alpha}V$.
By Lemma \ref{lem:maximal_vs_faces} we have 
\begin{equation*}
    V \cap \grade{n}{(\maxel{U})} = V \cap \faces{n}{+}U \cap \faces{n}{-}U \subseteq \faces{n}{+}V \cap \faces{n}{-}V
\end{equation*}
by the first part, and we conclude.
\end{proof}

\begin{lem} \label{lem:boundary_included_in_subset}
	Let $V \subseteq U$ be closed subsets of an oriented graded poset, $n \in \mathbb{N}$, and $\alpha \in \set{ + , - }$.
	If $\bound{n}{\alpha}U \subseteq V$, then $\bound{n}{\alpha}U \subseteq \bound{n}{\alpha}V$.
\end{lem}
\begin{proof}
	By Lemma \ref{lem:faces_intersection}, under the assumption that $\bound{n}{\alpha}U \subseteq V$, we have
	\[
		\faces{n}{\alpha}U = V \cap \faces{n}{\alpha}U \subseteq \faces{n}{\alpha}V
	\]
	and, for all $k < n$,
	\[
		\grade{k}{(\maxel{U})} = V \cap \grade{k}{(\maxel{U})} \subseteq \grade{k}{(\maxel{V})},
	\]
	and we conclude.
\end{proof}


\section{The category of oriented graded posets} \label{sec:category_ogpos}

\begin{guide}
	In this section, we define a notion of morphism of oriented graded posets.
	While for regular directed complexes it will be beneficial to extend this notion in more than one way, the choice we make here suits the relatively unstructured nature of general oriented graded posets, and while somewhat ``rigid'', it ensures the existence of those colimits that we need the most, namely, coproducts and gluings (pushouts of monomorphisms).
\end{guide}

\begin{dfn}[Morphism of oriented graded posets] \index{oriented graded poset!morphism} \index{morphism!of oriented graded posets}
Let $P, Q$ be oriented graded posets.
A \emph{morphism} $f\colon P \to Q$ is a function of their underlying sets which, for all $x \in P$ and $\alpha \in \set{ +, - }$, induces a bijection between $\faces{}{\alpha} x$ and $\faces{}{\alpha} f(x)$.
\end{dfn}

\begin{dfn}[The category $\ogpos$] \index{$\ogpos$}
We let $\ogpos$ denote the category whose objects are oriented graded posets and morphisms are morphisms of oriented graded posets.
\end{dfn}

\begin{lem} \label{lem:properties_of_morphisms}
Let $f\colon P \to Q$ be a morphism of oriented graded posets.
Then
\begin{enumerate}
	\item $f$ is order-preserving,
 	\item $f$ is closed,
	\item $f$ is dimension-preserving, that is, for all $x \in P$, $\dim{f(x)} = \dim{x}$.
\end{enumerate}
\end{lem}
\begin{proof}
Let $x, y \in P$ with $x \leq y$.
We proceed by induction on $\dim{y} - \dim{x}$, which is $\geq 0$ by Lemma \ref{lem:dimension_is_monotonic}.
If $\dim{y} = \dim{x}$, then $x = y$, so $f(x) = f(y)$.
Otherwise, by Lemma \ref{lem:if_locally_finite_height_then_covers} there exists $y'$ such that $x \leq y'$ and $y' \in \faces{}{} y$.
Then 
\[
	\dim{y'} - \dim{x} = (\dim{y} - 1) - \dim{x},
\] 
so $f(x) \leq f(y')$ by the inductive hypothesis, while $f(y') \in f(\faces{}{} y) = \faces{}{} f(y)$, hence $f(y') < f(y)$, by the definition of morphism.

Next, let $x \in P$ and $y \in \clset{f(x)}$, that is, $y \leq f(x)$.
We will prove that $y \in f(\clset{x})$ by induction on $\dim{f(x)} - \dim{y}$.
If $\dim{f(x)} = \dim{y}$ then $y = f(x)$ and we are done.
Otherwise, there exists $y'$ such that $y \leq y'$ and $y' \in \faces{}{} f(x)$.
Then $y' \in f(\faces{}{} x)$, that is, there exists $x' \in \faces{}{} x$ such that $y' = f(x')$.
Moreover, $\dim{f(x')} - \dim{y} = (\dim{f(x)} - 1) - \dim{y}$, so by the inductive hypothesis $y \in f(\clset{x'}) \subseteq f(\clset{x})$.
By Lemma \ref{lem:map_closed_iff_maps_lowersets_to_lowersets} this proves that $f$ is closed.

Finally, let $x \in P$; we will prove that $\dim{f(x)} = \dim{x}$ by induction on $\dim{x}$.
If $\dim{x} = 0$, then $\faces{}{} x = \varnothing$, so $\faces{}{} f(x) = \varnothing$.
It follows that there are no non-trivial paths in $\hasse{Q}$ starting from $f(x)$, hence, $\dim{f(x)} = 0$.
If $\dim{x} > 0$, the set $\faces{}{} x$ is non-empty and its elements have dimension $\dim{x} - 1$.
By the inductive hypothesis, for all $x' \in \faces{}{} x$, we have 
\[	\dim{f(x')} = \dim{x'} = \dim{x} - 1. \]
Since $f(x') \in \faces{}{}f(x)$, we conclude that $\dim{f(x)} = \dim{f(x')} + 1 = \dim{x}$.
\end{proof}

\begin{cor} \label{cor:underlying_poset_functor}
Forgetting the orientation determines a faithful functor
	\[ \fun{U}\colon \ogpos \to \posclos. \]
\end{cor}

\begin{cor} \label{cor:nskeleton_is_a_functor}
	Let $f\colon P \to Q$ be a morphism of oriented graded posets, $n \in \mathbb{N}$.
	Then
	\begin{enumerate}
		\item $f(\skel{n}{P}) \subseteq \skel{n}{Q}$,
		\item $f \mapsto \restr{f}{\skel{n}{P}}$ determines an endofunctor $\skel{n}{}$ on $\ogpos$.
	\end{enumerate}
\end{cor}

\begin{prop} \label{prop:hasse_graph_homomorphism}
Let $f\colon P \to Q$ be a morphism of oriented graded posets.
Then $f$ induces homomorphisms
\[	
	\hasse{f}\colon \hasse{P} \to \hasse{Q} \quad  \text{and} \quad 
	\hasseo{f}\colon \hasseo{P} \to \hasseo{Q}.
\]
These assignments determine functors $\hasse{}, \hasseo{}\colon \ogpos \to \gph$.
\end{prop}
\begin{proof}
	Let $x, y \in P$ and suppose there is an edge from $y$ to $x$ in $\hasse{P}$.
	Then $x \in \faces{}{}y$, hence $f(x) \in \faces{}{}f(y)$, hence there is an edge from $f(y)$ to $f(x)$ in $\hasse{Q}$.
	The case of oriented Hasse diagrams is similarly straightforward, as is functoriality of the assignments.
\end{proof}

\begin{dfn}[Flow preorder] \index{oriented graded poset!flow preorder}
	Let $P$ be an oriented graded poset.
	The \emph{flow preorder} $\precflow$ on $P$ is defined by $x \precflow y$ if and only there is a path from $x$ to $y$ in $\hasseo{P}$.
\end{dfn}

\begin{cor} \label{cor:morphisms_preserve_flow}
Let $f\colon P \to Q$ be a morphism of oriented graded posets, $x, y \in P$.
If $x \precflow y$, then $f(x) \precflow f(y)$.
\end{cor}
\begin{proof}
	By Proposition \ref{prop:hasse_graph_homomorphism}, $\hasseo{f}$ maps a path from $x$ to $y$ in $\hasseo{P}$ to a path from $f(x)$ to $f(y)$ in $\hasseo{Q}$.
\end{proof}

\begin{dfn}[Inclusion of oriented graded posets] \index{oriented graded poset!inclusion}
An \emph{inclusion} is an injective morphism of oriented graded posets.
\end{dfn}

\begin{lem} \label{lem:properties_of_inclusions}
Let $\imath\colon P \incl Q$ be an inclusion of oriented graded posets.
Then 
\begin{enumerate}
	\item $\imath$ is order-reflecting,
	\item $\imath$ reflects input and output faces, that is, if $\imath(x) \in \faces{}{\alpha}\imath(y)$, then $x \in \faces{}{\alpha}y$.
\end{enumerate}
\end{lem}
\begin{proof}
	The first fact is an immediate consequence of Lemma \ref{lem:closed_embedding_is_order_reflecting}, as $\fun{U}\imath$ is a closed embedding of posets.
	For the second one, suppose $\imath(x) \in \faces{}{\alpha}\imath(y)$.
	Then there exists a unique $x' \in \faces{}{\alpha}y$ such that $\imath(x) = \imath(x')$.
	Because $\imath$ is injective, $x = x'$, and we conclude.
\end{proof}

\begin{lem} \label{lem:characterisation_of_isomorphisms}
Let $f\colon P \to Q$ be a morphism of oriented graded posets.
The following are equivalent:
\begin{enumerate}[label=(\alph*)]
	\item $f$ is a surjective inclusion;
	\item $f$ is an isomorphism of oriented graded posets.
\end{enumerate}
\end{lem}
\begin{proof}
	If $f$ is an isomorphism of oriented graded posets, it induces an isomorphism of their underlying posets, hence a bijection on the underlying sets.

	Conversely, suppose that $f$ is a surjective inclusion.
	Then $\fun{U}f$ is a surjective closed embedding, hence an isomorphism of the underlying posets, which admits an inverse $\invrs{f}$.
	Since by Lemma \ref{lem:properties_of_inclusions} $f$ both preserves and reflects input and output faces, $\invrs{f}$ does too.
	We conclude that $\invrs{f}$ lifts to a morphism of oriented graded posets, inverse to $f$.
\end{proof}

\begin{prop} \label{prop:inclusions_preserve_faces}
Let $\imath\colon P \incl Q$ be an inclusion of oriented graded posets and $U \subseteq P$ a closed subset.
For all $n \in \mathbb{N}$ and $\alpha \in \set{ +, - }$,
\begin{enumerate}
	\item $\imath(\faces{n}{\alpha}U) = \faces{n}{\alpha} \imath(U)$, 
	\item $\imath(\grade{n}{ (\maxel{U}) }) = \grade{n}{ (\maxel{\imath(U)}) }$.
\end{enumerate}
\end{prop}
\begin{proof}
	Let $n \in \mathbb{N}$ and $\alpha \in \set{+, -}$.
	Given $x \in \imath(\faces{n}{\alpha} U)$, there is a unique $x' \in \faces{n}{\alpha} U$ such that $x = \imath(x')$.
	By Lemma \ref{lem:properties_of_inclusions}, $\dim{x} = \dim{x'} = n$.
	Suppose that there exists $y \in \cofaces{}{-\alpha}x \cap \imath(U)$.
	Then $y = \imath(y')$ for a unique $y' \in U$, and $\imath(x') \in \faces{}{-\alpha} \imath(y')$.
	By Lemma \ref{lem:properties_of_inclusions}, $x' \in \faces{}{-\alpha} y'$, that is, $y' \in \cofaces{}{-\alpha}x' \cap U$, a contradiction.
	It follows that $x \in \faces{n}{\alpha} \imath(U)$.
	The converse is analogous.
	
	Since $\imath$ is injective, direct images under $\imath$ commute with intersections, so by Lemma \ref{lem:maximal_vs_faces}
	\begin{align*}
		\imath(\grade{n}{ (\maxel{U}) }) & = \imath(\faces{n}{+} U \cap \faces{n}{-} U) = \imath(\faces{n}{+} U) \cap \imath(\faces{n}{-} U) = \\
						 & = \faces{n}{+} \imath(U) \cap \faces{n}{-} \imath(U) = \grade{n}{ (\maxel{\imath(U)}) }. \qedhere
	\end{align*}
\end{proof}

\begin{cor} \label{cor:inclusions_preserve_boundaries}
Let $\imath\colon P \incl Q$ be an inclusion of oriented graded posets and $U \subseteq P$ a closed subset.
For all $n \in \mathbb{N}$ and $\alpha \in \set{ +, - }$, 
\[
	\imath(\bound{n}{\alpha}U) = \bound{n}{\alpha}\imath(U).
\]
\end{cor}
\begin{proof}
	Follows from Proposition \ref{prop:inclusions_preserve_faces} together with the fact that $\imath$ is closed by Lemma \ref{lem:properties_of_morphisms}, so images under $\imath$ commute with closures.
\end{proof}

\begin{dfn}[Local embedding] \index{morphism!local embedding}
A morphism $f\colon P \to Q$ of oriented graded posets is a \emph{local embedding} if its underlying map of posets is a local embedding.
\end{dfn}

\begin{dfn}[The category $\ogposloc$] \index{$\ogposloc$}
	We let $\ogposloc$ denote the wide subcategory of $\ogpos$ whose morphisms are local embeddings.
\end{dfn}

\begin{exm}[A non-injective local embedding of regular directed complexes] \index[counterex]{A non-injective local embedding of regular directed complexes}
	If $P$ is the oriented face poset of the pasting diagram of two arrows
\[
	\begin{tikzcd}[sep=small]
	{{\scriptstyle 0}\;\bullet} && {{\scriptstyle 1}\;\bullet} && {{\scriptstyle 2}\;\bullet}
	\arrow["0", from=1-1, to=1-3]
	\arrow["1", from=1-3, to=1-5]
	\end{tikzcd} 
	\quad \quad
	\input{img/noninjective_domain.tex}
\]	
	and $Q$ the oriented face poset of a ``loop'' diagram with two arrows
\[
	\begin{tikzcd}[sep=small]
	{{\scriptstyle 0}\;\bullet} && {{\scriptstyle 1}\;\bullet}
	\arrow["0"', curve={height=12pt}, from=1-1, to=1-3]
	\arrow["1"', curve={height=12pt}, from=1-3, to=1-1]
	\end{tikzcd}
	\quad \quad
	\input{img/noninjective_codomain.tex}
\]
	then the morphism which identifies $(0, 0)$ and $(0, 2)$, that is,
	\[
		(n, k) \mapsto \begin{cases}
			(0, 0) & \text{if $(n, k) = (0, 2)$}, \\
			(n, k) & \text{otherwise},
		\end{cases}
	\]
	is a surjective local embedding $P \to Q$, but it is not an inclusion.
\end{exm}

\begin{exm}[A morphism of oriented graded posets which is not a local embedding] \index[counterex]{A morphism of oriented graded posets which is not a local embedding}
	It will be a non-trivial result (Corollary \ref{cor:morphisms_of_rdcpx_are_local_isomorphisms}) that \emph{every} morphism of regular directed complexes is a local embedding.
	To see that this is not obvious, we exhibit a simple counterexample to this property when one steps out of the class of regular directed complexes.
	Let $P$ and $Q$ be the oriented graded posets (with uniform orientation)
	\[
		\input{img/nonlocemb_1.tex} \quad \quad \quad \input{img/nonlocemb_2.tex}
	\]
	respectively.
	There is a surjective morphism $f\colon P \to Q$ which identifies $(0, 0)$ and $(0, 1)$, and sends $(n, k)$ to $(n, k)$ for all other $n, k$.
	This is not a local embedding, as it is evidently not an isomorphism on $\clset{(2, 0)} = P$.
\end{exm}

\begin{lem} \label{lem:local_embeddings_lift_orientation}
	Let $P$ be a poset, $Q$ an oriented graded poset, and let $f\colon P \to \fun{U}Q$ be a local embedding.
	Then there exists a unique orientation on $P$ such that $f$ lifts to a local embedding of oriented graded posets.
\end{lem}
\begin{proof}
	Since $\fun{U}Q$ is graded, by Lemma \ref{lem:local_embedding_into_graded} so is $P$.
	Let $x \in P$.
	By Lemma \ref{lem:local_embedding_preserves_faces}, $f$ induces a bijection between $\faces{}{}x$ and $\faces{}{}f(x)$.
	Then, for all $y \in \faces{}{}x$ and $\alpha \in \set{+, -}$, we let $y \in \faces{}{\alpha}x$ if and only if $f(y) \in \faces{}{\alpha}f(x)$.
	This defines an orientation on $P$ such that $f$ is, by construction, a local embedding of oriented graded posets.
	This is also the only possible orientation with this property, since with any other choice we would have $y \in \faces{}{\alpha}x$ but $f(y) \in \faces{}{-\alpha}f(x)$ for some $x, y \in P$ and $\alpha \in \set{ +, - }$.
\end{proof}

\begin{prop} \label{prop:pullbacks_of_inclusions}
	The category $\ogpos$ has pullbacks of inclusions, and they are both preserved and reflected by $\fun{U}\colon \ogpos \to \posclos$.
	Moreover, inclusions are stable under pullbacks.
\end{prop}
\begin{proof}
	Let $f\colon P \to Q$ be a morphism and $\imath\colon U \incl Q$ an inclusion of oriented graded posets.
	By Lemma \ref{lem:properties_of_morphisms} $\fun{U}f$ is a closed order-preserving map and $\fun{U}\imath$ a closed embedding of posets.
	By Lemma \ref{lem:pullbacks_of_closed_embeddings}, the pullback of $\fun{U}\imath$ along $\fun{U}f$ exists, and produces a span $f'\colon V \to \fun{U}U$, $j\colon V \incl \fun{U}P$ where $j$ is a closed embedding.
	By Lemma \ref{lem:local_embeddings_lift_orientation}, there is a unique orientation on $V$ such that $j$ lifts to an inclusion of oriented graded posets.
	Given $x \in V$, we then have $\imath \after \restr{f'}{\faces{}{}x} = f\after\restr{j}{\faces{}{}x}$; since $\imath$, $f$, and $j$ all preserve input and output faces, so does $f'$.
	
	This proves that the pullback square lifts to a square of morphisms in $\ogpos$.
	Given another span of morphisms $g\colon W \to P$, $h\colon W \to U$ such that $f\after g = \imath \after h$, it suffices to prove that the universal map $k\colon \fun{U}W \to \fun{U}V$ obtained in $\posclos$ lifts to a morphism of oriented graded posets.
	For all $x \in W$, we have $j \after \restr{k}{\faces{}{}x} = \restr{h}{\faces{}{}x}$.
	Since both $j$ and $h$ preserve input and output faces, we conclude that $k$ does too.
\end{proof}

\begin{lem} \label{lem:reflected_colimits_in_ogpos}
	Let $\fun{F}\colon \smcat{C} \to \ogpos$ be a diagram of inclusions of oriented graded posets, and suppose $\gamma$ is a colimit cone under $\fun{UF}$ whose components are all closed embeddings.
	Then there exists a unique colimit cone $\vec{\gamma}$ under $\fun{F}$ such that $\gamma = \fun{U}\vec{\gamma}$, whose components are all inclusions.
\end{lem}
\begin{proof}
	Let $P$ be the colimit of $\fun{UF}$, and let $x \in P$.
	By the construction of colimits in $\posclos$, there exist an object $c$ in $\smcat{C}$ and $y \in \fun{F}c$ such that $x = \gamma_c(y)$.
	Since $\gamma_c$ is a closed embedding, it restricts to an isomorphism between $\clset{y}$ and $\clset{x}$.
	There is a unique orientation on $\clset{x}$ that lifts this to an isomorphism of oriented graded posets.

	Suppose that there exist another $c'$ in $\smcat{C}$ and $y' \in \fun{F}c'$ such that $\gamma_{c'}(y') = x$.
	Then there exist a zig-zag of morphisms
\[\begin{tikzcd}
	c \equiv c_0 & {c_1} & {c_2} & \ldots & {c_{m-2}} & {c_{m-1}} & {c_m \equiv c'} 
	\arrow["{f_0}", from=1-2, to=1-1]
	\arrow["{f_1}"', from=1-2, to=1-3]
	\arrow[from=1-4, to=1-3]
	\arrow[from=1-4, to=1-5]
	\arrow["{f_{m-1}}", from=1-6, to=1-5]
	\arrow["{f_m}"', from=1-6, to=1-7]
\end{tikzcd}\]
in $\smcat{C}$, and a sequence $(y_i \in \fun{F}c_i)_{i=0}^m$ such that $y_0 = y$, $y_m = y'$, and for all $k \in \set{1, \ldots, \frac{m}{2}}$,
\[
	y_{2k-2} = \fun{F}f_{2k-2}(y_{2k-1}), \quad \quad y_{2k} = \fun{F}f_{2k-1}(y_{2k-1});
\]
note that, according to the construction of colimits, it is the images through the $\gamma_{c_i}$ of the two sides of each equation that need to be equal, but since all the $\gamma_{c_i}$ are injective, we can lift the equations to their domains.

Because by assumption all the $\fun{F}f_i$ are inclusions, they induce isomorphisms between all the $\clset{y_i}$.
It follows that the orientation transported to $\clset{x}$ from $\clset{y'}$ coincides with the orientation transported from $\clset{y}$.
Since $x$ was arbitrary, we conclude that we can give $P$ a unique orientation such that all the $\gamma_c$ are inclusions of oriented graded posets.
These form the cone $\vec{\gamma}$ under $\fun{F}$ in $\ogpos$.

It remains to show that $\vec{\gamma}$ is a colimit cone.
Let $\eta$ be another cone under $\fun{F}$ with tip $Q$.
Then the universal property of $\gamma$ produces a unique order-preserving map $f\colon \fun{U}P \to \fun{U}Q$; it suffices to show that $f$ lifts to a morphism of oriented graded posets.
By the same reasoning as before, for all $x \in P$, there exist $c$ in $\smcat{C}$ and $y \in \fun{F}c$ such that $x = \gamma_c(y)$.
Then $f \after \restr{\gamma_c}{\faces{}{}y} = \restr{\eta_c}{\faces{}{}y}$, and since both $\eta_c$ and $\gamma_c$ preserve input and output faces, so does $f$.
\end{proof}

\begin{prop} \label{prop:pushouts_in_ogpos}
The category $\ogpos$ has
\begin{enumerate}
    \item a strict initial object $\varnothing$,
    \item pushouts of inclusions along inclusions,
\end{enumerate}
preserved and reflected by $\fun{U}\colon \ogpos \to \posclos$.
Moreover,
\begin{enumerate}
	\item the pushout of an inclusion along an inclusion is an inclusion,
	\item a pushout square of inclusions is also a pullback square.
\end{enumerate}
\end{prop}
\begin{proof}
	The colimit cone exhibiting the empty poset as an initial object of $\posclos$ trivially satisfies the conditions of Lemma \ref{lem:reflected_colimits_in_ogpos}, so it lifts to $\ogpos$.
	It is a strict initial object in $\ogpos$ because it is in $\posclos$.

	Given a span of inclusions in $\ogpos$, its underlying span in $\posclos$ is a span of closed embeddings.
	By Lemma \ref{lem:pushouts_of_spans_of_closed_embeddings}, its pushout is both a pushout square and a pullback square of closed embeddings.
	By Proposition \ref{prop:pullbacks_of_inclusions} and Lemma \ref{lem:reflected_colimits_in_ogpos}, it lifts to a square of inclusions which is both a pushout and a pullback.
\end{proof}

\begin{prop} \label{prop:diagram_of_inclusions}
Every oriented graded poset $P$ is the colimit of the diagram of inclusions
\begin{equation*}
	\fun{U}P \to \ogpos, \quad \quad (x \leq y) \mapsto (\clset{{ x }} \incl \clset{{ y }}).
\end{equation*}
\end{prop}
\begin{proof}
	Consider the obvious cone whose component at $x$ is $\clset{x} \incl P$.
	The underlying cone in $\posclos$ is a colimit cone whose components are all closed embeddings.
	We conclude by Lemma \ref{lem:reflected_colimits_in_ogpos}.
\end{proof}

\begin{prop} \label{prop:factorisation_in_ogpos}
Let $f\colon P \to Q$ be a morphism of oriented graded posets.
Then $f$ factors as
\begin{enumerate}
	\item a surjective morphism $\widehat{f}\colon P \to f(P)$,
	\item followed by an inclusion $\imath\colon f(P) \incl Q$.
\end{enumerate}
This factorisation is unique up to unique isomorphism.
\end{prop}
\begin{proof}
By Lemma \ref{lem:properties_of_morphisms}, $f$ has an underlying closed order-preserving map of posets.
By Proposition \ref{prop:em_factorisation_in_posclos}, this factors up to unique isomorphism as a surjective closed order-preserving map $\widehat{f}\colon \fun{U}P \to f(P)$ followed by a closed embedding $\imath\colon f(P) \incl \fun{U}Q$.

By Lemma \ref{lem:local_embeddings_lift_orientation}, there is a unique orientation on $f(P)$ that makes $\imath$ an inclusion of oriented graded posets.
Then, for all $x \in P$ and $\alpha \in \set{+, -}$, we have 
\[
	\restr{\widehat{f}}{\faces{}{\alpha}x} = \invrs{(\restr{\imath}{\faces{}{\alpha}\widehat{f}(x)})}\after \restr{f}{\faces{}{\alpha}x},
\]
which is a bijection because both its factors are.
Thus $\widehat{f}$ is a morphism, and the factorisation lifts to $\ogpos$.
\end{proof}

\begin{cor} \label{cor:ofs_on_ogpos}
The classes of
\begin{enumerate}
	\item surjective morphisms,
	\item inclusions
\end{enumerate}
form an orthogonal factorisation system on $\ogpos$.
\end{cor}
\begin{proof}
Both classes are evidently closed under composition and contain all isomorphisms.
The statement then follows from Proposition \ref{prop:factorisation_in_ogpos}.
\end{proof}


\section{Oriented thin graded posets} \label{sec:oriented_thin}

\begin{guide}
	\emph{Thinness}, also known as the \emph{diamond property}, is a condition on posets associated with their order complex being a combinatorial manifold.
	In this section, we show that if thinness of the underlying poset is complemented by a constraint on orientations, then there is a canonical \cemph{augmented chain complex} associated to an oriented graded poset (Proposition \ref{prop:chain_complex_is_functorial}).
\end{guide}

\begin{dfn}[Positive least element] \index{oriented graded poset!positive least element}
	Let $P$ be an oriented graded poset, $\bot \in P$.
	We say that $\bot$ is a \emph{positive least element of $P$} if
	\begin{enumerate}
		\item $\bot$ is the least element of $P$,
		\item $\cofaces{}{}\bot = \cofaces{}{+}\bot$.
	\end{enumerate}
\end{dfn}

\begin{comm}
	There is, clearly, a dual notion of negative least element, that we will not take into consideration.
\end{comm}

\begin{dfn}[The category $\ogposbot$] \index{$\ogposbot$}
	We let $\ogposbot$ denote the full subcategory of $\ogpos$ on oriented graded posets with a positive least element.
\end{dfn}

\begin{lem} \label{lem:underlying_poset_of_ogp_with_positive_least_element}
	The restriction of $\fun{U}\colon \ogpos \to \posclos$ to $\ogposbot$ factors through the inclusion of subcategories $\posclosbot \incl \posclos$.
\end{lem}
\begin{proof}
	By definition, if $P$ has a positive least element, then $\fun{U}P$ has a least element.
	Let $f\colon P \to Q$ be a morphism in $\ogposbot$, let $\bot_P$ be the least element of $P$, and let $\bot_Q$ be the least element of $Q$.
	By Lemma \ref{lem:minimal_iff_dim_0} $\dim{\bot_P} = \dim{\bot_Q} = 0$.
	Let $x \in P$.
	Since $f$ is dimension-preserving, if $f(x) = \bot_Q$, then $\dim{x} = \dim{f(x)} = 0$.
	But $\bot_P$ is the only minimal element of $P$, so it is the only 0\nbd dimensional element of $P$.
	We conclude that $x = \bot_P$, so $f$ reflects the least element.
\end{proof}

\begin{prop} \label{prop:ogpos_equivalent_ogposbot}
	There exists a unique pair of functors
\[
		\augm{(-)}\colon \ogpos \to \ogposbot, \quad \quad
		\dimin{(-)}\colon \ogposbot \to \ogpos
\]
	such that 
	\begin{enumerate}
		\item the diagram of functors
			\begin{equation}
			\begin{tikzcd} \label{eq:diagram_of_augmentations}
			\ogpos && \ogposbot && \ogpos \\
			\posclos && \posclosbot && \posclos
			\arrow["{\augm{(-)}}", from=1-1, to=1-3]
			\arrow["{\fun{U}}", from=1-1, to=2-1]
			\arrow["{\augm{(-)}}", from=2-1, to=2-3]
			\arrow["{\fun{U}}", from=1-3, to=2-3]
			\arrow["{\dimin{(-)}}", from=1-3, to=1-5]
			\arrow["{\fun{U}}", from=1-5, to=2-5]
			\arrow["{\dimin{(-)}}", from=2-3, to=2-5]
			\end{tikzcd}
			\end{equation}
			commutes,
		\item for all oriented graded posets $P$, $x \in P$, and $\alpha \in \set{+, -}$,
			\[
				\cofaces{}{\alpha}\augm{x} = \set{\augm{y} \mid y \in \cofaces{}{\alpha}x}, 
			\]
		\item the two functors are inverse to each other up to natural isomorphism.
	\end{enumerate}
\end{prop}
\begin{proof}
	The requirement that (\ref{eq:diagram_of_augmentations}) commute fixes what the functors do on the underlying posets and closed maps.
	Let $P$ be an oriented graded poset.
	By Proposition \ref{prop:augmentation_preserves_graded}, $\augm{P}$ is graded, and by Lemma \ref{lem:cofaces_in_augmentation} the orientation on $\augm{P}$ is fixed uniquely by
	\begin{itemize}
		\item $\cofaces{}{\alpha}\augm{x} = \set{\augm{y} \mid y \in \cofaces{}{\alpha}x}$ for all $x \in P$ and $\alpha \in \set{+, -}$,
		\item $\cofaces{}{+}\bot = \cofaces{}{}\bot$, since $\bot$ must be a positive least element.
	\end{itemize}
	Finally, if $P$ is an oriented graded poset with a positive least element, then $\dimin{P}$ is a graded poset by Corollary \ref{cor:diminution_preserves_graded}, and the requirement that $\augm{\dimin{P}}$ be isomorphic to $P$ fixes the orientation on $\dimin{P}$ uniquely.
	
	It is straightforward to check that, with these orientations, $\augm{(-)}$ and $\dimin{(-)}$ are well-defined on morphisms.
	Then Lemma \ref{lem:characterisation_of_isomorphisms} implies that the natural isomorphisms that exhibit the equivalence between $\posclos$ and $\posclosbot$ lift to natural isomorphisms that exhibit an equivalence between $\ogpos$ and $\ogposbot$.
\end{proof}

\begin{dfn}[Augmentation of an oriented graded poset] \index{oriented graded poset!augmentation} \index{$\augm{P}$} \index{augmentation!of an oriented graded poset}
	Let $P$ be an oriented graded poset.
	The \emph{augmentation of $P$} is the oriented graded poset $\augm{P}$ obtained as the image of $P$ through the functor $\augm{(-)}\colon \ogpos \to \ogposbot$.
\end{dfn}

\begin{dfn}[Diminution of an oriented graded poset with positive least element] \index{oriented graded poset!diminution} \index{$\dimin{P}$} \index{diminution!of an oriented graded poset}
	Let $P$ be an oriented graded poset with a positive least element.
	The \emph{diminution of $P$} is the oriented graded poset $\dimin{P}$ obtained as the image of $P$ through the functor $\dimin{(-)}\colon \ogposbot \to \ogpos$.
\end{dfn}

\begin{dfn}[Thin graded poset] \index{graded poset!thin}
	Let $P$ be a graded poset with a least element.
	We say that $P$ is \emph{thin} if, for all $x, y \in P$ such that $x \leq y$ and $\codim{x}{y} = 2$, the interval $[x, y]$ has exactly 4 elements, that is, it is of the form
\[\begin{tikzcd}[column sep=small]
	& y \\
	{z_1} && {z_2} \\
	& x
	\arrow[from=1-2, to=2-1]
	\arrow[from=1-2, to=2-3]
	\arrow[from=2-1, to=3-2]
	\arrow[from=2-3, to=3-2]
\end{tikzcd}\]
	for exactly two elements $z_1, z_2$. 
\end{dfn}

\begin{comm}
The property of being thin is also known as the \emph{diamond property}, most commonly in the theory of abstract polytopes.
\end{comm}

\begin{dfn}[Oriented thin graded poset] \index{oriented graded poset!oriented thin}
	Let $P$ be an oriented graded poset with a positive least element.
	We say that $P$ is \emph{oriented thin} if, for all $x, y \in P$ such that $x \leq y$ and $\codim{x}{y} = 2$, the interval $[x, y]$ is of the form
\[\begin{tikzcd}[column sep=small]
	& y \\
	{z_1} && {z_2} \\
	& x
	\arrow["{\alpha}"', from=1-2, to=2-1]
	\arrow["{\beta}", from=1-2, to=2-3]
	\arrow["{\gamma}"', from=2-1, to=3-2]
	\arrow["{-\alpha\beta\gamma}", from=2-3, to=3-2]
\end{tikzcd}\]
for exactly two elements $z_1, z_2$, and for some $\alpha, \beta, \gamma \in \set{+, -}$.
\end{dfn}

\begin{rmk}
	Evidently, if $P$ is oriented thin, then its underlying graded poset is thin.
\end{rmk}

\begin{dfn}[The category $\otgpos$] \index{$\otgpos$}
	We let $\otgpos$ denote the full subcategory of $\ogposbot$ on the oriented thin graded posets.
\end{dfn}

\begin{comm}
	Oriented thinness is stated as a property of objects of $\ogposbot$, but via the equivalence of Proposition \ref{prop:ogpos_equivalent_ogposbot} it can be seen as a property of objects of $\ogpos$, translating via the pair $\augm{(-)}$, $\dimin{(-)}$.
\end{comm}

\begin{dfn}[Augmented chain complex] \index{chain complex} \index{chain complex!augmented} \index{augmented chain complex|see {chain complex}}
	An \emph{augmented chain complex} $C$ is a chain complex of abelian groups in non-negative degree
\[\begin{tikzcd}
	\ldots & {\chain{C}{n}} & {\chain{C}{n-1}} & \ldots & {\chain{C}{1}} & {\chain{C}{0}}
	\arrow["\der", from=1-1, to=1-2]
	\arrow["\der", from=1-2, to=1-3]
	\arrow["\der", from=1-3, to=1-4]
	\arrow["\der", from=1-4, to=1-5]
	\arrow["\der", from=1-5, to=1-6]
\end{tikzcd}\]
	together with a homomorphism $\eau\colon \chain{C}{0} \to \mathbb{Z}$ satisfying $\eau \after \der = 0$.
\end{dfn}

\begin{dfn}[Homomorphism of augmented chain complexes] \index{chain complex!homomorphism}
	Let $C$, $D$ be augmented chain complexes.
	A \emph{homomorphism} $f\colon C \to D$ is a sequence 
	\[
		(\chain{f}{n}\colon \chain{C}{n} \to \chain{D}{n})_{n \in \mathbb{N}}
	\]
	of homomorphisms of abelian groups such that the diagrams
\[\begin{tikzcd}
	{\chain{C}{n}} && {\chain{C}{n-1}} \\
	{\chain{D}{n}} && {\chain{D}{n-1}}
	\arrow["\der", from=1-1, to=1-3]
	\arrow["{\chain{f}{n}}", from=1-1, to=2-1]
	\arrow["{\chain{f}{n-1}}", from=1-3, to=2-3]
	\arrow["\der", from=2-1, to=2-3]
\end{tikzcd}
\quad \quad 
\begin{tikzcd}
	{\chain{C}{0}} \\
	{\chain{D}{0}} && {\mathbb{Z}}
	\arrow["{\chain{f}{0}}", from=1-1, to=2-1]
	\arrow["\eau", from=2-1, to=2-3]
	\arrow["\eau", curve={height=-12pt}, from=1-1, to=2-3]
\end{tikzcd}\]
	commute, for $n$ ranging in $\posnat$.
\end{dfn}

\begin{dfn}[The category $\chaug$] \index{$\chaug$}
	We let $\chaug$ denote the category whose objects are augmented chain complexes and morphisms are homomorphisms of augmented chain complexes.
\end{dfn}

\begin{dfn}[Augmented chain complex of an oriented thin graded poset] \index{chain complex!of an oriented thin graded poset} \index{$\freeab{P}$}
	Let $P$ be an oriented graded poset such that $\augm{P}$ is oriented thin.
	The \emph{augmented chain complex of $P$} is the diagram
\[\begin{tikzcd}
	\ldots & {\freeab{\grade{n}{P}}} & {\freeab{\grade{n-1}{P}}} & \ldots & {\freeab{\grade{1}{P}}} & {\freeab{\grade{0}{P}}}
	\arrow["\der", from=1-1, to=1-2]
	\arrow["\der", from=1-2, to=1-3]
	\arrow["\der", from=1-3, to=1-4]
	\arrow["\der", from=1-4, to=1-5]
	\arrow["\der", from=1-5, to=1-6]
\end{tikzcd}\]
of homomorphisms of abelian groups, where $\freeab{\grade{n}{P}}$ is the free abelian group on the set $\grade{n}{P}$ and, for each $n > 0$, the homomorphism $\der\colon \freeab{\grade{n}{P}} \to \freeab{\grade{n-1}{P}}$ is defined on the generators $x \in \grade{n}{P}$ by
\begin{equation} \label{eq:boundary_in_chain_complex}
	x \mapsto \sum_{y \in \faces{}{+}x} y - \sum_{y \in \faces{}{-}x} y,
\end{equation}
together with the homomorphism $\eau\colon \freeab{\grade{0}{P}} \to \mathbb{Z}$ defined on the generators $x \in \grade{0}{P}$ by $x \mapsto 1$.
\end{dfn}

\begin{lem} \label{lem:augmented_chain_complex_is_indeed_one}
	Let $P$ be an oriented graded poset such that $\augm{P}$ is oriented thin.
	Then the augmented chain complex of $P$ is well-defined as an augmented chain complex $\freeab{P}$.
\end{lem}
\begin{proof}
	It suffices to show that $\eau \after \der = 0$ and that $\der \after \der = 0\colon \freeab{\grade{n}{P}} \to \freeab{\grade{n-2}{P}}$ for all $n > 1$.
	Observe that there exist isomorphisms
	\[
		\mathbb{Z} \iso \freeab{\grade{0}{(\augm{P})}}, \quad \quad \freeab{\grade{n}{P}} \iso \freeab{\grade{n+1}{(\augm{P})}}
	\]
	for all $n \in \mathbb{N}$, through which both $\eau$ and all the $\der$ can be seen as defined by the expression (\ref{eq:boundary_in_chain_complex}) on elements of $\augm{P}$.
	Now, $\eau \after \der$ and $\der \after \der$ are defined, on elements of $\augm{P}$, by
	\[
		x \mapsto \sum_{y \in \faces{}{+}x} \sum_{z \in \faces{}{+}y} z + \sum_{y \in \faces{}{-}x} \sum_{z \in \faces{}{-}y} z - \sum_{y \in \faces{}{+}x}\sum_{z \in \faces{}{-}y} z - \sum_{y \in \faces{}{-}x} \sum_{z \in \faces{}{+}z} z.
	\]
	By oriented thinness of $\augm{P}$, any element $z$ that appears in this expression appears exactly twice with opposite signs, so the two cancel out, and the expression evaluates to 0.
\end{proof}

\begin{prop} \label{prop:chain_complex_is_functorial}
	Let $f\colon P \to Q$ be a morphism of oriented graded posets such that $\augm{P}$ and $\augm{Q}$ are oriented thin.
	Then the sequence of homomorphisms
	\begin{align*}
		\freeab{\grade{n}{f}}\colon \freeab{\grade{n}{P}} & \to \freeab{\grade{n}{Q}}, \\
		x \in \grade{n}{P} & \mapsto f(x) \in \grade{n}{Q}
	\end{align*}
	is a homomorphism $\freeab{f}\colon \freeab{P} \to \freeab{Q}$ of augmented chain complexes.
	This assignment determines a functor $\freeab{-}\colon \otgpos \to \chaug$.
\end{prop}
\begin{proof}
	Let $x \in \grade{0}{P}$.
	Then $\eau(\freeab{\grade{0}{f}}(x)) = \eau(f(x)) = 1 = \eau(x)$, so $\eau \after \freeab{\grade{0}{f}} = \eau$.
	Next, let $x \in \grade{n}{P}$ for $n > 0$.
	Then
	\begin{align*}
		\der (\freeab{\grade{n}{f}}(x)) 
		& = \der (f(x))
		= \;\; \sum_{\mathclap{y \in \faces{}{+}f(x)}} \; y \; - 
		\;\; \sum_{\mathclap{y \in \faces{}{-}f(x)}}\; y \; = \\
		& = \;\; \sum_{\mathclap{y' \in \faces{}{+}x}}\; f(y') \; - 
		\;\; \sum_{\mathclap{y' \in \faces{}{-}x}} \; f(y') \; = 
		\freeab{\grade{n-1}{f}}(\der x),
	\end{align*}
	using the fact that $f$ determines bijections between $\faces{}{\alpha}x$ and $\faces{}{\alpha}f(x)$ for all $\alpha \in \set{+, -}$.
	This proves that $\der \after \freeab{\grade{n}{f}} = \freeab{\grade{n-1}{f}} \after \der$, so 
	$(\freeab{\grade{n}{f}})_{n \in \mathbb{N}}$ is a homomorphism of augmented chain complexes.
	Functoriality is straightforward.
\end{proof}

\clearpage 
\thispagestyle{empty}

%% file: molecules.tex
\chapter{Molecules} \label{chap:molecules}
\thispagestyle{firstpage}

\begin{guide}
Not all oriented graded posets are the oriented face posets of shapes of pasting diagrams, and even when they are, it is not guaranteed that the shape can be univocally reconstructed from its oriented face poset.
In the introduction to the first chapter, we identified some conditions under which the latter ought to be possible, namely, that
\begin{enumerate}
	\item the underlying graded poset is the face poset of a regular cell complex, and
	\item both the input and the output boundary of each $n$\nbd cell are closed topological $(n-1)$\nbd balls.
\end{enumerate}
Moreover, the orientation of cells must be compatible with the whole diagram, and all its boundaries, admitting expressions as iterated composites in a strict $n$\nbd category.
In this chapter, we will define a class of oriented graded posets that satisfies these properties, and whose members we call \cemph{molecules}.
Molecules will be our combinatorial notion of ``shape of a pasting diagram'' throughout the book.

The class of molecules is defined inductively, as the smallest subclass containing $1$, the \cemph{point}, that is, the shape of a 0\nbd cell
\[
	\begin{tikzcd}[sep=small]
	{{\scriptstyle 0}\;\bullet}
	\end{tikzcd}
\]
and closed under two constructions.

In the first, we take two molecules $U$ and $V$ together with an isomorphism $\varphi$ between $\bound{k}{+}U$ and $\bound{k}{-}V$ for some natural number $k$, and form a new molecule $U \cp{k} V$ by \cemph{pasting} $U$ and $V$ along this isomorphism, which we express formally as a pushout in $\ogpos$.
While this appears to depend on $\varphi$, we will later prove that when $\varphi$ exists, it is unique.
This operation is a model of strict $n$\nbd categorical composition, whose equations it satisfies up to unique isomorphism, and on oriented face posets of pasting diagrams it behaves as you would expect: for example, if $\globe{2}$ is the 2\nbd dimensional \cemph{globe} (Section \ref{sec:globes_and_thetas})
\begin{equation} \label{eq:2globe}
\begin{tikzcd}[sep=small]
	{{\scriptstyle 0}\;\bullet} && {{\scriptstyle 1}\;\bullet}
	\arrow[""{name=0, anchor=center, inner sep=0}, "0"', curve={height=18pt}, from=1-1, to=1-3]
	\arrow[""{name=1, anchor=center, inner sep=0}, "1", curve={height=-18pt}, from=1-1, to=1-3]
	\arrow["0"', shorten <=5pt, shorten >=5pt, Rightarrow, from=0, to=1]
\end{tikzcd}
	\quad \quad \quad
	\input{img/globe2.tex}
\end{equation}
then $\globe{2} \cp{0} \globe{2}$ and $\globe{2} \cp{1} \globe{2}$ are both defined, and are, respectively, the oriented face posets
\begin{equation} \label{eq:globe2_cp0_globe2}
	\begin{tikzcd}[sep=small]
	{{\scriptstyle 0}\;\bullet} && {{\scriptstyle 1}\;\bullet} && {{\scriptstyle 2}\;\bullet}
	\arrow[""{name=0, anchor=center, inner sep=0}, "0"', curve={height=18pt}, from=1-1, to=1-3]
	\arrow[""{name=1, anchor=center, inner sep=0}, "2", curve={height=-18pt}, from=1-1, to=1-3]
	\arrow[""{name=2, anchor=center, inner sep=0}, "1"', curve={height=18pt}, from=1-3, to=1-5]
	\arrow[""{name=3, anchor=center, inner sep=0}, "3", curve={height=-18pt}, from=1-3, to=1-5]
	\arrow["0"', shorten <=5pt, shorten >=5pt, Rightarrow, from=0, to=1]
	\arrow["1"', shorten <=5pt, shorten >=5pt, Rightarrow, from=2, to=3]
\end{tikzcd}
	\quad \quad
	\input{img/globe2_cp0_globe2.tex}
	\quad ,
\end{equation}
\begin{equation} \label{eq:globe2_cp1_globe2}
\begin{tikzcd}[sep=small]
	{{\scriptstyle 0}\;\bullet} && {{\scriptstyle 1}\;\bullet}
	\arrow[""{name=0, anchor=center, inner sep=0}, "0"', curve={height=24pt}, from=1-1, to=1-3]
	\arrow[""{name=1, anchor=center, inner sep=0}, "2", curve={height=-24pt}, from=1-1, to=1-3]
	\arrow[""{name=2, anchor=center, inner sep=0}, "1"'{pos=0.3}, from=1-1, to=1-3]
	\arrow["0"', shorten <=3pt, shorten >=3pt, Rightarrow, from=0, to=2]
	\arrow["1"', shorten <=3pt, shorten >=3pt, Rightarrow, from=2, to=1]
\end{tikzcd}
	\quad \quad
	\input{img/globe2_cp1_globe2.tex}
	\quad .
\end{equation}
In the second construction, we take two molecules $U$ and $V$ of the same dimension $n$, together with an isomorphism $\varphi$ between $\bound{}{} U$ and $\bound{}{} V$, that restricts to isomorphisms between $\bound{}{\alpha}U$ and $\bound{}{\alpha} V$ for all $\alpha \in \set{ + , - }$ (that is, $U$ and $V$ have the ``same'' input and output boundary).
Moreover, we require that $U$ and $V$ are both \cemph{round}, which is a constraint on the intersections of their input and output boundaries, which conceptually corresponds to, and in a precise sense implies (Proposition \ref{prop:order_complex_of_round_molecule}), $U$ and $V$ being closed topological $n$\nbd balls.
We then form a new \cemph{atom} $U \celto V$, a molecule with no non-trivial pasting decomposition, by
\begin{enumerate} 
	\item first gluing $U$ and $V$ along $\varphi$, which is also a pushout in $\ogpos$ --- you should picture this as the construction of an $n$\nbd sphere from two $n$\nbd balls glued along their boundaries, 
	\item then adding a new $(n+1)$\nbd dimensional element whose input boundary is $U$ and output boundary is $V$ --- you should picture this as ``filling'' the $n$\nbd sphere to produce an $(n+1)$\nbd ball.
\end{enumerate}
Just as it was for pasting, $\varphi$ turns out to be unique when it exists, so the construction does not depend on the choice of an isomorphism.

For example, the \cemph{arrow} $\thearrow{}$ is obtained as $1 \celto 1$, and the 2\nbd dimensional globe as $\thearrow{} \celto \thearrow{}$.
We may also now express the oriented face poset (\ref{eq:example_ogposet}) of the shape (\ref{eq:example_shape}) as
\[
	( ( \thearrow{} \cp{0} \thearrow{} ) \celto \thearrow{} ) \cp{0} \thearrow{}.
\]
The first part of the chapter leads to the definition of molecules.
The second part is devoted to proving ``rigidity'' properties of molecules: most importantly, that when two molecules are isomorphic, they are isomorphic in a unique way (Corollary 
\ref{cor:molecule_unique_isomorphism}).
This will allow us to be relaxed about the distinction between isomorphism and equality of molecules.
\end{guide}


\section{Pastings and globularity} \label{sec:pastings}

\begin{guide}
	A characteristic property of cells in a strict $\omega$\nbd category is \cemph{globularity}, ensuring that the strict $\omega$\nbd category has an underlying $\omega$\nbd graph or globular set.
	This is a condition that binds together the input and output boundaries in different dimensions, implying that the iterated operation of taking lower-dimensional input or output boundaries can be reduced to a single instance.
	Correspondingly, molecules will satisfy a version of globularity, relative to their intrinsic notion of boundary.

	In this section, we define pastings of oriented graded posets, which will be one of the constructors for the class of molecules.
	The main result is Lemma \ref{lem:pasting_preserves_globularity}: pastings preserve the property of globularity.
\end{guide}

\begin{dfn}[Pasting construction] \index{oriented graded poset!pasting} \index{$U \cpiso{k}{\varphi} V$} \index{pasting!at the $k$-boundary}
Let $U$, $V$ be oriented graded posets, $k \in \mathbb{N}$, and let $\varphi\colon \bound{k}{+}U \incliso \bound{k}{-}V$ be an isomorphism.
The \emph{pasting of $U$ and $V$ at the $k$\nbd boundary along $\varphi$} is the oriented graded poset $U \cpiso{k}{\varphi} V$ obtained as the pushout
\[\begin{tikzcd}
	\bound{k}{+}U & \bound{k}{-}V & V \\
	U && U \cpiso{k}{\varphi} V
	\arrow["\varphi", hook, from=1-1, to=1-2]
	\arrow[hook, from=1-2, to=1-3]
	\arrow[hook', from=1-1, to=2-1]
	\arrow["\imath_U", hook, from=2-1, to=2-3]
	\arrow["\imath_V", hook', from=1-3, to=2-3]
	\arrow["\lrcorner"{anchor=center, pos=0.125, rotate=180}, draw=none, from=2-3, to=1-1]
\end{tikzcd}\]
in $\ogpos$.
\end{dfn}

\begin{dfn}[Splitting] \index{oriented graded poset!splitting}
Let $U$ be an oriented graded poset.
A \emph{splitting} of $U$ is an ordered pair of closed subsets $V$, $W \subseteq U$ such that $V \cup W = U$ and $V \cap W = \bound{k}{+}V = \bound{k}{-}W$ for some $k \in \mathbb{N}$.
Given such a splitting of $U$, we say that $U$ \emph{splits into $V \cup W$ along the $k$\nbd boundary}.
\end{dfn}

\begin{lem} \label{lem:pasting_produces_splitting}
Let $U$, $V$ be oriented graded posets, $k \in \mathbb{N}$, and suppose $U \cpiso{k}{\varphi} V$ is defined.
Then $U \cpiso{k}{\varphi} V$ splits into $\imath_U(U) \cup \imath_V(V)$ along the $k$\nbd boundary.
\end{lem}
\begin{proof}
	An immediate consequence of Corollary \ref{cor:inclusions_preserve_boundaries} together with Proposition \ref{prop:pushouts_in_ogpos}.
\end{proof}

\begin{comm}
We will regularly identify $U$ and $V$ with their isomorphic images, and say that $U \cpiso{k}{\varphi} V$ splits into $U \cup V$ along the $k$\nbd boundary.
\end{comm}

\begin{lem} \label{lem:pasting_top_boundary}
Let $U$, $V$ be oriented graded posets, $k \in \mathbb{N}$, and suppose $U \cpiso{k}{\varphi} V$ is defined.
Then
\begin{enumerate}
    \item $\bound{k}{-}(U \cpiso{k}{\varphi} V) = \bound{k}{-}U$,
    \item $\bound{k}{+}(U \cpiso{k}{\varphi} V) = \bound{k}{+}V$.
\end{enumerate}
\end{lem}
\begin{proof}
Since $\bound{k}{-}V = U \cap V \subseteq U$, we have $\faces{k}{-}V \subseteq U$ and $\grade{j}{(\maxel{V})} \subseteq U$ for all $j < k$.
It follows from Lemma \ref{lem:faces_of_union} that 
\[\grade{j}{( \maxel{(U \cup V)} )} \subseteq \grade{j}{ (\maxel{U}) }, \quad \quad 
\faces{k}{-}(U \cup V) \subseteq \faces{k}{-}U,\] 
so $\bound{k}{-}(U \cup V) \subseteq \bound{k}{-}U$.

Conversely, suppose $x \in \bound{k}{-}U$.
Then there exists $y$ such that $x \leq y$ and $y \in \faces{k}{-}U$ or $y \in \grade{j}{(\maxel{U})}$ for some $j < k$.

Suppose that $y \in \faces{k}{-}U$. 
If $y \notin V$ then $y \in \faces{k}{-}(U \cup V)$ by Lemma \ref{lem:faces_of_union}.
If $y \in V$ then $y \in \grade{k}{(U \cap V)}$ which by Lemma \ref{lem:maximal_in_boundary} is equal to $\faces{k}{-}V$. 
It follows that $y \in \faces{k}{-}U \cap \faces{k}{-}V$, and by Lemma \ref{lem:faces_of_union} $y \in \faces{k}{-}(U \cup V)$.

Suppose that $y \in \grade{j}{(\maxel{U})}$ with $j < k$.
By Lemma \ref{lem:maximal_in_boundary},
\begin{equation*}
    \grade{j}{(\maxel{U})} = \grade{j}{(\maxel{(\bound{k}{+}U)})} = \grade{j}{(\maxel{(\bound{k}{-}V)})} = \grade{j}{(\maxel{V})},
\end{equation*}
and by Lemma \ref{lem:faces_of_union} $y \in \grade{j}{(\maxel{(U \cup V)})}$.
In either case $x, y \in \bound{k}{-}(U \cup V)$, so $\bound{k}{-}(U \cup V) = \bound{k}{-}U$.
The proof that $\bound{k}{+}(U \cup V) = \bound{k}{+}V$ is dual.
\end{proof}

\begin{dfn}[Globularity] \index{oriented graded poset!globular}
Let $U$ be an oriented graded poset.
We say that $U$ is \emph{globular} if, for all $k, n \in \mathbb{N}$ and $\alpha, \beta \in \set{ +, - }$, if $k < n$ then
\begin{equation*}
    \bound{k}{\alpha}(\bound{n}{\beta}U) = \bound{k}{\alpha}U.
\end{equation*}
\end{dfn}

\begin{rmk}
	By Lemma \ref{lem:boundary_inclusion}, this is only non-trivial when $n < \dim{U}$.
\end{rmk}

\begin{exm}[An oriented graded poset which is not globular] \index[counterex]{An oriented graded poset which is not globular}
	As you would expect, all molecules will turn out to be globular, so it is easy to come up with examples.
	Instead, let $U$ be the oriented graded poset (not a molecule, nor a regular directed complex)
\[
	\input{img/nonglobular.tex} \; ;
\]
then both $\bound{1}{-}U = \clset{(1, 0)}$ and $\bound{1}{+}U = \clset{(1, 1)}$ are isomorphic to the arrow $\thearrow{}$, but $U$ is \emph{not} globular, since
\begin{align*}
	\bound{0}{-}(\bound{1}{-}U) & = \set{(0, 0)} \subsetneq \bound{0}{-} U = \set{(0, 0), (0, 2)} \supsetneq \bound{0}{-}(\bound{1}{+}U) = \set{(0, 2)}, \\
	\bound{0}{+}(\bound{1}{-}U) & = \set{(0, 1)} \subsetneq \bound{0}{-} U = \set{(0, 1), (0, 3)} \supsetneq \bound{0}{+}(\bound{1}{+}U) = \set{(0, 3)}.
\end{align*}
\end{exm}

\begin{lem} \label{lem:boundaries_of_globular_are_globular}
Let $U$ be a globular oriented graded poset, $n \in \mathbb{N}$, and $\beta \in \set{ +, - }$.
Then $\bound{n}{\alpha}U$ is globular.
\end{lem}
\begin{proof}
Let $k < m$ be natural numbers and $\alpha, \gamma \in \set{ +, - }$.
If $m < n$, using globularity of $U$ twice,
\begin{equation*}
    \bound{k}{\alpha}(\bound{m}{\gamma}(\bound{n}{\beta}U)) = \bound{k}{\alpha}(\bound{m}{\gamma}U) = \bound{k}{\alpha}U = \bound{k}{\alpha}(\bound{n}{\beta}U).
\end{equation*}
If $m \geq n$, by Lemma \ref{lem:boundary_inclusion} we have $\bound{m}{\gamma}(\bound{n}{\beta}U) = \bound{n}{\beta}U$, so
\begin{equation*}
    \bound{k}{\alpha}(\bound{m}{\gamma}(\bound{n}{\beta}U)) = \bound{k}{\alpha}(\bound{n}{\beta}U). \qedhere 
\end{equation*}
\end{proof}

\begin{lem} \label{lem:globular_finite_dim_boundary}
Let $U$ be a globular oriented graded poset, and suppose $n \eqdef \dim{U} < \infty$.
Then $\bound{}{}U = \bound{n-1}{}U$.
\end{lem}
\begin{proof}
	One inclusion is obvious.
	For the other, notice that for all $k < n-1$ and all $\alpha \in \set{+, -}$, we have $\bound{k}{\alpha}U = \bound{k}{\alpha}(\bound{n-1}{\alpha}U) \subseteq \bound{n-1}{\alpha}U$, hence $\bound{k}{}U \subseteq \bound{n-1}{}U$.
\end{proof}

\begin{dfn}[Input and output boundary, the globular case]
Let $U$ be a finite-dimensional globular oriented graded poset, $n \eqdef \dim{U}$.
For all $\alpha \in \set{+, -}$, we write $\bound{}{\alpha}U \eqdef \bound{n-1}{\alpha}U$.
\end{dfn}

\begin{rmk}
By Lemma \ref{lem:globular_finite_dim_boundary}, this makes it so $\bound{}{}U = \bound{}{-}U \cup \bound{}{+}U$.
\end{rmk}

\begin{lem} \label{lem:globular_interiors}
Let $U$ be a globular oriented graded poset, $x \in U$.
Then either
\begin{itemize}
	\item $x \in \inter{U}$, or 
	\item there exist $n < \dim{U}$ and $\alpha \in \set{+, -}$ such that $x \in \inter{\bound{n}{\alpha}U}$.
\end{itemize}
\end{lem}
\begin{proof}
	Either $x \in \inter{U}$ or $x \in \bound{}{}U$.
	In the first case, we are done.
	In the other, there exist $n < \dim{U}$ and $\alpha \in \set{+, -}$ such that $x \in \bound{n}{\alpha}U$.
	We may assume that $n = \dim{\bound{n}{\alpha}U}$; otherwise, supposing $m \eqdef \dim{\bound{n}{\alpha}U} < n$, we have
	\[
		\bound{m}{\alpha}U = \bound{m}{\alpha}(\bound{n}{\alpha}U) = \bound{n}{\alpha}U
	\]
	by globularity and Lemma \ref{lem:boundary_inclusion}.

	Now we proceed by induction on $n$.
	If $n = 0$, then $\bound{n}{\alpha}U = \inter{\bound{n}{\alpha}U}$ and we are done.
	Otherwise, either $x \in \inter{\bound{n}{\alpha}U}$ or there exist $k < n$ and $\beta \in \set{+, -}$ such that $x \in \bound{k}{\beta}(\bound{n}{\alpha}U) = \bound{k}{\beta}U$.
	In the first case, we are done, in the second we can apply the inductive hypothesis.
\end{proof}

\begin{lem} \label{lem:pasting_lower_boundary}
Let $U$, $V$ be globular oriented graded posets, $k \in \mathbb{N}$, and suppose $U \cpiso{k}{\varphi} V$ is defined.
For all $j < k$ and $\alpha \in \set{ +, - }$,
\begin{equation*}
    \bound{j}{\alpha}U = \bound{j}{\alpha}V = \bound{j}{\alpha}(U \cpiso{k}{\varphi} V).
\end{equation*}
\end{lem}
\begin{proof}
The first equality follows from globularity by
\begin{equation*}
    \bound{j}{\alpha}U = \bound{j}{\alpha}(\bound{k}{+}U) = \bound{j}{\alpha}(\bound{k}{-}V) = \bound{j}{\alpha}V.
\end{equation*}
From Corollary \ref{cor:boundary_of_union}, we also have $\bound{j}{\alpha}(U \cup V) \subseteq \bound{j}{\alpha}U = \bound{j}{\alpha}V$, so it suffices to prove the converse inclusion.

Let $x \in \bound{j}{\alpha}U$.
Then there exists $y$ such that $x \leq y$ and $y \in \faces{j}{\alpha}U$ or $y \in \grade{\ell}{(\maxel{U})}$ for some $\ell < j$.
Using Lemma \ref{lem:maximal_in_boundary} together with the fact that $\bound{j}{\alpha}U = \bound{j}{\alpha}V$, we get in the first case that $y \in \faces{j}{\alpha}V$ and in the second case that $y \in \grade{\ell}{(\maxel{V})}$.
We conclude by Lemma \ref{lem:faces_of_union}.
\end{proof}

\begin{lem} \label{lem:pasting_higher_boundary}
Let $U$, $V$ be globular oriented graded posets, $k \in \mathbb{N}$, and suppose $U \cpiso{k}{\varphi} V$ is defined.
For all $n > k$ and $\alpha \in \set{ +, - }$, the pasting $\bound{n}{\alpha}U \cpiso{k}{\varphi} \bound{n}{\alpha}V$ is defined and maps isomorphically onto $\bound{n}{\alpha}(U \cpiso{k}{\varphi} V)$.
\end{lem}
\begin{proof}
By globularity, $\bound{k}{+}(\bound{n}{\alpha}U) = \bound{k}{+}U$ and $\bound{k}{-}(\bound{n}{\alpha}V) = \bound{k}{-}V$, so $\varphi$ has the correct type to determine the pasting $\bound{n}{\alpha}U \cpiso{k}{\varphi} \bound{n}{\alpha}V$.

By Corollary \ref{cor:inclusions_preserve_boundaries}, the inclusions of $U$ and $V$ into $U \cpiso{k}{\varphi} V$ preserve boundaries, so by the universal property of $\bound{n}{\alpha}U \cpiso{k}{\varphi} \bound{n}{\alpha}V$ we get an inclusion
\begin{equation*}
    \bound{n}{\alpha}U \cpiso{k}{\varphi} \bound{n}{\alpha}V \incl U \cpiso{k}{\varphi} V.
\end{equation*}
It suffices then to show that its image is $\bound{n}{\alpha}(U \cpiso{k}{\varphi} V)$.
If we identify $U$ and $V$ with their isomorphic images in $U \cpiso{k}{\varphi} V$, this is equivalent to proving 
\[\bound{n}{\alpha}U \cup \bound{n}{\alpha}V \subseteq \bound{n}{\alpha}(U \cup V);\] 
the converse inclusion is given by Corollary \ref{cor:boundary_of_union}.

Let $x \in \bound{n}{\alpha}U$.
Then there exists $y$ such that $x \leq y$ and $y \in \faces{n}{\alpha}U$ or $y \in \grade{j}{(\maxel{U})}$ for some $j < n$.
If $y \in \faces{n}{\alpha}U$, since $U \cap V = \bound{k}{+}U = \bound{k}{-}V$ is at most $k$\nbd dimensional, by Lemma \ref{lem:faces_of_union}
\begin{equation*}
    \faces{n}{\alpha}(U \cup V) = \faces{n}{\alpha}U + \faces{n}{\alpha}V,
\end{equation*}
so $y \in \faces{n}{\alpha}(U \cup V)$.
Similarly, if $y \in \grade{j}{(\maxel{U})}$ and $k < j < n$,
\begin{equation*}
    \grade{j}{(\maxel{(U \cup V)})} = \grade{j}{(\maxel{U})} + \grade{j}{(\maxel{V})},
\end{equation*}
so $y \in \grade{j}{(\maxel{(U \cup V)})}$.
In either case $x, y \in \bound{n}{\alpha}(U \cup V)$.

Suppose then that $y \in \grade{j}{(\maxel{U})}$ with $j \leq k$.
By Lemma \ref{lem:maximal_vs_faces} we have $\grade{k}{(\maxel{U})} \subseteq \faces{k}{+}U$, so from Lemma \ref{lem:maximal_in_boundary} and $\bound{k}{+}U = \bound{k}{-}V$ we deduce that $y \in \faces{k}{-}V$ if $j = k$ and $y \in \grade{j}{(\maxel{V})}$ if $j < k$.
Applying Lemma \ref{lem:faces_of_union} once more, we deduce in the first case that $y \in \faces{k}{-}(U \cup V)$ and in the second case that $z \in \grade{j}{(\maxel{(U \cup V)})}$.
In either case, $x, y \in \bound{n}{\alpha}(U \cup V)$.

This proves that $\bound{n}{\alpha}U \subseteq \bound{n}{\alpha}(U \cup V)$; the proof that $\bound{n}{\alpha}V \subseteq \bound{n}{\alpha}(U \cup V)$ is symmetrical.
\end{proof}

\begin{lem} \label{lem:pasting_preserves_globularity}
Let $U$, $V$ be globular oriented graded posets, $k \in \mathbb{N}$, and suppose $U \cpiso{k}{\varphi} V$ is defined.
Then $U \cpiso{k}{\varphi} V$ is globular.
\end{lem}
\begin{proof}
Let $m, n \in \mathbb{N}$ such that $m < n$, and $\alpha, \beta \in \set{ +, - }$.
If $n < k$, by Lemma \ref{lem:pasting_lower_boundary}
\begin{equation*}
    \bound{m}{\alpha}(\bound{n}{\beta}(U \cpiso{k}{\varphi} V)) = \bound{m}{\alpha}(\bound{n}{\beta}U) = \bound{m}{\alpha}(U) = \bound{m}{\alpha}(U \cpiso{k}{\varphi} V).
\end{equation*}
If $n = k$, by Lemma \ref{lem:pasting_top_boundary} and Lemma \ref{lem:pasting_lower_boundary},
\begin{equation*}
    \bound{m}{\alpha}(\bound{n}{-}(U \cpiso{k}{\varphi} V)) = \bound{m}{\alpha}(\bound{n}{-}U) = \bound{m}{\alpha}(U) = \bound{m}{\alpha}(U \cpiso{k}{\varphi} V)
\end{equation*}
and
\begin{equation*}
    \bound{m}{\alpha}(\bound{n}{+}(U \cpiso{k}{\varphi} V)) = \bound{m}{\alpha}(\bound{n}{+}V) = \bound{m}{\alpha}(V) = \bound{m}{\alpha}(U \cpiso{k}{\varphi} V).
\end{equation*}
Finally, if $n > k$, by Lemma \ref{lem:pasting_higher_boundary} we have
\begin{equation*}
    \bound{m}{\alpha}(\bound{n}{\beta}(U \cpiso{k}{\varphi} V)) = 
    \bound{m}{\alpha}(\bound{n}{\beta}U \cpiso{k}{\varphi} \bound{n}{\beta}V),
\end{equation*}
and by Lemma \ref{lem:boundaries_of_globular_are_globular} $\bound{n}{\beta}U$ and $\bound{n}{\beta}V$ are globular.
If $m < k$ we use Lemma \ref{lem:pasting_lower_boundary} to obtain
\begin{equation*}
    \bound{m}{\alpha}(\bound{n}{\beta}U \cpiso{k}{\varphi} \bound{n}{\beta}V) =
    \bound{m}{\alpha}(\bound{n}{\beta}U) =
    \bound{m}{\alpha}U = \bound{m}{\alpha}(U \cpiso{k}{\varphi} V).
\end{equation*}
If $m = k$ we use Lemma \ref{lem:pasting_top_boundary} instead to obtain
\begin{equation*}
    \bound{m}{-}(\bound{n}{\beta}U \cpiso{k}{\varphi} \bound{n}{\beta}V) = 
    \bound{m}{-}(\bound{n}{\beta}U) = \bound{m}{-}U = \bound{m}{-}(U \cpiso{k}{\varphi} V)
\end{equation*}
and similarly
\begin{equation*}
    \bound{m}{+}(\bound{n}{\beta}U \cpiso{k}{\varphi} \bound{n}{\beta}V) = 
    \bound{m}{+}(\bound{n}{\beta}V) = \bound{m}{+}V = \bound{m}{+}(U \cpiso{k}{\varphi} V).
\end{equation*}
Finally, if $m > k$ we use Lemma \ref{lem:pasting_higher_boundary} once more to obtain
\begin{equation*}
    \bound{m}{\alpha}(\bound{n}{\beta}U \cpiso{k}{\varphi} \bound{n}{\beta}V) =
    \bound{m}{\alpha}(\bound{n}{\beta}U) \cpiso{k}{\varphi} \bound{m}{\alpha}(\bound{n}{\beta}V) =
    \bound{m}{\alpha}U \cpiso{k}{\varphi} \bound{m}{\alpha}V
\end{equation*}
and once more to obtain
\begin{equation*}
    \bound{m}{\alpha}U \cpiso{k}{\varphi} \bound{m}{\alpha}V = 
    \bound{m}{\alpha}(U \cpiso{k}{\varphi} V). \qedhere
\end{equation*}
\end{proof}


\section{Rewrites and roundness} \label{sec:rewrites}

\begin{guide}
	In this section, we introduce the second main constructor for the class of molecules: the \cemph{rewrite} construction.
	The name stems from higher-dimensional rewriting theory, where an $(n+1)$\nbd dimensional cell is seen as the embodiment of a rewrite step on $n$\nbd dimensional diagrams, that is, the operation of rewriting its input boundary into its output boundary.
	Indeed, the rewrite construction takes two $n$\nbd dimensional oriented graded posets with isomorphic boundaries, and forms an $(n+1)$\nbd dimensional oriented graded poset that has one of each as its input and output boundary.

	The construction is particularly well-behaved on \cemph{round} oriented graded posets, also defined here.
	Roundness is a specialisation of globularity, which in addition to having good topological properties (Proposition \ref{prop:order_complex_of_round_molecule}) also implies that an oriented graded poset admits a convenient disjoint partition into the interiors of its boundaries (Lemma \ref{lem:round_partition_into_interiors}).
	Unlike globularity, roundness is not in general preserved by pastings; however, crucially, it is preserved by the rewrite construction (Lemma \ref{lem:rewrite_preserves_roundness}).
\end{guide}

\begin{dfn}[Rewrite construction] \index{oriented graded poset!rewrite} \index{$U \celto^\varphi V$}
Let $U$, $V$ be globular oriented graded posets of the same finite dimension $n$, and suppose $\varphi\colon \bound{}{}U \incliso \bound{}{}V$ is an isomorphism restricting to isomorphisms $\varphi^\alpha\colon \bound{}{\alpha} U \incliso \bound{}{\alpha} V$ for each $\alpha \in \set{ +, - }$.
Construct the pushout
\[\begin{tikzcd}
	\bound{}{}U & \bound{}{}V & V \\
	U && \bound{}{}(U \celto^\varphi V)
	\arrow["\varphi", hook, from=1-1, to=1-2]
	\arrow[hook, from=1-2, to=1-3]
	\arrow[hook', from=1-1, to=2-1]
	\arrow[hook, from=2-1, to=2-3]
	\arrow[hook', from=1-3, to=2-3]
	\arrow["\lrcorner"{anchor=center, pos=0.125, rotate=180}, draw=none, from=2-3, to=1-1]
\end{tikzcd}\]
in $\ogpos$. 
The \emph{rewrite of $U$ into $V$ along $\varphi$} is the oriented graded poset $U \celto^\varphi V$ obtained by adjoining a single $(n+1)$\nbd dimensional element $\top$ to $\bound{}{}(U \celto^\varphi V)$, with
\begin{equation*}
    \faces{}{-}\top \eqdef U_n, \quad \quad \faces{}{+}\top \eqdef V_n.
\end{equation*}
\end{dfn}

\begin{comm}
By Corollary \ref{cor:inclusions_preserve_boundaries}, we can identify $U$ and $V$ with their isomorphic images in $U \celto^\varphi V$, in such a way that $U \celto^\varphi V$ is equal to $(U \cup V) + \set{ \top }$, with $U \cap V = \bound{}{}U = \bound{}{}V$.
\end{comm}

\begin{lem} \label{lem:boundaries_of_rewrite}
Let $U$, $V$ be oriented graded posets and suppose $U \celto^\varphi V$ is defined.
Then
\begin{enumerate}
    \item $\bound{}{-}(U \celto^\varphi V) = U$,
    \item $\bound{}{+}(U \celto^\varphi V) = V$.
\end{enumerate}
\end{lem}
\begin{proof}
Identifying $U$ and $V$ with their isomorphic images, we will prove that $\bound{}{-}(U \celto^\varphi V) = U$ and $\bound{}{+}(U \celto^\varphi V) = V$.
Let $n \eqdef \dim{U} = \dim{V}$.
By construction, we have $\faces{n}{-}(U \celto^\varphi V) = \grade{n}{U}$ and $\faces{n}{+}(U \celto^\varphi V) = \grade{n}{V}$.

For all $k < n$, we have $\grade{k}{ ( \maxel{(U \celto^\varphi V)} ) } = \grade{k}{ ( \maxel{(U \cup V)} ) }$.
We claim that this is equal to both $\grade{k}{ ( \maxel{U} ) }$ and $\grade{k}{ ( \maxel{V} ) }$.
For $k < n-1$,
\begin{equation*}
    \grade{k}{ ( \maxel{U} ) } = \grade{k}{ ( \maxel{\bound{}{\alpha}U} ) } = \grade{k}{ ( \maxel{\bound{}{\alpha}V} ) } = \grade{k}{ ( \maxel{V} ) }
\end{equation*}
by Lemma \ref{lem:maximal_in_boundary}.
For $k = n-1$, by Lemma \ref{lem:maximal_vs_faces}
\begin{equation*}
    \grade{n-1}{ ( \maxel{U} ) } = \faces{}{-}U \cap \faces{}{+}U = \faces{}{-}V \cap \faces{}{+}V = \grade{n-1}{ ( \maxel{V} ) }.
\end{equation*}
We then conclude by Lemma \ref{lem:faces_of_union}.
\end{proof}

\begin{lem} \label{lem:rewrite_preserves_globularity}
Let $U$, $V$ be oriented graded posets and suppose $U \celto^\varphi V$ is defined.
Then $U \celto^\varphi V$ is globular.
\end{lem}
\begin{proof}
For all $k < \dim{U} = \dim{V}$ and $\alpha \in \set{ +, - }$, we have
\begin{equation*}
    \bound{k}{\alpha}U = \bound{k}{\alpha}(\bound{}{\beta}U) = \bound{k}{\alpha}(\bound{}{\beta}V) = \bound{k}{\alpha}V
\end{equation*}
since $\bound{}{\beta}U = \bound{}{\beta}V$ and $U$, $V$ are globular.
It then suffices to show that, for all $k < \dim{U}$ and $\alpha \in \set{ +, - }$,
\begin{equation*}
    \bound{k}{\alpha}(U \celto^\varphi V) = \bound{k}{\alpha}U.
\end{equation*}
Indeed, suppose this holds, and let $k < n < \dim{(U \celto^\varphi V)}$ and $\alpha, \beta \in \set{ +, - }$.
If $n = \dim{U}$, then by Lemma \ref{lem:boundaries_of_rewrite}
\begin{equation*}
    \bound{k}{\alpha}(\bound{n}{-}(U \celto^\varphi V)) = \bound{k}{\alpha}U = \bound{k}{\alpha}(U \celto^\varphi V)
\end{equation*}
and similarly
\begin{equation*}
    \bound{k}{\alpha}(\bound{n}{+}(U \celto^\varphi V)) = \bound{k}{\alpha}V = \bound{k}{\alpha}U = \bound{k}{\alpha}(U \celto^\varphi V).
\end{equation*}
If $n < \dim{U}$, then
\begin{equation*}
    \bound{k}{\alpha}(\bound{n}{\beta}(U \celto^\varphi V)) = \bound{k}{\alpha}(\bound{n}{\beta}U) = \bound{k}{\alpha}U = \bound{k}{\alpha}(U \celto^\varphi V)
\end{equation*}
using the globularity of $U$.

Let then $k < \dim{U}$ and $\alpha \in \set{ +, - }$.
We have $\faces{k}{\alpha}(U \celto^\varphi V) = \faces{k}{\alpha}(U \cup V)$.
Since $\faces{k}{\alpha}U = \faces{k}{\alpha}V$, by Lemma \ref{cor:boundary_of_union} we have $\faces{k}{\alpha}(U \cup V) = \faces{k}{\alpha}U$.
Similarly, we prove that for all $j < k$ we have $\grade{j}{ ( \maxel{(U \cup V)} ) } = \grade{j}{ ( \maxel{U} ) }$.
It follows that $\bound{k}{\alpha}(U \celto^\varphi V) = \bound{k}{\alpha}U$.
\end{proof}

\begin{dfn}[Roundness] \index{oriented graded poset!round}
Let $U$ be an oriented graded poset.
We say that $U$ is \emph{round} if it is globular and, for all $n < \dim{U}$,
\begin{equation*}
    \bound{n}{-}U \cap \bound{n}{+}U = \bound{n-1}{}U.
\end{equation*}
\end{dfn}

\begin{lem} \label{lem:round_is_pure}
Let $U$ be a round oriented graded poset.
If $U$ is finite-dimensional, then $U$ is pure.
\end{lem}
\begin{proof}
We will prove the contrapositive.
Suppose that $U$ is not pure.
Then there exists a maximal element $x$ in $U$ with $k \eqdef \dim{x} < \dim{U}$.
By Lemma \ref{lem:maximal_vs_faces}, $x \in \bound{k}{-}U \cap \bound{k}{+}U$.
Then $\bound{k}{-}U \cap \bound{k}{+}U$ is $k$\nbd dimensional and cannot be equal to $\bound{k-1}{}U$, which is $(k-1)$\nbd dimensional.
It follows that $U$ is not round.
\end{proof}

\begin{lem} \label{lem:boundaries_of_round_and_globular}
Let $U$ be round, $n \in \mathbb{N}$, and $\alpha \in \set{ +, - }$.
Then $\bound{n}{\alpha}U$ is round.
\end{lem}
\begin{proof}
If $n \geq \dim{U}$ there is nothing to prove, so suppose $n < \dim{U}$.
By Lemma \ref{lem:boundaries_of_globular_are_globular}, $\bound{n}{\alpha}U$ is globular.
Let $k < \dim{(\bound{n}{\alpha}U)} \leq n$.
Then
\begin{equation*}
    \bound{k}{-}(\bound{n}{\alpha}U) \cap \bound{k}{+}(\bound{n}{\alpha}U) = \bound{k}{-}U \cap \bound{k}{+}U = \bound{k-1}{}U = \bound{k-1}{}(\bound{n}{\alpha}U)
\end{equation*}
using roundness of $U$.
\end{proof}

\begin{lem} \label{lem:round_partition_into_interiors}
	Let $U$ be a round oriented graded poset.
	Then $U$ is partitioned into
	\[
		\inter{U} + \sum_{k < \dim{U}} \left( \inter{\bound{k}{-}U} + \inter{\bound{k}{+}U} \right).
	\]
\end{lem}
\begin{proof}
Let $x \in U$.
Since $U$ is globular, by Lemma \ref{lem:globular_interiors} either $x \in \inter{U}$, or there exist $k < \dim{U}$ and $\alpha \in \set{+, -}$ such that $x \in \inter{\bound{k}{\alpha}U}$.
It then suffices to show that these subsets are all disjoint.
Let $\alpha \in \set{+, -}$ and $k < \dim{U}$.
Then
\[ \inter{U} \cap \inter{\bound{k}{\alpha}U} \subseteq \inter{U} \cap \bound{}{}U = \varnothing \]
since $\bound{k}{\alpha}U \subseteq \bound{}{}U$.
Let $\beta \in \set{+, -}$ and $j < \dim{U}$, and suppose $j \neq k$.
Assume without loss of generality that $j < k$.
Then 
\[
	\bound{j}{\beta}U = \bound{j}{\beta}(\bound{k-1}{\beta}(\bound{k}{\alpha}U)) \subseteq \bound{}{}(\bound{k}{\alpha}U)
\]
by globularity, so
\[
	\inter{\bound{j}{\beta}U} \cap \inter{\bound{k}{\alpha}U} \subseteq \bound{}{}(\bound{k}{\alpha}U) \cap \inter{\bound{k}{\alpha}U} = \varnothing.
\]
Finally, in the case $j = k$,
\[
	\bound{k}{-}U \cap \bound{k}{+}U = \bound{k-1}{}U = \bound{}{}(\bound{k}{+}U)
\]
by roundness, so
\[
	\inter{\bound{k}{-}U} \cap \inter{\bound{k}{+}U} \subseteq \bound{}{}(\bound{k}{+}U) \cap \inter{\bound{k}{+}U} = \varnothing.
\]
This concludes the proof.
\end{proof}

\begin{lem} \label{lem:rewrite_preserves_roundness}
Let $U$, $V$ be round and suppose $U \celto^\varphi V$ is defined.
Then $U \celto^\varphi V$ is round.
\end{lem}
\begin{proof}
Globularity follows from Lemma \ref{lem:rewrite_preserves_globularity}, so we only need to prove roundness.
Let $n \eqdef \dim{U} = \dim{V}$.
By Lemma \ref{lem:boundaries_of_rewrite}
\begin{equation*}
    \bound{}{-}(U \celto^\varphi V) \cap \bound{}{+}(U \celto^\varphi V) = U \cap V = \bound{}{}U = \bound{}{}V,
\end{equation*}
and by globularity $\bound{}{}U = \bound{}{}(\bound{}{-}(U \celto^\varphi V)) = \bound{n-1}{}(U \celto^\varphi V)$.
Finally, for $k < n$
\begin{equation*}
    \bound{k}{-}(U \celto^\varphi V) \cap \bound{k}{+}(U \celto^\varphi V) = \bound{k}{-}U \cap \bound{k}{+}U = \bound{k-1}{}U = \bound{k-1}{}(U \celto^\varphi V)
\end{equation*}
by globularity of $U \celto^\varphi V$ and roundness of $U$.
\end{proof}

\begin{exm}[A molecule which is pure but not round] \index[counterex]{A molecule which is pure but not round}
	Using Lemma \ref{lem:round_is_pure}, it is easy to find examples of oriented graded posets, or even molecules, which are not round: it suffices to pick one which is not pure.
	For example, (\ref{eq:example_ogposet}) is not pure, since it is 2\nbd dimensional but has a 1\nbd dimensional maximal element, so it is also not round.

	On the other hand, let $U$ be the molecule $\globe{2} \cp{0} \globe{2}$ from (\ref{eq:globe2_cp0_globe2}).
	Then $U$ is pure 2\nbd dimensional, but it is not round, since
\[ 
	\bound{0}{} U = \set{(0, 0), (0, 2)} \subsetneq 
	\bound{1}{-} U \cap \bound{1}{+} U = \set{(0, 0), (0, 1), (0, 2)}.
\]
\end{exm}


\section{The inductive definition of molecules} \label{sec:inductive_molecules}

\begin{guide}
	In this section, we finally reach the definition of \cemph{molecule} and of \cemph{atom}.
	We also prove some basic stability properties of the classes of molecules and atoms: all pastings of molecules are molecules (Lemma 
	\ref{lem:pastings_preserve_molecules}), all input and output boundaries of molecules are molecules (Lemma \ref{lem:molecules_are_globular}), and the lower set of each element of a molecule is an atom (Lemma \ref{lem:all_downsets_are_atoms}).
\end{guide}

\begin{dfn}[Point] \index{point}
The \emph{point} is the oriented graded poset $1$ with a single element and trivial orientation.
\end{dfn}

\begin{dfn}[Molecule] \index{molecule}
The class of \emph{molecules} is the inductive subclass of oriented graded posets closed under isomorphisms and generated by the following clauses.
\begin{enumerate}
    \item (\textit{Point}). The point is a molecule.
    \item (\textit{Paste}). Let $U$, $V$ be molecules, let $k < \min \set{ \dim U, \dim V }$, and let $\varphi\colon \bound{k}{+}U \incliso \bound{k}{-}V$ be an isomorphism.
	    Then $U \cpiso{k}{\varphi} V$ is a molecule.
    \item (\textit{Atom}). Let $U$, $V$ be \emph{round} molecules of the same finite dimension and let $\varphi\colon \bound{}{}U \incliso \bound{}{}V$ be an isomorphism restricting to $\varphi^\alpha\colon \bound{}{\alpha} U \incliso \bound{}{\alpha} V$ for each $\alpha \in \set{ +, - }$.
	    Then $U \celto^\varphi V$ is a molecule.
\end{enumerate}
\end{dfn}

\begin{comm} \label{comm:properties_isomorphism_invariant}
We will prove many properties of molecules by induction on their construction.
Since we only consider properties that are invariant under isomorphism, in these proofs we will not explicitly consider closure under isomorphisms.
\end{comm}

\begin{lem} \label{lem:molecules_are_finite}
	Let $U$ be a molecule.
	Then $\size{U}$ is finite.
\end{lem}
\begin{proof}
	By induction on the construction of $U$.
	If $U$ was produced by (\textit{Point}), then $\size{U} = 1$.
	If $U$ was produced by (\textit{Paste}), then it is equal to $V \cpiso{k}{\varphi} W$ for some molecules $V$, $W$ and $k < \min \set{ \dim V, \dim W }$.
	Then $\size{U} \leq \size{V} + \size{W}$, which is finite by the inductive hypothesis.
	If $U$ was produced by (\textit{Atom}), then it is of the form $V \celto^\varphi W$ for some molecules $V$, $W$.
	Then $\size{U} \leq \size{V} + \size{W} + 1$, which is finite by the inductive hypothesis.
\end{proof}

\begin{rmk}
	It follows evidently that all molecules have finite dimension.
\end{rmk}

\begin{lem} \label{lem:only_0_molecule}
Let $U$ be a molecule.
If $\dim{U} = 0$, then $U$ is isomorphic to the point.
\end{lem}
\begin{proof}
By induction on the construction of $U$.
If $U$ was produced by (\textit{Point}), then $U = 1$ and $\dim{U} = 0$.
If $U$ was produced by (\textit{Paste}), then it is equal to $V \cpiso{k}{\varphi} W$ where $V$, $W$ are molecules with $k < \min \set{ \dim V, \dim W }$.
Then $\dim{U} = \max\set{ \dim V, \dim W  } > k \geq 0$.
If $U$ was produced by (\textit{Atom}), then it is of the form $V \celto^\varphi W$, and $\dim{U} = \dim{V} + 1 = \dim{W} + 1 > 0$.
\end{proof}

\begin{lem} \label{lem:pastings_preserve_molecules}
Let $U$, $V$ be molecules, $k \in \mathbb{N}$, and $\varphi\colon \bound{k}{+}U \incliso \bound{k}{-}V$ an isomorphism.
Then $U \cpiso{k}{\varphi} V$ is a molecule.
\end{lem}
\begin{proof}
If $k < \min \set{ \dim U, \dim V }$, then this is an application of the (\textit{Paste}) constructor.
If $k \geq \dim U$, then $\bound{k}{+}U = U$ and $U \cp{k} V$ is isomorphic to $V$, which is a molecule by assumption.
Similarly, if $k \geq \dim V$, then $U \cp{k} V$ is isomorphic to $U$.
\end{proof}

\begin{lem} \label{lem:molecules_are_globular}
Let $U$ be a molecule, $n \in \mathbb{N}$, $\alpha \in \set{ +, - }$.
Then
\begin{enumerate}
    \item $U$ is globular,
    \item $\bound{n}{\alpha}U$ is a molecule,
    \item if $n \leq \dim{U}$, then $\dim{\bound{n}{\alpha}U} = n$.
\end{enumerate}
\end{lem}
\begin{proof}
By induction on the construction of $U$.
Suppose $U$ was produced by (\textit{Point}). 
Then $U$ is the point, it has no non-trivial boundaries, and is trivially globular.

Suppose $U$ was produced by (\textit{Paste}).
Then $U = V \cpiso{k}{\varphi} W$ for some molecules $V$, $W$.
By the inductive hypothesis, $V$ and $W$ are globular, and by Lemma \ref{lem:pasting_preserves_globularity} so is $U$.
We have $k < \min \set{ \dim{V}, \dim{W} }$.
If $n = k$, then by Lemma \ref{lem:pasting_top_boundary} $\bound{n}{-}U$ is equal to $\bound{n}{-}V$ and $\bound{n}{+}U$ to $\bound{n}{+}W$.
By the inductive hypothesis, both of these are $n$\nbd dimensional molecules.
If $n < k$, then by Lemma \ref{lem:pasting_lower_boundary} $\bound{n}{\alpha}U$ is equal to $\bound{n}{\alpha}V$, and again the inductive hypothesis applies.
If $n > k$, then by Lemma \ref{lem:pasting_higher_boundary} $\bound{n}{\alpha}U$ is equal to $\bound{n}{\alpha}V \cpiso{k}{\varphi} \bound{n}{\alpha}W$.
By the inductive hypothesis, $\bound{n}{\alpha}V$ and $\bound{n}{\alpha}W$ are molecules, and if $n < \dim{U} = \max \set{ \dim{V}, \dim{W} }$, at least one of them is $n$\nbd dimensional.

Finally, suppose $U$ was produced by (\textit{Atom}).
Then $U = V \celto^\varphi W$ for some round molecules $V$, $W$ of the same dimension.
By the inductive hypothesis, $V$ and $W$ are globular, and by Lemma \ref{lem:rewrite_preserves_globularity} so is $U$.
If $n \geq \dim{U}$, then $\bound{n}{\alpha}U = U$ is by assumption a molecule.
If $n = \dim{U}-1$, then by Lemma \ref{lem:boundaries_of_rewrite} $\bound{}{-}U$ is equal to $V$ and $\bound{}{+}U$ to $W$, both molecules of dimension $n$.
If $n < \dim{U}-1$, then $\bound{n}{\alpha}U = \bound{n}{\alpha}V = \bound{n}{\alpha}W$ by globularity, and the inductive hypothesis applies.
\end{proof}

\begin{dfn}[Atom] \index{molecule!atom} \index{atom}
An \emph{atom} is a molecule with a greatest element.
\end{dfn}

\begin{lem} \label{lem:atom_greatest_element}
Let $U$ be a molecule.
The following are equivalent:
\begin{enumerate}[label=(\alph*)]
    \item $U$ is an atom;
    \item the final constructor producing $U$ is (\textit{Point}) or (\textit{Atom}).
\end{enumerate}
\end{lem}
\begin{proof}
If $U$ was produced by (\textit{Point}), then $U$ is the point, which trivially has a greatest element.

If $U$ was produced by (\textit{Paste}), then $U$ splits into a union $V \cup W$, where $V \cap W = \bound{k}{+}V = \bound{k}{-}W$ and $k < \max\set{ \dim{V}, \dim{W} }$.
Then there exist elements $x_1 \in V$ and $x_2 \in W$ such that
\begin{enumerate}
    \item $x_1$ is maximal in $V$ and $x_2$ is maximal in $W$,
    \item $\dim{x_1} > k$ and $\dim{x_2} > k$.
\end{enumerate}
By Lemma \ref{lem:dimension_of_boundary}, $\dim{(V \cap W)} \leq k$, so neither $x_1$ nor $x_2$ are contained in $V \cap W$.
It follows that $x_1$ and $x_2$ are distinct maximal elements of $U$, so $U$ does not have a greatest element.

If $U$ was produced by (\textit{Atom}), then $U$ splits into $(U_{-} \cup U_{+}) + \set{ \top }$, where $U_{-}$ and $U_{+}$ are round molecules of dimension $n$, and $\faces{}{\alpha}\top = \grade{n}{(U_\alpha)}$ for each $\alpha \in \set{ +, - }$.
By Lemma \ref{lem:round_is_pure}, we have $U_\alpha = \clos{\grade{n}{(U_\alpha)}}$, so $U_\alpha = \bound{}{\alpha}\top \subseteq \clset{{ \top }}$.
It follows that all elements of $U$ are in the closure of $x$, that is, $x$ is the greatest element of $U$.
\end{proof}

\begin{cor} \label{cor:atoms_are_round}
All atoms are round.
\end{cor}
\begin{proof}
Let $U$ be an atom.
If it was produced by (\textit{Point}), it is trivially round.
If it was produced by (\textit{Atom}), it is round by Lemma \ref{lem:rewrite_preserves_roundness}.
\end{proof}

\begin{lem} \label{lem:all_downsets_are_atoms}
Let $U$ be a molecule, $x \in U$.
Then $\clset{{ x }}$ is an atom.
\end{lem}
\begin{proof}
By induction on the construction of $U$.
If $U$ was produced by (\textit{Point}), then $x$ must be the unique element of $U$ whose closure is $U$ itself.
If $U$ was produced by (\textit{Paste}), it splits into $V \cup W$, and $x \in V$ or $x \in W$; the inductive hypothesis applies.
If $U$ was produced by (\textit{Atom}), it is equal to $(V \cup W) + \set{ \top }$, and either $x \in V$ or $x \in W$, in which case the inductive hypothesis applies, or $x = \top$, and $\clset{{ x }} = U$ is an atom by definition.
\end{proof}

\begin{lem} \label{lem:molecules_are_connected}
Let $U$ be a molecule.
Then $U$ is connected, that is, 
\begin{enumerate}
    \item $U$ is non-empty,
    \item for all closed subsets $V$, $W \subseteq U$, if $U = V \cup W$ and $V \cap W = \varnothing$, then $V = \varnothing$ or $W = \varnothing$.
\end{enumerate}
\end{lem}
\begin{proof}
Let $V$, $W \subseteq U$ be closed subsets such that $U = V \cup W$.
We proceed by induction on the construction of $U$.
If $U$ was produced by (\textit{Point}) or by (\textit{Atom}), it is an atom by Lemma \ref{lem:atom_greatest_element}, so it has a greatest element $\top$.
Then either $\top \in V$ or $\top \in W$, implying $V = U$ and $W = U$, respectively, hence $V \cap W = W$ and $V \cap W = V$, respectively.

If $U$ was produced by (\textit{Paste}), it splits into $V' \cup W'$ along the $k$\nbd boundary for some molecules $V', W'$ and $k < \min \set{\dim{V'}, \dim{W'}}$.
We have that
\begin{equation*}
    V' = (V \cap V') \cup (W \cap V'), \quad \quad W' = (V \cap W') \cup (W \cap W').
\end{equation*}
By the inductive hypothesis, one of the parts in each union is empty.
Suppose $W \cap V' = \varnothing$ and $W \cap W' = \varnothing$.
Then $W \cap (V' \cup W') = W \cap U = W = \varnothing$.
Similarly, if $V \cap V' = \varnothing$ and $V \cap W' = \varnothing$, then $V = \varnothing$.

Suppose that $W \cap V' = \varnothing$ and $V \cap W' = \varnothing$, or that $V \cap V' = \varnothing$ and $W \cap W' = \varnothing$.
Then $\bound{k}{+}V' = V' \cap W' = V \cap V' \cap W \cap W' = \varnothing$, a contradiction.
\end{proof}

\begin{lem} \label{lem:maximal_0dim}
Let $U$ be a molecule.
Then $U$ has a maximal 0\nbd dimensional element if and only if $\dim{U} = 0$.
\end{lem}
\begin{proof}
If $\dim{U} = 0$, then $U$ is the point by Lemma \ref{lem:only_0_molecule}, hence has a greatest 0\nbd dimensional element.
Conversely, let $x$ be maximal and 0\nbd dimensional, and let $V \eqdef \clos{ ((\maxel{U}) \setminus \set{x}) }$.
Then $\set{x}$ is closed, $U = V \cup \set{x}$, and $V \cap \set{x} = \varnothing$.
By Lemma \ref{lem:molecules_are_connected}, $V = \varnothing$, so $U = \set{x}$.
\end{proof}


\section{Isomorphisms of molecules are unique} \label{sec:unique_iso}

\begin{guide}
	In this section, we prove that molecules have no non-trivial automorphisms (Proposition 
	\ref{prop:molecule_unique_automorphism}), hence that two isomorphic molecules are isomorphic in a unique way.
	We then show that pasting of molecules satisfies the equations of strict $\omega$\nbd categories up to unique isomorphism (Proposition 
	\ref{prop:associativity_of_pasting}, Proposition \ref{prop:unitality_of_pasting}, Proposition \ref{prop:interchange_of_pasting}).

	We derive the main result from properties of a graph-like structure --- a \emph{directed graph with open edges} --- associated to a molecule.
	This is, in fact, one way of formalising the \cemph{string diagram} associated to the top dimensions of a pasting diagram, and that is how we are going to interpret it.
\end{guide}

\begin{dfn}[Induced subgraph] \index{directed graph!subgraph!induced}
Let $\mathscr{G}$ be a directed graph and let $W \subseteq V_\mathscr{G}$.
The \emph{induced subgraph} of $\mathscr{G}$ on $W$ is the directed graph
\begin{equation*}
\restr{\mathscr{G}}{W} \eqdef
    \begin{tikzcd}
	{E'} & {W}
	\arrow["\restr{s}{E'}", shift left=1.5, from=1-1, to=1-2]
	\arrow["\restr{t}{E'}"', shift right=1.5, from=1-1, to=1-2]
    \end{tikzcd}
\end{equation*}
where $E' \eqdef \set{ e \in E_\mathscr{G} \mid s(e), t(e) \in W }$.
\end{dfn}

\begin{dfn}[Directed graph with open edges] \index{directed graph!with open edges}
A \emph{directed graph with open edges} is a directed graph
\begin{equation*}
\mathscr{G} \eqdef
    \begin{tikzcd}
	{E_\mathscr{G}} & {N_\mathscr{G} + W_\mathscr{G}}
	\arrow["s", shift left=1.5, from=1-1, to=1-2]
	\arrow["t"', shift right=1.5, from=1-1, to=1-2]
    \end{tikzcd}
\end{equation*}
with set of vertices bipartite into a set $N_\mathscr{G}$ of \emph{node vertices} and a set $W_\mathscr{G}$ of \emph{wire vertices}, satisfying the following properties:
\begin{enumerate}
    \item the bipartition $N_\mathscr{G} + W_\mathscr{G}$ exhibits $\mathscr{G}$ as a bipartite graph, that is, every edge connects a node vertex to a wire vertex or vice versa;
    \item each wire vertex is the source of at most one edge and the target of at most one edge.
\end{enumerate}
\end{dfn}

\begin{comm}
This is, up to inessential encoding details, the structure called an \emph{open graph} in \cite{dixon2013open} and simply a \emph{graph} in \cite{kock2016graphs}.
\end{comm}

\begin{dfn}[Boundary of a directed graph with open edges] \index{directed graph!with open edges!boundary} \index{boundary!of a directed graph with open edges}
Let $\mathscr{G}$ be a directed graph with open edges.
The \emph{input boundary} of $\mathscr{G}$ is the set
\begin{equation*}
    \faces{}{-}\mathscr{G} \eqdef \set{ x \in W_\mathscr{G} \mid \invrs{t}(x) = \varnothing }
\end{equation*}
and the \emph{output boundary} of $\mathscr{G}$ is the set
\begin{equation*}
    \faces{}{+}\mathscr{G} \eqdef \set{ x \in W_\mathscr{G} \mid \invrs{s}(x) = \varnothing }.
\end{equation*}
\end{dfn}

\begin{dfn}[Graph of a molecule] \index{molecule!graph} \index{$\graph{U}$}
Let $U$ be a molecule, $n \eqdef \dim{U}$.
The \emph{graph of $U$} is the directed graph
\begin{equation*}
\graph{U} \eqdef
    \begin{tikzcd}
	{ E_{\graph{U}} } & { N_{\graph{U}} + W_{\graph{U}} },
	\arrow["s", shift left=1.5, from=1-1, to=1-2]
	\arrow["t"', shift right=1.5, from=1-1, to=1-2]
    \end{tikzcd}
\end{equation*}
where
\begin{itemize}
    \item $E_{\graph{U}} \eqdef \set{ (x, y) \mid x \in U_n, y \in \faces{}{+}x } + \set{ (x, y) \mid y \in U_n, x \in \faces{}{-}y }$,
    \item $N_{\graph{U}} \eqdef \grade{n}{U}$,
    \item $W_{\graph{U}} \eqdef \grade{n-1}{U}$,
    \item $s\colon (x, y) \mapsto x$,
    \item $t\colon (x, y) \mapsto y$.
\end{itemize}
\end{dfn}

\begin{rmk}
	Note that, if we forget the separation of the vertex set into $N_{\graph{U}}$ and $W_{\graph{U}}$, then $\graph{U}$ is the induced subgraph of $\hasseo{U}$ on vertices of dimension $n$ and $n-1$.
\end{rmk}

\begin{prop} \label{prop:graph_of_molecule_properties}
Let $U$ be a molecule.
Then 
\begin{enumerate}
    \item $\graph{U}$ is a directed graph with open edges,
    \item $\graph{U}$ is acyclic,
    \item $\faces{}{\alpha}\graph{U} = \faces{}{\alpha}U$ for all $\alpha \in \set{ +, - }$.
\end{enumerate}
\end{prop}
\begin{proof}
The fact that $\faces{}{\alpha}\graph{U} = \faces{}{\alpha}U$ for all $\alpha \in \set{ +, - }$ is immediate from the definitions.
Moreover, $\graph{U}$ is bipartite by construction, so it suffices to check the other conditions.

We proceed by induction on the construction of $U$.
If $U$ was produced by (\textit{Point}) or by (\textit{Atom}), then by Lemma \ref{lem:atom_greatest_element} it has a greatest element $\top$.
In this case, $\graph{U}$ has a single edge $(x, \top)$ for each $x \in \faces{}{-}\top$ and a single edge $(\top, x)$ for each $x \in \faces{}{+}\top$.
Since $\faces{}{-}\top \cap \faces{}{+}\top = \varnothing$, the graph is acyclic.

If $U$ was produced by (\textit{Paste}), it is of the form $V \cpiso{k}{\varphi} W$.
Let $n \eqdef \dim{U}$.
If $k < n - 1$, then $\graph{U}$ is the disjoint union of the induced subgraphs on the vertices in the image of $V$ and $W$, respectively.
If $n = \dim{V} = \dim{W}$ we can conclude by the inductive hypothesis.
Otherwise, the inductive hypothesis applies to one of the components, while the other is a discrete graph with no node vertices, trivially satisfying the conditions of an acyclic directed graph with open edges.

If $k = n - 1$, observe first that necessarily $\dim{V} = \dim{W} = n$.
Then $\graph{U}$ is the union of $\graph{V}$ and $\graph{W}$, and their intersection consists of the wire vertices in $\faces{n-1}{+}V = \faces{n-1}{-}W$.
Let $x$ be a wire vertex.
If $x \in V \setminus W$ or $x \in V \setminus W$, it is the source of at most one edge and the target of at most one edge by the inductive hypothesis applied to $\graph{V}$ and $\graph{W}$.
If $x \in V \cap W$, then $x \in \faces{}{+}\graph{V}$, so it is the source of no edge of $\graph{V}$ and at most one edge of $\graph{W}$, and $x \in \faces{}{-}\graph{W}$, so it is the target of no edge of $\graph{W}$ and the source of at most one edge of $\graph{V}$.

Finally, suppose there is a cycle in $\graph{U}$.
Because $\graph{V}$ and $\graph{W}$ are separately acyclic, such a cycle needs to cross from $V$ to $W \setminus V$ and back.
However, a path entering $V$ from $W \setminus V$ must enter a wire vertex $y$ from a node vertex $x \in W$ such that $y \in \faces{}{+}x$.
But $\grade{n-1}{(V \cap W)} = \faces{}{-}W$, so this is impossible.
We conclude that $\graph{U}$ is acyclic.
\end{proof}

\begin{exm}[The graph and string diagram of a molecule] \index{string diagram} \label{exm:string_diagrams}
	Let $U$ be the oriented face poset of the pasting diagram
	\begin{equation} \label{eq:pasting_to_string}
		\begin{tikzcd}[sep=small]
	& \bullet \\
	\bullet &&& \bullet && \bullet \\
	&& \bullet
	\arrow["3", curve={height=-6pt}, from=2-1, to=1-2]
	\arrow[""{name=0, anchor=center, inner sep=0}, "5", curve={height=-12pt}, from=1-2, to=2-4]
	\arrow[""{name=1, anchor=center, inner sep=0}, "6", curve={height=-18pt}, from=2-4, to=2-6]
	\arrow[""{name=2, anchor=center, inner sep=0}, "2"', curve={height=18pt}, from=2-4, to=2-6]
	\arrow[""{name=3, anchor=center, inner sep=0}, "0"', curve={height=12pt}, from=2-1, to=3-3]
	\arrow["1"', curve={height=6pt}, from=3-3, to=2-4]
	\arrow["4", from=1-2, to=3-3]
	\arrow["0"', curve={height=-6pt}, shorten <=7pt, shorten >=3pt, Rightarrow, from=3, to=1-2]
	\arrow["2", shorten <=5pt, shorten >=5pt, Rightarrow, from=2, to=1]
	\arrow["1"', curve={height=6pt}, shorten <=3pt, shorten >=7pt, Rightarrow, from=3-3, to=0]
\end{tikzcd}\end{equation}
	where we left 0\nbd cells unlabelled.
	Then $\graph{U}$ is the directed graph with open edges
\[\begin{tikzcd}[column sep=tiny]
	&&& {{\scriptstyle (1, 5)}\;\bullet} \\
	&&& {{\scriptstyle (2, 1)}\;\bullet} && {{\scriptstyle (1, 6)}\;\bullet} \\
	{{\scriptstyle (1, 3)}\;\bullet} && {{\scriptstyle (1, 4)}\;\bullet} && {{\scriptstyle (1, 1)}\;\bullet} & {{\scriptstyle (2, 2)}\;\bullet} \\
				    & {{\scriptstyle (2, 0)}\;\bullet} &&&& {{\scriptstyle (1, 2)}\;\bullet} \\
	& {{\scriptstyle (1, 0)}\;\bullet}
	\arrow[from=5-2, to=4-2]
	\arrow[from=4-2, to=3-3]
	\arrow[from=4-2, to=3-1]
	\arrow[from=3-3, to=2-4]
	\arrow[from=3-5, to=2-4]
	\arrow[from=2-4, to=1-4]
	\arrow[from=4-6, to=3-6]
	\arrow[from=3-6, to=2-6]
\end{tikzcd}\]
	where we took advantage of the acyclicity of the graph to draw it in such a way that all edges point upward.

	Now, the idea is that wire vertices, which have at most one outgoing and one incoming edge, should really be seen as ``midpoints'' of a single \cemph{wire}, whose ``halves'' are the edges incident to the wire vertex.
	In order to make sure that every wire vertex is a midpoint, we frame the graph to a rectangular ``canvas'', by adding the missing half-wires to all wire vertices in $\faces{}{+}\graph{U}$ and in $\faces{}{-}{\graph{U}}$, and extending them to reach the top and bottom edges of the canvas.
	The result is something like
	\[
		\input{img/strdiag_example.tex}
	\]
	which is a typical \emph{string diagram} representation of (\ref{eq:pasting_to_string}).
\end{exm}

\begin{cor} \label{cor:codimension_1_elements}
Let $U$ be a molecule, $x \in U$ with $\codim{x}{U} = 1$, and $\alpha \in \set{+, -}$.
Then
\begin{enumerate}
    \item $x \in \maxel{U}$ if and only if $\size{\cofaces{}{}{x}} = 0$,
    \item $x \in \faces{}{\alpha}U \setminus \faces{}{-\alpha}{U}$ if and only if $\size{\cofaces{}{\alpha}x} = 1$ and $\size{\cofaces{}{-\alpha}x} = 0$,
    \item $x \notin \faces{}{}U$ if and only if $\size{\cofaces{}{+}x} = \size{\cofaces{}{-}x} = 1$.
\end{enumerate}
\end{cor}
\begin{proof}
By Proposition \ref{prop:graph_of_molecule_properties}, $\graph{U}$ is a directed graph with open edges, and by construction we can identify $\cofaces{}{-}x$ with $\invrs{s}(x)$ and $\cofaces{}{+}x$ with $\invrs{t}(x)$.
It follows that $\size{\cofaces{}{\alpha}x} \leq 1$.
The statement then follows from the isomorphism between $\faces{}{\alpha}U$ and $\faces{}{\alpha}\graph{U}$, combined with Lemma \ref{lem:maximal_vs_faces}.
\end{proof}

\begin{lem} \label{lem:path_in_graph_of_molecule}
Let $U$ be a molecule, $n \eqdef \dim{U} > 0$, and $x \in \grade{n}{U}$.
Then there exist $y_- \in \faces{}{-}U$ and $y_+ \in \faces{}{+}U$ such that there is a path from $y_-$ to $y_+$ passing through $x$ in $\graph{U}$.
\end{lem}
\begin{proof}
We construct a path $x = x_0 \to y_0 \to \ldots \to x_m \to y_+$ by successive extensions; the construction of a path from $y_-$ to $x$ is dual.
Suppose we have reached $x_i$.
By Lemma \ref{lem:all_downsets_are_atoms} $\clset{{ x_i }}$ is an atom, so $\bound{}{+}x_i$ is $(n-1)$\nbd dimensional and $\faces{}{+}x_i$ is non-empty.
Pick $y_i$ in $\faces{}{+}x_i$.
If $y_i$ has no input cofaces, then $y_i \in \faces{}{+}U$, so we can let $m \eqdef i$ and $y_+ \eqdef y_i$.
Otherwise, pick $x_{i+1} \in \cofaces{}{-}y_i$.
Since $\graph{U}$ is finite and acyclic by Proposition \ref{prop:graph_of_molecule_properties}, this procedure must terminate after a finite number of steps.
\end{proof}

\begin{prop} \label{prop:molecule_unique_automorphism}
Let $U$ be a molecule and $\imath\colon U \incliso U$ an automorphism.
Then $\imath$ is the identity.
\end{prop}
\begin{proof}
We proceed by induction on $n \eqdef \dim{U}$.
If $n = 0$, then $U = 1$ by Lemma \ref{lem:only_0_molecule}, and the only endomorphism of $1$ is the identity.

Suppose $n > 0$ and let $\alpha \in \set{ +, - }$.
By Proposition \ref{lem:molecules_are_globular}, $\bound{}{\alpha}U$ is a molecule of dimension $n - 1$, and $\imath(\bound{}{\alpha}U) = \bound{}{\alpha}U$.
By the inductive hypothesis, the restriction of $\imath$ to $\bound{}{\alpha}U$ is the identity.

Let $x \in \maxel{U}$, and suppose $\imath(x) = x$.
Then $\imath(\bound{}{\alpha}x) = \bound{}{\alpha}x$.
By Lemma \ref{lem:all_downsets_are_atoms}, $\clset{{ x }}$ is an atom, so $\bound{}{\alpha}x$ is a molecule of dimension strictly lower than $n$.
By the inductive hypothesis the restriction of $\imath$ to $\bound{}{\alpha}x$ is the identity.
Since $\clset{{ x }} = (\bound{}{-}x \cup \bound{}{+}x) + \set{ x }$, it follows that $\imath$ restricts to the identity on $\clset{{ x }}$.
Therefore, it suffices to prove that $\imath$ fixes all $x \in \maxel{U}$.

If $\dim{x} < n$, then $x \in \bound{}{\alpha}U$, and we have already proved $\imath(x) = x$.
Suppose then $\dim{x} = n$, and construct a path $y_- = y_0 \to x_0 \to \ldots \to y_m \to x_m = x$ in $\graph{U}$ as in Lemma \ref{lem:path_in_graph_of_molecule}.
Since $\imath$ preserves the covering relation and orientations, it maps this path to another path in $\graph{U}$.
We have $y_0 \in \bound{}{-}U$, so $\imath(y_0) = y_0$.
Suppose $\imath(y_i) = y_i$.
Since $y_i$ is a wire vertex in a directed graph with open edges, $x_i$ is the only node vertex with an edge from $y_i$, so necessarily $\imath(x_i) = x_i$.
If $i < m$, then $\imath$ is the identity on $\clset{{ x_i }}$, so $\imath(y_{i+1}) = y_{i+1}$.
Iterating until we reach $m$, we conclude.
\end{proof}

\begin{cor} \label{cor:molecule_unique_isomorphism}
Let $U$, $V$ be molecules.
If $U$ and $V$ are isomorphic, there exists a unique isomorphism $\varphi\colon U \incliso V$.
\end{cor}

\begin{comm} \index{$U \cp{k} V$}
It follows that, if $U$, $V$ are molecules, there is at most \emph{one} isomorphism $\varphi\colon \bound{k}{+}U \incliso \bound{k}{-}V$, so we can write
\begin{equation*}
    U \cp{k} V \eqdef U \cpiso{k}{\varphi} V,
\end{equation*}
and speak simply of \emph{the pasting of $U$ and $V$ at the $k$\nbd boundary}.
\end{comm}

\begin{cor} \label{cor:molecule_boundary_unique_isomorphism}
Let $U$, $V$ be round molecules, and suppose $\bound{}{\alpha}U$ and $\bound{}{\alpha}V$ are isomorphic for all $\alpha \in \set{ +, - }$.
Then there exists a unique isomorphism $\varphi\colon \bound{}{}U \incliso \bound{}{}V$ restricting to isomorphisms $\varphi^\alpha\colon \bound{}{\alpha} U \incliso \bound{}{\alpha} V$.
\end{cor}
\begin{proof}
By Corollary \ref{cor:molecule_unique_isomorphism}, the isomorphisms $\varphi^\alpha$ are uniquely defined.
They restrict to unique isomorphisms $\bound{}{\beta}(\bound{}{\alpha}U) \incliso \bound{}{\beta}(\bound{}{\alpha}V)$ for all $\beta \in \set{ +, - }$, which implies that the restrictions of $\varphi^-$ and $\varphi^+$ to $\bound{}{+}U \cap \bound{}{-}U = \bound{}{}(\bound{}{\alpha}U)$ are equal.
It follows that there is a unique extension of $\varphi^-, \varphi^+$ to a map $\varphi\colon \bound{}{}U \to \bound{}{}V$.
Since $V$ is also round, this map is injective, hence an isomorphism.
\end{proof}

\begin{comm} \index{$U \celto V$}
It follows that, if $U$, $V$ are round molecules, there is at most \emph{one} isomorphism $\varphi\colon \bound{}{}U \incliso \bound{}{}V$ restricting to $\varphi^\alpha\colon \bound{}{\alpha} U \incliso \bound{}{\alpha} V$ for all $\alpha \in \set{ +, - }$, so we can write

\begin{equation*}
    U \celto V \eqdef U \celto^\varphi V,
\end{equation*}
and speak simply of \emph{the rewrite of $U$ into $V$}.
\end{comm}

\begin{dfn}[Merger of a round molecule] \index{molecule!merger} \index{$\compos{U}$}
Let $U$ be a round molecule.
The \emph{merger} of $U$ is the atom $\compos{U} \eqdef \bound{}{-}U \celto \bound{}{+}U$.
\end{dfn}

\begin{exm}[The merger of a round molecule]
	Let $U$ be the oriented face poset of the pasting diagram
\[\begin{tikzcd}[column sep=small]
	\bullet && \bullet & \bullet \\
	& \bullet
	\arrow[curve={height=6pt}, from=1-1, to=2-2]
	\arrow[""{name=0, anchor=center, inner sep=0}, curve={height=12pt}, from=2-2, to=1-4]
	\arrow[""{name=1, anchor=center, inner sep=0}, curve={height=-12pt}, from=1-1, to=1-3]
	\arrow[from=1-3, to=1-4]
	\arrow[""{name=2, anchor=center, inner sep=0}, curve={height=-30pt}, from=1-1, to=1-4]
	\arrow[from=2-2, to=1-3]
	\arrow[shorten <=3pt, Rightarrow, from=0, to=1-3]
	\arrow[curve={height=-6pt}, shorten <=2pt, shorten >=5pt, Rightarrow, from=2-2, to=1]
	\arrow[shorten >=4pt, Rightarrow, from=1-3, to=2]
\end{tikzcd}\;,\]
	which is a round 2\nbd dimensional molecule.
	The merger $\compos{U}$ of $U$ is the oriented face poset of the pasting diagram
	\[\begin{tikzcd}[column sep=small]
	\bullet &&& \bullet \\
	& \bullet
	\arrow[curve={height=6pt}, from=1-1, to=2-2]
	\arrow[curve={height=12pt}, from=2-2, to=1-4]
	\arrow[""{name=0, anchor=center, inner sep=0}, curve={height=-24pt}, from=1-1, to=1-4]
	\arrow[curve={height=6pt}, shorten <=4pt, shorten >=8pt, Rightarrow, from=2-2, to=0]
\end{tikzcd} \;, \]
	obtained by replacing the interior of the first diagram with a single 2\nbd dimensional cell, while keeping its boundary fixed.
\end{exm}

\begin{lem} \label{lem:atom_merger_of_its_boundary}
Let $U$ be an atom, $\dim{U} > 0$.
Then $U$ is uniquely isomorphic to $\compos{U}$.
\end{lem}
\begin{proof}
By Lemma \ref{lem:atom_greatest_element}, $U$ was produced by (\textit{Atom}), so it is of the form $V \celto W$, and by Lemma \ref{lem:boundaries_of_rewrite} $\bound{}{-}U$ is isomorphic to $V$ and $\bound{}{+}U$ to $W$.
We conclude by Corollary \ref{cor:molecule_unique_isomorphism} and Corollary \ref{cor:molecule_boundary_unique_isomorphism}.
\end{proof}

\begin{prop} \label{prop:associativity_of_pasting}
Let $U$, $V$, $W$ be molecules and $k \in \mathbb{N}$ such that $U \cp{k} V$ and $V \cp{k} W$ are both defined.
Then $(U \cp{k} V) \cp{k} W$ and $U \cp{k} (V \cp{k} W)$ are both defined and uniquely isomorphic.
\end{prop}
\begin{proof}
By Lemma \ref{lem:pasting_top_boundary}, $\bound{k}{+}(U \cp{k} V)$ is isomorphic to $\bound{k}{+}V$, which is isomorphic to $\bound{k}{-}W$ since $V \cp{k} W$ is defined.
It follows that $(U \cp{k} V) \cp{k} W$ is defined.
A dual proof shows that $U \cp{k} (V \cp{k} W)$ is defined, and a routine argument shows that it satisfies the same universal property as $(U \cp{k} V) \cp{k} W$.
We conclude that the two are isomorphic, uniquely so by Corollary \ref{cor:molecule_unique_isomorphism}.
\end{proof}

\begin{comm} \label{comm:associativity_pasting_no_brackets}
Proposition \ref{prop:associativity_of_pasting} shows that pasting of molecules at the $k$\nbd boundary is associative up to unique isomorphism.
Thus, given any sequence $(\order{i}{U})_{i=1}^m$ of molecules and $k \in \mathbb{N}$ such that, for all $i \in \set{1, \ldots, m-1}$, the pasting $\order{i}{U} \cp{k} \order{i+1}{U}$ is defined, we may unambiguously write
\begin{equation*}
    \order{1}{U} \cp{k} \ldots \cp{k} \order{m}{U}
\end{equation*}
to stand for any binary bracketing up to unique isomorphism.
\end{comm}

\begin{prop} \label{prop:unitality_of_pasting}
Let $U$ be a molecule and $k \in \mathbb{N}$.
Then $U \cp{k} \bound{k}{+}U$ and $\bound{k}{-}U \cp{k} U$ are both defined and uniquely isomorphic to $U$.
\end{prop}
\begin{proof}
Since $\dim{\bound{k}{\alpha}U} \leq k$, we have $\bound{k}{\beta}(\bound{k}{\alpha}U) = \bound{k}{\alpha}U$ for all $\alpha, \beta \in \set{ +, - }$.
Moreover, by Lemma \ref{lem:molecules_are_globular} $\bound{k}{+}U$ is a molecule.
It follows that $U \cp{k} \bound{k}{+}U$ and $\bound{k}{-}U \cp{k} U$ are uniquely defined, and since the inclusion of $U$ in each of them is the pushout of an isomorphism, it is an isomorphism.
\end{proof}

\begin{prop} \label{prop:interchange_of_pasting}
Let $U, U', V, V'$ be molecules and $k < n \in \mathbb{N}$ such that $(U \cp{n} U') \cp{k} (V \cp{n} V')$ is defined.
Then $(U \cp{k} V) \cp{n} (U' \cp{k} V')$ is defined and uniquely isomorphic to $(U \cp{n} U') \cp{k} (V \cp{n} V')$.
\end{prop}
\begin{proof}
By Lemma \ref{lem:pasting_lower_boundary}, 
\begin{equation*}
    \bound{k}{+}(U \cp{n} U') = \bound{k}{+}U = \bound{k}{+}U', \quad \quad \bound{k}{-}(V \cp{n} V') = \bound{k}{-}V = \bound{k}{-}V'.
\end{equation*}
Since all of these are molecules, $U \cp{k} V$ and $U' \cp{k} V'$ are uniquely defined.
By Lemma \ref{lem:pasting_higher_boundary}, $\bound{n}{+}(U \cp{k} V)$ is isomorphic to $\bound{n}{+}U \cp{k} \bound{n}{+}V$, which, since $U \cp{n} U'$ and $V \cp{n} V'$ are defined, is isomorphic to $\bound{n}{-}U' \cp{k} \bound{n}{-}V'$, which is in turn isomorphic to $\bound{n}{-}(U' \cp{k} V')$.
It follows that $(U \cp{k} V) \cp{n} (U' \cp{k} V')$ is defined, and a routine argument shows that it satisfies the same universal property as $(U \cp{n} U') \cp{k} (V \cp{n} V')$.
\end{proof}

\clearpage 
\thispagestyle{empty}

%% file: layerings.tex
\chapter{Submolecules and layerings} \label{chap:layerings}
\thispagestyle{firstpage}

\begin{guide}
	In the paradigm of higher-dimensional rewriting, as originally formulated in the theory of polygraphs, a rewrite system consists of $(n+1)$\nbd dimensional cells, whose input and output boundaries are $n$\nbd dimensional diagrams.
	The basic computational step of such a rewrite system consists of matching the input $n$\nbd boundary of an $(n+1)$\nbd cell within a wider $n$\nbd dimensional pasting diagram, then ``applying'' the rewrite by substituting the output $n$\nbd boundary for the match.
	The rewrite is itself embodied by an $(n+1)$\nbd dimensional diagram, which can be seen as the result of gluing the $(n+1)$\nbd cell to the original $n$\nbd dimensional diagram at the location of the match.

	For example, a simple 2\nbd dimensional rewrite system, corresponding to a string rewrite system, may contain a 2\nbd dimensional cell
\[\begin{tikzcd}[sep=small]
	& \bullet \\
	\bullet && \bullet \\
	& \bullet
	\arrow["a", curve={height=-6pt}, from=2-1, to=1-2]
	\arrow["b", curve={height=-6pt}, from=1-2, to=2-3]
	\arrow["b"', curve={height=6pt}, from=2-1, to=3-2]
	\arrow["a"', curve={height=6pt}, from=3-2, to=2-3]
	\arrow[shorten <=3pt, shorten >=3pt, Rightarrow, from=3-2, to=1-2]
\end{tikzcd}\]
	embodying the substitution of the string $ab$ for the string $ba$.
	The input 1\nbd boundary of this cell can be matched twice in the string $babba$, represented by the 1\nbd dimensional pasting diagram
	\[\begin{tikzcd}
	\bullet & \bullet & \bullet & \bullet & \bullet & \bullet
	\arrow["b", from=1-1, to=1-2]
	\arrow["a", from=1-2, to=1-3]
	\arrow["b", from=1-3, to=1-4]
	\arrow["b", from=1-4, to=1-5]
	\arrow["a", from=1-5, to=1-6]
\end{tikzcd}\; ,\]
	and the corresponding rewrites are embodied by the 2\nbd dimensional pasting diagrams
\[\begin{tikzcd}[sep=small]
	& \bullet \\
	\bullet && \bullet & \bullet & \bullet & \bullet \\
	& \bullet
	\arrow["a", curve={height=-6pt}, from=2-1, to=1-2]
	\arrow["b", curve={height=-6pt}, from=1-2, to=2-3]
	\arrow["b"', curve={height=6pt}, from=2-1, to=3-2]
	\arrow["a"', curve={height=6pt}, from=3-2, to=2-3]
	\arrow[shorten <=3pt, shorten >=3pt, Rightarrow, from=3-2, to=1-2]
	\arrow["b", from=2-3, to=2-4]
	\arrow["b", from=2-4, to=2-5]
	\arrow["a", from=2-5, to=2-6]
\end{tikzcd}\;, \quad
\begin{tikzcd}[sep=small]
	&&&& \bullet \\
	\bullet & \bullet & \bullet & \bullet && \bullet \\
	&&&& \bullet
	\arrow["a", curve={height=-6pt}, from=2-4, to=1-5]
	\arrow["b", curve={height=-6pt}, from=1-5, to=2-6]
	\arrow["b"', curve={height=6pt}, from=2-4, to=3-5]
	\arrow["a"', curve={height=6pt}, from=3-5, to=2-6]
	\arrow[shorten <=3pt, shorten >=3pt, Rightarrow, from=3-5, to=1-5]
	\arrow["b", from=2-3, to=2-4]
	\arrow["a", from=2-2, to=2-3]
	\arrow["b", from=2-1, to=2-2]
\end{tikzcd}\;, \]
	whose output 1\nbd boundary is the result of each substitution.

	The ``gluing'' producing these diagrams is not, at the level of shapes, an instance of \emph{pasting} of molecules, as defined in the previous chapter: the match is not with the \emph{entire} 1\nbd dimensional boundary, but only a portion of it.
	The question is then: \emph{for what ``portions'' of a molecule is this operation well-defined, that is, it produces a well-formed molecule?}
	This has an evident answer in dimension 1 --- it is the ``linear subgraphs'' --- but looks increasingly complicated in higher dimensions.

	The general answer (Proposition \ref{prop:round_submolecule_substitution}) is: it is the images of molecules that can appear as factors in a ``pasting decomposition'' of a molecule, what we call its \cemph{submolecules}.
	However, this is far from a practical answer by itself: as we saw at the end of last chapter, pasting satisfies non-trivial equations up to isomorphism, so the problem of \emph{recognising submolecules} is not trivial.
	This chapter is largely devoted to the production of criteria for deciding this problem at least in the instance which is most relevant to higher-dimensional rewriting: \emph{round} submolecules of the same dimension as the ambient molecule, what we call \cemph{rewritable submolecules}.
	This culminates in Theorem \ref{thm:rewritable_submolecule_criterion}, which is the strongest fully general criterion that we know at the moment.

	Solving this problem involves understanding the ``space'' of possible pasting decompositions of a molecule.
	As it turns out, this can be largely reduced to understanding pasting decompositions of a special kind: those of the form
	\[
		U \simeq \order{1}{U} \cp{k} \ldots \cp{k} \order{m}{U},
	\]
	where the molecule $\order{i}{U}$ contains a \emph{single} maximal element of $U$ whose dimension is strictly larger than $k$, for each $i \in \set{1, \ldots, m}$.
	Such a decomposition is called a \cemph{$k$\nbd layering} of $U$.

	For example, if $U$ is the oriented face poset of the 2\nbd dimensional pasting diagram
\begin{equation} \label{eq:both_0_and1_layering}
	\begin{tikzcd}
	\bullet & \bullet & \bullet & \bullet
	\arrow[""{name=0, anchor=center, inner sep=0}, "0"', curve={height=12pt}, from=1-1, to=1-2]
	\arrow[""{name=1, anchor=center, inner sep=0}, "3", curve={height=-12pt}, from=1-1, to=1-2]
	\arrow["1", from=1-2, to=1-3]
	\arrow[""{name=2, anchor=center, inner sep=0}, "4", curve={height=-12pt}, from=1-3, to=1-4]
	\arrow[""{name=3, anchor=center, inner sep=0}, "2"', curve={height=12pt}, from=1-3, to=1-4]
	\arrow["1"', shorten <=3pt, shorten >=3pt, Rightarrow, from=3, to=2]
	\arrow["0", shorten <=3pt, shorten >=3pt, Rightarrow, from=0, to=1]
\end{tikzcd}\end{equation}
	then $U$ admits a single 0\nbd layering with layers
\begin{equation} \label{eq:exm_0_layering}
\begin{tikzcd}
	\bullet & \bullet & \bullet & \bullet & \bullet & \bullet
	\arrow[""{name=0, anchor=center, inner sep=0}, "0"', curve={height=12pt}, from=1-1, to=1-2]
	\arrow[""{name=1, anchor=center, inner sep=0}, "3", curve={height=-12pt}, from=1-1, to=1-2]
	\arrow[""{name=2, anchor=center, inner sep=0}, "4", curve={height=-12pt}, from=1-5, to=1-6]
	\arrow[""{name=3, anchor=center, inner sep=0}, "2"', curve={height=12pt}, from=1-5, to=1-6]
	\arrow["1", from=1-3, to=1-4]
	\arrow["1"', shorten <=3pt, shorten >=3pt, Rightarrow, from=3, to=2]
	\arrow["0", shorten <=3pt, shorten >=3pt, Rightarrow, from=0, to=1]
\end{tikzcd}\end{equation}
	and two 1\nbd layerings with layers
\begin{equation} \label{eq:exm_1_layering}
	\begin{tikzcd}
	\bullet & \bullet & \bullet & \bullet & \bullet & \bullet & \bullet & \bullet
	\arrow[""{name=0, anchor=center, inner sep=0}, "0"', curve={height=12pt}, from=1-1, to=1-2]
	\arrow[""{name=1, anchor=center, inner sep=0}, "3", curve={height=-12pt}, from=1-1, to=1-2]
	\arrow["2", from=1-3, to=1-4]
	\arrow["1", from=1-2, to=1-3]
	\arrow["3", from=1-5, to=1-6]
	\arrow["1", from=1-6, to=1-7]
	\arrow[""{name=2, anchor=center, inner sep=0}, "4", curve={height=-12pt}, from=1-7, to=1-8]
	\arrow[""{name=3, anchor=center, inner sep=0}, "2"', curve={height=12pt}, from=1-7, to=1-8]
	\arrow["0", shorten <=3pt, shorten >=3pt, Rightarrow, from=0, to=1]
	\arrow["1", shorten <=3pt, shorten >=3pt, Rightarrow, from=3, to=2]
\end{tikzcd}\;,\end{equation}
\begin{equation} \label{eq:exm_1_layering_bis}
	\begin{tikzcd}
	\bullet & \bullet & \bullet & \bullet & \bullet & \bullet & \bullet & \bullet
	\arrow[""{name=0, anchor=center, inner sep=0}, "0"', curve={height=12pt}, from=1-5, to=1-6]
	\arrow[""{name=1, anchor=center, inner sep=0}, "3", curve={height=-12pt}, from=1-5, to=1-6]
	\arrow["4", from=1-7, to=1-8]
	\arrow["1", from=1-6, to=1-7]
	\arrow["0", from=1-1, to=1-2]
	\arrow["1", from=1-2, to=1-3]
	\arrow[""{name=2, anchor=center, inner sep=0}, "4", curve={height=-12pt}, from=1-3, to=1-4]
	\arrow[""{name=3, anchor=center, inner sep=0}, "2"', curve={height=12pt}, from=1-3, to=1-4]
	\arrow["0", shorten <=3pt, shorten >=3pt, Rightarrow, from=0, to=1]
	\arrow["1", shorten <=3pt, shorten >=3pt, Rightarrow, from=3, to=2]
\end{tikzcd}\;,
\end{equation}
	respectively.
	Notice that, when $k = \dim{U} - 1$ (for example, $k = 1$ in this example), then $k$\nbd layerings have a rewriting-theoretic interpretation, as possible \emph{sequentialisations} of a diagram into sequences of individual rewrites.
	
	We will show that every molecule $U$ admits a $k$\nbd layering for some $k$ (Theorem 
	\ref{thm:molecules_admit_layerings}), and that the \emph{smallest} $k$ for which this happens always falls between two numerical values associated to $U$: its \cemph{frame dimension} and its \cemph{layering dimension} (Corollary 
	\ref{cor:least_dimension_with_layering}).

	Every $k$\nbd layering of $U$ determines a linear order on maximal elements of $U$ of dimension strictly greater than $k$ --- a \cemph{$k$\nbd ordering} --- and can be uniquely reconstructed from this linear order up to layer-wise isomorphism.
	For example, (\ref{eq:exm_0_layering}) corresponds to the 0\nbd ordering
	\[
		(2, 0), (1, 1), (2, 1)
	\]
	while (\ref{eq:exm_1_layering}) and (\ref{eq:exm_1_layering_bis}) correspond to the 1\nbd orderings
	\[
		(2, 0), (2, 1), \quad \quad (2, 1), (2, 0),
	\]
	respectively.
	The question of what $k$\nbd layerings are possible can thus be turned into the question of what $k$\nbd orderings determine valid $k$\nbd layerings.
	This turns out to be a complex question, for which we only have a neat answer in special cases, as we will see in Chapter \ref{chap:acyclic} --- all the more fascinating since it seems to be an intrinsically \emph{directed} problem, that does not have a counterpart in traditional combinatorial topology.
\end{guide}


\section{Submolecules and substitution} \label{sec:submolecules}

\begin{guide}
	In this section, we define \emph{submolecules} and prove some of their basic properties.
	We then define the operation of \cemph{substitution} of round molecules included into an oriented graded poset of the same dimension.
\end{guide}

\begin{dfn}[Submolecule inclusions] \index{submolecule!inclusion}
The class of \emph{submolecule inclusions} is the smallest subclass of inclusions of molecules such that
\begin{enumerate}
    \item all isomorphisms are submolecule inclusions,
    \item for all molecules $U$, $V$ and all $k \in \mathbb{N}$ such that the pasting $U \cp{k} V$ is defined, $U \incl (U \cp{k} V)$ and $V \incl (U \cp{k} V)$ are submolecule inclusions,
    \item the composite of two submolecule inclusions is a submolecule inclusion.
\end{enumerate}
\end{dfn}

\begin{dfn}[Submolecules] \index{molecule!submolecule|see {submolecule}} \index{submolecule} \index{$V \submol U$}
Let $U$ be a molecule and $V \subseteq U$ a closed subset.
We say that $V$ is a \emph{submolecule} of $U$, and write $V \submol U$, if the inclusion of $V$ into $U$ is a submolecule inclusion.
\end{dfn}

\begin{rmk}
An inclusion $\imath\colon V \incl U$ is a submolecule inclusion if and only if $\imath(V) \submol U$.
\end{rmk}

\begin{lem} \label{lem:boundary_is_submolecule}
Let $U$ be a molecule, $n \in \mathbb{N}$, $\alpha \in \set{ +, - }$.
Then $\bound{n}{\alpha}U \submol U$.
\end{lem}
\begin{proof}
By Lemma \ref{lem:molecules_are_globular}, $\bound{n}{\alpha} U$ is a molecule.
By Proposition \ref{prop:unitality_of_pasting}, the pastings $U \cp{n} \bound{n}{+}U$ and $\bound{n}{-}U \cp{n} U$ are both defined and uniquely isomorphic to $U$.
The inclusion of $\bound{n}{-}U$ into $U$ factors as the inclusion $\bound{n}{-}U \incl (\bound{n}{-}U \cp{n} U)$ followed by an isomorphism, and the inclusion of $\bound{n}{+}U$ factors as the inclusion $\bound{n}{+}U \incl (U \cp{n} \bound{n}{+}U)$ followed by an isomorphism.
\end{proof}

\begin{lem} \label{lem:downset_is_submolecule}
Let $U$ be a molecule, $x \in U$.
Then $\clset{{ x }} \submol U$.
\end{lem}
\begin{proof}
By Lemma \ref{lem:all_downsets_are_atoms}, $\clset{{ x }}$ is a molecule.
We proceed by induction on the construction of $U$.
If $U$ was produced by (\textit{Point}), then $x$ must be the unique element of $U$, so $\clset{{ x }} = U$.
If $U$ was produced by (\textit{Paste}), it splits into $V \cup W$ with $V$, $W \submol U$, and $x \in V$ or $x \in W$.
By the inductive hypothesis, $\clset{{ x }} \submol V$ or $\clset{{ x }} \submol W$.
If $U$ was produced by (\textit{Atom}), it is equal to $(V \cup W) + \set{ \top }$ with $V$, $W \submol U$ by Lemma \ref{lem:boundary_is_submolecule}, and either $x \in V$ or $x \in W$, in which case the inductive hypothesis applies, or $x = \top$, and $\clset{{ x }} = U$.
\end{proof}

\begin{lem} \label{lem:intersection_with_boundary}
Let $V \submol U$ be molecules, $n \in \mathbb{N}$, $\alpha \in \set{ +, - }$.
Then $V \cap \bound{n}{\alpha}U \subseteq \bound{n}{\alpha}V$.
\end{lem}
\begin{proof}
We proceed by induction on the construction of the inclusion $\imath\colon V \incl U$.
If $\imath$ is an isomorphism, that is, $V = U$, the statement is immediate.

Suppose $\imath$ is the inclusion of $V$ into a pasting.
Then there exist $W \submol U$ and $k \in \mathbb{N}$ such that $U$ splits into $V \cup W$ or $W \cup V$ along the $k$\nbd boundary.
Suppose without loss of generality that $U$ splits into $V \cup W$.
If $n > k$, by Lemma \ref{lem:pasting_higher_boundary} we have $\bound{n}{\alpha}U = \bound{n}{\alpha}V \cup \bound{n}{\alpha}W$ and by globularity $V \cap W = \bound{k}{+}V \subseteq \bound{n}{\alpha}V$.
Then
\begin{equation*}
    V \cap \bound{n}{\alpha}U = (V \cap \bound{n}{\alpha}V) \cup (V \cap \bound{n}{\alpha}W) \subseteq \bound{n}{\alpha}V.
\end{equation*}
If $n = k$, by Lemma \ref{lem:pasting_top_boundary} we have
\begin{align*}
    V \cap \bound{n}{-}U & = V \cap \bound{n}{-}V = \bound{n}{-}V, \\
    V \cap \bound{n}{+}U & = V \cap \bound{n}{+}W \subseteq V \cap W = \bound{n}{+}V.
\end{align*}
Finally, if $n < k$, by Lemma \ref{lem:pasting_lower_boundary} we have $\bound{n}{\alpha}U = \bound{n}{\alpha}V$.

Suppose $\imath$ is a composite of two submolecule inclusions, exhibiting $V \submol W$ and $W \submol U$.
Then
\begin{equation*}
    V \cap \bound{n}{\alpha}U = V \cap W \cap \bound{n}{\alpha}U \subseteq V \cap \bound{n}{\alpha}W \subseteq \bound{n}{\alpha}V
\end{equation*}
using the inductive hypothesis twice, and we conclude.
\end{proof}

\begin{comm} \label{comm:induction_on_submolecules} \index{submolecule!induction}
Let $\property{P}$ be a property of molecules such that, whenever $\mathrm{P}$ holds of a molecule $U$, then $\mathrm{P}$ holds of every submolecule $V \submol U$.
Because every proper submolecule of $U$ has strictly fewer elements than $U$, the submolecule relation on $U$ is well-founded, and its minimal elements are the 0\nbd dimensional one-element subsets $\set{x} \submol U$ for each $x \in \grade{0}{U}$.

If we want to prove that $\property{P}$ implies $\property{Q}$ for all molecules, we can then proceed by \emph{induction on submolecules}: assume that a molecule $U$ satisfies $\property{P}$, then
\begin{itemize}
    \item prove that $\set{x}$ satisfies $\property{Q}$ for all $x \in \grade{0}{U}$,
    \item prove that $U$ satisfies $\property{Q}$ under the assumption that every proper submolecule $V \submol U$ satisfies $\property{Q}$.
\end{itemize}
\end{comm}

\begin{lem} \label{lem:submolecule_rewrite}
Let $V$ be a molecule, $n < \dim{V}$, $\alpha \in \set{ +, - }$.
Consider a pushout diagram of the form
\[\begin{tikzcd}
	{\bound{n}{\alpha} V} && V \\
	U && {V \cup U}
	\arrow[hook, from=1-1, to=1-3]
	\arrow["\imath", hook', from=1-1, to=2-1]
	\arrow["{j_U}", hook, from=2-1, to=2-3]
	\arrow["{j_V}", hook', from=1-3, to=2-3]
	\arrow["\lrcorner"{anchor=center, pos=0.125, rotate=180}, draw=none, from=2-3, to=1-1]
\end{tikzcd}\]
in $\ogpos$.
If $\dim{U} = n$ and $\imath$ is a submolecule inclusion, then
\begin{enumerate}
    \item $V \cup U$ is a molecule,
    \item $j_U$ maps $U$ onto $\bound{n}{\alpha}(V \cup U)$,
    \item $j_V(V) \submol V \cup U$ and $j_V(\bound{n}{-\alpha}V) \submol \bound{n}{-\alpha}(V \cup U)$.
\end{enumerate}
\end{lem}
\begin{proof}
By induction on the construction of $\imath$.
If $\imath$ is an isomorphism, then $j_V$ is also an isomorphism, and all the statements are trivially satisfied.

Suppose $U$ is of the form $\bound{n}{\alpha}V \cp{k} W$ for some $k \in \mathbb{N}$, and $\imath$ is the inclusion of $\bound{n}{\alpha}V$ into the pasting.
Since $\dim{U} = n$, necessarily $\dim{W} \leq n$, so $\bound{n}{\alpha}W = W$ by Lemma \ref{lem:dimension_of_boundary}.
If $k \geq n$, then also $k \geq \dim{W}$, and in this case $\imath$ and $j_V$ are again isomorphisms.
Suppose that $k < n$.
Identifying $V$ with its image through $j_V$, $V \cup U$ splits into $V \cup W$ with
\begin{equation*}
    V \cap W = \bound{n}{\alpha}V \cap W = \bound{k}{-}W = \bound{k}{+}(\bound{n}{\alpha}V) = \bound{k}{+} V
\end{equation*}
where the final equation uses globularity of $V$.
This exhibits $V \cup U$ as $V \cp{k} W$, with $j_V$ the inclusion of $V$ into the pasting, and $j_U$ maps $\bound{n}{\alpha}V \cp{k} W$ onto $\bound{n}{\alpha}(V \cp{k} W)$ by Lemma \ref{lem:pasting_higher_boundary}.
By the same result, $\bound{n}{-\alpha}V \submol \bound{n}{-\alpha}(V \cp{k} W)$.
The case where $U$ is of the form $W \cp{k} \bound{n}{\alpha}V$ is dual.

By the pasting law for pushout squares, if the statement is true of two submolecule inclusions, it is also true of their composite.
\end{proof}

\begin{dfn}[Substitution] \index{submolecule!substitution} \index{$\subs{U}{W}{\imath(V)}$}
Let $U$ be an oriented graded poset, $V$, $W$ be round molecules with $\dim U = \dim V = \dim W$, $\imath\colon V \incl U$ an inclusion, and suppose that $V \celto W$ is defined.
Consider the pushout
\begin{equation} \label{eq:substitution}
	\begin{tikzcd}
	V & V \celto W \\
	U & U \cup (V \celto W)
	\arrow[hook, from=1-1, to=1-2]
	\arrow["\imath", hook', from=1-1, to=2-1]
	\arrow[hook, from=2-1, to=2-2]
	\arrow[hook', from=1-2, to=2-2]
	\arrow["\lrcorner"{anchor=center, pos=0.125, rotate=180}, draw=none, from=2-2, to=1-1]
\end{tikzcd}
\end{equation}
in $\ogpos$.
The \emph{substitution of $W$ for $\imath\colon V \incl U$} is the oriented graded poset
\[
    \subs{U}{W}{\imath(V)} \eqdef \bound{}{+} (U \cup (V \celto W)).
\]
\end{dfn}

\begin{comm}
When $\imath$ is the inclusion of a closed subset $V \subseteq U$, we will write simply $\subs{U}{W}{V}$.
\end{comm}

\begin{exm}[A substitution]
	Let $U$ and $W$ be $\thearrow{} \cp{0} \thearrow{}$, and let $V \incl U$ be the inclusion of $\thearrow{}$ as the second factor of $U$.
	Then the pushout square (\ref{eq:substitution}) looks like
\[
	\begin{tikzcd}[sep=small]
	&&&&&&& \bullet \\
	& \bullet & {} & \bullet &&& \bullet & {} & \bullet \\
	&&&&&&&&& {} \\
	\bullet &&&&& \bullet && \bullet \\
	& \bullet & {} & \bullet &&& \bullet & {} & \bullet
	\arrow[curve={height=6pt}, from=5-2, to=5-4]
	\arrow[curve={height=6pt}, from=4-1, to=5-2]
	\arrow[curve={height=6pt}, from=2-2, to=2-4]
	\arrow[""{name=0, anchor=center, inner sep=0}, curve={height=6pt}, from=2-7, to=2-9]
	\arrow[curve={height=-6pt}, from=2-7, to=1-8]
	\arrow[curve={height=-6pt}, from=1-8, to=2-9]
	\arrow[curve={height=6pt}, from=4-6, to=5-7]
	\arrow[""{name=1, anchor=center, inner sep=0}, curve={height=6pt}, from=5-7, to=5-9]
	\arrow[curve={height=-6pt}, from=5-7, to=4-8]
	\arrow[curve={height=-6pt}, from=4-8, to=5-9]
	\arrow[color=\mycolor, shorten <=11pt, shorten >=11pt, hook', from=2-3, to=5-3]
	\arrow[color=\mycolor, shorten <=29pt, shorten >=43pt, hook, from=2-3, to=2-8]
	\arrow[color=\mycolor, shorten <=11pt, shorten >=27pt, hook', from=2-8, to=5-8]
	\arrow[color=\mycolor, shorten <=29pt, shorten >=43pt, hook, from=5-3, to=5-8]
	\arrow[shorten <=5pt, Rightarrow, from=0, to=1-8]
	\arrow[shorten <=5pt, Rightarrow, from=1, to=4-8]
\end{tikzcd}
\]
	in terms of pasting diagrams.
	The substitution of $W$ for $V \incl U$ is the output boundary of the diagram in the bottom right corner, which is isomorphic to $\thearrow{} \cp{0} \thearrow{} \cp{0} \thearrow{}$.
\end{exm}

\begin{lem} \label{lem:full_substitution}
Let $U$, $V$ be round molecules such that $U \celto V$ is defined.
Then $\subs{U}{V}{U}$ is isomorphic to $V$.
\end{lem}
\begin{proof}
The inclusion of $U \celto V$ into $U \cup (U \celto V)$ is an isomorphism,
and we conclude by Lemma \ref{lem:boundaries_of_rewrite}.
\end{proof}

\begin{lem} \label{lem:substitution_direct_pushout}
Let $U$ be an oriented graded poset, $V$, $W$ round molecules, and $\imath\colon V \incl U$ an inclusion such that $\subs{U}{W}{\imath(V)}$ is defined.
Let $\varphi\colon \bound{}{}V \incliso \bound{}{}W$ be the isomorphism used in the construction of $V \celto W$.
Then $\subs{U}{W}{\imath(V)}$ can be constructed as the pushout
\begin{equation} \label{eq:substitution_direct_pushout}
\begin{tikzcd}
	\bound{}{}V & \bound{}{}W & W \\
	U \setminus \inter{\imath(V)} && \subs{U}{W}{\imath(V)}.
	\arrow["\varphi", hook, from=1-1, to=1-2]
	\arrow[hook, from=1-2, to=1-3]
	\arrow[hook', from=1-1, to=2-1]
	\arrow[hook, from=2-1, to=2-3]
	\arrow["j", hook', from=1-3, to=2-3]
	\arrow["\lrcorner"{anchor=center, pos=0.125, rotate=180}, draw=none, from=2-3, to=1-1]
\end{tikzcd}
\end{equation}
\end{lem}
\begin{proof}
We can safely identify $V$ with its image through $\imath$, and treat it as a closed subset of $U$.
First of all, observe that $U \setminus \inter{V}$ is the complement of the complement of a closed subset in a closed subset, so it is closed in $U$, and well-defined as an oriented graded poset.

Let $n \eqdef \dim U$, so $\dim (U \cup (V \celto W)) = \dim (V \celto W) = n + 1$.
Then
\begin{equation*}
    \faces{n}{+}(V \celto W) = \grade{n}{W}, \quad \quad \faces{n}{+}U = \grade{n}{U}
\end{equation*}
and since $U \cap (V \celto W) = V$, by Lemma \ref{lem:faces_of_union}
\begin{equation*}
	\faces{n}{+}(U \cup (V \celto W)) = \grade{n}{W} + (\grade{n}{U} \setminus \grade{n}{V}) = \grade{n}{W} + \grade{n}{(U \setminus \inter{V})},
\end{equation*}
while for all $k < n$
\begin{equation*}
	\grade{k}{(\maxel{(U \cup (V \celto W))})} = \grade{k}{(\maxel{U})} = \grade{k}{(\maxel{(U \setminus \inter{V}}))}
\end{equation*}
because both $V$ and $V \celto W$ are round, hence pure by Lemma \ref{lem:round_is_pure}, and do not contain any maximal elements of dimension $k$.

It follows that $\bound{}{+}(U \cup (V \celto W))$ is the union of $W$ and $U \setminus \inter{V}$, with intersection $\bound{}{}W = \bound{}{}V$.
\end{proof}

\begin{lem} \label{lem:revert_substitution}
Let $U$ be an oriented graded poset, $V$, $W$ round molecules, let $\imath\colon V \incl U$ be an inclusion such that the substitution $\subs{U}{W}{\imath(V)}$ is defined, and let $j\colon W \incl \subs{U}{W}{\imath(V)}$ be the right side of (\ref{eq:substitution_direct_pushout}).
Then $\subs{(\subs{U}{W}{\imath(V)})}{V}{j(W)}$ is defined and isomorphic to $U$.
\end{lem}
\begin{proof}
Since $W \celto V$ is defined whenever $V \celto W$ is defined, it follows that $\subs{(\subs{U}{W}{\imath(V)})}{V}{j(W)}$ is defined.
The isomorphism with $U$ is straightforward algebra of closed subsets using Lemma \ref{lem:substitution_direct_pushout} twice.
\end{proof}

\begin{lem} \label{lem:substitution_preserves_boundaries}
Let $U$, $V$, $W$ be molecules, $k < \dim{U}$, $\alpha \in \set{+, -}$, and let $\imath\colon V \incl U$ be a submolecule inclusion such that $\subs{U}{W}{\imath(V)}$ is defined.
Then $\bound{k}{\alpha}U$ is isomorphic to $\bound{k}{\alpha}(\subs{U}{W}{\imath(V)})$.
\end{lem}
\begin{proof}
By Lemma \ref{lem:submolecule_rewrite}, $U \cup (V \celto W)$ is a molecule and $U$ is isomorphic to its input boundary.
By globularity, $\bound{k}{\alpha}U$ is isomorphic to 
\begin{equation*}
    \bound{k}{\alpha}(\bound{}{+}(U \cup (V \celto W))) = \bound{k}{\alpha}(\subs{U}{W}{\imath(V)}). \qedhere
\end{equation*}
\end{proof}

\begin{lem} \label{lem:pasting_after_substitution}
Let $U, V, W, U', U''$ be molecules, $k < \dim{U}$, and $\imath\colon V \incl U$ be a submolecule inclusion such that $U \cp{k} U'$, $U'' \cp{k} U$, and $\subs{U}{W}{\imath(V)}$ are defined.
Then 
\begin{enumerate}
    \item $\subs{U}{W}{\imath(V)} \cp{k} U'$ and $U'' \cp{k} \subs{U}{W}{\imath(V)}$ are defined,
    \item if $\dim{U'} \leq \dim{U}$, then $\subs{(U \cp{k} U')}{W}{\imath_U(\imath(V))}$ is defined and isomorphic to $\subs{U}{W}{\imath(V)} \cp{k} U'$,
    \item if $\dim{U''} \leq \dim{U}$, then $\subs{(U'' \cp{k} U)}{W}{\imath_U(\imath(V))}$ is defined and isomorphic to $U'' \cp{k} \subs{U}{W}{\imath(V)}$.
\end{enumerate}
\end{lem}
\begin{proof}
Lemma \ref{lem:substitution_preserves_boundaries} implies that $\subs{U}{W}{\imath(V)} \cp{k} U'$ and $U'' \cp{k} \subs{U}{W}{\imath(V)}$ are defined.
The substitution $\subs{(U \cp{k} U')}{W}{\imath_U(\imath(V))}$ is then defined if and only if $\dim{(U \cp{k} U')} = \dim{U}$, equivalently, if and only if $\dim{U'} \leq \dim{U}$.
Similarly, $\subs{(U'' \cp{k} U)}{W}{\imath_U(\imath(V))}$ is defined if and only if $\dim{U''} \leq \dim{U}$.
The isomorphisms follow straightforwardly from the definitions using the pasting law for pushout squares.
\end{proof}

\begin{dfn}[Multiple substitution] \index{submolecule!substitution!multiple}
Let $U$ be an oriented graded poset, and let $(\order{i}{V}, \order{i}{W}, j_i\colon \order{i}{V} \incl U)_{i=1}^m$ be a finite family of triples of
\begin{itemize}
    \item round molecules $\order{i}{V}, \order{i}{W}$ such that $\dim{U} = \dim{\order{i}{V}} = \dim{\order{i}{W}}$ and $\order{i}{V} \celto \order{i}{W}$ is defined, together with
    \item an inclusion $j_i\colon \order{i}{V} \incl U$.
 \end{itemize}
 Suppose that, for all $k, \ell \in \set{1, \ldots, m}$, if $k \neq \ell$, then
 \begin{equation*}
    j_k(\order{k}{V}) \cap j_{\ell}(\order{\ell}{V}) = j_k(\bound{}{}\order{k}{V}) \cap j_{\ell}(\bound{}{}\order{\ell}{V}).
\end{equation*}
Then for all $k, \ell \in \set{1, \ldots, m}$, if $k \neq \ell$, the image of $j_\ell$ is a subset of $U \setminus (j_k(\order{k}{V}) \setminus j_k(\bound{}{}\order{k}{V}))$, which by Lemma \ref{lem:substitution_direct_pushout} can be included into the substitution $\subs{U}{\order{k}{W}}{j_k(\order{k}{V})}$.
We let $j_\ell$ also denote the composite inclusion $\order{\ell}{V} \incl \subs{U}{\order{k}{W}}{j_k(\order{k}{V})}$.

The \emph{multiple substitution of $(\order{i}{W})_{i=1}^m$ for $(j_i\colon \order{i}{V} \incl U)_{i=1}^m$} is defined recursively by
\begin{align*}
    \subs{U}{\order{i}{W}}{j_i(\order{i}{V})}_{i=1}^0 & \eqdef U, \\
    \subs{U}{\order{i}{W}}{j_i(\order{i}{V})}_{i=1}^k &  \eqdef 
    \subs{
        \subs{U}{\order{i}{W}}{j_i(\order{i}{V})}_{i=1}^{k-1}}
        {\order{k}{W}}{j_k(\order{k}{V})}, \quad k > 0.
\end{align*}
\end{dfn}


\section{Layerings} \label{sec:layerings}

\begin{guide}
	In this section, we present the definition and some fundamental facts about layerings of molecules, culminating in Theorem \ref{thm:molecules_admit_layerings}, asserting that a $k$\nbd layering exists if $k$ is at least the \emph{layering dimension} of a molecule.
	
	The results of this section justify the proof method of \cemph{induction on layering dimension}, which we will be using frequently in the rest of the book.
	We illustrate the utility of layering-based methods with a succinct classification proof for 1\nbd dimensional molecules (Lemma \ref{lem:only_1_molecules}).
\end{guide}

\begin{dfn}[Layering] \index{molecule!layering}
Let $U$ be a molecule, $-1 \leq k < \dim{U}$, and
\begin{equation*}
    m \eqdef \size{\bigcup_{i > k} \grade{i}{(\maxel{U})}}.
\end{equation*}
A \emph{$k$\nbd layering of $U$} is a sequence $(\order{i}{U})_{i=1}^m$ of molecules such that $U$ is isomorphic to
\begin{equation*}
    \order{1}{U} \cp{k} \ldots \cp{k} \order{m}{U}
\end{equation*}
and $\dim{\order{i}{U}} > k$ for all $i \in \set{ 1, \ldots, m }$.
\end{dfn}

\begin{comm} \label{comm:minus_one_layering}
For $k = -1$, it is implied that $m = 1$, and $U$ is an atom.
\end{comm}

\begin{comm}
We will regularly identify the molecules in a layering of $U$ with their isomorphic images in $U$.
\end{comm}

\begin{lem} \label{lem:layering_intersections}
Let $U$ be a molecule, $-1 \leq k < \dim{U}$, and $(\order{i}{U})_{i=1}^m$ a $k$\nbd layering of $U$.
For all $i < j \in \set{1, \ldots, m}$,
\begin{equation*}
    \order{i}{U} \cap \order{j}{U} = \bound{k}{+}\order{i}{U} \cap \bound{k}{-}\order{j}{U}.
\end{equation*}
\end{lem}
\begin{proof}
Let $i < j \in \set{1, \ldots, m}$, and
\begin{align*}
    V & \eqdef \order{1}{U} \cp{k} \ldots \cp{k} \order{i}{U}, \\
    W & \eqdef \bound{k}{+}\order{i}{U} \cp{k} \order{i+1}{U} \cp{k} \ldots \cp{k} \order{j-1}{U}, \\
    Z & \eqdef \order{j}{U} \cp{k} \ldots \cp{k} \order{m}{U}.
\end{align*}
Then $U$ splits into $V \cup (W \cp{k} Z)$ along the $k$\nbd boundary, so
\begin{equation*}
    \bound{k}{+}\order{i}{U} = \bound{k}{+}V = \bound{k}{-}(W \cp{k} Z) = V \cap (W \cp{k} Z).
\end{equation*}
Since $\order{i}{U} \subseteq V$ and $\order{j}{U} \subseteq (W \cp{k} Z)$, it follows that $\order{i}{U} \cap \order{j}{U} \subseteq \bound{k}{+}\order{i}{U}$.
Dually, from the fact that $U$ splits into $(V \cp{k} W) \cup Z$ along the $k$\nbd boundary, we derive $\order{i}{U} \cap \order{j}{U} \subseteq \bound{k}{-}\order{j}{U}$.
\end{proof}

\begin{lem} \label{lem:layering_basic_properties}
Let $U$ be a molecule, $k < \dim{U}$, and let $(\order{i}{U})_{i=1}^m$ be a $k$\nbd layering of $U$.
Then, for all $i \in \set{ 1, \ldots, m }$, $\order{i}{U}$ contains a single maximal element of dimension $> k$.
\end{lem}
\begin{proof}
Because $\dim{\order{i}{U}} > k$, each $\order{i}{U}$ contains at least one maximal element of dimension $> k$, and because
\begin{equation*}
    \dim{( \order{i}{U} \cap \order{j}{U} )} = \dim{(\bound{k}{+}\order{i}{U} \cap \bound{k}{-}\order{j}{U})} \leq k
\end{equation*}
by Lemma \ref{lem:layering_intersections}, no such maximal element is contained in two of them.
Since there are exactly $m$ maximal elements of dimension $> k$, it follows that each $\order{i}{U}$ contains exactly one of them.
\end{proof}

\begin{lem} \label{lem:higher_layerings}
Let $U$ be a molecule, $k \leq \ell < \dim{U}$.
If $U$ admits a $k$\nbd layering, then $U$ admits an $\ell$\nbd layering.
\end{lem}
\begin{proof}
Let $(\order{i}{U})_{i=1}^m$ be a $k$\nbd layering of $U$.
For all $i \in \set{1, \ldots, m}$, let
\begin{equation*}
    \order{i}{V} \eqdef \bound{\ell}{+}\order{1}{U} \cp{k} \ldots \cp{k} \bound{\ell}{+}\order{i-1}U \cp{k} \order{i}{U} \cp{k} \bound{\ell}{-}\order{i+1}{U} \cp{k} \ldots \cp{k} \bound{\ell}{-}\order{m}U.
\end{equation*}
By repeated applications of Proposition \ref{prop:interchange_of_pasting} followed by Proposition \ref{prop:unitality_of_pasting}, $U$ is isomorphic to
\begin{equation*}
    \order{1}{V} \cp{\ell} \ldots \cp{\ell} \order{m}{V}.
\end{equation*}
Restricting to the subsequence of $(\order{i}{V})_{i=1}^m$ on those $i \in \set{1,\ldots,m}$ such that $\dim{\order{i}{V}} > \ell$, which does not change the result by Proposition \ref{prop:unitality_of_pasting}, we obtain an $\ell$\nbd layering of $U$.
\end{proof}

\begin{dfn}[Layering dimension] \index{molecule!layering dimension} \index{$\lydim{U}$} \index{dimension!layering}
Let $U$ be a molecule.
The \emph{layering dimension} of $U$ is the integer
\begin{equation*}
    \lydim{U} \eqdef \min \set{ k \geq -1 \mid \size{\bigcup_{i > k+1} \grade{i}{(\maxel{U})}} \leq 1 }.
\end{equation*}
\end{dfn}

\begin{lem} \label{lem:layering_dimension_smaller_than_dimension}
Let $U$ be a molecule, $n \eqdef \dim{U}$.
Then
\begin{enumerate}
    \item $\lydim{U} \leq n - 1$,
    \item $\lydim{U} = n - 1$ if and only if $\size{\grade{n}{U}} > 1$.
\end{enumerate}
\end{lem}
\begin{proof}
We have 
\begin{equation*}
    \size{\bigcup_{i > n} \grade{i}{(\maxel{U})}} = \size{\varnothing} = 0,
\end{equation*}
so $\lydim{U} \leq n - 1$, with equality if and only if
\begin{equation*}
    \size{\bigcup_{i > n-1} \grade{i}{(\maxel{U})}} = \size{\grade{n}{(\maxel{U})}} = \size{\grade{n}{U}} > 1,
\end{equation*}
where we used Lemma \ref{lem:maximal_vs_faces}.
\end{proof}

\begin{lem} \label{lem:layering_dimension_atom}
Let $U$ be a molecule.
Then $\lydim{U} = -1$ if and only if $U$ is an atom.
\end{lem}
\begin{proof}
Suppose $\lydim{U} = -1$.
Either $\size{\bigcup_{i > 0} \grade{i}{(\maxel{U})}} = 0$, so $\dim{U} = 0$ and we conclude by Lemma \ref{lem:only_0_molecule}, or $1 = \size{\bigcup_{i > 0} \grade{i}{(\maxel{U})}} = \size{\maxel{U}}$ by Lemma \ref{lem:maximal_0dim}.
In either case, $U$ has a greatest element.
Conversely, if $U$ has a greatest element, $\size{\bigcup_{i > 0} \grade{i}{(\maxel{U})}} \leq \size{\maxel{U}} = 1$.
\end{proof}

\begin{lem} \label{lem:layering_dimensions_pasting}
Let $U$, $V$ be molecules and $k < \min \set{ \dim{U}, \dim{V} }$ such that $U \cp{k} V$ is defined.
Then
\begin{equation*}
    \lydim{(U \cp{k} V)} \geq \max \set{ \lydim{U}, \lydim{V}, k }.
\end{equation*}
\end{lem}
\begin{proof}
Identifying $U$ and $V$ with their isomorphic images, $U \cp{k} V$ splits into $U \cup V$ with $\dim{(U \cap V)} = \dim{\bound{k}{+}U} = k$.
By Lemma \ref{lem:faces_of_union}, for all $i > k$,
\begin{equation*}
    \grade{i}{ ( \maxel{(U \cp{k} V)} ) } = \grade{i}{(\maxel{U})} + \grade{i}{(\maxel{V})},
\end{equation*}
and since $k < \min \set{ \dim{U}, \dim{V} }$, both $U$ and $V$ have at least one maximal element of dimension strictly larger than $k$.
It follows that
\begin{equation*}
    \size{ \bigcup_{i > k} \grade{i}{ ( \maxel{ (U \cp{k} V) } ) } } = 
    \size{ \bigcup_{i > k} \grade{i}{ ( \maxel{U} ) } } + \size{ \bigcup_{i > k} \grade{i}{ ( \maxel{V} ) } } \geq 2,
\end{equation*}
so $k - 1 < \lydim{(U \cp{k} V)}$, that is, $k \leq \lydim{(U \cp{k} V)}$.
Furthermore, letting $n \eqdef \lydim{(U \cp{k} V)}$, since $n + 1 > k$,
\begin{equation*}
    \size{ \bigcup_{i > n+1} \grade{i}{ ( \maxel{U} ) } } + \size{ \bigcup_{i > n+1} \grade{i}{ ( \maxel{V} ) } } = \size{ \bigcup_{i > n+1} \grade{i}{ ( \maxel{ (U \cp{k} V) } ) } } \leq 1,
\end{equation*}
which implies that $\size{ \bigcup_{i > n+1} \grade{i}{ ( \maxel{U} ) } } \leq 1$ and $\size{ \bigcup_{i > n+1} \grade{i}{ ( \maxel{V} ) } } \leq 1$.
It follows that $\lydim{U} \leq \lydim{(U \cp{k} V)}$ and $\lydim{V} \leq \lydim{(U \cp{k} V)}$.
\end{proof}

\begin{lem} \label{lem:lydim_layering_properties}
Let $U$ be a molecule, $k \eqdef \lydim{U}$.
Suppose $k \geq 0$, and let $(\order{i}{U})_{i=1}^m$ be a $k$\nbd layering of $U$.
Then
\begin{enumerate}
    \item $m > 1$,
    \item for each $i \in \set{1, \ldots, m}$, $\lydim{\order{i}{U}} < k$,
    \item at most one of the $\order{i}{U}$ contains an element of dimension $> k + 1$.
\end{enumerate}
\end{lem}
\begin{proof}
By definition of $\lydim{U}$, if $k \geq 0$ and a $k$\nbd layering exists, then $m > 1$, for otherwise $k - 1 \leq \lydim{U}$, a contradiction.
Moreover, $U$ contains at most one element of dimension $> k + 1$, which can be contained at most in one of the $\order{i}{U}$.
Finally, by Lemma \ref{lem:layering_basic_properties}, we have $\size{ \bigcup_{j > k} \grade{j}{ ( \maxel{\order{i}{U}} ) } } = 1$, so $\lydim{\order{i}{U}} \leq k - 1 < k$.
\end{proof}

\begin{thm} \label{thm:molecules_admit_layerings}
Let $U$ be a molecule, $\lydim{U} \leq k < \dim{U}$.
Then $U$ admits a $k$\nbd layering.
\end{thm}
\begin{proof}
Let $k \eqdef \lydim{U}$.
If $k = -1$, then $U$ is an atom and admits the trivial layering $U = \order{1}{U}$.
If $k \geq 0$, by Lemma \ref{lem:layering_dimension_atom} $U$ is not an atom, so we can assume that $U$ was produced by (\textit{Paste}).
Then $U$ is equal to $V \cp{\ell} W$ for some molecules $V$, $W$ and $\ell < \min \set{ \dim{V}, \dim{W} }$.
By the inductive hypothesis, we have layerings
\begin{equation*}
    \order{1}{V} \cp{k_V} \ldots \cp{k_V} \order{m_V}{V}, \quad \quad \order{1}{W} \cp{k_W} \ldots \cp{k_W} \order{m_W}{W}
\end{equation*}
of $V$ and $W$, respectively, for $k_V \eqdef \lydim{V}$ and $k_W \eqdef \lydim{W}$.
Furthermore, by Lemma \ref{lem:layering_dimensions_pasting}, we know that $k \geq \max \set{ k_V, k_W, \ell }$.
Let
\begin{align*}
    n_V & \eqdef 
    \begin{cases}
        m_V & \text{if $k_V = k$,} \\
        1 & \text{if $k_V < k$ and $\dim{V} > k$,} \\
        0 & \text{if $k_V < \dim{V} < k$},
    \end{cases}
    \\
    n_W & \eqdef 
    \begin{cases}
        m_W & \text{if $k_W = k$,} \\
        1 & \text{if $k_W < k$ and $\dim{W} > k$,} \\
        0 & \text{if $k_W < \dim{W} < k$}.
    \end{cases}
\end{align*}
Notice that it can never be the case that $n_V = n_W = 0$.
We claim that we can decompose $V$ as
\begin{equation} \label{eq:padded_decomposition_1}
    \order{1}{\tilde{V}} \cp{k} \ldots \cp{k} \order{n_V}{\tilde{V}} \cp{k} 
    \underbrace{\bound{k}{+}V \cp{k} \ldots \cp{k} \bound{k}{+}V}_{\text{$n_W$ times}},
\end{equation}
where each $\order{i}{\tilde{V}}$ is a molecule containing exactly one maximal element of dimension $> k$.
If $k_V = k$, we let $\order{i}{\tilde{V}} \eqdef \order{i}{V}$ for all $i \in \set{1, \ldots, m_V}$.
If $k_V < k$, then $V$ contains at most one maximal element of dimension $> k_V + 1$, hence at most one maximal element of dimension $> k$.
If $\dim{V} > k$, it contains exactly one, and we let $\order{1}{\tilde{V}} \eqdef V$.
If $\dim{V} < k$, then $V = \bound{k}{+}V$.
By Proposition \ref{prop:unitality_of_pasting}, pasting copies of $\bound{k}{+}V$ does not change the result up to unique isomorphism.
Similarly, we can decompose $W$ as
\begin{equation} \label{eq:padded_decomposition_2}
    \underbrace{\bound{k}{-}W \cp{k} \ldots \cp{k} \bound{k}{-}W}_{\text{$n_V$ times}} \cp{k}
    \order{1}{\tilde{W}} \cp{k} \ldots \cp{k} \order{n_W}{\tilde{W}}
\end{equation}
where each $\order{i}{\tilde{W}}$ contains exactly one maximal element of dimension $> k$.

If $\ell = k$, since $\ell < \min \set{ \dim{V}, \dim{W} }$, we have $0 < \min \set{ n_V, n_W }$.
Then
\begin{equation*}
    \order{1}{\tilde{V}} \cp{k} \ldots \cp{k} \order{n_V}{\tilde{V}} \cp{k} 
    \order{1}{\tilde{W}} \cp{k} \ldots \cp{k} \order{n_W}{\tilde{W}}
\end{equation*}
is a $k$\nbd layering of $U$.
If $\ell < k$, let
\begin{equation*}
    \order{i}{U} \eqdef 
    \begin{cases}
        \order{i}{\tilde{V}} \cp{\ell} \bound{k}{-}W & \text{if $i \leq n_V$,} \\
        \bound{k}{+}V \cp{\ell} \order{i-n_V}{\tilde{W}} & \text{if $n_V < i \leq n_V+n_W$}.
    \end{cases}
\end{equation*}
Since $\dim{\bound{k}{-}V} = \dim{\bound{k}{+}W} = k$, each $\order{i}{U}$ still contains exactly one maximal element of dimension $> k$.
Plugging (\ref{eq:padded_decomposition_1}) and (\ref{eq:padded_decomposition_2}) in $V \cp{\ell} W$ and using Proposition \ref{prop:interchange_of_pasting} repeatedly, we deduce that $V \cp{\ell} W$ is isomorphic to
\begin{equation*}
    \order{1}{U} \cp{k} \ldots \cp{k} \order{n_V + n_W}{U},
\end{equation*}
which has the desired properties.
Necessarily, $n_V + n_W = m$.

For $\lydim{U} < k < \dim{U}$, the statement follows from Lemma \ref{lem:higher_layerings}.
\end{proof}

\begin{comm} \label{comm:induction_on_lydim}
Theorem \ref{thm:molecules_admit_layerings} in conjunction with Lemma \ref{lem:lydim_layering_properties} and Lemma \ref{lem:layering_dimension_atom} allows us to prove properties of molecules \emph{by induction on their layering dimension}.
That is, to prove that a property holds of all molecules $U$, it suffices to
\begin{itemize}
    \item prove that it holds when $\lydim{U} = -1$, that is, when $U$ is an atom,
    \item prove that it holds when $k \eqdef \lydim{U} \geq 0$, assuming that it holds of all the $(\order{i}{U})_{i=1}^m$ in a $k$\nbd layering of $U$.
\end{itemize}
\end{comm}

\begin{dfn}[Arrow] \index{arrow} \index{$\thearrow{}$}
The \emph{arrow} is the 1\nbd dimensional atom $\thearrow{} \eqdef (1 \celto 1)$.
\end{dfn}

\begin{lem} \label{lem:only_1_molecules} \index{$\thearrow{m}$}
Let $U$ be a 1\nbd dimensional molecule, $m \eqdef \size{\grade{1}{U}}$.
Then $U$ is isomorphic to $\thearrow{m} \eqdef \underbrace{\thearrow{} \cp{0} \ldots \cp{0} \thearrow{}}_{\text{$m$ times}}$.
\end{lem}
\begin{proof}
By Lemma \ref{lem:layering_dimension_smaller_than_dimension}, either $\lydim{U} = -1$ or $\lydim{U} = 0$.
In the first case, $U$ is an atom by Lemma \ref{lem:layering_dimension_atom}.
Because by Lemma \ref{lem:only_0_molecule} the point is the only 0\nbd dimensional molecule up to isomorphism, the arrow is the only 1\nbd dimensional atom, so $U$ is isomorphic to $\thearrow{}$.
In the second case, $U$ admits a 0\nbd layering $(\order{i}{U})_{i=1}^m$ by Theorem \ref{thm:molecules_admit_layerings}, and by Lemma \ref{lem:lydim_layering_properties}, for each $i \in \set{1, \ldots, m}$, necessarily $\lydim{\order{i}{U}} = -1$.
By the first part, $\order{i}{U}$ is isomorphic to $\thearrow{}$.
\end{proof}

\begin{lem} \label{lem:only_2_atoms} \index{$\disk{n}{m}$}
	Let $U$ be a $2$\nbd dimensional atom, $n \eqdef \size{\faces{}{-}U}$, $m \eqdef \size{\faces{}{+}U}$.
	Then $U$ is isomorphic to $\disk{n}{m} \eqdef \thearrow{n} \celto \thearrow{m}$.
\end{lem}
\begin{proof}
Immediate from Lemma \ref{lem:only_1_molecules}.
\end{proof}

\begin{comm}
	Of course, a direct inductive argument on the structure of molecules would also have sufficed to prove Lemma 
	\ref{lem:only_1_molecules}.
	Nevertheless, we use this as a simple illustration of how an argument based on layerings can make such proofs more concise.
\end{comm}


\section{Flow graphs and orderings} \label{sec:flow_graphs}

\begin{guide}
	In this section, we define two families of graphs associated to an oriented graded poset: its \cemph{$k$\nbd flow graphs} and its \cemph{maximal $k$\nbd flow graphs}.
	We then define a \emph{$k$\nbd ordering} of a molecule to be a topological sort of its maximal $k$\nbd flow graph, which exists if and only if the latter is acyclic.
	
	We prove that each $k$\nbd layering determines a unique $k$\nbd ordering (Proposition 
	\ref{prop:layerings_induce_orderings}) and give a criterion for when a $k$\nbd ordering corresponds to a $k$\nbd layering (Proposition 
	\ref{prop:layering_from_ordering}).

	Finally, we use the theory developed in these sections to prove that morphisms of molecules of the same dimension preserve interiors (Proposition \ref{prop:morphisms_of_molecules_preserve_interior}), an important step towards the proof of Theorem \ref{thm:morphisms_of_atoms_are_injective} in the next chapter.
\end{guide}

\begin{dfn}[Flow graph] \index{oriented graded poset!flow graph} \index{molecule!flow graph} \index{$\flow{k}{U}$} \index{flow graph!of an oriented graded poset}
Let $P$ be an oriented graded poset, $k \geq -1$.
The \emph{$k$\nbd flow graph of $P$} is the directed graph $\flow{k}{P}$ whose
\begin{itemize}
    \item set of vertices is $\bigcup_{i > k} \grade{i}{P}$, and
    \item set of edges is
	    \[
		    \set{(x, y) \mid \faces{k}{+}x \cap \faces{k}{-}y \neq \varnothing},
		\]
	with $s\colon (x, y) \mapsto x$ and $t\colon (x, y) \mapsto y$.
\end{itemize}
\end{dfn}

\begin{dfn}[Maximal flow graph] \index{oriented graded poset!flow graph!maximal} \index{molecule!flow graph!maximal} \index{$\maxflow{k}{U}$} \index{flow graph!maximal}
Let $P$ be a finite-dimensional oriented graded poset, $k \geq -1$.
The \emph{maximal $k$\nbd flow graph of $P$} is the induced subgraph $\maxflow{k}{P}$ of $\flow{k}{P}$ on the vertex set
\begin{equation*}
    \bigcup_{i > k} \grade{i}{(\maxel{P})} \subseteq \bigcup_{i > k} \grade{i}{P}.
\end{equation*}
\end{dfn}

\begin{rmk} \label{rmk:codim_1_maximal_flow_graph}
For $k \eqdef \dim{P} - 1$, $\flow{k}{P}$ and $\maxflow{k}{P}$ coincide.
\end{rmk}

\begin{dfn}[Topological sort] \index{directed graph!topological sort}
Let $\mathscr{G}$ be a directed acyclic graph with finite set of vertices, $m \eqdef \size{V_\mathscr{G}}$.
A \emph{topological sort of $\mathscr{G}$} is a linear ordering $(\order{i}{x})_{i=1}^m$ of $V_\mathscr{G}$ such that, for all edges $e \in E_\mathscr{G}$, if $s(e) = \order{i}{x}$ and $t(e) = \order{j}{x}$, then $i < j$.
\end{dfn}

\begin{dfn}[Ordering of a molecule] \index{molecule!ordering}
Let $U$ be a molecule, $k \geq -1$, and suppose $\maxflow{k}{U}$ is acyclic.
A \emph{$k$\nbd ordering of $U$} is a topological sort of $\maxflow{k}{U}$.
\end{dfn}

\begin{exm}[Flow graph and maximal flow graph of a molecule]
	Let $U$ be the oriented face poset of (\ref{eq:pasting_to_string}) from Example \ref{exm:string_diagrams}.
	The 0\nbd flow graph $\flow{0}{U}$ is
\[\begin{tikzcd}
	& {{\scriptstyle (1, 5)}\;\bullet} \\
	{{\scriptstyle (1, 3)}\;\bullet} & {{\scriptstyle (2, 1)}\;\bullet} && {{\scriptstyle (1, 6)}\;\bullet} \\
	& {{\scriptstyle (1, 4)}\;\bullet} && {{\scriptstyle (2, 2)}\;\bullet} \\
	& {{\scriptstyle (2, 0)}\;\bullet} & {{\scriptstyle (1, 1)}\;\bullet} & {{\scriptstyle (1, 2)}\;\bullet} \\
	& {{\scriptstyle (1, 0)}\;\bullet}
	\arrow[from=4-2, to=4-3]
	\arrow[from=5-2, to=4-3]
	\arrow[from=3-2, to=4-3]
	\arrow[from=2-1, to=3-2]
	\arrow[from=2-1, to=1-2]
	\arrow[from=1-2, to=3-4]
	\arrow[from=4-3, to=3-4]
	\arrow[from=2-1, to=2-2]
	\arrow[from=1-2, to=2-4]
	\arrow[from=1-2, to=4-4]
	\arrow[from=4-3, to=4-4]
	\arrow[from=4-3, to=2-4]
	\arrow[from=2-2, to=2-4]
	\arrow[from=2-2, to=3-4]
	\arrow[from=2-2, to=4-4]
\end{tikzcd}\]
	and the maximal 0\nbd flow graph $\maxflow{0}{U}$ is its induced subgraph
\[\begin{tikzcd}
	{{\scriptstyle (2, 1)}\;\bullet} && {{\scriptstyle (2, 2)}\;\bullet} \\
	{{\scriptstyle (2, 0)}\;\bullet}
	\arrow[from=1-1, to=1-3]
\end{tikzcd}\]
	while the 1\nbd flow graph $\flow{1}{U}$ is
\[\begin{tikzcd}
	{{\scriptstyle (2, 1)}\;\bullet} && {{\scriptstyle (2, 2)}\;\bullet} \\
	{{\scriptstyle (2, 0)}\;\bullet}
	\arrow[from=2-1, to=1-1]
\end{tikzcd}\]
	and it is equal to $\maxflow{1}{U}$, since every 2\nbd dimensional element of $U$ is maximal.
\end{exm}

\begin{lem} \label{lem:maxflow_acyclic_pasting}
Let $U$, $V$ be molecules and $k < \min \set {\dim{U}, \dim{V}}$ such that $U \cp{k} V$ is defined.
If $\maxflow{k}{U}$ and $\maxflow{k}{V}$ are acyclic, then $\maxflow{k}{(U \cp{k} V)}$ is acyclic.
\end{lem}
\begin{proof}
Suppose that $\maxflow{k}{U}$ and $\maxflow{k}{V}$ are acyclic.
We may identify $U$ and $V$ with their images in $U \cp{k} V$.
By Lemma \ref{lem:faces_of_union}, since $\dim{(U \cap V)} = k$,
\begin{equation*}
    \bigcup_{i > k} \grade{i}{ ( \maxel{ (U \cp{k} V) } ) }  =
    \bigcup_{i > k} \grade{i}{ ( \maxel{U} ) } + \bigcup_{i > k} \grade{i}{ ( \maxel{V} ) },
\end{equation*}
so $\maxflow{k}{U}$ and $\maxflow{k}{V}$ are isomorphic to the induced subgraphs of $\maxflow{k}{(U \cp{k} V)}$ on the vertices in $U$ and $V$, respectively.
It follows that a cycle in $\maxflow{k}{(U \cp{k} V)}$ cannot remain in $U$ or $V$, but has to visit vertices in both.
In particular, such a cycle has to go through an edge from $x \in V$ to $y \in U$, induced by the existence of $z \in \faces{k}{+}x \cap \faces{k}{-}y$.
But then $z \notin \bound{k}{-}V$ and $z \notin \bound{k}{+}U$, yet $z \in U \cap V$, a contradiction.
\end{proof}

\begin{prop} \label{prop:if_layering_then_ordering}
Let $U$ be a molecule, $k \geq -1$.
If $U$ admits a $k$\nbd layering, then $\maxflow{k}{U}$ is acyclic, and $U$ admits a $k$\nbd ordering.
\end{prop}
\begin{proof}
Let $(\order{i}{U})_{i=1}^m$ be a $k$\nbd layering of $U$.
For each $i \in \set{1, \ldots, m}$, the graph $\maxflow{k}{\order{i}{U}}$ is trivially acyclic by Lemma \ref{lem:layering_basic_properties}.
We conclude by applying Lemma \ref{lem:maxflow_acyclic_pasting} repeatedly.
\end{proof}

\begin{cor} \label{cor:flow_acyclic_in_codimension_1}
Let $U$ be a molecule, $n \eqdef \dim{U}$.
Then $\flow{n-1}{U}$ is acyclic.
\end{cor}
\begin{proof}
By Theorem \ref{thm:molecules_admit_layerings}, $U$ always admits an $(n-1)$\nbd layering.
We conclude by Proposition \ref{prop:if_layering_then_ordering} combined with the fact that $\flow{n-1}{U} = \maxflow{n-1}{U}$.
\end{proof}

\begin{dfn}[Sets of layerings and of orderings] \index{$\layerings{k}{U}$} \index{$\orderings{k}{U}$}
Let $U$ be a molecule, $k \geq -1$.
We let
\begin{align*}
    \layerings{k}{U} & \eqdef \set{\text{$k$\nbd layerings $(\order{i}{U})_{i=1}^m$ of $U$ up to layer-wise isomorphism}}, \\
    \orderings{k}{U} & \eqdef \set{\text{$k$\nbd orderings $(\order{i}{x})_{i=1}^m$ of $U$}}.
\end{align*}
\end{dfn}

\begin{prop} \label{prop:layerings_induce_orderings} \index{$\lto{k}{U}$}
Let $U$ be a molecule, $k \geq -1$.
For each $k$\nbd layering $(\order{i}{U})_{i=1}^m$ of $U$ and each $i \in \set{1, \ldots, m}$, let $\order{i}{x}$ be the only element of $\bigcup_{j > k}\grade{j}{(\maxel{U})}$ in the image of $\order{i}{U}$.
Then the assignment
\begin{equation} \label{eq:lay_to_ord}
    \lto{k}{U}\colon (\order{i}{U})_{i=1}^m \mapsto (\order{i}{x})_{i=1}^m
\end{equation}
determines an injective function $\layerings{k}{U} \incl \orderings{k}{U}$.
\end{prop}
\begin{proof}
By Lemma \ref{lem:layering_basic_properties}, the assignment $(\order{i}{U})_{i=1}^m \mapsto (\order{i}{x})_{i=1}^m$ is well-defined.
Let $i, j \in \set{1, \ldots, m}$, and suppose that there is an edge from $\order{i}{x}$ to $\order{j}{x}$ in $\maxflow{k}{U}$, that is, there exists $z \in \faces{k}{+}\order{i}{x} \cap \faces{k}{-}\order{j}{x}$.
By Proposition \ref{prop:if_layering_then_ordering}, $\maxflow{k}{U}$ is acyclic, so necessarily $i \neq j$.
If $j < i$, then $\order{j}{U} \cap \order{i}{U} \subseteq \bound{k}{+}\order{j}{U} \cap \bound{k}{-}\order{i}{U}$ by Lemma \ref{lem:layering_intersections}, contradicting the existence of $z$.
It follows that $i < j$, so $(\order{i}{x})_{i=1}^m$ is a $k$\nbd ordering of $U$.

Let $(\order{i}{V})_{i=1}^m$ be another $k$\nbd layering, and suppose it determines the same $k$\nbd ordering as $(\order{i}{U})_{i=1}^m$.
Then the image of both $\order{1}{U}$ and $\order{1}{V}$ in $U$ is
\begin{equation*}
\clos\set{\order{1}{x}} \cup \bound{}{-}U,
\end{equation*}
so $\order{1}{U}$ is isomorphic to $\order{1}{V}$.
If $m = 1$ we are done.
Otherwise, $(\order{i}{U})_{i=2}^m$ and $(\order{i}{V})_{i=2}^m$ are $k$\nbd layerings inducing the same $k$\nbd ordering on their image.
By recursion, we conclude that they are layer-wise isomorphic.
\end{proof}

\begin{lem} \label{lem:lto_bijection_above_layering}
Let $U$ be a molecule, $\ell \geq -1$.
If $U$ has an $\ell$\nbd layering, then for all $k > \ell$ the function $\lto{k}{U}\colon \layerings{k}{U} \incl \orderings{k}{U}$ is a bijection.
\end{lem}
\begin{proof}
Let $(\order{i}{U})_{i=1}^m$ be an $\ell$\nbd layering of $U$, and let $(\order{i}{x})_{i=1}^m$ be its image through $\lto{\ell}{U}$.
For $k > \ell$, let $(\order{i}{y})_{i=1}^p$ be a $k$\nbd ordering of $U$.
Then there exists a unique injection $\fun{j}\colon \set{1, \ldots, p} \incl \set{1, \ldots, m}$ such that $\order{i}{y} = \order{\fun{j}(i)}{x}$ for all $i \in \set{1, \ldots, p}$.
Let
\begin{align*}
	\order{i}{V} & \eqdef \bound{k}{\alpha(i,1)}\order{1}{U} \cp{\ell} \ldots \cp{\ell} \order{\fun{j}(i)}{U} \cp{\ell} \ldots \cp{\ell} \bound{k}{\alpha(i,m)}\order{m}{U}, \\
	\alpha(i,j) & \eqdef \begin{cases} + & \text{if $j = \fun{j}(i')$ for some $i' < i$,} \\ 
	- & \text{otherwise}.
	\end{cases}
\end{align*}
Applying Proposition \ref{prop:unitality_of_pasting} and Proposition \ref{prop:interchange_of_pasting} repeatedly, we find that $(\order{i}{V})_{i=1}^p$ is a $k$\nbd layering of $U$ and $(\order{i}{y})_{i=1}^p$ is its image through $\lto{k}{U}$.
This proves that $\lto{k}{U}$ is surjective, and we conclude by Proposition \ref{prop:layerings_induce_orderings}.
\end{proof}

\begin{dfn}[Rewrite steps] \index{molecule!rewrite steps}
Let $U$ be a molecule, $k \geq -1$, and let $(\order{i}{U})_{i=1}^m$ be a $k$\nbd layering of $U$.
The sequence $(\step{i}{U})_{i=0}^m$ of \emph{rewrite steps} associated with $(\order{i}{U})_{i=1}^m$ is defined recursively by
\begin{itemize}
    \item $\step{0}{U} \eqdef \bound{k}{-}U$,
    \item $\step{i}{U} \eqdef \bound{k}{+}\order{i}{U}$ for $i \in \set{1, \ldots, m}$.
\end{itemize}
\end{dfn}

\begin{lem} \label{lem:boundary_move}
Let $U$ be a molecule, $k \in \mathbb{N}$, and suppose
\begin{equation*}
    \bigcup_{i > k} \grade{i}{ ( \maxel{U} ) } = \set{x}.
\end{equation*}
Then, for all $\alpha \in \set { +, - }$,
\begin{enumerate}
    \item $\bound{k}{\alpha}x \submol \bound{k}{\alpha}U$, 
    \item $\bound{k}{\alpha}U$ is isomorphic to $\subs{\bound{k}{-\alpha}U}{\bound{k}{\alpha}x}{\bound{k}{-\alpha}x}$.
\end{enumerate}
\end{lem}
\begin{proof}
We proceed by induction on $\lydim{U}$.
If $\lydim{U} = -1$, then $U$ is an atom and equal to $\clset{{x}}$.
It follows that $\bound{k}{\alpha}x = \bound{k}{\alpha}U$, which is trivially a submolecule, and is isomorphic to $\subs{\bound{k}{-\alpha}U}{\bound{k}{\alpha}x}{\bound{k}{-\alpha}{x}}$ by Lemma \ref{lem:full_substitution}.

Suppose $\ell \eqdef \lydim{U} \geq 0$, and let $(\order{i}{U})_{i = 1}^m$ be an $\ell$\nbd layering of $U$.
Then $\ell \leq k-1 < k$ because $\size{\bigcup_{i > k} \grade{i}{ ( \maxel{U} ) }} = 1$.
By Lemma \ref{lem:pasting_higher_boundary}, $\bound{k}{\alpha}U$ is isomorphic to
\begin{equation*}
    \bound{k}{\alpha}\order{1}{U} \cp{\ell} \ldots \cp{\ell} \bound{k}{\alpha}\order{m}{U}.
\end{equation*}
Now $x$ is contained in a single $\order{i}{U}$.
By the inductive hypothesis, $\bound{k}{\alpha}x \submol \bound{k}{\alpha}\order{i}{U}$, and the latter is isomorphic to $\subs{\bound{k}{-\alpha}\order{i}{U}}{\bound{k}{\alpha}x}{\bound{k}{-\alpha}{x}}$.
We conclude by Lemma \ref{lem:pasting_after_substitution}.
\end{proof}

\begin{cor} \label{cor:rewrite_steps_substitution}
Let $U$ be a molecule, $k \geq -1$, $(\order{i}{U})_{i=1}^m$ a $k$\nbd layering of $U$, and $(\order{i}{x})_{i=1}^m$ the associated $k$\nbd ordering.
For all $i \in \set{1, \ldots, m}$, the $i$-th rewrite step $\step{i}{U}$ is isomorphic to $\subs{\step{i-1}{U}}{\bound{k}{+}\order{i}{x}}{\bound{k}{-}\order{i}{x}}$.
\end{cor}
\begin{proof}
By Proposition \ref{lem:layering_basic_properties}, $\order{i}{x}$ is the only element of dimension $> k$ in $\order{i}{U}$.
The result then follows from repeated application of Lemma \ref{lem:boundary_move}.
\end{proof}

\begin{exm}[Rewrite steps and string diagrams]
	When we draw a 3\nbd dimensional pasting diagram as a sequence of 2\nbd dimensional diagrams, we are informally using the rewrite steps associated to a 2\nbd layering of the diagram in order to understand its 3\nbd dimensional shape.
	Indeed, the ability to take layerings, and form the corresponding sequences of rewrite steps, is crucial to our ability to visualise higher-dimensional shapes, where our spatial intuition falters.
	This is especially powerful when coupled with the use of string diagrams as in Example 
	\ref{exm:string_diagrams}, which allows us to ``isolate'' the information of \emph{what codimension-1 faces are being rewritten at each step}.

	For example, in the higher algebraic theory whose models include \emph{monoidal bicategories}, there is a 4\nbd dimensional cell called the \emph{pentagonator}, a weak version of Mac Lane's pentagon equation.
	The oriented face poset $U$ of its shape has the oriented Hasse diagram
	\[
		\input{img/penta_hasse.tex}
	\]
	which is of course quite unintelligible.
	Let us try to gain a better understanding by using the aforementioned tools.

	The graph of $U$, as a string diagram, is
	\[
		\input{img/penta_strdiag4.tex}\;,
	\]
	so the greatest element of $U$ rewrites three 3\nbd dimensional faces into two 3\nbd dimensional faces.
	Since $U$ is an atom, to understand its shape it suffices to understand the round 3\nbd dimensional molecules $\bound{}{-}U$ and $\bound{}{+}U$.
	The graph of $\bound{}{-}U$ is
	\begin{equation} \label{eq:penta_input_strdiag}
		\input{img/penta_input_strdiag3.tex}
	\end{equation}
	which, seen as a rewrite sequence, contains three 3\nbd dimensional rewrites, each replacing two 2\nbd dimensional faces with two 2\nbd dimensional faces.
	This molecule has a unique 2\nbd layering, corresponding to the 2\nbd ordering $(3, 0), (3, 1), (3, 2)$.
	The rewrite steps associated with this 2\nbd layering have graphs
	\[
		\input{img/penta_input_step0.tex} \quad
		\input{img/penta_input_step1.tex}
	\]
	\[
		\input{img/penta_input_step2.tex} \quad
		\input{img/penta_input_step3.tex}
	\]
	which now have the ``correct'' dimension and can be directly understood as 2\nbd dimensional diagrams.
	One can follow the flow of 2\nbd dimensional elements as depicted in (\ref{eq:penta_input_strdiag}) across these rewrite steps: between the zeroth and the first step, $(3, 0)$ replaces $(2, 0)$ and $(2, 1)$ with $(2, 4)$ and $(2, 3)$; between the first and second step, $(3, 1)$ replaces $(2, 2)$ and $(2, 3)$ with $(2, 5)$ and $(2, 6)$; and so on.

	Meanwhile, the graph of $\bound{}{+}U$ is
	\begin{equation} \label{eq:penta_output_strdiag}
		\input{img/penta_output_strdiag3.tex}
	\end{equation}
	which contains only two 3\nbd dimensional rewrites, also replacing two 2\nbd dimensional faces with two 2\nbd dimensional faces.
	Notice that if we frame the graph of $\bound{}{-}U$ and the graph of $\bound{}{+}U$ with the same order of input and output wires, only one of the two admits a planar representation.
	
	Now, $\bound{}{+}U$ also has a unique 2\nbd layering, corresponding to the 2\nbd ordering $(3, 4), (3, 3)$, whose sequence of rewrite steps is
	\[
		\input{img/penta_input_step0.tex} \quad
		\input{img/penta_output_step1.tex}
	\]
	\[
		\input{img/penta_input_step3.tex}\;;
	\]
	notice that the zeroth and last rewrite step are equal to those of $\bound{}{-}U$, as expected since the two molecules have isomorphic boundaries.
	Once again, one can follow the flow of 2\nbd dimensional elements as depicted in (\ref{eq:penta_output_strdiag}) across these rewrite steps.

	Overall, we can now interpret $(4, 0)$ as a ``coherence'' cell between two different sequences of 3\nbd dimensional rewrites on 2\nbd dimensional diagrams (which semantically may be interpreted as ``rebracketing'' steps for a binary operation that is \emph{associative up to higher homotopy}).
	We briefly note that $U$ is also a dual of the \emph{oriented 4\nbd simplex}, see Section 
	\ref{sec:simplices}.
\end{exm}

\begin{prop} \label{prop:layering_from_ordering}
Let $U$ be a molecule, $k \geq -1$, and let $(\order{i}{x})_{i=1}^m$ be a $k$\nbd ordering of $U$.
Let
\begin{align*}
    \order{0}{U} & \eqdef \bound{k}{-}U, \\
    \order{i}{U} & \eqdef \bound{k}{+}\order{i-1}{U} \cup \clset{{\order{i}{x}}} \quad \text{for $i \in \set{1, \ldots, m}$}.
\end{align*}
The following are equivalent:
\begin{enumerate}[label=(\alph*)]
    \item $(\order{i}{U})_{i=1}^m$ is a $k$\nbd layering of $U$;
    \item for all $i \in \set{1, \ldots, m}$, $\bound{k}{-}\order{i}{x} \submol \bound{k}{-}\order{i}{U}$. \label{cond:inputs_are_submolecules}
\end{enumerate}
Moreover, for all $i \in \set{1, \ldots, m-1}$, if $\bound{k}{-}\order{i}{x} \submol \bound{k}{-}\order{i}{U}$, then $\order{i}{U}$ and $\bound{k}{+}\order{i}{U} = \bound{k}{-}\order{i+1}{U}$ are molecules.
\end{prop}
\begin{proof}
Suppose $(\order{i}{U})_{i=1}^m$ is a $k$\nbd layering.
Then, for all $i \in \set{1, \ldots, m}$, $\order{i}{U}$ is a molecule, and by Proposition \ref{lem:layering_basic_properties} $\order{i}{x}$ is the only element of dimension $> k$ in $\order{i}{U}$.
By Lemma \ref{lem:boundary_move}, $\bound{k}{-}\order{i}{x} \submol \bound{k}{-}\order{i}{U}$.

Conversely, it follows from Lemma \ref{lem:submolecule_rewrite} that for all $i \in \set{1, \ldots, m}$, if $\bound{k}{-}\order{i}{U}$ is a molecule and $\bound{k}{-}\order{i}{x} \submol \bound{k}{-}\order{i}{U}$, then $\order{i}{U}$ is a molecule, hence $\bound{k}{+}\order{i}{U}$ is a molecule.
Moreover, since $(\order{i}{x})_{i=1}^m$ is a $k$\nbd ordering, it is straightforward to prove that $\order{i}{U} \cap \order{i+1}{U} = \bound{k}{+}\order{i}{U} = \bound{k}{-}\order{i+1}{U}$ for all $i \in \set{1, \ldots, m-1}$.
Since $\bound{}{-}\order{1}{U} = \bound{}{-}U$ is a molecule, it follows by induction, assuming condition \ref{cond:inputs_are_submolecules}, that $\order{i}{U}$ is a molecule for all $i \in \set{1, \ldots, m}$.
This proves that $(\order{i}{U})_{i=1}^m$ is a $k$\nbd layering of $U$.
\end{proof}

\begin{lem} \label{lem:intersections_of_interiors_of_boundaries}
	Let $U$ be a molecule, $n \eqdef \dim{U}$, and $x, y \in \grade{n}{U}$.
	For all $\alpha \in \set{+, -}$, if $x \neq y$, then $\inter{\bound{}{\alpha}x} \cap \inter{\bound{}{\alpha}y} = \varnothing$.
\end{lem}
\begin{proof}
	By Theorem \ref{thm:molecules_admit_layerings}, $U$ has an $(n-1)$\nbd layering $(\order{i}{U})_{i=1}^m$; we will identify each layer with its image in $U$.
	Let $(\order{i}{x})_{i=1}^m$ be the associated $(n-1)$\nbd ordering.
	Necessarily, 
	\[
		\set{\order{i}{x} \mid i \in \set{1, \ldots, m}} = \grade{n}{(\maxel{U})} = \grade{n}{U}.
	\]
	If $x \neq y$, there exists a unique pair $i \neq j \in \set{1, \ldots, m}$ such that $x = \order{i}{x}$ and $y = \order{j}{x}$.
	Suppose without loss of generality that $i < j$ and that $\alpha = -$.
	Then
	\[
		\clset{x} \cap \clset{y} \subseteq \clset{x} \cap \order{i}{U} \cap \order{j}{U} \subseteq \clset{x} \cap \bound{}{+}\order{i}{U} \subseteq \bound{}{+}x
	\]
	using Lemma \ref{lem:layering_intersections} and Lemma \ref{lem:intersection_with_boundary}.
	But
	\[
		\inter{\bound{}{-}x} = \bound{}{-}x \setminus (\bound{}{-}x \cap \bound{}{+}x) = \bound{}{-}x \setminus \bound{}{+}x
	\]
	by roundness of $\clset{x}$.
	This proves that $\inter{\bound{}{-}x} \cap \clset{y} = \varnothing$.
	All other cases are analogous.
\end{proof}

\begin{lem} \label{lem:characterisation_of_interior}
	Let $U$ be a molecule, $n \eqdef \dim{U}$, and let $x \in U$.
	Then $x \in \inter{U}$ if and only if
	\begin{itemize}
		\item $\dim{x} = n$, or
		\item $\dim{x} < n$ and $x \in \inter{\bound{}{+}y_+} \cap \inter{\bound{}{-}y_-}$ for a pair $y_+, y_- \in \grade{n}{U}$.
	\end{itemize}
	In the latter case, the pair $y_+, y_-$ is unique with this property.
\end{lem}
\begin{proof}
	By Theorem \ref{thm:molecules_admit_layerings}, $U$ has an $(n-1)$\nbd layering $(\order{i}{U})_{i=1}^m$; we will identify each layer with its image in $U$.
	Let $(\order{i}{x})_{i=1}^m$ be the associated $(n-1)$\nbd ordering.
	
	Suppose that $x \in \inter{U}$.
	Then there exists at least one $i \in \set{1, \ldots, m}$ such that $x \leq \order{i}{x}$.
	If $x = \order{i}{x}$, we are done.
	
	Otherwise, we have $x \in \bound{}{\alpha}\order{i}{x}$ for some $\alpha \in \set{+, -}$.
	Suppose without loss of generality that $\alpha = +$; then $x \in \step{i}{U}$.
	We claim that there exists a smallest $j > i$ such that $x \in \inter{\bound{}{-}\order{j}{x}}$.
	If that was not the case, using Corollary \ref{cor:rewrite_steps_substitution} in conjunction with Lemma \ref{lem:substitution_direct_pushout}, we could derive that $x \in \step{j}{U}$ for all $j > i$; in particular, $x \in \step{m}{U} = \bound{}{+}U$, a contradiction.
	By a dual argument, there exists a greatest $k \leq i$ such that $x \in \inter{\bound{}{+}\order{k}{x}}$.
	We let $y_- \eqdef \order{j}{x}$ and $y_+ \eqdef \order{k}{x}$, which proves the existence of a pair with the required property.
	Uniqueness is then a consequence of Lemma \ref{lem:intersections_of_interiors_of_boundaries}.

	Conversely, suppose $x \in \bound{}{} U$, so $\dim{x} < n$.
	Then there exists $\alpha \in \set{+, -}$ such that $x \in \bound{}{\alpha}U$; without loss of generality, suppose $\alpha = -$.
	By Lemma \ref{lem:downset_is_submolecule} and Lemma \ref{lem:intersection_with_boundary}, for all $i \in \set{1, \ldots, m}$,
	\[
		\bound{}{+}\order{i}{x} \subseteq \bound{}{-}U \cap \bound{}{+}\order{i}{x} \subseteq \bound{}{-}\order{i}{x} \cap \bound{}{+}\order{i}{x} = \bound{}{}(\bound{}{+}\order{i}{x}),
	\]
	where in the last step we used roundness of $\clset{\order{i}{x}}$.
	Thus it is never the case that $x \in \inter{\bound{}{+}\order{i}{x}}$.
\end{proof}

\begin{rmk}
	Notice that, when $\codim{x}{U} = 1$, Lemma \ref{lem:characterisation_of_interior} is a special case of Corollary \ref{cor:codimension_1_elements}.
\end{rmk}

\begin{prop} \label{prop:morphisms_of_molecules_preserve_interior}
	Let $f\colon U \to V$ be a morphism of molecules and suppose $\dim{U} = \dim{V}$.
	Then $f(\inter{U}) \subseteq \inter{V}$.
\end{prop}
\begin{proof}
	We proceed by induction on $n \eqdef \dim{U} = \dim{V}$.
	If $n = 0$, then $U$ and $V$ are both the point and the statement is obvious.

	Suppose $n > 0$ and let $x \in \inter{U}$.
	By Lemma \ref{lem:characterisation_of_interior}, either $\dim{x} = n$ or there exists a pair $y_+, y_- \in \grade{n}{U}$ with $x \in \inter{\bound{}{+}y_+} \cap \inter{\bound{}{-}y_-}$.
	In the first case, $\dim{f(x)} = n$, so $f(x) \in \inter{V}$.
	In the second case, for each $\alpha \in \set{+, -}$,
	\[
		f(\bound{}{\alpha}y_\alpha) = f(\clos{\faces{}{\alpha}y_\alpha}) = \clos{\faces{}{\alpha}f(y_\alpha)} = \bound{}{\alpha}f(y_\alpha),
	\]
	and both $\bound{}{\alpha}y_\alpha$ and $\bound{}{\alpha}f(y_\alpha)$ are molecules of dimension $n-1$.
	By the inductive hypothesis applied to $\restr{f}{\bound{}{\alpha}y_\alpha}\colon \bound{}{\alpha}y_\alpha \to \bound{}{\alpha}f(y_\alpha)$, we have
	\[
		f(\inter{\bound{}{\alpha}y_\alpha}) \subseteq \inter{\bound{}{\alpha}f(y_\alpha)}
	\]
	hence $f(x) \subseteq \inter{\bound{}{+}f(y_+)} \cap \inter{\bound{}{-}f(y_-)}$.
	Using Lemma \ref{lem:characterisation_of_interior} once more, we conclude that $f(x) \in \inter{V}$.
\end{proof}


\section{Frame dimension} \label{sec:frame_dimension}

\begin{guide}
	Theorem \ref{thm:molecules_admit_layerings} puts an upper bound on the minimal $k$ such that a $k$\nbd layering of a molecule exists.
	In this section, we define the \emph{frame dimension} of a molecule, and prove that it is a \emph{lower} bound for the same value.

	One could hope, after a little experimenting, that this lower bound is always reached, that is, a layering in the frame dimension always exists.
	We will see in a later chapter that this is true up to dimension 3 (Theorem \ref{thm:dim3_frame_acyclic}), but fails in dimension 4.

	The main result that we prove in this section (Proposition 
	\ref{prop:intersection_of_maximal_elements}) is mainly of technical interest: it shows that there cannot be any overlap of the \emph{interiors} of the input $k$\nbd boundaries of two maximal elements of a molecule, as soon as $k$ is strictly smaller than the frame dimension; and similarly for the interiors of the output $k$\nbd boundaries.
\end{guide}

\begin{dfn}[Frame dimension] \index{molecule!frame dimension} \index{$\frdim{U}$} \index{dimension!frame}
Let $U$ be a molecule.
The \emph{frame dimension of $U$} is the integer
\begin{equation*}
    \frdim{U} \eqdef \dim{\bigcup \set{(\clset{{x}} \cap \clset{{y}}) \mid
    x, y \in \maxel{U}, x \neq y } }.
\end{equation*}
\end{dfn}

\begin{lem} \label{lem:frame_dimension_atom}
Let $U$ be a molecule.
Then $\frdim{U} = -1$ if and only if $U$ is an atom.
\end{lem}
\begin{proof}
Suppose that $\frdim{U} = -1$, and let $x \in \maxel{U}$.
Then, letting $V \eqdef \clos{ (( \maxel{U}) \setminus \set{x}) }$, we have $U = \clos \set{x} \cup V$ and $\clset{{x}} \cap V = \varnothing$.
By Lemma \ref{lem:molecules_are_connected}, $V = \varnothing$, so $x$ is the greatest element of $U$.
Conversely, if $U$ is an atom, there does not exist a pair of distinct elements of $\maxel{U}$, so $\frdim{U} = \dim{\varnothing} = -1$.
\end{proof}

\begin{lem} \label{lem:layering_vs_frame_dimension}
Let $U$ be a molecule.
Then $\frdim{U} \leq \lydim{U}$.
\end{lem}
\begin{proof}
Let $r \eqdef \frdim{U}$.
If $r = -1$, by Lemma \ref{lem:frame_dimension_atom} $U$ is an atom, and by Lemma \ref{lem:layering_dimension_atom} $\lydim{U}$ is also $-1$.
Suppose $r \geq 0$.
Then there exist distinct maximal elements $x, y \in U$ such that $\dim{(\clset{{x}} \cap \clset{{y}})} = r$.
Necessarily $r < \min \set{ \dim{x}, \dim{y} }$, so $x, y \in \bigcup_{i > r} \grade{i}{ (\maxel{U}) }$ and $\size{ \bigcup_{i > r} \grade{i}{ (\maxel{U}) } } \geq 2$.
It follows that $r - 1 < \lydim{U}$, that is, $r \leq \lydim{U}$.
\end{proof}

\begin{exm}[A molecule with unequal frame and layering dimensions] \index[counterex]{A molecule with unequal frame and layering dimensions}
	The 2\nbd dimensional oriented face poset $U$ of (\ref{eq:both_0_and1_layering}) satisfies
	\[
		\frdim{U} = 0, \quad \quad \lydim{U} = 1.
	\]
	Notice that there are no such examples in dimension 1: by the classification result Lemma \ref{lem:only_1_molecules}, a 1\nbd dimensional molecule is either the atom $\thearrow{}$, in which case both its frame and layering dimension are -1 by Lemma \ref{lem:layering_dimension_atom} and 
	Lemma \ref{lem:frame_dimension_atom}, or it is the molecule $\thearrow{m}$ for $m > 1$, in which case both its frame and layering dimension are 0.
\end{exm}

\begin{lem} \label{lem:layerings_only_above_frdim}
Let $U$ be a molecule, $k \geq -1$.
If $U$ admits a $k$\nbd layering, then $k \geq \frdim{U}$.
\end{lem}
\begin{proof}
Let $x, y$ be distinct maximal elements of $U$.
If $\min \set{\dim{x}, \dim{y}} \leq k$, then $\dim{(\clset{{x}} \cap \clset{{y}})} < k$.
Suppose that $k < \min \set{\dim{x}, \dim{y}}$, and let $(\order{i}{U})_{i=1}^m$ be a $k$\nbd layering of $U$.
By Lemma \ref{lem:layering_basic_properties} there exist $i \neq j$ such that $x \in \order{i}{U}$ and $y \in \order{j}{U}$.
By Lemma \ref{lem:layering_intersections}, there exists $\alpha \in \set{+, -}$ such that $\order{i}{U} \cap \order{j}{U} = \bound{k}{\alpha}\order{i}{U} \cap \bound{k}{-\alpha}\order{j}{U}$.
Then $\clset{{x}} \cap \clset{{y}} \subseteq \bound{k}{\alpha}\order{i}{U} \cap \bound{k}{-\alpha}\order{j}{U}$, so $\dim{(\clset{{x}} \cap \clset{{y}})} \leq k$.
\end{proof}

\begin{cor} \label{cor:least_dimension_with_layering}
Let $U$ be a molecule.
Then
\begin{equation*}
    \frdim{U} \leq \min \set{k \geq -1 \mid \text{$U$ admits a $k$\nbd layering}} \leq \lydim{U}.
\end{equation*}
\end{cor}
\begin{proof}
Follows from Lemma \ref{lem:layerings_only_above_frdim} and Theorem \ref{thm:molecules_admit_layerings}.
\end{proof}

\begin{lem} \label{lem:frdim_boundaries}
Let $U$ be a molecule, $x \in \maxel{U}$, $\alpha \in \set{+,-}$.
Suppose that $\dim{x} > \frdim{U}$.
Then, for all $k > \frdim{U}$,
\begin{equation*}
    \bound{k}{\alpha}x \subseteq \bound{k}{\alpha}U.
\end{equation*}
\end{lem}
\begin{proof}
If $k \geq \dim{x}$, then $\bound{k}{\alpha}x = \clset{{x}}$ and $x \in \bigcup_{j \leq k} \grade{j}{(\maxel{U})}$.
It follows that $x \in \bound{k}{\alpha}{U}$, so $\clset{{x}} \subseteq \bound{k}{\alpha}U$.

Suppose $k < \dim{x}$ and let $y \in \bound{k}{\alpha}x$.
By Corollary \ref{cor:atoms_are_round} and Lemma \ref{lem:round_is_pure}, $\bound{k}{\alpha}x$ is pure and $k$\nbd dimensional, so there exists $z \in \faces{k}{\alpha}x$ such that $y \in \clset{{z}}$.
Because $k > \frdim{U}$, $x$ is the only maximal element above $z$, so $\cofaces{}{}z \subseteq \clset{{x}}$.
It follows that $z \in \faces{k}{\alpha}U$, hence $y \in \bound{k}{\alpha}U$.
\end{proof}

\begin{prop} \label{prop:intersection_of_maximal_elements}
Let $U$ be a molecule, $x, y \in \maxel{U}$, $x \neq y$.
For all $k \geq \frdim{U}$,
\begin{equation*}
    \clset{{x}} \cap \clset{{y}} = (\bound{k}{-}x \cap \bound{k}{+}y) \cup (\bound{k}{+}x \cap \bound{k}{-}y).
\end{equation*}
\end{prop}
\begin{proof}
Suppose $k \geq \min{\set{\dim{x}, \dim{y}}}$.
Without loss of generality, suppose $k \geq \dim{x}$.
Since $x$ is maximal, $x \in \bound{k}{-}U \cap \bound{k}{+}U$.
Then
\begin{equation*}
    \clset{{x}} \cap \clset{{y}} \subseteq (\bound{k}{-}U \cap \bound{k}{+}U) \cap \clset{{y}} \subseteq \bound{k}{-}y \cap \bound{k}{+}y
\end{equation*}
by Lemma \ref{lem:intersection_with_boundary}.
Since $\clset{{x}} = \bound{k}{-}x = \bound{k}{+}x$, we are done.

We now proceed by induction on $k' \eqdef \min{\set{\dim{x}, \dim{y}}} - k \geq 0$, having established the base case $k' = 0$.
Suppose $k' > 0$, and consider an element $z \in \clset{{x}} \cap \clset{{y}}$.
By the inductive hypothesis, there exists $\alpha \in \set{+, -}$ such that $z \in \bound{k+1}{\alpha}x \cap \bound{k+1}{-\alpha}y$.
Since $k + 1 > k \geq \frdim{U}$, by Lemma \ref{lem:frdim_boundaries} we have $z \in \bound{k+1}{\alpha}U \cap \bound{k+1}{-\alpha}U$, so
\begin{align*}
    z \in \bound{k+1}{\alpha}x \cap \bound{k+1}{-\alpha}U & = \bound{k+1}{\alpha}x \cap \bound{k+1}{-\alpha}x = \bound{k}{}x, \\
    z \in \bound{k+1}{\alpha}U \cap \bound{k+1}{-\alpha}y & = \bound{k+1}{\alpha}y \cap \bound{k+1}{-\alpha}y = \bound{k}{}y,
\end{align*}
using Lemma \ref{lem:intersection_with_boundary} and roundness of $\clset{{x}}$ and $\clset{{y}}$.

Suppose there exists $\beta \in \set{ +, - }$ such that $z \in \bound{k}{\beta}x \cap \bound{k}{\beta}y$.
We want to prove that $z \in \bound{k}{-\beta}x$ or $z \in \bound{k}{-\beta}y$.
Without loss of generality, suppose $\beta = +$.
We have the following facts:
\begin{enumerate}
    \item $\bound{k}{+}x$ and $\bound{k}{+}y$ are pure and $k$\nbd dimensional,
    \item since $k + 1 \leq \min{\set{\dim{x}, \dim{y}}}$, both $\bound{k+1}{\alpha}x$ and $\bound{k+1}{\alpha}y$ are pure and $(k+1)$\nbd dimensional,
    \item $\bound{k}{+}x \subseteq \bound{k+1}{\alpha}x$ and $\bound{k}{+}y \subseteq \bound{k+1}{\alpha}y$ by globularity,
    \item $\bound{k+1}{\alpha}x \cup \bound{k+1}{\alpha}y \subseteq \bound{k+1}{\alpha}U$ which is a $(k+1)$\nbd dimensional molecule by Lemma \ref{lem:frdim_boundaries}.
\end{enumerate}
We will construct a sequence of pairs $(x_i, y_j)$ with the following properties:
\begin{enumerate}
    \item $x_i \in \faces{k+1}{\alpha}x$ and $y_j \in \faces{k+1}{\alpha}y$,
    \item $z \in \clset{{x_i}} \cap \clset{{y_j}}$,
    \item there exists $\gamma(i,j) \in \set{+, -}$ such that \begin{equation*}
        \clset{{x_i}} \cap \clset{{y_j}} = \bound{}{\gamma(i,j)}x_i \cap \bound{}{-\gamma(i,j)}y_j,
    \end{equation*}
    \item if $i' > i$, there exists a non-trivial path from $x_{i'}$ to $x_i$ in $\flow{k}{(\bound{k+1}{\alpha}x)}$,
    \item if $j' > j$, there exists a non-trivial path from $y_{j'}$ to $y_j$ in $\flow{k}{(\bound{k+1}{\alpha}y)}$.
\end{enumerate}
To start, since $\bound{k}{+}x$ and $\bound{k}{+}y$ are pure and $k$\nbd dimensional, there exist $x'_0 \in \faces{k}{+}x$ and $y'_0 \in \faces{k}{+}y$ such that $z \in \clset{{x'_0}} \cap \clset{{y'_0}}$.
Because $\bound{k+1}{\alpha}x$ and $\bound{k+1}{\alpha}y$ are also pure and $(k+1)$\nbd dimensional, we can pick
\begin{equation*}
    x_0 \in \faces{k+1}{\alpha}x \cap \cofaces{}{+}x'_0, \quad  y_0 \in \faces{k+1}{\alpha}y \cap \cofaces{}{+}y'_0.
\end{equation*}
By Lemma \ref{lem:layering_dimension_smaller_than_dimension} $\lydim{\bound{k+1}{\alpha}U} \leq k$, so by Theorem \ref{thm:molecules_admit_layerings} there exists a $k$\nbd layering $(\order{i}{V})_{i=1}^m$ of $\bound{k+1}{\alpha}U$.
By Lemma \ref{lem:layering_intersections} and Lemma \ref{lem:intersection_with_boundary}, there exist $i \neq j$ and $\gamma \in \set{+, -}$ such that $x_0 \in \order{i}{V}$, $y_0 \in \order{j}{V}$, and
\begin{equation*}
    \clset{{x_0}} \cap \clset{{y_0}} = (\clset{{x_0}} \cap \bound{k}{\gamma}\order{i}{V}) \cap (\clset{{y_0}} \cap \bound{k}{-\gamma}\order{j}{V}) = \bound{}{\gamma}x_0 \cap \bound{}{-\gamma}y_0.
\end{equation*}
We deduce that $z \in \bound{}{\gamma}x_0 \cap \bound{}{-\gamma}y_0$.

The condition on paths in $\flow{k}{(\bound{k+1}{\alpha}x)}$ and $\flow{k}{(\bound{k+1}{\alpha}y)}$ is trivially satisfied, so we have defined a pair $(x_0, y_0)$ satisfying all the properties with $\gamma(0, 0) \eqdef \gamma$.
This will form the base of the induction.

Suppose we have defined $(x_i, y_j)$ satisfying all the conditions.
If $\gamma(i, j) = +$, the next pair will be of the form $(x_i, y_{j+1})$; if $\gamma(i, j) = -$ it will be of the form $(x_{i+1}, y_j)$.
Suppose without loss of generality that $\gamma(i, j) = +$, that is,
\begin{equation*}
    \clset{{x_i}} \cap \clset{{y_j}} = \bound{}{+}x_i \cap \bound{}{-}y_j.
\end{equation*}
Since $\bound{k}{-}y_j$ is pure and $k$\nbd dimensional, there exists $w \in \faces{}{-}y_j$ such that $z \in \clset{{w}}$.
If $w \in \faces{k}{-}y$, we have proved that $z \in \bound{k}{-}y$.
Otherwise, since $w$ is a codimension-1 element of $\bound{k+1}{\alpha}y$, which is pure and $(k+1)$\nbd dimensional, by Corollary \ref{cor:codimension_1_elements} there exists a single $y_{j+1} \in \faces{k+1}{\alpha}y \cap \cofaces{}{+}w$.
Reasoning as before with Lemma \ref{lem:layering_intersections}, we find that
\begin{equation*}
    z \in \clset{{x_i}} \cap \clset{{y_{j+1}}} = \bound{}{\gamma'}x_i \cap \bound{}{-\gamma'}y_{j+1}
\end{equation*}
for some $\gamma' \in \set{+, -}$.
Moreover, since $w \in \faces{}{+}y_{j+1} \cap \faces{}{-}y_j$, there exists an edge from $y_{j+1}$ to $y_j$ in $\flow{k}{(\bound{k+1}{\alpha}y)}$.

The case $\gamma(i, j) = -$ is dual.
By Corollary \ref{cor:flow_acyclic_in_codimension_1}, paths in $\flow{k}{(\bound{k+1}{\alpha}x)}$ and $\flow{k}{(\bound{k+1}{\alpha}y)}$ are finite.
Eventually, then, we reach $x_i$ such that $\faces{}{-}x_i \subseteq \faces{k}{-}x$, proving that $z \in \bound{k}{-}x$; or we reach $y_j$ such that $\faces{}{-}y_j \subseteq \faces{k}{-}y$, proving that $z \in \bound{k}{-}y$.
\end{proof}


\section{Rewritable submolecules}

\begin{guide}
	In this section, using the theory developed so far, we produce a criterion to decide whether the inclusion of a round molecule into a molecule of the same dimension is a submolecule inclusion 
	(Theorem \ref{thm:rewritable_submolecule_criterion}).
	This amounts, in essence, to a decision algorithm, which is the one made explicit and analysed in \cite{hadzihasanovic2023higher}.
\end{guide}

\begin{dfn}[Rewritable submolecule] \index{submolecule!rewritable}
Let $U$ be a molecule.
A submolecule $V \submol U$ is \emph{rewritable} if $\dim{U} = \dim{V}$ and $V$ is round.
\end{dfn}

\begin{prop} \label{prop:round_submolecule_substitution}
Let $\imath\colon V \incl U$ be an inclusion of molecules such that $\dim{V} = \dim{U}$ and $V$ is round.
The following are equivalent:
\begin{enumerate}[label=(\alph*)]
    \item $\imath$ is a submolecule inclusion;
    \item for all molecules $W$ such that $V \celto W$ is defined, $\subs{U}{W}{\imath(V)}$ is a molecule and $j\colon W \incl \subs{U}{W}{\imath(V)}$ is a submolecule inclusion;
    \item $\subs{U}{\compos{V}}{\imath(V)}$ is a molecule.
\end{enumerate}
\end{prop}
\begin{proof}
If $\imath$ is a submolecule inclusion, by Lemma \ref{lem:submolecule_rewrite} $U \cup (V \celto W)$ and its output boundary $U[W/\imath(V)]$ are molecules, and the inclusion of $W$ into $U[W/\imath(V)]$ is a submolecule inclusion.

If $V$ is a round molecule, then $\compos{V}$ is an atom, which is round by Corollary \ref{cor:atoms_are_round}, and has boundaries isomorphic to those of $V$ by Lemma \ref{lem:boundaries_of_rewrite}.
By Corollary \ref{cor:molecule_boundary_unique_isomorphism}, $V \celto \compos{V}$ is defined, so the third condition is a special case of the second one.

Finally, suppose $\subs{U}{\compos{V}}{\imath(V)}$ is a molecule.
By Lemma \ref{lem:downset_is_submolecule}, since $\compos{V}$ is an atom, its inclusion $j$ into $\subs{U}{\compos{V}}{\imath(V)}$ is a submolecule inclusion.
Using Lemma \ref{lem:submolecule_rewrite} as in the first part, we deduce that $\subs{(\subs{U}{\compos{V}}{\imath(V)})}{V}{j(\compos{V})}$ is a molecule, and the inclusion of $V$ into it is a submolecule inclusion.
By Lemma \ref{lem:revert_substitution}, $\subs{(\subs{U}{\compos{V}}{\imath(V)})}{V}{j(\compos{V})}$ is isomorphic to $U$, and $\imath$ factors as this submolecule inclusion followed by an isomorphism.
\end{proof}

\begin{dfn}[Path-induced subgraph] \index{directed graph!subgraph!path-induced}
Let $\mathscr{G}$ be a directed graph and $W \subseteq V_\mathscr{G}$.
We say that $\restr{\mathscr{G}}{W}$ is \emph{path-induced} if, for all $x, y \in W$, every path from $x$ to $y$ in $\mathscr{G}$ is included in $\restr{\mathscr{G}}{W}$.
\end{dfn}

\begin{comm}
Path-induced subgraphs are also called \emph{convex subgraphs}, for example in \cite{bonchi2022string}.
\end{comm}

\begin{dfn}[Contraction of a connected subgraph] \index{directed graph!subgraph!contraction}
Let $\mathscr{G}$ be a directed graph and $W \subseteq V_\mathscr{G}$ such that $\restr{\mathscr{G}}{W}$ is connected.
The \emph{contraction of $\restr{\mathscr{G}}{W}$ in $\mathscr{G}$} is the graph minor $\mathscr{G}/(\restr{\mathscr{G}}{W})$ obtained by contracting every edge in $\restr{\mathscr{G}}{W}$.

Explicitly, the set of vertices of $\mathscr{G}/(\restr{\mathscr{G}}{W})$ is $(V_\mathscr{G} \setminus W) + \set{x_W}$, and for all pair of vertices $x, y$,
\begin{itemize}
    \item if $x, y \neq x_W$, there is an edge between $x$ and $y$ for each edge between $x$ and $y$ in $\mathscr{G}$,
    \item if $x = x_W$ and $y \neq x_W$, there is an edge from $x$ to $y$ for each pair of a vertex $z \in W$ and an edge from $x$ to $y$ in $\mathscr{G}$,
    \item if $x \neq x_W$ and $y = x_W$, there is an edge from $x$ to $y$ for each pair of a vertex $z \in W$ and an edge from $x$ to $z$ in $\mathscr{G}$,
    \item there are no edges from $x_W$ to $x_W$.
\end{itemize}
\end{dfn}

\begin{lem} \label{lem:connected_subgraph_conditions_path_induced}
Let $\mathscr{G}$ be a directed acyclic graph and $W \subseteq V_\mathscr{G}$ such that $\restr{\mathscr{G}}{W}$ is connected.
The following are equivalent:
\begin{enumerate}[label=(\alph*)]
    \item $\restr{\mathscr{G}}{W}$ is path-induced; \label{cond:path_induced}
    \item $\mathscr{G}/(\restr{\mathscr{G}}{W})$ is acyclic; \label{cond:acyclic_contraction}
    \item there is a topological sort of $\mathscr{G}$ in which vertices of $W$ are consecutive. \label{cond:consecutive_tsort}
\end{enumerate}
Moreover, under any of the equivalent conditions, there is a bijection between
\begin{itemize}
    \item topological sorts of $\mathscr{G}$ in which vertices of $W$ are consecutive,
    \item pairs of a topological sort of $\restr{\mathscr{G}}{W}$ and a topological sort of $\mathscr{G}/(\restr{\mathscr{G}}{W})$.
\end{itemize}
\end{lem}
\begin{proof}
We prove the contrapositive of the implication from \ref{cond:path_induced} to \ref{cond:acyclic_contraction}.
Suppose $\mathscr{G}/(\restr{\mathscr{G}}{W})$ has a cycle.
If the cycle does not pass through $x_W$, then it lifts to a cycle in $\mathscr{G}$, contradicting the assumption that $\mathscr{G}$ is acyclic.
It follows that the cycle contains a segment of the form $x_W \to x_1 \to \ldots \to x_m \to x_W$, where $m > 0$ and $x_i \neq x_W$ for all $i \in \set{1, \ldots, m}$.
Then there exist $y, z \in W$ and a path $y \to x_1 \to \ldots \to x_m \to z$ in $\mathscr{G}$, so $\restr{\mathscr{G}}{W}$ is not path-induced.

Next, suppose that $\mathscr{G}/(\restr{\mathscr{G}}{W})$ is acyclic.
Then both $\mathscr{G}/(\restr{\mathscr{G}}{W})$ and $\restr{\mathscr{G}}{W}$ are acyclic, so they admit topological sorts $(\order{i}{x})_{i=1}^m$ and $(\order{j}{y})_{j=1}^p$, respectively.
For exactly one $q \in \set{1, \ldots, m}$, $\order{i}{x} = x_W$.
We claim that 
\begin{equation*}
    ((\order{i}{x})_{i=1}^{q-1}, (\order{j}{y})_{j=1}^p, (\order{i}{x})_{i=q+1}^m)
\end{equation*}
is a topological sort of $\mathscr{G}$.
Indeed, for all edges from $x$ to $x'$ in $\mathscr{G}$,
\begin{itemize}
    \item if $x, x' \notin W$, then $x = \order{i}{x}$, $x' = \order{i'}{x}$ for some $i, i' \in \set{1, \ldots, m} \setminus \set{q}$, and there is an edge from $x$ to $x'$ in $\mathscr{G}/(\restr{\mathscr{G}}{W})$, so $i < i'$;
    \item if $x, x' \in W$, then $x = \order{j}{y}$, $x' = \order{j'}{y}$ for some $j, j' \in \set{1, \ldots, p}$, and there is an edge from $x$ to $x'$ in $\restr{\mathscr{G}}{W}$, so $j < j'$;
    \item if $x \in W$, $x' \notin W$, then $x = \order{j}{y}$, $x' = \order{i}{x}$ for some $i \in \set{1, \ldots, m} \setminus \set{q}$, $j \in \set{1, \ldots, p}$, and there is an edge from $x_W$ to $x'$ in $\mathscr{G}/(\restr{\mathscr{G}}{W})$, so $q < i$;
    \item if $x \notin W$, $x' \in W$, then $x = \order{i}{x}$, $x' = \order{j}{y}$ for some $i \in \set{1, \ldots, m} \setminus \set{q}$, $j \in \set{1, \ldots, p}$, and there is an edge from $x$ to $x_W$ in $\mathscr{G}/(\restr{\mathscr{G}}{W})$, so $i < q$.
\end{itemize}
This proves the implication from \ref{cond:acyclic_contraction} to \ref{cond:consecutive_tsort}.
Moreover, it defines an injection from pairs of a topological sort of $\restr{\mathscr{G}}{W}$ and a topological sort of $\mathscr{G}/(\restr{\mathscr{G}}{W})$ to topological sorts of $\mathscr{G}$ in which the vertices of $W$ are consecutive.
This will prove to be a bijection as soon as we have proven the converse implication.

Finally, we prove the contrapositive of the implication from \ref{cond:consecutive_tsort} to \ref{cond:path_induced}.
Suppose $\restr{\mathscr{G}}{W}$ is not path-induced, that is, there is a path $x \to x_1 \to \ldots \to x_m \to y$ in $\mathscr{G}$ such that $m > 0$, $x, y \in W$, and $x_i \notin W$ for all $i \in \set{1, \ldots, m}$.
It follows that the $x_i$ must come between $x$ and $y$ in every topological sort of $\mathscr{G}$, so the vertices of $W$ can never be consecutive.
\end{proof}

\begin{prop} \label{prop:round_molecule_connected_flowgraph}
Let $U$ be a molecule, $n \eqdef \dim{U}$.
If $U$ is round, then $\flow{n-1}{U}$ is connected.
\end{prop}
\begin{proof}
First of all, if $U$ is round, then it is pure, so the vertices of $\flow{n-1}{U}$ are the elements of $\grade{n}{U}$.
If $U$ is an atom, then $\flow{n-1}{U}$ consists of a single vertex and no edges, so it is trivially connected.
In particular this is true when $n = 0$ by Lemma \ref{lem:only_0_molecule}, so we can proceed by induction on $n$.

Suppose $n > 0$ and $\size{\grade{n}{U}} > 1$, which by Lemma \ref{lem:layering_dimension_smaller_than_dimension} implies $\lydim{U} = n - 1$.
Assume by way of contradiction that $\flow{n-1}{U}$ is not connected.
Then there is a bipartition $\grade{n}{U} = A + B$ such that there are no edges in $\flow{n-1}{U}$ between vertices in $A$ and vertices in $B$.
By Corollary \ref{cor:codimension_1_elements}, no element of codimension 1 in $U$ can be covered by two elements with the same orientation, so this implies that $\dim{(\clos{A} \cap \clos{B})} < n-1$.
Let 
\begin{equation*}
    A' \eqdef \set{x \in \faces{}{-}U \mid \cofaces{}{-}x \subseteq A}, \quad \quad
    B' \eqdef \set{x \in \faces{}{-}U \mid \cofaces{}{-}x \subseteq B}.
\end{equation*}
Then $A' + B'$ is a bipartition of $\faces{}{-}U$.
By Lemma \ref{lem:boundaries_of_round_and_globular}, $\bound{}{-}U$ is round, so by the inductive hypothesis $\flow{n-2}{(\bound{}{-}U)}$ is connected.
It follows that there exist $\alpha \in \set{+, -}$, $x \in A'$, $y \in B'$, and $z \in \grade{n-2}{U}$ such that $z \in \faces{}{\alpha}x \cap \faces{}{-\alpha}y$.
Then $z$ has two distinct cofaces in $\bound{}{-}U$, so by Corollary \ref{cor:codimension_1_elements} $z \notin \bound{}{}(\bound{}{-}U) = \bound{n-2}{}U$.
We claim that $z \in \bound{}{+}U$, contradicting the roundness of $U$.

By Theorem \ref{thm:molecules_admit_layerings}, there exists an $(n-1)$\nbd layering $(\order{i}{U})_{i=1}^m$ of $U$; we will identify the $\order{i}{U}$ with their isomorphic images in $U$.
Let $V_0 \eqdef \bound{}{-}U$ and $V_i \eqdef \bound{}{+}\order{i}{U}$ for each $i \in \set{1, \ldots, m}$.
We will prove that, for all $i \in \set{0, \ldots, m}$,
\begin{enumerate}
    \item $z \in V_i$,
    \item there exist $x_i \in \clos{A}$ and $y_i \in \clos{B}$ such that $\cofaces{}{\alpha} z \cap V_i = \set{x_i}$ and $\cofaces{}{-\alpha} z \cap V_i = \set{y_i}$.
\end{enumerate}
For $i = 0$, we have already established this with $x_0 \eqdef x$, $y_0 \eqdef y$.
Let $i \geq 0$, and assume this holds for $i - 1$.
By Lemma \ref{lem:layering_basic_properties}, there is a single $n$\nbd dimensional element $\order{i}{x}$ in $\order{i}{U}$, and by Lemma \ref{lem:boundary_move}
\begin{equation*}
    V_i = \subs{\bound{}{-}\order{i}{U}}{\bound{}{+}\order{i}{x}}{\bound{}{-}\order{i}{x}} = \subs{V_{i-1}}{\bound{}{+}\order{i}{x}}{\bound{}{-}\order{i}{x}}.
\end{equation*}
Suppose $\order{i}{x} \in A$.
Then $y_{i-1} \notin \clset{{\order{i}{x}}}$, so $y_{i-1} \in V_i$, and we let $y_i \eqdef y_{i-1}$.
If $x_{i-1} \notin \clset{{\order{i}{x}}}$ then also $x_{i-1} \in V_i$, and we let $x_i \eqdef x_{i-1}$.
Otherwise, $x_{i-1}$ is the only coface of $z$ in $\bound{}{-}\order{i}{x}$, so by Corollary \ref{cor:codimension_1_elements} we have $z \in \bound{}{\alpha}(\bound{}{-}\order{i}{x}) = \bound{}{\alpha}(\bound{}{+}\order{i}{x})$.
It follows that $z \in V_i$ and there exists a unique $x_i$ such that $\cofaces{}{\alpha}z \cap \bound{}{+}\order{i}{x} = \set{x_i}$.
The case $\order{i}{x} \in B$ is analogous.

Since $V_m = \bound{}{+}U$, we have proved that $z \in \bound{}{+}U$, a contradiction.
\end{proof}

\begin{exm}[A pure 3-dimensional molecule with connected 2-flow which is not round] \index[counterex]{A pure 3-dimensional molecule with connected 2-flow which is not round}
	The converse of Proposition \ref{prop:round_molecule_connected_flowgraph} does not hold in general --- not even when a molecule is \emph{pure}.
	Let $U$ be the 3\nbd dimensional molecule whose oriented Hasse diagram is
	\[
		\input{img/connflow_hasse.tex}
	\]
	which admits the 2\nbd layering $\order{1}{U} \cp{2} \order{2}{U}$, where
	\begin{align*}
		\order{1}{U} & \eqdef (\disk{1}{2} \celto (\disk{1}{2} \cp{1} (\globe{2} \cp{0} \globe{2}))) \cp{1} \disk{2}{1}, \\
		\order{2}{U} & \eqdef \disk{1}{2} \cp{1} (((\globe{2} \cp{0} \globe{2}) \cp{1} \disk{2}{1}) \celto \disk{2}{1}),
	\end{align*}
	corresponding to the sequence of rewrite steps
	\[
		\input{img/connflow_step0.tex} \quad
		\input{img/connflow_step1.tex} \quad
		\input{img/connflow_step2.tex}
	\]
	in string diagrams.
	Then $\flow{2}{U}$ is equal to
\[\begin{tikzcd}
	{{\scriptstyle (3, 0)}\;\bullet} && {{\scriptstyle (3, 1)}\;\bullet}
	\arrow[from=1-1, to=1-3]
\end{tikzcd}\]
	which is connected, and $U$ is pure, but $U$ is not round, since
	\[
		(0, 2) \in \bound{2}{+}U \cap \bound{2}{-}U, \quad \quad
		(0, 2) \not\in \bound{1}{}U.
	\]
\end{exm}

\begin{lem} \label{lem:flow_of_substitution_is_contraction_of_flow}
Let $\imath\colon V \incl U$ be an inclusion of molecules of the same dimension $n$, and suppose $V$ is round.
Then $\flow{n-1}{\subs{U}{\compos{V}}{\imath(V)}}$ is isomorphic to $\flow{n-1}{U}/\flow{n-1}{V}$.
\end{lem}
\begin{proof}
By Lemma \ref{lem:flow_under_inclusion} and Proposition \ref{prop:round_molecule_connected_flowgraph}, $\flow{n-1}{V}$ is a connected induced subgraph of $\flow{n-1}{U}$, so its contraction is well-defined.
Now, the vertices of $\flow{n-1}{\subs{U}{\compos{V}}{\imath(V)}}$ are either
\begin{itemize}
    \item $x \in \grade{n}{U} \setminus \grade{n}{V}$, or
    \item $x_V$ such that the image of $\compos{V}$ in $\subs{U}{\compos{V}}{\imath(V)}$ is $\clset{{x_V}}$.
\end{itemize}
Let $x, y$ be two vertices of $\flow{n-1}{\subs{U}{\compos{V}}{\imath(V)}}$.
\begin{itemize}
\item If $x, y \in \grade{n}{U} \setminus \grade{n}{V}$, then $\faces{}{+}x \cap \faces{}{-}y$ is the same in $\subs{U}{\compos{V}}{\imath(V)}$ as in $U$, so there is an edge from $x$ to $y$ in $\flow{n-1}{U}$ if and only if there is an edge in $\flow{n-1}{\subs{U}{\compos{V}}{\imath(V)}}$.
\item If $x = x_V$ then $\faces{}{+}x_V \cap \faces{}{-}y$ is in bijection with $\faces{}{+}V \cap \faces{}{-}y$ in $U$.
For all $z \in \faces{}{+}V$, since $V$ is pure and $n$\nbd dimensional, there exists $w \in \cofaces{}{+}z$.
If $\faces{}{+}x_V \cap \faces{}{-}y$ is non-empty, it follows that $\faces{}{+}z \cap \faces{}{-}y$ is non-empty in $U$ for some $z \in \grade{n}{\imath(V)}$.
Thus there exist $z \in \grade{n}{\imath(V)}$ and an edge from $z$ to $y$ in $\flow{n-1}{U}$.
\item Dually, if $y = x_V$, there is an edge from $x$ to $y$ in $\flow{n-1}{\subs{U}{\compos{V}}{\imath(V)}}$ if and only if there exist $z \in \grade{n}{\imath(V)}$ and an edge from $x$ to $z$ in $\flow{n-1}{U}$.
\item Finally, $\faces{}{+}V \cap \faces{}{-}V = \varnothing$ because $V$ is pure, so $\faces{}{+}x_V \cap \faces{}{-}x_V$ and there is no edge from $x_V$ to $x_V$.
\end{itemize}
It is then straightforward to establish an isomorphism with the explicit description of $\flow{n-1}{U}/\flow{n-1}{V}$.
\end{proof}

\begin{prop} \label{prop:round_submolecule_flow_path_induced}
Let $\imath\colon V \incl U$ be an inclusion of molecules such that $n \eqdef \dim{U} = \dim{V}$ and $V$ is round.
If $\imath$ is a submolecule inclusion, then $\flow{n-1}{V}$ is a path-induced subgraph of $\flow{n-1}{U}$.
\end{prop}
\begin{proof}
By Proposition \ref{prop:round_submolecule_substitution}, if $\imath$ is a submolecule inclusion then $\subs{U}{\compos{V}}{V}$ is a molecule.
By Corollary \ref{cor:flow_acyclic_in_codimension_1} $\flow{n-1}{\subs{U}{\compos{V}}{V}}$ is acyclic, and by Lemma \ref{lem:flow_of_substitution_is_contraction_of_flow} it is isomorphic to $\flow{n-1}{U}/\flow{n-1}{V}$.
It follows from Lemma \ref{lem:connected_subgraph_conditions_path_induced} that $\flow{n-1}{V}$ is a path-induced subgraph of $\flow{n-1}{U}$.
\end{proof}

\begin{lem} \label{lem:round_submolecules_from_layering}
Let $\imath\colon V \incl U$ be an inclusion of molecules of the same dimension $n$, suppose $V$ is round, and suppose $(\order{i}{y})_{i=1}^p$ is an $(n-1)$\nbd ordering induced by an $(n-1)$\nbd layering of $V$.
The following are equivalent:
\begin{enumerate}[label=(\alph*)]
    \item $\imath$ is a submolecule inclusion;
    \item there exist an $(n-1)$\nbd ordering $(\order{i}{x})_{i=1}^m$ induced by an $(n-1)$\nbd layering $(\order{i}{U})_{i=1}^m$ of $U$, and $q \in \set{1, \ldots, m}$ such that
    \begin{enumerate}[label=\arabic*.]
        \item $(\order{i}{x})_{i=q}^{p+q-1} = (\imath(\order{i}{y}))_{i=1}^p$,
        \item $\imath(\bound{}{-}V) \submol \bound{}{-}\order{q}{U}$.
    \end{enumerate}
\end{enumerate}
\end{lem}
\begin{proof}
Identify $V$ with its isomorphic image through $\imath$, and suppose that $\imath$ is a submolecule inclusion.
Then $\tilde{U} \eqdef \subs{U}{\compos{V}}{V}$ is a molecule by Proposition \ref{prop:round_submolecule_substitution}, and admits an $(n-1)$\nbd layering $(\order{i}{\tilde{U}})_{i=1}^{m-p+1}$ by Theorem \ref{thm:molecules_admit_layerings}.
Let $\clset{{x}}$ be the image of $\compos{V}$ in $\tilde{U}$; then $x \in \order{q}{\tilde{U}}$ for exactly one $q \in \set{1, \ldots, m-p+1}$.
Then $W \eqdef \subs{\order{q}{\tilde{U}}}{V}{\clset{{x}}}$ is defined, and by Lemma \ref{lem:pasting_after_substitution} combined with Lemma \ref{lem:revert_substitution}, $U$ is isomorphic to
\begin{equation*}
    \order{1}{\tilde{U}} \cp{n-1} \ldots \cp{n-1} \order{q-1}{\tilde{U}} \cp{n-1} W \cp{n-1} \order{q+1}{\tilde{U}} \cp{n-1} \ldots \cp{n-1} \order{m-p+1}{\tilde{U}}.
\end{equation*}
By Lemma \ref{lem:boundary_move}, $\bound{}{-}x \submol \bound{}{-}\order{q}{\tilde{U}}$, so by Lemma \ref{lem:substitution_preserves_boundaries} $\bound{}{-}V \submol \bound{}{-}W$.
We can apply the criterion of Proposition \ref{prop:layering_from_ordering} to deduce that $(\order{i}{y})_{i=1}^p$ is an $(n-1)$\nbd ordering of $W$ induced by an $(n-1)$\nbd layering $(\order{i}{W})_{i=1}^p$.
Letting 
\begin{equation*}
    (\order{i}{U})_{i=1}^m \eqdef ((\order{i}{\tilde{U}})_{i=1}^{q-1}, (\order{i}{W})_{i=1}^p, (\order{i}{\tilde{U}})_{i=q+1}^{m-p+1}),
\end{equation*}
produces an $(n-1)$\nbd layering of $U$, hence also an $(n-1)$\nbd ordering $(\order{i}{x})_{i=1}^m$ of $U$, with the property that $(\order{i}{x})_{i=q}^{p+q-1} = (\order{i}{y})_{i=1}^p$.

Conversely, let $(\order{i}{U})_{i=1}^m$ be an $(n-1)$\nbd layering of $U$ satisfying the properties in the statement, and let $W \submol U$ be the image of $\order{q}{U} \cp{n-1} \ldots \cp{n-1} \order{p+q-1}{U}$ in $U$.
Then $\grade{n}{W} = \grade{n}{V}$, so 
\begin{equation*}
    W = V \cup \bound{}{-}W.
\end{equation*}
Because $\bound{}{-}V \submol \bound{}{-}\order{q}{U} = \bound{}{-}W$, by Lemma \ref{lem:submolecule_rewrite} $V \submol W \submol U$.
\end{proof}

\begin{thm} \label{thm:rewritable_submolecule_criterion}
Let $\imath\colon V \incl U$ be an inclusion of molecules such that $n \eqdef \dim{U} = \dim{V}$ and $V$ is round, $m \eqdef \size{\grade{n}{U}}$, $p \eqdef \size{\grade{n}{V}}$.
The following are equivalent:
\begin{enumerate}[label=(\alph*)]
    \item $\imath$ is a submolecule inclusion;
    \item there is a topological sort $((\order{i}{x})_{i=1}^{q-1}, x_V, (\order{i}{x})_{i=q+1}^{m-p+1})$ of $\flow{n-1}{U}/\flow{n-1}{V}$ such that, letting
    \begin{align*}
    \order{0}{U} & \eqdef \bound{}{-}U, \\
    \order{q}{U} & \eqdef \bound{n-1}{+}\order{q-1}{U} \cup \imath(V), \\
    \order{i}{U} & \eqdef \bound{n-1}{+}\order{i-1}{U} \cup \clset{{\order{i}{x}}} \quad \text{for $i \neq q$},
    \end{align*}
    we have $\imath(\bound{}{-}V) \submol \bound{}{-}\order{q}{U}$ and $\bound{}{-}\order{i}{x} \submol \bound{}{-}\order{i}{U}$ for all $i \neq q$.
\end{enumerate}
\end{thm}
\begin{proof}
Identify $V$ with its isomorphic image through $\imath$, and suppose that $\imath$ is a submolecule inclusion.
Then $\tilde{U} \eqdef \subs{U}{\compos{V}}{V}$ is a molecule by Proposition \ref{prop:round_submolecule_substitution}, so it admits an $(n-1)$\nbd layering $(\order{i}{\tilde{U}})_{i=1}^{m-p+1}$, which induces an $(n-1)$\nbd ordering.
By Lemma \ref{lem:flow_of_substitution_is_contraction_of_flow}, this $(n-1)$\nbd ordering can be identified with a topological sort $((\order{i}{x})_{i=1}^{q-1}, x_V, (\order{i}{x})_{i=q+1}^{m-p+1})$ of $\flow{n-1}{U}/\flow{n-1}{V}$.
By Lemma \ref{lem:boundary_move}, we have $\bound{}{-}x_V \submol \bound{}{-}\order{q}{\tilde{U}}$ and $\bound{}{-}\order{i}{x} \submol \bound{}{-}\order{i}{\tilde{U}}$ for $i \neq q$.
By Lemma \ref{lem:pasting_after_substitution} combined with Lemma \ref{lem:revert_substitution}, letting $W \eqdef \subs{\order{q}{\tilde{U}}}{V}{\clset{{x_V}}}$, $U$ is isomorphic to
\begin{equation*}
    \order{1}{\tilde{U}} \cp{n-1} \ldots \cp{n-1} \order{q-1}{\tilde{U}} \cp{n-1} W \cp{n-1} \order{q+1}{\tilde{U}} \cp{n-1} \ldots \cp{n-1} \order{m-p+1}{\tilde{U}},
\end{equation*}
and $W$ is isomorphic to $\order{q}{U}$, while $\order{i}{\tilde{U}}$ is isomorphic to $\order{i}{U}$ for all $i \neq q$.
We conclude by Lemma \ref{lem:substitution_preserves_boundaries}.

Conversely, let $(\order{i}{y})_{i=1}^p$ be an $(n-1)$\nbd ordering induced by an $(n-1)$\nbd layering of $V$.
Then $((\order{i}{x})_{i=1}^{q-1}, (\order{i}{y})_{i=1}^p, (\order{i}{x})_{i=q+1}^{m-p+1})$ is an $(n-1)$\nbd ordering of $U$, and by the criterion of Proposition \ref{prop:layering_from_ordering} it is induced by an $(n-1)$\nbd layering. 
We conclude by Lemma \ref{lem:round_submolecules_from_layering}.
\end{proof}

\begin{exm}[An inclusion of a round molecule which is not a submolecule inclusion] \index[counterex]{An inclusion of a round molecule which is not a submolecule inclusion} \label{exm:non_submol}
	Let $U$ be the 3\nbd dimensional molecule whose oriented Hasse diagram is
	\[
		\input{img/nonsubmol_hasse.tex}
	\]
	which admits the 2\nbd layering $\order{1}{U} \cp{2} \order{2}{U} \cp{2} \order{3}{U}$, where
	\begin{align*}
		\order{1}{U} & \eqdef 
		\disk{1}{2} \cp{1} ((\globe{2} \celto (\globe{2} \cp{1} \globe{2})) \cp{0} \thearrow{}) \cp{1} \disk{2}{1}, \\
		\order{2}{U} & \eqdef
		\disk{1}{2} \cp{1} (\globe{2} \cp{0} \thearrow{}) \cp{1} (((\globe{2} \cp{0} \thearrow{}) \cp{1} \disk{2}{1}) \celto ((\thearrow{} \cp{0} \globe{2}) \cp{1} \disk{2}{1})), \\
		\order{3}{U} & \eqdef
		((\disk{1}{2} \cp{1} (\thearrow{} \cp{0} \globe{2})) \celto \disk{1}{2}) \cp{1} \disk{2}{1}
	\end{align*}
	corresponding to the 2\nbd ordering $(3, 1), (3, 2), (3, 0)$ and the sequence of rewrite steps
	\[
		\input{img/nonsubmol_step0.tex} \quad
		\input{img/nonsubmol_step1.tex}
	\]
	\[
		\input{img/nonsubmol_step2.tex} \quad
		\input{img/nonsubmol_step3.tex}
	\]
	in string diagrams.
	Then $V \eqdef \clset{(3, 0), (3, 1)} \subseteq U$ is a round 3\nbd dimensional molecule: it admits the 2\nbd layering $\order{1}{V} \cp{2} \order{1}{V}$, where
	\begin{align*}
		\order{1}{V} & \eqdef
		\disk{1}{2} \cp{1} ((\globe{2} \celto (\globe{2} \cp{1} \globe{2})) \cp{0} \globe{2}), \\
		\order{2}{V} & \eqdef
		((\disk{1}{2} \cp{1} (\globe{2} \cp{0} \globe{2})) \celto \disk{1}{2}) \cp{1} (\globe{2} \cp{0} \thearrow{}).
	\end{align*}
	However, the 2\nbd flow graph $\flow{2}{U}$ is
\[\begin{tikzcd}
	{{\scriptstyle (3, 1)}\;\bullet} && {{\scriptstyle (3, 0)}\;\bullet} \\
	& {{\scriptstyle (3, 2)}\;\bullet}
	\arrow[from=1-1, to=1-3]
	\arrow[from=1-1, to=2-2]
	\arrow[from=2-2, to=1-3]
\end{tikzcd}\]
	and its induced subgraph $\flow{2}{V}$ on $\set{(3, 0), (3, 1)}$ is not path-induced.
	The same ``non-convexity'' can also be detected in $\graph{U}$:
	\[
		\input{img/nonsubmol_graph.tex}
	\]
	It follows from Proposition \ref{prop:round_submolecule_flow_path_induced} that $V$ is not a submolecule of $U$.
	We will see, by Theorem \ref{thm:round_submolecule_dim2}, that dimension 3 is the lowest dimension in which this can happen.
\end{exm}

\clearpage 
\thispagestyle{empty}

%% file: omegacats.tex
\chapter{Diagrams in strict \omegatit-categories} \label{chap:omegacats}
\thispagestyle{firstpage}

\begin{guide}
	From the very start, we have informally thought of oriented graded posets as shapes of $n$\nbd categorical diagrams.
	It is time to make this connection precise, by defining a functor from the category of oriented graded posets to the category of \cemph{strict $\omega$\nbd categories} and \cemph{strict functors}.

	We have actually already done all the hard work in Chapter \ref{chap:molecules}, so it is now quite straightforward to prove that isomorphisms classes of molecules form a strict $\omega$\nbd category with pasting at the $k$\nbd boundary as $k$\nbd composition (Theorem 
	\ref{thm:molec_omega_category}).
	Having defined pastings universally as \emph{pushouts} now pays off, because those pushouts are still pushouts in the slice of $\ogpos$ over each oriented graded poset $P$, so it is almost trivial to show that \cemph{molecules over $P$} --- that is, isomorphism classes of morphisms $f\colon U \to P$, where $U$ is a molecule, in the slice of $\ogpos$ over $P$ --- also form a strict $\omega$\nbd category, naturally in $P$ (Theorem \ref{thm:molecin_ogp_omega_category}).

	This produces a strict $\omega$\nbd category $\molecin{P}$ associated to each oriented graded poset $P$, but functors out of it still lack one of the properties that we expected from combinatorial diagrams: namely, that they are uniquely determined by a labelling of the elements of $P$ with cells of the codomain.
	We finally restrict our attention to \cemph{regular directed complexes}, which are the oriented graded posets that are ``locally atoms'', that is, whose lower sets at each element are atoms.
	We prove that, when $P$ is a regular directed complex, a functor out of $\molecin{P}$ is uniquely determined by a labelling of $P$ (Proposition \ref{prop:combinatorial_diagrams_suffice}), so it deserves to be called a combinatorial diagram.
	This relies on further strengthening our earlier results on the rigidity of molecules, with the proof that \emph{every morphism} between atoms of the same dimension is, in fact, an isomorphism (Theorem 
	\ref{thm:morphisms_of_atoms_are_injective}).
	
	To conclude, when $U$ is a molecule, the fibrational definition of $\molecin{U}$ makes the \cemph{pasting theorem} almost a triviality: the isomorphism class of the identity of $U$ is a ``greatest'', that is, \emph{terminal} molecule over $U$ (Corollary \ref{cor:atom_only_one_topdim_cell}).
	Functors out of $\molecin{U}$ thus deserve to be called \cemph{pasting diagrams}, whose composite is the image of the greatest cell.
\end{guide}


\section{Fundamentals of strict \omegatit-categories} \label{sec:strict_omegacat}

\begin{guide}
	In this section, we fix notations and give some basic coordinates of the theory of strict $\omega$\nbd categories.
	We refer the reader to \cite{ara2023polygraphs} for a recent, thorough reference on the subject.
\end{guide}

\begin{dfn}[Reflexive $\omega$-graph] \index{reflexive $\omega$-graph} \index{cell}
	A \emph{reflexive $\omega$\nbd graph} is a set $X$, whose elements are called \emph{cells}, together with, for all $n \in \mathbb{N}$, operators
	\[
		\bound{n}{-}, \bound{n}{+}\colon X \to X
	\]
	called \emph{input} and \emph{output $n$\nbd boundary}, satisfying the following axioms.
	\begin{enumerate}
		\item (\textit{Finite dimension}). For all $t \in X$, there exists $n \in \mathbb{N}$ such that
			\[
				\bound{n}{-}t = \bound{n}{+}t = t.
			\]
		\item (\textit{Globularity}). For all $t \in X$, $k, n \in \mathbb{N}$, and $\alpha, \beta \in \set{ +, - }$,
			\[
				\bound{k}{\alpha}(\bound{n}{\beta}t) = 
				\begin{cases}
					\bound{k}{\alpha}t & \text{if $k < n$,} \\
					\bound{n}{\beta}t & \text{if $k \geq n$.}
				\end{cases}
			\]
	\end{enumerate}
\end{dfn}

\begin{comm}
	We will use variables $t, u, \ldots$ for cells in a reflexive $\omega$\nbd graph or strict $\omega$\nbd category.
\end{comm}

\begin{dfn}[Dimension of a cell] \index{cell!dimension} \index{$\dim{t}$} \index{dimension!of a cell}
	Let $t$ be a cell in a reflexive $\omega$\nbd graph.
	The \emph{dimension of $t$} is the natural number
	\[
		\dim{t} \eqdef \min \set{ n \in \mathbb{N} \mid \bound{n}{-}t = \bound{n}{+}t = t }.
	\]
\end{dfn}

\begin{comm}
	Reflexive $\omega$\nbd graphs are also known as \emph{reflexive globular sets}.
	Here we are using the so-called ``single-set'' definition.
	Alternatively, reflexive $\omega$\nbd graphs can defined as presheaves on a category whose objects are the natural numbers, hence as certain sequences of sets connected by functions; see for example \cite[Section 10.1]{leinster2004higher}.
\end{comm}

\begin{lem} \label{lem:cell_boundaries_properties}
	Let $t$ be a cell in a reflexive $\omega$\nbd graph, $n \in \mathbb{N}$, and $\alpha \in \set{+, -}$. 
	Then
	\begin{enumerate}
		\item $\dim \bound{n}{\alpha} t \leq \min \set{n, \dim{t}}$,
		\item if $n \geq \dim t$, then $\bound{n}{\alpha}t = t$.
	\end{enumerate}
\end{lem}

\begin{proof}
	By globularity, for all $\beta \in \set{+, -}$ we have $\bound{n}{\beta}(\bound{n}{\alpha}t) = \bound{n}{\alpha}t$, so by definition $\dim \bound{n}{\alpha} t \leq n$.
	If $n \leq k \eqdef \dim{t}$, this implies that $\dim \bound{n}{\alpha} t \leq \dim{t}$.
	If $n \geq k$, we have $\bound{n}{\alpha}t = \bound{n}{\alpha}(\bound{k}{\beta}t) = \bound{k}{\beta}t = t$.
\end{proof}

\begin{dfn}[Morphism of reflexive $\omega$-graphs] \index{reflexive $\omega$-graph!morphism} \index{morphism!of reflexive $\omega$-graphs} \index{$\romegagph$}
	Let $X, Y$ be reflexive $\omega$\nbd graphs.
	A \emph{morphism} $f\colon X \to Y$ is a function of their underlying sets which, for all $t \in X$, $n \in \mathbb{N}$, and $\alpha \in \set{+, -}$, satisfies
	\[
		f(\bound{n}{\alpha}t) = \bound{n}{\alpha}f(t).
	\]
	Reflexive $\omega$\nbd graphs and their morphisms form a category $\romegagph$.
\end{dfn}

\begin{lem} \label{lem:morphisms_of_omega_gph_decrease_dimension}
	Let $f\colon X \to Y$ be a morphism of reflexive $\omega$\nbd graphs, $t \in X$.
	Then $\dim{f(t)} \leq \dim{t}$.
\end{lem}
\begin{proof}
	Let $n \eqdef \dim{t}$.
	Then, for all $\alpha \in \set{+, -}$, $\bound{n}{\alpha}f(t) = f(\bound{n}{\alpha}t) = f(t)$.
\end{proof}

\begin{dfn}[Skeleta of a reflexive $\omega$-graph] \index{skeleton!of a reflexive $\omega$-graph} \index{$\skel{n}{X}$}
	Let $X$ be a reflexive $\omega$\nbd graph, $n \geq -1$.
	The \emph{$n$\nbd skeleton of $X$} is the reflexive $\omega$\nbd graph with set of cells
	\[
		\skel{n}{X} \eqdef \set{ t \in X \mid \dim{t} \leq n }
	\]
	and the boundary operators of $X$ restricted to this subset.
\end{dfn}

\begin{rmk}
	The boundary operators of $X$, restricted to $\skel{n}{X}$, have image in $\skel{n}{X}$ by Lemma \ref{lem:cell_boundaries_properties}.
\end{rmk}

\begin{dfn}[Reflexive $n$-graph] \index{reflexive $\omega$-graph!$n$-graph} \index{$\rngph{n}$}
	Let $n \in \mathbb{N}$.
	A \emph{reflexive $n$\nbd graph} is a reflexive $\omega$\nbd graph which is equal to its $n$\nbd skeleton, that is, has no cells of dimension strictly greater than $n$.
	We write $\rngph{n}$ for the full subcategory of $\romegagph$ on the reflexive $n$\nbd graphs.
\end{dfn}

\begin{dfn}[Composable pair of cells] \index{cell!composable pair} \index{$X \cpable{k} X$}
	Let $t, u$ be a pair of cells in a reflexive $\omega$\nbd graph, $k \in \mathbb{N}$.
	We say that $t$ and $u$ are \emph{$k$\nbd composable} if $\bound{k}{+}t = \bound{k}{-}u$.
	We write
	\[
		X \cpable{k} X \eqdef \set{ (t, u) \in X \times X \mid \bound{k}{+}t = \bound{k}{-}u }.
	\]
	for the set of $k$\nbd composable pairs of cells in $X$.
\end{dfn}

\begin{dfn}[Strict $\omega$-category] \index{strict $\omega$-category} \index{cell!composition} \index{$t \cp{k} u$}
	A \emph{strict $\omega$\nbd category} is a reflexive $\omega$\nbd graph $X$ together with, for all $k \in \mathbb{N}$, an operation
	\[
		- \cp{k} -\colon X \cpable{k} X \to X
	\]
	called \emph{$k$\nbd composition}, satisfying the following axioms.
	\begin{enumerate}
		\item (\textit{Compatibility with boundaries}).
			For all $k$\nbd composable pairs of cells $t, u$, all $n \in \mathbb{N}$, and $\alpha \in \set{+, -}$,
			\[
				\bound{n}{\alpha}(t \cp{k} u) =
				\begin{cases}
					\bound{n}{\alpha}t = \bound{n}{\alpha}u & \text{if $n < k$},\\
					\bound{k}{-}t & \text{if $n = k$, $\alpha = -$}, \\
					\bound{k}{+}u & \text{if $n = k$, $\alpha = +$}, \\
					\bound{n}{\alpha}t \cp{k} \bound{n}{\alpha}u & \text{if $n > k$}.
				\end{cases}
			\]
		\item (\textit{Associativity}).
			For all cells $t, u, v$ such that either side is defined,
			\[
				(t \cp{k} u) \cp{k} v = t \cp{k} (u \cp{k} v).
			\]
		\item (\textit{Unitality}).
			For all cells $t$,
			\[
				t \cp{k} \bound{k}{+}t = \bound{k}{-}t \cp{k} t = t.
			\]
		\item (\textit{Interchange}).
			For all cells $t, t', u, u'$ and $n > k$ such that the left-hand side is defined,
			\[
				(t \cp{n} t') \cp{k} (u \cp{n} u') = (t \cp{k} u) \cp{n} (t' \cp{k} u').
			\]
	\end{enumerate}
\end{dfn}

\begin{dfn}[Strict functor of strict $\omega$-categories] \index{strict $\omega$-category!strict functor|see {strict functor}} \index{$\omegacat$} \index{strict functor}
Let $X, Y$ be strict $\omega$-categories.
A \emph{strict functor} $f\colon X \to Y$ is a morphism of their underlying reflexive $\omega$\nbd graphs which, for all $k \in \mathbb{N}$ and $k$\nbd composable cells $t, u$ in $X$, satisfies
\[
	f(t \cp{k} u) = f(t) \cp{k} f(u).
\]
Strict $\omega$\nbd categories and strict functors form a category $\omegacat$.
\end{dfn}

\begin{prop} \label{prop:omegacat_complete_cocomplete}
	The category $\omegacat$ has all small limits and colimits.
\end{prop}
\begin{proof}
	See \cite[Corollary 14.2.5]{ara2023polygraphs}.
\end{proof}

\begin{lem} \label{lem:composition_with_lower_dim}
	Let $k \in \mathbb{N}$ and let $t, u$ be $k$\nbd composable cells in a strict $\omega$\nbd category.
	\begin{enumerate}
		\item If $\dim{u} \leq k$, then $t \cp{k} u = t$.
		\item If $\dim{t} \leq k$, then $t \cp{k} u = u$.
	\end{enumerate}
\end{lem}
\begin{proof}
	If $\dim{u} \leq k$, then by Lemma \ref{lem:cell_boundaries_properties} $u = \bound{k}{-}u$, which is equal to $\bound{k}{+}t$ because $t, u$ are $k$\nbd composable.
	Then $t \cp{k} u = t \cp{k} \bound{k}{+}t = t$ by unitality.
	The case $\dim{t} \leq k$ is analogous.
\end{proof}

\begin{lem} \label{lem:composition_does_not_increase_dim}
	Let $k \in \mathbb{N}$ and let $t, u$ be $k$\nbd composable cells in a strict $\omega$\nbd category.
	Then $\dim{t \cp{k} u} \leq \max \set{\dim{t}, \dim{u}}$.
\end{lem}
\begin{proof}
	Let $n \eqdef \max \set{\dim{t}, \dim{u}}$.
	If $n \leq k$, by Lemma \ref{lem:composition_with_lower_dim} we have $t \cp{k} u = t = u$, so $\dim{t \cp{k} u} = n$.
	If $n > k$, then for all $\alpha \in \set{+, -}$
	\[
		t \cp{k} u = \bound{n}{\alpha}t \cp{k} \bound{n}{\alpha}u = \bound{n}{\alpha}(t \cp{k} u)
	\]
	hence $\dim{t \cp{k} u} \leq n$.
\end{proof}

\begin{dfn}[Skeleta of a strict $\omega$-category] \index{skeleton!of a strict $\omega$-category} \index{$\skel{n}{X}$}
	Let $X$ be a strict $\omega$\nbd category, $n \geq -1$.
	The \emph{$n$\nbd skeleton of $X$} is the strict $\omega$\nbd category whose underlying reflexive $\omega$\nbd graph is $\skel{n}{X}$, together with the restriction of the composition operations of $X$ to $\skel{n}{X}$.
\end{dfn}

\begin{rmk}
	The $k$\nbd composition operation of $X$, restricted to pairs in $\skel{n}{X} \cpable{k} \skel{n}{X}$, has image in $\skel{n}{X}$ by Lemma \ref{lem:composition_does_not_increase_dim}.
\end{rmk}

\begin{rmk}
	Given a strict $\omega$\nbd category $X$, there is a natural cone of injective strict functors
\[\begin{tikzcd}
	{\varnothing \equiv\skel{-1}{X}} & {\skel{0}{X}} & {\skel{1}{X}} & \ldots & {\skel{n}{X}} & \ldots \\
	&&&&& X
	\arrow[hook, from=1-1, to=1-2]
	\arrow[hook, from=1-2, to=1-3]
	\arrow[hook, from=1-3, to=1-4]
	\arrow[hook, from=1-4, to=1-5]
	\arrow[hook, from=1-5, to=1-6]
	\arrow[hook, "{\imath_n}", from=1-5, to=2-6]
	\arrow[hook, "{\imath_1}", curve={height=6pt}, from=1-3, to=2-6]
	\arrow[hook, "{\imath_0}", curve={height=12pt}, from=1-2, to=2-6]
	\arrow[hook, "{\imath_{-1}}"', curve={height=18pt}, from=1-1, to=2-6]
\end{tikzcd}\]
	indexed by $n \geq -1$, exhibiting $X$ as the colimit of its skeleta in $\omegacat$.
\end{rmk}

\begin{dfn}[Strict $n$-category] \index{strict $\omega$-category!$n$-category} \index{$\ncat{n}$}
	Let $n \in \mathbb{N}$.
	A \emph{strict $n$\nbd category} is a strict $\omega$\nbd category which is equal to its $n$\nbd skeleton.
	We write $\ncat{n}$ for the full subcategory of $\omegacat$ on the strict $n$\nbd categories.
\end{dfn}

\begin{dfn}[Span of a set of cells] \index{strict $\omega$-category!span of a set of cells} 
	Let $X$ be a strict $\omega$\nbd category and $\gener{S}$ a set of cells in $X$.
	The \emph{span of $\gener{S}$} is the set of cells $\cspan{\gener{S}}$ inductively generated by the following clauses.
	\begin{enumerate}
		\item If $t \in \gener{S}$, then $t \in \cspan{\gener{S}}$.
		\item For all $k \in \mathbb{N}$, if $t, u \in \cspan{\gener{S}}$ are $k$\nbd composable, then $t \cp{k} u \in \cspan{\gener{S}}$.
	\end{enumerate}
\end{dfn}

\begin{dfn}[Generating set] \index{strict $\omega$-category!generating set}
	Let $X$ be a strict $\omega$\nbd category.
	A \emph{generating set for $X$} is a set $\gener{S}$ of cells such that $\cspan{\gener{S}}$ contains every cell in $X$.
\end{dfn}

\begin{lem} \label{lem:functors_equal_on_generating_set}
	Let $f,g\colon X \to Y$ be strict functors and let $\gener{S}$ be a generating set for $X$.
	If $f(t) = g(t)$ for all $t \in \gener{S}$, then $f = g$.
\end{lem}
\begin{proof}
	Let $t$ be a cell in $X$.
	Then $t \in \cspan{\gener{S}}$ because $\gener{S}$ is a generating set.
	We can thus proceed by structural induction on $\cspan{\gener{S}}$.
	If $t \in \gener{S}$, then $f(t) = g(t)$ by assumption.
	Otherwise, $t = t_1 \cp{k} t_2$ for some $k$\nbd composable pair $t_1$, $t_2$ such that, by the inductive hypothesis, $f(t_1) = g(t_1)$ and $f(t_2) = g(t_2)$.
	Then $f(t) = f(t_1) \cp{k} f(t_2) = g(t_1) \cp{k} g(t_2) = g(t)$ by the definition of strict functor.
\end{proof}

\begin{dfn}[Basis of a strict $\omega$-category] \index{strict $\omega$-category!basis} \index{basis!of a strict $\omega$-category}
	Let $X$ be a strict $\omega$\nbd category.
	A \emph{basis for $X$} is a generating set $\gener{S}$ which is minimal, in the sense that, if $\gener{T} \subseteq \gener{S}$ is also a generating set for $X$, then $\gener{T} = \gener{S}$.
\end{dfn}


\section{The \omegatit-category of molecules} \label{sec:omegacat_molecules}

\begin{guide}
	In this section, we prove that isomorphism classes of molecules form a strict $\omega$\nbd category, and derive that isomorphism classes of molecules \emph{over an oriented graded poset $P$} also do, naturally in $P$ (Proposition \ref{prop:molecin_is_a_functor}).
\end{guide}

\begin{dfn}[Isomorphism classes of molecules] \index{$\molec$} \index{$\atom$}
	For each oriented graded poset $P$, let $\isocl{P}$ denote its isomorphism class in $\ogpos$.
	We let
	\begin{align*}
		\molec & \eqdef \set{ \isocl{U} \mid \text{$U$ is a molecule} }, \\
		\atom & \eqdef \set{ \isocl{U} \mid \text{$U$ is an atom} } \subset \molec.
	\end{align*}
\end{dfn}

\begin{comm}
	Alternatively, $\molec$ can be defined as a set containing, for each isomorphism class of molecules, a chosen representative.

	In \cite{hadzihasanovic2023data} we described an algorithmic method to compute a standard representative for each isomorphism class of molecules.
\end{comm}

\begin{lem} \label{lem:molec_omega_graph}
	For all $n \in \mathbb{N}$ and $\alpha \in \set{+, -}$, let
	\begin{align*}
		\bound{n}{\alpha}\colon \molec & \to \molec, \\
		\isocl{U} & \mapsto \isocl{\bound{n}{\alpha}U}.
	\end{align*}
	Then $\molec$ together with these boundary operators is a reflexive $\omega$\nbd graph.
\end{lem}
\begin{proof}
	By Lemma \ref{lem:molecules_are_globular}, $\bound{n}{\alpha}U$ is a molecule, and by Corollary \ref{cor:inclusions_preserve_boundaries} if $U$ and $V$ are isomorphic, so are $\bound{n}{\alpha}U$ and $\bound{n}{\alpha}V$.
	This proves that $\bound{n}{\alpha}$ is well-defined as a function.
	The finite dimension axiom is implied by Corollary \ref{cor:dimension_from_boundary}.
	Finally, the globularity axiom is a consequence of Lemma \ref{lem:molecules_are_globular} together with Lemma \ref{lem:boundary_inclusion}.
\end{proof}

\begin{lem} \label{lem:dim_molecule_as_cell}
	Let $U$ be a molecule.
	The dimension of $\isocl{U}$ as a cell in $\molec$ is equal to the dimension of $U$ as an oriented graded poset.
\end{lem}
\begin{proof}
	Immediate from Corollary \ref{cor:dimension_from_boundary}.
\end{proof}

\begin{thm} \label{thm:molec_omega_category}
	For each $k \in \mathbb{N}$, let
	\begin{align*}
		- \cp{k} -\colon \molec \cpable{k} \molec & \to \molec, \\
		\isocl{U}, \isocl{V} & \mapsto \isocl{U \cp{k} V}.
	\end{align*}
	Then $\molec$ together with these composition operations is a strict $\omega$\nbd category.
\end{thm}
\begin{proof}
	First of all, $\isocl{U}$ and $\isocl{V}$ are $k$\nbd composable if and only if $\isocl{\bound{k}{+}U} = \isocl{\bound{k}{-}V}$, in which case there is a unique isomorphism between $\bound{k}{+}U$ and $\bound{k}{-}V$, hence $U \cp{k} V$ is defined.
	That this is independent of representatives is a consequence of Corollary $\ref{cor:inclusions_preserve_boundaries}$ and the universal property of pushouts.
	This proves that the composition operations are well-defined as functions.

	Then, compatibility with boundaries is a combined consequence of
	Lemma \ref{lem:pasting_lower_boundary}, Lemma \ref{lem:pasting_top_boundary}, and Lemma \ref{lem:pasting_higher_boundary}.
	Finally, associativity follows from Proposition \ref{prop:associativity_of_pasting}, unitality from Proposition \ref{prop:unitality_of_pasting}, and interchange from Proposition \ref{prop:interchange_of_pasting}.
\end{proof}

\begin{dfn}[The $\omega$-category of molecules] \index{strict $\omega$-category!of molecules}
	The \emph{$\omega$\nbd category of molecules} is the strict $\omega$\nbd category described in Theorem \ref{thm:molec_omega_category}.
\end{dfn}

\begin{prop} \label{prop:atoms_basis_for_molecules}
	The set $\atom$ is a basis for $\molec$.
\end{prop}
\begin{proof}
	Let $\isocl{U} \in \molec$.
	We prove that $\isocl{U} \in \cspan{\atom}$ by induction on $\lydim{U}$.
	If $\lydim{U} = -1$, then $U$ is an atom, so $\isocl{U} \in \atom \subseteq \cspan{\atom}$.
	Else, $k \eqdef \lydim{U} \geq 0$, and $U$ admits a $k$\nbd layering $(\order{i}{U})_{i=1}^m$.
	Then 
	\[
		\isocl{U} = \isocl{\order{1}{U}} \cp{k} \ldots \cp{k} \isocl{\order{m}{U}},
	\]
	and by the inductive hypothesis $\isocl{\order{i}{U}} \in \cspan{\atom}$ for all $i \in \set{1, \ldots, m}$.
	It follows that $\isocl{U} \in \cspan{\atom}$, which proves that $\atom$ is a generating set.

	Suppose $\gener{S} \subseteq \atom$ is such that $\cspan{\gener{S}} = \molec$, and let $\isocl{U} \in \atom$.
	By assumption $\isocl{U} \in \cspan{\gener{S}}$; we will prove that $\isocl{U} \in \gener{S}$ by structural induction.
	If $\isocl{U} \in \gener{S}$, then we are done.
	Otherwise, $\isocl{U} = \isocl{U_1} \cp{k} \isocl{U_2} = \isocl{U_1 \cp{k} U_2}$ for some $k$\nbd composable pair $\isocl{U_1}, \isocl{U_2}$, and by the inductive hypothesis $\isocl{U_1}$, $\isocl{U_2} \in \gener{S}$.
	By Lemma \ref{lem:composition_with_lower_dim}, if $k \geq \min \set{\dim{U_1}, \dim{U_2}}$, then $\isocl{U} = \isocl{U_1}$ or $\isocl{U} = \isocl{U_2}$, and we are done.
	If $k < \min \set{\dim{U_1}, \dim{U_2}}$, by Lemma \ref{lem:layering_dimensions_pasting} $\lydim{U} \geq k$, contradicting the fact that $U$ is an atom.
	We conclude that $\isocl{U} \in \gener{S}$, hence $\gener{S} = \atom$, and $\atom$ is a basis.
\end{proof}

\begin{dfn}[Molecules over an oriented graded poset] \index{$\molecin{P}$} \index{$\atomin{P}$}
	For each morphism $f\colon U \to P$ of oriented graded posets, let $\isocl{f}$ denote its isomorphism class in the slice category $\slice{\ogpos}{P}$.
	Given an oriented graded poset $P$, we let
	\begin{align*}
		\molecin{P} & \eqdef \set{ \isocl{f\colon U \to P} \mid \text{$U$ is a molecule} }, \\
		\atomin{P} & \eqdef \set{ \isocl{f\colon U \to P} \mid \text{$U$ is an atom} } \subseteq \molecin{P}.
	\end{align*}
	We call these \emph{molecules over $P$} and \emph{atoms over $P$}, respectively.
\end{dfn}

\begin{rmk}
	By definition of what an isomorphism is in a slice category, if $\isocl{f\colon U \to P} = \isocl{g\colon V \to P}$, then $\isocl{U} = \isocl{V}$.
\end{rmk}

\begin{lem} \label{lem:fibered_omega_graph}
	Let $P$ be an oriented graded poset.
	For all $k \in \mathbb{N}$ and $\alpha \in \set{+, -}$, let
	\begin{align*}
		\bound{k}{\alpha}\colon \molecin{P} & \to \molecin{P}, \\
		\isocl{f\colon U \to P} & \mapsto \isocl{\restr{f}{\bound{k}{\alpha}U}\colon \bound{k}{\alpha}U \to P}.
	\end{align*}
	Then $\molecin{P}$ together with these boundary operators is a reflexive $\omega$\nbd graph.
	In particular, if $n \eqdef \dim{P} < \infty$, then it is a reflexive $n$\nbd graph.
\end{lem}
\begin{proof}
	The fact that $\molecin{P}$ is a reflexive $\omega$\nbd graph is a straightforward fibred variant of Lemma \ref{lem:molec_omega_graph}.

	Suppose that $\dim{U} = n < \infty$.
	Then if a morphism $U \to X$ exists, necessarily $\dim{U} \leq n$.
	It follows from Lemma \ref{lem:dim_molecule_as_cell} that $\molecin{P}$ has no cells of dimension higher than $n$, so it is a reflexive $n$\nbd graph.
\end{proof}

\begin{dfn}[Pasting of molecules over an oriented graded poset]
	Let $k \in \mathbb{N}$ and suppose that $[f\colon U \to P]$, $[g\colon V \to P]$ are $k$\nbd composable molecules over $P$.
	The fact that $\bound{k}{+}[f] = \bound{k}{-}[g]$ translates to the existence of a diagram
	\[
	\begin{tikzcd}
		\bound{k}{+}U & \bound{k}{-}V & V \\
		U && P
		\arrow["\varphi", hook, from=1-1, to=1-2]
		\arrow[hook, from=1-2, to=1-3]
		\arrow[hook', from=1-1, to=2-1]
		\arrow["f", from=2-1, to=2-3]
		\arrow["g", from=1-3, to=2-3]
		\end{tikzcd}
	\]
	in $\ogpos$ for a necessarily unique isomorphism $\varphi$ of molecules.
	We define $f \cp{k} g\colon U \cp{k} V \to P$ to be the unique morphism induced by the universal property of the pushout diagram defining $U \cp{k} V$.
\end{dfn}

\begin{thm} \label{thm:molecin_ogp_omega_category}
	Let $P$ be an oriented graded poset and, for each $k \in \mathbb{N}$, let
	\begin{align*}
		- \cp{k} -\colon \molecin{P} \cpable{k} \molecin{P} & \to \molecin{P}, \\
		\isocl{f\colon U \to P}, \isocl{g\colon V \to P} & \mapsto \isocl{f \cp{k} g\colon U \cp{k} V \to P}.
	\end{align*}
	Then 
	\begin{enumerate}
		\item $\molecin{P}$ together with these composition operations is a strict $\omega$\nbd category,
		\item the set $\atomin{P}$ is a basis for $\molecin{P}$.
	\end{enumerate}
	In particular, if $n \eqdef \dim{P} < \infty$, then $\molecin{P}$ is a strict $n$\nbd category.
\end{thm}
\begin{proof}
	A simple fibred variant of Theorem \ref{thm:molec_omega_category} and Proposition \ref{prop:atoms_basis_for_molecules}.
\end{proof}

\begin{prop} \label{prop:molecin_is_a_functor} \index{strict functor!induced by a morphism}
	Let $f\colon P \to Q$ be a morphism of oriented graded posets.
	Then
	\begin{align*}
		\molecin{f}\colon \molecin{P} & \to \molecin{Q}, \\
		\isocl{g\colon U \to P} & \mapsto \isocl{f \after g\colon U \to Q}
	\end{align*}
	is a strict functor of strict $\omega$\nbd categories.
	This assignment determines a functor $\molecin{-}\colon \ogpos \to \omegacat$.
\end{prop}
\begin{proof}
	First of all, $\molecin{f}$ is well-defined as a function, for any isomorphism $g \iso g'$ in $\slice{\ogpos}{P}$ induces an isomorphism $f\circ g \iso f \circ g'$ in $\slice{\ogpos}{Q}$.
	Moreover, $\molecin{f}$ preserves boundary operators, since for all $n \in \mathbb{N}$ and $\alpha \in \set{ +, - }$
	\[
		f\circ (\restr{g}{\bound{n}{\alpha}U}) = \restr{(f \circ g)}{\bound{n}{\alpha}U}.
	\]
	Finally, given $k \in \mathbb{N}$ and a $k$\nbd composable pair $[g\colon U \to P]$ and $[h\colon V \to P]$,
	\[
		f \circ (g \cp{k} h) \quad \text{and} \quad (f \circ g) \cp{k} (f \circ h)
	\]
	both satisfy the property of universal morphisms from the pushout square defining $U \cp{k} V$, hence they are equal.
	This proves that $\molecin{f}$ preserves $k$\nbd composition operations, hence it is a strict functor.
	Functoriality of the assignment $f \mapsto \molecin{f}$ is straightforward.
\end{proof}

\begin{prop} \label{prop:skeleta_of_molecin}
	Let $P$ be an oriented graded poset, $n \geq -1$.
	Then $\skel{n}{\molecin{P}}$ is naturally isomorphic to $\molecin{\skel{n}{P}}$.
\end{prop}
\begin{proof}
	Let $\isocl{f\colon U \to P}$ be a cell in $\molecin{P}$.
	Since $f$ is dimension-preserving, $\dim{\isocl{f\colon U \to P}} \leq n$ if and only if $\dim{U} \leq n$ if and only if $f$ factors through the inclusion $\skel{n}{P} \incl P$.
	Naturality is straightforward.
\end{proof}


\section{Regular directed complexes} \label{sec:rdcpx}

\begin{guide}
	In this section, we finally give the definition of \emph{regular directed complex}, which has been anticipated since the beginning of the book.
	We prove that all regular directed complexes, augmented with a positive least element, are oriented thin (Proposition \ref{prop:oriented_diamond_rdc}).
	
	We then embark on the proof that all morphisms of regular directed complexes are local embeddings (Corollary \ref{cor:morphisms_of_rdcpx_are_local_isomorphisms}).
	We derive this from Theorem 
	\ref{thm:morphisms_of_atoms_are_injective}, which in conjunction with Proposition 
	\ref{prop:molecule_unique_automorphism} implies that there is at most one morphism between any pair of atoms of the same dimension.
	
	Having justified it beforehand, we conclude the chapter with a formal definition of \emph{diagram} and \emph{pasting diagram} in a strict $\omega$\nbd category based on the theory of regular directed complexes.
\end{guide}

\begin{dfn}[Regular directed complex] \index{regular directed complex}
A \emph{regular directed complex} is an oriented graded poset $P$ with the property that, for all $x \in P$, the closed subset $\clset{x}$ is an atom.
\end{dfn}

\begin{rmk}
By Lemma \ref{lem:all_downsets_are_atoms}, every molecule is a regular directed complex.
\end{rmk}

\begin{lem} \label{lem:if_local_embeds_into_rdcpx_then_is_rdcpx}
	Let $f\colon P \to Q$ be a local embedding of oriented graded posets.
	If $Q$ is a regular directed complex, then so is $P$.
\end{lem}
\begin{proof}
	For all $x \in P$, we have an isomorphism between $\clset{x}$ and $\clset{f(x)}$.
	Since the latter is an atom, so is the former.
\end{proof}

\begin{prop} \label{prop:oriented_diamond_rdc}
Let $P$ be a regular directed complex.
Then $\augm{P}$ is an oriented thin graded poset.
\end{prop}
\begin{proof}
	Let $x', y' \in \augm{P}$ with $\codim{x'}{y'} = 2$.
	Suppose first that $x' = \augm{x}$ for some $x \in P$.
	Then also $y' = \augm{y}$ for a unique $y \in P$, and the interval $[x', y']$ in $\augm{P}$ is isomorphic to the interval $[x, y]$ in $P$.
	We know that there exists at least one element $z_1$ such that $z_1 \in \faces{}{} y \cap \cofaces{}{} x$.
	Let $\alpha, \beta \in \set{ +, - }$ such that $z_1 \in \faces{}{\alpha}y \cap \cofaces{}{\beta}x$.
	Then $x \in \bound{}{\alpha}y$, which is a round molecule, and $\codim{x}{\bound{}{\alpha}y} = 1$.
	Since $z_1 \in \cofaces{}{\beta} x$, by Corollary \ref{cor:codimension_1_elements}, there are only two options: either
	\begin{itemize}
		\item $x \in \faces{}{\beta}(\bound{}{\alpha} y)$, in which case $\faces{}{\alpha} y \cap \cofaces{}{\beta} x = \set{z_1}$ and $\faces{}{\alpha} y \cap \cofaces{}{-\beta} x = \varnothing$, or
		\item $x \notin \faces{}{}(\bound{}{\alpha} y)$, in which case there exists a single other element $z_2$ such that $\faces{}{\alpha} y \cap \cofaces{}{\beta} x = \set{z_1}$ and $\faces{}{\alpha} y \cap \cofaces{}{-\beta} x = \set{z_2}$.
	\end{itemize}
	In the first case, by globularity of $\clset{y}$, we have $x \in \faces{}{\beta}(\bound{}{-\alpha} y)$, and by Corollary \ref{cor:codimension_1_elements} applied to $\bound{}{-\alpha} y$, we have $\cofaces{}{}x \cap \faces{}{-\alpha} y = \cofaces{}{\beta} x \cap \faces{}{-\alpha} y = \set{z_2}$ for a single element $z_2$.
	In the second case, since $\clset{y}$ is round, $x \notin \faces{}{}(\bound{}{\alpha}y)$ implies that $x \in \bound{}{\alpha}y \setminus \bound{}{-\alpha}y$, so no other elements of $\faces{}{}y$ covers $x$.

	We conclude that the interval $[x, y]$ is of the form
\[\begin{tikzcd}[column sep=small]
	& y \\
	{z_1} && {z_2} \\
	& x
	\arrow["{\alpha}"', from=1-2, to=2-1]
	\arrow["{-\alpha}", from=1-2, to=2-3]
	\arrow["{\beta}"', from=2-1, to=3-2]
	\arrow["{\beta}", from=2-3, to=3-2]
\end{tikzcd}
	\quad \text{or} \quad
\begin{tikzcd}[column sep=small]
	& y \\
	{z_1} && {z_2} \\
	& x
	\arrow["{\alpha}"', from=1-2, to=2-1]
	\arrow["{\alpha}", from=1-2, to=2-3]
	\arrow["{\beta}"', from=2-1, to=3-2]
	\arrow["{-\beta}", from=2-3, to=3-2]
\end{tikzcd}\]
	in the two cases, respectively.

	Now, suppose that $x' = \bot$.
	Then $y' = \augm{y}$ for some $y \in \grade{1}{P}$.
	By Lemma \ref{lem:only_1_molecules}, $\clset{y}$ is isomorphic to $\thearrow{}$, so there exists a unique pair $z_1, z_2 \in \grade{0}{P}$ such that $\faces{}{-}y = \set{z_1}$ and $\faces{}{+}y = \set{z_2}$.
	It follows that $[\bot, y']$ is of the form
\[\begin{tikzcd}[column sep=small]
	& {y'} \\
	{\augm{(z_1)}} && {\augm{(z_2)}} \\
	& {\bot}
	\arrow["-"', from=1-2, to=2-1]
	\arrow["+", from=1-2, to=2-3]
	\arrow["+"', from=2-1, to=3-2]
	\arrow["+", from=2-3, to=3-2]
\end{tikzcd}\]
	and we conclude.
\end{proof}

\begin{lem} \label{lem:morphisms_of_atoms_injective_codim_2}
	Let $f\colon U \to V$ be a morphism of atoms, $n \eqdef \dim{U} = \dim{V}$.
	For all $\alpha \in \set{ +, - }$,
	\begin{enumerate}
		\item $f$ is injective on $\grade{n-1}{U} = \faces{}{}U$ and on $\grade{n-2}{U}$,
		\item $f(\bound{}{\alpha}U) = \bound{}{\alpha}V$ and $f(\bound{n-2}{\alpha}U) = \bound{n-2}{\alpha}V$.
	\end{enumerate}
\end{lem}
\begin{proof}
	Let $\top_U, \top_V$ be the greatest elements of $U$ and $V$, respectively.
	Because $f$ is dimension-preserving and $\top_U$ and $\top_V$ are the only $n$\nbd dimensional elements in $U$ and $V$, respectively, $f(\top_U) = \top_V$.
	Moreover, $f$ is injective on $\faces{}{} \top_U = \grade{n-1}{U}$, so
	\[
		f(\bound{}{\alpha}U) = f(\clos{\faces{}{\alpha}\top_U}) = \clos{\faces{}{\alpha}f(\top_U)} = \bound{}{\alpha}V.
	\]

	Let $x, x' \in \grade{n-2}{U}$, and suppose $y \eqdef f(x) = f(x')$.
	By Proposition \ref{prop:oriented_diamond_rdc} applied to the intervals $[x, \top_U]$ and $[x', \top_U]$, we have 
	\[
		\cofaces{}{}x = \set{x_1, x_2}, \quad \cofaces{}{}x' = \set{x'_1, x'_2} 
	\]
	for some $x_1 \neq x_2, x'_1 \neq x'_2 \in \grade{n-1}{U}$, while applying it to $[y, \top_V]$ we have 		
	\[
		\cofaces{}{}y = \set{y_1, y_2}
	\]
	for some $y_1 \neq y_2 \in \grade{n-1}{V}$.
	Then $f(x_1), f(x_2), f(x'_1), f(x'_2) \in \set{y_1, y_2}$; since $f$ is injective on $\grade{n-1}{U}$, necessarily $\set{x_1, x_2} = \set{x'_1, x'_2}$.
	It follows that $x, x' \in \faces{}{}x_1$, and since $f$ is injective on $\faces{}{}x_1$, from $f(x) = f(x')$ we conclude that $x = x'$.
	This proves that $f$ is injective on $\grade{n-2}{U}$.

	As a consequence, for all $x \in \grade{n-2}{U}$ and all $\alpha \in \set{+, -}$, $f$ induces a bijection between $\cofaces{}{\alpha}x$ and $\cofaces{}{\alpha}f(x)$.
	Then $x \in \faces{n-2}{\alpha}U$ if and only if $\cofaces{}{}x = \cofaces{}{\alpha}x$ if and only if $\cofaces{}{}f(x) = \cofaces{}{\alpha}f(x)$ if and only if $f(x) \in \faces{n-2}{\alpha}V$.
	Since $\bound{n-2}{\alpha}U$ and $\bound{n-2}{\alpha}V$ are round, hence pure, this suffices to conclude that $f(\bound{n-2}{\alpha}U) = \bound{n-2}{\alpha}V$.
\end{proof}

\begin{lem} \label{lem:morphisms_of_atoms_preserve_boundaries}
	Let $f\colon U \to V$ be a morphism of atoms, $n \eqdef \dim{U} = \dim{V}$.
	For all $\alpha \in \set{ +, - }$ and $k < n$,
	\begin{enumerate}
		\item $f$ is injective on $\grade{k}{(\bound{k}{}U)} = \faces{k}{}U$ and on $\grade{k-1}{(\bound{k}{}U)}$,
		\item $f(\bound{k}{\alpha}U) = \bound{k}{\alpha}V$.
	\end{enumerate}
\end{lem}
\begin{proof}
	We proceed by induction on $(n - 1) - k$.
	For each $k$, we will prove that, for all $\alpha \in \set{+, -}$,
	\begin{enumerate}
		\item $f$ is injective on $\grade{k}{(\bound{k}{}U)} = \faces{k}{}U$ and on $\grade{k-1}{(\bound{k}{}U)}$,
		\item $f(\bound{k}{\alpha}U) = \bound{k}{\alpha}V$ and $f(\bound{k-1}{\alpha}U) = \bound{k-1}{\alpha}V$.
	\end{enumerate}
	These also imply that $f$ determines a bijection between $\faces{k-1}{\alpha}U$ and $\faces{k-1}{\alpha}V$.

	For $k = n-1$, the result is given by Lemma \ref{lem:morphisms_of_atoms_injective_codim_2}, since in this case $\grade{n-1}{(\bound{n-1}{}U)} = \grade{n-1}{U}$ and $\grade{n-2}{(\bound{n-1}{}U)} = \grade{n-2}{U}$.

	Suppose $k < n-1$, and let
	\[
		U' \eqdef \bound{k}{-}U \celto \bound{k}{+}U, \quad \quad V' \eqdef \bound{k}{-}V \celto \bound{k}{+}V;
	\]
	both of these are well-defined by roundness of $U$ and $V$.
	Moreover, we have unique isomorphisms
	\[
		\varphi\colon \bound{}{}U' \incliso \bound{k}{}U, \quad \quad
		\psi\colon \bound{k}{}V \incliso \bound{}{}V'
	\]
	restricting to the unique isomorphisms of molecules 
	\[
		\varphi^\alpha\colon \bound{}{\alpha}U' \incliso \bound{}{\alpha}U, \quad \quad
		\psi^\alpha\colon \bound{k}{\alpha}V \incliso \bound{}{\alpha}V'
	\]
	for each $\alpha \in \set{+, -}$.
	We define
	\[
		f'\colon U' \to V', \quad \quad 
		x \mapsto \begin{cases}
			\top & \text{if $x = \top$}, \\
			(\psi \after f \after \varphi)(x) & \text{if $x \in \bound{}{}U'$}.
		\end{cases}
	\]
	By the inductive hypothesis, $f(\bound{k}{\alpha}U) = \bound{k}{\alpha}V$, so this is well-defined as a function.
	Moreover, for all $x \in U'$ and $\alpha \in \set{+, -}$, we have
	\[
		\restr{f'}{\faces{}{\alpha}x} = \restr{(\psi \after f \after \varphi)}{\faces{}{\alpha}x}.
	\]
	Now, 
	\begin{itemize}
		\item if $x = \top$, then $\varphi$ maps $\faces{}{\alpha}x$ bijectively to $\faces{k}{\alpha}U$, which, by the inductive hypothesis, $f$ maps bijectively to $\faces{k}{\alpha}V$, which then $\psi$ maps bijectively to $\faces{}{\alpha}\top = \faces{}{\alpha}f(x)$;
		\item if $x \in \bound{}{}U'$, then $\varphi$ maps $\faces{}{\alpha}x$ bijectively to $\faces{}{\alpha}\varphi(x)$, which $f$ maps bijectively to $\faces{}{\alpha}(f \after \varphi)(x)$, which $\psi$ maps bijectively to $\faces{}{\alpha}(\psi \after f \after \varphi)(x)$.
	\end{itemize}
	This proves that $f'$ is a morphism of oriented graded posets, and in particular a morphism of $(k+1)$\nbd dimensional atoms.
	By Lemma \ref{lem:morphisms_of_atoms_injective_codim_2} applied to $f'$, 
	\begin{enumerate}
		\item $f'$ is injective on $\faces{}{}U'$ and on $\grade{k-1}{U'}$
		\item for all $\alpha \in \set{+, -}$, $f'(\bound{k-1}{\alpha}U') = \bound{k-1}{\alpha}V'$.
	\end{enumerate}
	Since $\varphi$ maps $\faces{}{}U'$ bijectively to $\faces{k}{}U$ and $\grade{k-1}{U'}$ to $\grade{k-1}{(\bound{k}{}U)}$, and $f'$ is equal to $\psi \after f \after \varphi$ on these sets, the first implies that $f$ is injective on $\faces{k}{}U$ and on $\grade{k-1}{(\bound{k}{}U)}$.
	Similarly, since $\varphi$ maps $\bound{k-1}{\alpha}U'$ bijectively to $\bound{k-1}{\alpha}U$ and $\psi$ maps $\bound{k-1}{\alpha}V'$ bijectively to $\bound{k-1}{\alpha}V$, the second implies that $f(\bound{k-1}{\alpha}U) = \bound{k-1}{\alpha}V$.
	This completes the inductive step and the proof.
\end{proof}

\begin{thm} \label{thm:morphisms_of_atoms_are_injective}
	Let $f\colon U \to V$ be a morphism of atoms of the same dimension.
	Then $f$ is an isomorphism.
\end{thm}
\begin{proof}
	Let $n \eqdef \dim{U} = \dim{V}$, and let $\top_U$, $\top_V$ be the greatest elements of $U$ and $V$, respectively.
	Because $\top_U$ and $\top_V$ are the only $n$\nbd dimensional elements of $U$ and $V$, respectively, and $f$ is dimension-preserving and closed,
	\[
		f(U) = f(\clset{\top_U}) = \clset{f(\top_U)} = \clset{\top_V} = V,
	\]
	so $f$ is surjective.
	It remains to prove that $f$ is injective.

	Let $x, x' \in U$ and suppose $y \eqdef f(x) = f(x')$.
	Because $f$ is dimension-preserving, $y = \top_V$ if and only if $x = x' = \top_U$, so we may assume $x, x' \in \bound{}{}U$ and $y \in \bound{}{}V$.
	Since $U$ is round, by Lemma \ref{lem:round_partition_into_interiors} there exist unique $k < n$ and $\alpha \in \set{ +, - }$ such that $x \in \inter{\bound{k}{\alpha}U}$.
	By Lemma \ref{lem:morphisms_of_atoms_preserve_boundaries}, $f$ restricts to 
	\[
		\restr{f}{\bound{k}{\alpha}U}\colon \bound{k}{\alpha}U \to \bound{k}{\alpha}V,
	\]
	which is a morphism of molecules of the same dimension.
	By Proposition \ref{prop:morphisms_of_molecules_preserve_interior}, $y = f(x) \in \inter{\bound{k}{\alpha}V}$.
	Applying the same reasoning to $x'$, we find unique $j < n$ and $\beta \in \set{ +, - }$ such that $x' \in \inter{\bound{j}{\beta}U}$, and deduce that $y \in \inter{\bound{j}{\beta}V}$.
	By Lemma \ref{lem:round_partition_into_interiors} applied to $V$, necessarily $j = k$ and $\alpha = \beta$.

	We have established that $x, x' \in \inter{\bound{k}{\alpha}U} \subseteq \bound{k}{}U$.
	Next, we proceed by induction on $\ell \eqdef k - \dim{x} = k - \dim{x'}$.
	If $\ell \in \set{0, 1}$, by Lemma \ref{lem:morphisms_of_atoms_preserve_boundaries} $f(x) = f(x')$ implies $x = x'$, so suppose $\ell > 1$.
	By Lemma \ref{lem:characterisation_of_interior} applied to the molecule $\bound{k}{\alpha}U$, there exist unique $x_+, x_-, x'_+, x'_- \in \grade{k}{(\bound{k}{\alpha}U)} = \faces{k}{\alpha}U$ such that
	\[
		x \in \inter{\bound{}{+}x_+} \cap \inter{\bound{}{-}x_-}, \quad \quad
		x' \in \inter{\bound{}{+}x'_+} \cap \inter{\bound{}{-}x'_-}.
	\]
	Reasoning as before with Lemma \ref{lem:morphisms_of_atoms_preserve_boundaries} and Proposition \ref{prop:morphisms_of_molecules_preserve_interior}, we deduce that
	\[
		y \in \inter{\bound{}{+}f(x_+)} \cap \inter{\bound{}{-}f(x_-)} \quad
		\text{but also} \quad
		y \in \inter{\bound{}{+}f(x'_+)} \cap \inter{\bound{}{-}f(x'_-)}.
	\]
	Applying Lemma \ref{lem:characterisation_of_interior} to $\bound{k}{\alpha}V$, we deduce that
	\[
		f(x_+) = f(x'_+), \quad \quad f(x_-) = f(x'_-),
	\]
	and since $f$ is injective on $\faces{k}{\alpha}U$ it follows that $x_+ = x'_+$ and $x_- = x'_-$.
	Then $x, x' \in \clset{x_+}$, and we may restrict ourselves to 
	\[
		\restr{f}{\clset{x_+}}\colon \clset{x_+} \to \clset{f(x_+)},
	\]
	also a morphism of atoms of the same dimension.
	Now $x, x' \in \inter{\bound{k-1}{+}x_+}$ and $k - 1 - \dim{x} < \ell$, so the inductive hypothesis applies, and we conclude.
\end{proof}

\begin{cor} \label{cor:atoms_over_regular_directed_complexes}
	Let $U$ be an atom, $P$ a regular directed complex, and $f\colon U \to P$ a morphism.
	Then $f$ is an inclusion.
\end{cor}
\begin{proof}
	Let $\top$ be the greatest element of $U$.
	By definition of regular directed complex, $f(U) = f(\clset{\top}) = \clset{f(\top)}$ is an atom, so by Theorem \ref{thm:morphisms_of_atoms_are_injective} $f$ restricts to an isomorphism with its image.
\end{proof}

\begin{cor} \label{cor:morphisms_of_rdcpx_are_local_isomorphisms}
	Let $f\colon P \to Q$ be a morphism of regular directed complexes.
	Then $f$ is a local embedding.
\end{cor}

\begin{cor} \label{cor:basis_of_omegacat_presented_by_rdcpx}
	Let $P$ be a regular directed complex.
	Then the set
	\[
		\set{ \isocl{ \clset{x} \incl P } \mid x \in P }
	\]
	is a basis for the $\omega$\nbd category $\molecin{P}$.
\end{cor}
\begin{proof}
	By Theorem \ref{thm:molecin_ogp_omega_category}, $\atomin{P}$ is a basis, and by Corollary \ref{cor:atoms_over_regular_directed_complexes} every morphism from an atom to $P$ is an inclusion, isomorphic to one of the form $\clset{x} \incl P$ for some $x \in P$.
\end{proof}

\begin{prop} \label{prop:molecule_over_atom_is_iso}
	Let $f\colon U \to V$ be a morphism of molecules of the same dimension.
	If $V$ is an atom, then $f$ is an isomorphism.
\end{prop}
\begin{proof}
	We proceed by induction on $\lydim{U}$.
	If $\lydim{U} = -1$, then $U$ is an atom, and the statement is the content of Theorem 
	\ref{thm:morphisms_of_atoms_are_injective}.

	Suppose $\lydim{U} = k \geq 0$.
	Then $U$ admits a $k$\nbd layering $(\order{i}{U})_{i=1}^m$.
	Let $i \in \set{1, \ldots, m}$ be such that $\dim{\order{i}{U}} = \dim{U}$, and write
	\begin{align*}
		W & \eqdef \order{1}{U} \cp{k} \ldots \cp{k} \order{i-1}{U} \cp{k} \bound{k}{-}\order{i}{U}, \\
		W' & \eqdef \bound{k}{+}\order{i}{U} \cp{k} \order{i+1}{U} \cp{k} \ldots \cp{k} \order{m}{U},
	\end{align*}
	so that we may identify $U$ with $W \cp{k} \order{i}{U} \cp{k} W'$.
	By the inductive hypothesis, $\restr{f}{\order{i}{U}}\colon \order{i}{U} \to V$ is an isomorphism, so by Corollary \ref{cor:inclusions_preserve_boundaries}
	\begin{align*}
		f(\bound{k}{-}\order{i}{U}) & = f(\bound{k}{+}W) = \bound{k}{-}V, \\
		f(\bound{k}{+}\order{i}{U}) & = f(\bound{k}{-}W') = \bound{k}{+}V.
	\end{align*}
	Suppose that $\dim{W} > k$.
	Then there exist $x \in \faces{k}{+}W = \faces{k}{-}\order{i}{U}$ and $y \in \cofaces{}{+}x \cap W$.
	It follows that $f(y) \in \cofaces{}{+}f(x)$, but $f(x) \in \faces{k}{-}V$, a contradiction.
	Then $\dim{W} \leq k$.
	Dually, we prove that $\dim{W'} \leq k$, so $i = m = 1$, contradicting $\lydim{U} = k$.
\end{proof}

\begin{cor} \label{cor:atom_only_one_topdim_cell}
	Let $U$ be an atom, $n \eqdef \dim{U}$.
	Then $\isocl{\idd{U}\colon U \to U}$ is the only $n$\nbd dimensional cell of $\molecin{U}$.
\end{cor}

\begin{dfn}[Diagram in a strict $\omega$-category] \index{diagram!in a strict $\omega$-category} \index{shape!of a diagram}
	Let $X$ be a strict $\omega$\nbd category and $P$ a regular directed complex.
	A \emph{diagram of shape $P$ in $X$} is a strict functor $d\colon \molecin{P} \to X$.
\end{dfn}

\begin{dfn}[Combinatorial diagram] \index{diagram!combinatorial} \index{$\ldiag{d}$}
	Let $X$ be a strict $\omega$\nbd category, $P$ a regular directed complex, and $d$ a diagram of shape $P$ in $X$.
	The \emph{combinatorial diagram} associated with $d$ is the function
	\begin{align*}
		\ldiag{d}\colon P &\to X, \\
		x & \mapsto d\isocl{\clset{x} \incl P},
	\end{align*}
	where $P$ and $X$ should be read as the underlying set of $P$ and the set of cells of $X$, respectively.
\end{dfn}

\begin{prop} \label{prop:combinatorial_diagrams_suffice}
	Let $X$ be a strict $\omega$\nbd category, $P$ a regular directed complex, and $d, d'$ two diagrams of shape $P$ in $X$.
	If $\ldiag{d} = \ldiag{d'}$, then $d = d'$.
\end{prop}
\begin{proof}
	By Corollary \ref{cor:basis_of_omegacat_presented_by_rdcpx}, if $\ldiag{d} = \ldiag{d'}$, then $d(t) = d'(t)$ for all cells $t$ in a basis of $\molecin{P}$, hence $d = d'$ by Lemma \ref{lem:functors_equal_on_generating_set}.
\end{proof}

\begin{dfn}[Pasting diagram] \index{diagram!pasting diagram}
	Let $X$ be a strict $\omega$\nbd category.
	A \emph{pasting diagram in $X$} is a diagram in $X$ whose shape is a molecule.
\end{dfn}

\begin{dfn}[Boundary of a pasting diagram] \index{diagram!boundary}
	Let $X$ be a strict $\omega$\nbd category, $d$ be a pasting diagram of shape $U$ in $X$, $n \in \mathbb{N}$, and $\alpha \in \set{+, -}$.
	We let $\bound{n}{\alpha}d \eqdef \restr{d}{\bound{n}{\alpha}U}$, a pasting diagram of shape $\bound{n}{\alpha}U$ in $X$.
	We also let $\bound{n}{}d \eqdef \restr{d}{\bound{n}{}U}$, a diagram of shape $\bound{n}{}U$, and omit $n$ when equal to $\dim{U} - 1$.
\end{dfn}

\begin{dfn}[Composite of a pasting diagram] \index{diagram!composite}
	Let $X$ be a strict $\omega$\nbd category and $d$ a pasting diagram of shape $U$ in $X$.
	The \emph{composite} of $d$ is the cell $d\isocl{\idd{U}\colon U \to U}$.
\end{dfn}

\begin{exm}[Commutative diagrams] \index{diagram!commutative} \label{exm:commutative_diagram}
	The commonly used informal notion of \emph{commutative diagram} in a category is largely subsumed by the notion of diagram in a strict $\omega$\nbd category.
	
	Indeed, a small category $\smcat{C}$ can be seen as a strict 1\nbd category.
	Given a ``commutative diagram in $\smcat{C}$'', we can add a 2\nbd cell, with either direction, between any two paths in its underlying graph that are supposed to ``commute''.
	If the oriented face poset $P$ of this diagram shape is a 2\nbd dimensional regular directed complex, then the commutative diagram is represented by a strict functor $\molecin{P} \to \smcat{C}$, that is, a diagram of shape $P$ in $\smcat{C}$.
	This is because any 2\nbd dimensional cell in $\molecin{P}$ has to be mapped to a lower\nbd dimensional cell in $\smcat{C}$, so the composites of its input and output 1\nbd boundaries must be equal in $\smcat{C}$, which is precisely the meaning of ``commutativity''.

	For example, a commutative square
	\[\begin{tikzcd}
	c && c' \\
	d && d'
	\arrow["g", from=1-1, to=1-3]
	\arrow["f", from=1-1, to=2-1]
	\arrow["f'", from=1-3, to=2-3]
	\arrow["h", from=2-1, to=2-3]
\end{tikzcd}\]
	in $\smcat{C}$ corresponds to a diagram whose shape is the oriented face poset of
	\[\begin{tikzcd}
	\bullet && \bullet \\
	\bullet && \bullet
	\arrow[from=1-1, to=1-3]
	\arrow[from=1-1, to=2-1]
	\arrow[from=1-3, to=2-3]
	\arrow[from=2-1, to=2-3]
	\arrow[shorten <=5pt, shorten >=5pt, Rightarrow, from=2-1, to=1-3]
\end{tikzcd}\]
	which is a 2\nbd dimensional atom.
	A \emph{non-necessarily commutative} square, on the other hand, is a diagram whose shape is the oriented face poset of
\[\begin{tikzcd}
	\bullet && \bullet \\
	\bullet && \bullet
	\arrow[from=1-1, to=1-3]
	\arrow[from=1-1, to=2-1]
	\arrow[from=1-3, to=2-3]
	\arrow[from=2-1, to=2-3]
\end{tikzcd}\]
	which is a 1\nbd dimensional regular directed complex.
\end{exm}

\clearpage
\thispagestyle{empty}

%% file: maps.tex
\chapter{Maps and comaps} \label{chap:maps}
\thispagestyle{firstpage}

\begin{guide}
	As should become more and more apparent, regular directed complexes are a very convenient class of \emph{objects} for the purposes of higher-categorical combinatorics.
	However, as seen at the end of last chapter, their \emph{morphisms} as oriented graded posets are very rigid: they are local embeddings.
	The purpose of this chapter is to explore more expressive notions of morphism.

	As objects at a crossroads of topology and higher category theory, regular directed complexes need to balance the pull of the two sides.
	Thus, we want morphisms between regular directed complexes $P$ and $Q$ to be interpreted as strict functors between $\molecin{P}$ and $\molecin{Q}$, but we do not want \emph{all} strict functors, since not all strict functors are sound for geometric realisation.
	Key to topological soundness, as we will see in Chapter \ref{chap:geometric}, is that morphisms have an underlying order-preserving map of posets, so the question becomes: \emph{what order-preserving maps induce strict functors between $\molecin{P}$ and $\molecin{Q}$?}

	We will answer this question twice: first \emph{covariantly}, and then \emph{contravariantly}.
	The answers determine the classes of \cemph{maps} and of \cemph{comaps} of regular directed complexes, respectively.

	Maps $p\colon P \to Q$ must ``push forward'' molecules over $P$ to molecules over $Q$, which are represented by local embeddings.
	However, this cannot be achieved by taking the \emph{direct image} of $p$ as an order-preserving map, since even the direct image of a molecule through a local embedding is not necessarily a molecule.
	Instead, we have to consider a different orthogonal factorisation system on order-preserving maps, whose right class is precisely the local embeddings, and whose left class is the class of \cemph{final maps}.
	This is the posetal version of the \emph{comprehensive factorisation system} of categories and functors.

	Injective maps of regular directed complexes coincide with their inclusions as oriented graded posets (Lemma \ref{lem:injective_maps_are_inclusions}), and more generally dimension-preserving maps of regular directed complexes coincide with the local embeddings and with their morphisms as oriented graded posets (Proposition 
	\ref{prop:characterisation_of_morphisms_among_maps}).
	What is new is that maps are allowed to \emph{decrease} the dimension of elements; in particular, the point is a terminal object in the category of regular directed complexes and maps, which it certainly was not in $\ogpos$.
	As we will see in Chapter \ref{chap:special}, this has the consequence that several \emph{shape categories} used in models of $(\infty, n)$\nbd categories, including both \emph{coface} and \emph{codegeneracy} maps, appear as subcategories of regular directed complexes and maps.

	Comaps of regular directed complexes have a more straightforward definition, as there is no restriction to ``pulling back'' molecules over $Q$ using ordinary inverse images.
	Dually to maps, comaps are allowed to \emph{increase} the dimension of elements; in fact, the only comaps that are also maps are isomorphisms (Proposition 
	\ref{prop:maps_comaps_trivial_intersection}).
	Duals of comaps can be seen as combinatorial \cemph{subdivisions} of regular cell complexes, preserving the dimension, boundary, and interior of each cell.

	As we hoped, maps and comaps of regular directed complexes admit a natural interpretation as strict functors of strict $\omega$\nbd categories (Theorem \ref{thm:molecin_is_a_functor_on_maps} and Theorem \ref{thm:molecin_is_a_functor_on_comaps}).
	These interpretations are \emph{pseudomonic} functors, which roughly means that a regular directed complex can be reconstructed up to isomorphism from the strict $\omega$\nbd category that it presents (Proposition \ref{prop:pseudomonicity_of_molecin_maps} and Proposition 
	\ref{prop:pseudomonicity_of_molecin_comaps}).
	Finally, both maps and comaps also admit functorial interpretations as homomorphisms of augmented chain complexes (Proposition \ref{prop:chain_complex_is_functorial_on_maps} and Proposition 
	\ref{prop:chain_complex_is_functorial_on_comaps}), with the same variance as the interpretations in strict $\omega$\nbd categories.
\end{guide}


\section{Pushforwards and pullbacks} \label{sec:pushforwards}

\begin{guide}
	In this section we provide some order-theoretic complements to Chapter \ref{chap:order}.
	We describe the orthogonal factorisation system on $\posclos$ whose left class is closed final maps, and right class is local embeddings, as well as its interaction with the usual epi--mono factorisation system, giving rise to a ternary factorisation system.
	Then, we define the \cemph{pushforward} of a local embedding along a closed order-preserving map, as well as the \cemph{pullback} of a local embedding along a not necessarily closed order-preserving map, which will play a role in the functorial action of maps and comaps, respectively, on molecules over a regular directed complex.
\end{guide}

\begin{dfn}[Final map] \index{map!final}
Let $f\colon P \to Q$ be an order-preserving map of posets.
We say that $f$ is \emph{final} if, for all $y \in Q$,
\begin{enumerate}
	\item there exists $x \in P$ such that $y \leq f(x)$, and
	\item for all $x, x' \in P$, if $y \leq f(x)$ and $y \leq f(x')$, then there exists a zig-zag
\[
	x \leq x_1 \geq x_2 \leq \ldots \geq x_{m-1} \leq x_m \geq x'
\]
	in $P$ such that $y \leq f(x_i)$ for all $i \in \set{1, \ldots, m}$.
\end{enumerate}
\end{dfn}

\begin{lem} \label{lem:final_closed_maps_are_surjective}
	Let $f\colon P \to Q$ be an order-preserving map of posets.
	If $f$ is closed and final, then $f$ is surjective.
\end{lem}
\begin{proof}
	Let $y \in Q$.
	Since $f$ is final, there exists $x \in P$ such that $y \leq f(x)$, so $y \in \clset{f(x)}$.
	Since $f$ is closed, $y \in f(\clset{x})$.
\end{proof}

\begin{prop} \label{prop:comprehensive_factorisation_system}
Let $f\colon P \to Q$ be a closed order-preserving map of posets.
Then $f$ factors as
\begin{enumerate}
	\item a closed final map $f_\clas{F}\colon P \to \pfw{f}P$,
	\item followed by a local embedding $f_\clas{L}\colon \pfw{f}P \to Q$.
\end{enumerate}
This factorisation is unique up to unique isomorphism.
\end{prop}
\begin{proof}
As shown in \cite{street1973comprehensive}, every functor of categories factors uniquely up to unique isomorphism as a \emph{final functor}, followed by a \emph{discrete fibration}.
This factorisation restricts to posets, identified with categories whose hom-sets have at most one element; the two classes then correspond precisely to final maps and local embeddings.

Explicitly, we can construct $\pfw{f}P$ as a poset of equivalences classes $\isocl{x, y}$ of pairs of $x \in P$ and $y \leq f(x)$.
The equivalence relation is generated by
\[
	\text{$\isocl{x, y} = \isocl{x', y}$ if $x \leq x'$},
\]
or, more explicitly, $\isocl{x, y} = \isocl{x', y}$ if and only if there exists a zig-zag
\[
	x \leq x_1 \geq x_2 \leq \ldots \geq x_{m-1} \leq x_m \geq x'
\]
in $P$ such that $y \leq f(x_i)$ for all $i \in \set{1, \ldots, m}$.
The partial order is defined by $\isocl{x, y} \leq \isocl{x', y'}$ if and only if $\isocl{x, y} = \isocl{x', y}$ and $y \leq y'$ in $Q$.
The maps $f_\clas{F}$ and $f_\clas{L}$ are defined by
\begin{align*}
	f_\clas{F}\colon P \to \pfw{f}P, \quad \quad  & x \mapsto \isocl{x, f(x)}, \\
	f_\clas{L}\colon \pfw{f}P \to Q, \quad \quad  & \isocl{x, y} \mapsto y.
\end{align*}
We only need to show that $f_\clas{F}$ is closed when $f$ is; as a local embedding, $f_u$ is always closed.
Suppose $\isocl{x, y} \leq f_\clas{F}(x') = \isocl{x', f(x')}$.
Then $\isocl{x, y} = \isocl{x', y}$, so $y \leq f(x')$.
Because $f$ is closed, there exists $y' \leq x'$ such that $f(y') = y$.
Then $f_\clas{F}(y') = \isocl{y', y} = \isocl{x', y} = \isocl{x, y}$, proving that $f_\clas{F}$ is closed.
\end{proof}

\begin{rmk}
	Notice that $\pfw{f}P$ is the colimit in $\posclos$ of the diagram of closed embeddings
	\[
		P \to \posclos, \quad \quad (x \leq y) \mapsto (\clset{f(x)} \incl \clset{f(y)}),
	\]
	exhibited by the cone mapping $\clset{f(x)}$ isomorphically onto $\clset{\isocl{x, f(x)}}$.
\end{rmk}

\begin{cor} \label{cor:ofs_final_local_embedding}
The classes of
\begin{enumerate}
	\item closed final maps,
	\item local embeddings
\end{enumerate}
form an orthogonal factorisation system on $\posclos$.
\end{cor}
\begin{proof}
It follows from the category-theoretic generalisation that both classes are closed under composition and contain all isomorphisms.
The statement then follows from Proposition \ref{prop:comprehensive_factorisation_system}.
\end{proof}

\begin{cor} \label{cor:ternary_factorisation_system}
The classes of
\begin{enumerate}
	\item closed final maps,
	\item surjective local embeddings,
	\item closed embeddings
\end{enumerate}
form a ternary factorisation system on $\posclos$.
\end{cor}
\begin{proof}
Follows from Corollary \ref{cor:em_ofs_on_posclos}, Corollary \ref{cor:ofs_final_local_embedding}, and the fact that closed embeddings are local embeddings.
\end{proof}

\begin{rmk}
	We note that this fact is specific to $\posclos$, as monomorphisms in $\poscat$ are not necessarily local embeddings.
\end{rmk}

\begin{dfn}[Pushforward of a local embedding] \index{map!local embedding!pushforward} \index{$\pfw{p}f$}
	Let $p\colon P \to Q$ be a closed order-preserving map of posets and let $f\colon U \to P$ be a local embedding.
	The \emph{pushforward of $f$ through $p$} is the local embedding
	\[
		\pfw{p}f\colon \pfw{(p \after f)}U \to Q
	\]
	obtained from the factorisation of $p \after f$ as a closed final map followed by a local embedding.
\end{dfn}

\begin{dfn}[Pushforward of a closed subset] \index{poset!subsets!pushforward} \index{$\pfw{p}U$}
	Let $U \subseteq P$ be a closed subset of a poset and let $p\colon P \to Q$ be a closed order-preserving map.
	The \emph{pushforward of $U$ along $p$} is the local embedding
	\[	
		(\restr{p}{U})_\clas{L}\colon \pfw{p}U \to Q
	\]
	where $\pfw{p}U \eqdef \pfw{(\restr{p}{U})}U$.
	
	Given a closed subset $V \subseteq Q$, we write $\pfw{p}U = V$ if
	\begin{enumerate}
		\item the canonical surjective local embedding $\pfw{p}U \to p(U)$ is an isomorphism,
		\item $p(U) = V$.
	\end{enumerate}
	Notice that the first condition is equivalent to $\restr{p}{U}\colon U \to p(U)$ being final.
\end{dfn}

\begin{lem} \label{lem:pushforward_of_atom}
	Let $p\colon P \to Q$ be a closed order-preserving map of posets, $x \in P$.
	Then
	\[
		\pfw{p}\clset{x} = \clset{p(x)}.
	\]
\end{lem}
\begin{proof}
	It suffices to show that the surjective closed order-preserving map $\restr{p}{\clset{x}}\colon \clset{x} \to \clset{p(x)}$ is final.
	But this is necessarily the case, since for all $y, y' \in \clset{x}$ we have the zig-zag $y \leq x \geq y'$.
\end{proof}

\begin{lem} \label{lem:local_embeddings_pullback_stable}
	Let $p\colon P \to Q$ be an order-preserving map of posets and let $f\colon V \to Q$ be a local embedding.
	Then the pullback $\pb{p}f\colon \pb{p}V \to P$ of $f$ along $p$ in $\poscat$ is a local embedding.
\end{lem}
\begin{proof}
	By construction of pullbacks in $\poscat$, the elements of $\pb{p}V$ are pairs $(x, y)$ of $x \in P$ and $y \in V$ such that $p(x) = f(y)$, with the partial order
	\[
		\text{$(x', y') \leq (x, y)$  if and only if  $x' \leq x$ and $y' \leq y$}.
	\]
	Let $(x, y) \in \pb{p}V$, and let $x' \leq \pb{p}f(x, y) = x$.
	Then $p(x') \leq p(x) = f(y)$.
	Because $f$ is a local embedding, there exists a unique $y' \leq y$ such that $f(y') = p(x')$.
	But then $(x', y') \in \clset{(x, y)}$ and $\pb{p}f(x', y') = x'$, which proves that $\pb{p}f$ is closed.
	Moreover, $(x', y')$ is the unique lift of $x'$ to $\clset{(x, y)}$ because $y'$ is the unique lift of $p(x')$ to $\clset{y}$.
	This proves that $\restr{\pb{p}f}{\clset{(x, y)}}$ is a closed embedding, and we conclude that $\pb{p}f$ is a local embedding.
\end{proof}

\begin{dfn}[Pullback of a local embedding] \index{map!local embedding!pullback} \index{$\pb{p}f$}
	Let $p\colon P \to Q$ be an order-preserving map of posets and let $f\colon V \to Q$ be a local embedding.
	The \emph{pullback of $f$ along $p$} is the local embedding
	\[
		\pb{p}f\colon \pb{p}V \to P
	\]
	obtained by pulling $f$ back along $p$ in $\poscat$.
\end{dfn}


\section{Maps of regular directed complexes} \label{sec:maps}

\begin{guide}
	In this section, we define the category $\rdcpxmap$ of regular directed complexes and maps.
	We characterise its subcategory $\rdcpx$ whose morphisms are local embeddings, which coincides with the full subcategory of $\ogpos$ on the regular directed complexes (Proposition 
	\ref{prop:characterisation_of_morphisms_among_maps}).
	We look at limits and colimits, and show that $\rdcpxmap$ admits a ternary factorisation system, induced by the one on $\posclos$ (Corollary \ref{cor:ternary_factorisation_system_rdcpxmap}).
	Then, we define a functor $\molecin{-}\colon \rdcpxmap \to \omegacat$, which coincides with the one defined on $\ogpos$ on their joint subcategory $\rdcpx$, as well as a functor $\freeab{-}\colon \rdcpxmap \to \chaug$.
	Finally, we prove that fibres of maps of regular directed complexes are regular directed complexes (Proposition \ref{prop:fibres_of_rdcpx_are_rdcpx}).
\end{guide}

\begin{dfn}[Map of regular directed complexes] \index{regular directed complex!map} \index{map!of regular directed complexes}
Let $P$, $Q$ be regular directed complexes.
A \emph{map} $p\colon P \to Q$ is a closed order-preserving map of their underlying posets which, for all $x \in P$, $n \in \mathbb{N}$, and $\alpha \in \set{+, -}$, satisfies
\[
	\pfw{p}\bound{n}{\alpha}x = \bound{n}{\alpha}p(x).
\]
\end{dfn}

\begin{rmk}
	This condition is equivalent to $p(\bound{n}{\alpha}x) = \bound{n}{\alpha}p(x)$ and the surjective map $\restr{p}{\bound{n}{\alpha}x}\colon \bound{n}{\alpha}x \to p(\bound{n}{\alpha}x)$ being final.
\end{rmk}

\begin{rmk} \label{rmk:properties_of_maps}
	By Lemma \ref{lem:closed_map_of_graded_posets_dim_non_increasing}, a map of regular directed complexes is dimension-non-increasing.
\end{rmk}

\begin{lem} \label{lem:if_dim_decreased_faces_map_to_same_element}
	Let $p\colon P \to Q$ be a map of regular directed complexes, $x \in P$.
	If $\dim{p(x)} < \dim{x}$, then for all $\alpha \in \set{+, -}$ there exists $x^\alpha \in \faces{}{\alpha}x$ such that $p(x^\alpha) = p(x)$.
\end{lem}
\begin{proof}
	Let $n \eqdef \dim{x} > 0$.
	By Lemma \ref{lem:boundary_inclusion} we have, for all $\alpha \in \set{+, -}$,
	\[
		p(\bound{n-1}{\alpha}x) = \bound{n-1}{\alpha}p(x) = \clset{p(x)}.
	\]
	Since $p$ is order-preserving and $\bound{n-1}{\alpha}x = \clos \faces{}{\alpha}x$, there must exist $x^\alpha \in \faces{}{\alpha}x$ such that $p(x^\alpha) = p(x)$.
\end{proof}

\begin{lem} \label{lem:element_with_decreased_dimension}
	Let $p\colon P \to Q$ be a map of regular directed complexes, $x \in P$.
	Then there exists $x' \leq x$ such that
	\begin{enumerate}
		\item $p(x') = p(x)$,
		\item $\dim{p(x)} = \dim{x'}$.
	\end{enumerate}
\end{lem}
\begin{proof}
	By induction on $\dim{x} - \dim{p(x)}$, which is $\geq 0$ by Remark \ref{rmk:properties_of_maps}.
	If $\dim{x} = \dim{p(x)}$, then the statement is true for $x' \eqdef x$.
	Let $n \eqdef \dim{x}$ and suppose $n > \dim{p(x)}$.
	By Lemma \ref{lem:if_dim_decreased_faces_map_to_same_element} there exists $x' \in \faces{}{}x$ such that $p(x') = p(x)$.
	Since $\dim{x'} - \dim{p(x)} < \dim{x} - \dim{p(x)}$, the inductive hypothesis applies, and we conclude.
\end{proof}

\begin{dfn}[The category $\rdcpxmap$] \index{$\rdcpxmap$}
There is a category $\rdcpxmap$ whose objects are regular directed complexes and morphisms are maps of regular directed complexes. 
\end{dfn}

\begin{comm}
	The downward arrow in $\rdcpxmap$ is meant to remind that morphisms in this category are allowed to \emph{decrease} the dimension.
\end{comm}

\begin{prop} \label{prop:forgetful_from_rdcpxmap_to_pos}
	Forgetting the orientation determines a faithful functor
	\[ \fun{U}\colon \rdcpxmap \to \posclos. \]
\end{prop}
\begin{proof}
	By definition.
\end{proof}

\begin{prop} \label{prop:terminal_object_of_rdcpxmap}
	The point $1$ is a terminal object in $\rdcpxmap$.
\end{prop}
\begin{proof}
	For each regular directed complex $P$, there exists a unique function $p\colon P \to 1$.
	For all $x \in P$, $n \in \mathbb{N}$, and $\alpha \in \set{+, -}$, obviously $p(\bound{n}{\alpha}x) = \bound{n}{\alpha}p(x)$, and $\restr{p}{\bound{n}{\alpha}x}$ is final because $\bound{n}{\alpha}x$ is connected, by Lemma \ref{lem:molecules_are_connected}.
\end{proof}

\begin{lem} \label{lem:morphisms_are_maps}
	Let $P$, $Q$ be regular directed complexes and $f\colon P \to Q$ a morphism of oriented graded posets.
	Then $f$ is a map.
\end{lem}
\begin{proof}
	Let $x \in P$.
	Since both $P$ and $Q$ are regular directed complexes, the restriction $\restr{f}{\clset{x}}\colon \clset{x} \to \clset{f(x)}$ is a morphism of atoms of the same dimension.
	By Theorem \ref{thm:morphisms_of_atoms_are_injective}, it is an isomorphism, so by Corollary \ref{cor:inclusions_preserve_boundaries} it restricts to an isomorphism on all input and output boundaries.
\end{proof}

\begin{lem} \label{lem:injective_maps_are_inclusions}
	Let $p\colon P \to Q$ be a map of regular directed complexes.
	The following are equivalent:
	\begin{enumerate}[label=(\alph*)]
		\item $p$ is injective;
		\item $p$ is an inclusion, that is, an injective morphism of oriented graded posets.
	\end{enumerate}
\end{lem}
\begin{proof}
	If $p$ is an injective morphism, it is an injective map by Lemma \ref{lem:morphisms_are_maps}.
	Conversely, suppose $p$ is an injective map.
	Then $\fun{U}p$ is a closed embedding, so by 
	Lemma \ref{lem:local_embedding_preserves_faces}, for all $x \in P$, it induces a bijection between $\faces{}{}x$ and $\faces{}{}p(x)$.
	For all $\alpha \in \set{ +, - }$ this, coupled with the fact that $p(\bound{}{\alpha}x) = \bound{}{\alpha}p(x)$, implies that $p$ induces a bijection between $\faces{}{\alpha}x$ and $\faces{}{\alpha}p(x)$, that is, $p$ is a morphism.
\end{proof}

\begin{comm}
	Lemma \ref{lem:injective_maps_are_inclusions} allows us to not distinguish between injective maps and injective morphisms of regular directed complexes, and call both of them inclusions.
	In particular, isomorphisms of regular directed complexes in $\rdcpxmap$ are the same as their isomorphisms in $\ogpos$, that is, by Lemma \ref{lem:characterisation_of_isomorphisms}, the surjective inclusions.
\end{comm}

\begin{exm}[A surjective, non-injective map between atoms of the same dimension] \index[counterex]{A surjective, non-injective map between atoms of the same dimension} \label{exm:noninjective_map_atom}
	Consider the atom $\disk{2}{1}$ with the numbering of elements
	\[
\begin{tikzcd}[sep=small]
	{{\scriptstyle 0}\;\bullet} && {{\scriptstyle 2}\;\bullet} \\
	& {{\scriptstyle 1}\;\bullet}
	\arrow[""{name=0, anchor=center, inner sep=0}, "2", curve={height=-18pt}, from=1-1, to=1-3]
	\arrow["0"', curve={height=6pt}, from=1-1, to=2-2]
	\arrow["1"', curve={height=6pt}, from=2-2, to=1-3]
	\arrow["0", shorten <=3pt, shorten >=6pt, Rightarrow, from=2-2, to=0]
\end{tikzcd}
\]
	as well as the 2\nbd globe $\globe{2}$ with the same numbering as in (\ref{eq:2globe}).
	There is a surjective map $p\colon \disk{2}{1} \to \globe{2}$ defined by
\begin{align*}
	(0, 0) & \mapsto (0, 0), \quad 
	       & (0, 1), (0, 2) & \mapsto (0, 1), \\ 
	(1, 0) & \mapsto (1, 0), \quad 
		       & (1, 1) & \mapsto (0, 1), \quad 
		       & (1, 2) & \mapsto (1, 1), \\
	(2, 0) & \mapsto (2, 0), && &&
\end{align*}
	which can be seen as ``collapsing'' the cell $(1, 1)$ in the input boundary of $\disk{2}{1}$.
	There is also a dual map collapsing $(1, 0)$.

	This shows that Theorem \ref{thm:morphisms_of_atoms_are_injective} does not extend to surjective maps.
	Notice that $p$ is also an example of a surjective map of atoms which does not admit a section, since the section would have to be an injective map, hence an inclusion by Lemma \ref{lem:injective_maps_are_inclusions}, and the only inclusions of atoms of the same dimension are isomorphisms.
	
	An immediate consequence is that the full subcategory of $\rdcpxmap$ on the atoms is not a \emph{Reedy category}.
	There are restricted classes of maps of atoms which do form a Reedy category --- for example, maps which are posetal Grothendieck fibrations --- but we will not study them in the book.
\end{exm}

\begin{dfn}[Final map of regular directed complexes] \index{map!final}
	A map of regular directed complexes is \emph{final} if its underlying map of posets is final.
\end{dfn}

\begin{rmk}
	If $p$ is a final map of regular directed complexes, $\fun{U}p$ is a closed final map of posets, so by Lemma \ref{lem:final_closed_maps_are_surjective} $p$ is surjective.
\end{rmk}

\begin{lem} \label{lem:final_map_injective_in_topdim}
	Let $p\colon P \to Q$ be a final map of regular directed complexes, $n \eqdef \dim{P} < \infty$, and let $x, x' \in \grade{n}{P}$.
	If $p(x) = p(x')$ and $\dim{p(x)} = n$, then $x = x'$.
\end{lem}
\begin{proof}
	Since $p$ is final, there exists a zig-zag
	\[
		x \leq x_1 \geq x_2 \leq \ldots \geq x_{m-1} \leq x_m \geq x'
	\]
	in $P$ such that $p(x) = p(x') \leq p(x_i)$ for all $i \in \set{1, \ldots, m}$.
	By Lemma \ref{lem:dimension_is_monotonic}, it follows that $n = \dim{p(x)} \leq \dim{p(x_i)}$, which, since $p$ is dimension-non-increasing and $\dim{x_i} \leq n$, implies that $\dim{x_i} = \dim{p(x_i)} = n$ for all $i \in \set{1, \ldots, m}$. 
	We conclude that $x = x_i = x'$ for all $i \in \set{1, \ldots, m}$.
\end{proof}

\begin{lem} \label{lem:two_faces_only_in_dim_of_image}
	Let $p\colon P \to Q$ be a map of regular directed complexes, let $x \in P$, and $k \eqdef \dim {p(x)}$.
	Then, for all $\alpha \in \set{+, -}$, there exists a unique $x^\alpha \in \faces{k}{\alpha}x$ such that $p(x) = p(x^\alpha)$.
\end{lem}
\begin{proof}
	Let $\alpha \in \set{+, -}$.
	We have $p(\bound{k}{\alpha}x) = \bound{k}{\alpha}p(x) = \clset{p(x)}$, so there exists $x^\alpha \in \bound{k}{\alpha}x$ such that $p(x^\alpha) = p(x)$, and because $p$ is dimension-non-increasing, necessarily $x^\alpha \in \grade{k}{(\bound{k}{\alpha}x)} = \faces{k}{\alpha}x$.
	Moreover, $\restr{p}{\bound{k}{\alpha}x}\colon \bound{k}{\alpha}x \to \clset{p(x)}$ is a final map of $k$\nbd dimensional regular directed complexes and $\grade{k}{(\clset{p(x)})} = \set{p(x)}$.
	We conclude by Lemma \ref{lem:final_map_injective_in_topdim}.
\end{proof}

\begin{dfn}[Local embedding of regular directed complexes] \index{map!local embedding}
	A map of regular directed complexes is a \emph{local embedding} if its underlying map of posets is a local embedding.
\end{dfn}

\begin{prop} \label{prop:characterisation_of_morphisms_among_maps} 
	Let $p\colon P \to Q$ be a map of regular directed complexes.
	The following are equivalent:
	\begin{enumerate}[label=(\alph*)]
		\item $p$ is a morphism of oriented graded posets;
		\item $p$ is a local embedding;
		\item $p$ is dimension-preserving.
	\end{enumerate}
\end{prop}
\begin{proof}
	If $p$ is a morphism of oriented graded posets, it is a local embedding by Corollary \ref{cor:morphisms_of_rdcpx_are_local_isomorphisms}.
	Moreover, a local embedding of graded posets is evidently dimension-preserving.
	
	Suppose that $p$ is a dimension-preserving map and let $x \in P$ and $\alpha \in \set{+, -}$.
	If $\dim{x} = 0$, then $\faces{}{\alpha}x$ and $\faces{}{\alpha}p(x)$ are both empty, so $p$ trivially induces a bijection.
	Suppose $n \eqdef \dim{x} > 0$.
	Then 
	\[
		p(\faces{}{\alpha}x) = p(\grade{n-1}{(\bound{n-1}{\alpha}x)}) = \grade{n-1}{(\bound{n-1}{\alpha}p(x))} = \faces{}{\alpha}p(x).
	\]
	Let $y, y' \in \faces{}{\alpha}x$ and suppose $p(y) = p(y')$.
	Because $\restr{p}{\bound{n-1}{\alpha}x}\colon \bound{n-1}{\alpha}x \to \bound{n-1}{\alpha}p(x)$ is final, it follows from Lemma \ref{lem:final_map_injective_in_topdim} that $y = y'$.
	We conclude that $p$ determines a bijection between $\faces{}{\alpha}x$ and $\faces{}{\alpha}p(x)$, that is, $p$ is a morphism of oriented graded posets.
\end{proof}

\begin{dfn}[The category $\rdcpx$] \index{$\rdcpx$}
	We let $\rdcpx$ denote the wide subcategory of $\rdcpxmap$ whose morphisms are local embeddings.
\end{dfn}

\begin{comm}
	The equality sign in $\rdcpx$ is meant to remind that morphisms in this category must \emph{preserve} the dimension.
\end{comm}

\begin{cor} \label{cor:rdcpx_as_full_subcategory}
	The category $\rdcpx$ can be identified with the full subcategory of $\ogpos$ whose objects are the regular directed complexes.
	Moreover, the diagram of functors
\[\begin{tikzcd}
	\rdcpx && \ogpos \\
	\rdcpxmap && \posclos
	\arrow[hook, from=1-1, to=1-3]
	\arrow[hook', from=1-1, to=2-1]
	\arrow["{\fun{U}}", from=2-1, to=2-3]
	\arrow["{\fun{U}}", from=1-3, to=2-3]
\end{tikzcd}\]
	commutes.
\end{cor}

\begin{rmk}
	In fact, the inclusion $\rdcpx \incl \ogpos$ factors through the inclusion $\ogposloc \incl \ogpos$.
\end{rmk}

\begin{prop} \label{prop:functor_from_rdcpx_to_otgpos}
	There is a unique functor $\augm{(-)}\colon \rdcpx \to \otgpos$ such that the diagram of functors
\[\begin{tikzcd}
	\rdcpx && \otgpos \\
	\ogpos && \ogposbot
	\arrow["{\augm{(-)}}", hook, from=1-1, to=1-3]
	\arrow[hook', from=1-1, to=2-1]
	\arrow["{\augm{(-)}}", from=2-1, to=2-3]
	\arrow[hook', from=1-3, to=2-3]
\end{tikzcd}\]
	commutes.
	Moreover, this functor is full and faithful.
\end{prop}
\begin{proof}
	The functor exists by Proposition \ref{prop:oriented_diamond_rdc}, and is full and faithful by Proposition \ref{prop:ogpos_equivalent_ogposbot} and Corollary \ref{cor:rdcpx_as_full_subcategory}.
	Uniqueness is enforced by the requirement of strict commutativity.
\end{proof}

\begin{prop} \label{prop:pullbacks_of_inclusions_rdcpxmap}
	The category $\rdcpxmap$ has pullbacks of inclusions, and they are both preserved and reflected by $\fun{U}\colon \rdcpxmap \to \posclos$ and by $\rdcpx \incl \rdcpxmap$.
	Moreover, inclusions are stable under pullback.
\end{prop}
\begin{proof}
	A straightforward variant of the proof of Proposition \ref{prop:pullbacks_of_inclusions}, using Lemma \ref{lem:if_local_embeds_into_rdcpx_then_is_rdcpx} to deduce that the pullback is a regular directed complex.
\end{proof}

\begin{lem} \label{lem:colimits_of_inclusions_in_rdcpx}
	Let $\fun{F}\colon C \to \rdcpx$ be a diagram of inclusions of regular directed complexes, and suppose $\gamma$ is a colimit cone under $\fun{UF}$ whose components are all closed embeddings.
	Then
	\begin{enumerate}
		\item there exists a unique colimit cone $\vec{\gamma}$ under $\fun{F}$ such that $\gamma = \fun{U}\vec{\gamma}$, whose components are all inclusions;
		\item this colimit is preserved and reflected by $\rdcpx \incl \rdcpxmap$.
	\end{enumerate}
\end{lem}
\begin{proof}
	By Lemma \ref{lem:reflected_colimits_in_ogpos}, the colimit cone lifts to $\ogpos$.
	Since $\rdcpx$ is a full subcategory of $\ogpos$, this is a colimit cone in $\rdcpx$ as long as its image lies in $\rdcpx$.
	Let $P$ be the tip of the cone and $x \in P$.
	Then there exist an object $c$ in $C$ and $y \in \fun{F}c$ such that $x = \vec{\gamma}_c(y)$.
	Since $\vec{\gamma}_c$ is an inclusion, it restricts to an isomorphism between $\clset{y}$ and $\clset{x}$.
	Because $\fun{F}c$ is a regular directed complex, $\clset{y}$ is an atom, hence $\clset{x}$ is an atom.

	Finally, to show that $\rdcpx \incl \rdcpxmap$ preserves this colimit, consider a cone $\eta$ under $\fun{F}$ in $\rdcpxmap$ with tip $Q$.
	Then $\fun{U}$ maps this cone onto a cone under $\fun{UF}$ in $\posclos$.
	Let $f\colon \fun{U}P \to \fun{U}Q$ be the closed order-preserving map induced by the universal property of $\gamma$, and let $x \in P$.
	As seen before, there exist an object $c$ in $C$ and $y \in \fun{F}c$ such that $x = \vec{\gamma}_c(y)$.
	Then 
	\[
		\restr{f}{\clset{x}} = \restr{\fun{U}\eta_c}{\clset{y}} \after \invrs{(\restr{\gamma_c}{\clset{y}})},
	\]
	and both factors on the right-hand side lift to maps of regular directed complexes.
	It follows that $f$ lifts to a map of regular directed complexes.
\end{proof}

\begin{prop} \label{prop:pushouts_in_rdcpx}
The category $\rdcpx$ has
\begin{enumerate}
    \item a strict initial object $\varnothing$,
    \item pushouts of inclusions along inclusions,
\end{enumerate}
preserved and reflected by $\fun{U}\colon \rdcpx \to \posclos$ and by the subcategory inclusion $\rdcpx \incl \rdcpxmap$.
Moreover,
\begin{enumerate}
	\item the pushout of an inclusion along an inclusion is an inclusion,
	\item a pushout square of inclusions is also a pullback square.
\end{enumerate}
\end{prop}
\begin{proof}
	Straightforward variant of Proposition \ref{prop:pushouts_in_ogpos} using Lemma \ref{lem:colimits_of_inclusions_in_rdcpx} and Proposition \ref{prop:pullbacks_of_inclusions_rdcpxmap}.
\end{proof}

\begin{prop} \label{prop:em_factorisation_in_rdcpxmap}
Let $p\colon P \to Q$ be a map of regular directed complexes.
Then $p$ factors as
\begin{enumerate}
	\item a surjective map $\widehat{p}\colon P \to p(P)$,
	\item followed by an inclusion $\imath\colon p(P) \incl Q$.
\end{enumerate}
This factorisation is unique up to unique isomorphism.
\end{prop}
\begin{proof}
	Straightforward variant of Proposition \ref{prop:factorisation_in_ogpos}.
\end{proof}

\begin{cor} \label{cor:em_ofs_on_rdcpxmap}
The classes of
\begin{enumerate}
	\item surjective maps,
	\item inclusions
\end{enumerate}
form an orthogonal factorisation system on $\rdcpxmap$.
\end{cor}

\begin{prop} \label{prop:comprehensive_factorisation_system_on_rdcpxmap}
Let $p\colon P \to Q$ be a map of regular directed complexes.
Then $p$ factors as
\begin{enumerate}
	\item a final map $p_\clas{F}\colon P \to \pfw{p}P$,
	\item followed by a local embedding $p_\clas{L}\colon \pfw{p}P \to Q$.
\end{enumerate}
This factorisation is unique up to unique isomorphism.
\end{prop}
\begin{proof}
	By Proposition \ref{prop:comprehensive_factorisation_system}, the underlying closed order-preserving map of $p$ factors as a closed final map $p_\clas{F}\colon \fun{U}P \to \pfw{p}P$ followed by a local embedding $p_\clas{L}\colon \pfw{p}P \to \fun{U}Q$.
	By Lemma \ref{lem:if_local_embeds_into_rdcpx_then_is_rdcpx} and Proposition 
	\ref{prop:characterisation_of_morphisms_among_maps}, the latter lifts uniquely to a local embedding of regular directed complexes, so it suffices to show that $p_\clas{F}$ lifts to a map of regular directed complexes.
	For all $x \in P$,
	\[
		\restr{p_\clas{F}}{\clset{x}} = 
		\invrs{(\restr{p_\clas{L}}{\clset{p_\clas{F}(x)}})} \after \restr{p}{\clset{x}},
	\]
	and both factors on the right-hand side lift to maps of regular directed complexes.
	It follows that $p_\clas{F}$ lifts to a map of regular directed complexes.
\end{proof}

\begin{cor} \label{cor:ofs_final_local_embedding_rdcpxmap}
The classes of
\begin{enumerate}
	\item final maps,
	\item local embeddings
\end{enumerate}
form an orthogonal factorisation system on $\rdcpxmap$.
\end{cor}

\begin{cor} \label{cor:ternary_factorisation_system_rdcpxmap}
The classes of
\begin{enumerate}
	\item final maps,
	\item surjective local embeddings,
	\item inclusions
\end{enumerate}
form a ternary factorisation system on $\rdcpxmap$.
\end{cor}

\begin{prop} \label{prop:final_map_preserves_molecules}
	Let $p\colon U \to V$ be a final map of regular directed complexes, and suppose $U$ is a molecule.
	Then
	\begin{enumerate}
		\item $V$ is a molecule,
		\item for all $n \in \mathbb{N}$ and $\alpha \in \set{+, -}$, $\pfw{p}\bound{n}{\alpha}U = \bound{n}{\alpha}V$,
		\item for all $k \in \mathbb{N}$, if $U$ splits into submolecules $W \cup W'$ along the $k$\nbd boundary, then $V$ splits into submolecules $\pfw{p}W \cup \pfw{p}W'$ along the $k$\nbd boundary.
	\end{enumerate}
\end{prop}
\begin{proof}
	We proceed by induction on the construction of $U$.
	If $U$ was produced by (\textit{Point}) or (\textit{Atom}), it is an atom by Lemma
	\ref{lem:atom_greatest_element}, so it has a greatest element $\top$.
	Since $p$ is final, it is surjective, so $V = p(\clset{\top}) = \clset{p(\top)}$.
	Because $V$ is a regular directed complex, it follows that $V$ is an atom, and
	\[
		\pfw{p}\bound{n}{\alpha}U = \pfw{p}\bound{n}{\alpha}\top = \bound{n}{\alpha}p(\top) = \bound{n}{\alpha}V
	\]
	by definition of map.
	
	Suppose that $U$ was produced by (\textit{Paste}).
	Then $U$ is obtained from a pushout of the form
\[\begin{tikzcd}
	{\overline{W}} && {W'} \\
	W && U
	\arrow[hook', from=1-1, to=2-1]
	\arrow[hook', from=1-3, to=2-3]
	\arrow[hook, from=2-1, to=2-3]
	\arrow[hook, from=1-1, to=1-3]
	\arrow["\lrcorner"{anchor=center, pos=0.125, rotate=180}, draw=none, from=2-3, to=1-1]
\end{tikzcd}\]
	with $\overline{W} = \bound{k}{+}W = \bound{k}{-}W'$ for some $k < \min \set{\dim{W}, \dim{W'}}$ and molecules $W, W'$.
	Factorising the restrictions of $p$ as final maps followed by local embeddings, we obtain a commutative diagram
\[\begin{tikzcd}
	&&& {\pfw{p}\overline{W}} && {\pfw{p}W'} \\
	{\overline{W}} && {W'} & {\pfw{p}W} \\
	W && U &&&& V
	\arrow[hook', from=2-1, to=3-1]
	\arrow[hook', from=2-3, to=3-3]
	\arrow[hook, from=3-1, to=3-3]
	\arrow[hook, from=2-1, to=2-3]
	\arrow["\lrcorner"{anchor=center, pos=0.125, rotate=180}, draw=none, from=3-3, to=2-1]
	\arrow["p", from=3-3, to=3-7]
	\arrow["{\overline{f}}", near start, curve={height=-6pt}, from=2-1, to=1-4]
	\arrow["f'", near start, curve={height=-6pt}, from=2-3, to=1-6]
	\arrow["f", near start, crossing over, curve={height=-6pt}, from=3-1, to=2-4]
	\arrow["j", hook, from=1-4, to=1-6]
	\arrow["\ell", from=2-4, to=3-7]
	\arrow["\ell'", curve={height=-6pt}, from=1-6, to=3-7]
	\arrow["\imath"', crossing over, near start, hook', from=1-4, to=2-4]
\end{tikzcd}\]
where, by the inductive hypothesis, $\pfw{p}W$ and $\pfw{p}W'$ are molecules, $\imath$ and $j$ are inclusions, and
\begin{align*}
	\pfw{p}\overline{W} & = \pfw{p}(\bound{k}{+}W) = \bound{k}{+}\pfw{p}W = \\
			    & = \pfw{p}(\bound{k}{-}W') = \bound{k}{-}\pfw{p}W'.
\end{align*}
By Proposition \ref{prop:pushouts_in_rdcpx}, the pushout of $\imath$ and $j$ exists in $\rdcpx$ and constructs the molecule $\pfw{p}W \cp{k} \pfw{p}W'$.
Since $\imath, j, \ell, \ell'$ are all local embeddings, the universal property of the pushout $\pfw{p}W \cp{k} \pfw{p}W'$ in $\rdcpx$ produces a unique local embedding
\[
	\varphi\colon \pfw{p}W \cp{k} \pfw{p}W' \to V
\]
such that the diagram
\[\begin{tikzcd}
	&&& {\pfw{p}\overline{W}} && {\pfw{p}W'} \\
	{\overline{W}} && {W'} & {\pfw{p}W} && {\pfw{p}W \cp{k} \pfw{p}W'} \\
	W && U &&&& V	
	\arrow[hook', from=2-1, to=3-1]
	\arrow[hook', from=2-3, to=3-3]
	\arrow[hook, from=3-1, to=3-3]
	\arrow[hook, from=2-1, to=2-3]
	\arrow["\lrcorner"{anchor=center, pos=0.125, rotate=180}, draw=none, from=3-3, to=2-1]
	\arrow["p", from=3-3, to=3-7]
	\arrow["{\overline{f}}", near start, curve={height=-6pt}, from=2-1, to=1-4]
	\arrow["f'", near start, curve={height=-6pt}, from=2-3, to=1-6]
	\arrow["f", near start, crossing over, curve={height=-6pt}, from=3-1, to=2-4]
	\arrow["j", hook, from=1-4, to=1-6]
	\arrow["\ell", from=2-4, to=3-7]
	\arrow["\ell'", curve={height=-6pt}, from=1-6, to=3-7]
	\arrow["\imath"', crossing over, near start, hook', from=1-4, to=2-4]
	\arrow[hook', from=1-6, to=2-6]
	\arrow[hook, from=2-4, to=2-6]
	\arrow["\varphi", from=2-6, to=3-7]
	\arrow["\lrcorner"{anchor=center, pos=0.125, rotate=180}, draw=none, from=2-6, to=1-4]
\end{tikzcd}\]
commutes.
Finally, the universal property of the pushout diagram exhibiting $U$ as a pasting in $\rdcpxmap$ produces a unique map
\[
	p'\colon U \to \pfw{p}W \cp{k} \pfw{p}W'
\]
such that the diagram
\[\begin{tikzcd}
	&&& {\pfw{p}\overline{W}} && {\pfw{p}W'} \\
	{\overline{W}} && {W'} & {\pfw{p}W} && {\pfw{p}W \cp{k} \pfw{p}W'} \\
	W && U &&&& V	
	\arrow[hook', from=2-1, to=3-1]
	\arrow[hook', from=2-3, to=3-3]
	\arrow[hook, from=3-1, to=3-3]
	\arrow[hook, from=2-1, to=2-3]
	\arrow["\lrcorner"{anchor=center, pos=0.125, rotate=180}, draw=none, from=3-3, to=2-1]
	\arrow["p", from=3-3, to=3-7]
	\arrow["{\overline{f}}", near start, curve={height=-6pt}, from=2-1, to=1-4]
	\arrow["f'", near start, curve={height=-6pt}, from=2-3, to=1-6]
	\arrow["f", near start, crossing over, curve={height=-6pt}, from=3-1, to=2-4]
	\arrow["j", hook, from=1-4, to=1-6]
	\arrow["\imath"', crossing over, near start, hook', from=1-4, to=2-4]
	\arrow[hook', from=1-6, to=2-6]
	\arrow[hook, from=2-4, to=2-6]
	\arrow["\varphi", from=2-6, to=3-7]
	\arrow["\lrcorner"{anchor=center, pos=0.125, rotate=180}, draw=none, from=2-6, to=1-4]
	\arrow["{p'}", from=3-3, to=2-6]
\end{tikzcd}\]
commutes.
Because $f, f', \overline{f}$ are all final maps, by general properties of left classes of orthogonal factorisation systems \cite[Proposition 2.1.1.(c)]{freyd1972categories}, we know that $p'$ is a final map.
Since $p = \varphi \after p'$, with $p$ and $p'$ final and $\varphi$ a local embedding, we conclude that $\varphi$ is an isomorphism.
It follows that $V$ is a molecule, splitting into $\pfw{p}W \cup \pfw{p}W'$ along the $k$\nbd boundary.

It only remains to show that, for all $n \in \mathbb{N}$ and $\alpha \in \set{+, -}$, $\pfw{p}\bound{n}{\alpha}U = \bound{n}{\alpha}V$.
For $n \leq k$, $\bound{n}{\alpha}U$ is contained in $\bound{n}{\alpha}W$ or $\bound{n}{\alpha}W'$, so we already know this from the inductive hypothesis.
For $n > k$, by Lemma \ref{lem:pasting_higher_boundary} $\bound{n}{\alpha}U$ splits into $\bound{n}{\alpha}W \cup \bound{n}{\alpha}W'$ along the $k$\nbd boundary.
Then, by what we have already proved, $\pfw{p}\bound{n}{\alpha}U$ splits into $\pfw{p}\bound{n}{\alpha}W \cup \pfw{p}\bound{n}{\alpha}W'$ along the $k$\nbd boundary, so it is uniquely isomorphic to $\bound{n}{\alpha}(\pfw{p}W \cp{k} \pfw{p}W')$, which is uniquely isomorphic to $\bound{n}{\alpha}V$.
\end{proof}

\begin{exm}[The ternary factorisation of a map]
	Let $U$ and $P$ be the oriented face posets of the 1\nbd dimensional diagrams
	\[\begin{tikzcd}
	{{\scriptstyle 0}\;\bullet} & {{\scriptstyle 1}\;\bullet} & {{\scriptstyle 2}\;\bullet} & {{\scriptstyle 3}\;\bullet}
	\arrow["0", from=1-1, to=1-2]
	\arrow["1", from=1-2, to=1-3]
	\arrow["2", from=1-3, to=1-4]
\end{tikzcd}\;, \quad \quad
\begin{tikzcd}
	{{\scriptstyle 0}\;\bullet} & {{\scriptstyle 1}\;\bullet} & {{\scriptstyle 2}\;\bullet}
	\arrow["0"', curve={height=12pt}, from=1-1, to=1-2]
	\arrow["2", from=1-2, to=1-3]
	\arrow["1"', curve={height=12pt}, from=1-2, to=1-1]
\end{tikzcd} \;, \]
	respectively, and let $p\colon U \to P$ be the map defined by
\begin{align*}
	(0, 0), (0, 3) & \mapsto (0, 0), \quad 
	       & (0, 1), (0, 2) & \mapsto (0, 1), \\ 
	(1, 0) & \mapsto (1, 0), \quad 
		& (1, 1) & \mapsto (0, 1), \quad
		       & (1, 2) & \mapsto (1, 1).
\end{align*}
	Let $\pfw{p}{U}$ and $p(U)$ be the oriented face posets of the diagrams
	\[\begin{tikzcd}
	{{\scriptstyle 0}\;\bullet} & {{\scriptstyle 1}\;\bullet} & {{\scriptstyle 2}\;\bullet}
	\arrow["0", from=1-1, to=1-2]
	\arrow["1", from=1-2, to=1-3]
\end{tikzcd}\;, \quad \quad 
	\begin{tikzcd}
	{{\scriptstyle 0}\;\bullet} & {{\scriptstyle 1}\;\bullet}
	\arrow["0"', curve={height=12pt}, from=1-1, to=1-2]
	\arrow["1"', curve={height=12pt}, from=1-2, to=1-1]
\end{tikzcd} \;, \]
	respectively.
	Then $p$ factors as
	\begin{enumerate}
		\item the final map $p_{\clas{F}}\colon U \to \pfw{p}{U}$ defined by
\begin{align*}
	(0, 0) & \mapsto (0, 0), \quad 
	       & (0, 1), (0, 2) & \mapsto (0, 1), 
	       & (0, 3) & \mapsto (0, 2), \\ 
	(1, 0) & \mapsto (1, 0), \quad 
		& (1, 1) & \mapsto (0, 1), \quad
		       & (1, 2) & \mapsto (1, 1),
\end{align*}
		which can be seen as the ``purely collapsing part'' of $p$, which does not identify any cells unless they are in the boundary of a collapsed higher-dimensional cell,
	\item followed by the surjective local embedding $\widehat{p}\colon \pfw{p}{U} \to p(U)$ defined by
\begin{align*}
	(0, 0), (0, 2) & \mapsto (0, 0), \quad 
	       & (0, 1) & \mapsto (0, 1), \\ 
	(1, 0) & \mapsto (1, 0), \quad 
		& (1, 1) & \mapsto (1, 1),
\end{align*}
	which can be seen as the part of $p$ which ``rigidly identifies'' some cells,
\item followed by the evident inclusion $p(U) \incl P$.
\end{enumerate}
	Since $U$ is a molecule, we see that, compatibly with Proposition 
	\ref{prop:final_map_preserves_molecules}, the pushforward $\pfw{p}{U}$ is also a molecule, whereas the direct image $p(U)$ is no longer a molecule.
\end{exm}

\begin{thm} \label{thm:molecin_is_a_functor_on_maps} \index{strict functor!induced by a map}
	Let $p\colon P \to Q$ be a map of regular directed complexes.
  	Then there exists a strict functor of strict $\omega$\nbd categories defined by
	\begin{align*}
		\pfw{p} \equiv \molecin{p}\colon \molecin{P} & \to \molecin{Q}, \\
			\isocl{ f\colon U \to P } & \mapsto \isocl{ \pfw{p}f\colon \pfw{(p \after f)}U \to Q }.
	\end{align*}
	This assignment determines a functor $\molecin{-}\colon \rdcpxmap \to \omegacat$, such that the square of functors
	\begin{equation} \label{eq:molecin_morphisms_and_maps}
\begin{tikzcd}
	\rdcpx && \ogpos \\
	\rdcpxmap && \omegacat
	\arrow[hook, from=1-1, to=1-3]
	\arrow[hook', from=1-1, to=2-1]
	\arrow["{\molecin{-}}"', from=2-1, to=2-3]
	\arrow["{\molecin{-}}", from=1-3, to=2-3]
\end{tikzcd}
	\end{equation}
	commutes.
\end{thm}
\begin{proof}
	The fact that the assignment is functorial and respects isomorphism classes is an immediate consequence of the functoriality of orthogonal factorisations.
	The fact that $\pfw{(p \after f)}U$ is a molecule whenever $U$ is a molecule, as well as compatibility with boundary operators and $k$\nbd composition operations, all follow straightforwardly from Proposition \ref{prop:final_map_preserves_molecules}.
	
	Finally, commutativity of the diagram (\ref{eq:molecin_morphisms_and_maps}) follows from the fact that, when $p$ is already a local embedding, $p \after f = (p \after f) \after \idd{U}$ is a factorisation of $p \after f$ as a final map followed by a local embedding.
	Thus the assignment restricts to the one from Proposition \ref{prop:molecin_is_a_functor} on the local embeddings, that is, on the maps in the image of $\rdcpx$.
\end{proof}

\begin{rmk} \label{rmk:characterisation_of_maps_by_functors}
	As a kind of converse to this result, observe that if $P$, $Q$ are regular directed complexes, and $p\colon P \to Q$ is an order-preserving map of their underlying posets such that
	\[
		\pfw{p}\colon \isocl{f\colon U \to P} \mapsto \isocl{\pfw{p}f\colon \pfw{(p \after f)}U \to Q}
	\]
	determines a strict functor $\molecin{P} \to \molecin{Q}$, then by virtue of $\pfw{p}$ being a morphism of reflexive $\omega$\nbd graphs, and the inclusions $\isocl{\clset{x} \incl P}$ being cells in $\molecin{P}$, it must be the case that for all $x \in P$, $n \in \mathbb{N}$, and $\alpha \in \set{+, -}$,
\[
	\pfw{p}\bound{n}{\alpha}x = \bound{n}{\alpha}p(x),
\]
that is, $p$ is a map.
Thus maps of regular directed complexes are characterised among functions with an underlying order-preserving map of posets by the property that they induce strict functors by pushforward.
\end{rmk}

\begin{prop} \label{prop:pseudomonicity_of_molecin_maps}
	The functor $\molecin{-}\colon \rdcpxmap \to \omegacat$ is pseudomonic, that is,
	\begin{enumerate}
		\item it is faithful,
		\item it is full on isomorphisms,
		\item it reflects isomorphisms.
	\end{enumerate}
\end{prop}
\begin{proof}
	Let $P$, $Q$ be regular directed complexes, let $p, p'\colon P \to Q$ be maps, and suppose that $\pfw{p} = \pfw{p'}$.
	Then for all $x \in P$, we have 
	\[
		\isocl{\clset{p(x)} \incl Q} = \pfw{p}\isocl{\clset{x} \incl P} = \pfw{p'}\isocl{\clset{x} \incl P} = \isocl{\clset{p'(x)} \incl Q},
	\]
	so $\clset{p(x)} = \clset{p'(x)}$ as closed subsets of $Q$.
	It follows that $p(x) = p'(x)$, hence $p = p'$.

	Next, let $\varphi\colon \molecin{P} \to \molecin{Q}$ be an isomorphism of strict $\omega$\nbd categories with inverse $\invrs{\varphi}$.
	Then, in particular, both $\varphi$ and its inverse preserve the dimension of cells.
	Consider $x \in P$ and let $\isocl{f\colon U \to Q} \eqdef \varphi\isocl{\clset{x} \incl P}$.
	We claim that $U$ must be an atom.
	Assume for the sake of contradiction that $U$ is not an atom, or, equivalently by Lemma 
	\ref{lem:layering_dimension_atom}, that $k \eqdef \lydim{U} \geq 0$.
	Then $U$ admits a $k$\nbd layering $(\order{i}{U})_{i=1}^m$, inducing a decomposition
	\[
		\isocl{f} = \isocl{\order{1}{f}} \cp{k} \ldots \cp{k} \isocl{\order{m}{f}}
	\]
	of $\isocl{f}$, with each factor having dimension $> k$.
	Then
	\[
		\isocl{\clset{x} \incl P} = \invrs{\varphi}\isocl{f} = \invrs{\varphi}\isocl{\order{1}{f}} \cp{k} \ldots \cp{k} \invrs{\varphi}\isocl{\order{m}{f}},
	\]
	and each factor has dimension $> k$.
	By Lemma \ref{lem:layering_dimensions_pasting}, this would imply that $\lydim{\clset{x}} \geq k$, a contradiction.

	By Corollary \ref{cor:atoms_over_regular_directed_complexes}, we conclude that $\isocl{f} = \isocl{\clset{y} \incl Q}$ for a unique $y \in Q$, and we let $\widehat{\varphi}(x) \eqdef y$.
	Such an assignment for each $x \in P$ defines a function $\widehat{\varphi}\colon P \to Q$, which is easily determined to be an isomorphism of the underlying posets of $P$ and $Q$; it suffices to show that it is a map of regular directed complexes.
	We claim that, for all $x \in P$ and $U \submol \clset{x}$, we have
\[
	\varphi \isocl{U \incl P} = \isocl{\widehat{\varphi}(U) \incl Q}.
\]
	We proceed by induction on submolecules of $\clset{x}$.
	If $U = \clset{y}$ for some $y \leq x$, then the statement is true by construction of $\widehat{\varphi}$.
	In particular, this is true when $\dim{y} = 0$, which takes care of the base case.
	Suppose that $U$ is not an atom.
	Then there exist proper submolecules $V, W \submol U$ and $k < \min \set{\dim{V}, \dim{W}}$ such that $U$ splits into $V \cup W$ along the $k$\nbd boundary.
	Then
	\begin{align*}
	\varphi \isocl{U \incl P} & = \varphi \isocl{V \incl P} \cp{k} \varphi \isocl{W \incl P} = \isocl{\widehat{\varphi}(V) \incl Q} \cp{k} \isocl{\widehat{\varphi}(W) \incl Q} = \\
				  & = \isocl{\widehat{\varphi}(V) \cup \widehat{\varphi}(W) \incl Q} = \isocl{\widehat{\varphi}(V \cup W) \incl Q} = \isocl{\widehat{\varphi}(U) \incl Q},
	\end{align*}
	using the fact that $\varphi$ respects $k$\nbd composition, the inductive hypothesis applied to $V$ and $W$, and the fact that direct images preserve unions.
	This concludes the inductive proof.
	Let $x \in P$, $n \in \mathbb{N}$, and $\alpha \in \set{+, -}$.
	Since $\bound{n}{\alpha}x \submol \clset{x}$ by Lemma \ref{lem:boundary_is_submolecule}, we have in particular
	\[
		\isocl{\widehat{\varphi}(\bound{n}{\alpha}x) \incl Q} = \varphi \bound{n}{\alpha}\isocl{\clset{x} \incl P} = \bound{n}{\alpha}\varphi\isocl{\clset{x} \incl P} = \isocl{\bound{n}{\alpha}\widehat{\varphi}(x) \incl Q}.
	\]
	This proves that $\widehat{\varphi}(\bound{n}{\alpha}x) = \bound{n}{\alpha}\widehat{\varphi}(x)$, and since $\widehat{\varphi}$ is injective, all its restrictions are final maps onto their image.
	We conclude that $\widehat{\varphi}$ is a map of regular directed complexes such that $\varphi = \pfw{\widehat{\varphi}}$.

	Finally, observe that we not only determined that $\varphi$ lifts to a map, but that it lifts to a bijective map, hence, by Lemma \ref{lem:characterisation_of_isomorphisms} and Lemma 
	\ref{lem:injective_maps_are_inclusions}, to an isomorphism of regular directed complexes.
	This implies that $\molecin{-}$ reflects isomorphisms.
\end{proof}

\begin{prop} \label{prop:chain_complex_is_functorial_on_maps}
	Let $p\colon P \to Q$ be a map of regular directed complexes.
	Then the sequence of homomorphisms
	\begin{align*}
		\freeab{\grade{n}{p}}\colon \freeab{\grade{n}{P}} & \to \freeab{\grade{n}{Q}}, \\
		x \in \grade{n}{P} & \mapsto \begin{cases}
			p(x) 
			& \text{if $\dim{p(x)} = n$}, \\
			0
			& \text{if $\dim{p(x)} < n$}
		\end{cases}
	\end{align*}
	is a homomorphism $\freeab{p}\colon \freeab{P} \to \freeab{Q}$ of augmented chain complexes.
	This assignment determines a functor $\freeab{-}\colon \rdcpxmap \to \chaug$, such that the diagram of functors
\begin{equation} \label{eq:map_functor_to_chaug}
\begin{tikzcd}
	\rdcpx && \otgpos \\
	\rdcpxmap && \chaug
	\arrow["{\augm{(-)}}", hook, from=1-1, to=1-3]
	\arrow[hook', from=1-1, to=2-1]
	\arrow["{\freeab{-}}", from=2-1, to=2-3]
	\arrow["{\freeab{-}}", from=1-3, to=2-3]
\end{tikzcd}\end{equation}
	commutes.
\end{prop}
\begin{proof}
	Let $x \in \grade{0}{P}$.
	Then $\dim{p(x)} = 0$, so $\freeab{\grade{0}{p}}(x) = p(x)$, and $\eau(p(x)) = 1$, so $\eau \after \freeab{\grade{0}{p}} = \eau$.
	Next, let $x \in \grade{n}{P}$ for $n > 0$.
	We distinguish three cases.
	First, suppose that $\dim{p(x)} = n$.
	Then
	\[
		\der (\freeab{\grade{n}{p}}(x)) = \der (p(x)) = 
		\; \sum_{\mathclap{y \in \faces{}{+}p(x)}}\; y \; - 
		\;\; \sum_{\mathclap{y \in \faces{}{-}p(x)}}\; y\,.
	\]
	For all $\alpha \in \set{+, -}$, since $\restr{p}{\bound{n-1}{\alpha}x}$ is final onto its image, by Lemma \ref{lem:final_map_injective_in_topdim} $\restr{p}{\faces{}{\alpha}x}$ is injective.
	Moreover, for all $y' \in \faces{}{\alpha}x$, either $\dim{p(y')} = n-1$ and $p(y') \in \faces{}{\alpha}p(x)$, or $\dim{p(y')} < n - 1$ and then $\freeab{\grade{n-1}{p}}(y') = 0$.
	It follows that
	\[
		\der (\freeab{\grade{n}{p}}(x)) = \sum_{y' \in \faces{}{+}x} \freeab{\grade{n-1}{p}}(y') - \sum_{y' \in \faces{}{-}x} \freeab{\grade{n-1}{p}}(y') = \freeab{\grade{n-1}{p}}(\der x).
	\]
	Next, suppose that $\dim{p(x)} = n - 1$.
	Then by Lemma \ref{lem:two_faces_only_in_dim_of_image} there exists a unique $y^\alpha \in \faces{}{\alpha}x$ such that $p(y^\alpha) = p(x)$.
	For all other $y' \in \faces{}{\alpha}x$, we have $\dim{p(y')} < n - 1$, hence $\freeab{\grade{n-1}{p}}(y') = 0$.
	Then
	\begin{align*}
		\der (\freeab{\grade{n}{p}}(x)) & = \der (0) = 0 = 
		p(y^+) - p(y^-) = \\
						& = \sum_{y' \in \faces{}{+}x} \freeab{\grade{n-1}{p}}(y') - \sum_{y' \in \faces{}{+}x} \freeab{\grade{n-1}{p}}(y') = \freeab{\grade{n-1}{p}}(\der x).
	\end{align*}
	Finally, if $\dim{p(x)} < n - 1$, then also $\dim{p(y')} < n - 1$ for all $y' \in \faces{}{}x$, so
	\[
		\der (\freeab{\grade{n}{p}}(x)) = 0 = \freeab{\grade{n-1}{p}}(\der x).
	\]
	This proves that $\der \after \freeab{\grade{n}{p}} = \freeab{\grade{n-1}{p}} \after \der$, so	$(\freeab{\grade{n}{p}})_{n \in \mathbb{N}}$ is a homomorphism of augmented chain complexes.
	Functoriality and commutativity of (\ref{eq:map_functor_to_chaug}) are straightforward.
\end{proof}

\begin{dfn}[Collapsing dimension] \index{map!collapsing dimension} \index{dimension!collapsing} \index{$\cdim{p}x$}
	Let $p\colon P \to Q$ be a map of regular directed complexes, $x \in P$.
	The \emph{collapsing dimension of $x$ through $p$} is the natural number $\cdim{p}x \eqdef \dim x - \dim p(x)$.
\end{dfn}

\begin{lem} \label{lem:grading_of_fibres_of_maps}
	Let $p\colon P \to Q$ be a map of regular directed complexes, $y \in Q$.
	Then the fibre $\invrs{p}y$ is graded by $x \mapsto \cdim{p}x$.
\end{lem}
\begin{proof}
	Because $P$ has locally finite height, so does $\invrs{p}y$.
	To prove that $\invrs{p}y$ is graded, it then suffices to show that
	\begin{enumerate}
		\item $x \in \invrs{p}y$ is minimal if and only if $\cdim{p}x = 0$,
		\item if $x' \prec x$ in $\invrs{p}y$, then $\cdim{p}x' = \cdim{p}x - 1$.
	\end{enumerate}
	Suppose that $\cdim{p}x = 0$.
	Then $x$ is minimal because $p$ is dimension-non-increasing.
	Conversely, suppose that $\cdim{p}x > 0$.
	Then, by Lemma \ref{lem:if_dim_decreased_faces_map_to_same_element}, there exists $x' \in \faces{}{}x$ such that $p(x) = p(x')$, so $x$ is not minimal in $\invrs{p}y$.
	For the second point, suppose that $x' \prec x$ in the fibre, and suppose that $x' < x'' \leq x$.
	Then $y = p(x') \leq p(x'') \leq p(x) = y$, that is, $x''$ is in the fibre, which implies that $x'' = x$.
	It follows that $x' \prec x$ in $P$, so $\dim{x'} = \dim{x} - 1$ and $\cdim{p}{x'} = \cdim{p}{x} - 1$.
\end{proof}

\begin{lem} \label{lem:boundaries_of_fibre_of_top_element}
	Let $p\colon U \to V$ be a map of regular directed complexes such that $V$ is an atom with greatest element $\top$, and let $m \eqdef \dim{V}$.
	Then, for all $k \in \mathbb{N}$ and $\alpha \in \set{+, -}$, we have $\bound{k}{\alpha}\invrs{p}\top = \invrs{p}\top \cap \bound{k+m}{\alpha}U$.
\end{lem}
\begin{proof}
	By Lemma \ref{lem:grading_of_fibres_of_maps}, $\invrs{p}\top$ is a graded poset and it inherits an orientation from $U$, so $\bound{k}{\alpha}\invrs{p}\top$ is well-defined.
	Suppose that $x \in \faces{k}{\alpha}\invrs{p}\top$.
	Then $\cdim{p}x = \dim{x} - \dim{\top} = \dim{x} - m = k$, so $\dim{x} = k+m$.
	Because $y \in \cofaces{}{}x$ implies $\top = p(x) \leq p(y)$, so $p(y) = \top$, the set of cofaces of $x$ in $U$ is equal to the set of cofaces of $x$ in $\invrs{p}\top$, so the latter being empty implies that $x \in \faces{k+m}{\alpha}U$.
	The converse is analogous, and proves that $\faces{k}{\alpha}\invrs{p}\top = \invrs{p}\top \cap \faces{k+m}{\alpha}U$.
	Similarly, we derive $\grade{k}{(\maxel{\invrs{p}\top})} = \invrs{p}\top \cap \grade{k+m}{(\maxel{U})}$.
	Finally, for all subsets $V \subseteq \invrs{p}\top$, the closure of $V$ in $\invrs{p}\top$ is equal to $\invrs{p}\top \cap \clos V$, from which we conclude.
\end{proof}

\begin{lem} \label{lem:fibre_of_top_element_preserves_globularity_roundness}
	Let $p\colon U \to V$ be a map of regular directed complexes such that $V$ is an atom with greatest element $\top$.
	Then
	\begin{enumerate}
		\item if $U$ is globular, then $\invrs{p}\top$ is globular,
		\item if $U$ is round, then $\invrs{p}\top$ is round.
	\end{enumerate}
\end{lem}
\begin{proof}
	Immediate from Lemma \ref{lem:boundaries_of_fibre_of_top_element}.
\end{proof}

\begin{lem} \label{lem:fibre_of_top_element_is_molecule}
	Let $p\colon U \to V$ be a map of regular directed complexes such that $U$ is a molecule and $V$ is an atom with greatest element $\top$.
	Then $\invrs{p}\top$ is either empty, or it is a molecule.
\end{lem}
\begin{proof}
	By Proposition \ref{prop:comprehensive_factorisation_system_on_rdcpxmap}, we can factorise $p$ as a final map $p_\mathscr{F}\colon U \to p_*U$ followed by a local embedding $p_\mathscr{L}\colon p_*U \to V$ whose image is $p(U)$, and $p_*U$ is a molecule by Proposition \ref{prop:final_map_preserves_molecules}.
	We have $\dim{p_*U} \leq \dim{V}$.
	If $\dim{p_*U} < \dim{V}$, then $\top \notin p(U)$, and $\invrs{p}\top$ is empty.
	If $\dim{p_*U} = \dim{V}$, then by Proposition \ref{prop:molecule_over_atom_is_iso} combined with Proposition \ref{prop:characterisation_of_morphisms_among_maps}, $p_\mathscr{L}$ is an isomorphism.
	Thus, if $\invrs{p}\top$ is not empty, $p$ is final and preserves all boundaries.

	Assuming $p$ final, let $n \eqdef \dim{U}$ and $m \eqdef \dim{V}$; we will prove by induction on $k \eqdef n - m \geq 0$ that, in this case, $\invrs{p}\top$ is a molecule of dimension $k$.
	Suppose that $k = 0$.
	Then by Lemma \ref{lem:final_map_injective_in_topdim} there exists a unique $x \in \invrs{p}\top$, so $\invrs{p}\top$ is a point, that is, a 0\nbd dimensional molecule.
	Suppose that $k > 0$; we proceed by induction on submolecules.
	If $U$ is an atom, for each $\alpha \in \set{+, -}$, the map $p^\alpha \eqdef \restr{p}{\bound{}{\alpha}U}\colon \bound{}{\alpha}U \to V$ is final, and by the inductive hypothesis $\invrs{(p^\alpha)}\top = \invrs{p}\top \cap \bound{}{\alpha}U = \bound{}{\alpha}\invrs{p}\top$ is a $(k-1)$\nbd dimensional molecule, which by Lemma \ref{lem:fibre_of_top_element_preserves_globularity_roundness} is round.
	Moreover,
	\[
		\bound{}{+}\invrs{p}\top \cap \bound{}{-}\invrs{p}\top =
		\invrs{p}\top \cap \bound{}{+}U \cap \bound{}{-}U =
		\invrs{p}\top \cap \bound{n-2}{}U = \bound{k-2}{}\invrs{p}\top
	\]
	which suffices to prove that $\invrs{p}\top$ is a $k$\nbd dimensional atom.

	Finally, suppose that $U$ splits into proper submolecules $U_1 \cp{\ell} U_2$, and for each $i \in \set{1, 2}$ let $p_i \eqdef \restr{p}{U_i}\colon U_i \to V$.
	By the inductive hypothesis, for each $i \in \set{1, 2}$, $\invrs{p_i}\top$ is either empty or a molecule, and in the latter case $p_i$ is final.
	If either is empty, then $\invrs{p}\top$ is equal to the other, and we are done.
	Suppose that $\ell < m$; we claim that either $\invrs{p_1}{\top}$ or $\invrs{p_2}\top$ is empty.
	Suppose by way of contradiction that $\top \in p(U_1) \cap p(U_2)$.
	Then $p(U_1) = p(U_2) = V$, so
	\[
		\bound{\ell}{+}V = p(\bound{\ell}{+}U_1) = p(\bound{\ell}{-}U_2) = \bound{\ell}{-}V,
	\]
	contradicting $\dim V = m > \ell$ by Corollary \ref{cor:dimension_from_boundary}.
	It remains to consider the case where $\ell \geq m$, and both $\invrs{p_1}\top$ and $\invrs{p_2}\top$ are non-empty.
	Then, $\invrs{p}\top = \invrs{p_1}\top \cup \invrs{p_2}\top$, and
	\begin{align*}
		\invrs{p_1}\top \cap \invrs{p_2}\top =
		\invrs{p}\top \cap U_1 \cap U_2
		&= \invrs{p}\top \cap \bound{\ell}{+}U_1 =
		\bound{\ell-m}{+}\invrs{p_1}\top \\
		&= \invrs{p}\top \cap \bound{\ell}{-}U_2 =
		\bound{\ell-m}{-}\invrs{p_2}\top.
	\end{align*}
	We conclude that $\invrs{p}\top$ splits into $\invrs{p_1}\top \cp{\ell-m} \invrs{p_2}\top$.
\end{proof}

\begin{prop} \label{prop:fibres_of_rdcpx_are_rdcpx}
	Let $p\colon P \to Q$ be a map of regular directed complexes and $y \in Q$.
	Then the fibre $\invrs{p}y$ is a regular directed complex.
\end{prop}
\begin{proof}
	Let $x \in \invrs{p}y$, and let $p_x \eqdef \restr{p}{\clset{x}}\colon \clset{x} \to \clset{y}$.
	Then $\clset{x} \cap \invrs{p}y$ is isomorphic to $\invrs{p_x}y$, which is an atom by Lemma \ref{lem:fibre_of_top_element_is_molecule}.
\end{proof}


\section{Comaps of regular directed complexes} \label{sec:comaps}

\begin{guide}
	In this section, we define comaps of regular directed complexes and prove that they compose (Corollary \ref{cor:comaps_compose}), which allows us to define a category $\rdcpxcomap$ of regular directed complexes and comaps.
	We prove that the ``intersection'' of maps and comaps consists precisely of isomorphisms (Proposition \ref{prop:maps_comaps_trivial_intersection}), which also reassures us that the notion of isomorphism of regular directed complexes remains stable across the different categories that they form.
	Finally, we show that comaps naturally determine both strict functors of strict $\omega$\nbd categories and homomorphisms of augmented chain complexes, compatibly with the constructions already defined on the overlaps, but \emph{contravariantly} instead of covariantly.
\end{guide}

\begin{dfn}[Comap of regular directed complexes] \index{map!comap} \index{regular directed complex!comap} \index{comap}
Let $P$, $Q$ be regular directed complexes.
A \emph{comap} $c\colon P \to Q$ is an order-preserving map of their underlying posets such that, for all $y \in Q$, $n \in \mathbb{N}$, and $\alpha \in \set{+, -}$,
\begin{enumerate}
	\item $\invrs{c}{\clset{y}}$ is a molecule,
	\item $\invrs{c}{\bound{n}{\alpha}y} = \bound{n}{\alpha}\invrs{c}\clset{y}$.
\end{enumerate}
\end{dfn}

\begin{lem} \label{lem:comaps_increase_dimension}
Let $c\colon P \to Q$ be a comap of regular directed complexes.
Then $c$ is dimension-non-decreasing.
\end{lem}
\begin{proof}
Let $x \in P$ and let $U \eqdef \invrs{c}{\clset{x}}$ and $n \eqdef \dim{c(x)}$.
For all $\alpha \in \set{+, -}$,
\[
	\bound{n}{\alpha}U = \invrs{c}\bound{n}{\alpha}c(x) = \invrs{c}\clset{c(x)} = U,
\]
so by Corollary \ref{cor:dimension_from_boundary} $\dim{U} \leq n$.
Since $x \in U$, necessarily $\dim{x} \leq n$.
\end{proof}

\begin{prop} \label{prop:comaps_inverse_image_preserves_molecules}
Let $c\colon U \to V$ be a comap of regular directed complexes, and suppose $V$ is a molecule.
Then
\begin{enumerate}
	\item $U$ is a molecule,
	\item $\dim{U} = \dim{V}$,
	\item for all $n \in \mathbb{N}$ and $\alpha \in \set{+, -}$, $\invrs{c}\bound{n}{\alpha}V = \bound{n}{\alpha}U$,
	\item for all $k \in \mathbb{N}$, if $V$ splits into submolecules $W \cup W'$ along the $k$\nbd boundary, then $U$ splits into submolecules $\invrs{c}W \cup \invrs{c}W'$ along the $k$\nbd boundary.
\end{enumerate}
\end{prop}
\begin{proof}
We proceed by induction on the construction of $V$.
If $V$ was produced by (\textit{Point}), then $V$ is the point.
Since $c$ is a comap, $U = \invrs{c}V$ is a molecule, and since $c$ is dimension-non-increasing, $\dim{U} = 0$.
It follows from Lemma \ref{lem:only_0_molecule} that $U$ is also the point.
All statements are then evident.

Suppose that $V$ was produced by (\textit{Paste}).
Then $V$ splits into submolecules $W \cup W'$ along the $k$\nbd boundary.
By the inductive hypothesis, $\invrs{c}W$ and $\invrs{c}W'$ are molecules with $\dim{W} = \dim{\invrs{c}W}$ and $\dim{W'} = \dim{\invrs{c}{W}}$.
Moreover,
\begin{align*}
	\invrs{c}W \cap \invrs{c}W' & = \invrs{c}(W \cap W') = \begin{cases}
		\invrs{c}\bound{k}{+}W = \bound{k}{+}\invrs{c}W, \\
		\invrs{c}\bound{k}{-}W' = \bound{k}{-}\invrs{c}W'
	\end{cases} \\
		\invrs{c}W \cup \invrs{c}W' & = \invrs{c}(W \cup W') = \invrs{c}V = U
\end{align*}
using the inductive hypothesis and the fact that inverse images preserve both unions and intersections.
It follows that $U$ is a molecule, splitting into submolecules $\invrs{c}W$ and $\invrs{c}W'$ along the $k$\nbd boundary.
Moreover,
\[
	\dim{U} = \max \set{\dim{\invrs{c}W}, \dim{\invrs{c}W'}} = \max \set{\dim{W}, \dim{W'}} = \dim{V},
\]
while the fact that $\invrs{c}\bound{n}{\alpha}V = \bound{n}{\alpha}U$ follows from 
\[
	\bound{n}{\alpha}U = \bound{n}{\alpha}\invrs{c}W \cup \bound{n}{\alpha}\invrs{c}W' = \invrs{c}(\bound{n}{\alpha}W \cup \bound{n}{\alpha}W') = \invrs{c}\bound{n}{\alpha}V
\]
for $n > k$, and from the inductive hypothesis for $n \leq k$.

Finally, suppose that $V$ was produced by (\textit{Atom}).
Then $V$ has a greatest element $\top$.
It follows that $U = \invrs{c}V = \invrs{c}\clset{\top}$ is a molecule, and for all $n \in \mathbb{N}$ and $\alpha \in \set{+, -}$,
\[
	\invrs{c}\bound{n}{\alpha}V = \invrs{c}\bound{n}{\alpha}\top = \bound{n}{\alpha}\invrs{c}\clset{\top} = \bound{n}{\alpha}U
\]
by definition of comap.

Let $n \eqdef \dim{V}$; it only remains to show that $\dim{U} = n$.
Because $c$ is dimension-non-decreasing, necessarily $\dim{U} \leq n$.
Moreover, by the inductive hypothesis, for all $\alpha \in \set{+, -}$, $\bound{n-1}{\alpha}U$ is an $(n-1)$\nbd dimensional molecule equal to $\invrs{c}\bound{n-1}{\alpha}V$, so $n-1 \leq \dim{U} \leq n$.
Suppose for the sake of contradiction that $\dim{U} = n-1$.
Then $U = \bound{n-1}{-}U = \bound{n-1}{+}U$, so
\[
	U = \bound{n-1}{-}U \cap \bound{n-1}{+}U = \invrs{c}\bound{n-1}{-}V \cap \invrs{c}\bound{n-1}{+}V = \invrs{c}\bound{n-2}{}V = \bound{n-2}{}U,
\]
using the inductive hypothesis, the fact that inverse images preserve intersections, and that $V$ is round by Corollary \ref{cor:atoms_are_round}.
Then $\dim{U} \leq n - 2$, a contradiction.
We conclude that $\dim{U} = n$.
\end{proof}

\begin{cor} \label{cor:comaps_compose}
	Let $c\colon P \to Q$ and $d\colon Q \to R$ be comaps of regular directed complexes.
	Then $d \after c\colon P \to R$ is a comap.
\end{cor}
\begin{proof}
	Straightforward consequence of Proposition \ref{prop:comaps_inverse_image_preserves_molecules}.
\end{proof}

\begin{lem} \label{lem:comap_inverse_image_boundaries}
	Let $c\colon U \to V$ be a comap of molecules, $n \in \mathbb{N}$, and $\alpha \in \set{+, -}$.
	Then $\faces{n}{\alpha}U = \grade{n}{(\invrs{c}\faces{n}{\alpha}V)}$.
\end{lem}
\begin{proof}
	By Proposition \ref{prop:comaps_inverse_image_preserves_molecules} and 
	Lemma \ref{lem:maximal_in_boundary}, we have 
	\[
		\faces{n}{\alpha}U = \grade{n}{(\bound{n}{\alpha}U)} = 
		\grade{n}{(\invrs{c}\bound{n}{\alpha}V)}.
	\]
	Let $x \in \faces{n}{\alpha}U$.
	By Lemma \ref{lem:comaps_increase_dimension}, we have $n \leq \dim{c(x)}$, but $c(x) \in \bound{n}{\alpha}U$ which is at most $n$\nbd dimensional, so $\dim{c(x)} = n$.
	We conclude that $c(x) \in \faces{n}{\alpha}V$.
\end{proof}

\begin{dfn}[The category $\rdcpxcomap$] \index{$\rdcpxcomap$}
There is a category $\rdcpxcomap$ whose objects are regular directed complexes and morphisms are comaps of regular directed complexes.
\end{dfn}

\begin{comm}
	The upward arrow in $\rdcpxcomap$ is meant to remind that morphisms in this category are allowed to \emph{increase} the dimension.
\end{comm}

\begin{prop} \label{prop:forgetful_from_rdcpxcomap_to_pos}
	Forgetting the orientation determines a faithful functor 
	\[ \fun{U}\colon \rdcpxcomap \to \poscat. \]
\end{prop}
\begin{proof}
	By definition.
\end{proof}

\begin{lem} \label{lem:comaps_inverse_image_and_roundness}
	Let $c\colon U \to V$ be a comap of molecules.
	If $V$ is round, then $U$ is round.
\end{lem}
\begin{proof}
	By Proposition \ref{prop:comaps_inverse_image_preserves_molecules}, $\dim{U} = \dim{V}$.
	Suppose that $V$ is round. 
	Then, for all $n < \dim{U}$,
	\begin{align*}
		\bound{n}{-}U \cap \bound{n}{+}U & = \invrs{c}\bound{n}{-}V \cap \invrs{c}\bound{n}{+}V = \invrs{c}(\bound{n}{-}V \cap \bound{n}{+}V) = \\
						 & = \invrs{c}(\bound{n-1}{}V) = \invrs{c}\bound{n-1}{-}V \cup \invrs{c}\bound{n-1}{+}V = \bound{n-1}{}U
	\end{align*}
	since inverse images preserve unions and intersections.
\end{proof}

\begin{lem} \label{lem:comaps_inverse_image_interiors}
	Let $c\colon U \to V$ be a comap of molecules, $n \in \mathbb{N}$, $\alpha \in \set{+, -}$.
	Then $\invrs{c}\inter{\bound{n}{\alpha}V} = \inter{\bound{n}{\alpha}U}$.
\end{lem}
\begin{proof}
	Follows immediately from Proposition \ref{prop:comaps_inverse_image_preserves_molecules} and the fact that inverse images preserve set differences.
\end{proof}

\begin{lem} \label{lem:comaps_are_surjective}
Let $c\colon P \to Q$ be a comap of regular directed complexes.
Then $c$ is surjective.
\end{lem}
\begin{proof}
Let $y \in Q$.
By Proposition \ref{prop:comaps_inverse_image_preserves_molecules}, $U \eqdef \invrs{c}{\clset{y}}$ is a molecule with $\dim{U} = \dim{\clset{y}} = \dim{y}$.
Let $x \in U$ be an element of dimension $\dim{U}$.
Then $c(x) \in \clset{y}$ and $\dim{y} = \dim{x} \leq \dim{c(x)}$, so necessarily $c(x) = y$.
\end{proof}

\begin{exm}[A surjective, non-injective comap between atoms of the same dimension] \index[counterex]{A surjective, non-injective comap between atoms of the same dimension}
	Consider the atoms $\disk{2}{1}$ and $\globe{2}$ exactly as in Example \ref{exm:noninjective_map_atom}.
	There is a comap $c\colon \disk{2}{1} \to \globe{2}$ defined by
\begin{align*}
	(0, 0) & \mapsto (0, 0), \quad 
	       & (0, 1) & \mapsto (1, 0), \quad
	       & (0, 2) & \mapsto (0, 1) \\ 
	(1, 0) & \mapsto (1, 0), \quad 
		       & (1, 1) & \mapsto (1, 0), \quad 
		       & (1, 2) & \mapsto (1, 1), \\
	(2, 0) & \mapsto (2, 0), && &&
\end{align*}
	which can be seen as ``merging'' the interior of the input boundary of $\disk{2}{1}$ into a single cell.
	We can also see the formal dual of $c$ as exhibiting a \emph{subdivision} of the input boundary of $\globe{2}$, which is an atom, into a non-atomic 1\nbd dimensional molecule.
	
	By Proposition \ref{prop:comaps_inverse_image_preserves_molecules}, Lemma 
	\ref{lem:comaps_inverse_image_and_roundness}, Lemma 
	\ref{lem:comaps_inverse_image_interiors}, all comaps are of this sort: the inverse image of each atom is a round molecule, whose boundaries are inverse images of the boundaries of the atom, and this allows us to interpret the comap as dual to a subdivision.
	Again, we find that Theorem \ref{thm:morphisms_of_atoms_are_injective} cannot be extended: there is a rich combinatorics of comaps between atoms of the same dimension.

	We note that, by Proposition \ref{prop:comap_to_globe}, $c$ is the \emph{only} comap from $\disk{2}{1}$ to $\globe{2}$. 
\end{exm}

\begin{prop} \label{prop:maps_comaps_trivial_intersection}
Let $P$, $Q$ be regular directed complexes and $f\colon P \to Q$ be an order-preserving map of their underlying posets.
The following are equivalent:
\begin{enumerate}[label=(\alph*)]
	\item $f$ is both a map and a comap of regular directed complexes;
	\item $f$ is an invertible map of regular directed complexes;
	\item $f$ is an invertible comap of regular directed complexes;
	\item $f$ is an isomorphism of oriented graded posets.
\end{enumerate}
\end{prop}
\begin{proof}
	The equivalence between invertible maps and isomorphisms of oriented graded posets follows from Lemma \ref{lem:morphisms_are_maps} and Lemma \ref{lem:injective_maps_are_inclusions}.

	Suppose that $f$ is an invertible map with inverse $\invrs{f}$.
	Then, for all $y \in Q$, $n \in \mathbb{N}$, and $\alpha \in \set{+, -}$, we have that $\invrs{f}\clset{y} = \clset{\invrs{f}(y)}$ is an atom, and that $\invrs{f}\bound{n}{\alpha}y = \bound{n}{\alpha}\invrs{f}y$.
	It follows that $f$ is a comap, and by symmetry so is $\invrs{f}$.

	Conversely, suppose that $f$ is an invertible comap with inverse $\invrs{f}$.
	Then, for all $x \in P$, $n \in \mathbb{N}$, and $\alpha \in \set{+, -}$, we have that 
	\[
		f(\bound{n}{\alpha}x) = \invrs{(\invrs{f})}(\bound{n}{\alpha}x) = \bound{n}{\alpha}\invrs{(\invrs{f})}\clset{x} = \bound{n}{\alpha}f(x).
	\]
	Moreover, $\restr{f}{\bound{n}{\alpha}x}$ is a bijection with its image, so in particular it is final.
	This proves that $f$ is a map, and dually so is $\invrs{f}$.

	We have proved that invertible maps, invertible comaps, and isomorphisms of oriented graded posets all coincide.
	This also implies that an isomorphism of oriented graded posets is both a map and a comap.

	To conclude, suppose that $f$ is both a map and a comap.
	By Remark \ref{rmk:properties_of_maps} and Lemma \ref{lem:comaps_increase_dimension}, $f$ is both dimension-non-increasing and dimension-non-decreasing, so it is dimension-preserving.
	By Proposition \ref{prop:characterisation_of_morphisms_among_maps}, it is a morphism of oriented graded posets, and by Lemma \ref{lem:comaps_are_surjective} it is also surjective.
	
	It remains to show that $f$ is also injective.
	Let $x, x' \in P$ and suppose that $f(x) = f(x')$.
	Let $V \eqdef \clset{f(x)}$ and $U \eqdef \invrs{f}V$, so that $x, x' \in U$.
	Since $f$ is a comap, by Proposition \ref{prop:comaps_inverse_image_preserves_molecules}, $U$ is a molecule of the same dimension as $V$.
	By Proposition \ref{prop:molecule_over_atom_is_iso}, because $V$ is an atom, the restriction $\restr{f}{U}\colon U \to V$ is an isomorphism.
	We conclude that $x = x'$, and $f$ is an isomorphism of oriented graded posets.
\end{proof}

\begin{lem} \label{lem:pullback_of_local_embedding_along_comap}
	Let $c\colon P \to Q$ be a comap and $f\colon V \to Q$ be a local embedding of regular directed complexes.
	Consider the pullback
\[\begin{tikzcd}
	{\pb{c}V} && {\fun{U}V} \\
	{\fun{U}P} && {\fun{U}Q}
	\arrow["{\fun{U}f}", from=1-3, to=2-3]
	\arrow["{\fun{U}c}", from=2-1, to=2-3]
	\arrow["d", from=1-1, to=1-3]
	\arrow["{\pb{c}f}", from=1-1, to=2-1]
	\arrow["\lrcorner"{anchor=center, pos=0.125}, draw=none, from=1-1, to=2-3]
\end{tikzcd}\]
	of their underlying order-preserving maps in $\poscat$.
	Then
	\begin{enumerate}
		\item $\pb{c}f$ lifts to a local embedding of regular directed complexes,
		\item $d$ lifts to a comap of regular directed complexes.
	\end{enumerate}
\end{lem}
\begin{proof}
	By Lemma \ref{lem:local_embeddings_pullback_stable}, $\pb{c}f$ is a local embedding, so by Lemma \ref{lem:local_embeddings_lift_orientation} and Lemma \ref{lem:if_local_embeds_into_rdcpx_then_is_rdcpx} it lifts to a local embedding of regular directed complexes.
	It then suffices to show that $d$ lifts to a comap.

	Let $y \in V$, and consider $\invrs{d}\clset{y}$.
	By the pasting law for pullbacks, the restriction
	\[
		\restr{\pb{c}f}{\invrs{d}\clset{y}}\colon \invrs{d}\clset{y} \to P
	\]
	is equal, up to unique isomorphism, to the pullback of the restriction 
	\[
		\restr{f}{\clset{y}}\colon \clset{y} \to Q
	\]
	along $c$.
	Since $f$ is a local embedding, $\restr{f}{\clset{y}}$ is a closed embedding.
	By Lemma \ref{lem:pullbacks_of_closed_embeddings}, closed embeddings are stable under pullbacks in $\poscat$, so $\restr{\pb{c}f}{\invrs{d}\clset{y}}$ is a closed embedding with image $\invrs{c}\clset{f(y)}$.
	Because $c$ is a comap, the latter is a molecule, so $\invrs{d}\clset{y}$ is a molecule.

	By the same argument, for all $n \in \mathbb{N}$ and $\alpha \in \set{+, -}$, we find that
	\[
		\restr{\pb{c}f}{\invrs{d}\bound{n}{\alpha}y}\colon \invrs{d}\bound{n}{\alpha}y \to P
	\]
	is a closed embedding with image $\invrs{c}\bound{n}{\alpha}f(y) = \bound{n}{\alpha}\invrs{c}\clset{f(y)}$.
	It is therefore equal, up to unique isomorphism, to the restriction of $\restr{\pb{c}f}{\invrs{d}\clset{y}}$ to $\bound{n}{\alpha}\invrs{d}\clset{y}$.
	We conclude that $\invrs{d}\bound{n}{\alpha}y = \bound{n}{\alpha}\invrs{d}\clset{y}$, hence $d$ lifts to a comap of regular directed complexes.
\end{proof}

\begin{dfn}[The groupoid $\rdcpxiso$] \index{$\rdcpxiso$}
	We let $\rdcpxiso$ denote the groupoid whose objects are regular directed complexes and morphisms are isomorphisms.
\end{dfn}

\begin{rmk}
	By Proposition \ref{prop:maps_comaps_trivial_intersection}, $\rdcpxiso$ may be identified with the core groupoid of any of the categories $\rdcpx$, $\rdcpxmap$, or $\rdcpxcomap$.	
\end{rmk}

\begin{thm} \label{thm:molecin_is_a_functor_on_comaps} \index{strict functor!induced by a comap}
	Let $c\colon P \to Q$ be a comap of regular directed complexes.
	Then there exists a strict functor of strict $\omega$\nbd categories defined by
	\begin{align*}
		\pb{c} \equiv \molecin{\pb{c}}\colon \molecin{Q} & \to \molecin{P}, \\
			\isocl{ f\colon U \to Q } & \mapsto \isocl{ \pb{c}f\colon \pb{c}U \to P }.
	\end{align*}
	This assignment determines a functor $\molecin{\pb{-}}\colon \opp{\rdcpxcomap} \to \omegacat$, such that the diagram of functors
\begin{equation} \label{eq:molecin_maps_and_comaps}
\begin{tikzcd}
	{\opp{\rdcpxiso}} & \rdcpxiso && \rdcpxmap \\
	{\opp{\rdcpxcomap}} &&& \omegacat
	\arrow["{\invrs{(-)}}", from=1-1, to=1-2]
	\arrow[hook, from=1-2, to=1-4]
	\arrow[hook', from=1-1, to=2-1]
	\arrow["{\molecin{\pb{-}}}", from=2-1, to=2-4]
	\arrow["{\molecin{-}}", from=1-4, to=2-4]
\end{tikzcd}
\end{equation}
	commutes.
\end{thm}
\begin{proof}
	The fact that the assignment respects isomorphism classes follows from the universal property of pullbacks, and functoriality from the pasting law for pullbacks.
The fact that $\pb{c}U$ is a molecule whenever $U$ is a molecule, as well as compatibility with boundary operators and $k$\nbd composition operations, follow straightforwardly from Proposition \ref{prop:comaps_inverse_image_preserves_molecules}, using Lemma \ref{lem:pullback_of_local_embedding_along_comap} to deduce that pullback produces a comap $\pb{c}U \to U$.
Finally, commutativity of diagram (\ref{eq:molecin_maps_and_comaps}) follows from the fact that pullback along an isomorphism is equal, up to unique isomorphism, to composition with its inverse.
\end{proof}

\begin{rmk}
	A dual of Remark \ref{rmk:characterisation_of_maps_by_functors} holds about comaps: they are exactly characterised, among functions with an underlying order-preserving map of posets, by the property that they induce strict functors between $\omega$\nbd categories of molecules by pullback.
\end{rmk}

\begin{prop} \label{prop:pseudomonicity_of_molecin_comaps}
	The functor $\molecin{\pb{-}}\colon \opp{\rdcpxcomap} \to \omegacat$ is pseudomonic, that is,
	\begin{enumerate}
		\item it is faithful,
		\item it is full on isomorphisms,
		\item it reflects isomorphisms.
	\end{enumerate}
\end{prop}
\begin{proof}
	Let $P$, $Q$ be regular directed complexes and $\varphi\colon \molecin{P} \to \molecin{Q}$ an isomorphism of strict $\omega$\nbd categories with inverse $\psi$.
	By Proposition \ref{prop:pseudomonicity_of_molecin_maps}, $\psi$ lifts to an invertible map $\widehat{\psi}\colon Q \to P$, which by Proposition \ref{prop:maps_comaps_trivial_intersection} is also an invertible comap.
	Then, by commutativity of the diagram (\ref{eq:molecin_maps_and_comaps}), we have $\varphi = \pb{\widehat{\psi}}$, proving that $\molecin{-}$ is full on isomorphism and that it reflects them.

	Let $c, d\colon P \to Q$ be comaps and suppose that $\pb{c} = \pb{d}$.
	For all $x \in P$,
	\begin{align*}
		\isocl{\invrs{c}\clset{c(x)} \incl P} & = \pb{c}\isocl{\clset{c(x)} \incl Q} = \pb{d}\isocl{\clset{c(x)} \incl Q} = \\ 
						      & = \isocl{\invrs{d}\clset{c(x)} \incl P},
	\end{align*}
	so $d(x) \in \clset{c(x)}$, that is, $d(x) \leq c(x)$.
	Dually, we prove $c(x) \leq d(x)$, so $c(x) = d(x)$.
	We conclude that $c = d$ and that $\molecin{\pb{-}}$ is faithful.
\end{proof}

\begin{exm}[A strict functor which is neither the image of a map nor of a comap] \index[counterex]{A strict functor which is neither the image of a map nor of a comap} \label{exm:not_map_nor_comap}
	Consider the 2\nbd dimensional molecules $\globe{2}$ and $\globe{2} \cp{0} \globe{2}$, as in the introduction to Chapter \ref{chap:molecules}.
	There is a strict functor $f\colon \molecin{\globe{2}} \to \molecin{(\globe{2} \cp{0} \globe{2})}$ which is uniquely determined by the assignment
	\[
		\isocl{\idd{\globe{2}}} \mapsto \isocl{\idd{(\globe{2} \cp{0} \globe{2})}}.
	\]
	We claim that there is no map $p\colon \globe{2} \to \globe{2} \cp{0} \globe{2}$ such that $f = \molecin{p}$.
	Since $\globe{2}$ is an atom, by Lemma 
	\ref{lem:pushforward_of_atom} $\pfw{p}\globe{2}$ is an atom for all maps $p$, so $\pfw{p}\idd{\globe{2}}$ cannot be isomorphic to $\idd{(\globe{2} \cp{0} \globe{2})}$, whose domain is not an atom.

	Likewise, there is no comap $c\colon \globe{2} \cp{0} \globe{2} \to \globe{2}$ such that $f = \molecin{\pb{c}}$:
	by Lemma \ref{lem:comaps_inverse_image_and_roundness}, since $\globe{2}$ is round, the domain of any comap whose codomain is $\globe{2}$ must be round, and $\globe{2} \cp{0} \globe{2}$ is not round.

	On the other hand, let $U$ be the oriented face poset of
\begin{equation} \label{eq:between_map_and_comap}
		\begin{tikzcd}[sep=small]
	&& {{\scriptstyle 3}\;\bullet} \\
	{{\scriptstyle 0}\;\bullet} &&& {{\scriptstyle 2}\;\bullet} \\
	& {{\scriptstyle 1}\;\bullet}
	\arrow["0"', curve={height=6pt}, from=2-1, to=3-2]
	\arrow["2", from=3-2, to=1-3]
	\arrow[""{name=0, anchor=center, inner sep=0}, "3", curve={height=-12pt}, from=2-1, to=1-3]
	\arrow[""{name=1, anchor=center, inner sep=0}, "1"', curve={height=12pt}, from=3-2, to=2-4]
	\arrow["4", curve={height=-6pt}, from=1-3, to=2-4]
	\arrow["0"', curve={height=-6pt}, shorten >=7pt, Rightarrow, from=3-2, to=0]
	\arrow["1"', curve={height=6pt}, shorten <=7pt, Rightarrow, from=1, to=1-3]
\end{tikzcd}\;,\end{equation}
	which is a 2\nbd dimensional round molecule.
	There is a surjective map $p\colon U \to \globe{2} \cp{0} \globe{2}$ which ``collapses'' the cell $(1, 2)$, defined by
\begin{align*}
	(0, 0) & \mapsto (0, 0),
	       & (0, 1), (0, 3) & \mapsto (0, 1),
	       & (0, 2) & \mapsto (0, 2) \\ 
	(1, 0) & \mapsto (1, 0),
		       & (1, 1) & \mapsto (1, 1),
		       & (1, 2) & \mapsto (0, 1), \\ 
	(1, 3) & \mapsto (1, 2),
		       & (1, 4) & \mapsto (1, 3), \\
	(2, 0) & \mapsto (2, 0), 
		& (2, 1) & \mapsto (2, 1),
\end{align*}
	as well as a unique comap $c\colon U \to \globe{2}$ defined by
\begin{align*}
	(0, 0) & \mapsto (0, 0),
	       & (0, 1) & \mapsto (1, 0),
	       & (0, 2) & \mapsto (0, 2), \\
	(0, 3) & \mapsto (1, 1), \\
	(1, 0), (1, 1) & \mapsto (1, 0),
		       & (1, 2) & \mapsto (2, 0),
		       & (1, 3), (1, 4) & \mapsto (1, 1), \\ 
	(2, 0), (2, 1) & \mapsto (2, 0).
\end{align*}
	Then $f = \molecin{p} \after \molecin{\pb{c}}$.
	This leaves open the possibility that maps and comaps jointly generate the full subcategory of $\omegacat$ on strict $\omega$\nbd categories of the form $\molecin{P}$ where $P$ is a regular directed complex.
	We do not know a counterexample, but have not attempted a proof either.
\end{exm}

\begin{lem} \label{lem:chain_complex_boundary_of_molecule}
	Let $P$ be a regular directed complex and let $U \subseteq P$ be a molecule, $n \eqdef \dim{U} > 0$.
	Then, in the augmented chain complex $\freeab{P}$,
	\[
		\der \left( \sum_{x \in \grade{n}{U}} x \right) = \sum_{y \in \faces{}{+}U} y - \sum_{y \in \faces{}{-}U} y.
	\]
\end{lem}
\begin{proof}
	By linearity of $\der$, we have
	\[
		\der \left( \sum_{x \in \grade{n}{U}} x \right) = 
			\sum_{x \in \grade{n}{U}} \sum_{y \in \faces{}{+}x} y - \sum_{x \in \grade{n}{U}} \sum_{y \in \faces{}{-}x} y.
	\]
	Let $y \in \grade{n-1}{U}$.
	By Corollary \ref{cor:codimension_1_elements}, if $y \in \maxel{U}$, then $y$ does not appear in this sum, if $y \in \faces{}{\alpha}U \setminus \faces{}{-\alpha}U$, then $y$ appears exactly once with sign $\alpha$, and if $y \not\in \faces{}{}U$, then $y$ appears exactly twice with opposite signs, so the two cancel out.
	Thus we are left with 
	\[
		\sum_{y \in \faces{}{+}U} y - \sum_{y \in \faces{}{-}U} y,
	\]
	as claimed.
\end{proof}

\begin{prop} \label{prop:chain_complex_is_functorial_on_comaps}
	Let $c\colon P \to Q$ be a comap of regular directed complexes.
	Then the sequence of homomorphisms
	\begin{align*}
		\freeab{\grade{n}{\pb{c}}}\colon \freeab{\grade{n}{Q}} & \to \freeab{\grade{n}{P}}, \\
		y \in \grade{n}{Q} & \mapsto 
		\sum_{x \in \grade{n}{\invrs{c}(y)}} x
	\end{align*}
	is a homomorphism $\freeab{\pb{c}}\colon \freeab{Q} \to \freeab{P}$ of augmented chain complexes.
	This assignment determines a functor $\freeab{\pb{-}}\colon \opp{\rdcpxcomap} \to \chaug$, such that the diagram of functors
\begin{equation} \label{eq:comap_functor_to_chaug}
\begin{tikzcd}
	{\opp{\rdcpxiso}} & \rdcpxiso && \rdcpxmap \\
	{\opp{\rdcpxcomap}} &&& \chaug
	\arrow["{\invrs{(-)}}", from=1-1, to=1-2]
	\arrow[hook, from=1-2, to=1-4]
	\arrow[hook', from=1-1, to=2-1]
	\arrow["{\freeab{\pb{-}}}", from=2-1, to=2-4]
	\arrow["{\freeab{-}}", from=1-4, to=2-4]
\end{tikzcd}
\end{equation}
	commutes.
\end{prop}
\begin{proof}
	Let $y \in \grade{0}{Q}$.
	Then by Proposition \ref{prop:comaps_inverse_image_preserves_molecules} $\invrs{c}{\clset{y}}$ is a 0\nbd dimensional molecule, so by Lemma \ref{lem:only_0_molecule} it is equal to $\set{x}$ for some $x \in \grade{0}{P}$.
	It follows that $\freeab{\grade{0}{\pb{c}}}(y) = x$, and $\eau (y) = \eau (x) = 1$, so $\eau \after \freeab{\grade{0}{\pb{c}}} = \eau$.

	Next, let $y \in \grade{n}{Q}$ for $n > 0$.
	Then $\invrs{c}{\clset{y}}$ is an $n$\nbd dimensional molecule.
	Moreover, as special cases of 
	Lemma \ref{lem:comap_inverse_image_boundaries}, for all $\alpha \in \set{+, -}$,
	\[ 
		\grade{n}{(\invrs{c}{\clset{y}})} = \grade{n}{\invrs{c}(y)}, \quad \quad
		\faces{}{\alpha}\invrs{c}\clset{y} = \grade{n-1}{(\invrs{c}\faces{}{\alpha}y)}.
	\]
	By Lemma \ref{lem:chain_complex_boundary_of_molecule},
	\begin{align*}
		\der(\freeab{\grade{n}{\pb{c}}}(y)) & = \quad  
		\; \sum_{\mathclap{z \in \faces{}{+}\invrs{c}\clset{y}}} \; z \quad -
		\quad \; \sum_{\mathclap{z \in \faces{}{-}\invrs{c}\clset{y}}} \; z \quad = \\
		& = \quad 
		\;\; \sum_{\mathclap{z \in \grade{n-1}{({\invrs{c}\faces{}{+}y)}}}}\; z \quad -
		\quad \;\; \sum_{\mathclap{z \in \grade{n-1}{(\invrs{c}\faces{}{-}y)}}} \; z \quad = \\
		& = \sum_{z' \in \faces{}{+}y} \freeab{\grade{n-1}{\pb{c}}}(z') -
		\sum_{z' \in \faces{}{-}y} \freeab{\grade{n-1}{\pb{c}}}(z') = 
		\freeab{\grade{n-1}{\pb{c}}}( \der (y) ).
	\end{align*}
	This proves that $\der \after \freeab{\grade{n}{\pb{c}}} = \freeab{\grade{n-1}{\pb{c}}} \after \der$, so $(\freeab{\grade{n}{\pb{c}}})_{n \in \mathbb{N}}$ is a homomorphism of augmented chain complexes.
	Functoriality and commutativity of (\ref{eq:comap_functor_to_chaug}) are both straightforward checks on the definitions.
\end{proof}

\clearpage
\thispagestyle{empty}

%% file: constructions.tex
\chapter{Constructions and operations} \label{chap:constructions}
\thispagestyle{firstpage}

\begin{guide}
	One of the most remarkable aspects of higher-categorical combinatorics is just how many topological constructions have ``directed'' analogues, and that these may rightfully be seen as \emph{more refined} version of their undirected counterparts, due to the combinatorial possibilities afforded by direction-reversing dualities, so much that one could be justified in thinking that the usual, undirected topology is but a shadow of a more fundamental directed topology.
	At the same time, these aspects are among the least developed and understood, in part due to the relative immaturity of higher category theory, and the fact that the most commonly used models of higher categories, at the time of writing, seem ill-suited to these constructions.

	The clearest example is \cemph{Gray products}, also known as Crans--Gray products, Gray tensor products, or just tensor products, which are the directed counterpart of cartesian products of spaces.
	Their importance is evident due to their role in representing higher-dimensional lax and oplax cells between functors --- similar to the role of products, in particular \cemph{cylinders}, in representing higher homotopies between spaces --- but they are also difficult to work with in most models, with the exception of those based on \emph{cubical} combinatorics.
	A similar discourse applies to \cemph{joins}, which, particularly in the special case of \cemph{cones}, play a fundamental role in the theory of lax and oplax limits and colimits, yet are only a natural construction in \emph{simplicial} models of higher categories.
	To complete the picture, \cemph{suspensions} are only really adapted to \emph{globular} combinatorics.

	The reason, as we will see in Chapter \ref{chap:special}, is that particular classes of diagram shapes such as oriented cubes, oriented simplices, and globes are only closed under one each of these constructions, and are in fact \emph{characterised} by this property. 
	The good news is that regular directed complexes are closed under all these constructions.
	In fact, Gray products and joins of regular directed complexes have a remarkably simple definition, in contrast to, for example, strict $\omega$\nbd categories, suggesting that this is the ``right'' framework to study these constructions. 

	What is less simple, for example, is proving that \emph{molecules} are closed under Gray products (Proposition \ref{prop:gray_product_of_molecules}).
	The reason is that, while Gray products preserve the pushout diagrams that define a pasting of molecules, the resulting pushout diagram does \emph{not} define a pasting except in trivial circumstances.
	To address this, we introduce a notion of \cemph{generalised pasting}, which provides a criterion for when gluing two molecules at a \emph{portion} of their boundaries produces a molecule.
	This is an original technical innovation and provides a neater understanding of why ``pasting diagrams are closed under Gray products'' than the earlier literature.

	Subsequently, we also prove that molecules are closed under suspensions (Proposition 
	\ref{prop:suspension_of_molecules}) and joins (Proposition \ref{prop:join_of_molecules}).
	We complete the survey with a section on direction-reversing \cemph{duals}, which are the only operation that is ``intrinsically directed'' since it operates only on the orientation and not on the underlying poset.

	For each of these constructions, we also prove compatibility with an analogous construction on augmented chain complexes.
	Compatibility with strict $\omega$\nbd categories is a thornier business, since the definition of Gray products and joins of strict $\omega$\nbd categories relies on Steiner's theory of directed chain complexes, so we postpone any comparison --- except for the easier case of suspensions and duals --- to Chapter \ref{chap:steiner}.
\end{guide}


\section{Generalised pastings} \label{sec:generalised_pastings}

\begin{guide}
	In this section, we define generalised pastings, and prove that they preserve molecules (Lemma \ref{lem:generalised_pasting}.
	We then focus on a special case, \cemph{pasting at a submolecule}, which will feature in Chapter \ref{chap:special} as one of the constructors for the class of positive opetopes.
\end{guide}

\begin{dfn}[Generalised pasting of molecules] \index{molecule!generalised pasting} \index{$U \gencp{k} V$} \index{pasting!generalised}
	Let $U, V$ be molecules, $k \in \mathbb{N}$, and let
	\[\begin{tikzcd}
		{U \cap V} && V \\
		U && {U \cup V}
		\arrow[hook, from=1-1, to=1-3]
		\arrow[hook', from=1-1, to=2-1]
		\arrow[hook, from=2-1, to=2-3]
		\arrow[hook', from=1-3, to=2-3]
		\arrow["\lrcorner"{anchor=center, pos=0.125, rotate=180}, draw=none, from=2-3, to=1-1]
	\end{tikzcd}\]
	be a pushout diagram of inclusions.
	We say that $U \cup V$ is a \emph{generalised pasting of $U$ and $V$ at the $k$\nbd boundary}, and write $U \gencp{k} V$ for $U \cup V$, if
	\begin{enumerate}
		\item $U \cap V \submol \bound{k}{+}U$ and $U \cap V \submol \bound{k}{-}V$,
		\item $\bound{k}{-}(U \cup V)$ and $\bound{k}{+}(U \cup V)$ are molecules,
		\item $\bound{k}{-}U \submol \bound{k}{-}(U \cup V)$ and $\bound{k}{+}V \submol \bound{k}{+}(U \cup V)$.
	\end{enumerate}
\end{dfn}

\begin{rmk}
	Unlike $U \cp{k} V$, which is uniquely defined up to unique isomorphism given molecules $U$, $V$ and $k$, the generalised pasting $U \gencp{k} V$ is \emph{not} uniquely defined by these data, so $- \gencp{k} -$ should not be read as an operation; one needs the entire span $U \cap V \incl U, V$ in order to specify it.
	This is why we speak of \emph{a} generalised pasting.
\end{rmk}

\begin{lem} \label{lem:pasting_is_generalised_pasting_in_all_higher_dim}
	Let $U, V$ be molecules, $k \in \mathbb{N}$, and suppose $U \cp{k} V$ is defined.
	Then, for all $m \geq k$,
	\[
		U \cp{k} V = U \gencp{m} V.
	\]
\end{lem}
\begin{proof}
	Since $U \cp{k} V$ is a molecule, all its boundaries are molecules.
	Moreover, for all $m \geq k$,
	\begin{enumerate}
		\item $U \cap V = \bound{k}{+}U \submol \bound{m}{+}U$ and $U \cap V = \bound{k}{-}V \submol \bound{m}{-}V$,
		\item $\bound{k}{-}U = \bound{k}{-}(U \cp{k} V)$, $\bound{k}{+}V = \bound{k}{+}(U \cp{k} V)$, while for all $m > k$, $\alpha \in \set{+, -}$ we have $\bound{m}{\alpha}U, \bound{m}{\alpha}V \submol \bound{m}{\alpha}U \cp{k} \bound{m}{\alpha}V = \bound{m}{\alpha}(U \cp{k} V)$,
	\end{enumerate}
	and we conclude.
\end{proof}

\begin{lem} \label{lem:generalised_pasting}
	Let $U, V$ be molecules, $k \in \mathbb{N}$, and let $U \gencp{k} V$ be a generalised pasting of $U$ and $V$.
	Then
	\begin{enumerate}
		\item $U \cup \bound{k}{-}(U \gencp{k} V)$ and $V \cup \bound{k}{+}(U \gencp{k} V)$ are molecules,
		\item $U \submol U \cup \bound{k}{-}(U \gencp{k} V)$ and $V \submol V \cup \bound{k}{+}(U \gencp{k} V)$,
		\item $U \gencp{k} V$ is a molecule isomorphic to 
			\[
				(U \cup \bound{k}{-}(U \gencp{k} V)) \cp{k} (V \cup \bound{k}{+}(U \gencp{k} V)).
			\]
	\end{enumerate}
\end{lem}
\begin{proof}
Consider the pushout diagrams
\[\begin{tikzcd}
	\bound{k}{-}U && \bound{k}{-}(U \gencp{k} V) \\
	U && U \cup \bound{k}{-}(U \gencp{k} V)
	\arrow[hook, from=1-1, to=1-3]
	\arrow[hook', from=1-1, to=2-1]
	\arrow[hook, from=2-1, to=2-3]
	\arrow[hook', from=1-3, to=2-3]
	\arrow["\lrcorner"{anchor=center, pos=0.125, rotate=180}, draw=none, from=2-3, to=1-1]
\end{tikzcd}
\; \text{and} \;
\begin{tikzcd}
	\bound{k}{+}V && \bound{k}{+}(U \gencp{k} V) \\
	V && V \cup \bound{k}{+}(U \gencp{k} V)
	\arrow[hook, from=1-1, to=1-3]
	\arrow[hook', from=1-1, to=2-1]
	\arrow[hook, from=2-1, to=2-3]
	\arrow[hook', from=1-3, to=2-3]
	\arrow["\lrcorner"{anchor=center, pos=0.125, rotate=180}, draw=none, from=2-3, to=1-1]
\end{tikzcd}
\]
in $\ogpos$.
By assumption, they both satisfy the conditions of Lemma \ref{lem:submolecule_rewrite}, from which it follows that
\begin{enumerate}
	\item $U \cup \bound{k}{-}(U \gencp{k} V)$ and $V \cup \bound{k}{+}(U \gencp{k} V)$ are molecules,
	\item $U \submol (U \cup \bound{k}{-}(U \gencp{k} V))$ and $V \submol (V \cup \bound{k}{+}(U \gencp{k} V))$,
	\item $\bound{k}{-}(U \cup \bound{k}{-}(U \gencp{k} V)) = \bound{k}{-}(U \gencp{k} V)$,
	\item $\bound{k}{+}(V \cup \bound{k}{+}(U \gencp{k} V)) = \bound{k}{+}(U \gencp{k} V)$
	\item $\bound{k}{+}U \submol \bound{k}{+}(U \cup \bound{k}{-}(U \gencp{k} V))$ and $\bound{k}{-}V \submol \bound{k}{-}(V \cup \bound{k}{+}(U \gencp{k} V))$.
\end{enumerate}
We claim that 
\[
	\bound{k}{+}(U \cup \bound{k}{-}(U \gencp{k} V)) = \bound{k}{+}U \cup \bound{k}{-}V = \bound{k}{-}(V \cup \bound{k}{+}(U \gencp{k} V))
\]
as subsets of $U \gencp{k} V$.
By Lemma \ref{lem:faces_of_union}, $\faces{k}{+}(U \cup \bound{k}{-}(U \gencp{k} V))$ splits into
\[
	\faces{k}{+}U \cap \faces{k}{-}(U \gencp{k} V) + \faces{k}{+}U \setminus \bound{k}{-}(U \gencp{k} V) + \faces{k}{-}(U \gencp{k} V) \setminus U.
\]
Then
\[
	\faces{k}{+}U \cap \faces{k}{-}(U \gencp{k} V) = \faces{k}{+}U \cap \faces{k}{-}U = \grade{k}{(\maxel{U})}
\]
by Lemma \ref{lem:maximal_vs_faces}, while
\[
	\faces{k}{+}U \setminus \bound{k}{-}(U \gencp{k} V) = \faces{k}{+}U \setminus \faces{k}{-}U
\]
so the first two terms are jointly equal to $\faces{k}{+}U$.
Finally,
\[
	\faces{k}{-}(U \gencp{k} V) \setminus U = \faces{k}{-}V \setminus U.
\]
Since $\faces{k}{-}V \cap U = \grade{k}{(U \cap V)} \subseteq \faces{k}{+}U$, we have
\[
	\faces{k}{+}U + \faces{k}{-}V \setminus U = \faces{k}{+}U \cup \faces{k}{-}V.
\]
By Lemma \ref{lem:faces_of_union} again, for $j < k$, we have that $\grade{j}{(\maxel{(U \cup \bound{k}{-}(U \gencp{k} V))})}$ splits into
\[
	\grade{j}{(\maxel{U})} \cap \grade{j}{(\maxel{(U \gencp{k} V)})} + \grade{j}{(\maxel{U})} \setminus \bound{k}{-}(U \gencp{k} V) + \grade{j}{(\maxel{(U \gencp{k} V)})} \setminus U.
\]
Here
\[
	\grade{j}{(\maxel{U})} \cap \grade{j}{(\maxel{(U \gencp{k} V)})} = \grade{j}{(\maxel{U})} \cap \grade{j}{(\maxel{V})} + \grade{j}{(\maxel{U})} \setminus V,
\]
while
\[
	\grade{j}{(\maxel{U})} \setminus \bound{k}{-}(U \gencp{k} V) = \varnothing 
\]
since $\grade{j}{(\maxel{U})} \subseteq \bound{k}{-}U \subseteq \bound{k}{-}(U \gencp{k} V)$.
Finally,
\[
	\grade{j}{(\maxel{(U \gencp{k} V)})} \setminus U = \grade{j}{(\maxel{V})} \setminus U,
\]
and we conclude that $\grade{j}{(\maxel{(U \cup \bound{k}{-}(U \gencp{k} V))})} = \grade{j}{(\maxel{(U \gencp{k} V)})}$.
Since this is included in $\grade{j}{(\maxel{U})} \cup \grade{j}{(\maxel{V})}$, we have immediately
\[
	\bound{k}{+}(U \cup \bound{k}{-}(U \gencp{k} V)) \subseteq \bound{k}{+}U \cup \bound{k}{-}V.
\]
Conversely, we already know that $\bound{k}{+}U \subseteq \bound{k}{+}(U \cup \bound{k}{-}(U \gencp{k} V))$.
Suppose $x \in \bound{k}{-}V$.
Then there exists $y$ such that $x \leq y$ and $y \in \faces{k}{-}V$ or $y \in \grade{j}{(\maxel{V})}$ for some $j < k$.
In the first case, we are done.
In the second case, either $y \in \grade{j}{(\maxel{(U \gencp{k} V)})}$ and we are done, or there exists $z \in U$ such that $y < z$.
But then $x \in U \cap V \subseteq \bound{k}{+}U$, and we conclude that
\[
	\bound{k}{+}(U \cup \bound{k}{-}(U \gencp{k} V)) = \bound{k}{+}U \cup \bound{k}{-}V.
\]
A dual argument proves that
\[
	\bound{k}{-}(V \cup \bound{k}{+}(U \gencp{k} V)) = \bound{k}{+}U \cup \bound{k}{-}V.
\]
It follows that $U \gencp{k} V$ is isomorphic to the pasting
\[
	(U \cup \bound{k}{-}(U \gencp{k} V)) \cp{k} (V \cup \bound{k}{+}(U \gencp{k} V)),
\]
and we conclude.
\end{proof}

\begin{exm}[A generalised pasting] \label{exm:generalised_pasting}
	Let $U$ be the oriented face poset of (\ref{eq:between_map_and_comap}) from Example 
	\ref{exm:not_map_nor_comap}.
	Then $U$ is a generalised pasting of the form 
	\[
		\disk{1}{2} \gencp{1} \disk{2}{1}, 
	\]
	where $\disk{1}{2} \simeq \clset{(2, 1)}$ and $\disk{2}{1} \simeq \clset{(2, 0)}$.
\end{exm}

\begin{dfn}[Pasting at a submolecule] \index{molecule!pasting!at a submolecule} \index{$U \cpsub{k,W}{+} V, U \cpsub{k,W}{-} V$} \index{submolecule!pasting at} \index{pasting!at a submolecule}
	Let $U$, $V$ be molecules, $k \in \mathbb{N}$, $\alpha \in \set{+, -}$, and $W \submol \bound{k}{\alpha}V$ such that $W$ is isomorphic to $\bound{k}{-\alpha}U$.
	The \emph{pasting of $U$ at the submolecule $W \submol \bound{k}{\alpha}V$} is the oriented graded poset $U \cpsub{k,W}{\alpha}V$ obtained as the pushout
\[\begin{tikzcd}
	\bound{k}{-\alpha}U & W & V \\
	U && U \cpsub{k,W}{\alpha}V
	\arrow["\sim", hook, from=1-1, to=1-2]
	\arrow[hook, from=1-2, to=1-3]
	\arrow[hook', from=1-1, to=2-1]
	\arrow[hook, from=2-1, to=2-3]
	\arrow[hook', from=1-3, to=2-3]
	\arrow["\lrcorner"{anchor=center, pos=0.125, rotate=180}, draw=none, from=2-3, to=1-1]
\end{tikzcd}\]
	in $\ogpos$.
	We omit the index $k$ when $k = \dim{V} - 1$.
\end{dfn}

\begin{lem} \label{lem:pasting_at_submolecule}
	Let $U$, $V$ be molecules, $k \in \mathbb{N}$, $\alpha \in \set{+, -}$, and $W \submol \bound{k}{\alpha}V$ such that $U \cpsub{k,W}{\alpha}V$ is defined.
	Then
	\begin{enumerate}
		\item $U \cup \bound{k}{\alpha}V$ is a molecule,
		\item $U \cpsub{k,W}{\alpha}V$ is a molecule, isomorphic to
			\[
				\begin{cases}
					(U \cup \bound{k}{\alpha}V) \cp{k} V & \text{if $\alpha = -$,} \\
					V \cp{k} (U \cup \bound{k}{\alpha} V) & \text{if $\alpha = +$.}
				\end{cases}
			\]
	\end{enumerate}
\end{lem}
\begin{proof}
	Suppose that $\alpha = -$; the case $\alpha = +$ is dual.
	The pushout diagram
\[\begin{tikzcd}
	\bound{k}{+}U & W & \bound{k}{-}V \\
	U && U \cup \bound{k}{-}V
	\arrow["\sim", hook, from=1-1, to=1-2]
	\arrow[hook, from=1-2, to=1-3]
	\arrow[hook', from=1-1, to=2-1]
	\arrow[hook, from=2-1, to=2-3]
	\arrow[hook', from=1-3, to=2-3]
	\arrow["\lrcorner"{anchor=center, pos=0.125, rotate=180}, draw=none, from=2-3, to=1-1]
\end{tikzcd}\]
satisfies the conditions of Lemma \ref{lem:submolecule_rewrite}, so $U \cup \bound{k}{-}V$ is a molecule with
	\[
		\bound{k}{-}U \submol \bound{k}{-}(U \cup \bound{k}{-}V), \quad \quad
		\bound{k}{+}(U \cup \bound{k}{-}V) = \bound{k}{+}(\bound{k}{-}V) = \bound{k}{-}V.
	\]
	It follows that $(U \cup \bound{k}{-}V) \cp{k} V$ is defined and isomorphic to $U \cpsub{k,W}{-} V$.
\end{proof}

\begin{rmk}
	As a consequence of Lemma \ref{lem:pasting_at_submolecule}, $U \cpsub{k,W}{\alpha} V$ is, in particular, a generalised pasting at the $k$\nbd boundary.
\end{rmk}

\begin{lem} \label{lem:round_pasting_at_submolecule}
	Let $U$, $V$ be molecules, $\dim{U} = \dim{V}$, $\alpha \in \set{+, -}$, and $W \submol \bound{}{\alpha}V$ such that $U \cpsub{W}{\alpha} V$ is defined.
	If $V$ is round, then $U \cpsub{W}{\alpha} V$ is round.
\end{lem}
\begin{proof}
	Suppose without loss of generality that $\alpha = -$.
	By Lemma \ref{lem:pasting_at_submolecule},
	\begin{align*}
		\bound{}{+}(U \cpsub{W}{-} V) & = \bound{}{+}V, \\
		\bound{}{-}(U \cpsub{W}{-} V) & = \bound{}{-}(U \cup \bound{}{-}V).
	\end{align*}
	Since $U \cpsub{W}{-} V$ is globular, for all $k < \dim{V} - 1$ and $\beta \in \set{+, -}$,
	\[
		\bound{k}{\beta}(U \cpsub{W}{-} V) = \bound{k}{\beta}(\bound{}{-}(U \cpsub{W}{-} V)) = \bound{k}{\beta}(\bound{}{+}V) = \bound{k}{\beta}V,
	\]
	and by Corollary \ref{cor:boundary_of_union} and the fact that $U \cap V \subseteq \bound{}{-}V$,
	\[
		\bound{}{+}V \cap \bound{}{-}(U \cup \bound{}{-}V) \subseteq (\bound{}{+}V \cap \bound{}{-}U) \cup (\bound{}{+}V \cap \bound{}{-}V) = \bound{}{+}V \cap \bound{}{-}V.
	\]
	Roundness of $U \cpsub{W}{-} V$ then follows immediately from roundness of $V$.
\end{proof}

\begin{lem} \label{lem:boundary_of_round_pasting_at_submolecule}
	Let $U$, $V$ be molecules, $\dim{U} = \dim{V}$, $\alpha \in \set{+, -}$, and $W \submol \bound{}{\alpha}V$ such that $U \cpsub{W}{\alpha} V$ is defined.
	If $U$ is round, then $\bound{}{\alpha}(U \cpsub{W}{\alpha} V)$ is isomorphic to $\subs{\bound{}{\alpha}V}{\bound{}{\alpha}U}{W}$.
\end{lem}
\begin{proof}
	By Lemma \ref{lem:pasting_at_submolecule}, $\bound{}{\alpha}(U \cpsub{W}{\alpha}V)$ is equal to $\bound{}{\alpha}(U \cup \bound{}{\alpha}V)$.
	Since $U$ is round, $\compos{U}$ is defined and $\bound{}{\beta}\compos{U}$ is isomorphic to $\bound{}{\beta}U$ for all $\beta \in \set{+, -}$.
	Then the pushout
\[\begin{tikzcd}
	\bound{}{-\alpha}\compos{U} & W & \bound{}{\alpha}V \\
	\compos{U} && \compos{U} \cup \bound{}{\alpha}V
	\arrow["\sim", hook, from=1-1, to=1-2]
	\arrow[hook, from=1-2, to=1-3]
	\arrow[hook', from=1-1, to=2-1]
	\arrow[hook, from=2-1, to=2-3]
	\arrow[hook', from=1-3, to=2-3]
	\arrow["\lrcorner"{anchor=center, pos=0.125, rotate=180}, draw=none, from=2-3, to=1-1]
\end{tikzcd}\]
	is defined, and $\bound{}{\alpha}(U \cup \bound{}{\alpha}V)$ is isomorphic to $\bound{}{\alpha}(\compos{U} \cup \bound{}{\alpha}V)$.
	By definition, up to a duality, this is $\subs{\bound{}{\alpha}V}{\bound{}{\alpha}\compos{U}}{W}$, which is isomorphic to $\subs{\bound{}{\alpha}V}{\bound{}{\alpha}U}{W}$.
\end{proof}

\begin{exm}[A pasting at a submolecule]
	The generalised pasting of Example \ref{exm:generalised_pasting} is \emph{not} a pasting at a submolecule.
	On the other hand, let $U$ be the oriented face poset of
\[\begin{tikzcd}[column sep=small]
	& \bullet \\
	\bullet &&& \bullet
	\arrow[""{name=0, anchor=center, inner sep=0}, "0"', curve={height=18pt}, from=2-1, to=2-4]
	\arrow["1", curve={height=-6pt}, from=2-1, to=1-2]
	\arrow[""{name=1, anchor=center, inner sep=0}, "2"', curve={height=12pt}, from=1-2, to=2-4]
	\arrow[""{name=2, anchor=center, inner sep=0}, "3", curve={height=-12pt}, from=1-2, to=2-4]
	\arrow["0", curve={height=-6pt}, shorten <=6pt, Rightarrow, from=0, to=1-2]
	\arrow["1"', shorten <=3pt, shorten >=3pt, Rightarrow, from=1, to=2]
\end{tikzcd}\]
	which is a 2\nbd dimensional molecule.
	Then $U$ is isomorphic to $\globe{2} \cpsub{W}{+} \disk{1}{2}$, where $W \eqdef \clset{(1, 2)} \submol \bound{}{+}\disk{1}{2} \simeq \bound{}{+}(2, 0)$.
\end{exm}


\section{Gray products} \label{sec:gray}

\begin{guide}
	In this section, we define the Gray product of oriented graded posets.
	We give expressions for the boundaries of a Gray product in terms of the boundaries of its factors, and work our way towards the proof that the Gray product of two molecules is a molecule (Proposition 
	\ref{prop:gray_product_of_molecules}), which implies that the Gray product of two regular directed complexes is a regular directed complex (Corollary \ref{cor:gray_products_of_rdcpx}).

	We show that Gray products are part of monoidal structures on $\ogpos$, $\rdcpxmap$, and $\rdcpxcomap$, compatibly on their overlaps.
	Finally, we recall the definition of the tensor product of augmented chain complexes, and prove that the $\freeab{-}$ construction maps Gray products to tensor products.
\end{guide}

\begin{dfn}[Gray product of oriented graded posets] \index{oriented graded poset!Gray product} \index{regular directed complex!Gray product} \index{$P \gray Q$} \index{product!Gray|see {Gray product}} \index{Gray product!of oriented graded posets}
	Let $P, Q$ be oriented graded posets.
	The \emph{Gray product} of $P$ and $Q$ is the oriented graded poset $P \gray Q$ whose
	\begin{itemize}
		\item underlying graded poset is the product $P \times Q$ of the underlying posets,
		\item orientation is defined, for all $(x, y) \in P \times Q$ and $\alpha \in \set{+, -}$, by
			\[
				\faces{}{\alpha}(x, y) = \faces{}{\alpha}x \times \set{y} + \set{x} \times \faces{}{(-)^{\dim{x}}\alpha}y.
			\]
	\end{itemize}
\end{dfn}

\begin{exm}[A non-symmetric Gray product of molecules] \index[counterex]{A non-symmetric Gray product of molecules} \label{exm:nonsymmetry_gray}
	As soon as both factors are at least 1\nbd dimensional, the Gray product ceases to be symmetric, in general.
	For example, $\thearrow{} \gray \thearrow{2}$ and $\thearrow{2} \gray \thearrow{}$ are isomorphic to the oriented face posets of
	\[\begin{tikzcd}[sep=small]
	& \bullet \\
	\bullet && \bullet \\
	& \bullet && \bullet \\
	&& \bullet
	\arrow[from=4-3, to=3-4]
	\arrow[from=2-1, to=1-2]
	\arrow[from=3-2, to=2-3]
	\arrow[shorten <=3pt, shorten >=3pt, Rightarrow, from=4-3, to=2-3]
	\arrow[shorten <=6pt, shorten >=3pt, Rightarrow, from=3-2, to=1-2]
	\arrow[from=2-1, to=3-2]
	\arrow[from=3-2, to=4-3]
	\arrow[from=2-3, to=3-4]
	\arrow[from=1-2, to=2-3]
\end{tikzcd} \quad \quad \text{and} \quad \quad
\begin{tikzcd}[sep=small]
	&& \bullet \\
	& \bullet && \bullet \\
	\bullet && \bullet \\
	& \bullet
	\arrow[from=3-1, to=4-2]
	\arrow[from=2-2, to=3-3]
	\arrow[from=3-1, to=2-2]
	\arrow[from=2-2, to=1-3]
	\arrow[from=3-3, to=2-4]
	\arrow[from=4-2, to=3-3]
	\arrow[from=1-3, to=2-4]
	\arrow[shorten <=3pt, shorten >=3pt, Rightarrow, from=4-2, to=2-2]
	\arrow[shorten <=3pt, shorten >=3pt, Rightarrow, from=3-3, to=1-3]
\end{tikzcd}\]
	respectively, which are not isomorphic to each other.
\end{exm}

\begin{prop} \label{prop:gray_products_monoidal_ogpos}
	There is a unique monoidal structure $(\ogpos, \gray, 1)$ on $\ogpos$ such that
	\begin{enumerate}
		\item $P \gray Q$ is the Gray product of $P$ and $Q$,
		\item $\fun{U}\colon (\ogpos, \gray, 1) \to (\posclos, \times, 1)$ is a strict monoidal functor.
	\end{enumerate}
\end{prop}
\begin{proof}
	The requirements fix the monoidal structure uniquely, so it suffices to show that it is well-defined.
	
	First of all, suppose $f\colon P \to P'$ and $g\colon Q \to Q'$ are morphisms of oriented graded posets, and let $(x, y) \in P \gray Q$ and $\alpha, \beta \in \set{+, -}$.
	Then
	\begin{enumerate}
		\item $f$ maps $\faces{}{\alpha}x$ bijectively onto $\faces{}{\alpha}f(x)$, 
		\item $g$ maps $\faces{}{\beta}y$ bijectively onto $\faces{}{\beta}g(y)$,
		\item $\dim{x} = \dim{f(x)}$,
	\end{enumerate}
	so $f \times g$ determines a bijection between
	\[
		\faces{}{\alpha}(x, y) = \faces{}{\alpha}x \times \set{y} + \set{x} \times \faces{}{(-)^{\dim{x}}\alpha}y
	\]
	and
	\[
		\faces{}{\alpha}(f(x), g(y)) = \faces{}{\alpha}f(x) \times \set{g(y)} + \set{f(x)} \times \faces{}{(-)^{\dim{f(x)}}\alpha}g(y).
	\]
	This proves that $f \times g$ lifts to a morphism $f \gray g\colon P \gray Q \to P' \gray Q'$ of oriented graded posets.

	It then suffices to show that the structural isomorphisms of $(\posclos, \times, 1)$ lift to isomorphisms of oriented graded posets.
	Let $((x, y), z) \in (P \gray Q) \gray R$ and $\alpha \in \set{+, -}$.
	Then
	\[
	\faces{}{\alpha}((x, y), z) = \faces{}{\alpha}(x, y) \times \set{z} + \set{(x, y)} \times \faces{}{(-)^{\dim{(x, y)}}\alpha}z 
\]
	which is equal to
	\[
		(\faces{}{\alpha}x \times \set{y}) \times \set{z} + (\set{x} \times \faces{}{(-)^{\dim{x}}\alpha}y ) \times \set{z} + (\set{x} \times \set{y}) \times \faces{}{(-)^{\dim{x} + \dim{y}}\alpha}z,
	\]
	using Lemma \ref{lem:product_of_posets_faces} to rewrite $\dim{(x, y)}$ as $\dim{x} + \dim{y}$.
	The associator $(P \times Q) \times R \iso P \times (Q \times R)$ maps this to
	\[
		\faces{}{\alpha}x \times (\set{y} \times \set{z}) + \set{x} \times (\faces{}{(-)^{\dim{x}}\alpha}y \times \set{z}) + \set{x} \times (\set{y} \times \faces{}{(-)^{\dim{x} + \dim{y}}\alpha}z)
	\]
	which is equal to
	\[
		\faces{}{\alpha}x \times \set{(y, z)} + \set{x} \times \faces{}{(-)^{\dim{x}}}(y, z) = \faces{}{\alpha}(x, (y, z)).
	\]
	Thus the associator lifts to an isomorphism of oriented graded posets.
	Finally, let $(x, *) \in P \gray 1$.
	Since $\faces{}{}* = \varnothing$,
	\[
		\faces{}{\alpha}(x, *) = \faces{}{\alpha}x \times \set{*}
	\]
	which the right unitor $P \times 1 \iso P$ maps to $\faces{}{\alpha}x$.
	Thus the right unitor lifts to an isomorphism of oriented graded posets.
	A similar proof, using the fact that $(-)^{\dim{*}} = (-)^0 = +$, shows the same for the left unitor, and we conclude.
\end{proof}

\begin{prop} \label{prop:gray_product_augm_oriented_thin}
	Let $P$, $Q$ be oriented graded posets such that $\augm{P}$, $\augm{Q}$ are oriented thin.
	Then $\augm{(P \gray Q)}$ is oriented thin.
\end{prop}
\begin{proof}
	Let $z, w$ be elements of $\augm{(P \gray Q)}$ such that $z \leq w$ and $\codim{z}{w} = 2$.
	If $z = \bot$, then $w = \augm{(x, y)}$ for a pair of elements $x \in P$, $y \in Q$ with $\dim{x} + \dim{y} = 1$.
	If $\dim{x} = 1$ and $\dim{y} = 0$, then $[z, w]$ is isomorphic to the interval $[\bot, x]$ in $\augm{P}$, and if $\dim{x} = 0$ and $\dim{y} = 1$, then $[z, w]$ is isomorphic to the interval $[\bot, y]$ in $\augm{Q}$, so we conclude by oriented thinness of $\augm{P}$ and $\augm{Q}$.

	Suppose $\dim{z} > 0$.
	Then $z = \augm{(x', y')}$ and $w = \augm{(x, y)}$ for some elements $(x', y'), (x, y) \in P \gray Q$ with $\codim{x',y'}{(x, y)} = 2$, and the interval $[z, w]$ is isomorphic to the interval $[(x', y'), (x, y)]$ in $P \gray Q$.
	Let $k \eqdef \dim{x}$.
	Then either $y = y'$ and $\codim{x'}{x} = 2$, in which case the interval $[(x', y'), (x, y)]$ has the form
\[\begin{tikzcd}[column sep=small]
	& {(x, y)} \\
	{(z_1, y)} && {(z_2, y)} \\
	& {(x', y)}
	\arrow["{\alpha}"', from=1-2, to=2-1]
	\arrow["{\beta}", from=1-2, to=2-3]
	\arrow["{\gamma}"', from=2-1, to=3-2]
	\arrow["{-\alpha\beta\gamma}", from=2-3, to=3-2]
\end{tikzcd}\]
	by oriented thinness of $\augm{P}$, or $x = x'$ and $\codim{y'}{y} = 2$, in which case the interval has the form
\[\begin{tikzcd}[column sep=small]
	& {(x, y)} \\
	{(x, z'_1)} && {(x, z'_2)} \\
	& {(x, y')}
	\arrow["{(-)^{k}\alpha}"', from=1-2, to=2-1]
	\arrow["{(-)^{k}\beta}", from=1-2, to=2-3]
	\arrow["{(-)^{k}\gamma}"', from=2-1, to=3-2]
	\arrow["{-(-)^{k}\alpha\beta\gamma}", from=2-3, to=3-2]
\end{tikzcd}\]
	by oriented thinness of $\augm{Q}$, or $\codim{x'}{x} = \codim{y'}{y} = 1$, in which case the interval has the form
\[\begin{tikzcd}[column sep=small]
	& {(x, y)} \\
	{(x', y)} && {(x, y')} \\
	& {(x', y')}
	\arrow["{\alpha}"', from=1-2, to=2-1]
	\arrow["{(-)^{k}\beta}", from=1-2, to=2-3]
	\arrow["{(-)^{k-1}\beta}"', from=2-1, to=3-2]
	\arrow["{\alpha}", from=2-3, to=3-2]
\end{tikzcd}\]
	for the unique $\alpha, \beta \in \set{+, -}$ such that $x' \in \faces{}{\alpha}x$ and $y' \in \faces{}{\beta}y$.
	The defining condition of oriented thinness holds in all three cases.
\end{proof}

\begin{cor} \label{cor:gray_monoidal_on_otgpos}
	The monoidal structure $(\ogposbot, \augm{\gray}, \augm{1})$ obtained by transporting $(\ogpos, \gray, 1)$ along the equivalence $\augm{(-)}$ restricts to a monoidal structure on $\otgpos$.
\end{cor}

\begin{lem} \label{lem:gray_product_preserves_inclusions}
	Let $\imath\colon P \incl P'$ and $j\colon Q \incl Q'$ be inclusions of oriented graded posets.
	Then $\imath \gray j\colon P \gray Q \to P' \gray Q'$ is an inclusion.
\end{lem}
\begin{proof}
	Follows from Lemma \ref{lem:preservation_of_closed_embeddings}.
\end{proof}

\begin{rmk}
	Since the Gray product preserves inclusions, we can let it act on closed subsets $U \subseteq P$, $V \subseteq Q$, producing closed subsets $U \gray V \subseteq P \gray Q$.
\end{rmk}

\begin{lem} \label{lem:gray_products_preserve_some_colimits}
	Let $P$ be an oriented graded poset, let $\fun{F}$ be a diagram of inclusions in $\ogpos$, and let $\gamma$ be a colimit cone under $\fun{F}$ whose components are all inclusions.
	Then 
	\begin{enumerate}
		\item $P \gray \gamma$ is a colimit cone under $P \gray \fun{F}$,
		\item $\gamma \gray P$ is a colimit cone under $\fun{F} \gray P$.
	\end{enumerate}
\end{lem}
\begin{proof}
	By Lemma \ref{lem:reflected_colimits_in_ogpos}, $\fun{U}\gamma$ is a colimit cone in $\posclos$ whose components are all closed embeddings.
	Now, $P \gray \gamma$ has an underlying cone $\fun{U}P \times \fun{U}\gamma$.
	By Lemma \ref{lem:product_preserves_colimits}, $\fun{U}P \times \fun{U}\gamma$ is a colimit cone in $\posclos$ whose components are all closed embeddings by Lemma \ref{lem:preservation_of_closed_embeddings}.
	By Lemma \ref{lem:reflected_colimits_in_ogpos}, we conclude that $P \gray \gamma$ is a colimit cone in $\ogpos$ whose components are all inclusions.
	A symmetrical argument proves the same for $\gamma \gray P$.
\end{proof}

\begin{lem} \label{lem:faces_of_gray_product}
	Let $U, V$ be closed subsets of oriented graded posets, $n \in \mathbb{N}$, and $\alpha \in \set{+, -}$.
	Then
	\begin{enumerate}
		\item $\faces{n}{\alpha}(U \gray V) = \displaystyle\sum_{k = 0}^n \faces{k}{\alpha}U \times \faces{n-k}{(-)^k \alpha}V$,
		\item $\grade{n}{(\maxel{(U \gray V)})} = \displaystyle\sum_{k = 0}^n \grade{k}{(\maxel{U})} \times \grade{n-k}{(\maxel{V})}$.
	\end{enumerate}
\end{lem}
\begin{proof}
	Let $(x, y) \in \grade{n}{U}$.
	Then $\dim{(x, y)} = \dim{x} + \dim{y} = n$, so letting $k \eqdef \dim{x}$, we have $\dim{y} = n-k$.
	Now,
	\[
		\cofaces{}{-\alpha}(x, y) = \cofaces{}{-\alpha}x \times \set{y} + \set{x} \times \cofaces{}{-(-)^{k}\alpha}y,
	\]
	so $\cofaces{}{-\alpha}(x, y) \cap (U \gray V) = \varnothing$ if and only if
	\[
		\cofaces{}{-\alpha}x \cap U = \varnothing \quad \text{and} \quad \cofaces{}{-(-)^{k} \alpha}y \cap V = \varnothing.
	\]
	It follows that $(x, y) \in \faces{n}{\alpha}(U \gray V)$ if and only if $(x, y) \in \faces{k}{\alpha}U \times \faces{n-k}{(-)^k \alpha}V$.
	The proof for maximal elements is analogous, and simpler.
\end{proof}

\begin{cor} \label{cor:boundaries_of_gray_product}
	Let $U, V$ be closed subsets of oriented graded posets, $n \in \mathbb{N}$, and $\alpha \in \set{+, -}$.
	Then
	\begin{equation} \label{eq:boundaries_of_gray}
		\bound{n}{\alpha}(U \gray V) = \bigcup_{k=0}^n \bound{k}{\alpha}U \gray \bound{n-k}{(-)^k \alpha}V.
	\end{equation}
\end{cor}

\begin{rmk}
	Note that we can safely let the union range over all $k \in \mathbb{N}$, since for $k < 0$ or $k > n$ one of the factors is empty.
\end{rmk}

\begin{lem} \label{lem:gray_preserves_globular}
	Let $U$, $V$ be globular oriented graded posets.
	Then $U \gray V$ is globular.
\end{lem}
\begin{proof}
	Let $k, n \in \mathbb{N}$ and $\alpha \in \set{+, -}$ with $k < n$.
	Then, for all $j \in \set{0, \ldots, k}$, we have
	\[
		\bound{k-j}{(-)^j\alpha}V = \bound{k-j}{(-)^j\alpha}(\bound{n-j}{(-)^j\alpha}V) \subseteq \bound{n-j}{(-)^j\alpha}V,
	\]
	and consequently, using equation (\ref{eq:boundaries_of_gray}), $\bound{k}{\alpha}(U \gray V) \subseteq \bound{n}{\alpha}(U \gray V)$.
	It follows from Lemma \ref{lem:boundary_included_in_subset} that
	\[
		\bound{k}{\alpha}(U \gray V) \subseteq \bound{k}{\alpha}(\bound{n}{\alpha}(U \gray V)).
	\]
	Conversely, let $(x, y)$ be a maximal element in $\bound{k}{\alpha}(\bound{n}{\alpha}(U \gray V))$, and let
	\[
		i \eqdef \min \set{ j \in \set{0, \ldots, n} \mid x \in \bound{j}{\alpha}U },
	\]
	which is well-defined by equation (\ref{eq:boundaries_of_gray}).
	Then $(x, y) \in \bound{i}{\alpha}U \gray \bound{n-i}{(-)^i\alpha}V$, and
	by Lemma \ref{lem:faces_intersection} $(x, y)$ must be a maximal element in
	\[
		\bound{k}{\alpha}(\bound{i}{\alpha}U \gray \bound{n-i}{(-)^i\alpha}V) = \bigcup_{j=0}^k \bound{j}{\alpha}(\bound{i}{\alpha}U) \gray \bound{k-j}{(-)^j\alpha}(\bound{n-i}{(-)^i\alpha}V).
	\]
	Since by construction $x \notin \bound{j}{\alpha}U = \bound{j}{\alpha}(\bound{i}{\alpha}U)$ for all $j < i$, in fact
	\[
		(x, y) \in \bigcup_{j=i}^k \bound{j}{\alpha}U \gray \bound{k-j}{(-)^j\alpha}V \subseteq \bound{k}{\alpha}(U \gray V),
	\]
	where we used the fact that, by globularity of $V$,
	\[
		\bound{k-j}{(-)^j\alpha}(\bound{n-i}{(-)^i\alpha}V) = \bound{k-j}{(-)^j\alpha}V
	\]
	for all $j \geq i$, since $k < n$ and $j \geq i$ imply that $k - j < n-i$.
	This proves that
	\[
		\maxel{(\bound{k}{\alpha}(\bound{n}{\alpha}(U \gray V)))} \subseteq \bound{k}{\alpha}(U \gray V),
	\]
	so using Lemma \ref{lem:closure_of_maximal} we conclude that $\bound{k}{\alpha}(\bound{n}{\alpha}(U \gray V)) = \bound{k}{\alpha}(U \gray V)$.

	Next, we have, for all $j \in \set{0, \ldots, n-1}$, that
	\[
		\bound{j}{\alpha}U = \bound{j}{\alpha}(\bound{j+1}{-\alpha}U) \subseteq \bound{j+1}{-\alpha}U.
	\]
	It follows that
	\begin{align*}
		\bound{n-1}{\alpha}(U \gray V) & = \bigcup_{j=1}^n\bound{j-1}{\alpha}U \gray \bound{n-j}{(-)^{j-1}\alpha}V \subseteq \\
					       & \subseteq \bigcup_{j=0}^n \bound{j}{-\alpha}U \gray \bound{n-j}{-(-)^j\alpha}V = \bound{n}{-\alpha}(U \gray V),
	\end{align*}
	and by Lemma \ref{lem:boundary_included_in_subset} we have
	\[
		\bound{n-1}{\alpha}(U \gray V) \subseteq \bound{n-1}{\alpha}(\bound{n}{-\alpha}(U \gray V)).
	\]
	Conversely, let $(x, y)$ be a maximal element in $\bound{n-1}{\alpha}(\bound{n}{-\alpha}(U \gray V))$, and let
	\[
		i \eqdef \max \set{ j \in \set{0, \ldots, n} \mid y \in \bound{n-j}{-(-)^{j}\alpha}V }.
	\]
	Proceeding dually to the first part of the proof, we deduce that $(x, y)$ is maximal in $\bound{n-1}{\alpha}(\bound{i}{-\alpha}U \gray \bound{n-i}{-(-)^i\alpha}V)$, and, from that,
	\[
		(x, y) \in \bound{n-1}{\alpha}(U \gray V)
	\]
	hence $\bound{n-1}{\alpha}(\bound{n}{-\alpha}(U \gray V)) = \bound{n-1}{\alpha}(U \gray V)$.
	Finally,
	\[
		\bound{k}{\alpha}(\bound{n}{-\alpha}(U \gray V)) = \bound{k}{\alpha}(\bound{n-1}{\alpha}(\bound{n}{-\alpha}(U \gray V))) = \bound{k}{\alpha}(\bound{n-1}{\alpha}(U \gray V)) = \bound{k}{\alpha}(U \gray V)
	\]
	and we conclude.
\end{proof}

\begin{lem} \label{lem:boundaries_of_globular_gray_product}
	Let $U$, $V$ be globular oriented graded posets, $n \in \mathbb{N}$, and $\alpha \in \set{+, -}$.
	Then
	\[
		\bound{n}{\alpha}(U \gray V) = \quad \quad \bigcup_{\mathclap{k = \max \set{0, n - \dim{V}}}}^{\mathclap{\min \set{n, \dim{U}}}} \;\; \bound{k}{\alpha}U \gray \bound{n-k}{(-)^k \alpha}V = \bound{n}{\alpha}(\bound{n}{\alpha}U \gray \bound{n}{\alpha}V).
	\]
\end{lem}
\begin{proof}
	Let $m \eqdef \dim{V}$ and suppose that $0 \leq k < n - m$, so $\dim{V} < n-k$.
	Then $\bound{n-k}{(-)^k \alpha}V = V$, while by globularity of $U$
	\[
		\bound{k}{\alpha}U = \bound{k}{\alpha}(\bound{n-m}{\alpha}U) \subseteq \bound{n-m}{\alpha}U,
	\]
	so 
	\[
		\bound{k}{\alpha}U \gray \bound{n-k}{(-)^k \alpha}V \subseteq \bound{n-m}{\alpha}U \gray \bound{m}{(-)^{n-m}\alpha}V, 
	\]
	and we can omit the first $n - \dim{V}$ terms in the union (\ref{eq:boundaries_of_gray}).
	Similarly, let $p \eqdef \dim{U}$ and suppose $p < k \leq n$.
	Then $\bound{k}{\alpha}U = U$, while by globularity of $V$
	\[
		\bound{n-k}{(-)^k\alpha}V = \bound{n-k}{(-)^k\alpha}(\bound{n-p}{(-)^p\alpha}V) \subseteq \bound{n-p}{(-)^p\alpha}V,
	\]
	so
	\[
		\bound{k}{\alpha}U \gray \bound{n-k}{(-)^k \alpha}V \subseteq \bound{p}{\alpha}U \gray \bound{n-p}{(-)^p\alpha}V, 
	\]
	and we can omit the last $n - \dim{U}$ terms in the union (\ref{eq:boundaries_of_gray}).
	This proves one equation.
	For the other, it suffices to observe that, for all $k \in \set{0, \ldots, n}$,
	\[
		\bound{k}{\alpha}U = \bound{k}{\alpha}(\bound{n}{\alpha}U), \quad \bound{n-k}{(-)^k\alpha}V = \bound{n-k}{(-)^k\alpha}(\bound{n}{\alpha}V). \qedhere
	\]
\end{proof}

\begin{lem} \label{lem:gray_preserves_round}
	Let $U$, $V$ be round oriented graded posets.
	Then $U \gray V$ is round.
\end{lem}
\begin{proof}
	By Lemma \ref{lem:gray_preserves_globular}, $U \gray V$ is globular.
	Let $n < \dim{U \gray V}$.
	Then
	\begin{align*}
		\bound{n}{-}(U \gray V) \cap \bound{n}{+}(U \gray V) & = 
		\bigcup_{i=0}^n (\bound{i}{-}U \gray \bound{n-i}{-(-)^i}V) \cap
		\bigcup_{j=0}^n (\bound{j}{+}U \gray \bound{n-j}{(-)^j}V) = \\
								     & =
		\bigcup_{i, j=0}^n (\bound{i}{-}U \cap \bound{j}{+}U) \gray
		(\bound{n-i}{-(-)^i}V \cap \bound{n-j}{(-)^j}V)
	\end{align*}
	by elementary properties of unions, intersections, and products.
	We can split this union into unions over pairs $(i, j)$ with $i < j$, $i > j$, and $i = j$, respectively.
	When $i < j$, by globularity of $U$ and $V$, we have
	\[
		\bound{i}{-}U \cap \bound{j}{+}U = \bound{i}{-}U, \quad \quad \bound{n-i}{-(-)^i}V \cap \bound{n-j}{(-)^j}V = \bound{n-j}{(-)^j}V,
	\]
	and also $\bound{n-j}{(-)^j}V \subseteq \bound{n-i-1}{-(-)^{i}}V$ whenever $j > i+1$, so
	\begin{equation} \label{eq:round_gray_first_half}
		\bigcup_{i < j} (\bound{i}{-}U \cap \bound{j}{+}U) \gray
		(\bound{n-i}{-(-)^i}V \cap \bound{n-j}{(-)^j}V) = \bigcup_{i=0}^n \bound{i}{-}U \gray \bound{n-i-1}{-(-)^i}V,
	\end{equation}
	which is equal to $\bound{n-1}{-}(U \gray V)$.
	By a dual argument,
	\begin{equation} \label{eq:round_gray_second_half}
		\bigcup_{i > j} (\bound{i}{-}U \cap \bound{j}{+}U) \gray
		(\bound{n-i}{-(-)^i}V \cap \bound{n-j}{(-)^j}V) = \bigcup_{i=0}^n \bound{i}{+}U \gray \bound{n-i-1}{(-)^i}V,
	\end{equation}
	which is equal to $\bound{n-1}{+}(U \gray V)$.
	Finally, by roundness of $U$ and $V$, we have
	\begin{align*}
		\bound{i}{-}U \cap \bound{i}{+}U & = \begin{cases}
			\bound{i-1}{}U & \text{if $i < \dim{U}$,} \\
			U & \text{if $i \geq \dim{U}$},
		\end{cases} \\
		\bound{n-i}{-(-)^i}V \cap \bound{n-i}{(-)^i}V & = \begin{cases}
			\bound{n-i-1}{}V & \text{if $n - i < \dim{V}$,} \\
			V & \text{if $n - i \geq \dim{V}$}.
		\end{cases}
	\end{align*}
	Thus $\displaystyle\bigcup_{i =0}^n (\bound{i}{-}U \cap \bound{i}{+}U) \gray (\bound{n-i}{-(-)^i}V \cap \bound{n-i}{(-)^i}V)$ is equal to 
	\[
		\bigcup_{i=0}^{n-\dim{V}} \bound{i-1}{}U \gray V \cup
		\bigcup_{i=n-\dim{V}+1}^{\dim{U}-1} \bound{i-1}{}U \gray \bound{n-i-1}{}V \cup
		\bigcup_{i=\dim{U}}^n U \gray \bound{n-i-1}{}V,
	\]
	whose every term is included either in (\ref{eq:round_gray_first_half}) or in (\ref{eq:round_gray_second_half}).
	We conclude that
	\[
		\bound{n}{-}(U \gray V) \cap \bound{n}{+}(U \gray V) = \bound{n-1}{-}(U \gray V) \cup \bound{n-1}{+}(U \gray V) = \bound{n-1}{}(U \gray V),
	\]
	that is, $U \gray V$ is round.
\end{proof}

\begin{lem} \label{lem:splitting_boundary_of_globular_gray}
	Let $U$, $V$ be globular oriented graded posets, $n, j \in \mathbb{N}$, and $\alpha \in \set{+, -}$.
	Then
	\begin{enumerate}
		\item $\bound{n}{\alpha}(U \gray V) = \bound{n}{\alpha}(\bound{j}{\alpha}U \gray V) \cup \bound{n}{\alpha}(U \gray \bound{n-j-1}{(-)^{j+1}\alpha}V)$,
		\item $\bound{n}{\alpha}(\bound{j}{\alpha}U \gray V) \cap \bound{n}{\alpha}(U \gray \bound{n-j-1}{(-)^{j+1}\alpha}V) = \bound{j}{\alpha}U \gray \bound{n-j-1}{(-)^{j+1}\alpha}V$.
	\end{enumerate}
\end{lem}
\begin{proof}
	By Lemma \ref{lem:boundaries_of_globular_gray_product}, we have
	\[
		\bound{n}{\alpha}(\bound{j}{\alpha}U \gray V) = \bigcup_{k=0}^j \bound{k}{\alpha}(\bound{j}{\alpha}U) \gray \bound{n-k}{(-)^k\alpha}V = \bigcup_{k=0}^j \bound{k}{\alpha}U \gray \bound{n-k}{(-)^k\alpha}V,
	\]
	while
	\[
		\bound{n}{\alpha}(U \gray \bound{n-j-1}{(-)^{j+1}\alpha}V) = \bigcup_{k=j+1}^n \bound{k}{\alpha}U \gray \bound{n-k}{(-)^k\alpha}(\bound{n-j-1}{(-)^{j+1}\alpha}V) = \bigcup_{k=j+1}^n \bound{k}{\alpha}U \gray \bound{n-k}{(-)^k\alpha}V,
	\]
	and the union of the two is $\bound{n}{\alpha}(U \gray V)$.
	Their intersection is contained in
	\[
		(\bound{j}{\alpha}U \gray V) \cap (U \gray \bound{n-j-1}{(-)^{j+1}\alpha}V) = \bound{j}{\alpha}U \gray \bound{n-j-1}{(-)^{j+1}\alpha}V.
	\]
	But $\bound{j}{\alpha}U \subseteq \bound{j+1}{\alpha}U$ and $\bound{n-j-1}{(-)^{j+1}\alpha}V \subseteq \bound{n-j}{(-)^{j}\alpha}V$ by globularity, so
	\begin{align*}
		\bound{j}{\alpha}U \gray \bound{n-j-1}{(-)^{j+1}\alpha}V & \subseteq \bound{j}{\alpha}U \gray \bound{n-j}{(-)^j\alpha}V \subseteq \bound{n}{\alpha}(\bound{j}{\alpha}U \gray V), \\
		\bound{j}{\alpha}U \gray \bound{n-j-1}{(-)^{j+1}\alpha}V & \subseteq \bound{j+1}{\alpha}U \gray \bound{n-j-1}{(-)^{j+1}\alpha}V \subseteq \bound{n}{\alpha}(U \gray \bound{n-j-1}{(-)^{j+1}\alpha}V),
	\end{align*}
	and we conclude.
\end{proof}

\begin{prop} \label{prop:gray_product_of_molecules}
	Let $U, V$ be molecules.
	Then
	\begin{enumerate}
		\item $U \gray V$ is a molecule,
		\item for all $n, j \in \mathbb{N}$,
			\begin{align*}
				\bound{n}{-}(U \gray V) & = 
				\bound{n}{-}(\bound{j}{-}U \gray V) \gencp{n-1} \bound{n}{-}(U \gray \bound{n-j-1}{(-)^{j}}V), \\
				\bound{n}{+}(U \gray V) & = 
				\bound{n}{+}(U \gray \bound{n-j-1}{(-)^{j+1}}V) \gencp{n-1} \bound{n}{+}(\bound{j}{+}U \gray V)
			\end{align*}
		\item if $U$ splits into $W \cup W'$ along the $k$\nbd boundary, then
			\[
				U \gray V = (W \gray V) \gencp{k+\dim{V}} (W' \gray V),
			\]
			and if $V$ splits into $W \cup W'$ along the $k$\nbd boundary, then
			\[
				U \gray V = \begin{cases}
					(U \gray W) \gencp{k+\dim{U}} (U \gray W') & \text{if $\dim{U}$ is even}, \\
					(U \gray W') \gencp{k+\dim{U}} (U \gray W) & \text{if $\dim{U}$ is odd}.
				\end{cases}
			\]
	\end{enumerate}
\end{prop}
\begin{proof}
	We proceed by double induction on submolecules of $U$ and $V$.
	We have that $\set{x} \gray V$ is isomorphic to $V$ for all $x \in \grade{0}{U}$ and $U \gray \set{y}$ is isomorphic to $U$ for all $y \in \grade{0}{V}$, and in this case all statements follow from $U$ and $V$ being molecules.
	Now, suppose that the statement holds for all proper submolecules of $U$ or $V$.
	We will prove by recursion on $n$ that $\bound{n}{\alpha}(U \gray V)$ is a molecule and satisfies the statement for all $\alpha \in \set{+, -}$ and $n \in \set{0, \ldots, \dim{U \gray V} - 1}$.
	We have
	\[
		\bound{0}{\alpha}(U \gray V) = \bound{0}{\alpha}U \gray \bound{0}{\alpha}V
	\]
	which is evidently a 0\nbd dimensional molecule.
	Suppose $n > 0$, and let $j$ be such that $n - \dim{V} \leq j < \dim{U}$; this is always possible since $n < \dim{U} + \dim{V}$.
	By Lemma \ref{lem:splitting_boundary_of_globular_gray}, we have
	\begin{align*}
		\bound{n}{\alpha}(U \gray V) & = \bound{n}{\alpha}(\bound{j}{\alpha}U \gray V) \cup \bound{n}{\alpha}(U \gray \bound{n-j-1}{(-)^{j+1}\alpha}V), \\
		\bound{j}{\alpha}U \gray \bound{n-j-1}{(-)^{j+1}\alpha}V & = \bound{n}{\alpha}(\bound{j}{\alpha}U \gray V) \cap \bound{n}{\alpha}(U \gray \bound{n-j-1}{(-)^{j+1}\alpha}V)
	\end{align*}
	and both $\bound{j}{\alpha}U \submol U$ and $\bound{n-j-1}{(-)^{j+1}\alpha}V \submol V$ are proper submolecules.
	By the inductive hypothesis on submolecules, $\bound{n}{\alpha}(\bound{j}{\alpha}U \gray V)$ and $\bound{n}{\alpha}(U \gray \bound{n-j-1}{(-)^{j+1}\alpha}V)$ are both molecules.
	Moreover, for all $\beta \in \set{+, -}$,
	\begin{enumerate}
		\item $\bound{n-1}{\beta}(U \gray V) = \bound{n-1}{\beta}(\bound{n}{\alpha}(U \gray V))$ by Lemma \ref{lem:gray_preserves_globular},
		\item $\bound{n-1}{\beta}(U \gray V)$ is a molecule by the inductive hypothesis on $n$,
		\item $\bound{n-1}{\alpha}(\bound{j}{\alpha}U \gray V) \submol \bound{n-1}{\alpha}(U \gray V)$ and $\bound{n-1}{-\alpha}(U \gray \bound{n-j-1}{(-)^{j+1}\alpha}V) \submol \bound{n-1}{-\alpha}(U \gray V)$ by the inductive hypothesis on $n$,
		\item $\bound{j}{\alpha}U \gray \bound{n-j-1}{(-)^{j+1}\alpha}V = \bound{n-1}{-\alpha}(\bound{j}{\alpha}U \gray \bound{n-j-1}{(-)^{j+1}\alpha}V) \submol \bound{n-1}{-\alpha}(\bound{j}{\alpha}U \gray V)$, and also
			$\bound{j}{\alpha}U \gray \bound{n-j-1}{(-)^{j+1}\alpha}V = \bound{n-1}{\alpha}(\bound{j}{\alpha}U \gray \bound{n-j-1}{(-)^{j+1}\alpha}V) \submol \bound{n-1}{\alpha}(U \gray \bound{n-j-1}{(-)^{j+1}\alpha}V)$, by the inductive hypothesis on submolecules.
	\end{enumerate}
	Thus the assumptions of Lemma \ref{lem:generalised_pasting} are satisfied for the pushout diagram
\[\begin{tikzcd}
	\bound{j}{\alpha}U \gray \bound{n-j-1}{(-)^{j+1}\alpha}V && \bound{n}{\alpha}(U \gray \bound{n-j-1}{(-)^{j+1}\alpha}V) \\
	\bound{n}{\alpha}(\bound{j}{\alpha}U \gray V) && \bound{n}{\alpha}(U \gray V),
	\arrow[hook, from=1-1, to=1-3]
	\arrow[hook', from=1-1, to=2-1]
	\arrow[hook, from=2-1, to=2-3]
	\arrow[hook', from=1-3, to=2-3]
	\arrow["\lrcorner"{anchor=center, pos=0.125, rotate=180}, draw=none, from=2-3, to=1-1]
\end{tikzcd}\]
which proves that $\bound{n}{\alpha}(U \gray V)$ is a molecule and satisfies the statement for all $n - \dim{V} \leq j < \dim{U}$.
	The remaining cases with $j < n - \dim{V}$ and $j \geq \dim{U}$ are either trivial, or follow by the inductive hypothesis on submolecules.

	If $U$ and $V$ are both atoms, then this, together with the fact that $U \gray V$ is round by Corollary \ref{cor:atoms_are_round} and Lemma \ref{lem:gray_preserves_round}, is sufficient to prove that $U \gray V$ is an atom, isomorphic to $\bound{}{-}(U \gray V) \celto \bound{}{+}(U \gray V)$.
	Otherwise, suppose that $U$ splits into proper submolecules $W \cup W'$ along the $k$\nbd boundary.
	Then
	\begin{align*}
		U \gray V & = (W \gray V) \cup (W' \gray V), \\
		(W \gray V) \cap (W' \gray V) & = \bound{k}{+}W \gray V = \bound{k}{-}W' \gray V.
	\end{align*}
	Letting $m \eqdef \dim{V}$, we have that
	\begin{enumerate}
		\item $\bound{k}{+}W \gray V = \bound{k+m}{+}(\bound{k}{+}W \gray V) \submol \bound{k+m}{+}(W \gray V)$ and, symmetrically,
			$\bound{k}{-}W' \gray V = \bound{k+m}{-}(\bound{k}{-}W' \gray V) \submol \bound{k+m}{-}(W' \gray V)$ by the inductive hypothesis on submolecules,
		\item $\bound{k+m}{-}(U \gray V)$ and $\bound{k+m}{+}(U \gray V)$ are molecules by the first part of the proof,
	\end{enumerate}
	so to apply Lemma \ref{lem:generalised_pasting} it suffices to show that
	\[
		\bound{k+m}{-}(W \gray V) \submol \bound{k+m}{-}(U \gray V) 
		\quad \text{and} \quad
		\bound{k+m}{+}(W' \gray V) \submol \bound{k+m}{+}(U \gray V).
	\]
	We will deduce $\bound{k+m}{-}(W \gray V) \submol \bound{k+m}{-}(U \gray V)$ from the following statement: for all $j \in \set{1, \ldots, m}$,
	\begin{align*}
		\bound{k+m}{-}(U & \gray \bound{m-j}{(-)^{k+j-1}}V) = \\
				 & \bound{k+m}{-}(W \gray \bound{m-j}{(-)^{k+j-1}}V) \gencp{k+m-j} \bound{k+m}{-}(W' \gray \bound{m-j}{(-)^{k+j-1}}V).
	\end{align*}
	We will prove this by backward recursion on $j$.
	For $j = m$,
	\begin{align*}
		\bound{k+m}{-}(U & \gray \bound{0}{(-)^{k+m-1}}V) = \bound{k+m}{-}U \gray \bound{0}{(-)^{k+m-1}}V = \\
								& = (\bound{k+m}{-}W \cp{k} \bound{k+m}{-}W') \gray \bound{0}{(-)^{k+m-1}}V = \\
								& = \bound{k+m}{-}(W \gray \bound{0}{(-)^{k+m-1}}V) \cp{k} \bound{k+m}{-}(W' \gray \bound{0}{(-)^{k+m-1}}V),
	\end{align*}
	which has the required form by Lemma \ref{lem:pasting_is_generalised_pasting_in_all_higher_dim}.
	Let $0 < j < m$.
	Then
	\[
		\bound{k+m}{-}(U \gray \bound{m-j}{(-)^{k+j-1}}V) = (\bound{k+j}{-}U \gray \bound{m-j}{(-)^{k+j-1}}V) \gencp{k+m-1} \bound{k+m}{-}(U \gray \bound{m-j-1}{(-)^{k+j}}V)
	\]
	by the first part of the proof.
Now, $\bound{k+j}{-}U = \bound{k+j}{-}W \cp{k} \bound{k+j}{-}W'$, and since $\bound{m-j}{(-)^{k+j-1}}V$ is a proper submolecule of $V$, we can use the inductive hypothesis on submolecules to rewrite $\bound{k+j}{-}U \gray \bound{m-j}{(-)^{k+j-1}}V$ as
	\[
		(\bound{k+j}{-}W \gray \bound{m-j}{(-)^{k+j-1}}V) \gencp{k+m-j} (\bound{k+j}{-}W' \gray \bound{m-j}{(-)^{k+j-1}}V).
	\]
	By the inductive hypothesis on $j$, $\bound{k+m}{-}(U \gray \bound{m-j-1}{(-)^{k+j}}V)$ is equal to
	\[
		\bound{k+m}{-}(W \gray \bound{m-j-1}{(-)^{k+j}}V) \gencp{k+m-j-1} \bound{k+m}{-}(W' \gray \bound{m-j-1}{(-)^{k+j}}V).
	\]
	Now, both
	\begin{align*}
		(\bound{k+j}{-}W \gray \bound{m-j}{(-)^{k+j-1}}V) & \gencp{k+m-1} \bound{k+m}{-}(W \gray \bound{m-j-1}{(-)^{k+j}}V), \\
		(\bound{k+j}{-}W' \gray \bound{m-j}{(-)^{k+j-1}}V) & \gencp{k+m-1} \bound{k+m}{-}(W' \gray \bound{m-j-1}{(-)^{k+j}}V)
	\end{align*}
	are defined and equal to $\bound{k+m}{-}(W \gray \bound{m-j}{(-)^{k+j-1}}V)$ and to $\bound{k+m}{-}(W' \gray \bound{m-j}{(-)^{k+j-1}}V)$, respectively.
	Expanding the generalised pastings as in Lemma \ref{lem:generalised_pasting} and using Proposition \ref{prop:interchange_of_pasting}, we deduce that $\bound{k+m}{-}(U \gray \bound{m-j}{(-)^{k+j-1}}V)$ is equal to
	\[
		\bound{k+m}{-}(W \gray \bound{m-j}{(-)^{k+j-1}}V) \gencp{k+m-j} \bound{k+m}{-}(W' \gray \bound{m-j}{(-)^{k+j-1}}V),
	\]
	which completes the inductive step on $j$.
	Finally, we have
	\[
		\bound{k+m}{-}(U \gray V) = (\bound{k}{-}U \gray V) \gencp{k+m-1} \bound{k+m}{-}(U \gray \bound{}{(-)^{k}}V).
	\]
	But $\bound{k}{-}U = \bound{k}{-}W$, so this is equal to
	\begin{align*}
		(\bound{k}{-}W \gray V) & \gencp{k+m-1} (\bound{k+m}{-}(W \gray \bound{}{(-)^{k}}V) \gencp{k+m-1} \bound{k+m}{-}(W' \gray \bound{}{(-)^{k}}V)) = \\
					& = \bound{k+m}{-}(W \gray V) \gencp{k+m-1} \bound{k+m}{-}(W' \gray \bound{}{(-)^{k}}V),
	\end{align*}
	expanding the generalised pastings and using Proposition \ref{prop:associativity_of_pasting}.

	This proves that $\bound{k+m}{-}(W \gray V) \submol \bound{k+m}{-}(U \gray V)$, and a dual proof determines that $\bound{k+m}{+}(W' \gray V) \submol \bound{k+m}{+}(U \gray V)$.
	By Lemma \ref{lem:generalised_pasting}, we conclude that
	\[
		U \gray V = (W \gray V) \gencp{k+m} (W' \gray V).
	\]
	The proof when $V$ splits into proper submolecules is entirely analogous, with a little extra care about sign flips when $U$ is odd-dimensional.
\end{proof}

\begin{cor} \label{cor:gray_products_of_rdcpx}
	Let $P$, $Q$ be regular directed complexes.
	Then $P \gray Q$ is a regular directed complex.
\end{cor}
\begin{proof}
	Let $(x, y) \in P \gray Q$.
	Then $\clset{(x, y)}$ is isomorphic to $\clset{x} \gray \clset{y}$, which is an atom by Proposition \ref{prop:gray_product_of_molecules}.
\end{proof}

\begin{cor} \label{cor:monoidal_structure_on_rdcpx}
	The monoidal structure $(\ogpos, \gray, 1)$ restricts to a monoidal structure on $\rdcpx$ and on $\rdcpxiso$.
\end{cor}

\begin{prop} \label{prop:gray_product_of_maps}
	There is a unique monoidal structure $(\rdcpxmap, \gray, 1)$ on $\rdcpxmap$ such that both
	\begin{enumerate}
		\item $(\rdcpx, \gray, 1) \incl (\rdcpxmap, \gray, 1)$ and
		\item $\fun{U}\colon (\rdcpxmap, \gray, 1) \to (\posclos, \times, 1)$
	\end{enumerate}
	are strict monoidal functors.
\end{prop}
\begin{proof}
	Since the monoidal structure extends the one on $\rdcpx$, it is the Gray product on objects, and the requirement that $\fun{U}$ be strict monoidal determines uniquely what it does on maps, so it suffices to show that if $p\colon P \to P'$ and $q\colon Q \to Q'$ are maps of regular directed complexes, then $p \times q$ lifts to a map $p \gray q\colon P \gray Q \to P' \gray Q'$ of regular directed complexes.

	Let $(x, y) \in P \gray Q$, $n \in \mathbb{N}$, and $\alpha \in \set{+, -}$.
	Then
	\begin{align*}
		\bound{n}{\alpha}(p(x), q(y)) & = \bigcup_{k \in \mathbb{N}} \bound{k}{\alpha}p(x) \gray \bound{n-k}{(-)^k\alpha}q(y) = \bigcup_{k \in \mathbb{N}} p(\bound{k}{\alpha}x) \gray q(\bound{n-k}{(-)^k\alpha}y) = \\
					      & = (p \times q) \left( \bigcup_{k \in \mathbb{N}} \bound{k}{\alpha}x \gray \bound{n-k}{(-)^k\alpha}y \right) = (p \times q)\bound{n}{\alpha}(x, y).
	\end{align*}
	Now, let $(u, v), (u', v') \in \bound{n}{\alpha}(x, y)$ and suppose that $(w, z) \leq (p(u), q(v))$ and $(w, z) \leq (p(u'), q(v'))$.
	We have 
	\[
		(u, v) \in \bound{i}{\alpha}x \gray \bound{n-i}{(-)^i\alpha}y, \quad (u', v') \in \bound{j}{\alpha}x \gray \bound{n-j}{(-)^j\alpha}y
	\]
	for some $i, j \in \set{0, \ldots, n}$.
	Suppose without loss of generality that $i \leq j$.
	Then $u, u' \in \bound{j}{\alpha}x$, $w \leq p(u)$ and $w \leq p(u')$.
	Because $\restr{p}{\bound{j}{\alpha}x}$ is final onto its image, there is a zig-zag
	\[
		u \leq u_1 \geq \ldots \leq u_m \geq u'
	\]
	in $\bound{j}{\alpha}x$ such that $w \leq p(u_i)$ for all $i \in \set{1, \ldots, m}$.
	This induces a zig-zag
	\[
		(u, v') \leq (u_1, v') \geq \ldots \leq (u_m, v') \geq (u', v')
	\]
	in $\bound{j}{\alpha}x \gray \bound{n-j}{(-)^j\alpha}y \subseteq \bound{n}{\alpha}(x, y)$ with $(w, z) \leq (p(u_i), q(v'))$ for all $i \in \set{1, \ldots, m}$.
	Now, since $i \leq j$, we also have $v, v' \in \bound{n-i}{(-)^i\alpha}y$, $z \leq q(v)$, and $z \leq q(v')$.
	Because $\restr{q}{\bound{n-i}{(-)^i\alpha}y}$ is final onto its image, there is a zig-zag
	\[
		v \leq v_1 \geq \ldots \leq v_\ell \geq v'
	\]
	in $\bound{n-i}{(-)^i\alpha}y$ such that $z \leq q(v_i)$ for all $i \in \set{1, \ldots, \ell}$.
	Then
	\[
		(u, v) \leq (u, v_1) \geq \ldots \leq (u, v_\ell) \geq (u, v')
	\]
	is a zig-zag in $\bound{i}{\alpha}x \gray \bound{n-i}{(-)^i\alpha}y \subseteq \bound{n}{\alpha}(x, y)$ with $(w, z) \leq (p(u), q(v_i))$ for all $i \in \set{1, \ldots, \ell}$.
	Concatenating the two zig-zags, we deduce that $\restr{(p \times q)}{\bound{n}{\alpha}(x, y)}$ is final onto its image, and $p \times q$ lifts to a map.
\end{proof}

\begin{rmk}
	By Proposition \ref{prop:terminal_object_of_rdcpxmap}, $(\rdcpxmap, \gray, 1)$ is \emph{semicartesian} monoidal, that is, the monoidal unit is the terminal object.
	It follows that, for all regular directed complexes $P, Q$, there are natural projection maps
	\[
		P \gray Q \to P, \quad \quad P \gray Q \to Q
	\]
	obtained by composing $\idd{P} \gray \varepsilon\colon P \gray Q \to P \gray 1$ with a right unitor and $\varepsilon \gray \idd{Q}\colon P \gray Q \to 1 \gray Q$ with a left unitor, where $\varepsilon$ is the unique map to the terminal object.
\end{rmk}

\begin{prop} \label{prop:gray_product_of_comaps}
	There is a unique monoidal structure $(\rdcpxcomap, \gray, 1)$ on $\rdcpxcomap$ such that both
	\begin{enumerate}
		\item $(\rdcpxiso, \gray, 1) \incl (\rdcpxcomap, \gray, 1)$ and
		\item $\fun{U}\colon (\rdcpxcomap, \gray, 1) \to (\poscat, \times, 1)$
	\end{enumerate}
	are strict monoidal functors.
\end{prop}
\begin{proof}
	As in the case of maps, the requirements fix the monoidal structure uniquely, so it suffices to show that, if $c\colon P \to P'$ and $d\colon Q \to Q'$ are comaps, then $c \times d$ lifts to a comap $c \gray d\colon P \gray Q \to P' \gray Q'$.

	Let $(x, y) \in P' \gray Q'$.
	Then 
	\[
		\invrs{(c \times d)}\clset{(x, y)} = \invrs{c}\clset{x} \gray \invrs{d}\clset{y},
	\]
	which is a molecule by Proposition \ref{prop:gray_product_of_molecules}.
	Moreover,
	\begin{align*}
		\bound{n}{\alpha} & \invrs{(c \times d)}\clset{(x, y)} = \bigcup_{k\in \mathbb{N}} \bound{k}{\alpha}\invrs{c}\clset{x} \gray \bound{n-k}{(-)^k\alpha}\invrs{d}\clset{y} = \\
				  & = \bigcup_{k \in \mathbb{N}} \invrs{c}(\bound{k}{\alpha}x) \gray \invrs{d}(\bound{n-k}{(-)^{k}\alpha}y) = \\
		& = \invrs{(c \times d)}\left( \bigcup_{k \in \mathbb{N}} \bound{k}{\alpha}x \gray \bound{n-k}{(-)^k\alpha}y \right) = \invrs{(c \times d)}\bound{n}{\alpha}(x, y)
	\end{align*}
	using the definition of comap and basic properties of inverse images and products.
	It follows that $c \times d$ lifts to a comap.
\end{proof}

\begin{dfn}[Tensor product of augmented chain complexes] \index{chain complex!tensor product} \index{$C \otimes D$}
	Let $C$, $D$ be augmented chain complexes.
	The \emph{tensor product of $C$ and $D$} is the augmented chain complex $C \otimes D$ defined by
	\[
		\grade{n}{(C \otimes D)} \eqdef \bigoplus_{k=0}^n \grade{k}{C} \otimes \grade{n-k}{D}
	\]
	for all $n \in \mathbb{N}$, with
	\begin{align*}
		\der\colon \grade{n}{(C \otimes D)} & \to \grade{n-1}{(C \otimes D)}, \\
		x \otimes y & 
		\mapsto \der (x) \otimes y + (-)^k x \otimes \der (y)
	\end{align*}
	for each $n > 0$, $k \leq n$, $x \in \grade{k}{C}$, and $y \in \grade{n-k}{D}$, and 
	\begin{align*}
		\eau\colon \grade{0}{(C \otimes D)} & \to \mathbb{Z}, \\
		x \otimes y & \mapsto \eau (x) \eau(y)
	\end{align*}
	for all $x \in \grade{0}{C}$ and $y \in \grade{0}{D}$.
	The tensor product extends to a monoidal structure on $\chaug$, whose unit is the chain complex $\mathbb{Z}$ with
	\[
		\grade{n}{\mathbb{Z}} \eqdef \begin{cases}
			\mathbb{Z} & \text{if $n = 0$}, \\
			0 & \text{if $n > 0$},
		\end{cases}
		\quad \quad \eau \eqdef \idd{\mathbb{Z}}\colon \mathbb{Z} \to \mathbb{Z}.
	\]
\end{dfn}

\begin{prop} \label{prop:gray_to_tensor_product}
	Let $P$, $Q$ be oriented graded posets such that $\augm{P}$, $\augm{Q}$ are oriented thin.
	Then $\freeab{P} \otimes \freeab{Q}$ is naturally isomorphic to $\freeab{(P \gray Q)}$.
\end{prop}
\begin{proof}
	By Proposition \ref{prop:product_of_graded_is_graded}, we have, for all $n \in \mathbb{N}$,
	\[
		\grade{n}{(P \times Q)} = \sum_{k=0}^n \grade{k}{P} \times \grade{n-k}{Q},
	\]
	and there is a natural isomorphism of abelian groups
	\begin{align*}
		\grade{n}{\varphi}\colon \grade{n}{\freeab{(P \gray Q)}} \equiv \freeab{\left( \sum_{k=0}^n \grade{k}{P} \times \grade{n-k}{Q} \right)} &\to
		\bigoplus_{k=0}^n \freeab{\grade{k}{P}} \otimes \freeab{\grade{n-k}{Q}}
		\equiv \grade{n}{(\freeab{P} \otimes \freeab{Q})}, \\
		(x, y) \in \grade{k}{P} \times \grade{n-k}{Q} & \mapsto x \otimes y.
	\end{align*}
	Then, for all $n > 0$, $k \leq n$, $x \in \grade{k}{P}$, and $y \in \grade{n-k}{Q}$,
	\begin{align*}
		& \der \grade{n}{\varphi}(x, y) = \der (x \otimes y) = \der (x) \otimes y + (-)^k x \otimes \der(y) = \\
		& \; = \sum_{x' \in \faces{}{+}x} x' \otimes y - \sum_{x' \in \faces{}{-}x} x' \otimes y + (-)^k\left( \sum_{y' \in \faces{}{+}y} x \otimes y' - \sum_{y' \in \faces{}{-}y} x \otimes y'\right) = \\
		& \; = \quad \; \sum_{\mathclap{(x', y') \in \faces{}{+}(x, y)}} \; x' \otimes y'
		\; - \qquad \sum_{\mathclap{(x', y') \in \faces{}{-}(x, y)}} \; x' \otimes y' 
		\; = \; \grade{n-1}{\varphi} (\der(x, y)),
	\end{align*}
	while for all $x \in \grade{0}{P}$ and $y \in \grade{0}{Q}$,
	\[
		\eau \grade{0}{\varphi}(x, y) = \eau(x \otimes y) = \eau(x)\eau(y) = 1 = \eau(x, y)
	\]
	which proves that $(\grade{n}{\varphi})_{n \in \mathbb{N}}$ is an isomorphism of augmented chain complexes.
	Naturality is straightforward.
\end{proof}

\begin{cor} \label{cor:monoidal_functors_to_chaug}
	The functors 
	\[
		\freeab{-}\colon \otgpos \to \chaug, \quad \freeab{-}\colon \rdcpxmap \to \chaug, \quad \freeab{\pb{-}}\colon \opp{\rdcpxcomap} \to \chaug
	\]
	lift to strong monoidal functors
	\begin{align*}
		&\freeab{-} \colon (\otgpos, \augm{\gray}, \augm{1}) \to (\chaug, \otimes, \mathbb{Z}), \\
		&\freeab{-} \colon (\rdcpxmap, \gray, 1) \to (\chaug, \otimes, \mathbb{Z}), \\
		&\freeab{\pb{-}} \colon (\opp{\rdcpxcomap}, \gray, 1) \to (\chaug, \otimes, \mathbb{Z}).
	\end{align*}
\end{cor}
\begin{proof}
	It suffices to check that naturality of the isomorphism of Proposition \ref{prop:gray_to_tensor_product} extends to maps and comaps, and that this natural isomorphism, as well as the evident isomorphism $\freeab{1} \iso \mathbb{Z}$, is compatible with unitors and associators.
	All of these are straightforward checks.
\end{proof}

\begin{exm}[Cylinders] \label{exm:cylinders} \index{cylinder}
	The \emph{cylinder} on a space $X$ is its cartesian product with the topological interval $I$: due to the symmetry of cartesian products, there is no essential difference between $I \times X$ and $X \times I$.
	
	On the other hand, since Gray products are not symmetric, there is both a \emph{left cylinder} $\thearrow{} \gray P$ and a \emph{right cylinder} $P \gray \thearrow{}$ on a regular directed complex $P$.
	We have seen a 2\nbd dimensional instance of both in Example \ref{exm:nonsymmetry_gray}, so here we consider a 3\nbd dimensional example.

	The left cylinder $\thearrow{} \gray \disk{2}{1}$ is a 3\nbd dimensional atom whose input and output boundaries are the oriented face posets of
	\[\begin{tikzcd}[sep=small]
	& \bullet \\
	&&&& \bullet \\
	\bullet \\
	& \bullet && \bullet
	\arrow[curve={height=6pt}, from=3-1, to=4-2]
	\arrow[curve={height=6pt}, from=4-2, to=4-4]
	\arrow[from=4-4, to=2-5]
	\arrow[from=3-1, to=1-2]
	\arrow[curve={height=-12pt}, from=1-2, to=2-5]
	\arrow[""{name=0, anchor=center, inner sep=0}, curve={height=-12pt}, from=3-1, to=4-4]
	\arrow[shorten <=8pt, curve={height=6pt}, Rightarrow, from=4-4, to=1-2]
	\arrow[shorten >=4pt, Rightarrow, from=4-2, to=0]
	\end{tikzcd} \quad \quad \text{and} \quad \quad 
\begin{tikzcd}[sep=small]
	& \bullet \\
	&& \bullet && \bullet \\
	\bullet \\
	& \bullet && \bullet
	\arrow[curve={height=6pt}, from=3-1, to=4-2]
	\arrow[curve={height=6pt}, from=4-2, to=4-4]
	\arrow[from=4-4, to=2-5]
	\arrow[from=3-1, to=1-2]
	\arrow[""{name=0, anchor=center, inner sep=0}, curve={height=-12pt}, from=1-2, to=2-5]
	\arrow[curve={height=6pt}, from=1-2, to=2-3]
	\arrow[curve={height=6pt}, from=2-3, to=2-5]
	\arrow[from=4-2, to=2-3]
	\arrow[curve={height=-6pt}, Rightarrow, from=4-4, to=2-3]
	\arrow[shorten <=3pt, shorten >=6pt, curve={height=-6pt}, Rightarrow, from=4-2, to=1-2]
	\arrow[shorten >=4pt, Rightarrow, from=2-3, to=0]
\end{tikzcd}\]
respectively, while the input and output boundaries of the right cylinder $\disk{2}{1} \gray \thearrow{}$ are the oriented face posets of
\[\begin{tikzcd}[sep=small]
	&&& \bullet \\
	\bullet && \bullet \\
	&&&& \bullet \\
	& \bullet && \bullet
	\arrow[curve={height=6pt}, from=4-4, to=3-5]
	\arrow[curve={height=6pt}, from=4-2, to=4-4]
	\arrow[from=2-1, to=4-2]
	\arrow[from=1-4, to=3-5]
	\arrow[""{name=0, anchor=center, inner sep=0}, curve={height=-12pt}, from=2-1, to=1-4]
	\arrow[curve={height=6pt}, from=2-3, to=1-4]
	\arrow[curve={height=6pt}, from=2-1, to=2-3]
	\arrow[from=2-3, to=4-4]
	\arrow[curve={height=6pt}, Rightarrow, from=4-2, to=2-3]
	\arrow[shorten <=3pt, shorten >=6pt, curve={height=6pt}, Rightarrow, from=4-4, to=1-4]
	\arrow[shorten >=4pt, Rightarrow, from=2-3, to=0]
\end{tikzcd} \quad \quad \text{and} \quad \quad
\begin{tikzcd}[sep=small]
	&&& \bullet \\
	\bullet \\
	&&&& \bullet \\
	& \bullet && \bullet
	\arrow[curve={height=6pt}, from=4-4, to=3-5]
	\arrow[curve={height=6pt}, from=4-2, to=4-4]
	\arrow[from=2-1, to=4-2]
	\arrow[from=1-4, to=3-5]
	\arrow[curve={height=-12pt}, from=2-1, to=1-4]
	\arrow[""{name=0, anchor=center, inner sep=0}, curve={height=-12pt}, from=4-2, to=3-5]
	\arrow[shorten <=8pt, curve={height=-6pt}, Rightarrow, from=4-2, to=1-4]
	\arrow[shorten >=4pt, Rightarrow, from=4-4, to=0]
\end{tikzcd}\]
	respectively.
	Notice that the two cylinders happen to be related by a duality reversing the direction of 1\nbd cells and 3\nbd cells: this is a consequence of the symmetries $\disk{2}{1} \simeq \opp{\disk{2}{1}}$ and $\thearrow{} \simeq \opp{\thearrow{}}$, coupled with the general fact that $P \gray Q \simeq \opp{(\opp{Q} \gray \opp{P})}$ as one case of Proposition \ref{prop:duals_and_gray_products}.
\end{exm}


\section{Suspensions} \label{sec:suspension}

\begin{guide}
	In this section, we study the directed counterpart of the (\emph{unreduced}, or \emph{two-point}) suspension of spaces.
	We prove that the suspension of oriented graded posets preserves molecules (Proposition 
	\ref{prop:suspension_of_molecules}), hence regular directed complexes (Corollary 
	\ref{cor:suspension_of_rdcpx}), and determines endofunctors on $\ogpos$, $\rdcpxmap$, and $\rdcpxcomap$, compatibly on their overlaps.
	Finally, we recall the definition of the suspension of a strict $\omega$\nbd category, as well as an analogous construction on augmented chain complexes, and prove that these constructions are compatible with $\molecin{-}$ and $\freeab{-}$, respectively.
\end{guide}

\begin{dfn}[Suspension of an oriented graded poset] \index{oriented graded poset!suspension} \index{regular directed complex!suspension} \index{$\sus{P}$} \index{suspension!of an oriented graded poset}
	Let $P$ be an oriented graded poset.
	The \emph{suspension of $P$} is the oriented graded poset $\sus{P}$ whose
	\begin{itemize}
		\item underlying set is 
			\[
				\set{\sus{x} \mid x \in P} + \set{\bot^+, \bot^-},
			\]
		\item partial order and orientation are defined by
			\[
			\cofaces{}{\alpha}x \eqdef 
			\begin{cases}
				\set{\sus{y} \mid y \in \cofaces{}{\alpha}x'} & 
				\text{if $x = \sus{x'}$, $x' \in P$}, \\
				\set{\sus{y} \mid y \in \grade{0}{P}} &
				\text{if $x = \bot^\alpha$}, \\
				\varnothing &
				\text{if $x = \bot^{-\alpha}$},
			\end{cases}
		\]
			for all $x \in \sus{P}$ and $\alpha \in \set{+, -}$.
	\end{itemize}
\end{dfn}

\begin{comm}
	The $n$\nbd categorical suspension is commonly denoted by $\fun{\Sigma}$, for example in \cite{steiner1993algebra} or \cite{ara2020joint}.
	In topology, however, this symbol usually denotes the \emph{reduced} suspension of a pointed space, whereas the operation we define here is related to the plain unreduced suspension, usually denoted by $\fun{S}$.
\end{comm}

\begin{lem} \label{lem:suspension_basic_properties}
Let $P$ be an oriented graded poset.
Then
\begin{enumerate}
	\item $\sus{P}$ is well-defined as an oriented graded poset,
	\item for all $x \in \sus{P}$, $\dim{x}$ is equal to
	\[
		\begin{cases}
		\dim{x'} + 1 & \text{if $x = \sus{x'}$, $x' \in P$,} \\
		0 & \text{if $x = \bot^\alpha$, $\alpha \in \set{+, -}$},
		\end{cases}
	\]
	\item for all $x \in \sus{P}$ and $\alpha \in \set{+, -}$, $\faces{}{\alpha}x$ is equal to
	\[
	\begin{cases}
		\set{\sus{y} \mid y \in \faces{}{\alpha}x'} & \text{if $x = \sus{x'}$, $\dim{x'} > 0$,} \\
		\set{\bot^\alpha} & \text{if $\dim{x} = 1$}, \\
		\varnothing & \text{if $\dim{x} = 0$}.
	\end{cases}
	\]
\end{enumerate}
\end{lem}
\begin{proof}
	By construction, $\bot^+$ and $\bot^-$ are not cofaces of any element, so they are minimal, and their dimension is well-defined and equal to $0$.
	Suppose $x = \sus{x'}$ for some $x' \in P$; we proceed by induction on $\dim{x'}$.
	If $\dim{x'} = 0$, then $\faces{}{\alpha}\sus{x'} = \set{\bot^\alpha}$ for all $\alpha \in \set{+, -}$.
	Thus all faces of $x$ have dimension 0, and $\dim{x} = 1$.
	Suppose $\dim{x'} > 0$.
	Then 
	\[
		\faces{}{\alpha}\sus{x'} = \set{\sus{y'} \mid y' \in {\faces{}{\alpha}x'}},
	\]
	and by the inductive hypothesis all its elements have dimension $\dim{x'}$.
	It follows that the dimension of $\sus{x'}$ is well-defined and equal to $\dim{x'} + 1$.
\end{proof}

\begin{lem} \label{lem:suspension_is_a_functor}
	Let $f\colon P \to Q$ be a morphism of oriented graded posets.
	Then
	\begin{align*}
		\sus{f}\colon \sus{P} & \to \sus{Q}, \\
		x & \mapsto \begin{cases}
			\sus{f(x')} & \text{if $x = \sus{x'}$, $x' \in P$, } \\
			\bot^\alpha & \text{if $x = \bot^\alpha$, $\alpha \in \set{+, -}$},
		\end{cases}
	\end{align*}
	is a morphism of oriented graded posets.
	This assignment determines a faithful endofunctor $\sus{}$ on $\ogpos$.
\end{lem}
\begin{proof}
	Let $x \in \sus{P}$ and $\alpha \in \set{+, -}$.
	The fact that $\sus{f}$ induces a bijection between $\faces{}{\alpha}x$ and $\faces{}{\alpha}f(x)$ follows by a straightforward case distinction from Lemma \ref{lem:suspension_basic_properties}.
	Functoriality and faithfulness are equally straightforward.
\end{proof}

\begin{prop} \label{prop:suspension_of_oriented_thin}
	Let $P$ be an oriented graded poset such that $\augm{P}$ is oriented thin. Then $\augm{(\sus{P})}$ is oriented thin.
\end{prop}
\begin{proof}
	Let $x, y \in \augm{(\sus{P})}$ such that $x \leq y$ and $\codim{x}{y} = 2$.
	If $x = \bot$, then $y = \augm{(\sus{y'})}$ for some $y' \in \grade{0}{P}$, so by construction the interval $[x, y]$ is of the form
\[\begin{tikzcd}
	& {y} \\
	{\augm{(\bot^+)}} && {\augm{(\bot^-)}} \\
	& {\bot}
	\arrow["+"', from=1-2, to=2-1]
	\arrow["-", from=1-2, to=2-3]
	\arrow["+"', from=2-1, to=3-2]
	\arrow["+", from=2-3, to=3-2]
\end{tikzcd}\]
	in $\augm{(\sus{P})}$.
	Otherwise, $x = \augm{x'}$ and $y = \augm{y'}$ for some $x', y' \in \sus{P}$, and the interval $[x, y]$ is isomorphic to the interval $[x', y']$ in $\sus{P}$.
	If $x' = \bot^\alpha$ for some $\alpha \in \set{+, -}$, then $[x', y']$ is, by construction, isomorphic to the interval $[\bot, y']$ in $\augm{P}$ at the level of the underlying posets.
	Since $\augm{P}$ is oriented thin, $[x', y']$ is of the form
\[\begin{tikzcd}
	& {y'} \\
	{z_1} && {z_2} \\
	& {\bot^\alpha}
	\arrow["+"', from=1-2, to=2-1]
	\arrow["-", from=1-2, to=2-3]
	\arrow["\alpha"', from=2-1, to=3-2]
	\arrow["\alpha", from=2-3, to=3-2]
\end{tikzcd}\]
	for some elements $z_1, z_2 \in \sus{P}$.
	Finally, if $\dim{x'} > 1$, then $x' = \sus{x''}$ and $y' = \sus{y''}$ for some $x'', y'' \in P$, and the interval $[x', y']$ is isomorphic to the interval $[x'', y'']$ in $P$.
	We conclude by oriented thinness of $\augm{P}$.
\end{proof}

\begin{lem} \label{lem:suspension_preserves_inclusions}
	Let $\imath\colon P \incl Q$ be an inclusion of oriented graded posets.
	Then $\sus{\imath}\colon \sus{P} \to \sus{Q}$ is an inclusion.
\end{lem}
\begin{proof}
	By inspection of the definition.
\end{proof}

\begin{rmk}
	Since $\sus{}$ preserves inclusions, we can let it act on closed subsets $U \subseteq P$, producing closed subsets $\sus{U} \subseteq \sus{P}$.
\end{rmk}

\begin{lem} \label{lem:suspension_preserves_connected_colimits}
	Let $\fun{F}\colon \smcat{C} \to \ogpos$ be a connected diagram of inclusions of oriented graded posets, and let $\gamma$ be a colimit cone under $\fun{F}$ whose components are all inclusions.
	Then $\sus{\gamma}$ is a colimit cone under $\sus{\fun{F}}$.
\end{lem}
\begin{proof}
	By Lemma \ref{lem:suspension_preserves_inclusions}, $\sus{\fun{F}}$ is a diagram of inclusions and the components of $\sus{\gamma}$ are all inclusions, so by Lemma
	\ref{lem:colimits_in_posclos} and Lemma 
	\ref{lem:reflected_colimits_in_ogpos} it suffices to prove that the underlying cone in $\poscat$ is a colimit cone.
	Let $P$ be the tip of $\gamma$ and let $\eta$ be a cone under $\fun{U}\sus{\fun{F}}$ with tip $Q$.
	The functions 
	\begin{align*}
		\imath_c\colon \fun{U}\fun{F}c & \to \fun{U}\sus{\fun{F}c}, \\
		x & \mapsto \sus{x}
	\end{align*}
	are order-preserving injections of posets for each object $c$ in $\smcat{C}$, so the restrictions of the components $\eta_c$ along these injections form a cone under $\fun{UF}$ in $\poscat$.
	By the universal property of $\fun{U}\gamma$, there is a unique order-preserving map $f\colon \fun{U}P \to Q$ such that $f \after \fun{U}\gamma_c = \eta_c \after \imath_c$ for all objects $c$ in $\smcat{C}$.
	Let $\overbar{c}$ be an arbitrary object of $\smcat{C}$; note that $\smcat{C}$ is non-empty, since it is connected.
	We let
	\begin{align*}
		f'\colon \fun{U}\sus{P} & \to Q, \\
		x & \mapsto \begin{cases}
			f(x') & \text{if $x = \sus{x'}$, $x' \in P$}, \\
			\eta_{\overbar{c}}(\bot^\alpha) & 
			\text{if $x = \bot^\alpha$, $\alpha \in \set{+, -}$}.
		\end{cases}
	\end{align*}
	We claim that $f'$ satisfies $f' \after \fun{U}\sus{\gamma_c} = \eta_c$ for all objects $c$ in $\smcat{C}$.
	This is true by construction when $c = \overbar{c}$.
	For any other $c$, since $\smcat{C}$ is connected, there exists a zig-zag of morphisms
\[\begin{tikzcd}
	\overbar{c} \equiv c_0 & {c_1} & {c_2} & \ldots & {c_{m-2}} & {c_{m-1}} & {c_m \equiv c} 
	\arrow["{f_0}", from=1-2, to=1-1]
	\arrow["{f_1}"', from=1-2, to=1-3]
	\arrow[from=1-4, to=1-3]
	\arrow[from=1-4, to=1-5]
	\arrow["{f_{m-1}}", from=1-6, to=1-5]
	\arrow["{f_m}"', from=1-6, to=1-7]
\end{tikzcd}\]
and $\sus{\fun{F}f_i}(\bot^\alpha) = \bot^\alpha$ for all $i \in \set{1, \ldots, m}$ and $\alpha \in \set{+, -}$.
Because the $\sus{\fun{F}f_i}$ are all injective and $\eta$ is a cone under $\fun{U}\sus{\fun{F}}$, it follows that
	\[
		\eta_{c_i}(\bot^\alpha) = \eta_{c_j}(\bot^\alpha)
	\]
	for all $i, j \in \set{0, \ldots, m}$ and $\alpha \in \set{+, -}$, hence $\eta_c(\bot^\alpha) = \eta_{\overbar{c}}(\bot^\alpha)$ for all $\alpha \in \set{+, -}$.
	This proves that $f' \after \fun{U}\sus{\gamma_c} = \eta_c$ for all objects $c$ in $\smcat{C}$, and $f'$ is evidently unique with this property.
	Thus $\fun{U}\sus{\gamma}$ is a colimit cone in $\poscat$ and $\sus{\gamma}$ is a colimit cone in $\ogpos$.
\end{proof}

\begin{rmk}
	Note that it is necessary to assume $\smcat{C}$ is connected in Lemma \ref{lem:suspension_preserves_connected_colimits}: $\sus{}$ does not preserve coproducts nor the initial object.
\end{rmk}

\begin{lem} \label{lem:suspension_of_union_and_intersection}
	Let $U, V$ be closed subsets of an oriented graded poset.
	Then
	\begin{enumerate}
		\item $\sus{(U \cup V)}	= \sus{U} \cup \sus{V}$,
		\item $\sus{(U \cap V)} = \sus{U} \cap \sus{V}$.
	\end{enumerate}
\end{lem}
\begin{proof}
	We have
	\begin{align*}
		& \sus{(U \cup V)} = \set{\bot^+, \bot^-} + \set{\sus{x} \mid x \in U \cup V} = \\
				 & = (\set{\bot^+, \bot^-} + \set{\sus{x} \mid x \in U}) \cup (\set{\bot^+, \bot^-} + \set{\sus{x} \mid x \in V}) = \sus{U} \cup \sus{V}.
	\end{align*}
	The case of intersections is similar.
\end{proof}

\begin{lem} \label{lem:suspension_of_closed_subset}
	Let $U \subseteq P$ be a non-empty closed subset of an oriented graded poset.
	Then
	\[
		\sus{U} = \clset{\sus{x} \mid x \in U}.
	\]
\end{lem}
\begin{proof}
	We have $\set{\sus{x} \mid x \in U} \subseteq \sus{U}$ and $\sus{U}$ is closed, so one inclusion follows from Lemma \ref{lem:closure_is_monotonic}.
	For the other, it suffices to show that 
	\[\set{\bot^+, \bot^-} \subseteq \clset{\sus{x} \mid x \in U}.\]
	But if $U$ is non-empty, then there exists $x \in \grade{0}{U}$, and $\set{\bot^+, \bot^-} = \faces{}{}\sus{x}$.
\end{proof}

\begin{comm}
	Observe that Lemma \ref{lem:suspension_of_closed_subset} fails when $U$ is empty, since $\sus{\varnothing} = \set{\bot^+, \bot^-} \neq \clset{\varnothing} = \varnothing$.
\end{comm}

\begin{lem} \label{lem:faces_of_suspension}
	Let $U$ be a non-empty closed subset of an oriented graded poset, $n \in \mathbb{N}$, $\alpha \in \set{+, -}$.
	Then
	\[
		\faces{n}{\alpha}\sus{U} = \begin{cases}
			\set{\bot^\alpha} & \text{if $n = 0$}, \\
			\set{\sus{x} \mid x \in \faces{n-1}{\alpha}U} & \text{if $n > 0$}.
		\end{cases}
	\]
\end{lem}
\begin{proof}
	Let $x \in \grade{n}{\sus{U}}$.
	If $n = 0$, then $x = \bot^\beta$ for some $\beta \in \set{+, -}$.
	Then 
	\[
		\cofaces{}{-\alpha}x \cap \sus{U} = \begin{cases}
			\set{\sus{y} \mid y \in \grade{0}{U}} & \text{if $\beta = -\alpha$}, \\
			\varnothing & \text{if $\beta = \alpha$}.
		\end{cases}
	\]
	Since $U$ is non-empty, $\grade{0}{U}$ is non-empty, so $x \in \faces{0}{\alpha}\sus{U}$ if and only if $x = \bot^\alpha$.
	Now suppose $n > 0$.
	Then $x = \sus{x'}$ for some $x' \in \grade{n-1}{U}$, and
	\[
		\cofaces{}{-\alpha}x \cap \sus{U} = \set{\sus{y} \mid y \in \cofaces{}{-\alpha}x' \cap U},
	\]
	which is empty if and only if $x' \in \faces{n-1}{\alpha}U$.
\end{proof}

\begin{cor} \label{cor:boundary_of_suspension}
	Let $U$ be a non-empty closed subset of an oriented graded poset, $n \in \mathbb{N}$, $\alpha \in \set{+, -}$. Then
	\[
		\bound{n}{\alpha}\sus{U} = \begin{cases}
			\set{\bot^\alpha} & \text{if $n = 0$}, \\
			\sus{\bound{n-1}{\alpha}U} & \text{if $n > 0$}.
		\end{cases}
	\]
\end{cor}
\begin{proof}
	Follows from Lemma \ref{lem:suspension_of_closed_subset}, Lemma \ref{lem:faces_of_suspension}, and Lemma \ref{lem:maximal_vs_faces}.
\end{proof}

\begin{lem} \label{lem:suspension_preserves_globular_and_round}
	Let $U$ be a non-empty oriented graded poset.
	Then
	\begin{enumerate}
		\item if $U$ is globular, then $\sus{U}$ is globular,
		\item if $U$ is round, then $\sus{U}$ is round.
	\end{enumerate}
\end{lem}
\begin{proof}
	Suppose $U$ is globular, let $k, n \in \mathbb{N}$, $\alpha, \beta \in \set{+, -}$, and suppose $k < n$.
	By Corollary \ref{cor:boundary_of_suspension},
	\[
		\bound{k}{\alpha}(\bound{n}{\beta}\sus{U}) = \bound{k}{\alpha}(\sus{\bound{n-1}{\beta}U}) = \begin{cases}
			\sus{\bound{k-1}{\alpha}(\bound{n-1}{\beta}U)} & \text{if $k > 0$}, \\
			\set{\bot^\alpha} & \text{if $k = 0$}.
		\end{cases}
	\]
	If $k = 0$, this is always equal to $\bound{k}{\alpha}\sus{U}$.
	If $k > 0$, by globularity of $U$, this is equal to $\sus{\bound{k-1}{\alpha}U}$, which is equal to $\bound{k}{\alpha}\sus{U}$.

	Next, suppose that $U$ is round, and let $n < \dim{U}$.
	For $n = 0$, we have
	\[
		\bound{0}{+}\sus{U} \cap \bound{0}{-}\sus{V} = \set{\bot^+} \cap \set{\bot^-} = \varnothing = \bound{-1}{}\sus{U}.
	\]
	By roundness of $U$ and Lemma \ref{lem:suspension_of_union_and_intersection}, for $n = 1$ we have
	\begin{align*}
		\bound{1}{+}\sus{U} \cap \bound{1}{-}\sus{U} & = \sus{\bound{0}{+}U} \cap \sus{\bound{0}{-}U} = \sus{(\bound{0}{+}U \cap \bound{0}{-}U)} = \\
							     & = \sus{\varnothing} = \set{\bot^+, \bot^-} = \bound{0}{}\sus{U}
	\end{align*}
	and for $n > 1$ we have
	\begin{align*}
		\bound{n}{+}\sus{U} \cap \bound{n}{-}\sus{U} & = \sus{\bound{n-1}{+}U} \cap \sus{\bound{n-1}{-}U} = \sus{(\bound{n-1}{+}U \cap \bound{n-1}{-}U)} = \\ 
		& = \sus{\bound{n-2}{}U} = \bound{n-1}{}\sus{U}.
	\end{align*}
	This proves that $\sus{U}$ is round.
\end{proof}

\begin{prop} \label{prop:suspension_of_molecules}
	Let $U$ be a molecule.
	Then
	\begin{enumerate}
		\item $\sus{U}$ is a molecule,
		\item if $U$ is isomorphic to $V \cp{k} W$, then $\sus{U}$ is isomorphic to $\sus{V} \cp{k+1} \sus{W}$,
		\item if $U$ is isomorphic to $V \celto W$, then $\sus{U}$ is isomorphic to $\sus{V} \celto \sus{W}$.
	\end{enumerate}
\end{prop}
\begin{proof}
We proceed by induction on the construction of $U$.
If $U$ was produced by (\textit{Point}), then $U = 1$ and $\sus{U}$ is isomorphic to $1 \celto 1$ by inspection.

If $U$ was produced by (\textit{Paste}), then it is of the form $V \cp{k} W$ for some molecules $V, W$ and $k < \min \set{\dim{V}, \dim{W}}$.
Let 
\[\begin{tikzcd}
	\bound{k}{+}V & \bound{k}{-}W & W \\
	V && U
	\arrow["\sim", hook, from=1-1, to=1-2]
	\arrow[hook, from=1-2, to=1-3]
	\arrow[hook', from=1-1, to=2-1]
	\arrow[hook, from=2-1, to=2-3]
	\arrow[hook', from=1-3, to=2-3]
	\arrow["\lrcorner"{anchor=center, pos=0.125, rotate=180}, draw=none, from=2-3, to=1-1]
\end{tikzcd}\]
be the pushout diagram exhibiting $U$ as $V \cp{k} W$.
Since pushouts are connected colimits, by Lemma \ref{lem:suspension_preserves_connected_colimits} together with Corollary \ref{cor:boundary_of_suspension}, $\sus{}$ maps this diagram onto a pushout diagram
\[\begin{tikzcd}
	\bound{k+1}{+}\sus{V} & \bound{k+1}{-}\sus{W} & \sus{W} \\
	\sus{V} && \sus{U}
	\arrow["\sim", hook, from=1-1, to=1-2]
	\arrow[hook, from=1-2, to=1-3]
	\arrow[hook', from=1-1, to=2-1]
	\arrow[hook, from=2-1, to=2-3]
	\arrow[hook', from=1-3, to=2-3]
	\arrow["\lrcorner"{anchor=center, pos=0.125, rotate=180}, draw=none, from=2-3, to=1-1]
\end{tikzcd}\]
where $\sus{V}$ and $\sus{W}$ are molecules by the inductive hypothesis.
This diagram exhibits $\sus{U}$ as $\sus{V} \cp{k+1} \sus{W}$.

Finally, if $U$ was produced by (\textit{Atom}), then it is of the form $V \celto W$ for some round molecules $V, W$ of the same dimension $n \eqdef \dim{U} - 1$.
Let 
\[\begin{tikzcd}
	\bound{}{}V & \bound{}{}W & W \\
	V && \bound{}{}U
	\arrow["\sim", hook, from=1-1, to=1-2]
	\arrow[hook, from=1-2, to=1-3]
	\arrow[hook', from=1-1, to=2-1]
	\arrow[hook, from=2-1, to=2-3]
	\arrow[hook', from=1-3, to=2-3]
	\arrow["\lrcorner"{anchor=center, pos=0.125, rotate=180}, draw=none, from=2-3, to=1-1]
\end{tikzcd}\]
be the pushout diagram exhibiting $\bound{}{}U$ as a gluing of $V$ and $W$.
Using Lemma \ref{lem:suspension_preserves_connected_colimits} and Corollary \ref{cor:boundary_of_suspension} once more, $\sus{}$ maps this diagram onto a pushout diagram
\[\begin{tikzcd}
	\bound{}{}\sus{V} & \bound{}{}\sus{W} & \sus{W} \\
	\sus{V} && \bound{}{}\sus{U}
	\arrow["\sim", hook, from=1-1, to=1-2]
	\arrow[hook, from=1-2, to=1-3]
	\arrow[hook', from=1-1, to=2-1]
	\arrow[hook, from=2-1, to=2-3]
	\arrow[hook', from=1-3, to=2-3]
	\arrow["\lrcorner"{anchor=center, pos=0.125, rotate=180}, draw=none, from=2-3, to=1-1]
\end{tikzcd}\]
where $\sus{V}$ and $\sus{W}$ are both round molecules of the same dimension by the inductive hypothesis and Lemma \ref{lem:suspension_preserves_globular_and_round}.
This exhibits $\bound{}{}\sus{U}$ as $\bound{}{}(\sus{V} \celto \sus{W})$, hence exhibits $\sus{U}$ as $\sus{V} \celto \sus{W}$.
\end{proof}

\begin{cor} \label{cor:suspension_of_rdcpx}
	Let $P$ be a regular directed complex.
	Then $\sus{P}$ is a regular directed complex.
\end{cor}
\begin{proof}
	Let $x \in \sus{P}$.
	Either $x = \bot^\alpha$ for some $\alpha \in \set{+, -}$, so $\clset{x} = \set{\bot^\alpha}$ is a 0\nbd dimensional atom, or $x = \sus{x'}$ for some $x' \in P$, in which case 
	\[\clset{x} = \sus{\clset{x'}}\] 
	is an atom by Proposition \ref{prop:suspension_of_molecules}.
\end{proof}

\begin{prop} \label{prop:suspension_is_a_functor_on_maps}
	Let $p\colon P \to Q$ be a map of regular directed complexes.
	Then
	\begin{align*}
		\sus{p}\colon \sus{P} & \to \sus{Q}, \\
		x & \mapsto \begin{cases}
			\sus{p(x')} & \text{if $x = \sus{x'}$, $x' \in P$, } \\
			\bot^\alpha & \text{if $x = \bot^\alpha$, $\alpha \in \set{+, -}$},
		\end{cases}
	\end{align*}
	is a map of regular directed complexes.
	This assignment determines a faithful endofunctor $\sus{}$ on $\rdcpxmap$, such that the diagram of functors
\begin{equation} \label{eq:suspension_functors_commute}
	\begin{tikzcd}
	\ogpos && \rdcpx && \rdcpxmap \\
	\ogpos && \rdcpx && \rdcpxmap
	\arrow[hook', from=1-3, to=1-1]
	\arrow[hook, from=1-3, to=1-5]
	\arrow["{\sus{}}", from=1-3, to=2-3]
	\arrow["{\sus{}}", from=1-1, to=2-1]
	\arrow["{\sus{}}", from=1-5, to=2-5]
	\arrow[hook, from=2-3, to=2-5]
	\arrow[hook', from=2-3, to=2-1]
\end{tikzcd}
\end{equation}
	commutes.
\end{prop}
\begin{proof}
	By Corollary \ref{cor:suspension_of_rdcpx}, $\sus{P}$ and $\sus{Q}$ are regular directed complexes.
	Functoriality, faithfulness, and commutativity of (\ref{eq:suspension_functors_commute}) are evident by inspection of the definitions, so it suffices to show that $\sus{p}$ is a map.
	Let $x \in \sus{P}$, $n \in \mathbb{N}$, and $\alpha \in \set{+, -}$.
	If $x = \bot^\beta$ for some $\beta \in \set{+, -}$, then
	\[
		\sus{p}(\bound{n}{\alpha}x) = \sus{p}(\set{\bot^\beta}) = \set{\bot^\beta} = \bound{n}{\alpha}\sus{p}(x),
	\]
	while if $x = \sus{x'}$ for some $x' \in P$, then for $n > 0$
	\[
		\sus{p}(\bound{n}{\alpha}x) = \sus{p}(\sus{\bound{n-1}{\alpha}x'}) = \sus{(p(\bound{n-1}{\alpha}x'))} = \sus{\bound{n-1}{\alpha}x'} = \bound{n}{\alpha}x
	\]
	while for $n = 0$
	\[
		\sus{p}(\bound{0}{\alpha}x) = \sus{p}(\set{\bot^\alpha}) = \set{\bot^\alpha} = \bound{0}{\alpha}\sus{p}(x)
	\]
	where in both cases we used Corollary \ref{cor:boundary_of_suspension}.
	Finally, let $u, v \in \bound{n}{\alpha}x$ and suppose $z \leq \sus{p}(u)$ and $z \leq \sus{p}(v)$.
	If $z = \sus{z'}$ for some $z \in P$, then necessarily also $u = \sus{u'}$, $v = \sus{v'}$, and $x = \sus{x'}$ for some $x' \in P$ and $u', v' \in \bound{n-1}{\alpha}x'$ such that $z' \leq u'$ and $z' \leq v'$.
	Since $\restr{p}{\bound{n-1}{\alpha}x'}$ is final onto its image, there is a zig-zag
	\[
		u' \leq u_1 \geq \ldots \leq u_m \geq v'
	\]
	in $\bound{n-1}{\alpha}x'$ such that $z' \leq p(u_i)$ for all $i \in \set{1, \ldots, m}$.
	This induces a zig-zag
	\[
		u \leq \sus{u_1} \geq \ldots \leq \sus{u_m} \geq v
	\]
	in $\bound{n}{\alpha}x$ such that $z \leq \sus{p}(\sus{u_i})$ for all $i \in \set{1, \ldots, m}$.
	If, instead, $z = \bot^\beta$ for some $\beta \in \set{+, -}$, then
	\[
		u \geq \bot^\beta \leq v
	\]
	is a zig-zag in $\bound{n}{\alpha}x$ such that $z \leq \sus{p}(\bot^\beta)$.
	This proves that $\restr{\sus{p}}{\bound{n}{\alpha}x}$ is final onto its image, completing the proof that $\sus{p}$ is a map.
\end{proof}

\begin{prop} \label{prop:suspension_is_a_functor_on_comaps}
	Let $c\colon P \to Q$ be a comap of regular directed complexes.
	Then
	\begin{align*}
		\sus{c}\colon \sus{P} & \to \sus{Q}, \\
		x & \mapsto \begin{cases}
			\sus{c(x')} & \text{if $x = \sus{x'}$, $x' \in P$, } \\
			\bot^\alpha & \text{if $x = \bot^\alpha$, $\alpha \in \set{+, -}$},
		\end{cases}
	\end{align*}
	is a comap of regular directed complexes.
	This assignment determines a faithful endofunctor $\sus{}$ on $\rdcpxcomap$.
\end{prop}
\begin{proof}
	By Corollary \ref{cor:suspension_of_rdcpx}, both $\sus{P}$ and $\sus{Q}$ are regular directed complexes.
	Functoriality and faithfulness are straightforward, so it suffices to prove that $\sus{c}$ is a comap.
	Let $x \in \sus{Q}$, $n \in \mathbb{N}$, and $\alpha \in \set{+, -}$.
	If $x = \bot^\beta$ for some $\beta \in \set{+, -}$, then
	\[
		\invrs{(\sus{c})}\bound{n}{\alpha}x = \invrs{(\sus{c})}\set{\bot^\beta} = \set{\bot^\beta} = \bound{n}{\alpha}\invrs{(\sus{c})}\clset{x},
	\]
	which is a 0\nbd dimensional atom.
	If $x = \sus{x'}$ for some $x' \in Q$, then for $n > 0$
	\[
		\invrs{(\sus{c})}\bound{n}{\alpha}x = \invrs{(\sus{c})}\sus{\bound{n-1}{\alpha}x'} = \sus{(\bound{n-1}{\alpha}\invrs{c}{\clset{x'}})} = \bound{n}{\alpha}\invrs{(\sus{c})}\clset{x},
	\]
	which is a molecule because $\invrs{c}{\clset{x'}}$ is, while for $n = 0$
	\[
		\invrs{(\sus{c})}\bound{0}{\alpha}x = \invrs{(\sus{c})}\set{\bot^\alpha} = \set{\bot^\alpha} = \bound{0}{\alpha}\invrs{(\sus{c})}\clset{x},
	\]
	which again is a 0\nbd dimensional atom.
	This proves that $\sus{c}$ is a comap.
\end{proof}

\begin{dfn}[Suspension of a strict $\omega$-category] \index{strict $\omega$-category!suspension} \index{$\sus{X}$} \index{suspension!of a strict $\omega$-category}
	Let $X$ be a strict $\omega$-category.
	The \emph{suspension of $X$} is the strict $\omega$-category $\sus{X}$ whose set of cells is
	\[
		\set{\sus{t} \mid t \in X} + \set{\bot^+, \bot^-},
	\]
	with the boundary operators defined, for all $t \in \sus{X}$, $n \in \mathbb{N}$, $\alpha \in \set{+, -}$, by
	\[
		\bound{n}{\alpha}t \eqdef \begin{cases}
			\sus{\bound{n-1}{\alpha}t'} & \text{if $n > 0$, $t = \sus{t'}$, $t' \in X$}, \\
			\bot^\alpha & \text{if $n = 0$, $t = \sus{t'}$, $t' \in X$}, \\
			t & \text{if $t \in \set{\bot^+, \bot^-}$},
		\end{cases}
	\]
	and the $k$\nbd composition operations defined, for all $k \in \mathbb{N}$ and $k$\nbd composable pairs $t, u$ in $\sus{X}$, by
	\[
		t \cp{k} u \eqdef \begin{cases}
			\sus{(t' \cp{k-1} u')} & \text{if $t = \sus{t'}$, $u = \sus{u'}$, $t', u' \in X$}, \\
			t & \text{if $u \in \set{\bot^+, \bot^-}$}, \\
			u & \text{if $t \in \set{\bot^+, \bot^-}$}.
		\end{cases}
	\]
	Suspension extends to an endofunctor $\sus{}\colon \omegacat \to \omegacat$.
\end{dfn}

\begin{comm}
	This is the $\omega$\nbd categorical suspension as defined, for example, in \cite[Section 2.2]{ozornova2023quillen}.
	As noted there, this is not exactly the same as the definition in 
	\cite[Section B.6]{ara2020joint}, where suspended cells are also reversed.
\end{comm}

\begin{prop} \label{prop:compatibility_with_omegacat_suspension}
	Let $P$ be a regular directed complex.
	Then $\sus{\molecin{P}}$ is naturally isomorphic to $\molecin{\sus{P}}$.
\end{prop}
\begin{proof}
	Let $\varphi\colon \sus{\molecin{P}} \to \molecin{\sus{P}}$ be defined by
	\begin{align*}
		\bot^\alpha & \mapsto \isocl{\set{\bot^\alpha} \incl \sus{P}}, \\ 
		\sus{\isocl{f\colon U \to P}} & \mapsto \isocl{\sus{f}\colon \sus{U} \to \sus{P}}
	\end{align*}
	for all $\alpha \in \set{+, -}$ and cells $\isocl{f\colon U \to P}$ in $\molecin{P}$.
	Then $\varphi$ is a morphism of reflexive $\omega$\nbd graphs by 
	Corollary \ref{cor:boundary_of_suspension} and a strict functor of strict $\omega$\nbd categories by Proposition \ref{prop:suspension_of_molecules}.
	It is also injective due to faithfulness of $\sus{}$ on $\ogpos$, so it suffices to show that it is surjective.
	By Corollary \ref{cor:basis_of_omegacat_presented_by_rdcpx}, the set $\gener{S} \eqdef \set{\isocl{\clset{x} \incl \sus{P}} \mid x \in \sus{P}}$ is a basis for $\molecin{\sus{P}}$, so we may proceed by structural induction on $\cspan{\gener{S}}$.
	If $t \in \gener{S}$, then by definition of $\sus{P}$ it is either $\isocl{\set{\bot^\alpha} \incl \sus{P}} = \varphi(\bot^\alpha)$ for some $\alpha \in \set{+, -}$, or it is $\isocl{\clset{\sus{x'}} \incl \sus{P}} = \varphi{\isocl{\clset{x'} \incl P}}$ for some $x' \in P$.
	
	Suppose that $t = u \cp{k} v$; we may assume $k < \min \set{\dim{u}, \dim{v}}$, for otherwise $t = u$ or $t = v$.
	By the inductive hypothesis, $u$ and $v$ are in the image of $\varphi$, and since they are not $0$\nbd dimensional, they must be of the form $\isocl{\sus{f}\colon \sus{U} \to \sus{P}}$ and $\isocl{\sus{g}\colon \sus{V} \to \sus{P}}$, respectively.
	Then $\bound{0}{+}u = \isocl{\set{\bot^+} \incl \sus{P}}$ is not equal to $\bound{0}{-}v = \isocl{\set{\bot^-} \incl \sus{P}}$, so $k > 0$.
	By Corollary \ref{cor:boundary_of_suspension}, since $\isocl{\sus{f}}$ and $\isocl{\sus{g}}$ are $k$\nbd composable, it follows that $\isocl{f}$ and $\isocl{g}$ are $(k-1)$\nbd composable, so by Proposition \ref{prop:suspension_of_molecules} $t = \isocl{\sus{(f \cp{k-1} g)}} = \varphi(\sus{\isocl{f \cp{k-1} g}})$.
	This proves that $\varphi$ is surjective, so it is an isomorphism of strict $\omega$\nbd categories.
	Naturality is straightforward.
\end{proof}

\begin{comm}
	To avoid stating multiple variants of each naturality result, when the category is unspecified, \emph{naturally isomorphic} should be interpreted as \emph{naturally} with respect to any of the notions of morphisms that we have considered, when well-typed.
	For example, when restricted to regular directed complexes, the natural isomorphism of Proposition \ref{prop:compatibility_with_omegacat_suspension} is also natural over maps and comaps.
\end{comm}

\begin{dfn}[Suspension of an augmented chain complex] \index{chain complex!suspension} \index{$\sus{C}$} \index{suspension!of an augmented chain complex}
	Let $C$ be an augmented chain complex.
	The \emph{suspension of $C$} is the augmented chain complex $\sus{C}$ defined by
	\[
		\grade{n}{\sus{C}} \eqdef \begin{cases}
				\freeab{\set{\bot^+, \bot^-}} & \text{if $n = 0$,} \\
				\grade{n-1}{C} & \text{if $n > 0$,}
			\end{cases}
	\]
	with $\der\colon \grade{n}{\sus{C}} \to \grade{n-1}{\sus{C}}$ defined by
	\[
		\begin{cases}
			x \mapsto \eau(x)(\bot^+ - \bot^-) & \text{if $n = 1$,} \\
			\der\colon \grade{n-1}{C} \to \grade{n-2}{C} & \text{if $n > 1$,}
		\end{cases}
	\]
	and $\eau\colon \grade{0}{\sus{C}} \to \mathbb{Z}$ defined by $\bot^\alpha \mapsto 1$ for all $\alpha \in \set{+, -}$.
	The suspension extends to an endofunctor $\sus{}$ on $\chaug$.
\end{dfn}

\begin{prop} \label{prop:suspension_of_otgpos_and_chaug}
	Let $P$ be an oriented graded poset such that $\augm{P}$ is oriented thin.
	Then $\freeab{(\sus{P})}$ is naturally isomorphic to $\sus{(\freeab{P})}$.
\end{prop}
\begin{proof}
	By Proposition \ref{prop:suspension_of_oriented_thin}, $\augm{(\sus{P})}$ is oriented thin, so $\freeab{(\sus{P})}$ is well-defined.
	For each $n \in \mathbb{N}$, $\grade{n}{\freeab{(\sus{P})}}$ and $\grade{n}{\sus{(\freeab{P})}}$ are by construction free abelian groups on isomorphic sets.
	The fact that the evident isomorphisms determine an isomorphism of augmented chain complexes is a straightforward check, as is naturality.
\end{proof}


\section{Joins} \label{sec:join}

\begin{guide}
	In this section, we define the join of oriented graded posets.
	This is defined in analogy with the cellular join of posets, with the augmentation and Gray product of oriented graded posets playing the role of their undirected counterparts.
	However, in order to prove that the join of two molecules is a molecule (Proposition 
	\ref{prop:join_of_molecules}), we rely on a relation between the join of two oriented graded posets and the Gray product of their \emph{suspensions} (Lemma 
	\ref{lem:joins_and_gray_of_suspension}).
	
	Finally, we prove that the join is part of a monoidal structure on $\ogpos$, $\rdcpxmap$, and $\rdcpxcomap$, compatibly on their overlaps, and that it is compatible with an analogous construction on augmented chain complexes.
\end{guide}

\begin{prop} \label{prop:gray_product_ogposbot}
	The monoidal structure $(\ogpos, \gray, 1)$ restricts to a monoidal structure on $\ogposbot$.
\end{prop}
\begin{proof}
	Let $P$, $Q$ be oriented graded posets with positive least elements $\bot_P$ and $\bot_Q$, respectively.
	Then $(\bot_P, \bot_Q)$ is the least element of $P \gray Q$.
	Moreover,
	\begin{align*}
		\cofaces{}{}(\bot_P, \bot_Q) & = 
		\cofaces{}{}\bot_P \times \set{\bot_Q} + \set{\bot_P} \times \cofaces{}{}\bot_Q = \\
					     & =
		\cofaces{}{+}\bot_P \times \set{\bot_Q} + \set{\bot_P} \times \cofaces{}{+}\bot_Q = \cofaces{}{+}(\bot_P, \bot_Q),
	\end{align*}
	which proves that $(\bot_P, \bot_Q)$ is a positive least element.
\end{proof}

\begin{dfn}[Join of oriented graded posets] \index{oriented graded poset!join} \index{regular directed complex!join} \index{$P \join Q$} \index{join!of oriented graded posets}
	Let $P$, $Q$ be oriented graded posets.
	The \emph{join of $P$ and $Q$} is the oriented graded poset $P \join Q \eqdef \dimin{(\augm{P}  \gray \augm{Q})}$.
\end{dfn}

\begin{exm}[A non-symmetric join of molecules] \index[counterex]{A non-symmetric join of molecules} \label{exm:nonsymmetry_join}
	Since the join is defined in terms of the Gray product, which is non-symmetric, it can be expected to also be non-symmetric, and that is indeed the case.
	For example, $1 \join \thearrow{2}$ and $\thearrow{2} \join 1$ are isomorphic to the oriented face posets of 
\[\begin{tikzcd}[column sep=small]
	& \bullet & \bullet \\
	\bullet &&& \bullet
	\arrow[curve={height=-6pt}, from=2-1, to=1-2]
	\arrow[from=1-2, to=1-3]
	\arrow[""{name=0, anchor=center, inner sep=0}, curve={height=6pt}, from=2-1, to=1-3]
	\arrow[curve={height=-6pt}, from=1-3, to=2-4]
	\arrow[""{name=1, anchor=center, inner sep=0}, curve={height=12pt}, from=2-1, to=2-4]
	\arrow[curve={height=6pt}, shorten <=5pt, Rightarrow, from=1, to=1-3]
	\arrow[shorten <=2pt, Rightarrow, from=0, to=1-2]
\end{tikzcd} \quad \quad \text{and} \quad \quad 
\begin{tikzcd}[column sep=small]
	& \bullet & \bullet \\
	\bullet &&& \bullet
	\arrow[curve={height=-6pt}, from=1-3, to=2-4]
	\arrow[from=1-2, to=1-3]
	\arrow[""{name=0, anchor=center, inner sep=0}, curve={height=6pt}, from=1-2, to=2-4]
	\arrow[curve={height=-6pt}, from=2-1, to=1-2]
	\arrow[""{name=1, anchor=center, inner sep=0}, curve={height=12pt}, from=2-1, to=2-4]
	\arrow[curve={height=-6pt}, shorten <=5pt, Rightarrow, from=1, to=1-2]
	\arrow[shorten <=2pt, Rightarrow, from=0, to=1-3]
\end{tikzcd}
\]
	respectively, which are not isomorphic to each other.
\end{exm}

\begin{prop} \label{prop:join_monoidal_structure}
	The join of oriented graded posets extends to a unique monoidal structure $(\ogpos, \join, \varnothing)$ on $\ogpos$ such that
	\begin{enumerate}
		\item $\dimin{(-)}\colon (\ogposbot, \gray, 1) \to (\ogpos, \join, \varnothing)$ is a strong monoidal functor,
		\item $\fun{U}\colon (\ogpos, \join, \varnothing) \to (\posclos, \join, \varnothing)$
		is a strict monoidal functor.
	\end{enumerate}
\end{prop}
\begin{proof}
	The requirement that $\fun{U}$ be strict monoidal fixes the monoidal structure uniquely on the underlying posets and maps, while $\augm{(-)}$ being strong monoidal fixes the orientations.
	That this is well-defined follows from Proposition \ref{prop:ogpos_equivalent_ogposbot} and Proposition \ref{prop:join_monoidal_structure_posclos}.
\end{proof}

\begin{comm}
	We will adopt the same notation for elements of the join of two oriented graded posets as for the cellular join of their underlying posets.
\end{comm}

\begin{lem} \label{lem:join_faces}
	Let $P$, $Q$ be oriented graded posets, $z \in P \join Q$.
	Then, for all $\alpha \in \set{+, -}$, $\faces{}{\alpha}z$ is equal to
	\[
	\begin{cases}
		\set{\inj{x'} \mid x' \in \faces{}{\alpha}x} & 
		\text{if $z = \inj{x}$, $x \in P$,} \\
		\set{\inr{y'} \mid y' \in \faces{}{\alpha}y} & 
		\text{if $z = \inr{y}$, $y \in Q$,} \\
		\set{\inr{y}} + \set{x \join y' \mid y' \in \faces{}{-}y} & 
		\text{if $z = x \join y$, $x \in \grade{0}{P}$, $\alpha = +$,} \\
		\set{x' \join y \mid x' \in \faces{}{\alpha}x} + \set{\inj{x}} & 
		\text{if $z = x \join y$, $y \in \grade{0}{Q}$, $\alpha = -(-)^{\dim{x}}$,} \\
		\set{x' \join y \mid x' \in \faces{}{\alpha}x} + & \text{} \\ 
		\quad + \set{x \join y' \mid y' \in \faces{}{-(-)^{\dim{x}}\alpha}y} &
		\text{if $z = x \join y$, $x \in P$, $y \in Q$, otherwise.}
	\end{cases}
\]
\end{lem}
\begin{proof}
	By a simple case distinction, based on the definition of the orientation on $\augm{P} \gray \augm{Q}$.
\end{proof}

\begin{prop} \label{prop:gray_product_oriented_thin}
	Let $P$, $Q$ be oriented thin graded posets.
	Then $P \gray Q$ is oriented thin.
\end{prop}
\begin{proof}
	An easier variant of the proof of Proposition \ref{prop:gray_product_augm_oriented_thin}.
\end{proof}

\begin{cor} \label{cor:gray_product_on_otgpos}
	The monoidal structure $(\ogposbot, \gray, 1)$ restricts to a monoidal structure on $\otgpos$.
\end{cor}

\begin{cor} \label{cor:join_of_oriented_thin}
	Let $P$, $Q$ be oriented graded posets such that $\augm{P}$ and $\augm{Q}$ are oriented thin.
	Then $\augm{(P \join Q)}$ is oriented thin.
\end{cor}

\begin{lem} \label{lem:join_preserves_inclusions}
	Let $\imath\colon P \incl P'$ and $j\colon Q \incl Q'$ be inclusions of oriented graded posets.
	Then $\imath \join j\colon P \join Q \to P' \join Q'$ is an inclusion.
\end{lem}
\begin{proof}
	Follows from Lemma \ref{lem:preservation_of_closed_embeddings}.
\end{proof}

\begin{rmk}
	Since the join preserves inclusions, we can let it act on closed subsets $U \subseteq P$, $V \subseteq Q$, producing closed subsets $U \join V \subseteq P \join Q$.
\end{rmk}

\begin{rmk}
	The monoidal structure $(\ogpos, \join, \varnothing)$ is \emph{semicocartesian} monoidal, that is, the monoidal unit is the initial object.
	Moreover, the unique morphism from the initial object is always an inclusion.
	It follows that, for all oriented graded posets $P$, $Q$, there are natural inclusions
	\[
		P \incl P \join Q, \quad \quad Q \incl P \join Q
	\]
	obtained by composing a right unitor with $P \join \varnothing \incl P \join Q$ and a left unitor with $\varnothing \join Q \incl P \join Q$.
\end{rmk}

\begin{lem} \label{lem:join_preserves_connected_colimits_ogpos}
	Let $P$ be an oriented graded poset, let $\fun{F}$ be a connected diagram of inclusions in $\ogpos$, and let $\gamma$ be a colimit cone under $\fun{F}$ whose components are all inclusions.
	Then
	\begin{enumerate}
		\item $P \join \gamma$ is a colimit cone under $P \join \fun{F}$,
		\item $\gamma \join P$ is a colimit cone under $\fun{F} \join P$.
	\end{enumerate}
\end{lem}
\begin{proof}
	By Lemma \ref{lem:reflected_colimits_in_ogpos}, $\fun{U}\gamma$ is a colimit cone in $\posclos$ whose components are all closed embeddings.
	Now, $P \join \gamma$ has an underlying cone $\fun{U}P \join \fun{U}\gamma$.
	By Lemma \ref{lem:join_preserves_connected_colimits}, $\fun{U}P \join \fun{U}\gamma$ is a colimit cone in $\posclos$ whose components are all closed embeddings by Lemma \ref{lem:preservation_of_closed_embeddings}.
	By Lemma \ref{lem:reflected_colimits_in_ogpos}, we conclude that $P \join \gamma$ is a colimit cone in $\ogpos$ whose components are all inclusions.
	A symmetrical argument proves the same for $\gamma \join P$.
\end{proof}

\begin{dfn}[Injection of join into Gray product of suspensions]
	Let $P$, $Q$ be oriented graded posets.
	We define an injective function
	\begin{align*}
		s\colon P \join Q & \to \sus{P} \gray \sus{Q}, \\
		z & \mapsto \begin{cases}
			(\sus{x}, \bot^+) & \text{if $z = \inj{x}$, $x \in P$}, \\
			(\bot^+, \sus{y}) & \text{if $z = \inr{y}$, $y \in Q$}, \\
			(\sus{x}, \sus{y}) & \text{if $z = x \join y$, $x \in P$, $y \in Q$}.
		\end{cases}
	\end{align*}
\end{dfn}

\begin{lem} \label{lem:joins_and_gray_of_suspension}
	Let $P$, $Q$ be oriented graded posets.
	Then
	\begin{enumerate}
		\item $s\colon P \join Q \to \sus{P} \gray \sus{Q}$ is an order-preserving map of the underlying posets,
		\item for all $z \in P \join Q$ and $\alpha \in \set{+, -}$, $s$ induces a bijection between $\cofaces{}{\alpha}z$ and $\cofaces{}{\alpha}s(z)$,
		\item for all closed subsets $U \subseteq \sus{P} \gray \sus{Q}$, $n \in \mathbb{N}$, and $\alpha \in \set{+, -}$,
			\[
				\faces{n}{\alpha}\invrs{s}U = \invrs{s}\faces{n+1}{\alpha}U, \quad \quad
				\bound{n}{\alpha}\invrs{s}U = \invrs{s}\bound{n+1}{\alpha}U.
			\]
	\end{enumerate}
\end{lem}
\begin{proof}
	A simple case distinction shows that, for all $z \in P \join Q$ and $\alpha \in \set{+, -}$, 
	\begin{enumerate}
		\item $\dim{s(z)} = \dim{z} + 1$, and
		\item $s$ induces a bijection between $\cofaces{}{\alpha}z$ and $\cofaces{}{\alpha} s(z)$.
	\end{enumerate}
	By Lemma \ref{lem:hasse_diagram_paths}, this suffices to prove that $s$ is order-preserving.
	Given a closed subset $U \subseteq \sus{P} \gray \sus{Q}$, it follows by Lemma 
	\ref{lem:inverse_images_of_closed} that $\invrs{s}U \subseteq P \join Q$ is closed.
	Let $n \in \mathbb{N}$, $\alpha \in \set{+, -}$, and let $z \in \grade{n}{(P \join Q)}$.
	Because 
	\[
		\cofaces{}{-\alpha}s(z) = s(\cofaces{}{-\alpha}z),
	\]
	we have $\cofaces{}{-\alpha}z \cap \invrs{s}{U} = \varnothing$ if and only if $\cofaces{}{-\alpha}s(z) \cap U = \varnothing$, so $z \in \faces{n}{\alpha}\invrs{s}{U}$ if and only if $s(z) \in \faces{n+1}{\alpha}U$.	
	It follows that $\faces{n}{\alpha}\invrs{s}U = \invrs{s}\faces{n+1}{\alpha}U$.
	By Lemma \ref{lem:maximal_vs_faces} and the fact that inverse images of order-preserving maps preserve closures, unions, and intersections, we conclude that $\bound{n}{\alpha}\invrs{s}{U} = \invrs{s}\bound{n+1}{\alpha}U$.
\end{proof}

\begin{lem} \label{lem:boundaries_of_join}
	Let $U$, $V$ be non-empty closed subsets of oriented graded posets, $n \in \mathbb{N}$, $\alpha \in \set{+, -}$.
	Then $\bound{n}{\alpha}(U \join V)$ is equal to
	\[
	\begin{cases}
		\displaystyle \varnothing \join \bound{n}{+}V \cup \bigcup_{k=1}^n \bound{k-1}{+}U \join \bound{n-k}{(-)^{k}}V &
		\text{if $n$ is even, $\alpha = +$}, \\
		\displaystyle \bigcup_{k=1}^n \bound{k-1}{-}U \join \bound{n-k}{-(-)^{k}}V \cup \bound{n}{-}U \join \varnothing &
		\text{if $n$ is even, $\alpha = -$}, \\
		\displaystyle \varnothing \join \bound{n}{+}V \cup \bigcup_{k=1}^n \bound{k-1}{+}U \join \bound{n-k}{(-)^k}V \cup \bound{n}{+}U \join \varnothing &
		\text{if $n$ is odd, $\alpha = +$}, \\
		\displaystyle \bigcup_{k=1}^n \bound{k-1}{-}U \join \bound{n-k}{-(-)^k}V &
		\text{if $n$ is odd, $\alpha = -$}.
	\end{cases}
	\]	
\end{lem}
\begin{proof}
	Consider the injection $s\colon U \join V \to \sus{U} \gray \sus{V}$.
	By Lemma \ref{lem:joins_and_gray_of_suspension},
	\[
		\bound{n}{\alpha}(U \join V) = \bound{n}{\alpha}\invrs{s}(\sus{U} \gray \sus{V}) = \invrs{s}\bound{n+1}{\alpha}(\sus{U} \gray \sus{V}).
	\]
	By Corollary \ref{cor:boundaries_of_gray_product}, this is equal to
	\[
		\invrs{s}\left( \bigcup_{k=0}^{n+1} \bound{k}{\alpha}\sus{U} \gray \bound{n+1-k}{(-)^k\alpha}\sus{V} \right),
	\]
	which by Corollary \ref{cor:boundary_of_suspension} we can rewrite as
	\[
		\invrs{s}\left(
		\set{\bot^\alpha} \gray \sus{\bound{n}{\alpha}V} \cup
		\bigcup_{k=1}^n \sus{\bound{k-1}{\alpha}U} \gray \sus{\bound{n-k}{(-)^k\alpha}V} \cup
		\sus{\bound{n}{\alpha}U} \gray \set{\bot^{(-)^{n+1}\alpha}} \right).
	\]
	Since inverse images are compatible with unions, we have
	\[
		\invrs{s}\left( 
		\bigcup_{k=1}^n \sus{\bound{k-1}{\alpha}U} \gray \sus{\bound{n-k}{(-)^k\alpha}V}
		\right) = \bigcup_{k=1}^n \bound{k-1}{\alpha}U \join \bound{n-k}{(-)^k\alpha}V,
	\]
	while
	\begin{align*}
		\invrs{s}(\set{\bot^\alpha} \gray \sus{\bound{n}{\alpha}V}) & = 
		\begin{cases}
			\varnothing \join \bound{n}{+}V &
			\text{if $\alpha = +$}, \\
			\varnothing &
			\text{if $\alpha = -$},
		\end{cases} \\
		\invrs{s}(\sus{\bound{n}{\alpha}U} \gray \set{\bot^{(-)^{n+1}\alpha}}) & =
		\begin{cases}
			\bound{n}{\alpha}U \join \varnothing &
			\text{if $\alpha = (-)^{n+1}$}, \\
			\varnothing &
			\text{if $\alpha = -(-)^{n+1}$},
		\end{cases}
	\end{align*}
	and we conclude.
\end{proof}

\begin{lem} \label{lem:join_preserves_globular_and_round}
	Let $U$, $V$ be oriented graded posets.
	Then
	\begin{enumerate}
		\item if $U$ and $V$ are globular, then $U \join V$ is globular,
		\item if $U$ and $V$ are round, then $U \join V$ is round.
	\end{enumerate}
\end{lem}
\begin{proof}
	Let $s\colon U \join V \to \sus{U} \gray \sus{V}$ be the injection of Lemma 
	\ref{lem:joins_and_gray_of_suspension}, $n \in \mathbb{N}$, and $\alpha \in \set{+, -}$.
	Then
	\begin{equation} \label{eq:inverse_image_of_gray_sus}
		\bound{n}{\alpha}(U \join V) = \bound{n}{\alpha}\invrs{s}(\sus{U} \gray \sus{V}) = \invrs{s}\bound{n+1}{\alpha}(\sus{U} \gray \sus{V}).
	\end{equation}
	Suppose that $U$ and $V$ are globular.
	Then $\sus{U} \gray \sus{V}$ is globular by Lemma 
	\ref{lem:gray_preserves_globular} and Lemma \ref{lem:suspension_preserves_globular_and_round}.
	Globularity of $U \join V$ immediately follows from (\ref{eq:inverse_image_of_gray_sus}).
	
	Similarly, if $U$ and $V$ are round, then $\sus{U} \gray \sus{V}$ is round by Lemma
	\ref{lem:gray_preserves_round} and Lemma \ref{lem:suspension_preserves_globular_and_round}.
	Roundness of $U \join V$ then follows from (\ref{eq:inverse_image_of_gray_sus}) and the compatibility of inverse images with intersections.
\end{proof}

\begin{lem} \label{lem:inverse_image_of_gray_of_suspension}
	Let $U$, $V$ be molecules and consider $s\colon U \join V \to \sus{U} \gray \sus{V}$.
	For all $W \submol \sus{U} \gray \sus{V}$,
	\begin{enumerate}
		\item if $\invrs{s}W \neq \varnothing$, then $\invrs{s}W$ is a molecule,
		\item if $W' \submol W$ and $\invrs{s}W' \neq \varnothing$, then $\invrs{s}W' \submol \invrs{s}W$.
	\end{enumerate}
\end{lem}
\begin{proof}
	By Proposition \ref{prop:gray_product_of_molecules} and Proposition 
	\ref{prop:suspension_of_molecules}, $\sus{U} \gray \sus{V}$ is a molecule.
	We proceed by induction on submolecules $W \submol \sus{U} \gray \sus{V}$.
	If $\dim{W} = 0$, then $\invrs{s}W = \varnothing$, so we may assume $\dim{W} > 0$.

	Suppose that $W$ is an atom.
	Then $W = \clset{(x, y)}$ for some $x \in \sus{U}$ and $y \in \sus{V}$, and since $\invrs{s}W \neq \varnothing$, necessarily $(x, y) = s(z)$ for some $z \in (U \join V)$.
	It follows that $\invrs{s}W = \clset{z}$.
	If $z = \inj{x'}$ or $z = \inr{y'}$ for some $x' \in \grade{0}{U}$ and $y' \in \grade{0}{V}$, then $\clset{z} = \set{z}$, which is a 0\nbd dimensional atom.
	Otherwise, by Lemma \ref{lem:join_faces}, $z$ has at least one input and one output face, hence $\invrs{s}\bound{}{\alpha}W \neq \varnothing$ for all $\alpha \in \set{+, -}$.
	By the inductive hypothesis, $\invrs{s}\bound{}{\alpha}W$ is a molecule, and by Lemma 
	\ref{lem:joins_and_gray_of_suspension} it is equal to $\bound{}{\alpha}\invrs{s}W$.
	The same result, combined with the fact that inverse images preserve intersections, implies that $\invrs{s}W$ is round, which suffices to prove that $\invrs{s}W$ is an atom, isomorphic to $\bound{}{-}\invrs{s}W \celto \bound{}{+}\invrs{s}{W}$.
	Moreover, any proper submolecule of $W$ is a submolecule of its input or output boundary, so the inductive hypothesis applies.

	Next, suppose that $W$ splits into proper submolecules $W' \cup W''$ along the $k$\nbd boundary.
	If $\invrs{s}W' = \varnothing$, then $\invrs{s}W = \invrs{s}W''$, and dually if $\invrs{s}W'' = \varnothing$, then $\invrs{s}W = \invrs{s}W'$, and the claim follows from the inductive hypothesis.
	Suppose, instead, that $\invrs{s}W', \invrs{s}W'' \neq \varnothing$.
	We claim that, in this case, $k > 0$.
	Indeed, $W'$ must contain an element of the form $(x, \bot^+)$ or $(\bot^+, y)$, so it also contains $(\bot^+, \bot^+)$.
	But because $\cofaces{}{-}(\bot^+, \bot^+) = \varnothing$, necessarily $\bound{0}{+}W' = \set{(\bot^+, \bot^+)}$.
	Then $\bound{0}{+}W' = \bound{0}{-}W''$ would imply $W'' = \set{(\bot^+, \bot^+)}$, so $W = W' \cup W'' = W'$, contradicting the assumption that $W', W''$ are proper submolecules.
	Because $k > 0$, by Lemma \ref{lem:joins_and_gray_of_suspension} we have
	\[
		\invrs{s}W' \cap \invrs{s}W'' = \invrs{s}(W' \cap W'') = 
		\begin{cases}
			\invrs{s}\bound{k}{+}W' = \bound{k-1}{+}\invrs{s}W', \\
			\invrs{s}\bound{k}{-}W'' = \bound{k-1}{-}\invrs{s}W'',
		\end{cases}
	\]
	and by the inductive hypothesis both $\invrs{s}W'$ and $\invrs{s}W''$ are molecules.
	We conclude that $\invrs{s}W$ is a molecule, isomorphic to $\invrs{s}W' \cp{k-1} \invrs{s}W''$.
	This completes the induction.
\end{proof}

\begin{prop} \label{prop:join_of_molecules}
	Let $U$, $V$ be molecules.
	Then
	\begin{enumerate}
		\item $U \join V$ is a molecule,
		\item if $U' \submol U$ and $V' \submol V$, then $U' \join V' \submol U \join V$.
	\end{enumerate}
\end{prop}
\begin{proof}
	Follows immediately from Lemma \ref{lem:inverse_image_of_gray_of_suspension}, since $U \join V = \invrs{s}(\sus{U} \gray \sus{V})$ and, if $U' \submol U$ and $V' \submol V$, then $\sus{U'} \gray \sus{V'} \submol \sus{U} \gray \sus{V}$.
\end{proof}

\begin{cor} \label{cor:join_preserves_rdcpx}
	Let $P$, $Q$ be regular directed complexes.
	Then $P \join Q$ is a regular directed complex.
\end{cor}
\begin{proof}
	Let $z \in P \join Q$.
	If $z = \inj{x}$ for some $x \in P$, then $\clset{z}$ is the image of the inclusion $\clset{x} \join \varnothing \incl P \join Q$, and $\clset{x} \join \varnothing$ is isomorphic to the atom $\clset{x}$.
	Dually, if $z = \inr{y}$ for some $y \in P$, then $\clset{z}$ is isomorphic to the atom $\clset{y}$.
	Finally, if $z = x \join y$ for some $x \in P$ and $y \in Q$, then $\clset{z}$ is isomorphic to $\clset{x} \join \clset{y}$, which is an atom by Proposition \ref{prop:join_of_molecules}.
\end{proof}

\begin{cor} \label{cor:join_restricts_to_rdcpx}
	The monoidal structure $(\ogpos, \join, \varnothing)$ restricts to a monoidal structure on $\rdcpx$ and on $\rdcpxiso$.
\end{cor}

\begin{prop} \label{prop:join_monoidal_on_maps}
	There is a unique monoidal structure $(\rdcpxmap, \join, \varnothing)$ on $\rdcpxmap$ such that both
	\begin{enumerate}
		\item $(\rdcpx, \join, \varnothing) \incl (\rdcpxmap, \join, \varnothing)$ and
		\item $\fun{U}\colon (\rdcpxmap, \join, \varnothing) \to (\posclos, \join, \varnothing)$
	\end{enumerate}
	are strict monoidal functors.
\end{prop}
\begin{proof}
	Since the monoidal structure extends the one on $\rdcpx$, it is the join on objects, and the requirement that $\fun{U}$ be strict monoidal determines uniquely its action on maps, so it suffices to show that if $p\colon P \to P'$ and $q\colon Q \to Q'$ are maps of regular directed complexes, then the join of the underlying closed order-preserving maps lifts to a map $p\join q\colon P \join Q \to P' \join Q'$ of regular directed complexes.

	Consider the injections $s\colon P \join Q \to \sus{P} \gray \sus{Q}$ and $s'\colon P' \join Q' \to \sus{P'} \gray \sus{Q'}$.
	For all $z \in P \join Q$,
	\[
		s'((p \join q)(z)) = (\sus{p} \gray \sus{q})(s(z)),
	\]
	and $\sus{p} \gray \sus{q}$ is a map by Proposition \ref{prop:gray_product_of_maps} and Proposition \ref{prop:suspension_is_a_functor_on_maps}.
	Then, for all $n \in \mathbb{N}$ and $\alpha \in \set{+, -}$,
	\begin{align*}
		\bound{n}{\alpha}(p \join q)(z) & = \bound{n}{\alpha}\invrs{s'}\clos{(\sus{p} \gray \sus{q})(s(z))} = \invrs{s'}(\sus{p} \gray \sus{q})(\bound{n+1}{\alpha}s(z)) = \\
						& = \invrs{s'}(\sus{p} \gray \sus{q})(\clos{ s(\bound{n}{\alpha}z) }) = (p \join q)(\bound{n}{\alpha}z),
	\end{align*}
	and the fact that $\restr{(p \join q)}{\bound{n}{\alpha}z}$ is final onto its image follows by Lemma \ref{lem:joins_and_gray_of_suspension} from the fact that $\restr{(\sus{p} \gray \sus{q})}{\bound{n+1}{\alpha}s(z)}$ is final onto its image.
\end{proof}

\begin{prop} \label{prop:join_monoidal_on_comaps}
	There is a unique monoidal structure $(\rdcpxcomap, \join, \varnothing)$ on $\rdcpxcomap$ such that both
	\begin{enumerate}
		\item $(\rdcpxiso, \join, \varnothing) \incl (\rdcpxcomap, \join, \varnothing)$ and
		\item $\fun{U}\colon (\rdcpxcomap, \join, \varnothing) \to (\poscat, \join, \varnothing)$
	\end{enumerate}
	are strict monoidal functors.
\end{prop}
\begin{proof}
	As in the case of maps, the requirements fix the monoidal structure uniquely, so it suffices to prove that, if $c\colon P \to P'$ and $d\colon Q \to Q'$ are comaps, then the join of the underlying order-preserving maps lifts to a comap $c \join d\colon P \join Q \to P' \join Q'$.
	Consider the injections $s\colon P \join Q \to \sus{P} \gray \sus{Q}$ and $s'\colon P' \join Q' \to \sus{P'} \gray \sus{Q'}$, and let $z \in P' \join Q'$.
	Then
	\[
		\invrs{(c \join d)}\clset{z} = \invrs{s}\invrs{(\sus{c} \gray \sus{d})}\clset{s'(z)}.
	\]
	Now, $\invrs{(\sus{c} \gray \sus{d})}\clset{s'(z)}$ is a molecule by Proposition \ref{prop:gray_product_of_comaps} and Proposition \ref{prop:suspension_is_a_functor_on_comaps}, and in fact a submolecule of $\sus{P} \gray \sus{Q}$ by Proposition \ref{prop:comaps_inverse_image_preserves_molecules} and 
	Lemma \ref{lem:downset_is_submolecule}.
	It follows from Lemma \ref{lem:inverse_image_of_gray_of_suspension} that $\invrs{(c \join d)}\clset{z}$ is a molecule.
	
	Finally, let $n \in \mathbb{N}$ and $\alpha \in \set{+, -}$.
	By Lemma \ref{lem:joins_and_gray_of_suspension},
	\begin{align*}
		\invrs{(c \join d)}\bound{n}{\alpha}z & = \invrs{(c \join d)}\invrs{s'}\bound{n+1}{\alpha}s'(z) = \invrs{s}\invrs{(\sus{c} \gray \sus{d})}\bound{n+1}{\alpha}s'(z) = \\
						      & = \bound{n}{\alpha}\invrs{s}\invrs{(\sus{c} \gray \sus{d})}\clset{s'(z)} = \bound{n}{\alpha}\invrs{(c \join d)}\clset{z},
	\end{align*}
	and we conclude that $c \join d$ is a comap.
\end{proof}

\begin{dfn}[Join of augmented chain complexes] \index{chain complex!join} \index{join!of augmented chain complexes}
	Let $C$, $D$ be augmented chain complexes.
	The \emph{join of $C$ and $D$} is the augmented chain complex $C \join D$ with
	\[
		\grade{n}{(C \join D)} \eqdef \grade{n}{D} \oplus \left( \bigoplus_{k=0}^{n-1} \grade{k}{C} \otimes \grade{n-1-k}{D} \right) \oplus \grade{n}{C}
	\]
	for all $n \in \mathbb{N}$, together with the homomorphisms $\der, \eau$ defined as follows.
	For each $x \in \grade{n}{C}$ and $y \in \grade{m}{D}$, let
	\[
		\inj{x} \eqdef x \in \grade{n}{(C \join D)}, \quad 
		\inr{y} \eqdef y \in \grade{m}{(C \join D)}, \quad
		x \join y \eqdef x \otimes y \in \grade{n+m+1}{(C \join D)}.
	\]
	Then $\der\colon \grade{n}{(C \join D)} \to \grade{n-1}{(C \join D)}$ is defined, for each $n > 0$, by
	\begin{align*}
		\inj{x} & \mapsto \inj{\der (x)}, \\
		\inr{y} & \mapsto \inr{\der (y)}, \\
		x \join y & \mapsto
		\begin{cases}
			\eau (x) (\inr{y}) - \eau (y) (\inj{x}) 
			& \text{if $n = 1$, $x \in \grade{0}{C}$, $y \in \grade{0}{D}$,} \\
			\eau (x) (\inr{y}) - x \join \der(y) 
			& \text{if $n > 1$, $x \in \grade{0}{C}$, $y \in \grade{n-1}{C}$}, \\
			\der (x) \join y + (-)^{n} \eau(y) (\inj{x})
			& \text{if $n > 1$, $x \in \grade{n-1}{C}$, $y \in \grade{0}{D}$}, \\
				\der(x) \join y + (-)^{k+1} x \join \der(y)
			& \text{if $x \in \grade{k}{C}$, $y \in \grade{n-1-k}{D}$, $0 < k < n-1$},
		\end{cases}
	\end{align*}
	while $\eau\colon \grade{0}{(C \join D)} \to \mathbb{Z}$ is defined by
	\[
		\inj{x} \mapsto \eau(x), \quad \quad \inr{y} \mapsto \eau(y).
	\]
	The join extends to a monoidal structure on $\chaug$, whose unit is the chain complex $0$ which is $0$ in every degree, with the unique homomorphism $\eau\colon 0 \to \mathbb{Z}$.
\end{dfn}

\begin{comm}
	This definition is based on \cite[\S 6.5]{ara2020joint}.
	As detailed there, the join can be seen as arising from treating $C$ and $D$ as (non-augmented) chain complexes shifted by one degree, with $\eau$ as the final $\der$ homomorphism; then taking their tensor product as chain complexes; and, finally, shifting degrees back by one, and treating the last $\der$ homomorphism of the tensor product, which will have $\mathbb{Z} \otimes \mathbb{Z} \simeq \mathbb{Z}$ as codomain, as the $\eau$ homomorphism of an augmented chain complex.
\end{comm}

\begin{prop} \label{prop:join_compatible_with_chaug}
	Let $P$, $Q$ be oriented graded posets such that $\augm{P}$, $\augm{Q}$ are oriented thin.
	Then $\freeab{P} \join \freeab{Q}$ is naturally isomorphic to $\freeab{(P \join Q)}$.
\end{prop}
\begin{proof}
	The augmented chain complex $\freeab{(P \join Q)}$ is well-defined by 
	Corollary \ref{cor:join_of_oriented_thin}.
	Moreover, there are evident isomorphisms
	\begin{align*}
		\grade{n}{\varphi}\colon \grade{n}{\freeab{(P \join Q)}} & \to \grade{n}{(\freeab{P} \join \freeab{Q})}, \\
		z \in \grade{n}{(P \join Q)} & \mapsto
		\begin{cases}
			\inj{x} & \text{if $z = \inj{x}$, $x \in P$}, \\
			\inr{y} & \text{if $z = \inr{y}$, $y \in Q$}, \\
			x \join y & \text{if $z = x \join y$, $x \in P$, $y \in Q$}.
		\end{cases}
	\end{align*}
	The fact that these determine an isomorphism of augmented chain complexes is then a straightforward check using Lemma \ref{lem:join_faces}, as is naturality.
\end{proof}

\begin{cor} \label{cor:join_monoidal_functors_to_chaug}
	The functors 
	\[
		\freeab{-}\colon \otgpos \to \chaug, \quad \freeab{-}\colon \rdcpxmap \to \chaug, \quad \freeab{\pb{-}}\colon \opp{\rdcpxcomap} \to \chaug
	\]
	lift to strong monoidal functors
	\begin{align*}
		&\freeab{-}\colon (\otgpos, \otimes, 1) \to (\chaug, \join, 0), \\
		&\freeab{-}\colon (\rdcpxmap, \join, \varnothing) \to (\chaug, \join, 0), \\
		&\freeab{\pb{-}}\colon (\opp{\rdcpxcomap}, \join, \varnothing) \to (\chaug, \join, 0).
	\end{align*}
\end{cor}
\begin{proof}
	It suffices to check that naturality of the isomorphism of Proposition \ref{prop:join_compatible_with_chaug} extends to maps and comaps, and that this natural isomorphism, as well as the evident isomorphism $\freeab{\varnothing} \iso 0$, is compatible with unitors and associators.
	All of these are straightforward checks.
\end{proof}

\begin{exm}[Cones] \label{exm:cones} \index{cone}
	One way of constructing the \emph{cone} on a space $X$ is to take the join of $X$ and a point.
	As was the case with cylinders, for a regular directed complex $P$ this splits into two options, the \emph{left cone} $1 \join P$ and the \emph{right cone} $P \join 1$.
	In Example \ref{exm:nonsymmetry_join} we have seen both a left and a right cone in dimension 2, so here we consider a 3\nbd dimensional example.
	
	The left cone $1 \join \disk{2}{1}$ is a 3\nbd dimensional atom whose input and output boundaries are the oriented face posets of
\[
	\begin{tikzcd}[sep=small]
	& \bullet \\
	\\
	\bullet &&&& \bullet
	\arrow[""{name=0, anchor=center, inner sep=0}, curve={height=6pt}, from=3-1, to=3-5]
	\arrow[curve={height=-6pt}, from=3-1, to=1-2]
	\arrow[curve={height=-12pt}, from=1-2, to=3-5]
	\arrow[shorten <=6pt, Rightarrow, from=0, to=1-2]
\end{tikzcd}
	\quad \quad \text{and} \quad \quad
	\begin{tikzcd}[sep=small]
	& \bullet \\
	&& \bullet \\
	\bullet &&&& \bullet
	\arrow[""{name=0, anchor=center, inner sep=0}, curve={height=6pt}, from=3-1, to=3-5]
	\arrow[curve={height=-6pt}, from=3-1, to=1-2]
	\arrow[from=1-2, to=2-3]
	\arrow[""{name=1, anchor=center, inner sep=0}, from=3-1, to=2-3]
	\arrow[from=2-3, to=3-5]
	\arrow[""{name=2, anchor=center, inner sep=0}, curve={height=-24pt}, from=1-2, to=3-5]
	\arrow[shorten <=4pt, Rightarrow, from=0, to=2-3]
	\arrow[curve={height=-6pt}, shorten <=3pt, Rightarrow, from=1, to=1-2]
	\arrow[shorten >=4pt, Rightarrow, from=2-3, to=2]
\end{tikzcd}\]
respectively, while the input and output boundaries of the right cone $\disk{2}{1} \join 1$ are the oriented face posets of 
\[\begin{tikzcd}[sep=small]
	&&& \bullet \\
	&& \bullet \\
	\bullet &&&& \bullet
	\arrow[""{name=0, anchor=center, inner sep=0}, curve={height=6pt}, from=3-1, to=3-5]
	\arrow[curve={height=-6pt}, from=1-4, to=3-5]
	\arrow[""{name=1, anchor=center, inner sep=0}, curve={height=-24pt}, from=3-1, to=1-4]
	\arrow[from=2-3, to=1-4]
	\arrow[""{name=2, anchor=center, inner sep=0}, from=2-3, to=3-5]
	\arrow[from=3-1, to=2-3]
	\arrow[shorten <=4pt, Rightarrow, from=0, to=2-3]
	\arrow[shorten >=4pt, Rightarrow, from=2-3, to=1]
	\arrow[curve={height=6pt}, shorten <=3pt, Rightarrow, from=2, to=1-4]
\end{tikzcd}
	\quad \quad \text{and} \quad \quad 
\begin{tikzcd}[sep=small]
	&&& \bullet \\
	\\
	\bullet &&&& \bullet
	\arrow[""{name=0, anchor=center, inner sep=0}, curve={height=6pt}, from=3-1, to=3-5]
	\arrow[curve={height=-6pt}, from=1-4, to=3-5]
	\arrow[curve={height=-12pt}, from=3-1, to=1-4]
	\arrow[shorten <=6pt, Rightarrow, from=0, to=1-4]
\end{tikzcd}
\]
	respectively.
	Once again, the two are related by a duality reversing the direction of odd-dimensional cells, a consequence of the symmetry of $\disk{2}{1}$ and the point under this duality and of the fact that $P \join Q \simeq \opp{(\opp{Q} \join \opp{P})}$ by Proposition 
	\ref{prop:odd_duals_and_joins}.
\end{exm}


\section{Duals} \label{sec:duals}

\begin{guide}
	In this section, we define the duals of an oriented graded poset, a family of operators, indexed by sets $J$ of strictly positive natural numbers, which exchange the input and output faces of an $n$\nbd dimensional element whenever $n \in J$, while keeping the underlying poset fixed up to isomorphism.

	We prove that each dual of a molecule is a molecule (Proposition \ref{prop:dual_preserves_molecules}), hence each dual of a regular directed complex is a regular directed complex (Corollary 
	\ref{cor:dual_preserves_rdcpx}).
	Each dual determines an involutive endofunctor on $\ogpos$, $\rdcpxmap$, and $\rdcpxcomap$, compatibly on the overlaps, and in fact the family of all duals determines a representation of the abelian group of subsets of $\posnat$ with symmetric set difference as multiplication.

	We prove that duals of oriented graded posets are compatible with the analogous constructions on strict $\omega$\nbd categories and augmented chain complexes via the $\molecin{-}$ and $\freeab{-}$ functors.
	Finally, we consider the interaction of duals with Gray products, suspensions, and joins.
	These single out three duals in particular as more interesting than the rest: the \cemph{odd dual} $\opp{-}$, which reverses Gray products and joins; the \cemph{even dual} $\coo{-}$, which reverses Gray products; and the \cemph{total dual} $\optot{-}$, which preserves Gray products.
\end{guide}

\begin{dfn}[Duals of an oriented graded poset] \index{oriented graded poset!dual} \index{$\dual{J}{P}$} \index{dual!of an oriented graded poset}
	Let $P$ be an oriented graded poset, $J \subseteq \posnat$.
	The \emph{$J$\nbd dual of $P$} is the oriented graded poset $\dual{J}{P}$ whose
	\begin{itemize}
		\item underlying set is
			\[
				\set{ \dual{J}{x} \mid x \in P },
			\]
		\item partial order and orientation are defined by
			\[
				\faces{}{\alpha} \dual{J}{x} \eqdef
				\begin{cases}
					\set{ \dual{J}{y} \mid y \in \faces{}{-\alpha}x } &
					\text{if $\dim{x} \in J$}, \\
					\set{ \dual{J}{y} \mid y \in \faces{}{\alpha}x } &
					\text{if $\dim{x} \not\in J$}
				\end{cases}
			\]
		for all $x \in P$ and $\alpha \in \set{+,-}$.
	\end{itemize}
\end{dfn}

\begin{lem} \label{lem:duals_have_same_underlying poset}
	Let $P$ be an oriented graded poset, $J \subseteq \posnat$.
	Then
	\begin{align*}
		\dual{J}{}\colon \fun{U}P & \to \fun{U}\dual{J}{P}, \\
			x & \mapsto \dual{J}{x}
	\end{align*}
	is an isomorphism of posets.
\end{lem}
\begin{proof}
	The function is evidently a bijection at the level of underlying sets.
	Moreover, for all $x, y \in P$, we have $y \in \faces{}{}x$ if and only if $\dual{J}{y} \in \faces{}{}\dual{J}{x}$, so $\dual{J}{}$ preserves and reflects the covering relation.
	It follows that $\dual{J}{}$ is an isomorphism of posets.
\end{proof}

\begin{lem} \label{lem:dual_is_functor_on_ogpos}
	Let $f\colon P \to Q$ be a morphism of oriented graded posets, $J \subseteq \posnat$.
	Then
	\begin{align*}
		\dual{J}{f}\colon \dual{J}{P} & \to \dual{J}{Q}, \\
			\dual{J}{x} & \mapsto \dual{J}{f(x)}
	\end{align*}
	is a morphism of oriented graded posets.
	This assignment determines an endofunctor $\dual{J}{}$ on $\ogpos$, such that the diagram of functors
	\begin{equation} \label{eq:underlying_of_duals}
\begin{tikzcd}
	\ogpos && \ogpos \\
	& \posclos
	\arrow["{\dual{J}{}}", from=1-1, to=1-3]
	\arrow["{\fun{U}}"', from=1-1, to=2-2]
	\arrow["{\fun{U}}", from=1-3, to=2-2]
\end{tikzcd}\end{equation}
	commutes up to natural isomorphism.
\end{lem}
\begin{proof}
	Let $x \in P$ and $\alpha \in \set{+, -}$.
	Then $\dual{J}{f}$ determines a bijection between $\faces{}{\alpha}\dual{J}{x}$ and $\faces{}{\alpha}\dual{J}{f(x)}$ because $f$ determines a bijection between $\faces{}{\beta}x$ and $\faces{}{\beta}f(x)$ for $\beta \eqdef -\alpha$ if $\dim{x} \in J$ and $\beta \eqdef \alpha$ if $\dim{x} \not\in J$.
	This proves that $\dual{J}{f}$ is a morphism of oriented graded posets.
	Functoriality is evident, and commutativity of (\ref{eq:underlying_of_duals}) is a consequence of Lemma \ref{lem:duals_have_same_underlying poset}.
\end{proof}

\begin{prop} \label{prop:dual_of_oriented_thin}
	Let $P$ be an oriented graded poset such that $\augm{P}$ is oriented thin, $J \subseteq \posnat$.
	Then $\augm{(\dual{J}{P})}$ is oriented thin.
\end{prop}
\begin{proof}
	By Lemma \ref{lem:duals_have_same_underlying poset}, $\augm{(\dual{J}{P})}$ and $\augm{P}$ have isomorphic underlying posets.
	Given an interval $[x, y]$ in $\augm{P}$ with $\codim{x}{y} = 2$, which by oriented thinness is of the form
\[\begin{tikzcd}
	& {y} \\
	{z_1} && {z_2} \\
	& {x}
	\arrow["\alpha"', from=1-2, to=2-1]
	\arrow["\beta", from=1-2, to=2-3]
	\arrow["\gamma"', from=2-1, to=3-2]
	\arrow["-\alpha\beta\gamma", from=2-3, to=3-2]
\end{tikzcd}\]
	the corresponding interval $[x', y']$ in $\augm{(\dual{J}{P})}$ is of the form
\[\begin{tikzcd}
	& {y'} \\
	{z'_1} && {z'_2} \\
	& {x'}
	\arrow["(-)^i\alpha"', from=1-2, to=2-1]
	\arrow["(-)^i\beta", from=1-2, to=2-3]
	\arrow["(-)^j\gamma"', from=2-1, to=3-2]
	\arrow["-(-)^j\alpha\beta\gamma", from=2-3, to=3-2]
\end{tikzcd}\]
	with $i, j \in \set{0, 1}$ depending on $J$.
	For any choice of $i, j$, this satisfies the defining condition of oriented thinness.
\end{proof}

\begin{lem} \label{lem:dual_preserves_inclusions}
	Let $\imath\colon P \incl Q$ be an inclusion of oriented graded posets, $J \subseteq \posnat$.
	Then $\dual{J}{\imath}\colon \dual{J}{P} \to \dual{J}{Q}$ is an inclusion.
\end{lem}
\begin{proof}
	Immediate from the fact that $\fun{U}\dual{J}{}$ is naturally isomorphic to $\fun{U}$, and a morphism is an inclusion if and only if its underlying map is injective.
\end{proof}

\begin{rmk}
	Since $\dual{J}{}$ preserves inclusions, we can let it act on closed subsets $U \subseteq P$, producing closed subsets $\dual{J}{U} \subseteq \dual{J}{P}$.
\end{rmk}

\begin{lem} \label{lem:dual_preserves_reflected_limits_colimits}
	Let $\gamma$ be a limit or colimit cone in $\ogpos$ that is preserved and reflected by $\fun{U}\colon \ogpos \to \posclos$, and let $J \subseteq \posnat$.
	Then $\dual{J}{\gamma}$ is a limit or colimit cone in $\ogpos$.
\end{lem}
\begin{proof}
	Immediate from the fact that $\fun{U}\dual{J}{}$ is naturally isomorphic to $\fun{U}$.
\end{proof}

\begin{lem} \label{lem:dual_over_nothing}
	Let $P$ be an oriented graded poset.
	Then $P$ is naturally isomorphic to $\dual{\varnothing}{P}$.
\end{lem}
\begin{proof}
	Straightforward.
\end{proof}

\begin{lem} \label{lem:iterated_duals}
	Let $P$ be an oriented graded poset, $I, J \subseteq \posnat$, and let $K \eqdef I \setminus J + J \setminus I$.
	Then
	\begin{align*}
		\mu_{I,J}\colon \dual{I}{\dual{J}{P}} & \to \dual{K}{P}, \\
			\dual{I}{\dual{J}{x}} & \mapsto \dual{K}{x}
	\end{align*}
	is a natural isomorphism of oriented graded posets.
\end{lem}
\begin{proof}
	By Lemma \ref{lem:duals_have_same_underlying poset}, $\mu_{I, J}$ is an isomorphism of the underlying posets.
	Let $x \in P$, $\alpha \in \set{+, -}$.
	By a simple case distinction, we have that $\faces{}{\alpha}\dual{I}{\dual{J}{x}}$ is in bijection with $\faces{}{-\alpha}x$ if $\dim{x} \in I \setminus J$ or $\dim{x} \in J \setminus I$, and in bijection with $\faces{}{\alpha}x$ otherwise.
	The same is true of $\faces{}{\alpha}\dual{K}{x}$.
	Naturality is evident.
\end{proof}

\begin{rmk}
	We can see the combination of Lemma \ref{lem:dual_over_nothing} and Lemma \ref{lem:iterated_duals} as the statement that the family of endofunctors $\dual{J}{}$ is a representation in the monoidal category of endofunctors of $\ogpos$ of the group whose
	\begin{itemize}
		\item elements are subsets $J \subseteq \posnat$,
		\item multiplication is symmetric set difference, with $\varnothing$ as unit.
	\end{itemize}
	Since every element of this group has order 2, we have in particular that $\dual{J}{\dual{J}{}}$ is naturally isomorphic to $\bigid{\ogpos}$ for all $J \subseteq \posnat$.
\end{rmk}

\begin{lem} \label{lem:faces_of_dual}
	Let $U$ be a closed subset of an oriented graded poset, $n \in \mathbb{N}$, $\alpha \in \set{+, -}$, and $J \subseteq \posnat$.
	Then
	\[
		\faces{n}{\alpha}\dual{J}{U} = 
		\begin{cases}
			\set{\dual{J}{x} \mid x \in \faces{n}{-\alpha}U} &
				\text{if $n + 1 \in J$,} \\
			\set{\dual{J}{x} \mid x \in \faces{n}{\alpha}U} &
				\text{if $n + 1 \not \in J$.}
		\end{cases}
	\]
\end{lem}
\begin{proof}
	First of all, $x \in \grade{n}{U}$ if and only if $\dual{J}{x} \in \grade{n}{(\dual{J}{U})}$.
	Let $\beta \eqdef -\alpha$ if $n + 1 \in J$ and $\beta \eqdef \alpha$ otherwise.
	Then $\cofaces{}{-\beta}x \cap U = \varnothing$ if and only if, for all $y \in \grade{n+1}{U}$, $x \notin \faces{}{-\beta}y$, if and only if, for all $\dual{J}{y} \in \grade{n+1}{(\dual{J}{U})}$, $\dual{J}{x} \notin \faces{}{-\alpha}\dual{J}{y}$, if and only if $\cofaces{}{-\alpha}\dual{J}{x} \cap \dual{J}{U} = \varnothing$.
\end{proof}

\begin{cor} \label{cor:boundaries_of_dual}
	Let $U$ be a closed subset of an oriented graded poset, $n \in \mathbb{N}$, $\alpha \in \set{+, -}$, and $J \subseteq \posnat$.
	Then
	\[
		\bound{n}{\alpha}\dual{J}{U} =
		\begin{cases}
			\dual{J}{\bound{n}{-\alpha}U} &
			\text{if $n + 1 \in J$,} \\
			\dual{J}{\bound{n}{\alpha}U} &
			\text{if $n + 1 \not \in J$.}
		\end{cases}
	\]
\end{cor}

\begin{cor} \label{cor:dual_preserves_globular_and_round}
	Let $U$ be an oriented graded poset, $J \subseteq \posnat$.
	Then
	\begin{enumerate}
		\item if $U$ is globular, then $\dual{J}{U}$ is globular,
		\item if $U$ is round, then $\dual{J}{U}$ is round.
	\end{enumerate}
\end{cor}

\begin{prop} \label{prop:dual_preserves_molecules}
	Let $U$ be a molecule, $J \subseteq \posnat$.
	Then
	\begin{enumerate}
		\item $\dual{J}{U}$ is a molecule,
		\item if $U$ is isomorphic to $V \cp{k} W$, then $\dual{J}{U}$ is isomorphic to
			\[\begin{cases}
				\dual{J}{W} \cp{k} \dual{J}{V} & \text{if $k + 1 \in J$,} \\
				\dual{J}{V} \cp{k} \dual{J}{W} & \text{if $k + 1 \not\in J$,}
			\end{cases}\]
		\item if $U$ is isomorphic to $V \celto W$, then $\dual{J}{U}$ is isomorphic to
			\[\begin{cases}
				\dual{J}{W} \celto \dual{J}{V} & \text{if $\dim{U} \in J$,} \\
				\dual{J}{V} \celto \dual{J}{W} & \text{if $\dim{U} \not\in J$.}
			\end{cases}\]
	\end{enumerate}
\end{prop}
\begin{proof}
We proceed by induction on the construction of $U$.
If $U$ was produced by (\textit{Point}), then $U = 1$ and $\dual{J}{1}$ is isomorphic to $1$.

If $U$ was produced by (\textit{Paste}), then it is of the form $V \cp{k} W$ for some molecules $V, W$ and $k < \min \set{\dim{V}, \dim{W}}$.
Let 
\[\begin{tikzcd}
	\bound{k}{+}V & \bound{k}{-}W & W \\
	V && U
	\arrow["\sim", hook, from=1-1, to=1-2]
	\arrow[hook, from=1-2, to=1-3]
	\arrow[hook', from=1-1, to=2-1]
	\arrow[hook, from=2-1, to=2-3]
	\arrow[hook', from=1-3, to=2-3]
	\arrow["\lrcorner"{anchor=center, pos=0.125, rotate=180}, draw=none, from=2-3, to=1-1]
\end{tikzcd}\]
be the pushout diagram exhibiting $U$ as $V \cp{k} W$.
By Lemma \ref{lem:dual_preserves_reflected_limits_colimits} and Corollary \ref{cor:boundaries_of_dual}, $\dual{J}{}$ takes this to a pushout diagram
\[\begin{tikzcd}
	\bound{k}{\alpha}\dual{J}{V} & \bound{k}{-\alpha}\dual{J}{W} & \dual{J}{W} \\
	\dual{J}{V} && \dual{J}{U}
	\arrow["\sim", hook, from=1-1, to=1-2]
	\arrow[hook, from=1-2, to=1-3]
	\arrow[hook', from=1-1, to=2-1]
	\arrow[hook, from=2-1, to=2-3]
	\arrow[hook', from=1-3, to=2-3]
	\arrow["\lrcorner"{anchor=center, pos=0.125, rotate=180}, draw=none, from=2-3, to=1-1]
\end{tikzcd}\]
where $\dual{J}{V}$ and $\dual{J}{W}$ are molecules by the inductive hypothesis, and $\alpha \eqdef -$ if $k + 1 \in J$, $\alpha \eqdef +$ if $k + 1 \not \in J$.
Depending on the case, this diagram exhibits $\dual{J}{U}$ as $\dual{J}{W} \cp{k} \dual{J}{V}$ or $\dual{J}{V} \cp{k} \dual{J}{W}$.

Finally, if $U$ was produced by (\textit{Atom}), then it is of the form $V \celto W$ for some round molecules $V, W$ of the same dimension.
Let 
\[\begin{tikzcd}
	\bound{}{}V & \bound{}{}W & W \\
	V && \bound{}{}U
	\arrow["\sim", hook, from=1-1, to=1-2]
	\arrow[hook, from=1-2, to=1-3]
	\arrow[hook', from=1-1, to=2-1]
	\arrow[hook, from=2-1, to=2-3]
	\arrow[hook', from=1-3, to=2-3]
	\arrow["\lrcorner"{anchor=center, pos=0.125, rotate=180}, draw=none, from=2-3, to=1-1]
\end{tikzcd}\]
be the pushout diagram exhibiting $\bound{}{}U$ as a gluing of $V$ and $W$.
Then $\dual{J}{W}$ takes this to a pushout diagram
\[\begin{tikzcd}
	\bound{}{}\dual{J}{V} & \bound{}{}\dual{J}{W} & \dual{J}{W} \\
	\dual{J}{V} && \bound{}{}\dual{J}{U}
	\arrow["\sim", hook, from=1-1, to=1-2]
	\arrow[hook, from=1-2, to=1-3]
	\arrow[hook', from=1-1, to=2-1]
	\arrow[hook, from=2-1, to=2-3]
	\arrow[hook', from=1-3, to=2-3]
	\arrow["\lrcorner"{anchor=center, pos=0.125, rotate=180}, draw=none, from=2-3, to=1-1]
\end{tikzcd}\]
where $\dual{J}{V}$ and $\dual{J}{W}$ are both round molecules of the same dimension by the inductive hypothesis and Corollary \ref{cor:dual_preserves_globular_and_round}.
This determines an isomorphism between $\bound{}{}\dual{J}{U}$ and both 
\[
	\bound{}{}(\dual{J}{W} \celto \dual{J}{V}) \quad \text{and} \quad \bound{}{}(\dual{J}{V} \celto \dual{J}{W}).
\]
Let $\top$ be the greatest element of $U$.
For all $\alpha \in \set{+, -}$, $\faces{}{\alpha}\dual{J}{\top}$ is isomorphic to $\faces{}{-\alpha}\top$ if $\dim{U} = \dim{\top} \in J$, and to $\faces{}{\alpha}\top$ otherwise.
Depending on the case, we can extend the isomorphism of boundaries to an isomorphism between $\dual{J}{U}$ and $\dual{J}{W} \celto \dual{J}{V}$, or between $\dual{J}{U}$ and $\dual{J}{V} \celto \dual{J}{W}$.
\end{proof}

\begin{cor} \label{cor:dual_preserves_rdcpx}
	Let $P$ be a regular directed complex, $J \subseteq \posnat$.
	Then $\dual{J}{P}$ is a regular directed complex.
\end{cor}
\begin{proof} 
	Every element of $\dual{J}{P}$ is of the form $\dual{J}{x}$ for some $x \in P$, and by Proposition \ref{prop:dual_preserves_molecules} $\clset{\dual{J}{x}} = \dual{J}{\clset{x}}$ is an atom.
\end{proof}

\begin{cor} \label{cor:dual_of_layering}
	Let $U$ be a molecule, $k \geq -1$, $J \subseteq \posnat$, and let $(\order{i}{U})_{i=1}^m$ be a $k$\nbd layering of $U$.
	Then
	\begin{itemize}
		\item if $k + 1 \in J$, then $(\dual{J}{\order{m+1-i}{U}})_{i=1}^m$ is a $k$\nbd layering of $\dual{J}{U}$,
		\item if $k + 1 \not\in J$, then $(\dual{J}{\order{i}{U}})_{i=1}^m$ is a $k$\nbd layering of $\dual{J}{U}$.	
	\end{itemize}
\end{cor}

\begin{prop} \label{prop:dual_functorial_on_maps}
	Let $p\colon P \to Q$ be a map of regular directed complexes, $J \subseteq \posnat$.
	Then
	\begin{align*}
		\dual{J}{p}\colon \dual{J}{P} & \to \dual{J}{Q}, \\
			\dual{J}{x} & \mapsto \dual{J}{p(x)}
	\end{align*}
	is a map of regular directed complexes.
	This assignment determines an endofunctor $\dual{J}{}$ on $\rdcpxmap$, such that the diagram of functors
\begin{equation} \label{eq:dual_functors_commute}
	\begin{tikzcd}
	\ogpos && \rdcpx && \rdcpxmap \\
	\ogpos && \rdcpx && \rdcpxmap
	\arrow[hook', from=1-3, to=1-1]
	\arrow[hook, from=1-3, to=1-5]
	\arrow["{\dual{J}{}}", from=1-3, to=2-3]
	\arrow["{\dual{J}{}}", from=1-1, to=2-1]
	\arrow["{\dual{J}{}}", from=1-5, to=2-5]
	\arrow[hook, from=2-3, to=2-5]
	\arrow[hook', from=2-3, to=2-1]
\end{tikzcd}
\end{equation}
	commutes.
\end{prop}
\begin{proof}
	Let $x \in P$, $n \in \mathbb{N}$, and $\alpha \in \set{+, -}$.
	By Corollary \ref{cor:boundaries_of_dual},
	\begin{align*}
		\dual{J}{p}(\bound{n}{\alpha}\dual{J}{x})
		& = \dual{J}{p}(\dual{J}{\bound{n}{\beta}x})
		= \dual{J}{(p(\bound{n}{\beta}x))} = \\
		& = \dual{J}{\bound{n}{\beta}p(x)}
		= \bound{n}{\alpha}\dual{J}{p(x)},
	\end{align*}
	where $\beta = -\alpha$ if $n+1 \in J$ and $\beta = \alpha$ otherwise.
	Moreover, $\restr{p}{\bound{n}{\alpha}x}$ and $\restr{\dual{J}{p}}{\bound{n}{\alpha}\dual{J}{x}}$ have the same underlying map of posets, so one is final onto its image if and only if the other is.
	This proves that $\dual{J}{p}$ is a map of regular directed complexes.
	Functoriality of $\dual{J}{}$ and commutativity of (\ref{eq:dual_functors_commute}) are straightforward.
\end{proof}

\begin{prop} \label{prop:dual_functorial_on_comaps}
	Let $c\colon P \to Q$ be a comap of regular directed complexes, $J \subseteq \posnat$.
	Then
	\begin{align*}
		\dual{J}{c}\colon \dual{J}{P} & \to \dual{J}{Q}, \\
			\dual{J}{x} & \mapsto \dual{J}{c(x)}
	\end{align*}
	is a comap of regular directed complexes.
	This assignment determines an endofunctor $\dual{J}{}$ on $\rdcpxcomap$.
\end{prop}
\begin{proof}
	Let $y \in Q$, $n \in \mathbb{N}$, and $\alpha \in \set{+, -}$.
	We have
	\[
		\invrs{(\dual{J}{c})}\clset{\dual{J}{y}} = \dual{J}{(\invrs{c}\clset{y})},
	\]
	which by Proposition \ref{prop:dual_preserves_molecules} is a molecule because $\invrs{c}\clset{y}$ is a molecule.
	Moreover, by Corollary \ref{cor:boundaries_of_dual},
	\begin{align*}
		\invrs{(\dual{J}{c})}(\bound{n}{\alpha}\dual{J}{y})
		& = \invrs{(\dual{J}{c})}(\dual{J}{\bound{n}{\beta}y})
		= \dual{J}{(\invrs{c}(\bound{n}{\beta}y))} = \\
		& = \dual{J}{(\bound{n}{\beta}\invrs{c}\clset{y})}
		= \bound{n}{\alpha}\dual{J}{(\invrs{c}\clset{y})},
	\end{align*}
	where $\beta = -\alpha$ if $n+1 \in J$ and $\beta = \alpha$ otherwise.
	This proves that $\dual{J}{c}$ is a comap.
	Functoriality is straightforward.
\end{proof}

\begin{prop} \label{prop:suspension_and_dual}
	Let $P$ be an oriented graded poset, $J \subseteq \posnat$, and let $\sus{J} \eqdef \set{n+1 \mid n \in J}$.
	Then
	\begin{align*}
		\varphi\colon \sus{\dual{J}{P}} & \to \dual{\sus{J}}{\sus{P}}, \\
		x & \mapsto \begin{cases}
			\dual{\sus{J}}{\sus{x'}} &
			\text{if $x = \sus{\dual{J}{x'}}$, $x' \in P$,} \\
			\dual{\sus{J}}{\bot^\alpha} &
			\text{if $x = \bot^\alpha$, $\alpha \in \set{+, -}$}
		\end{cases}
	\end{align*}
	is a natural isomorphism of oriented graded posets.
\end{prop}
\begin{proof}
	The function is evidently an isomorphism of the underlying posets, so it suffices to show it is compatible with the orientations.
	Let $x \in \sus{\dual{J}{P}}$ and $\alpha \in \set{+, -}$.
	We may only consider the case $x = \sus{\dual{J}{x'}}$ for some $x' \in P$.
	If $\dim{x'} = 0$, then $\dim{x'} \not\in J$, and consequently $\dim{\sus{x'}} = \dim{x'} + 1 \not\in \sus{J}$.
	It follows that $\faces{}{\alpha}x = \set{\bot^\alpha}$ and $\faces{}{\alpha}\varphi(x) = \dual{\sus{J}}{\bot^\alpha}$.
	If $\dim{x'} > 0$, then by Lemma \ref{lem:faces_of_suspension} and the definition of $J$\nbd dual,
	\[
		\faces{}{\alpha}x = \set{\sus{\dual{J}{y'}} \mid \dual{J}{y'} \in \faces{}{\alpha}\dual{J}{x'}} = \set{\sus{\dual{J}{y'}} \mid y' \in \faces{}{\beta}x'}
	\]
	with $\beta = -\alpha$ if $\dim{x'} \in J$ and $\beta = \alpha$ otherwise.
	Meanwhile,
	\[
		\faces{}{\alpha}\varphi(x) = \set{\dual{\sus{J}}{y} \mid y \in \faces{}{\beta}\sus{x'}} = \set{\dual{\sus{J}}{\sus{y'}} \mid y' \in \faces{}{\beta}x'}
	\]
	with $\beta = -\alpha$ if $\dim{\sus{x'}} = \dim{x'} + 1 \in \sus{J}$ and $\beta = \alpha$ otherwise.
	By definition of $\sus{J}$, $\varphi$ determines a bijection between these sets.
\end{proof}

\begin{prop} \label{prop:augmentation_and_dual}
	Let $P$ be an oriented graded poset, $J \subseteq \posnat$, and let $\sus{J} \eqdef \set{n+1 \mid n \in J}$.
	Then
	\begin{align*}
		\varphi\colon \augm{(\dual{J}{P})} & \to \dual{\sus{J}}{(\augm{P})}, \\
		x & \mapsto \begin{cases}
			\dual{\sus{J}}{(\augm{x'})} &
			\text{if $x = \augm{(\dual{J}{x'})}$, $x' \in P$,} \\
			\dual{\sus{J}}{\bot} &
			\text{if $x = \bot$}
		\end{cases}
	\end{align*}
	is a natural isomorphism of oriented graded posets.
\end{prop}
\begin{proof}
	A straightforward variation on the proof of Proposition \ref{prop:suspension_and_dual}.
\end{proof}

\begin{dfn}[Odd dual of an oriented graded poset] \index{oriented graded poset!dual!odd} \index{$\opp{P}$} \index{dual!odd}
	Let $P$ be an oriented graded poset.
	The \emph{odd dual of $P$} is the oriented graded poset $\opp{P} \eqdef \dual{J}{P}$ for $J \eqdef \set{2n + 1 \mid n \in \mathbb{N}}$.
\end{dfn}

\begin{dfn}[Even dual of an oriented graded poset] \index{oriented graded poset!dual!even} \index{$\coo{P}$} \index{dual!even}
	Let $P$ be an oriented graded poset.
	The \emph{even dual of $P$} is the oriented graded poset $\coo{P} \eqdef \dual{J}{P}$ for $J \eqdef \set{2n + 2 \mid n \in \mathbb{N}}$.
\end{dfn}

\begin{dfn}[Total dual of an oriented graded poset] \index{oriented graded poset!dual!total} \index{$\optot{P}$} \index{dual!total}
	Let $P$ be an oriented graded poset.
	The \emph{total dual of $P$} is the oriented graded poset $\optot{P} \eqdef \dual{J}{P}$ for $J \eqdef \posnat$.
\end{dfn}

\begin{comm}
	We extend the notations $\opp{-}$, $\coo{-}$, and $\optot{-}$ to the elements of the odd, even, and total dual of an oriented graded poset.
\end{comm}

\begin{rmk}
	As special instances of Lemma \ref{lem:iterated_duals}, we have natural isomorphisms
	\[
		\opp{(\coo{P})} \iso \optot{P}, \quad \quad \coo{(\opp{P})} \iso \optot{P}.
	\]
	Moreover, as special instances of Proposition \ref{prop:suspension_and_dual} and of Proposition \ref{prop:augmentation_and_dual}, respectively, we have natural isomorphisms
	\[
		\sus{(\opp{P})} \iso \coo{(\sus{P})}, \quad \quad \augm{(\opp{P})} \iso \coo{(\augm{P})}.
	\]
\end{rmk}

\begin{lem} \label{lem:odd_and_even_dual_faces}
	Let $x \in P$, $\alpha \in \set{+, -}$.
	Then
	\begin{align*}
	\faces{}{\alpha}\opp{x} & = \set{\opp{y} \mid y \in \faces{}{(-)^{\dim{x}}\alpha}x}, \\
	\faces{}{\alpha}\coo{x} & = \set{\coo{y} \mid y \in \faces{}{(-)^{\dim{x}+1}\alpha}x}.
	\end{align*}
\end{lem}
\begin{proof}
	Follows from the definition and the observation that $\dim{x}$ is odd if and only if $(-)^{\dim{x}} = -$ and even if and only if $(-)^{\dim{x}+1} = -$.
\end{proof}

\begin{exm}[The duals of a 2-dimensional molecule]
	Up to isomorphism, every dual of a 2\nbd dimensional regular directed complex $P$ is either $P$, $\opp{P}$, $\coo{P}$, or $\optot{P}$.

	Let $U$ be the oriented face poset of our initial example (\ref{eq:example_shape}).
	Then $U$, $\opp{U}$, $\coo{U}$, and $\optot{U}$ are the oriented face posets of
\[\begin{tikzcd}[sep=small]
	\bullet && \bullet && \bullet \\
	& \bullet
	\arrow[curve={height=6pt}, from=1-1, to=2-2]
	\arrow[curve={height=6pt}, from=2-2, to=1-3]
	\arrow[""{name=0, anchor=center, inner sep=0}, curve={height=-18pt}, from=1-1, to=1-3]
	\arrow[from=1-3, to=1-5]
	\arrow[shorten <=3pt, shorten >=6pt, Rightarrow, from=2-2, to=0]
\end{tikzcd}\;, \quad \quad
\begin{tikzcd}[sep=small]
	\bullet && \bullet && \bullet \\
	&&& \bullet
	\arrow[curve={height=6pt}, from=2-4, to=1-5]
	\arrow[curve={height=6pt}, from=1-3, to=2-4]
	\arrow[""{name=0, anchor=center, inner sep=0}, curve={height=-18pt}, from=1-3, to=1-5]
	\arrow[from=1-1, to=1-3]
	\arrow[shorten <=3pt, shorten >=6pt, Rightarrow, from=2-4, to=0]
\end{tikzcd}\;, \]
\[
\begin{tikzcd}[sep=small]
	& \bullet \\
	\bullet && \bullet && \bullet
	\arrow[curve={height=-6pt}, from=2-1, to=1-2]
	\arrow[curve={height=-6pt}, from=1-2, to=2-3]
	\arrow[""{name=0, anchor=center, inner sep=0}, curve={height=18pt}, from=2-1, to=2-3]
	\arrow[from=2-3, to=2-5]
	\arrow[shorten <=6pt, shorten >=3pt, Rightarrow, from=0, to=1-2]
\end{tikzcd} \quad \text{and} \quad
	\begin{tikzcd}[sep=small]
	&&& \bullet \\
	\bullet && \bullet && \bullet
	\arrow[curve={height=-6pt}, from=1-4, to=2-5]
	\arrow[curve={height=-6pt}, from=2-3, to=1-4]
	\arrow[""{name=0, anchor=center, inner sep=0}, curve={height=18pt}, from=2-3, to=2-5]
	\arrow[from=2-1, to=2-3]
	\arrow[shorten <=6pt, shorten >=3pt, Rightarrow, from=0, to=1-4]
\end{tikzcd}\]
respectively.
\end{exm}

\begin{prop} \label{prop:duals_and_gray_products}
	Let $P$, $Q$ be oriented graded posets.
	Then
	\begin{align*}
		\opp{(P \gray Q)} &\to \opp{Q} \gray \opp{P}, 
		& \coo{(P \gray Q)} &\to \coo{Q} \gray \coo{P}, 
		& \optot{(P \gray Q)} &\to \optot{P} \gray \optot{Q}, \\
		\opp{(x, y)} &\mapsto (\opp{y}, \opp{x}), 
		& \coo{(x, y)} &\mapsto (\coo{y}, \coo{x}), 
		& \optot{(x, y)} &\mapsto (\optot{x}, \optot{y})
	\end{align*}
	are natural isomorphisms of oriented graded posets.
\end{prop}
\begin{proof}
	Clearly these are isomorphisms of the underlying posets, so it suffices to show that they are compatible with orientations.
	Let $x \in P$, $y \in Q$, and $\alpha \in \set{+, -}$.
	By Lemma \ref{lem:odd_and_even_dual_faces},
	\[
	\faces{}{\alpha}\opp{(x, y)} = \set{\opp{(x', y')} \mid (x', y') \in \faces{}{(-)^{\dim{x} + \dim{y}}\alpha}(x, y)},
	\]
	which by the definition of Gray product is equal to
	\[
		\set{\opp{(x', y)} \mid x' \in \faces{}{(-)^{\dim{x} + \dim{y}}\alpha}x} +
		\set{\opp{(x, y')} \mid y' \in \faces{}{(-)^{\dim{y}}\alpha}y}
	\]
	because $(-)^{\dim{x} + \dim{y}}(-)^{\dim{x}} = (-)^{\dim{y}}$.
	Using Lemma \ref{lem:odd_and_even_dual_faces} again, we see that $\opp{(x, y)} \mapsto (\opp{y}, \opp{x})$ maps this to
	\[
		\set{\opp{y}} \times \faces{}{(-)^{\dim{y}}\alpha}\opp{x} + \faces{}{\alpha}\opp{y} \times \set{\opp{x}} = \faces{}{\alpha}(\opp{y}, \opp{x}).
	\]
	The case of even duals is completely analogous.
	Finally, the case of total duals follows from the cases of even and odd duals by the natural isomorphism $\coo{(\opp{P})} \iso \optot{P}$.
\end{proof}

\begin{prop} \label{prop:odd_duals_and_joins}
	Let $P$, $Q$ be oriented graded posets.
	Then
	\begin{align*}
		\opp{(P \join Q)} &\to \opp{Q} \join \opp{P}, \\
		\opp{z} &\mapsto \begin{cases}
			\inr{\opp{x}} & \text{if $z = \inj{x}$, $x \in P$}, \\
			\inj{\opp{y}} & \text{if $z = \inr{y}$, $y \in Q$}, \\
			\opp{y} \join \opp{x} & \text{if $z = x \join y$, $x \in P$, $y \in Q$}
		\end{cases}
	\end{align*}
	is a natural isomorphism of oriented graded posets.
\end{prop}
\begin{proof}
	From Proposition \ref{prop:augmentation_and_dual}, we have a natural isomorphism 
	\[
		\augm{(\opp{(P \join Q)})} \iso \coo{(\augm{(P \join Q)})} = \coo{(\augm{P} \gray \augm{Q})},
	\]
	and by Proposition \ref{prop:duals_and_gray_products} we have a natural isomorphism
	\[
		\coo{(\augm{P} \gray \augm{Q})} \iso \coo{(\augm{Q})} \gray \coo{(\augm{P})}.
	\]
	Finally, from Proposition \ref{prop:augmentation_and_dual} again, we obtain a natural isomorphism
	\[
		\coo{(\augm{Q})} \gray \coo{(\augm{P})}	\iso \augm{(\opp{Q})} \gray \augm{(\opp{P})} = \augm{(\opp{Q} \join \opp{P})}.
	\]
	Transporting along the equivalence $\dimin{(-)}$, we get a natural isomorphism
	\[
		\opp{(P \join Q)} \iso \opp{Q} \join \opp{P}
	\]
	which, when made explicit, corresponds to the one in the statement.
\end{proof}

\begin{dfn}[Duals of a strict $\omega$-category] \index{strict $\omega$-category!dual} \index{$\dual{J}{X}$} \index{dual!of a strict $\omega$-category}
	Let $X$ be a strict $\omega$\nbd category, $J \subseteq \posnat$.
	The \emph{$J$-dual of $X$} is the strict $\omega$-category $\dual{J}{X}$ whose set of cells is
	\[
		\set{\dual{J}{t} \mid t \in X}
	\]
	with the boundary operators defined, for all $\dual{J}{t} \in \dual{J}{X}$, $n \in \mathbb{N}$, $\alpha \in \set{+, -}$, by
	\[
		\bound{n}{\alpha}\dual{J}{t} \eqdef
		\begin{cases}
			\dual{J}{\bound{n}{-\alpha}t} & \text{if $n + 1 \in J$}, \\
			\dual{J}{\bound{n}{\alpha}t} & \text{if $n + 1 \not\in J$}, 
		\end{cases}
	\]
	and the $k$\nbd composition operations defined, for all $k \in \mathbb{N}$ and $k$\nbd composable pairs $\dual{J}{t}, \dual{J}{u}$ in $\dual{J}{X}$, by
	\[
		\dual{J}{t} \cp{k} \dual{J}{u} \eqdef
		\begin{cases}
			\dual{J}{(u \cp{k} t)} & \text{if $k + 1 \in J$}, \\
			\dual{J}{(t \cp{k} u)} & \text{if $k + 1 \not\in J$}.
		\end{cases}
	\]
	The $J$\nbd dual extends to an endofunctor $\dual{J}\colon \omegacat \to \omegacat$.
\end{dfn}

\begin{prop} \label{prop:dual_omegacat_and_ogpos}
	Let $P$ be an oriented graded poset.
	Then $\dual{J}{\molecin{P}}$ is naturally isomorphic to $\molecin{\dual{J}{P}}$.	
\end{prop}
\begin{proof}
	Let $\varphi\colon \dual{J}{\molecin{P}} \to \molecin{\dual{J}{P}}$ be defined by
	\[
		\dual{J}{\isocl{f\colon U \to P}} \mapsto \isocl{\dual{J}{f}\colon \dual{J}{U} \to \dual{J}{P}}
	\]
	for each $\isocl{f\colon U \to P}$ in $\molecin{P}$.
	Then $\varphi$ is a morphism of reflexive $\omega$\nbd graphs by Corollary 
	\ref{cor:boundaries_of_dual} and a strict functor of strict $\omega$\nbd categories by Proposition \ref{prop:dual_preserves_molecules}.
	Moreover, $\varphi$ is both injective and surjective by virtue of $\dual{J}{}$ being an involution on $\ogpos$ up to natural isomorphism.
	We conclude that $\varphi$ is an isomorphism of strict $\omega$\nbd categories.
	Naturality is straightforward.
\end{proof}

\begin{dfn}[Duals of an augmented chain complex] \index{chain complex!dual} \index{dual!of an augmented chain complex}
	Let $C$ be an augmented chain complex, $J \subseteq \posnat$.
	The \emph{$J$-dual of $C$} is the augmented chain complex $\dual{J}{C}$ with
	\[
		\grade{n}{\dual{J}{C}} \eqdef \grade{n}{C}
	\]
	for all $n \in \mathbb{N}$, with $\der\colon \grade{n}{\dual{J}{C}} \to \grade{n-1}{\dual{J}{C}}$ defined by
	\[
		\begin{cases}
			-\der\colon \grade{n}{C} \to \grade{n-1}{C} 
			& \text{if $n \in J$,} \\
			\der\colon \grade{n}{C} \to \grade{n-1}{C}
			& \text{if $n \not\in J$}
		\end{cases}
	\]
	for each $n > 0$, and $\eau\colon \grade{0}{\dual{J}{C}} \to \mathbb{Z}$ equal to $\eau\colon \grade{0}{C} \to \mathbb{Z}$.
	The $J$\nbd dual extends to an endofunctor $\dual{J}{}\colon \chaug \to \chaug$.
\end{dfn}

\begin{prop} \label{prop:dual_chaug_and_ogtpos}
	Let $P$ be an oriented graded poset such that $\augm{P}$ is oriented thin.
	Then $\freeab{(\dual{J}{P})}$ is naturally isomorphic to $\dual{J}{(\freeab{P})}$.
\end{prop}
\begin{proof}
	By Proposition \ref{prop:dual_of_oriented_thin}, $\augm{(\dual{J}{P})}$ is oriented thin, so $\freeab{(\dual{J}{P})}$ is well-defined.
	For all $n \in \mathbb{N}$, the abelian groups $\grade{n}{\freeab{(\dual{J}{P})}}$ and $\grade{n}{\dual{J}{(\freeab{P})}}$ are, by construction, free on isomorphic sets.
	That the evident isomorphisms commute with $\der, \eau$ is a straightforward check of the definitions.
\end{proof}

\clearpage
\thispagestyle{empty}

%% file: acyclic.tex
\chapter{Acyclicity} \label{chap:acyclic}
\thispagestyle{firstpage}

\begin{guide}
	The recourse to some form of \cemph{acyclicity} condition is so common in the literature on higher-categorical diagrams, that one could think it is intrinsic to the theory.
	As should be apparent, since we made it this far without discussing it, this is not the case.
	There are several reasons why the imposition of acyclicity is widespread, some more valid than others.

	One reason is that, up to dimension 2, all molecules are acyclic in the strongest possible sense (Proposition \ref{prop:dim2_acyclic}).
	Moreover, as we will see in Chapter \ref{chap:special}, all globes, thetas, oriented simplices, oriented cubes, and positive opetopes are acyclic in this strong sense.
	So many higher category theorists may not have encountered a non-acyclic pasting diagram shape, or at least not \emph{realise} that they have.
	On the other hand, it is likely that most category theorists have, in some form, seen the pasting diagrams
\begin{equation} \label{eq:triangle_equation}
	\begin{tikzcd}[sep=small]
	&& {\smcat{C}} &&&& {\smcat{C}} \\
	{\smcat{C}} &&& {\smcat{D}} && {\smcat{C}} &&& {\smcat{D}} \\
	& {\smcat{D}} &&&&&& {\smcat{D}}
	\arrow["{\fun{R}}"', curve={height=6pt}, from=2-1, to=3-2]
	\arrow["{\fun{L}}", from=3-2, to=1-3]
	\arrow[""{name=0, anchor=center, inner sep=0}, "{\idd{\smcat{C}}}", curve={height=-12pt}, from=2-1, to=1-3]
	\arrow[""{name=1, anchor=center, inner sep=0}, "{\idd{\smcat{D}}}"', curve={height=12pt}, from=3-2, to=2-4]
	\arrow["{\fun{R}}", curve={height=-6pt}, from=1-3, to=2-4]
	\arrow["{\idd{\smcat{C}}}", curve={height=-6pt}, from=2-6, to=1-7]
	\arrow["{\fun{R}}"', from=1-7, to=3-8]
	\arrow[""{name=2, anchor=center, inner sep=0}, "{\fun{R}}", curve={height=-12pt}, from=1-7, to=2-9]
	\arrow[""{name=3, anchor=center, inner sep=0}, "{\fun{R}}"', curve={height=12pt}, from=2-6, to=3-8]
	\arrow["{\idd{\smcat{D}}}"', curve={height=6pt}, from=3-8, to=2-9]
	\arrow["\varepsilon"', curve={height=-6pt}, shorten >=7pt, Rightarrow, from=3-2, to=0]
	\arrow["\eta"', curve={height=6pt}, shorten <=7pt, Rightarrow, from=1, to=1-3]
	\arrow["{\idd{\fun{R}}}", curve={height=-6pt}, shorten <=7pt, Rightarrow, from=3, to=1-7]
	\arrow["{\idd{\fun{R}}}", curve={height=6pt}, shorten >=7pt, Rightarrow, from=3-8, to=2]
\end{tikzcd}\end{equation}
	as the two sides of one of the \emph{triangle equations} of adjunctions.
	Putting a 3\nbd cell between these two sides --- either because we are dealing with a lax adjunction or pseudoadjunction of 2\nbd categories, or simply because we want to express their equation as a commutative diagram in the style of Example \ref{exm:commutative_diagram} --- leads to a non-acyclic pasting diagram shape: the interior arrows ``loop around'' from the input to the output 1\nbd boundary.

	Pseudoadjunctions are likely to be the first higher algebraic structure that a category theorist encounters (possibly second to pseudomonoids), which shows that pasting diagrams with cycles are likely to arise ``in nature'' and not only in exotic contexts.

	Another reason for acyclicity, as discussed in the book introduction, is wishing for \emph{subsets} of cells of a diagram shape to form a strict $\omega$\nbd category, which is akin to wishing for linear subgraphs rather than paths in a directed graph to form a category, and we have already successfully overcome.
	Another reason, in ``analytic'' approaches that focus on \emph{characterisation} of well-formed pasting diagrams, is that an acyclicity condition is like a wide-bladed scythe, cutting away many undesirable examples in one swoop, at the price of losing many good ones as well.

	There is, on the other hand, one reason which is more valid than the others: some form of acyclicity seems to play a role in whether the strict $\omega$\nbd category presented by a diagram shape is or is not a \cemph{polygraph}, that is, ``freely generated'' in the appropriate sense.
	The interest in this condition comes, sometimes, from wanting a stronger form of pasting theorem, which will not only guarantee that a pasting diagram has a unique composite --- something that we have already achieved with molecules, without acyclicity --- but that the pasting diagram shapes \emph{classify valid expressions in the syntax of strict $\omega$\nbd categories, up to the equations of strict $\omega$\nbd categories}.

	Indeed, we will see that it is \emph{not} always the case that $\molecin{P}$ is a polygraph, when $P$ is a regular directed complex.
	In the light of Chapter \ref{chap:geometric}, one could see this as a demonstration that the algebra of strict $\omega$\nbd categories is even more star-crossed than we already knew: not only is it not sound for general homotopical algebra, it also fails to be complete for relations that exist in the pasting of topological cells.
	On the other hand, it is still useful for many reasons to investigate conditions under which $\molecin{P}$ is a polygraph, after all.

	Building on earlier work by Steiner, we identify a most general, very technical condition on $P$, called ``having \cemph{frame-acyclic} molecules'', which ensures that $\molecin{P}$ is a polygraph (Theorem \ref{thm:frame_acyclic_presents_polygraphs}).
	While this has to do with acyclicity, it is actually weak enough that \emph{all} regular directed complexes up to dimension 3 automatically satisfy it, even ones that exhibit all sorts of direct loops (Theorem \ref{thm:dim3_frame_acyclic}).
	
	Frame-acyclicity is tied to the existence of $k$\nbd layerings when $k$ is equal to the frame dimension, that is, the least possible dimension in which a layering can exist (Theorem 
	\ref{thm:frame_acyclic_has_frame_layerings}).
	As seen in Chapter \ref{chap:layerings}, this is connected to the problem of recognising rewritable submolecules, and we can use frame-acyclicity in order to provide a simplified, fully satisfactory criterion for recognising rewritable submolecule up to dimension 3 (Theorem \ref{thm:round_submolecule_dim3}).

	Finally, in Section \ref{sec:stronger_acyclic}, we revisit the stronger acyclicity conditions, and contextualise them as \emph{sufficient conditions} for having frame-acyclic molecules, which are in general easier to check.
	We study the stability of these properties under the constructions that we considered in Chapter \ref{chap:constructions}, as well as their implications on morphisms of oriented graded posets and regular directed complexes.
\end{guide}


\section{Frame-acyclic molecules} \label{sec:frame_acyclic}

\begin{guide}
	In this section, we define frame-acyclicity of a molecule, and prove that it is equivalent to the existence of layerings in the frame dimension \emph{across all submolecules} of the molecule (Theorem 
	\ref{thm:frame_acyclic_has_frame_layerings} and Corollary 
	\ref{cor:frame_acyclicity_equivalent_conditions}).
	We then prove some consequences of frame-acyclicity, for which we do not have an independent proof: there exists a path between any two elements in the oriented Hasse diagram of a frame-acyclic molecule (Proposition \ref{prop:hasseo_connected_frame_acyclic}), as well as a path in one of another family of graphs, the \cemph{extended flow graphs} (Proposition 
	\ref{prop:extflow_connected_frame_acyclic}).
\end{guide}

\begin{dfn}[Frame-acyclic molecule] \index{molecule!acyclic!frame-}
Let $U$ be a molecule.
We say that $U$ is \emph{frame-acyclic} if for all submolecules $V \submol U$, if $r \eqdef \frdim{V}$, then $\maxflow{r}{V}$ is acyclic.
\end{dfn}

\begin{lem} \label{lem:frame_acyclic_submolecule}
Let $U$ be a molecule, $V \submol U$.
If $U$ is frame-acyclic, then $V$ is frame-acyclic.
\end{lem}
\begin{proof}
Straightforward.
\end{proof}

\begin{lem} \label{lem:frdim_layerings_bijection_with_orderings}
Let $U$ be a molecule.
Suppose that for all submolecules $V \submol U$, if $r \eqdef \frdim{V}$, then $V$ admits an $r$\nbd layering.
Then for all $k \geq \frdim{U}$ the function $\lto{k}{U}\colon \layerings{k}{U} \incl \orderings{k}{U}$ is a bijection.
\end{lem}
\begin{proof}
Let $r \eqdef \frdim{U}$.
By assumption, there exists an $r$\nbd layering of $U$, so by Lemma \ref{lem:lto_bijection_above_layering} it suffices to show that $\lto{r}{U}$ is a bijection.

Given two $r$\nbd orderings $(\order{i}{x})_{i=1}^m$ and $(\order{i}{y})_{i=1}^m$, there exists a unique permutation $\sigma$ such that $\order{i}{x} = \order{\sigma(i)}{y}$ for all $i \in \set{1, \ldots, m}$.
Let $\fun{d}((\order{i}{x})_{i=1}^m, (\order{i}{y})_{i=1}^m)$ be the number of pairs $(j, j')$ such that $j < j'$ but $\sigma(j') < \sigma(j)$.
Under the assumption that $(\order{i}{x})_{i=1}^m$ is in the image of $\lto{r}{U}$, we will prove that $(\order{i}{y})_{i=1}^m$ is also in the image of $\lto{r}{U}$ by induction on $\fun{d}((\order{i}{x})_{i=1}^m, (\order{i}{y})_{i=1}^m)$.
Since the image of $\lto{r}{U}$ is not empty, this will suffice to prove that $\lto{r}{U}$ is surjective, hence bijective by Proposition \ref{prop:layerings_induce_orderings}.

If $\fun{d}((\order{i}{x})_{i=1}^m, (\order{i}{y})_{i=1}^m) = 0$, then $\order{i}{x} = \order{i}{y}$ for all $i \in \set{1, \ldots, m}$, and there is nothing left to prove.

Suppose $\fun{d}((\order{i}{x})_{i=1}^m, (\order{i}{y})_{i=1}^m) > 0$.
Then there exists $j < m$ such that $\sigma(j+1) < \sigma(j)$.
Suppose $(\order{i}{x})_{i=1}^m$ is the image of the $r$\nbd layering $(\order{i}{U})_{i=1}^m$.
Let $V \submol U$ be the image of $\order{j}{U} \cp{r} \order{j+1}{U}$ in $U$, and let 
\begin{equation*}
    z_1 \eqdef \order{j}{x} = \order{\sigma(j)}{y}, \quad z_2 \eqdef \order{j+1}{x} = \order{\sigma(j+1)}{y}.
\end{equation*}
Because $z_1$ comes before $z_2$ in one $r$\nbd ordering, but after in another, there can be no edge between them in $\maxflow{r}{U}$, so
\begin{equation*}
    \dim{ (\clset{{z_1}} \cap \clset{{z_2}}) } < r.
\end{equation*}
Since $z_1, z_2$ are the only maximal elements of dimension $> r$ in $V$, we deduce that $\ell \eqdef \frdim{V} < r$.
By assumption, there exists an $\ell$\nbd layering of $V$.
In particular, there exist molecules $\order{1}{V}, \order{2}{V}$ such that
\begin{enumerate}
    \item $z_i$ is in the image of $\order{i}{V}$ for all $i \in \set{1, 2}$, and
    \item $V$ is isomorphic to $\order{1}{V} \cp{\ell} \order{2}{V}$ or to $\order{2}{V} \cp{\ell} \order{1}{V}$.
\end{enumerate}
Without loss of generality suppose that $V$ is isomorphic to $\order{1}{V} \cp{\ell} \order{2}{V}$.
By Proposition \ref{prop:unitality_of_pasting} and Proposition \ref{prop:interchange_of_pasting}, letting
\begin{align*}
    \order{j}{\tilde{U}} & \eqdef \bound{r}{-}\order{1}{V} \cp{\ell} \order{2}{V}, \\
    \order{j+1}{\tilde{U}} & \eqdef \order{1}{V} \cp{\ell} \bound{r}{+} \order{2}{V},
\end{align*}
we have that $V$ is isomorphic to $\order{j}{\tilde{U}} \cp{r} \order{j+1}{\tilde{U}}$.
Letting $\order{i}{\tilde{U}} \eqdef \order{i}{U}$ for $i \notin \set{j, j+1}$, we have that $(\order{i}{\tilde{U}})_{i=1}^m$ is an $r$\nbd layering of $U$, and
\begin{equation*}
    \lto{r}{U}\colon (\order{i}{\tilde{U}})_{i=1}^m \mapsto (\order{i}{\tilde{x}})_{i=1}^m = (\order{1}{x}, \ldots, \order{j+1}{x}, \order{j}{x}, \ldots, \order{m}{x}).
\end{equation*}
Then $\fun{d}((\order{i}{\tilde{x}})_{i=1}^m, (\order{i}{y})_{i=1}^m) < \fun{d}((\order{i}{x})_{i=1}^m, (\order{i}{y})_{i=1}^m)$ and $(\order{i}{\tilde{x}})_{i=1}^m$ is in the image of $\lto{r}{U}$.
We conclude by the inductive hypothesis.
\end{proof}

\begin{thm} \label{thm:frame_acyclic_has_frame_layerings}
Let $U$ be a molecule, $r \eqdef \frdim{U}$.
If $U$ is frame-acyclic, then $U$ admits an $r$\nbd layering.
\end{thm}
\begin{proof}
By Lemma \ref{lem:frame_acyclic_submolecule}, we can proceed by induction on submolecules.
For all $x \in \grade{0}{U}$, we have $\frdim{\set{x}} = -1$, and $\set{x}$ admits the trivial $(-1)$\nbd layering, which proves the base case.

We construct a finite plane tree of submolecules $\order{j_1, \ldots, j_p}{U} \submol U$, as follows:
\begin{itemize}
    \item the root is $\order{}{U} \eqdef U$;
    \item if $\lydim{\order{j_1, \ldots, j_p}{U}} \leq r$, then we let $\lydim{\order{j_1, \ldots, j_p}{U}}$ be a leaf;
    \item if $k \eqdef \lydim{\order{j_1, \ldots, j_p}{U}} > r$, then we pick a $k$\nbd layering $(\order{i}{V})_{i=1}^q$ of $\order{j_1, \ldots, j_p}{U}$, which is possible by Theorem \ref{thm:molecules_admit_layerings}, and for each $i \in \set{1, \ldots, q}$, we let the image of $\order{i}{V}$ be a child $\order{j_1, \ldots, j_p, i}{U}$ of $\order{j_1, \ldots, j_p}{U}$.
\end{itemize}
By Lemma \ref{lem:lydim_layering_properties}, the layering dimension of the children of a node is strictly smaller than that of the node, so the procedure terminates.

Fix an $r$\nbd ordering $(\order{i}{x})_{i=1}^m$ of $U$; this is possible because $\maxflow{r}{U}$ is acyclic.
Let $V \eqdef \order{j_1, \ldots, j_p}{U}$ be a node of the tree.
We have
\begin{equation*}
    \bigcup_{j > r} \grade{j}{(\maxel{V})} = 
    \sum_{i = 1}^m \bigcup_{j > r} 
    \left(\grade{j}{(\maxel{V})} \cap \clset{{\order{i}{x}}}\right) \eqqcolon \sum_{i = 1}^m \order{i}{M};
\end{equation*}
the $\order{i}{M}$ form a partition because $\frdim{U} = r$, so every element of dimension $> r$ is in the closure of $\order{i}{x}$ for a unique $i \in \set{1, \ldots, m}$.
We claim that $V$ is isomorphic to
\begin{equation*}
    \order{1}{V} \cp{r} \ldots \cp{r} \order{m}{V}
\end{equation*}
for some molecules $(\order{i}{V})_{i=1}^m$ such that, for each $i \in \set{1, \ldots, m}$, identifying $\order{i}{V}$ with its image in $V$, we have
\begin{equation*}
    \bigcup_{j > r} \grade{j}{(\maxel{\order{i}{V}})} = 
    \order{i}{M}.
\end{equation*}
We will prove this by backward induction on the tree $\order{j_1, \ldots, j_p}{U}$.

Suppose $V$ is a leaf, so $\lydim{V} \leq r$.
Then $V$ admits an $r$\nbd layering.
For each $i \in \set{1, \ldots, m}$, fix a topological sort $(\order{i, j}{y})_{j=1}^{p_i}$ of the induced subgraph $\restr{\maxflow{r}{V}}{\order{i}{M}}$.
We claim that $((\order{i,j}{y})_{j=1}^{p_i})_{i=1}^m$ is an $r$\nbd ordering of $V$.

Suppose there is an edge from $x$ to $x'$ in $\maxflow{r}{V}$.
Then $x \in \order{i}{M}$, $x' \in \order{i'}{M}$ for a unique pair $i, i' \in \set{1, \ldots, m}$.
If $i = i'$, then $x = \order{i, j}{y}$ and $x' = \order{i, j'}{y}$ for some $j, j' \in \set{1, \ldots, p_i}$, and $j < j'$ because $(\order{i, j}{y})_{j=1}^{p_i}$ is a topological sort of $\restr{\maxflow{r}{V}}{\order{i}{M}}$.
If $i \neq i'$, then there exists 
\begin{equation*}
    z \in \faces{r}{+}x \cap \faces{r}{-}x' \subseteq \clset{{\order{i}{x}}} \cap \clset{{\order{i'}{x}}}.
\end{equation*}
Since $\bound{r}{\alpha}\order{i}{x}$ and $\bound{r}{\alpha}\order{i'}{x}$ is pure and $r$\nbd dimensional for all $\alpha \in \set{+, -}$, by Proposition \ref{prop:intersection_of_maximal_elements}
\begin{equation*}
    z \in (\faces{r}{+}\order{i}{x} \cap \faces{r}{-}\order{i'}{x}) \cup (\faces{r}{-}\order{i}{x} \cap \faces{r}{+}\order{i'}{x}),
\end{equation*}
and by Lemma \ref{lem:faces_intersection} $\faces{r}{-}\order{i}{x} \cap \clset{{x}} \subseteq \faces{r}{-}x$ which is disjoint from $\faces{r}{+}x$, so $z \in \faces{r}{+}\order{i}{x} \cap \faces{r}{-}\order{i'}{x}$.
It follows that there is an edge from $\order{i}{x}$ to $\order{i'}{x}$ in $\maxflow{r}{U}$, so $i < i'$ because $(\order{i}{x})_{i=1}^m$ is a topological sort of $\maxflow{r}{U}$.
This proves that $((\order{i,j}{y})_{j=1}^{p_i})_{i=1}^m$ is an $r$\nbd ordering of $V$.

Let $W \submol V$, $\ell \eqdef \frdim{W}$.
If $V \neq U$ or $W \neq U$, then $W$ admits an $\ell$\nbd layering by the inductive hypothesis on proper submolecules of $U$.
If $W = V = U$ then $\ell = r$ and $W$ admits an $\ell$\nbd layering by Theorem \ref{thm:molecules_admit_layerings}.
In either case, $V$ satisfies the conditions of Lemma \ref{lem:frdim_layerings_bijection_with_orderings}, and since $r \geq \lydim{V} \geq \frdim{V}$, every $r$\nbd ordering of $V$ comes from an $r$\nbd layering of $V$.

It follows that $((\order{i,j}{y})_{j=1}^{p_i})_{i=1}^m$ comes from an $r$\nbd layering $((\order{i,j}{W})_{j=1}^{p_i})_{i=1}^m$, and we can define
\begin{equation*}
    \order{i}{V} \eqdef \order{i,1}{W} \cp{r} \ldots \cp{r} \order{i, p_i}{W}
\end{equation*}
for each $i \in \set{1, \ldots, m}$, satisfying the desired condition.
    
Now, suppose that $V$ is not a leaf, so $k \eqdef \lydim{V} > r$, and $V$ has children $(\order{j}{W})_{j=1}^q$ forming a $k$\nbd layering of $V$.
By the inductive hypothesis, each of the $\order{j}{W}$ has a decomposition
\begin{equation*}
    \order{j, 1}{W} \cp{r} \ldots \cp{r} \order{j, m}{W}
\end{equation*}
such that the maximal elements of dimension $> r$ in the image of $\order{j, i}{W}$ are contained in $\clset{{\order{i}{x}}}$.
Then, for each $i \in \set{1, \ldots, m}$ and $j, j' \in \set{1, \ldots, q}$,
\begin{equation*}
    \order{j, i}{W} \cap \order{j'}{W} \subseteq \order{j', i}{W},
\end{equation*}
so $\order{i}{V} \eqdef \order{1, i}{W} \cp{k} \ldots \cp{k} \order{q, i}{W}$ is defined.
Using Proposition \ref{prop:interchange_of_pasting} repeatedly, we conclude that $V$ is isomorphic to $\order{1}{V} \cp{r} \ldots \cp{r} \order{m}{V}$.

This concludes the induction on the tree $\order{j_1, \ldots, j_p}{U}$.
In particular, for the root $\order{}{U} = U$, the decomposition $\order{1}{U} \cp{r} \ldots \cp{r} \order{m}{U}$ satisfies
\begin{equation*}
    \bigcup_{j > r} \grade{j}{(\maxel{\order{i}{U}})} = 
    \set{\order{i}{x}},
\end{equation*}
that is, $(\order{i}{U})_{i=1}^m$ is an $r$\nbd layering of $U$.
\end{proof}

\begin{cor} \label{cor:frame_acyclicity_equivalent_conditions}
Let $U$ be a molecule.
The following are equivalent:
\begin{enumerate}[label=(\alph*)]
    \item $U$ is frame-acyclic; \label{cond:frame_acyclic}
    \item for all $V \submol U$ and all $\frdim{V} \leq k < \dim{V}$, $V$ admits a $k$\nbd layering; \label{cond:frdim_layerings}
    \item for all $V \submol U$ and all $\frdim{V} \leq k < \dim{V}$, the sets $\layerings{k}{V}$ and $\orderings{k}{V}$ are non-empty and equinumerous. \label{cond:frdim_bijection}
\end{enumerate}
\end{cor}
\begin{proof}
The implication from \ref{cond:frame_acyclic} to \ref{cond:frdim_layerings} is a consequence of Theorem \ref{thm:frame_acyclic_has_frame_layerings} together with \ref{lem:frame_acyclic_submolecule} and Lemma \ref{lem:higher_layerings}.
The implication from \ref{cond:frdim_layerings} to \ref{cond:frdim_bijection} is Lemma \ref{lem:frdim_layerings_bijection_with_orderings}.
Finally, the implication from \ref{cond:frdim_bijection} to \ref{cond:frame_acyclic} follows from Proposition \ref{prop:if_layering_then_ordering}.
\end{proof}

\begin{dfn}[Extended flow graph] \index{oriented graded poset!flow graph!extended} \index{molecule!flow graph!extended} \index{$\extflow{k}{U}$}
Let $P$ be an oriented graded poset, $k \geq -1$.
The \emph{extended $k$\nbd flow graph of $P$} is the bipartite directed graph $\extflow{k}{P}$ whose
\begin{itemize}
	\item set of vertices is 
		\[
			P = \bigcup_{i \leq k} \grade{i}{P} + \bigcup_{i > k} \grade{i}{P},
		\]
	\item set of edges is $E_- + E_+$, where
		\begin{align*}
			E_- & \eqdef \set{ (y, x) \mid y \in \bigcup_{i \leq k}\grade{i}{P}, x \in \bigcup_{i > k} \grade{i}{P}, y \in \inter{\bound{k}{-}x} }, \\
			E_+ & \eqdef \set{ (y, x) \mid y \in \bigcup_{i > k}\grade{i}{P}, x \in \bigcup_{i \leq k}\grade{i}{P}, x \in \inter{\bound{k}{+}y} },
		\end{align*}
	with $s\colon (y, x) \mapsto y$ and $t\colon (y, x) \mapsto x$.
\end{itemize}
\end{dfn}

\begin{lem} \label{lem:path_in_flow_induces_path_in_extflow}
	Let $P$ be an oriented graded poset, $k \geq -1$, and suppose $x, y \in \bigcup_{i > k}\grade{i}{P}$.
	If there exists a path from $x$ to $y$ in $\flow{k}{P}$, then there exists a path from $x$ to $y$ in $\extflow{k}{P}$.
\end{lem}
\begin{proof}
	Consider a path $x = x_0 \to x_1 \to \ldots \to x_m \to y$ from $x$ to $y$ in $\flow{k}{P}$.
	By definition of the $k$\nbd flow graph, for all $i \in \set{1, \ldots, m}$, there exists $z_i \in \faces{k}{+}x_{i-1} \cap \faces{k}{-}x_i$.
	By definition of the extended $k$\nbd flow graph, there exist edges $x_{i-1} \to z_i$ and $z_i \to x_i$ in $\extflow{k}{P}$.
	Concatenating all the two-step paths $x_{i-1} \to z_i \to x_i$, we obtain a path from $x$ to $y$ in $\extflow{k}{P}$.
\end{proof}

\begin{lem} \label{lem:maximal_elements_flow_connected_frame_acyclic}
	Let $U$ be a frame-acyclic molecule, $x, y \in \maxel{U}$, and let $V \submol U$ be a minimal submolecule containing $\set{x, y}$.
	For $r \eqdef \frdim{V}$, either
\begin{itemize}
	\item there is a path from $x$ to $y$ in $\flow{r}{U}$, or
	\item there is a path from $y$ to $x$ in $\flow{r}{U}$.
\end{itemize}
\end{lem}

\begin{comm}
	By \emph{minimal} submolecule we mean that, if $W \submol V$ and $\set{x, y} \subseteq W$, then $W = V$.
	This is well-defined because $\set{x, y} \subseteq U$ and the submolecule relation is well-founded.
\end{comm}

\begin{proof}
By Theorem \ref{thm:frame_acyclic_has_frame_layerings}, $V$ admits an $r$\nbd layering $(\order{i}{V})_{i=1}^m$ inducing an $r$\nbd ordering $(\order{i}{z})_{i=1}^m$; we will identify each layer with its image in $V$.
By the minimality assumption on $V$, and since $x, y$ are maximal, necessarily $x = \order{1}{z}$ and $y = \order{m}{z}$, or the other way around.
If $m = 1$, then $x = y$ and we are done.

Otherwise, suppose without loss of generality that $x = \order{1}{z}$ and $y = \order{m}{z}$.
We claim that there is a path from $x$ to $y$ in $\maxflow{r}{V}$.
Suppose by way of contradiction that there is no such path.
Then it is possible to pick an $r$\nbd ordering $(\order{i}{w})_{i=1}^m$ of $V$ in which $\order{1}{z}$ and $\order{m}{z}$ are consecutive, say $\order{1}{z} = \order{j}{w}$ and $\order{m}{z} = \order{j+1}{w}$.
By Corollary \ref{cor:frame_acyclicity_equivalent_conditions}, this is induced by an $r$\nbd layering $(\order{i}{W})_{i=1}^m$.
Then $W \eqdef \order{j}{W} \cp{r} \order{j+1}{W}$ is a submolecule of $V$ containing $\set{x, y}$, so by minimality it is equal to $V$.
However, $\order{1}{z}$ and $\order{m}{z}$ are the only maximal elements of $W$ of dimension $> r$, yet $\clset{\order{1}{z}} \cap \clset{\order{m}{z}}$, which is equal to $\bound{r}{+}\order{1}{z} \cap \bound{r}{-}\order{m}{z}$ by Lemma \ref{lem:layering_intersections}, cannot contain any $r$\nbd dimensional elements.
We conclude that $\frdim{W} < r$, a contradiction.

Thus there is a path from $\order{1}{z}$ to $\order{m}{z}$ in $\maxflow{r}{V}$, which induces a path in $\flow{r}{U}$ by Lemma \ref{lem:flow_under_inclusion}.
\end{proof}

\begin{prop} \label{prop:extflow_connected_frame_acyclic}
Let $U$ be a frame-acyclic molecule, $x, y \in U$.
Then there exists $k \geq -1$ such that either
\begin{itemize}
	\item there is a path from $x$ to $y$ in $\extflow{k}{U}$, or
	\item there is a path from $y$ to $x$ in $\extflow{k}{U}$.
\end{itemize}
\end{prop}
\begin{proof}
Let $V \submol U$ be a minimal submolecule containing $\set{x, y}$, $r \eqdef \frdim{V}$.
Suppose first that $r = -1$; by Lemma \ref{lem:frame_dimension_atom}, $V$ is an atom, equal to $\clset{z}$ for some $z \in U$.
Let $n \eqdef \dim{z}$.
If $n = 0$, then $z = x = y$, and there is a trivial path between $x$ and $y$, so suppose $n > 0$.
If $x, y \in \bound{}{\alpha}z$ for some $\alpha \in \set{+, -}$, by Lemma \ref{lem:boundary_is_submolecule} $\bound{}{\alpha}z$ is a proper submolecule of $\clset{z}$ containing $\set{x, y}$, contradicting the minimality assumption.
It follows that, necessarily,
\[
	x \in \inter{\bound{}{\alpha}z} \cup \set{z}, \quad y \in \inter{\bound{}{-\alpha}z} \cup \set{z}
\]
for some $\alpha \in \set{+, -}$.
Then if $\alpha = -$ there is a path from $x$ to $y$, and if $\alpha = +$ there is a path from $y$ to $x$ in $\extflow{n-1}{U}$.

Suppose that $r \geq 0$.
By Theorem \ref{thm:frame_acyclic_has_frame_layerings}, $V$ admits an $r$\nbd layering $(\order{i}{V})_{i=1}^m$ inducing an $r$\nbd ordering $(\order{i}{z})_{i=1}^m$; we will identify each layer with its image in $V$.
By the minimality assumption, either $x \in \order{1}{V} \setminus \bound{r}{+}\order{1}{V}$ and $y \in \order{m}{V} \setminus \bound{r}{-}\order{m}{V}$, or the other way around.
Assume without loss of generality that it is the first way.
Then $x \in \clset{\order{1}{z}} \setminus \bound{r}{+}\order{1}{z}$ and $y \in \clset{\order{m}{z}} \setminus \bound{r}{-}\order{m}{z}$.

Necessarily, $V$ is also a minimal submolecule containing $\order{1}{z}$ and $\order{m}{z}$.
By Lemma \ref{lem:maximal_elements_flow_connected_frame_acyclic}, there is a path from $\order{1}{z}$ to $\order{m}{z}$ in $\flow{r}{U}$, which induces a path in $\extflow{r}{U}$ by Lemma 
\ref{lem:path_in_flow_induces_path_in_extflow}.
It then suffices to show that there is a path from $x$ to $\order{1}{z}$ and a path from $\order{m}{z}$ to $y$.

We claim that $\dim{\order{1}{z}} = r+1$.
Indeed, if $\dim{\order{1}{z}} > r+1$, then there exists $\alpha \in \set{ +, - }$ such that $x \in \bound{}{\alpha}\order{1}{z}$, which is a submolecule of $\bound{}{\alpha}\order{1}{V}$ by
Lemma \ref{lem:boundary_move}.
Then
\[
	W \eqdef \bound{}{-}\order{1}{V} \cp{r} (\order{2}{V} \cp{r} \ldots \cp{r} \order{m}{V})
\]
is a proper submolecule of $V$ containing $\set{x, y}$, contradicting minimality.
It follows that $\dim{\order{1}{z}} = r+1$, hence 
\[
	\clset{\order{1}{z}} \setminus \bound{r}{+}\order{1}{z} = \set{\order{1}{z}} \cup \inter{\bound{r}{-}\order{1}{z}},
\]
so either $x = \order{1}{z}$ or $x \in \inter{\bound{r}{-}\order{1}{z}}$, and there is an edge from $x$ to $\order{1}{z}$ in $\extflow{r}{U}$.
A dual argument produces a path from $\order{m}{z}$ to $y$, and we conclude.
\end{proof}

\begin{lem} \label{lem:path_to_faces_in_hasseo}
Let $P$ be a regular directed complex, $x \in P$, and $k \leq \dim{x}$.
Then
\begin{enumerate}
    \item for all $y \in \faces{k}{-}x$, there is a path from $y$ to $x$ in $\hasseo{P}$,
    \item for all $y \in \faces{k}{+}x$, there is a path from $x$ to $y$ in $\hasseo{P}$.
\end{enumerate}
\end{lem}
\begin{proof}
Let $\alpha \in \set{+, -}$.
We proceed by induction on $\dim{x} - k$.
If $\dim{x} = k$, then $\faces{k}{\alpha}x = \set{x}$, and there trivially exists a path from $x$ to $x$.
Suppose $\dim{x} > k$ and let $y \in \faces{k}{\alpha}x$.
By globularity, $y \in \faces{}{\alpha}(\bound{k+1}{\alpha}x)$.
By Corollary \ref{cor:atoms_are_round}, $\clset{{x}}$ is round, and by Lemma \ref{lem:boundaries_of_round_and_globular} so is $\bound{k+1}{\alpha}x$.
In particular $\bound{k+1}{\alpha}x$ is pure, so $y$ is not maximal in it, and there exists $z \in \cofaces{}{\alpha}y \cap \faces{k+1}{\alpha}x$.
Then
\begin{itemize}
    \item if $\alpha = -$, there is an edge from $y$ to $z$, and by the inductive hypothesis there is a path from $z$ to $x$ in $\hasseo{P}$;
    \item if $\alpha = +$, there is an edge from $z$ to $y$, and by the inductive hypothesis there is a path from $x$ to $z$ in $\hasseo{P}$. \qedhere
\end{itemize}
\end{proof}

\begin{lem} \label{lem:interior_elements_have_both_input_and_output_cofaces}
Let $U$ be a molecule, $x \in \inter{U}$.
Then either $x \in \maxel{U}$, or for all $\alpha \in \set{+, -}$ there exists $y \in \cofaces{}{\alpha}x \cap \inter{U}$.
\end{lem}
\begin{proof}
	Suppose $x \notin \maxel{U}$, let $\alpha \in \set{+, -}$, and assume by way of contradiction that $\cofaces{}{\alpha}x \cap U = \varnothing$.
	Then, for $k \eqdef \dim{x}$, we have $x \in \faces{k}{-\alpha}U \subseteq \bound{k}{}U$, which is contained in $\bound{}{}U$ by globularity of $U$, a contradiction.
	This implies that there exists $y \in \cofaces{}{\alpha}x \cap U$.
	If $y \in \bound{}{}U$, then $x \in \bound{}{}U$ since it is a closed subset, a contradiction.
	We conclude that $y \in \cofaces{}{\alpha}x \cap \inter{U}$.
\end{proof}

\begin{lem} \label{lem:paths_to_interiors_in_hasseo}
	Let $P$ be a regular directed complex, $x \in P$, and $k \leq \dim{x}$.
	Then
	\begin{enumerate}
		\item for all $y \in \inter{\bound{k}{-}x}$, there is a path from $y$ to $x$ in $\hasseo{P}$,
		\item for all $y \in \inter{\bound{k}{+}x}$, there is a path from $x$ to $y$ in $\hasseo{P}$.
	\end{enumerate}
\end{lem}
\begin{proof}
	Let $y \in \inter{\bound{k}{-}x}$; we proceed by induction on $k - \dim{y}$.
	Because $\clset{x}$ is an atom, $\bound{k}{-}x$ is a round molecule, hence pure and $k$\nbd dimensional, and its set of maximal elements is $\faces{k}{-}x$.
	If $\dim{y} = k$, then $y \in \faces{k}{-}x$, and by Lemma \ref{lem:path_to_faces_in_hasseo} there is a path from $y$ to $x$ in $\hasseo{P}$. 
	Otherwise, $y$ is not maximal in $\bound{k}{-}x$, so by Lemma \ref{lem:interior_elements_have_both_input_and_output_cofaces} there exists $y' \in \cofaces{}{+}y \cap \inter{\bound{k}{-}x}$.
	This means that there is an edge from $y$ to $y'$ in $\hasseo{P}$, and since $k - \dim{y'} = k - \dim{y} - 1$, the inductive hypothesis gives us a path from $y'$ to $x$.
	Concatenating the two produces a path from $y$ to $x$.
	The case $y \in \inter{\bound{k}{+}x}$ is dual.
\end{proof}

\begin{lem} \label{lem:path_in_extflow_induces_path_in_hasseo}
	Let $P$ be a regular directed complex, $k \geq -1$, and $x, y \in P$.
	If there exists a path from $x$ to $y$ in $\extflow{k}{P}$, then there exists a path from $x$ to $y$ in $\hasseo{P}$.
\end{lem}
\begin{proof}
	It suffices to consider the case where there is an edge $x \to y$ in $\extflow{k}{P}$.
	Then either $\dim{y} > k$ and $x \in \inter{\bound{k}{-}y}$ or $\dim{x} > k$ and $y \in \inter{\bound{k}{+}x}$.
	In both cases, by Lemma \ref{lem:paths_to_interiors_in_hasseo}, there is a path from $x$ to $y$ in $\hasseo{P}$.
\end{proof}

\begin{prop} \label{prop:hasseo_connected_frame_acyclic}
Let $U$ be a frame-acyclic molecule, $x, y \in U$.
Either
\begin{itemize}
	\item there is a path from $x$ to $y$ in $\hasseo{U}$, or
	\item there is a path from $y$ to $x$ in $\hasseo{U}$.
\end{itemize}
\end{prop}
\begin{proof}
	Follows from Proposition \ref{prop:extflow_connected_frame_acyclic} and Lemma \ref{lem:path_in_extflow_induces_path_in_hasseo}.
\end{proof}

\begin{comm}
	We do not know, at the moment, whether Proposition \ref{prop:hasseo_connected_frame_acyclic} can be generalised from frame-acyclic molecules to all molecules.
\end{comm}


\section{Presenting polygraphs} \label{sec:polygraphs}

\begin{guide}
	We say that an oriented graded poset $P$ \emph{has frame-acyclic molecules} if all molecules over $P$ are frame-acyclic.
	In this section, after recalling the definition of polygraph, we prove that if $P$ has frame-acyclic molecules, then $\molecin{P}$ is a polygraph, freely generated by the atoms over $P$ (Theorem 
	\ref{thm:frame_acyclic_presents_polygraphs}).
	We then exhibit a 4\nbd dimensional counterexample to this property, which also serves as a counterexample to the extension of other statements beyond the frame-acyclic case.
\end{guide}

\begin{dfn}[Cellular extension of a strict $\omega$-category] \index{strict $\omega$-category!cellular extension}
	Let $X$ be a strict $\omega$\nbd category.
	A \emph{cellular extension of $X$} is a strict $\omega$\nbd category $X_\gener{S}$ together with a pushout diagram
\[
	\begin{tikzcd}
		{\coprod_{e \in \gener{S}} \molecin{\bound{}{}U_e}} &&& {\coprod_{e \in \gener{S}} \molecin{U_e}} \\
		X &&& {X_\gener{S}}
	\arrow["{(\bound{}{}e)_{e \in \gener{S}}}", from=1-1, to=2-1]
	\arrow["{(e)_{e \in \gener{S}}}", from=1-4, to=2-4]
	\arrow[hook, from=2-1, to=2-4]
	\arrow["{\coprod_{e \in \gener{S}}\molecin{\imath_e}}", hook, from=1-1, to=1-4]
	\arrow["\lrcorner"{anchor=center, pos=0.125, rotate=180}, draw=none, from=2-4, to=1-1]
\end{tikzcd}\]
in $\omegacat$, where, for each $e \in \gener{S}$, $U_e$ is an atom and $\imath_e \colon \bound{}{}U_e \incl U_e$ is the inclusion of its boundary.
\end{dfn}

\begin{comm}
	This is a non-standard definition of cellular extension; the standard definition essentially amounts to requiring that each $U_e$ be a \emph{globe}, see Section \ref{sec:globes_and_thetas}.
	However, as a consequence of Lemma 
	\ref{lem:atom_is_cellular_extension_of_boundary}, combined with the pasting law for pushout squares and the fact that coproducts of pushout squares are pushout squares, we can turn every cellular extension in this generalised sense into a cellular extension in the restricted sense, without changing $X_\gener{S}$ nor, up to a bijection, the set $\gener{S}$.
\end{comm}

\begin{rmk}
	The functor $X \incl X_\gener{S}$ in a cellular extension is always injective, as shown in \cite[Section 4]{makkai2005word}.
\end{rmk}

\begin{dfn}[Polygraph] \index{strict $\omega$-category!polygraph|see {polygraph}} \index{polygraph} \index{$\cwcom{X}{S}$}
	A \emph{polygraph} is a strict $\omega$-category $X$ together with, for each $n \in \mathbb{N}$, a pushout diagram
\[
	\begin{tikzcd}
		{\coprod_{e \in \grade{n}{\gener{S}}} \molecin{\bound{}{}U_e}} &&& {\coprod_{e \in \grade{n}{\gener{S}}} \molecin{U_e}} \\
		\skel{n-1}{X} &&& \skel{n}{X}
	\arrow["{(\bound{}{}e)_{e \in \grade{n}{\gener{S}}}}", from=1-1, to=2-1]
	\arrow["{(e)_{e \in \grade{n}{\gener{S}}}}", from=1-4, to=2-4]
	\arrow[hook, from=2-1, to=2-4]
	\arrow["{\coprod_{e \in \grade{n}{\gener{S}}}\molecin{\imath_e}}", hook, from=1-1, to=1-4]
	\arrow["\lrcorner"{anchor=center, pos=0.125, rotate=180}, draw=none, from=2-4, to=1-1]
\end{tikzcd}\]
in $\omegacat$, exhibiting $\skel{n}{X}$ as a cellular extension of $\skel{n-1}{X}$, such that $U_e$ is an $n$\nbd dimensional atom for all $e \in \grade{n}{\gener{S}}$.
The set
\[
	\gener{S} \eqdef \sum_{n \in \mathbb{N}} \set{e\isocl{\idd{U_e}} \mid e \in \grade{n}{\gener{S}}}
\]
is called the set of \emph{generating cells} of the polygraph.
We write $\cwcom{X}{S}$ for a polygraph $X$ with set $\gener{S}$ of generating cells.
\end{dfn}

\begin{comm}
	There is a notorious terminology divide between the ``French school'', where the term \emph{polygraph} originated after \cite{burroni1993higher}, and the ``Australian school'', where polygraphs are more commonly known as \emph{computads} after \cite{street1976limits}.
	While \emph{computad} has temporal precedence, \emph{polygraph} has arguably become more established after recent activity culminating in the book \cite{ara2023polygraphs}.
\end{comm}

\begin{lem} \label{lem:generating_cells_are_generating}
	Let $\cwcom{X}{S}$ be a polygraph.
	Then $\gener{S}$ is a basis for $X$.
\end{lem}
\begin{proof}
	The fact that $\gener{S}$ is a generating set and its minimality are consequences of \cite[Proposition 15.1.8 and Lemma 16.6.2]{ara2023polygraphs}, respectively.
\end{proof}

\begin{dfn}[Map of polygraphs] \index{polygraph!map} \index{map!of polygraphs}
	Let $\cwcom{X}{S}$ and $\cwcom{Y}{T}$ be polygraphs.
	A \emph{map} $f\colon \cwcom{X}{S} \to \cwcom{Y}{T}$ is a pair of
	\begin{enumerate}
		\item a strict functor $f\colon X \to Y$,
		\item a function $\fpos{f}\colon \gener{S} \to \gener{T}$,
	\end{enumerate}
	such that, for all $t \in \gener{S}$, $f(t) = \fpos{f}(t)$.
\end{dfn}

\begin{dfn}[The category $\polmap$] \index{$\polmap$}
	We let $\polmap$ denote the category of polygraphs and maps of polygraphs.
\end{dfn}

\begin{comm}
	A more restrictive notion of morphism of polygraphs, which is the one considered in 
	\cite[Section 16.2]{ara2023polygraphs}, requires $\fpos{f}$ to be dimension-preserving.
\end{comm}

\begin{prop} \label{prop:forgetful_from_polmap_is_pseudomonic}
	The forgetful functor $\fun{U}\colon \polmap \to \omegacat$ is pseudomonic, that is,
	\begin{enumerate}
		\item it is faithful,
		\item it is full on isomorphisms,
		\item it reflects isomorphisms.
	\end{enumerate}
\end{prop}
\begin{proof}
	Faithfulness is a consequence of 
	Lemma \ref{lem:functors_equal_on_generating_set} combined with Lemma 
	\ref{lem:generating_cells_are_generating}.
	Fullness and reflection of isomorphisms are \cite[Proposition 16.6.3]{ara2023polygraphs}.
\end{proof}

\begin{comm}
	As a consequence of Proposition \ref{prop:forgetful_from_polmap_is_pseudomonic}, a strict $\omega$\nbd category admits at most one structure of polygraph up to unique isomorphism in $\polmap$; that is, a structure of polygraph is a ``property-like structure''.
	Thus one is justified in saying that a strict $\omega$\nbd category \emph{is} or \emph{is not} a polygraph.
\end{comm}

\begin{dfn}[Oriented graded poset with frame-acyclic molecules] \index{oriented graded poset!with frame-acyclic molecules}
	Let $P$ be an oriented graded poset.
	We say that $P$ \emph{has frame-acyclic molecules} if, for all molecules $U$, if there exists a morphism $f\colon U \to P$, then $U$ is frame-acyclic.
\end{dfn}

\begin{lem} \label{lem:frame_acyclic_cellular_extension}
	Let $P$ be an oriented graded poset, $n \in \mathbb{N}$, and let $\grade{n}{\gener{S}}$ be a set containing one pasting diagram
	\[
		e \equiv \molecin{e}\colon \molecin{U_e} \to \skel{n}{\molecin{P}}
	\]
	for each $\isocl{e\colon U_e \to P}$ in $\atomin{P}$ such that $\dim{U_e} = n$.
	If $\skel{n}{P}$ has frame-acyclic molecules, then
\[
	\begin{tikzcd}
		{\coprod_{e \in \grade{n}{\gener{S}}} \molecin{\bound{}{}U_e}} &&& {\coprod_{e \in \grade{n}{\gener{S}}} \molecin{U_e}} \\
		\skel{n-1}{\molecin{P}} &&& \skel{n}{\molecin{P}}
		\arrow["{(\bound{}{}e)_{e \in \grade{n}{\gener{S}}}}", from=1-1, to=2-1]
		\arrow["{(e)_{e \in \grade{n}{\gener{S}}}}", from=1-4, to=2-4]
		\arrow[hook, from=2-1, to=2-4]
		\arrow["{\coprod_{e \in \grade{n}{\gener{S}}}\molecin{\imath_e}}", hook, from=1-1, to=1-4]
		\arrow["\lrcorner"{anchor=center, pos=0.125, rotate=180}, draw=none, from=2-4, to=1-1]
\end{tikzcd}\]
	is a pushout diagram in $\omegacat$, exhibiting $\skel{n}{\molecin{P}}$ as a cellular extension of $\skel{n-1}{\molecin{P}}$.
\end{lem}
\begin{proof}
	Let $X$ be a strict $\omega$\nbd category and let
\[
	\begin{tikzcd}
		{\coprod_{e \in \grade{n}{\gener{S}}} \molecin{\bound{}{}U_e}} &&& {\coprod_{e \in \grade{n}{\gener{S}}} \molecin{U_e}} \\
		\skel{n-1}{\molecin{P}} &&& X
		\arrow["{(\bound{}{}e)_{e \in \grade{n}{\gener{S}}}}", from=1-1, to=2-1]
		\arrow["\ell", from=1-4, to=2-4]
		\arrow["h", from=2-1, to=2-4]
		\arrow["{\coprod_{e \in \grade{n}{\gener{S}}}\molecin{\imath_e}}", hook, from=1-1, to=1-4]
\end{tikzcd}\]
	be a commutative diagram of strict functors.
	We define $\overbar{h}\colon \skel{n}{\molecin{P}} \to X$ as follows on cells $\isocl{f\colon U \to P}$ in $\skel{n}{\molecin{P}}$.
	If $\dim{U} < n$, then we let 
	\[	\overbar{h}\isocl{f} \eqdef h\isocl{f}. \]
	Suppose $\dim{U} = n$; we proceed by induction on $\lydim{U}$.
	If $\lydim{U} = -1$, then by Lemma \ref{lem:layering_dimension_atom} $U$ is an atom, so there exists a unique $\molecin{e} \in \grade{n}{\gener{S}}$ such that $\isocl{f} = \isocl{e}$, and we let
	\[
		\overbar{h}\isocl{f} \eqdef \ell \isocl{\idd{U_e}}.
	\]
	If $\lydim{U} = k \geq 0$, then by Theorem \ref{thm:molecules_admit_layerings} $U$ admits a $k$\nbd layering $(\order{i}{U})_{i=1}^m$, and each layer $\order{i}{U}$ has strictly lower layering dimension.
	Then we let
	\[
		\overbar{h}\isocl{f} \eqdef \overbar{h}\isocl{\restr{f}{\order{1}{U}}} \cp{k} \ldots \cp{k} \overbar{h}\isocl{\restr{f}{\order{m}{U}}}.
	\]
	By construction, if $\overbar{h}$ is well-defined, then it is a strict functor satisfying $\overbar{h} \after (e)_{e \in \grade{n}{\gener{S}}} = \ell$ and restricting to $h$ on $\skel{n-1}{\molecin{P}}$.
	Moreover, let $h'$ be another strict functor with the same property.
	Then $h'$ agrees with $\overbar{h}$ on every cell in $\atomin{\skel{n}{P}}$, which is a basis of $\skel{n}{\molecin{P}}$ by 
	Theorem \ref{thm:molecin_ogp_omega_category} and 
	Proposition \ref{prop:skeleta_of_molecin}.
	It follows from Lemma \ref{lem:functors_equal_on_generating_set} that $h' = \overbar{h}$.
	It only remains to show that $\overbar{h}$ is well-defined, that is, it is independent of the choice of a $k$\nbd layering of $U$ when $\dim{U} = n$ and $k \eqdef \lydim{U} \geq 0$.

	We may assume, inductively, that $\overbar{h}$ is well-defined on all cells $\isocl{g\colon V \to P}$ such that $\dim{V} < n$ or $\lydim{V} < k$.
	Let $(\order{i}{U})_{i=1}^m$ and $(\order{i}{V})_{i=1}^m$ be two $k$\nbd layerings of $U$ and let $(\order{i}{x})_{i=1}^m$, $(\order{i}{y})_{i=1}^m$ be the induced $k$\nbd orderings.
	We now proceed as in the proof of Lemma 
	\ref{lem:frdim_layerings_bijection_with_orderings}, letting $\sigma$ be the unique permutation such that $\order{i}{x} = \order{\sigma(i)}{y}$ for all $i \in \set{1, \ldots, m}$, letting
	\[
		d \eqdef \fun{d}((\order{i}{x})_{i=1}^m, (\order{i}{y})_{i=1}^m)
	\]
	be the number of pairs $(j, j')$ such that $j < j'$ but $\sigma(j') < \sigma(j)$, and proceeding by induction on $d$.
	If $d = 0$, then the two layerings are equal up to layer-wise isomorphism.
	If $d > 0$, then there exists $j < m$ such that $\sigma(j + 1) < \sigma(j)$, and we let $W \submol U$ be the image of $\order{j}{U} \cp{k} \order{j+1}{U}$ in $U$.
	Then $W$ contains exactly two elements 
	\[
		z_1 \eqdef \order{j}{x} = \order{\sigma(j)}{y}, \quad z_2 \eqdef \order{j+1}{x} = \order{\sigma(j+1)}{y}
	\]
	of dimension $> k$, yet there can be no edge between them in $\maxflow{k}{U}$, from which we deduce that $r \eqdef \frdim{W} < k$.
	By assumption, $W$ is frame-acyclic, so by Theorem 
	\ref{thm:frame_acyclic_has_frame_layerings} there exists an $r$\nbd layering of $W$, hence also a pair of molecules $\order{1}{W}$, $\order{2}{W}$, each containing a single element of dimension $> k$, such that $W$ is isomorphic to $\order{1}{W} \cp{r} \order{2}{W}$.
	We may assume, without loss of generality, that $z_1$ is in the image of $\order{1}{W}$ and $z_2$ in the image of $\order{2}{W}$.
	We then have
	\begin{align*}
		& \overbar{h}\isocl{\restr{f}{\order{j}{U}}} \cp{k} 
		\overbar{h}\isocl{\restr{f}{\order{j+1}{U}}} = \\ 
		& \quad = \left( \overbar{h}\isocl{\restr{f}{\order{1}{W}}} \cp{r} 
		\overbar{h}\isocl{\restr{f}{\bound{k}{-}\order{2}{W}}} \right) \cp{k}
		\left( \overbar{h}\isocl{\restr{f}{\bound{k}{+}\order{1}{W}}} \cp{r}
		\overbar{h}\isocl{\restr{f}{\order{2}{W}}} \right),
	\end{align*}
	which by interchange and unitality in $X$ is equal to
	\begin{align*}
		& \overbar{h}\isocl{\restr{f}{\order{1}{W}}} \cp{r} 
		\overbar{h}\isocl{\restr{f}{\order{2}{W}}} = \\ 
		& \quad = \left( \overbar{h}\isocl{\restr{f}{\bound{k}{-}\order{1}{W}}} \cp{r} 
		\overbar{h}\isocl{\restr{f}{\order{2}{W}}} \right) \cp{k}
		\left( \overbar{h}\isocl{\restr{f}{\order{1}{W}}} \cp{r}
		\overbar{h}\isocl{\restr{f}{\bound{k}{+}\order{2}{W}}} \right) = \\
		& \quad = \overbar{h}\isocl{\restr{f}{\order{j}{\tilde{U}}}} \cp{k} \overbar{h}\isocl{\restr{f}{\order{j+1}{\tilde{U}}}},
	\end{align*}
	where we let
	\begin{align*}
    		\order{j}{\tilde{U}} & \eqdef \bound{k}{-}\order{1}{W} \cp{r} \order{2}{W}, \\
		\order{j+1}{\tilde{U}} & \eqdef \order{1}{W} \cp{r} \bound{k}{+} \order{2}{W}.
	\end{align*}
	Notice that all the $n$\nbd dimensional cells in this calculation involve molecules whose layering dimension is $< k$, so $\overbar{h}$ is well-defined on each of them.
	Letting $\order{i}{\tilde{U}} \eqdef \order{i}{U}$ for all $i \not\in \set{j, j+1}$, we have that 
	\begin{enumerate}
		\item $(\order{i}{\tilde{U}})_{i=1}^m$ is a $k$\nbd layering of $U$,
		\item the definition of $\overbar{h}\isocl{f}$ using $(\order{i}{U})_{i=1}^m$ is equal to the one using $(\order{i}{\tilde{U}})$, and 
		\item the induced $k$\nbd ordering $(\order{i}{\tilde{x}})_{i=1}^m \eqdef (\order{1}{x}, \ldots, \order{j+1}{x}, \order{j}{x}, \ldots, \order{m}{x})$ satisfies $\fun{d}((\order{i}{\tilde{x}})_{i=1}^m, (\order{i}{y})_{i=1}^m) < d$,
	\end{enumerate}
	so, by the inductive hypothesis on $d$, the definition of $\overbar{h}\isocl{f}$ using $(\order{i}{\tilde{U}})_{i=1}^m$ is equal to the definition using $(\order{i}{V})_{i=1}^m$.
	We conclude that $\overbar{h}\isocl{f}$ is well-defined, which completes the proof.
\end{proof}

\begin{thm} \label{thm:frame_acyclic_presents_polygraphs}
	Let $P$ be an oriented graded poset with frame-acyclic molecules.
	Then $\molecin{P}$ is a polygraph whose set of generating cells is $\atomin{P}$.
\end{thm}
\begin{proof}
	If $P$ has frame-acyclic molecules, then $\skel{n}{P}$ has frame-acyclic molecules for all $n \in \mathbb{N}$.
	The statement then follows from Lemma \ref{lem:frame_acyclic_cellular_extension}.
\end{proof}

\begin{comm}
	In fact, by the roundness property of atoms, if $\molecin{P}$ is a polygraph, then it is a \emph{regular polygraph} in the sense of \cite{henry2018regular}.
\end{comm}

\begin{dfn}[The category $\rdcpxmapfa$] \index{$\rdcpxmapfa$}
	We let $\rdcpxmapfa$ denote the full subcategory of $\rdcpxmap$ on regular directed complexes with frame-acyclic molecules.
\end{dfn}

\begin{prop} \label{cor:functor_to_polmap}
	Let $p\colon P \to Q$ be a map of regular directed complexes with frame-acyclic molecules.
	Then $\pfw{p} \equiv \molecin{p}\colon \molecin{P} \to \molecin{Q}$ together with the function
	\begin{align*}
		\fpos{\pfw{p}}\colon \atomin{P} &\to \atomin{Q}, \\
		\isocl{\clset{x} \incl P} & \mapsto \isocl{\clset{p(x)} \incl Q}
	\end{align*}
	is a map of polygraphs.
	This determines a functor $\molecin{-}\colon \rdcpxmapfa \to \polmap$, such that the diagram of functors
\[\begin{tikzcd}
	\rdcpxmapfa && \polmap \\
	& \omegacat
	\arrow["{\molecin{-}}", from=1-1, to=1-3]
	\arrow["{\molecin{-}}"', from=1-1, to=2-2]
	\arrow["{\fun{U}}", from=1-3, to=2-2]
\end{tikzcd}\]
	commutes.
\end{prop}
\begin{proof}
	Follows from Theorem \ref{thm:frame_acyclic_has_frame_layerings}, Corollary 
	\ref{cor:basis_of_omegacat_presented_by_rdcpx}, and Lemma 
	\ref{lem:pushforward_of_atom}.
\end{proof}

\begin{exm}[A molecule which is not frame-acyclic] \index[counterex]{A molecule which is not frame-acyclic} \label{exm:not_frame_acyclic}
	We construct a 4\nbd dimensional molecule $\overbar{U}$ which is not frame-acyclic.
	This example is based on \cite[Section 8]{steiner1993algebra}, and has appeared, modulo minor tweaks, as \cite[Example 126]{hadzihasanovic2023higher}.

	We start by constructing a 3\nbd dimensional molecule $U$ which is isomorphic both to the input and to the output boundary of $\overbar{U}$.
	This molecule is isomorphic to $\order{1}{U} \cp{2} \order{2}{U}$ (this is \emph{not} a layering), where
	\begin{align*}
		\order{1}{U} & \eqdef 
		\disk{1}{2} \cp{1} ((\globe{2} \celto (\globe{2} \cp{1} \globe{2}))
		\cp{0} (\globe{2} \celto (\globe{2} \cp{1} \globe{2}))) \cp{1} \disk{2}{1}, \\
		\order{2}{U} & \eqdef 
		((\disk{1}{2} \cp{1} (\globe{2} \cp{0} \globe{2})) \celto \disk{1}{2})
		\cp{1} (((\globe{2} \cp{0} \globe{2}) \cp{1} \disk{2}{1}) \celto \disk{2}{1}).
	\end{align*}
	The graph of $U$ is
	\[
		\input{img/frameac_bound_graph.tex}
	\]
	and there are several possible 2\nbd orderings, all of which determine 2\nbd layerings.
	For example, the 2\nbd ordering $((3, 2), (3, 1), (3, 0), (3, 3))$ corresponds to the sequence of rewrite steps
	\[
		\input{img/frameac_bound_step0.tex} \quad
		\input{img/frameac_bound_step1.tex} \quad
		\input{img/frameac_bound_step2.tex}
	\]
	\[
		\input{img/frameac_bound_step3.tex} \quad
		\input{img/frameac_bound_step4.tex} \quad
	\]
	where we omitted the wire labels.
	Now, if we let
	\[
		V \eqdef \clset{(3, 1), (3, 0)}, \quad \quad W \eqdef \clset{(3, 2), (3, 3)},
	\]
	then both $V$ and $W$ are rewritable submolecules of $U$ that only overlap in their boundaries, but they are \emph{not} ``simultaneously'' rewritable, that is, the multiple substitution of $V$ and $W$ does not necessarily produce a molecule: to see this, observe that the inclusion of $V$ into $\subs{U}{\compos{W}}{W}$ is essentially a variant of the inclusion of Example \ref{exm:non_submol}.
	
	We now construct $\overbar{U}$ as the 4\nbd dimensional molecule corresponding to the following pair of pairs of rewrites, performed in either order:
	\begin{itemize}
		\item rewrite $V$ to $\compos{V}$, then revert,
		\item rewrite $W$ to $\compos{W}$, then revert.
	\end{itemize}
	After each pair of rewrites we obtain a molecule isomorphic to $U$ by Lemma 
	\ref{lem:revert_substitution}, which justifies performing the two pairs in either order.

	The oriented Hasse diagram of $\overbar{U}$ is
	\[
		\input{img/frameac_hasse.tex}
	\]
	and its graph is
	\[
		\input{img/frameac_graph.tex}
	\]
	which has two connected components, reflecting the fact that the two pairs of rewrites operate on submolecules that do not share any 3\nbd dimensional elements.
	The sequence of rewrite steps where we rewrite $V$ first corresponds to the 3\nbd layering $(\order{i}{\overbar{U}})_{i=1}^4$ inducing the 3\nbd ordering $((4, 0), (4, 2), (4, 1), (4, 3))$.
	Let 
	\[
		\overbar{V} \eqdef \order{2}{\overbar{U}} \cp{3} \order{3}{\overbar{U}} = \clset{(4, 2), (4, 1)}.
	\]
	Then $\overbar{V}$ is a submolecule of $\overbar{U}$, and $\frdim{\overbar{V}} = 2$ because $\clset{(4, 1)}$ and $\clset{(4, 2)}$ do not share any 3\nbd dimensional elements.
	However,
	\[
		(2, 7) \in \faces{2}{+}(4, 2) \cap \faces{2}{-}(4, 1), \quad \quad
		(2, 5) \in \faces{2}{+}(4, 1) \cap \faces{2}{-}(4, 2),
	\]
	leading to a cycle in $\maxflow{2}{\overbar{V}}$.
	This proves that $\overbar{V}$ and $\overbar{U}$ are not frame-acyclic.
\end{exm}

\begin{exm}[An ordering which is not induced by a layering] \index[counterex]{An ordering which is not induced by a layering}
	Let $\overbar{U}$ be the non frame-acyclic molecule of Example \ref{exm:not_frame_acyclic}.
	We already know, by Corollary \ref{cor:frame_acyclicity_equivalent_conditions}, that some submolecule of $\overbar{U}$ must admit an ordering which does not come from a layering.
	In fact, we will show something stronger.
	The 3\nbd flow graph of $\overbar{U}$ is 
\[\begin{tikzcd}[sep=small]
	{{\scriptstyle (4, 0)}\;\bullet} && {{\scriptstyle (4, 2)}\;\bullet} \\
	{{\scriptstyle (4, 1)}\;\bullet} && {{\scriptstyle (4, 3)}\;\bullet}
	\arrow[from=1-1, to=1-3]
	\arrow[from=2-1, to=2-3]
\end{tikzcd}\]
	which admits the topological sort $((4, 1), (4, 0), (4, 2), (4, 3))$.
	However, this 3\nbd ordering is \emph{not} induced by a 3\nbd layering: while the first rewrite is well-defined, the resulting 3\nbd dimensional molecule is isomorphic to $\subs{U}{\compos{W}}{W}$, and the inclusion of $\bound{}{-}(4, 0) = V$ into it is not a submolecule inclusion, so it does not pass the test of Proposition \ref{prop:layering_from_ordering}.

	In fact, there are only two valid 3\nbd layerings of $\overbar{U}$: the first induces the 3\nbd ordering $((4, 0), (4, 2), (4, 1), (4, 3))$, corresponding to the sequence
	\[
		\input{img/frameac_step0.tex} \quad
		\input{img/frameac_step1.tex}
	\]
	\[
		\input{img/frameac_step2.tex} \quad
		\input{img/frameac_step3.tex}
	\]
	\[
		\input{img/frameac_step4.tex}
	\]
	while the second induces the 3\nbd ordering $((4, 1), (4, 3), (4, 0), (4, 2))$, corresponding to the sequence with the same initial and final step, but with
	\[
		\input{img/frameac_bis_step1.tex} \quad
		\input{img/frameac_bis_step2.tex}
	\]
	\[
		\input{img/frameac_bis_step3.tex}
	\]
	as intermediate steps.
	This shows not only that not every 3\nbd ordering of $\overbar{U}$ is induced by a 3\nbd layering, but that \emph{there is no graph whose vertices are the elements of $\maxel{\overbar{U}}$, and whose topological sorts correspond to the 3\nbd layerings of $\overbar{U}$}, since no graph on 4 vertices admits exactly 2 topological sorts: that is, we cannot fix the problem simply by adding some extra edges to maximal flow graphs.
\end{exm}

\begin{exm}[A regular directed complex $P$ such that $\molecin{P}$ is not a polygraph] \index[counterex]{A regular directed complex $P$ such that $\molecin{P}$ is not a polygraph}
	Once more, we can take $P$ to be the non frame-acyclic molecule $\overbar{U}$ of Example \ref{exm:not_frame_acyclic} (this fact was suggested to us by F.~Loubaton).

	Let $\order{1}{V} \cp{3} \order{2}{V}$ and $\order{1}{W} \cp{3} \order{2}{W}$ be the decompositions of $\overbar{U}$ induced by the 3\nbd orderings
	\[
		((4, 0), (4, 2), (4, 1), (4, 3)), \quad \quad ((4, 1), (4, 3), (4, 0), (4, 2)),
	\]
	respectively; that is,
	\begin{align*}
		\grade{4}{(\order{1}{V})} & = \grade{4}{(\order{2}{W})} = \set{(4, 0), (4, 2)}, \\
		\grade{4}{(\order{2}{V})} & = \grade{4}{(\order{1}{W})} = \set{(4, 1), (4, 3)}.
	\end{align*}
	Then, in the strict 4\nbd category $\molecin{\overbar{U}}$,
	\[
		\isocl{\order{1}{V} \incl \overbar{U}} \cp{3} 
		\isocl{\order{2}{V} \incl \overbar{U}} = 
		\isocl{\order{1}{W} \incl \overbar{U}} \cp{3}
		\isocl{\order{2}{W} \incl \overbar{U}}
	\]
	since both sides are equal to $\isocl{\idd{\overbar{U}}}$, yet we claim that this identity is not provable from the equations of strict $\omega$\nbd categories, for any pasting decompositions of $\order{1}{V}$, $\order{2}{V}$, $\order{1}{W}$, and $\order{2}{W}$ into atoms.

	Indeed, in any decomposition of the left-hand side, $\clset{(4, 2)}$ must appear to the left of $\clset{(4, 1)}$, while in any decomposition of the right-hand side, it must appear to the right.
	An examination of the axioms of strict $\omega$\nbd categories shows that the only equation that can invert the order of two factors is the interchange equation, which, to be applicable, would require $\clset{(4, 2)}$ and $\clset{(4, 1)}$ to appear as factors in opposite sides of a pasting of the form $V' \cp{2} W'$.
	But this is not possible due to the existence of a cycle $(4, 2) \to (4, 1) \to (4, 2)$ in $\flow{2}{\overbar{U}}$, as observed in Example \ref{exm:not_frame_acyclic}.
	We conclude that $\molecin{\overbar{U}}$ is not a polygraph.
\end{exm}

\begin{comm}
	The fact that there exist non frame-acyclic molecules $U$, such that $\molecin{U}$ is not a polygraph, implies that the axioms of strict $\omega$\nbd categories are \emph{not} complete for the equations satisfied by the pasting of molecules in dimension $> 3$.
	It is an open question whether there is a finitary axiomatisation of all the equations satisfied by pasting.

	Since pasting of molecules is ``topologically sound'' by the results of Chapter 
	\ref{chap:geometric}, this leads to the somewhat counterintuitive realisation that, while the equations of strict $\omega$\nbd categories are \emph{too strict} as an axiomatisation of \emph{composition} in an $\infty$\nbd groupoid, they are \emph{not strong enough} as an axiomatisation of the relations that exist in the pasting of topological cells.
	The paradox is resolved by the fact that composition and pasting are fundamentally different operations: the first is a ``closed'' operation, which needs to output an object of the same type as its factors, while the second outputs an object of a different type.
	The necessity of weakness is tied to the demand that composition be closed: compare the ``associative-up-to-homotopy'' composition of paths in the fundamental groupoid of a space with the strictly associative composition of \emph{Moore paths}, which is really a \emph{pasting} of paths.
\end{comm}


\section{Stronger acyclicity conditions} \label{sec:stronger_acyclic}

\begin{guide}
	In this section, we study three acyclicity conditions of decreasing strength.

	The first condition, called simply \cemph{acyclicity}, has the strongest consequences for a generic oriented graded poset $P$: it implies that $P$ has frame-acyclic molecules (Proposition 
	\ref{prop:acyclic_has_frame_acyclic_molecules}) and that every morphism from a molecule to $P$ is an inclusion (Proposition \ref{prop:morphism_from_molecule_to_acyclic_is_injective}).
	Acyclic oriented graded posets are closed under pasting (Lemma 
	\ref{lem:acyclic_closed_under_pasting}), Gray products (Proposition 
	\ref{prop:gray_product_of_acyclic}), and joins (Proposition 
	\ref{prop:join_of_acyclic}), but not under duals, except for total duals (Proposition 
	\ref{prop:acyclic_stable_under_total_dual}).

	The second and the third condition, called \cemph{strong dimension-wise acyclicity} and \cemph{dimension-wise acyclicity}, imply that $P$ has frame-acyclic molecules when $P$ is a regular directed complex.
	Strongly dimension-wise acyclic and dimension-wise acyclic oriented graded posets are closed under all duals (Proposition \ref{prop:dw_acyclic_stable_under_all_duals}), but not, in general, under Gray products and joins.
	Moreover, strong dimension-wise acyclicity of an oriented graded poset $P$ implies that all \emph{local embeddings} of a molecule over $P$ are inclusions (Proposition 
	\ref{prop:molecule_over_strongly_dimensionwise_acyclic}).

	All acyclicity conditions are stable under suspensions (Proposition 
	\ref{prop:suspension_of_acyclic}).
\end{guide}

\begin{dfn}[Acyclic oriented graded poset] \index{oriented graded poset!acyclic} \index{molecule!acyclic}
Let $P$ be an oriented graded poset.
We say that $P$ is \emph{acyclic} if $\hasseo{P}$ is acyclic.
\end{dfn}

\begin{dfn}[Strongly dimension-wise acyclic oriented graded poset] \index{oriented graded poset!acyclic!strongly dimension-wise} \index{molecule!acyclic!strongly dimension-wise}
Let $P$ be an oriented graded poset.
We say that $P$ is \emph{strongly dimension-wise acyclic} if, for all $k \geq -1$, $\extflow{k}{P}$ is acyclic.
\end{dfn}

\begin{dfn}[Dimension-wise acyclic oriented graded poset] \index{oriented graded poset!acyclic!dimension-wise} \index{molecule!acyclic!dimension-wise}
Let $P$ be an oriented graded poset.
We say that $P$ is \emph{dimension-wise acyclic} if, for all $k \geq -1$, $\flow{k}{P}$ is acyclic.
\end{dfn}

\begin{rmk}
	Note that the condition that $\extflow{k}{P}$ or $\flow{k}{P}$ be acyclic is trivial when $k = -1$ or when $k \geq \dim{P}$.
\end{rmk}

\begin{comm} \label{comm:acyclicity_in_steiner}
Acyclicity and strong dimension-wise acyclicity correspond, essentially, to \emph{total loop-freeness} and \emph{loop-freeness} in \cite{steiner1993algebra}.
Theorem \ref{thm:frame_acyclic_has_frame_layerings}, on the other hand, shows that frame-acyclicity corresponds to what Steiner calls being \emph{split}.
We use different terminology both because we find it more informative, and because both our definitions and the underlying framework differ slightly but significantly.
\end{comm}

\begin{prop} \label{prop:acyclicity_implications}
Let $P$ be a regular directed complex.
Then
\begin{enumerate}
	\item if $P$ is acyclic, then it is strongly dimension-wise acyclic,
	\item if $P$ is strongly dimension-wise acyclic, then it is dimension-wise acyclic.
\end{enumerate}
\end{prop}
\begin{proof}
For all $k \geq -1$, by Lemma \ref{lem:path_in_flow_induces_path_in_extflow} a cycle in $\flow{k}{U}$ induces a cycle in $\extflow{k}{U}$.
By Lemma \ref{lem:path_in_extflow_induces_path_in_hasseo}, a cycle in $\extflow{k}{U}$ induces a cycle in $\hasseo{U}$.
\end{proof}

\begin{exm}[A molecule which is strongly dimension-wise acyclic but not acyclic] \index[counterex]{A molecule which is strongly dimension-wise acyclic but not acyclic} \label{exm:non_acyclic}
	We prove that the first implication in Proposition \ref{prop:acyclicity_implications} is strict.
	Let $U$ be a 3\nbd dimensional atom whose input and output boundaries are the oriented face posets of
	\[\begin{tikzcd}[column sep=small]
	& {{\scriptstyle 3}\; \bullet} \\
	{{\scriptstyle 0}\; \bullet} &&& {{\scriptstyle 2}\; \bullet} \\
	&& {{\scriptstyle 1}\; \bullet}
	\arrow["0"', curve={height=12pt}, from=2-1, to=3-3]
	\arrow["1"', curve={height=6pt}, from=3-3, to=2-4]
	\arrow["2", curve={height=-6pt}, from=2-1, to=1-2]
	\arrow["3", curve={height=-12pt}, from=1-2, to=2-4]
	\arrow["0", curve={height=6pt}, Rightarrow, from=3-3, to=1-2]
\end{tikzcd} \quad \text{and} \quad
\begin{tikzcd}[column sep=small]
	& {{\scriptstyle 3}\; \bullet} \\
	{{\scriptstyle 0}\; \bullet} &&& {{\scriptstyle 2}\; \bullet} \\
	&& {{\scriptstyle 1}\; \bullet}
	\arrow[""{name=0, anchor=center, inner sep=0}, "0"', curve={height=12pt}, from=2-1, to=3-3]
	\arrow["1"', curve={height=6pt}, from=3-3, to=2-4]
	\arrow["2", curve={height=-6pt}, from=2-1, to=1-2]
	\arrow[""{name=1, anchor=center, inner sep=0}, "3", curve={height=-12pt}, from=1-2, to=2-4]
	\arrow["4", from=1-2, to=3-3]
	\arrow["1"', curve={height=-6pt}, shorten <=7pt, Rightarrow, from=0, to=1-2]
	\arrow["2"', curve={height=6pt}, shorten >=7pt, Rightarrow, from=3-3, to=1]
\end{tikzcd}\]
	respectively.
	Then the extended 0\nbd flow graph $\extflow{0}{U}$ is
	\[\begin{tikzcd}[sep=small]
	& {{\scriptstyle (1, 2)}\;\bullet} & {{\scriptstyle (0, 3)}\;\bullet} &&& {{\scriptstyle (1, 3)}\;\bullet} \\
	& {{\scriptstyle (2, 0)}\;\bullet} &&&& {{\scriptstyle (2, 2)}\;\bullet} \\
	{{\scriptstyle (0, 0)}\;\bullet} &&& {{\scriptstyle (1, 4)}\;\bullet} &&& {{\scriptstyle (0, 2)}\;\bullet} \\
	& {{\scriptstyle (2, 1)}\;\bullet} &&&& {{\scriptstyle (3, 0)}\;\bullet} \\
	& {{\scriptstyle (1, 0)}\;\bullet} &&& {{\scriptstyle (0, 1)}\;\bullet} & {{\scriptstyle (1, 1)}\;\bullet}
	\arrow[from=5-2, to=5-5]
	\arrow[from=5-5, to=5-6]
	\arrow[from=5-6, to=3-7]
	\arrow[from=1-6, to=3-7]
	\arrow[from=2-6, to=3-7]
	\arrow[from=2-2, to=3-7]
	\arrow[from=3-4, to=5-5]
	\arrow[from=4-2, to=5-5]
	\arrow[from=1-3, to=3-4]
	\arrow[from=1-3, to=2-6]
	\arrow[from=1-3, to=1-6]
	\arrow[from=1-2, to=1-3]
	\arrow[from=3-1, to=1-2]
	\arrow[from=3-1, to=2-2]
	\arrow[from=3-1, to=4-2]
	\arrow[from=3-1, to=5-2]
	\arrow[from=3-1, to=4-6]
	\arrow[from=4-6, to=3-7]
\end{tikzcd}\]
	while the extended 1\nbd flow graph $\extflow{1}{U}$ is
\[\begin{tikzcd}[column sep=small]
	& {{\scriptstyle (1, 1)}\;\bullet} && {{\scriptstyle (2, 0)}\;\bullet} && {{\scriptstyle (1, 2)}\;\bullet} \\
	{{\scriptstyle (0, 0)}\;\bullet} & {{\scriptstyle (0, 1)}\;\bullet} && {{\scriptstyle (3, 0)}\;\bullet} && {{\scriptstyle (0, 3)}\;\bullet} & {{\scriptstyle (0, 2)}\;\bullet} \\
	& {{\scriptstyle (1, 0)}\;\bullet} & {{\scriptstyle (2, 1)}\;\bullet} & {{\scriptstyle (1, 4)}\;\bullet} & {{\scriptstyle (2, 2)}\;\bullet} & {{\scriptstyle (1, 3)}\;\bullet}
	\arrow[from=3-2, to=1-4]
	\arrow[from=1-4, to=1-6]
	\arrow[from=1-4, to=2-6]
	\arrow[from=1-4, to=3-6]
	\arrow[from=2-2, to=1-4]
	\arrow[from=1-2, to=1-4]
	\arrow[from=3-2, to=3-3]
	\arrow[from=3-3, to=1-6]
	\arrow[from=3-3, to=2-6]
	\arrow[from=3-3, to=3-4]
	\arrow[from=3-4, to=3-5]
	\arrow[from=2-2, to=3-5]
	\arrow[from=1-2, to=3-5]
	\arrow[from=3-5, to=3-6]
	\arrow[from=1-2, to=2-4]
	\arrow[from=2-2, to=2-4]
	\arrow[from=3-2, to=2-4]
	\arrow[from=2-4, to=1-6]
	\arrow[from=2-4, to=2-6]
	\arrow[from=2-4, to=3-6]
\end{tikzcd}\]
	and the extended 2\nbd flow graph $\extflow{2}{U}$ is
\[\begin{tikzcd}[sep=small]
	&& {{\scriptstyle (2, 2)}\;\bullet} & {{\scriptstyle (0, 0)}\;\bullet} & {{\scriptstyle (0, 3)}\;\bullet} & {{\scriptstyle (1, 2)}\;\bullet} \\
	{{\scriptstyle (2, 0)}\;\bullet} & {{\scriptstyle (3, 0)}\;\bullet} & {{\scriptstyle (1, 4)}\;\bullet} & {{\scriptstyle (0, 1)}\;\bullet} & {{\scriptstyle (1, 0)}\;\bullet} & {{\scriptstyle (1, 3)}\;\bullet} \\
	&& {{\scriptstyle (2, 1)}\;\bullet} & {{\scriptstyle (0, 2)}\;\bullet} & {{\scriptstyle (1, 1)}\;\bullet}
	\arrow[from=2-1, to=2-2]
	\arrow[from=2-2, to=3-3]
	\arrow[from=2-2, to=2-3]
	\arrow[from=2-2, to=1-3]
\end{tikzcd}\]
	all of which are acyclic.
	All other extended flow graphs are discrete, so $U$ is strongly dimension-wise acyclic.
	However, the oriented Hasse diagram
	\[
		\input{img/nonacy_hasse.tex}
	\]
	contains the cycle
	\[
		(0, 1) \to (1, 1) \to (2, 0) \to (3, 0) \to (2, 1) \to (1, 4) \to (0, 1)
	\]
	so $U$ is not acyclic.
\end{exm}

\begin{exm}[A regular directed complex which is dimension-wise acyclic but not strongly dimension-wise acyclic] \index[counterex]{A regular directed complex which is dimension-wise acyclic but not strongly dimension-wise acyclic} \label{exm:nonsdw}
	We prove that the second implication in Proposition 
	\ref{prop:acyclicity_implications} is also strict.
	Let $U$ be the oriented face poset of
\[\begin{tikzcd}
	& {{\color{\mycolor} \bullet}} \\
	\bullet && \bullet \\
	& {{\color{\mycolor} \bullet}}
	\arrow[curve={height=6pt}, from=2-1, to=3-2]
	\arrow[""{name=0, anchor=center, inner sep=0}, from=2-1, to=2-3]
	\arrow[curve={height=6pt}, from=3-2, to=2-3]
	\arrow[curve={height=-6pt}, from=2-1, to=1-2]
	\arrow[curve={height=-6pt}, from=1-2, to=2-3]
	\arrow[shorten >=3pt, Rightarrow, from=3-2, to=0]
	\arrow[shorten <=3pt, Rightarrow, from=0, to=1-2]
\end{tikzcd}\]
	which is a 2\nbd dimensional molecule, and let $P$ be the result of identifying the two 0\nbd dimensional cells marked in a different colour.
	This is a regular directed complex whose oriented Hasse diagram is
	\[
		\input{img/nonsdw_hasse.tex} 
	\]
	and whose 0\nbd flow graph $\flow{0}{P}$ and 1\nbd flow graph $\flow{1}{P}$ are
	\[\begin{tikzcd}[sep=tiny]
	{{\scriptstyle (1, 3)}\;\bullet} && {{\scriptstyle (1, 4)}\;\bullet} \\
	{{\scriptstyle (1, 0)}\;\bullet} && {{\scriptstyle (1, 1)}\;\bullet} \\
	{{\scriptstyle (2, 0)}\;\bullet} & {{\scriptstyle (1, 2)}\;\bullet} & {{\scriptstyle (2, 1)}\;\bullet}
	\arrow[from=2-1, to=1-3]
	\arrow[from=2-1, to=2-3]
	\arrow[from=1-1, to=2-3]
	\arrow[from=1-1, to=1-3]
	\end{tikzcd} \quad \text{and} \quad \begin{tikzcd}[sep=tiny]
	{{\scriptstyle (2, 0)}\;\bullet} && {{\scriptstyle (2, 1)}\;\bullet}
	\arrow[from=1-1, to=1-3]
\end{tikzcd}\]
	respectively.
	All other flow graphs are discrete, so $P$ is dimension-wise acyclic.
	However, the extended 1\nbd flow graph $\extflow{1}{P}$ is
	\[\begin{tikzcd}[sep=small]
	{{\scriptstyle (1, 0)}\;\bullet} && {{\scriptstyle (0, 1)}\;\bullet} && {{\scriptstyle (1, 4)}\;\bullet} \\
	{{\scriptstyle (1, 1)}\;\bullet} & {{\scriptstyle (2, 0)}\;\bullet} & {{\scriptstyle (1, 2)}\;\bullet} & {{\scriptstyle (2, 1)}\;\bullet} & {{\scriptstyle (1, 3)}\;\bullet} \\
	& {{\scriptstyle (0, 0)}\;\bullet} && {{\scriptstyle (0, 2)}\;\bullet}
	\arrow[from=2-2, to=2-3]
	\arrow[from=2-3, to=2-4]
	\arrow[from=2-4, to=1-3]
	\arrow[from=1-3, to=2-2]
	\arrow[from=2-4, to=1-5]
	\arrow[from=2-4, to=2-5]
	\arrow[from=2-1, to=2-2]
	\arrow[from=1-1, to=2-2]
\end{tikzcd}\]
	which contains the cycle
	\[
		(0, 1) \to (2, 0) \to (1, 2) \to (2, 1) \to (0, 1),
	\]
	so $P$ is not strongly dimension-wise acyclic.
\end{exm}

\begin{lem} \label{lem:flow_under_inclusion}
Let $\imath\colon P \incl Q$ be an inclusion of oriented graded posets, $k \geq -1$.
Then 
\begin{enumerate}
	\item $\flow{k}{P}$ is isomorphic to the induced subgraph of $\flow{k}{Q}$, and
	\item $\extflow{k}{P}$ is isomorphic to the induced subgraph of $\extflow{k}{Q}$
\end{enumerate}
on the vertices in the image of $\imath$.
\end{lem}
\begin{proof}
	Follows by definition from Proposition \ref{prop:inclusions_preserve_faces} and Corollary \ref{cor:inclusions_preserve_boundaries}.
\end{proof}

\begin{lem} \label{lem:if_flow_acyclic_then_maxflow_acyclic}
Let $U$ be a molecule, $k \geq -1$.
If $\flow{k}{U}$ is acyclic, then $\maxflow{k}{U}$ is acyclic, and $U$ admits a $k$\nbd ordering.
\end{lem}
\begin{proof}
Every induced subgraph of a directed acyclic graph is acyclic, and every directed acyclic graph admits a $k$\nbd ordering.
\end{proof}

\begin{prop} \label{prop:dimensionwise_implies_frame_acyclic}
Let $U$ be a molecule.
If $U$ is dimension-wise acyclic, then $U$ is frame-acyclic.
\end{prop}
\begin{proof}
Let $V \submol U$ be a submolecule inclusion, $r \eqdef \frdim{V}$.
If $U$ is dimension-wise acyclic, then $\flow{r}{U}$ is acyclic.
Then by Lemma \ref{lem:flow_under_inclusion} so is $\flow{r}{V}$, and by Lemma \ref{lem:if_flow_acyclic_then_maxflow_acyclic} so is $\maxflow{r}{V}$.
\end{proof}

\begin{exm}[A molecule which is frame-acyclic but not dimension-wise acyclic] \index[counterex]{A molecule which is frame-acyclic but not dimension-wise acyclic} \label{exm:non_dw_acy}
	We prove that the implication of Proposition \ref{prop:dimensionwise_implies_frame_acyclic} is strict.
	Let $U$ be a 3\nbd dimensional atom whose input and output boundaries are the oriented face posets of the two sides of (\ref{eq:triangle_equation}), that is,
\[
\begin{tikzcd}[column sep=small]
	&& {{\scriptstyle 3}\; \bullet} \\
	{{\scriptstyle 0}\; \bullet} &&& {{\scriptstyle 2}\; \bullet} \\
	& {{\scriptstyle 1}\; \bullet}
	\arrow[""{name=0, anchor=center, inner sep=0}, "1"', curve={height=12pt}, from=3-2, to=2-4]
	\arrow["0"', curve={height=6pt}, from=2-1, to=3-2]
	\arrow["4", curve={height=-6pt}, from=1-3, to=2-4]
	\arrow[""{name=1, anchor=center, inner sep=0}, "3", curve={height=-12pt}, from=2-1, to=1-3]
	\arrow["2", from=3-2, to=1-3]
	\arrow["1", curve={height=6pt}, shorten <=7pt, Rightarrow, from=0, to=1-3]
	\arrow["0", curve={height=-6pt}, shorten >=7pt, Rightarrow, from=3-2, to=1]
\end{tikzcd}
	\quad \text{and} \quad
	\begin{tikzcd}[column sep=small]
	& {{\scriptstyle 3}\; \bullet} \\
	{{\scriptstyle 0}\; \bullet} &&& {{\scriptstyle 2}\; \bullet} \\
	&& {{\scriptstyle 1}\; \bullet}
	\arrow[""{name=0, anchor=center, inner sep=0}, "0"', curve={height=12pt}, from=2-1, to=3-3]
	\arrow["1"', curve={height=6pt}, from=3-3, to=2-4]
	\arrow["3", curve={height=-6pt}, from=2-1, to=1-2]
	\arrow[""{name=1, anchor=center, inner sep=0}, "4", curve={height=-12pt}, from=1-2, to=2-4]
	\arrow["5", from=1-2, to=3-3]
	\arrow["2"', curve={height=-6pt}, shorten <=7pt, Rightarrow, from=0, to=1-2]
	\arrow["3"', curve={height=6pt}, shorten >=7pt, Rightarrow, from=3-3, to=1]
\end{tikzcd}\]
	respectively.
	Then $\flow{0}{U}$ contains the cycle
	\[
		(1, 2) \to (1, 5) \to (1, 2),
	\]
	so $U$ is not dimension-wise acyclic.
	However, since $U$ is 3\nbd dimensional, it is frame-acyclic as a consequence of Theorem 
	\ref{thm:dim3_frame_acyclic}.
\end{exm}

\begin{prop} \label{prop:acyclic_flow_is_linear_order}
Let $U$ be an acyclic molecule.
Then the flow preorder $\precflow$ is a linear order on $U$.
\end{prop}
\begin{proof}
	Let $x, y \in U$.
	By Proposition \ref{prop:acyclicity_implications} combined with Proposition \ref{prop:dimensionwise_implies_frame_acyclic}, $U$ is frame-acyclic, so by Proposition \ref{prop:hasseo_connected_frame_acyclic} either $x \precflow y$ or $y \precflow x$.
	Since $\hasseo{U}$ is acyclic, if both hold then $x = y$.
\end{proof}

\begin{lem} \label{lem:morphism_into_acyclic}
Let $f\colon P \to Q$ be a morphism of oriented graded posets.
If $Q$ is acyclic, then $P$ is acyclic.
\end{lem}
\begin{proof}
Suppose that there is a cycle in $\hasseo{P}$.
By Proposition \ref{prop:hasse_graph_homomorphism}, $\hasseo{f}$ maps it onto a cycle in $\hasseo{Q}$.
\end{proof}

\begin{prop} \label{prop:acyclic_has_frame_acyclic_molecules}
	Let $P$ be an acyclic oriented graded poset.
	Then $P$ has frame-acyclic molecules.
\end{prop}
\begin{proof}
	Let $U$ be a molecule and $f\colon U \to P$ be a morphism.
	By Lemma \ref{lem:morphism_into_acyclic}, $U$ is acyclic, so by Proposition
	\ref{prop:acyclicity_implications} combined with 
	Proposition \ref{prop:dimensionwise_implies_frame_acyclic}, it is frame-acyclic.
\end{proof}

\begin{cor} \label{cor:acyclic_presents_polygraphs}
	Let $P$ be an acyclic oriented graded poset.
	Then $\molecin{P}$ is a polygraph.
\end{cor}

\begin{prop} \label{prop:morphism_from_molecule_to_acyclic_is_injective}
Let $U$ be a molecule, $P$ an acyclic oriented graded poset, and $f\colon U \to P$ a morphism.
Then $f$ is an inclusion.
\end{prop}
\begin{proof}
	Let $x, y \in U$ and suppose that $f(x) = f(y)$.
	By Lemma \ref{lem:morphism_into_acyclic}, $U$ is acyclic, so by Proposition \ref{prop:acyclic_flow_is_linear_order} there exists a path from $x$ to $y$ or a path from $y$ to $x$ in $\hasseo{U}$.
	Then $\hasseo{f}$ maps this onto a cycle in $\hasseo{P}$, a contradiction, unless $x = y$ and the path is constant.
	We conclude that $f$ is injective.
\end{proof}

\begin{cor} \label{cor:acyclic_omega_cat_basis}
Let $P$ be an acyclic oriented graded poset.
Then
\begin{align*}
	\molecin{P} & = \set{ \isocl{ U \incl P } \mid U \subseteq P, \text{$U$ is a molecule} }, \\
	\atomin{P} & = \set{ \isocl{ \clset{x} \incl P } \mid x \in P, \text{$\clset{x}$ is an atom} }.
\end{align*}
\end{cor}
\begin{proof}
	By Proposition \ref{prop:morphism_from_molecule_to_acyclic_is_injective}, every morphism from a molecule to $P$ is an inclusion, equivalent to a subset inclusion $U \incl P$ for some closed subset $U \subseteq P$.
	In particular, every morphism from an atom to $P$ is equivalent to the inclusion $\clset{x} \incl P$ for some $x \in P$.
\end{proof}

\begin{prop} \label{prop:extflow_graph_homomorphism}
	Let $f\colon P \to Q$ be a local embedding of oriented graded posets.
	Then, for all $k \geq -1$, $f$ induces homomorphisms
	\[
		\flow{k}{f}\colon \flow{k}{P} \to \flow{k}{Q} \quad \text{and} \quad
		\extflow{k}{f}\colon \extflow{k}{P} \to \extflow{k}{Q}.
	\]
	These assignments determine functors $\flow{k}, \extflow{k}\colon \ogposloc \to \gph$.
\end{prop}
\begin{proof}
	For all $x \in P$, the restriction $\restr{f}{\clset{x}}$ is an inclusion.
	By Proposition \ref{prop:inclusions_preserve_faces}, for all $\alpha \in \set{+, -}$, if $y \in \faces{k}{\alpha}x$, then $f(y) \in \faces{k}{\alpha}f(x)$.
	Thus, if there is an edge between $x$ and $y$ in $\flow{k}{P}$, then there is an edge between $f(x)$ and $f(y)$ in $\flow{k}{Q}$.
	The case of extended $k$\nbd flow graphs is similar, using Corollary 
	\ref{cor:inclusions_preserve_boundaries}, and functoriality is straightforward.
\end{proof}

\begin{cor} \label{cor:local_embedding_into_dimensionwise_acyclic}
Let $f\colon P \to Q$ be a local embedding of oriented graded posets.
Then
\begin{enumerate}
	\item if $Q$ is dimension-wise acyclic, then so is $P$,
	\item if $Q$ is strongly dimension-wise acyclic, then so is $P$.
\end{enumerate}
\end{cor}
\begin{proof}
	Let $k \geq -1$ and suppose that there is a cycle in $\flow{k}{P}$ or $\extflow{k}{P}$.
	By Proposition \ref{prop:extflow_graph_homomorphism}, $\flow{k}{f}$ or $\extflow{k}{f}$ maps it onto a cycle in $\flow{k}{Q}$ or $\extflow{k}{Q}$.
\end{proof}

\begin{prop} \label{prop:dw_acyclic_rdcpx_has_frame_acyclic_molecules}
	Let $P$ be a dimension-wise acyclic regular directed complex.
	Then $P$ has frame-acyclic molecules.
\end{prop}
\begin{proof}
	Let $U$ be a molecule and $f\colon U \to P$ a morphism.
	By Corollary \ref{cor:morphisms_of_rdcpx_are_local_isomorphisms}, $f$ is a local embedding, so by Corollary \ref{cor:local_embedding_into_dimensionwise_acyclic} $U$ is dimension-wise acyclic.
	It follows from Proposition \ref{prop:dimensionwise_implies_frame_acyclic} that $U$ is frame-acyclic.
\end{proof}

\begin{cor} \label{cor:dw_acyclic_rdcpx_presents_polygraphs}
	Let $P$ be a dimension-wise acyclic regular directed complex.
	Then $\molecin{P}$ is a polygraph.
\end{cor}

\begin{prop} \label{prop:molecule_over_strongly_dimensionwise_acyclic}
Let $U$ be a molecule, $P$ a strongly dimension-wise acyclic oriented graded poset, and $f\colon U \to P$ a local embedding.
Then $f$ is an inclusion.
\end{prop}
\begin{proof}
Let $x, y \in U$ and suppose that $f(x) = f(y)$.
By Corollary \ref{cor:local_embedding_into_dimensionwise_acyclic}, $U$ is strongly dimension-wise acyclic.
It follows from Proposition \ref{prop:acyclicity_implications} and \ref{prop:dimensionwise_implies_frame_acyclic} that $U$ is frame acyclic, so by Proposition \ref{prop:extflow_connected_frame_acyclic} there exists $k \geq -1$ such that there is a path from $x$ to $y$ or a path from $y$ to $x$ in $\extflow{k}{U}$.
Then by Proposition \ref{prop:extflow_graph_homomorphism} $\extflow{k}{f}$ maps this onto a cycle in $\extflow{k}{P}$, a contradiction, unless $x = y$ and the path is constant.
We conclude that $f$ is injective.
\end{proof}

\begin{exm}[A local embedding of a molecule into a dimension-wise acyclic regular directed complex which is not an inclusion] \index[counterex]{A local embedding of a molecule into a dimension-wise acyclic regular directed complex which is not an inclusion}
	We prove that Proposition 
	\ref{prop:molecule_over_strongly_dimensionwise_acyclic} does not extend to dimension-wise acyclic regular directed complexes.
	
	Let $f\colon U \to P$ be the canonical quotient map defining the regular directed complex $P$ in Example \ref{exm:nonsdw}.
	Then $U$ is a molecule, $P$ is dimension-wise acyclic, and $f$ is a local embedding, but evidently $f$ is not injective.
\end{exm}

\begin{cor} \label{cor:strongly_dimensionwise_acyclic_omega_cat}
Let $P$ be a strongly dimension-wise acyclic regular directed complex.
Then
\[
	\molecin{P} = \set{ \isocl{ U \incl P } \mid U \subseteq P, \text{$U$ is a molecule} }. 
\]
\end{cor}

\begin{lem} \label{lem:acyclic_closed_under_pasting}
	Let $U, V$ be molecules and $k < \min \set{\dim{U}, \dim{V}}$ such that $U \cp{k} V$ is defined.
	If $U$ and $V$ are acyclic, then $U \cp{k} V$ is acyclic.
\end{lem}
\begin{proof}
	We will identify $U$ and $V$ with their isomorphic images in $U \cp{k} V$.
	Consider any path in $\hasseo{(U \cp{k} V)}$.
	If the path passes through an edge $x \to y$ such that $y \in U \setminus V$ but $x \in V$, since $V$ is closed, necessarily $x \in \faces{}{-}y$, so $x \in U \cap V$.
	Similarly, if $x \in U \setminus V$ and $y \in V$, then $y \in \faces{}{+}x$, so $y \in U \cap V$.
	It follows that any path has to pass through $U \cap V$ whenever it enters and leaves $U \setminus V$, and dually whenever it enters and leaves $V \setminus U$.

	Because $U$ and $V$ are acyclic, any cycle in $\hasseo{(U \cp{k} V)}$ has to pass through elements both in $V \setminus U$ and $U \setminus V$.
	We may then represent the cycle as
	\[
		x_\bullet \equiv (x_0 \to x_1 \to \ldots \to x_m = x_0)
	\]
	in such a way that $x_0 \in U \cap V$.
	Given such a representation, we let $S(x_\bullet)$ be the set
	\[
		\set{(i, j) \mid 
		\text{$i < j$, $x_i, x_j \in U \cap V$, $x_k \in (U \setminus V) \cup (V \setminus U)$ for all $i < k < j$} }
	\]
	of maximal segments of $x_\bullet$ which are entirely included in $U \setminus V$ or $V \setminus U$, marked with their endpoints in $U \cap V$.

	Suppose $x_\bullet$ is a cycle such that $\size{S(x_\bullet)}$ is minimal.
	Necessarily, $\size{S(x_\bullet)} > 1$, for otherwise the cycle is entirely contained in $U$ or in $V$.
	Let $(i, j) \in S(x_\bullet)$, and suppose without loss of generality that the segment $x_i \to \ldots \to x_j$ is contained in $U$.
	Then $x_i, x_j \in \bound{k}{+}U$, which is an acyclic molecule.
	By Proposition \ref{prop:hasseo_connected_frame_acyclic}, there exists a path from $x_i$ to $x_j$ or a path from $x_j$ to $x_i$ in $\hasseo{(\bound{k}{+}U)}$.
	In the latter case, there would be a cycle in $\hasseo{U}$, a contradiction.
	In the first case, we can replace the segment $x_i \to \ldots \to x_j$ of $x_\bullet$ with a path entirely contained in $U \cap V$.
	This produces a new cycle $x'_\bullet$ with $\size{S(x'_\bullet)} = \size{S(x_\bullet)} - 1$, contradicting minimality of $x_\bullet$.
\end{proof}

\begin{exm}[A pasting of strongly dimension-wise acyclic molecules which is not dimension-wise acyclic] \index[counterex]{A pasting of strongly dimension-wise acyclic molecules which is not dimension-wise acyclic}
	We prove that Lemma \ref{lem:acyclic_closed_under_pasting} does not extend to weaker acyclicity conditions.
	Let $V$ be a 3\nbd dimensional atom whose input and output boundaries are the oriented face posets of
\[\begin{tikzcd}[sep=small]
	&& \bullet \\
	\bullet &&& \bullet \\
	& \bullet
	\arrow[""{name=0, anchor=center, inner sep=0}, curve={height=12pt}, from=3-2, to=2-4]
	\arrow[curve={height=6pt}, from=2-1, to=3-2]
	\arrow[curve={height=-6pt}, from=1-3, to=2-4]
	\arrow[""{name=1, anchor=center, inner sep=0}, curve={height=-12pt}, from=2-1, to=1-3]
	\arrow[from=3-2, to=1-3]
	\arrow[curve={height=6pt}, shorten <=7pt, Rightarrow, from=0, to=1-3]
	\arrow[curve={height=-6pt}, shorten >=7pt, Rightarrow, from=3-2, to=1]
\end{tikzcd}
\quad \text{and} \quad 
\begin{tikzcd}[sep=small]
	&& \bullet \\
	\bullet &&& \bullet \\
	& \bullet
	\arrow[curve={height=12pt}, from=3-2, to=2-4]
	\arrow[curve={height=6pt}, from=2-1, to=3-2]
	\arrow[curve={height=-6pt}, from=1-3, to=2-4]
	\arrow[curve={height=-12pt}, from=2-1, to=1-3]
	\arrow[curve={height=-6pt}, Rightarrow, from=3-2, to=1-3]
\end{tikzcd}
\]
	respectively, and let $U$ be the 3\nbd dimensional atom from Example 
	\ref{exm:non_acyclic}.
	Then both $V$ and $U$ are strongly dimension-wise acyclic.
	However, the boundary of the pasting $V \cp{2} U$ is isomorphic to the boundary of the 3\nbd dimensional atom from Example \ref{exm:non_dw_acy}, which contains a cycle in its 0\nbd flow graph.
	We conclude that $V \cp{2} U$ is not dimension-wise acyclic.
\end{exm}

\begin{prop} \label{prop:gray_product_of_acyclic}
	Let $P$, $Q$ be acyclic oriented graded posets.
	Then $P \gray Q$ is acyclic.
\end{prop}
\begin{proof}
We will prove the contrapositive.
Suppose that $\hasseo{(P \gray Q)}$ has a cycle
\[
	(x_0, y_0) \to (x_1, y_1) \to \ldots \to (x_{m-1}, y_{m-1}) \to (x_m, y_m) = (x_0, y_0).
\]
Let 
\[
	p \eqdef \size{\set{i \in \set{1, \ldots, m} \mid x_i \neq x_{i-1}}},
\]
and define recursively
\[
	\fun{j}(i) \eqdef \begin{cases}
		0 & \text{for $i = 0$}, \\
		\min \set{j > \fun{j}(i-1) \mid x_j \neq x_{\fun{j}(i-1)}} & \text{for $i \in \set{1, \ldots, p}$}.
	\end{cases}
\]
Suppose that $p > 0$.
By definition of the orientation in $P \gray Q$, for all $i \in \set{1, \ldots, p}$, since there is an edge
\[
	(x_{\fun{j}(i) - 1}, y_{\fun{j}(i) - 1}) \to (x_{\fun{j}(i)}, y_{\fun{j}(i)})
\]
in $\hasseo{(P \gray Q)}$ and $x_{\fun{j}(i)} \neq x_{\fun{j}(i) - 1} = x_{\fun{j}(i-1)}$, necessarily $y \eqdef y_{\fun{j}(i) - 1} = y_{\fun{j}(i)}$ and
\[
	x_{\fun{j}(i)} \in \faces{}{+}x_{\fun{j}(i-1)} \quad \text{or} \quad x_{\fun{j}(i-1)} \in \faces{}{-}x_{\fun{j}(i)}.
\]
It follows that
\[
	x_{\fun{j}(0)} \to x_{\fun{j}(1)} \to \ldots \to x_{\fun{j(p)}} = x_{\fun{j}(0)}
\]
is a cycle in $\hasseo{P}$.
Finally, suppose that $p = 0$.
Then $x_i = x_0$ for all $i \in \set{1, \ldots, m}$.
It follows that, if $\dim{x_0}$ is even, then
\[
	y_0 \to y_1 \to \ldots \to y_{m-1} \to y_m = y_0
\]
is a cycle in $\hasseo{Q}$, and if $\dim{x_0}$ is odd, then
\[
	y_m \to y_{m-1} \to \ldots \to y_1 \to y_0 = y_m
\]
is a cycle in $\hasseo{Q}$.
\end{proof}

\begin{exm}[A Gray product of strongly dimension-wise acyclic molecules which is not dimension-wise acyclic] \index[counterex]{A Gray product of strongly dimension-wise acyclic molecules which is not dimension-wise acyclic} \label{exm:dw_acyclic_not_gray_stable}
	We prove that Proposition \ref{prop:gray_product_of_acyclic} does not extend to weaker acyclicity conditions.
	Let $U$ be a 3\nbd dimensional atom whose input and output boundary are the oriented face posets of
\[\begin{tikzcd}[column sep=small]
	{{\scriptstyle 0}\; \bullet} &&&& {{\scriptstyle 2}\; \bullet} \\
	&& {{\scriptstyle 1}\; \bullet}
	\arrow["0"', curve={height=12pt}, from=1-1, to=2-3]
	\arrow[""{name=0, anchor=center, inner sep=0}, "3", curve={height=-12pt}, from=1-1, to=1-5]
	\arrow[""{name=1, anchor=center, inner sep=0}, "1"', curve={height=12pt}, from=2-3, to=1-5]
	\arrow[""{name=2, anchor=center, inner sep=0}, "2", curve={height=-12pt}, from=2-3, to=1-5]
	\arrow["0", curve={height=-6pt}, shorten >=5pt, Rightarrow, from=2-3, to=0]
	\arrow["1"', shorten <=3pt, shorten >=3pt, Rightarrow, from=1, to=2]
\end{tikzcd}
\quad \text{and} \quad
	\begin{tikzcd}[column sep=small]
	{{\scriptstyle 0}\; \bullet} &&&& {{\scriptstyle 2}\; \bullet} \\
	&& {{\scriptstyle 1}\; \bullet}
	\arrow[""{name=0, anchor=center, inner sep=0}, "0"', curve={height=12pt}, from=1-1, to=2-3]
	\arrow[""{name=1, anchor=center, inner sep=0}, "3", curve={height=-12pt}, from=1-1, to=1-5]
	\arrow[""{name=2, anchor=center, inner sep=0}, "4", curve={height=-12pt}, from=1-1, to=2-3]
	\arrow["1"', curve={height=12pt}, from=2-3, to=1-5]
	\arrow["2", shorten <=3pt, shorten >=3pt, Rightarrow, from=0, to=2]
	\arrow["3"', curve={height=6pt}, shorten >=5pt, Rightarrow, from=2-3, to=1]
\end{tikzcd}\]
	respectively.
	Then $U$ is strongly dimension-wise acyclic.
	However, in $U \gray U$, writing $x \gray y$ instead of $(x, y)$ for better readability, we have
	\begin{align*}
		(0, 1) \gray (2, 2) & \in
		\faces{}{+}((0, 1) \gray (3, 0)) \cap \faces{}{-}((1, 1) \gray (2, 2)), \\
		(1, 1) \gray (1, 0) & \in
		\faces{}{+}((1, 1) \gray (2, 2)) \cap \faces{}{-}((2, 1) \gray (1, 0)), \\
		(2, 1) \gray (0, 1) & \in
		\faces{}{+}((2, 1) \gray (1, 0)) \cap \faces{}{-}((3, 0) \gray (0, 1)), \\
		(2, 2) \gray (0, 1) & \in
		\faces{}{+}((3, 0) \gray (0, 1)) \cap \faces{}{-}((2, 2) \gray (1, 2)), \\
		(1, 4) \gray (1, 2) & \in
		\faces{}{+}((2, 2) \gray (1, 2)) \cap \faces{}{-}((1, 4) \gray (2, 1)), \\
		(0, 1) \gray (2, 1) & \in
		\faces{}{+}((1, 4) \gray (2, 1)) \cap \faces{}{-}((0, 1) \gray (3, 0)).
	\end{align*}
	These relations determine a cycle in $\flow{2}{(U \gray U)}$.
	This proves that $U \gray U$ is not dimension-wise acyclic.
\end{exm}

\begin{lem} \label{lem:flow_graph_of_suspension}
	Let $P$ be an oriented graded poset, $k \in \mathbb{N}$.
	Then
	\begin{enumerate}
		\item $x \mapsto \sus{x}$ induces an isomorphism of directed graphs $\flow{k}{P} \iso \flow{k+1}{\sus{P}}$,
		\item $x \mapsto \sus{x}$ induces an embedding of directed graphs $\extflow{k}{P} \incl \extflow{k+1}{\sus{P}}$, whose complement is the discrete graph on $\{\bot^-, \bot^+\}$.
	\end{enumerate}
\end{lem}
\begin{proof}
	Follows from Lemma \ref{lem:faces_of_suspension} and Corollary 
	\ref{cor:boundary_of_suspension}, together with the observation that $\bound{0}{}\sus{x} = \set{\bot^+, \bot^-}$ for all $x \in P$, so $\bot^+$ and $\bot^-$ can never appear in the interior of $\bound{k+1}{+}x$ or $\bound{k+1}{-}x$.
\end{proof}

\begin{prop} \label{prop:suspension_of_acyclic}
	Let $P$ be an oriented graded poset.
	Then
	\begin{enumerate}
		\item if $P$ is acyclic, then $\sus{P}$ is acyclic,
		\item if $P$ is strongly dimension-wise acyclic, then $\sus{P}$ is strongly dimension-wise acyclic,
		\item if $P$ is dimension-wise acyclic, then $\sus{P}$ is dimension-wise acyclic.
	\end{enumerate}
\end{prop}
\begin{proof}
The strongly dimension-wise acyclic and dimension-wise acyclic cases follow from Lemma 
\ref{lem:flow_graph_of_suspension}, together with the observation that $\flow{0}{\sus{P}}$ is always a discrete graph, and that $\extflow{0}{\sus{P}}$ has 
\begin{itemize}
	\item an edge from $\bot^-$ to every element of the form $\sus{x}$,
	\item an edge from every element of the form $\sus{x}$ to $\bot^+$,
\end{itemize}
and no other edges.
For the acyclic case, let
\[
	x_0 \to x_1 \to \ldots \to x_{m-1} \to x_m
\]
be a path in $\hasseo{\sus{P}}$.
Suppose that $x_i = \bot^-$ for some $i \in \set{1, \ldots, m}$.
Then necessarily $i = 0$, because $\cofaces{}{+}\bot^- = \varnothing$.
Dually, if $x_i = \bot^+$, then $i = m$.
Thus any path in $\hasseo{\sus{P}}$ is of the form
\[
	(\bot^- \to) \; \sus{x'_1} \to \ldots \to \sus{x'_{m-1}} \; (\to \bot^+)
\]
where $x'_1 \to \ldots \to x'_{m-1}$ is a path in $\hasseo{P}$, and only contains a cycle if the latter does.
\end{proof}

\begin{prop} \label{prop:join_of_acyclic}
	Let $P$, $Q$ be acyclic oriented graded posets.
	Then $P \join Q$ is acyclic.
\end{prop}
\begin{proof}
	Let $s\colon P \join Q \to \sus{P} \gray \sus{Q}$ be the injection of Lemma 
	\ref{lem:joins_and_gray_of_suspension}.
	Then every path
	\[
		z_0 \to z_1 \to \ldots \to z_{m-1} \to z_m
	\]
	in $\hasseo{(P \join Q)}$ induces a path
	\[
		s(z_0) \to s(z_1) \to \ldots \to s(z_{m-1}) \to s(z_m)
	\]
	in $\hasseo{(\sus{P} \gray \sus{Q})}$.
	Since $\sus{P} \gray \sus{Q}$ is acyclic by Proposition \ref{prop:gray_product_of_acyclic} and Proposition \ref{prop:suspension_of_acyclic}, it follows that $P \join Q$ is acyclic.
\end{proof}

\begin{dfn}[The categories $\rdcpxmapac$ and $\rdcpxcomapac$] \index{$\rdcpxmapac$} \index{$\rdcpxcomapac$}
	We let $\rdcpxmapac$ and $\rdcpxcomapac$ denote the full subcategories of $\rdcpxmap$ and $\rdcpxcomap$, respectively, on the acyclic regular directed complexes.
\end{dfn}

\begin{cor} \label{cor:monoidal_structures_on_acyclic}
	The monoidal structures
	\[
		(\rdcpxmap, \gray, 1), \quad \quad (\rdcpxmap, \join, \varnothing)
	\]
	restrict to monoidal structures on $\rdcpxmapac$, and the monoidal structures
	\[
		(\rdcpxcomap, \gray, 1) \quad \quad (\rdcpxcomap, \join, \varnothing)
	\]
	restrict to monoidal structures on $\rdcpxcomapac$.
\end{cor}

\begin{dfn}[Converse of a directed graph] \index{directed graph!converse} \index{$\optot{\mathscr{G}}$}
	Let $\mathscr{G}$ be a directed graph.
	The \emph{converse of $\mathscr{G}$} is the directed graph $\optot{\mathscr{G}}$ with
	\begin{itemize}
		\item the same sets of vertices and edges as $\mathscr{G}$,
		\item source and target functions swapped with respect to $\mathscr{G}$.
	\end{itemize}
\end{dfn}

\begin{lem} \label{lem:graph_acyclic_iff_converse_acyclic}
	A directed graph is acyclic if and only if its converse is acyclic.
\end{lem}
\begin{proof}
	Let $\mathscr{G}$ be a directed graph.
	Then paths $x_0 \to \ldots \to x_m$ in $\mathscr{G}$ correspond biunivocally to paths $x_m \to \ldots \to x_0$ in $\optot{\mathscr{G}}$.
\end{proof}

\begin{lem} \label{lem:flow_graphs_under_dual}
	Let $P$ be an oriented graded poset, $J \subseteq \posnat$, $k \geq -1$, and consider the bijection $\dual{J}{}\colon x \mapsto \dual{J}{x}$ between the underlying sets of $P$ and $\dual{J}{P}$.
	Then
	\begin{enumerate}
		\item if $k + 1 \in J$, then $\dual{J}{}$ induces isomorphisms of directed graphs
			\[
				\optot{(\maxflow{k}{P})} \iso \maxflow{k}{\dual{J}{P}}, \quad \optot{(\flow{k}{P})} \iso \flow{k}{\dual{J}{P}}, \quad \optot{(\extflow{k}{P})} \iso \extflow{k}{\dual{J}{P}},
			\]
		\item if $k + 1 \not\in J$, then $\dual{J}{}$ induces isomorphisms of directed graphs
			\[
				\maxflow{k}{P} \iso \maxflow{k}{\dual{J}{P}}, \quad \flow{k}{P} \iso \flow{k}{\dual{J}{P}}, \quad \quad \extflow{k}{P} \iso \extflow{k}{\dual{J}{P}}.
			\]
	\end{enumerate}
\end{lem}
\begin{proof}
	Follows immediately from Lemma \ref{lem:faces_of_globes} and Corollary \ref{cor:boundaries_of_dual} by the definitions of $k$\nbd flow, maximal $k$\nbd flow, and extended $k$\nbd flow graph.
\end{proof}

\begin{prop} \label{prop:dw_acyclic_stable_under_all_duals}
	Let $P$ be an oriented graded poset, $J \subseteq \posnat$.
	Then
	\begin{enumerate}
		\item if $P$ is frame-acyclic, then so is $\dual{J}{P}$,
		\item if $P$ is dimension-wise acyclic, then so is $\dual{J}{P}$,
		\item if $P$ is strongly dimension-wise acyclic, then so is $\dual{J}{P}$.
	\end{enumerate}
\end{prop}
\begin{proof}
	Follows from Lemma \ref{lem:graph_acyclic_iff_converse_acyclic} combined with Lemma 
	\ref{lem:flow_graphs_under_dual}.
\end{proof}

\begin{exm}[A join of strongly dimension-wise acyclic molecules which is not dimension-wise acyclic] \index[counterex]{A join of strongly dimension-wise acyclic molecules which is not dimension-wise acyclic} \label{exm:dw_acyclic_not_join_stable}
	We show that Proposition \ref{prop:join_of_acyclic} does not extend to weaker acyclicity conditions.
	
	Let $U$ be the same 3\nbd dimensional atom as in Example 
	\ref{exm:dw_acyclic_not_gray_stable}.
	Since $U$ is strongly dimension-wise acyclic, by Proposition 
	\ref{prop:dw_acyclic_stable_under_all_duals} so is its total dual $\optot{U}$.
	Using the isomorphism between $\augm{(U \join \optot{U})}$ and $\augm{U} \gray \augm{(\optot{U})}$, since the total dual counteracts the orientation reversal on faces of the second factor due to dimensions being raised by 1, we see that the cycle in $\flow{2}{(U \gray U)}$ maps to a cycle
	\begin{align*}
		(0, 1) \join \optot{(3, 0)} & \to (1, 1) \join \optot{(2, 2)} \to (2, 1) \join \optot{(1, 0)} \to (3, 0) \join \optot{(0, 1)} \to \\
					    & \to (2, 2) \join \optot{(1, 2)} \to (1, 4) \join \optot{(2, 1)} \to (0, 1) \join \optot{(3, 0)}
	\end{align*}
	in $\flow{3}{(U \join \optot{U})}$.
	This proves that $U \join \optot{U}$ is not dimension-wise acyclic.
\end{exm}

\begin{lem} \label{lem:total_dual_hasseo}
	Let $P$ be an oriented graded poset.
	Then the bijection $x \mapsto \optot{x}$ induces an isomorphism $\optot{(\hasseo{P})} \iso \hasseo{(\optot{P})}$ of directed graphs.
\end{lem}
\begin{proof}
	For all $x, y \in P$ and $\alpha \in \set{+, -}$, we have $\optot{x} \in \faces{}{\alpha}\optot{y}$ if and only if $x \in \faces{}{-\alpha}y$.
	The claim then follows by the definition of oriented Hasse diagram.
\end{proof}

\begin{prop} \label{prop:acyclic_stable_under_total_dual}
	Let $P$ be an acyclic oriented graded poset.
	Then $\optot{P}$ is acyclic.
\end{prop}
\begin{proof}
	Follows from Lemma \ref{lem:graph_acyclic_iff_converse_acyclic} combined with Lemma
	\ref{lem:total_dual_hasseo}.
\end{proof}

\begin{exm}[An acyclic molecule with a dual which is not acyclic] \index[counterex]{An acyclic molecule with a dual which is not acyclic}
	Let $U$ be the 3\nbd dimensional atom of Example \ref{exm:non_acyclic}.
	Then $U$ is not acyclic, but $\opp{U}$ is acyclic.
	Since $U$ is isomorphic to $\opp{(\opp{U})}$, we conclude that acyclicity is not in general stable under duals, except for the total dual.
\end{exm}


\section{In low dimensions} \label{sec:in_low_dim}

\begin{guide}
	In this section, we study molecules and regular directed complexes up to dimension 3, which covers most practical uses of diagrammatic reasoning.
	We show that all molecules up to dimension 2 are acyclic (Proposition 
	\ref{prop:dim2_acyclic}), and in fact admit two ``orthogonal'' partial orders on their elements, a \cemph{horizontal order} and \cemph{vertical order}, which jointly form a linear order.
	This reflects the intuition that, given a pair of cells in a 2\nbd dimensional pasting diagram shape, the first is either \emph{below}, or \emph{above}, or \emph{to the left}, or \emph{to the right} of the other, and implies that recognition of rewritable submolecules is trivial in dimension 2 (Theorem 
	\ref{thm:round_submolecule_dim2}).

	Then, we prove that molecules up to dimension 3 are frame-acyclic (Theorem 
	\ref{thm:dim3_frame_acyclic}), which implies that every oriented graded poset up to dimension 3 has frame-acyclic molecules and presents a polygraph (Corollary \ref{cor:dim3_has_frame_acyclic_molecules}).
	We use this to formulate a simple criterion for the rewritable submolecule problem in dimension 3 (Theorem \ref{thm:round_submolecule_dim3}).
\end{guide}

\begin{lem} \label{lem:1_molecule_unique_ordering}
	Let $U$ be a molecule, $\dim{U} \leq 1$.
	Then
	\begin{enumerate}
	    \item $U$ is round,
	    \item $U$ is acyclic,
	    \item $\hasseo{U}$ is a linear graph with $\size{U}$ vertices,
	    \item $\flow{0}{U}$ is a linear graph with $\size{\grade{1}{U}}$ vertices,
	    \item $U$ admits a unique 0\nbd ordering.
	\end{enumerate}
\end{lem}
\begin{proof}
All straightforward checks using Lemma \ref{lem:only_0_molecule} and Lemma \ref{lem:only_1_molecules}.
\end{proof}

\begin{prop} \label{prop:in_dim1_all_inclusions_are_submolecule}
	Let $\imath\colon V \incl U$ be an inclusion of $1$\nbd dimensional molecules.
	Then $\imath$ is a submolecule inclusion.
\end{prop}
\begin{proof}
By Lemma \ref{lem:flow_under_inclusion} $\flow{0}{V}$ is an induced subgraph of $\flow{0}{U}$.
By Lemma  \ref{lem:1_molecule_unique_ordering} both of them are linear graphs, and an induced subgraph of a linear graph is a linear graph if and only if its vertices are consecutive in the ambient graph.
All other conditions of Lemma \ref{lem:round_submolecules_from_layering} are trivially satisfied.
\end{proof}

\begin{prop} \label{prop:dim2_acyclic}
	Let $U$ be a molecule, $\dim{U} \leq 2$.
	Then $U$ is acyclic.
\end{prop}
\begin{proof}
If $\dim{U} < 2$ this is part of Lemma \ref{lem:1_molecule_unique_ordering}, and if $U$ is a 2\nbd dimensional atom it can be checked directly using Lemma \ref{lem:only_2_atoms}.
The statement then follows by an easy induction from Lemma \ref{lem:acyclic_closed_under_pasting}.
\end{proof}

\begin{dfn}[Horizontal and vertical order] \index{order!horizontal} \index{order!vertical}
	Let $U$ be a molecule, $\dim{U} \leq 2$.
	The \emph{horizontal order} $\prech$ and the \emph{vertical order} $\precv$ on $U$ are defined by
\begin{itemize}
	\item $x \prech y$ if and only if there is a path from $x$ to $y$ in $\extflow{0}{U}$,
	\item $x \precv y$ if and only if there is a path from $x$ to $y$ in $\extflow{1}{U}$.
\end{itemize}
\end{dfn}

\begin{rmk}
It follows from Proposition \ref{prop:dim2_acyclic} together with Proposition \ref{prop:acyclicity_implications} that $\prech$ and $\precv$ are partial orders.
\end{rmk}

\begin{lem} \label{lem:dim2_horizontal_or_vertical_path}
	Let $U$ be a molecule, $\dim{U} \leq 2$.
	Then
	\begin{enumerate}
	    \item the union of $\prech$ and $\precv$ is a linear order on $\grade{1}{U}$,
	    \item the intersection of $\prech$ and $\precv$ is the identity relation on $\grade{1}{U}$.
	\end{enumerate}
\end{lem}
\begin{proof}
	By Proposition \ref{prop:dim2_acyclic} together with Proposition \ref{prop:acyclicity_implications} and Proposition \ref{prop:dimensionwise_implies_frame_acyclic}, $U$ is acyclic, strongly dimension-wise acyclic and frame-acyclic.

	Let $x, y \in U$, and suppose $x \neq y$.
	By Proposition \ref{prop:extflow_connected_frame_acyclic}, there exists $k \geq -1$ such that there is a path from $x$ to $y$ or a path from $y$ to $x$ in $\extflow{k}{U}$.
	Necessarily $k \in \set{0, 1}$, because otherwise $\extflow{k}{U}$ is a discrete graph.
	Suppose, without loss of generality, that $k = 0$ and that there is a path from $x$ to $y$.
	Then $x \prech y$.
	Because $\extflow{0}{U}$ is acyclic, it cannot be the case that $y \prech x$.
	It cannot be the case that $y \precv x$ either, for otherwise, by Lemma \ref{lem:path_in_extflow_induces_path_in_hasseo}, the paths exhibiting $x \prech y$ and $y \precv x$ would induce a cycle in $\hasseo{U}$.
	This proves that the union of $\prech$ and $\precv$ is a linear order.

	Suppose that $x \prech y$ and $x \precv y$; we will prove that $x = y$.
	If $\dim{U} < 2$, then $x \precv y$ immediately implies $x = y$, so suppose $\dim{U} = 2$.
	Then $U$ admits a 1\nbd layering with an associated 1\nbd ordering $(\order{i}{x})_{i=1}^m$.
	Let 
	\[
		\order{0}{V} \eqdef \bound{}{-}U, \quad \quad \order{i}{V} \eqdef \order{i-1}{V} \cup \clset{\order{i}{x}}
	\]
	for $i \in \set{1, \ldots, m}$, the images of increasing initial segments of the 1\nbd layering.
	We let $\order{i}{\prech}$ and $\order{i}{\precv}$ be the horizontal and vertical orders on $\order{i}{V}$, which are subsets of those on $U$, and proceed by induction on $i \in \set{0, \ldots, m}$.
	Since $\dim{\order{0}{V}} = 1$, we have already proved the base case.

	Let $i > 0$ and suppose $x \order{i}{\prech} y$ and $x \order{i}{\precv} y$.
	If $y \in \order{i-1}{V}$, then any path to $y$ in $\extflow{1}{\order{i}{V}}$ has to be contained in $\order{i-1}{V}$ because $\bound{}{+}\order{i}{x} \subseteq \bound{}{+}\order{i}{V}$.
	Any path to $y$ in $\extflow{0}{\order{i}{V}}$ can, at worst, contain segments which enter $\clset{\order{i}{x}}$ from $\bound{0}{-}\order{i}{x}$ and leave it from $\bound{0}{+}\order{i}{x}$.
	Any such segment can be replaced with a path traversing $\bound{}{-}\order{i}{x}$, producing a path contained in $\order{i-1}{V}$.
	It follows that $x \order{i-1}{\prech} y$ and $x \order{i-1}{\precv} y$, so $x = y$ by the inductive hypothesis.

	Suppose that $y \notin \order{i-1}{V}$, that is, $y \in \set{\order{i}{x}} \cup \inter{\bound{}{+}\order{i}{x}}$.
	Then any non-trivial path to $y$ in $\extflow{1}{\order{i}{V}}$ consists of a path contained in $\order{i-1}{V}$ to some $z \in \inter{\bound{}{-}\order{i}{x}}$, followed by the path $z \to \order{i}{x}$, and possibly a single step $\order{i}{x} \to y$.
	Then either $x = y$, or $x = \order{i}{x}$, or $x \in \order{i-1}{V}$ and $x \order{i-1}{\precv} z$.
	In the latter case, a path from $x$ to $y$ in $\extflow{0}{\order{i}{V}}$ has to enter $\clset{\order{i}{x}}$ from $\bound{0}{-}\order{i}{x}$, and then cannot leave it, for otherwise it would have to re-enter from $\bound{0}{-}\order{i}{x}$ creating a cycle.
	It follows that such a path consists of a path contained in $\order{i-1}{V}$ from $x$ to $\bound{0}{-}\order{i}{x}$, followed by a segment contained in $\bound{}{+}\order{i}{x}$.
	Since there is a path from $\bound{0}{-}\order{i}{x}$ to any element of $\bound{}{-}\order{i}{x}$, it follows that $x \order{i-1}{\prech} z$, so by the inductive hypothesis $x = z$.
	But then the path from $x$ to $z$ through $\bound{0}{-}\order{i}{x}$ is a cycle, a contradiction.
	The case $x = \order{i}{x}$ is dealt with in a similar way.
	Since $U = \order{m}{V}$, we conclude.
\end{proof}

\begin{lem} \label{lem:dim2_inclusion_is_path_induced}
	Let $\imath\colon V \incl U$ be an inclusion of molecules, $\dim{U} \leq 2$.
	Then $\flow{1}{V}$ is a path-induced subgraph of $\flow{1}{U}$.
\end{lem}
\begin{proof}
	Both $U$ and $\imath(V)$ are molecules of dimension $\leq 2$.
	Let $\prech$, $\precv$ be the horizontal and vertical order on $U$, and $\prech^V$, $\precv^V$ those on $\imath(V)$, which are subsets of those on $U$.

	Suppose by way of contradiction that $\flow{1}{V}$ is not path-induced. 
	Then there exists a path $x_0 \to \ldots \to x_m$ in $\flow{1}{U}$ such that $m > 1$, $x_0, x_m \in \imath(V)$, and $x_i \notin \imath(V)$ for all $i \in \set{1, \ldots, m-1}$.
	By definition, there exist 1\nbd dimensional elements $y_i \in \faces{}{+}x_{i-1} \cap \faces{}{-}x_{i}$ for all $i \in \set{1, \ldots, m}$.
	Then $y_1 \precv y_m$ and $y_1 \neq y_m$.
	By Corollary \ref{cor:codimension_1_elements}, $x_{i-1}$ and $x_i$ are the only cofaces of $y_i$ for all $i \in \set{1, \ldots, m}$.
	Necessarily, then, $y_1 \in \faces{}{+}\imath(V)$ and $y_m \in \faces{}{-}\imath(V)$, 
so it is not possible that $y_1 \precv^V y_m$.
	Then by Lemma \ref{lem:dim2_horizontal_or_vertical_path} applied to $\imath(V)$, one of $y_1 \prech^V y_m$, $y_m \prech^V y_1$, or $y_m \precv^V y_1$ must hold.
	Combined with $y_1 \precv y_m$, each of these implies $y_1 = y_m$ by Lemma \ref{lem:dim2_horizontal_or_vertical_path} applied to $U$, a contradiction.
\end{proof}

\begin{thm} \label{thm:round_submolecule_dim2}
	Let $\imath\colon V \incl U$ be an inclusion of molecules such that $\dim{U} = \dim{V} = 2$ and $V$ is round.
	Then $\imath$ is a submolecule inclusion.
\end{thm}
\begin{proof}
	By Lemma \ref{lem:dim2_inclusion_is_path_induced} combined with Lemma \ref{lem:connected_subgraph_conditions_path_induced}, there exists a 1\nbd ordering $(\order{i}{x})_{i=1}^m$ of $U$ in which the elements of $\imath(V)$ are consecutive, that is, $\order{i}{x} \in \imath(V)$ if and only if $p \leq i \leq q$ for some $p, q \in \set{1, \ldots, m}$.

	By Proposition \ref{prop:dim2_acyclic} $U$ is acyclic, so by Proposition \ref{prop:acyclicity_implications} and Proposition \ref{prop:dimensionwise_implies_frame_acyclic} it is frame-acyclic, and by Corollary \ref{cor:frame_acyclicity_equivalent_conditions} the 1\nbd ordering comes from a 1\nbd layering $(\order{i}{U})_{i=1}^m$ such that $\imath(\bound{}{-}V) \subseteq \bound{}{-}\order{p}{U}$.
	Since both are 1\nbd dimensional molecules, by Proposition \ref{prop:in_dim1_all_inclusions_are_submolecule} $\imath(\bound{}{-}V) \submol \bound{}{-}\order{p}{U}$, and we conclude by Lemma \ref{lem:round_submolecules_from_layering}.
\end{proof}

\begin{lem} \label{lem:dim2_multiple_substitutions}
	Let $U$ be a molecule, $\dim{U} \leq 2$, and suppose that 
	\[
		(\order{i}{V}, \order{i}{W}, j_i\colon \order{i}{V} \incl U)_{i=1}^m
	\]
	is a family of triples such that the multiple substitution $\subs{U}{\order{i}{W}}{j_i(\order{i}{V})}_{i=1}^m$ is defined.
	Then $\subs{U}{\order{i}{W}}{j_i(\order{i}{V})}_{i=1}^m$ is a molecule.
\end{lem}
\begin{proof}
	By assumption $\subs{U}{\order{i}{W}}{j_i(\order{i}{V})}_{i=1}^0 = U$ is a molecule.
	For $k > 0$, supposing that $\subs{U}{\order{i}{W}}{j_i(\order{i}{V})}_{i=1}^{k-1}$ is a molecule, 
\begin{equation*}
    j_k\colon \order{k}{V} \incl \subs{U}{\order{i}{W}}{j_i(\order{i}{V})}_{i=1}^{k-1}
\end{equation*}
	is the inclusion of a round molecule into a molecule of the same dimension $\leq 2$.
	By Proposition \ref{prop:in_dim1_all_inclusions_are_submolecule} and Theorem \ref{thm:round_submolecule_dim2}, it is a submolecule inclusion.
	We conclude by Proposition \ref{prop:round_submolecule_substitution}.
\end{proof}

\begin{comm}
	The admissibility of multiple substitution on molecules is a combinatorial analogue of the topological \emph{domain replacement condition} from \cite[Definition 3.9]{power1991pasting}.
	Thus Lemma \ref{lem:dim2_multiple_substitutions} implies that the domain replacement condition holds automatically for 2\nbd dimensional molecules, something that was also observed by Power in relation to 2\nbd dimensional pasting schemes.
\end{comm}

\begin{thm} \label{thm:dim3_frame_acyclic}
	Let $U$ be a molecule, $\dim{U} \leq 3$.
	Then $U$ is frame-acyclic.
\end{thm}
\begin{proof}
	Since $V \submol U$ implies $\dim{V} \leq \dim{U}$, it suffices to show that, if $r \eqdef \frdim{U}$, then $\maxflow{r}{U}$ is acyclic.
	If $\dim{U} < 3$, this is implied by Proposition \ref{prop:dim2_acyclic} together with Proposition \ref{prop:acyclicity_implications} and Proposition \ref{prop:dimensionwise_implies_frame_acyclic}.
	If $\dim{U} = 3$ and $r = 2$, this is implied by Corollary \ref{cor:flow_acyclic_in_codimension_1} together with Remark \ref{rmk:codim_1_maximal_flow_graph}.

	It suffices then to consider the case $\dim{U} = 3$ and $r < 2$.
	Let $(\order{i}{x})_{i=1}^m$ be any 2\nbd ordering of $U$, which in this case is simply an enumeration of its 3\nbd dimensional elements.
	By Lemma \ref{lem:frdim_boundaries}, for all $i \in \set{1, \ldots, m}$ we have
\[
	\bound{}{-}\order{i}{x} \subseteq \bound{}{-}U,
\]
	and since $\bound{}{-}\order{i}{x}$ is a round 2\nbd dimensional molecule and $\bound{}{-}U$ is a 2\nbd dimensional molecule, by Theorem \ref{thm:round_submolecule_dim2} we have in fact $\bound{}{-}\order{i}{x} \submol \bound{}{-}U$.

	By Proposition \ref{prop:intersection_of_maximal_elements}, for all $i \neq j \in \set{1, \ldots, m}$,
\[
    \bound{}{-}\order{i}{x} \cap \bound{}{-}\order{j}{x} \subseteq \clset{{\order{i}{x}}} \cap \clset{{\order{j}{x}}} \subseteq \bound{r}{}\order{i}{x} \cap \bound{r}{}\order{j}{x}.
\]
	It follows from Lemma \ref{lem:dim2_multiple_substitutions} that the multiple substitution
\[
    \tilde{U} \eqdef \subs{\bound{}{-}U}{\compos{\bound{}{-}\order{i}{x}}}{\bound{}{-}\order{i}{x}}_{i=1}^m
\]
	is defined, and is a 2\nbd dimensional molecule with the same frame dimension as $U$.
	We claim that $\maxflow{r}{\tilde{U}}$ is isomorphic to $\maxflow{r}{U}$:
	\begin{itemize}
		\item every maximal element of dimension $< 3$ in $U$ is maximal in $\bound{}{-}U$ and not affected by the multiple substitution, so it appears unchanged in $\tilde{U}$,
		\item every maximal element $\order{i}{x}$ of dimension 3 in $U$ has a corresponding maximal element of dimension 2 in $\tilde{U}$, the image of $\compos{\bound{}{-}\order{i}{x}}$, with the same lower\nbd dimensional boundaries.
	\end{itemize}
	Because $\tilde{U}$ is 2\nbd dimensional, it is acyclic, so $\maxflow{r}{\tilde{U}}$ is acyclic, and consequently $\maxflow{r}{U}$ is acyclic.
\end{proof}

\begin{cor} \label{cor:dim3_has_frame_acyclic_molecules}
	Let $P$ be an oriented graded poset, $\dim{P} \leq 3$.
	Then
	\begin{enumerate}
		\item $P$ has frame-acyclic molecules,
		\item $\molecin{P}$ is a polygraph.
	\end{enumerate}
\end{cor}

\begin{thm} \label{thm:round_submolecule_dim3}
	Let $\imath\colon V \incl U$ be an inclusion of molecules such that $\dim{U} = \dim{V} = 3$ and $V$ is round.
	The following are equivalent:
	\begin{enumerate}[label=(\alph*)]
	    \item $\imath$ is a submolecule inclusion;
	    \item $\flow{2}{V}$ is a path-induced subgraph of $\flow{2}{U}$.
	\end{enumerate}
\end{thm}
\begin{proof}
	One implication is Proposition \ref{prop:round_submolecule_flow_path_induced}, so we only need to prove the converse.
	By Lemma \ref{lem:connected_subgraph_conditions_path_induced}, if $\flow{2}{V}$ is path-induced, then there exists a 2\nbd ordering $(\order{i}{x})_{i=1}^m$ of $U$ in which the elements of $\imath(V)$ are consecutive, that is, $\order{i}{x} \in \imath(V)$ if and only if $p \leq i \leq q$ for some $p, q \in \set{1, \ldots, m}$.

	By Theorem \ref{thm:dim3_frame_acyclic}, $U$ is frame-acyclic, so by Corollary \ref{cor:frame_acyclicity_equivalent_conditions} the 2\nbd ordering comes from a 2\nbd layering $(\order{i}{U})_{i=1}^m$ such that $\imath(\bound{}{-}V) \subseteq \bound{}{-}\order{p}{U}$.
	Since both are 2\nbd dimensional molecules and $\bound{}{-}V$ is round, we have $\imath(\bound{}{-}V) \submol \bound{}{-}\order{p}{U}$ by Theorem \ref{thm:round_submolecule_dim2}, and we conclude by Lemma \ref{lem:round_submolecules_from_layering}.
\end{proof}

\clearpage
\thispagestyle{empty}

%% file: special.tex
\chapter{Special shapes} \label{chap:special}
\thispagestyle{firstpage}

\begin{guide}
	Every model of higher categories must, at least, specify how $n$\nbd cells relate to $(n+1)$\nbd cells in a higher category.
	This is typically modelled by an algebra of \emph{face operators}, restricting an $(n+1)$\nbd cell to its $n$\nbd dimensional faces.
	Dually, this algebra specifies the ``shape'' of an $(n+1)$\nbd cell, the face operators being dual to \cemph{coface} inclusions of $n$\nbd dimensional shapes into an $(n+1)$\nbd dimensional shape.
	As a consequence, a higher category, within a particular model, has an underlying presheaf on some \cemph{shape category}, whose objects are shapes of $n$\nbd cells for each $n \in \mathbb{N}$ and morphisms are at least coface inclusions (but often include \cemph{codegeneracy} maps as well).
	
	In all the examples of which we are aware, the shapes admit an interpretation as shapes of higher-categorical pasting diagrams, which is then leveraged in order to include strict $\omega$\nbd categories into the model via a nerve construction.
	Thus, we can study these shape categories as special cases of the general theory of higher-categorical pasting diagrams.

	In most cases, the choice of shapes is ``minimalistic'' --- ideally, one shape per dimension --- so one looks at \emph{uniform} families of shapes, such that all $k$\nbd dimensional faces of an $n$\nbd dimensional shape have the same shape, for all $k$ and $n \geq k$.
	There are exactly three such families up to duality, all of them starting with the point: the \cemph{globes}
\[\begin{tikzcd}[column sep=scriptsize]
	\bullet & \bullet & \bullet & \bullet & \bullet & \bullet & \bullet & \bullet
	\arrow[from=1-1, to=1-2]
	\arrow[""{name=0, anchor=center, inner sep=0}, curve={height=12pt}, from=1-3, to=1-4]
	\arrow[""{name=1, anchor=center, inner sep=0}, curve={height=12pt}, from=1-7, to=1-8]
	\arrow[""{name=2, anchor=center, inner sep=0}, curve={height=-12pt}, from=1-7, to=1-8]
	\arrow[shorten <=2pt, shorten >=2pt, Rightarrow, from=1-6, to=1-7]
	\arrow[""{name=3, anchor=center, inner sep=0}, curve={height=12pt}, from=1-5, to=1-6]
	\arrow[""{name=4, anchor=center, inner sep=0}, curve={height=-12pt}, from=1-5, to=1-6]
	\arrow[""{name=5, anchor=center, inner sep=0}, curve={height=-12pt}, from=1-3, to=1-4]
	\arrow[shorten <=3pt, shorten >=3pt, Rightarrow, from=1, to=2]
	\arrow[shorten <=3pt, shorten >=3pt, Rightarrow, from=3, to=4]
	\arrow[shorten <=3pt, shorten >=3pt, Rightarrow, from=0, to=5]
\end{tikzcd}\]
	the \cemph{oriented simplices}
\[\begin{tikzcd}[column sep=small]
	&&& \bullet &&& \bullet & \bullet &&& \bullet & \bullet \\
	\bullet & \bullet & \bullet && \bullet & \bullet &&& \bullet & \bullet &&& \bullet
	\arrow[from=2-1, to=2-2]
	\arrow[""{name=0, anchor=center, inner sep=0}, curve={height=0pt}, from=2-3, to=2-5]
	\arrow[curve={height=-0pt}, from=2-3, to=1-4]
	\arrow[curve={height=-0pt}, from=1-4, to=2-5]
	\arrow[""{name=1, anchor=center, inner sep=0}, curve={height=0pt}, from=2-6, to=2-9]
	\arrow[curve={height=-0pt}, from=2-6, to=1-7]
	\arrow[from=1-7, to=1-8]
	\arrow[""{name=2, anchor=center, inner sep=0}, curve={height=-0pt}, from=1-8, to=2-9]
	\arrow[""{name=3, anchor=center, inner sep=0}, curve={height=0pt}, from=2-6, to=1-8]
	\arrow[""{name=4, anchor=center, inner sep=0}, curve={height=-0pt}, from=2-10, to=1-11]
	\arrow[from=1-11, to=1-12]
	\arrow[curve={height=-0pt}, from=1-12, to=2-13]
	\arrow[""{name=5, anchor=center, inner sep=0}, curve={height=0pt}, from=2-10, to=2-13]
	\arrow[""{name=6, anchor=center, inner sep=0}, curve={height=0pt}, from=1-11, to=2-13]
	\arrow[shorten <=4pt, Rightarrow, from=0, to=1-4]
	\arrow[curve={height=0pt}, shorten <=5pt, Rightarrow, from=1, to=1-8]
	\arrow[shorten <=2pt, Rightarrow, from=3, to=1-7]
	\arrow[shorten <=17pt, shorten >=17pt, Rightarrow, from=2, to=4]
	\arrow[curve={height=-0pt}, shorten <=5pt, Rightarrow, from=5, to=1-11]
	\arrow[shorten <=2pt, Rightarrow, from=6, to=1-12]
\end{tikzcd}\]
	and the \cemph{oriented cubes}
	\[\begin{tikzcd}[sep=small]
	&&& \bullet &&& \bullet & \bullet &&& \bullet & \bullet \\
	\bullet & \bullet & \bullet && \bullet & \bullet & \bullet && \bullet & \bullet && \bullet & \bullet \\
	&&& \bullet &&& \bullet & \bullet &&& \bullet & \bullet
	\arrow[from=2-1, to=2-2]
	\arrow[from=2-3, to=3-4]
	\arrow[from=3-4, to=2-5]
	\arrow[from=2-3, to=1-4]
	\arrow[from=1-4, to=2-5]
	\arrow[shorten <=3pt, shorten >=3pt, Rightarrow, from=3-4, to=1-4]
	\arrow[from=2-6, to=3-7]
	\arrow[from=3-7, to=3-8]
	\arrow[from=3-8, to=2-9]
	\arrow[from=2-6, to=1-7]
	\arrow[from=1-7, to=1-8]
	\arrow[from=1-8, to=2-9]
	\arrow[from=2-10, to=3-11]
	\arrow[from=3-11, to=3-12]
	\arrow[from=3-12, to=2-13]
	\arrow[from=2-10, to=1-11]
	\arrow[from=1-11, to=1-12]
	\arrow[from=1-12, to=2-13]
	\arrow[from=2-7, to=3-8]
	\arrow[from=2-7, to=1-8]
	\arrow[from=2-6, to=2-7]
	\arrow[from=1-11, to=2-12]
	\arrow[from=3-11, to=2-12]
	\arrow[from=2-12, to=2-13]
	\arrow[shorten <=3pt, shorten >=3pt, Rightarrow, from=3-11, to=1-11]
	\arrow[Rightarrow, from=3-12, to=2-12]
	\arrow[Rightarrow, from=2-12, to=1-12]
	\arrow[Rightarrow, from=3-7, to=2-7]
	\arrow[Rightarrow, from=2-7, to=1-7]
	\arrow[shorten <=3pt, shorten >=3pt, Rightarrow, from=3-8, to=1-8]
	\arrow[Rightarrow, from=2-9, to=2-10]
\end{tikzcd}\]
	which are directed versions of \emph{ditopes}, \emph{simplices}, and \emph{hypercubes}, three uniform families of polytopes.

	We will see that the existence and uniformity of these three families is a consequence of being inductively generated by three constructions of oriented graded posets that we studied in Chapter \ref{chap:constructions}: globes by \emph{suspensions}, oriented simplices by \emph{cones}, and oriented cubes by \emph{cylinders}.

	Globes are the simplest possible shape in each dimension.
	However, since they have only one input and one output face in each dimension, they are not very well-adapted to expressing \emph{composition}, which involves mediation between a diagram of multiple cells and a single cell, its composite.
	For this reason, globular models of higher categories tend not to work with presheaves over globes, but over \emph{pasting diagrams of globes}.
	These pasting diagrams are known by a multitude of names, but what remains stable is that their category is labelled $\thetacat$, so in this book we will call them \cemph{thetas}.
	Thetas form the shape category of the \emph{complete Segal $\grade{n}{\thetacat}$\nbd space} model of $(\infty, n)$\nbd categories \cite{rezk2010cartesian}, among others.
	They admit a convenient combinatorial representation as \cemph{plane trees} (Proposition 
	\ref{prop:thetas_are_encoded_by_trees}), and can also be characterised as an inductive subclass of molecules, generated by suspensions and pastings at the 0\nbd boundary.

	Oriented simplices and cubes do not suffer from the same limitation as globes: it is possible to express certain compositions of $n$\nbd dimensional simplices or cubes through a single $(n+1)$\nbd dimensional simplex or cube.
	Consequently, oriented simplices and cubes form, directly, the shape categories of the \emph{complicial} \cite{verity2008complicial} and \emph{comical} \cite{campion2020cubical} models of $(\infty, n)$\nbd categories, respectively.
	We will show that the simplex category and the category of cubes with connections, two common shape categories, can be represented faithfully as subcategories of $\rdcpxmap$, and that \emph{cubical compositions} can be modelled as comaps in $\rdcpxcomap$.

	There is one more class of shapes that enjoyed a significant popularity in the early days of higher category theory: this is the \cemph{opetopes} \cite{baez1998higher}, also known as \emph{multitopes} \cite{hermida2000weak}.
	The opetopes are ``many-to-one'' cell shapes, which can be seen as a suitable combinatorics for the algebra of higher-dimensional composition: an $(n+1)$\nbd dimensional opetope can express the composition of any number of $n$\nbd dimensional opetopes in its input boundary.
	Because of the ``many-to-one'' constraint, the composable opetopic diagrams are \emph{tree-like}, and indeed opetopes can be classified by certain higher-dimensional tree-like structures called \cemph{zoom complexes}.

	There is one serious obstacle in modelling opetopes in our framework: the underlying cell complex structure of an opetope is not always regular.
	This is due to the fact that, for $n > 1$, there are \emph{nullary} $n$\nbd dimensional opetopes whose input boundary is degenerate, that is, has dimension strictly lower than $n - 1$.
	The first such opetope has the shape
	\[\begin{tikzcd}[sep=scriptsize]
	{} \\
	\bullet
	\arrow[loop, looseness=7, out=135, in=45]
	\arrow[Rightarrow, from=2-1, to=1-1]
\end{tikzcd}\]
	which describes a non-regular cell structure on a topological disk.
	On the other hand, it has already been observed that admitting such shapes leads to serious complications both on the combinatorial and on the topological side, which has led M.~Zawadowski to introduce \cemph{positive opetopes} as a better-behaved class leading to a better-behaved shape category \cite{zawadowski2017positive}.
	Positivity bars nullary shapes, and on opetopes it is essentially equivalent to regularity.
	We are thus able to characterise positive opetopes as yet another inductive subclass of molecules.

	One could write an entire book on the combinatorics of each of these classes of shapes, so the aim of this chapter is not to be exhaustive.
	Instead, we try to give a feeling of how it can be fruitful to treat these classes \emph{as subclasses of molecules}.
	Indeed, we are able to give sleek proofs of non-trivial facts about each of these classes using some general machinery that applies to all molecules: for example, that they are classes of \emph{acyclic} molecules, that a positive opetope can be reconstructed uniquely from its zoom complex (Proposition 
	\ref{prop:positive_opetope_reconstruction_zoom_complex}), and that recognition of rewritable submolecules is trivial for the boundaries of positive opetopes (Proposition 
	\ref{prop:inclusions_of_positive_opetope_trees_are_submolecule}).
\end{guide}


\section{Globes and thetas} \label{sec:globes_and_thetas}

\begin{guide}
	Globes are the very standard of what the shape of an $n$\nbd categorical cell is, so to a certain extent they are unavoidable in developing a higher category theory; every other cell shape must be able to encode a globe shape.

	We define globes as the class of oriented graded posets generated by the point under suspension, and deduce that all globes are acyclic molecules.
	We show that suspension can also be replaced by the rewrite construction as a constructor for globes (Corollary 
	\ref{cor:globes_as_generated_by_point_and_atom}).
	We classify the elements of globes and their maps, and prove that the $n$\nbd globe is terminal in the full subcategory of $\rdcpxcomap$ on round $n$\nbd dimensional molecules (Proposition 
	\ref{prop:comap_to_globe}).
	This formalises the intuition that ``all round molecules are subdivisions of globes''.
	
	Then, we focus on thetas, which we define as generated by the point under suspension and pasting at the 0\nbd boundary (a ``wedge sum'').
	We reprove the correspondence between thetas and plane trees (Proposition 
	\ref{prop:thetas_are_encoded_by_trees}), and finally characterise thetas as being precisely the molecules whose every atom is a globe (Proposition \ref{prop:thetas_are_molecules_of_globes}).
\end{guide}

\begin{dfn}[Globe] \index{globe} \index{$\globe{n}$}
	The class of \emph{globes} is the inductive subclass of oriented graded posets closed under isomorphisms and generated by the following clauses.
	\begin{enumerate}
		\item (\textit{Point}). The point is a globe.
		\item (\textit{Suspend}). If $U$ is a globe, then $\sus{U}$ is a globe.
	\end{enumerate}
	We let
	\[
		\globe{0} \eqdef 1, \quad \quad \globe{n} \eqdef \sus{\globe{n-1}}
	\]
	for $n > 0$, and call $\globe{n}$ the \emph{$n$\nbd globe}.
\end{dfn}

\begin{lem} \label{lem:basic_properties_of_globes}
	Let $U$ be a globe, $n \eqdef \dim{U}$.
	Then
	\begin{enumerate}
		\item $U$ is an atom,
		\item $U$ is uniquely isomorphic to $\globe{n}$,
		\item $U$ is acyclic.
	\end{enumerate}
\end{lem}
\begin{proof}
By induction on the construction of $U$.
If $U$ was produced by (\textit{Point}), then $U$ is a 0\nbd dimensional atom, isomorphic to $\globe{0}$, and is obviously acyclic.
If $U$ was produced by (\textit{Suspend}), then it is equal to $\sus{V}$ for some globe $V$.
By the inductive hypothesis, $V$ is an atom, isomorphic to $\globe{n}$ for $n \eqdef \dim{V}$, and is acyclic.
Then $U$ is isomorphic to $\sus{\globe{n}} = \globe{n+1}$, it is an atom by Proposition 
\ref{prop:suspension_of_molecules}, and is acyclic by Proposition \ref{prop:suspension_of_acyclic}.
\end{proof}

\begin{lem} \label{lem:globes_via_rewrite}
	Let $U$ be a globe.
	Then $U \celto U$ is defined and isomorphic to $\sus{U}$.
\end{lem}
\begin{proof}
By Lemma \ref{lem:basic_properties_of_globes}, we can assume that $U$ is $\globe{n}$ for $n \eqdef \dim{U}$, and proceed by induction on $n$.
If $n = 0$, then $\globe{0} \celto \globe{0}$ is the arrow $\thearrow{}$, and we can explicitly construct an isomorphism with $\globe{1}$.
For $n > 0$, we have $\globe{n} = \sus{\globe{n-1}}$, and by the inductive hypothesis the latter is isomorphic to $\globe{n-1} \celto \globe{n-1}$.
It then follows from Proposition \ref{prop:suspension_of_molecules} that $\sus{\globe{n}}$ is isomorphic to $\sus{\globe{n-1}} \celto \sus{\globe{n-1}}$, that is, to $\globe{n} \celto \globe{n}$.
\end{proof}

\begin{cor} \label{cor:globes_as_generated_by_point_and_atom}
	The class of globes is the inductive subclass of oriented graded posets closed under isomorphisms and generated by the clauses
	\begin{enumerate}
		\item (\textit{Point}). The point is a globe.
		\item (\textit{Atom}). If $U$ is a globe, then $U \celto U$ is a globe.
	\end{enumerate}
\end{cor}
\begin{proof}
	By Lemma \ref{lem:globes_via_rewrite}, (\textit{Atom}) is equivalent to (\textit{Suspend}) on globes.
\end{proof}

\begin{lem} \label{lem:boundaries_of_globes}
	Let $U$ be a globe, $n \in \mathbb{N}$, $\alpha \in \set{+, -}$, and $x \in U$.
	Then 
	\begin{enumerate}
		\item $\bound{n}{\alpha}U$ is a globe,
		\item $\clset{x}$ is a globe.
	\end{enumerate}
\end{lem}
\begin{proof}
	We proceed by induction on $m \eqdef \dim{U}$.
	If $m = 0$, then $U$ is the point, $\bound{n}{\alpha}U = U$, and $x$ is the only element of $U$, so $\clset{x} = U$.
	If $m > 0$, then by Corollary \ref{cor:globes_as_generated_by_point_and_atom} $U$ is isomorphic to $V \celto V$ for some globe $V$ of dimension $m-1$.
	If $n \geq m$, $\bound{n}{\alpha}U = U$, and if $n < m$, by globularity of $U$ and $V$ and Lemma \ref{lem:boundaries_of_rewrite},
	\[
		\bound{n}{\alpha}U = \bound{n}{\alpha}(\bound{m-1}{\alpha}U) = \bound{n}{\alpha}(\bound{m-1}{\alpha}V) = \bound{n}{\alpha}V,
	\]
	which is a globe by the inductive hypothesis.
	Moreover, $x$ is either the greatest element of $U$, in which case $\clset{x} = U$, or $x \in \bound{}{}U$, so $x$ is in the image of one of the copies of $V$ and $\clset{x}$ is a globe by the inductive hypothesis.
\end{proof}

\begin{lem} \label{lem:duals_of_globes}
	Let $U$ be a globe, $J \subseteq \posnat$.
	Then $\dual{J}{U}$ is a globe isomorphic to $U$.
\end{lem}
\begin{proof}
	By induction on the construction of $U$ using Corollary \ref{cor:globes_as_generated_by_point_and_atom}.
	If $U$ was produced by (\textit{Point}), this is obvious.
	If $U$ was produced by (\textit{Atom}), then it is isomorphic to $V \celto V$ for some globe $V$.
	By Proposition \ref{prop:dual_preserves_molecules}, $\dual{J}{U}$ is isomorphic to $\dual{J}{V} \celto \dual{J}{V}$, which by the inductive hypothesis is isomorphic to $V \celto V$, hence to $U$.
\end{proof}

\begin{dfn}[Elements of globes] \index{globe!elements}
	By recursion on $n \in \mathbb{N}$, we introduce the following notation for elements of the $n$\nbd globe.
	We let $0$ be the unique element of the $0$\nbd globe.
	For all $n > 0$ and $\alpha \in \set{+, -}$, we let
	\begin{align*}
		0^\alpha & \eqdef \bot^\alpha, \\
		k^\alpha & \eqdef \sus{{(k-1)}^\alpha} \quad \text{for $0 < k < n$}, \\
		n & \eqdef \sus{(n-1)}.
	\end{align*}
\end{dfn}

\begin{lem} \label{lem:faces_of_globes}
	Let $k, n \in \mathbb{N}$, $\alpha \in \set{+, -}$, and suppose $k < n$.
	Then
	\begin{enumerate}
		\item $\faces{n}{\alpha}\globe{n} = \grade{n}{(\globe{n})} = \set{n}$,
		\item $\faces{k}{\alpha}\globe{n} = \set{k^\alpha}$ and $\grade{k}{(\globe{n})} = \set{k^+, k^-}$.
	\end{enumerate}
\end{lem}
\begin{proof}
	Straightforward induction using Lemma \ref{lem:faces_of_suspension}.
\end{proof}

\begin{lem} \label{lem:unique_surjective_map_between_globes}
	Let $k, n \in \mathbb{N}$, $k \leq n$.
	Then there exists a unique surjective map $\tau\colon \globe{n} \to \globe{k}$.
\end{lem}
\begin{proof}
	The unique map $\varepsilon\colon \globe{n-k} \to \globe{0}$ to the terminal object is surjective.
	Then
	\[
		\tau \eqdef \underbrace{\sus{}\ldots\sus{}}_{\text{$k$ times}}\varepsilon\colon \globe{n} \to \globe{k}
	\]
	is a surjective map by Proposition \ref{prop:suspension_is_a_functor_on_maps}.
	Given another surjective map $p\colon \globe{n} \to \globe{k}$, we must have $p(n) = k$.
	Let $j < n$ and $\alpha \in \set{+, -}$.
	For $j \geq k$, a recursive application of Lemma \ref{lem:if_dim_decreased_faces_map_to_same_element} forces $p(j^\alpha) = k$, while for $j < k$, since $p(\bound{j}{\alpha}n) = \bound{j}{\alpha}k$, necessarily $p(j^\alpha) = j^\alpha$.
\end{proof}

\begin{prop} \label{prop:factorisation_of_maps_between_globes}
	Let $n, m \in \mathbb{N}$, let $p\colon \globe{n} \to \globe{m}$ be a map, and let $k \eqdef \dim{p(\globe{n})}$.
	Then $p$ factors uniquely as
	\begin{enumerate}
		\item the unique surjective map $\tau\colon \globe{n} \to \globe{k}$,
		\item followed by an inclusion $\imath\colon \globe{k} \incl \globe{m}$.
	\end{enumerate}
\end{prop}
\begin{proof}
	By Lemma \ref{lem:boundaries_of_globes}, the image of $p$ is a $k$\nbd globe.
	The statement then follows from Proposition \ref{prop:em_factorisation_in_rdcpxmap} together with Lemma \ref{lem:unique_surjective_map_between_globes}.
\end{proof}

\begin{prop} \label{prop:comap_to_globe}
	Let $U$ be a regular directed complex, $n \in \mathbb{N}$.
	The following are equivalent:
	\begin{enumerate}[label=(\alph*)]
		\item there exists a comap $c\colon U \to \globe{n}$;
		\item $U$ is a round molecule and $\dim{U} = n$.
	\end{enumerate}
	When either equivalent condition holds, the comap $c$ is unique.
\end{prop}
\begin{proof}
	Suppose that a comap $c\colon U \to \globe{n}$ exists.
	By Proposition \ref{prop:comaps_inverse_image_preserves_molecules}, $U$ is a molecule.
	Since $\globe{n}$ is an atom, hence round, by Lemma \ref{lem:comaps_inverse_image_and_roundness} $U$ is also round.
	Conversely, let
	\begin{align*}
		c\colon U & \to \globe{n}, \\
			x & \mapsto \begin{cases}
				n & \text{if $x \in \inter{U}$,} \\
				k^\alpha & \text{if $x \in \inter{\bound{k}{\alpha}U}$, $k < n$, $\alpha \in \set{+, -}$};
			\end{cases}
	\end{align*}
	this is well-defined as a function by Lemma \ref{lem:round_partition_into_interiors}, and it is straightforward to check that it determines a comap.
	Uniqueness then follows from Lemma \ref{lem:comaps_inverse_image_interiors}, since $\inter{\globe{n}} = \set{n}$ and $\inter{\bound{k}{\alpha}\globe{n}} = \set{k^\alpha}$ by Lemma \ref{lem:faces_of_globes}.
\end{proof}

\begin{lem} \label{lem:atom_is_cellular_extension_of_boundary}
	Let $U$ be an atom, $n \eqdef \dim{U}$, and let $c\colon U \to \globe{n}$ be the unique comap to the $n$\nbd globe.
	Then the diagram
\begin{equation} \label{eq:atom_as_cellular_extension_pushout}
	\begin{tikzcd}
		{\molecin{\bound{}{}\globe{n}}} && {\molecin{\globe{n}}} \\
	{\molecin{\bound{}{}U}} && {\molecin{U}}
	\arrow[hook, from=1-1, to=1-3]
	\arrow["{\pb{(\restr{c}{\bound{}{}U})}}", from=1-1, to=2-1]
	\arrow["{\pb{c}}", from=1-3, to=2-3]
	\arrow[hook, from=2-1, to=2-3]
	\arrow["\lrcorner"{anchor=center, pos=0.125, rotate=180}, draw=none, from=2-3, to=1-1]
\end{tikzcd}\end{equation}
	is a pushout square in $\omegacat$.
\end{lem}
\begin{proof}
	Let $X$ be a strict $\omega$\nbd category and let
\begin{equation} \label{eq:atom_as_cellular_extension_square}
		\begin{tikzcd}
	{\molecin{\bound{}{}\globe{n}}} && {\molecin{\globe{n}}} \\
	{\molecin{\bound{}{}U}} && {X}
	\arrow[hook, from=1-1, to=1-3]
	\arrow["{\pb{(\restr{c}{\bound{}{}U})}}", from=1-1, to=2-1]
	\arrow["d", from=1-3, to=2-3]
	\arrow["f", from=2-1, to=2-3]
\end{tikzcd}\end{equation}
	be a commutative square of strict functors.
	Then $d$ classifies a cell $u \eqdef d\isocl{\idd{\globe{n}}}$ of $X$ while $f$ restricts to two parallel pasting diagrams of shape $\bound{}{\alpha}U$ for $\alpha \in \set{+, -}$, whose composites we denote by $u^+$ and $u^-$.
	Then commutativity of (\ref{eq:atom_as_cellular_extension_square}) means that $\bound{n-1}{\alpha}u = u^\alpha$ for all $\alpha \in \set{+, -}$.

	Now, $\bound{}{}U = \skel{n-1}{U}$, so $\molecin{\bound{}{}U}$ is isomorphic to $\skel{n-1}{\molecin{U}}$ by Proposition \ref{prop:skeleta_of_molecin}.
	Moreover, by Corollary \ref{cor:atom_only_one_topdim_cell}, $\isocl{\idd{U}}$ is the only $n$\nbd dimensional cell of $\molecin{U}$.
	Then
	\begin{align*}
		f'\colon \molecin{U} & \to X, \\
		t & \mapsto \begin{cases}
			f(t) & \text{if $\dim{t} < n$}, \\
			u & \text{if $t = \isocl{\idd{U}}$}
		\end{cases}
	\end{align*}
	is well-defined as a strict functor, and it is the unique strict functor that extends $f$ and has the property that $d = f' \after \pb{c}$, since this forces
	\[
		f'\isocl{\idd{U}} = f'(\pb{c}\isocl{\idd{\globe{n}}}) = d\isocl{\idd{\globe{n}}} = u.
	\]
	This proves that (\ref{eq:atom_as_cellular_extension_pushout}) is a pushout.
\end{proof}

\begin{dfn}[Theta] \index{theta}
	The class of \emph{thetas} is the inductive subclass of oriented graded posets closed under isomorphisms and generated by the following clauses.
	\begin{enumerate}
		\item (\textit{Point}). The point is a theta.
		\item (\textit{Suspend}). If $U$ is a theta, then $\sus{U}$ is a theta.
		\item (\textit{Wedge}). If $U, V$ are thetas and $0 < \min \set{\dim{U}, \dim{V}}$, then $U \cp{0} V$ is a theta.
	\end{enumerate}
\end{dfn}

\begin{rmk}
	Since the clauses generating globes are a subset of those generating thetas, every globe is a theta.
\end{rmk}

\begin{comm}
Thetas, or the $\omega$\nbd categories that they present, are known by a variety of names in the literature.
They were called \emph{Batanin cells} in \cite{joyal1997disks} after \cite{batanin1998monoidal}, \emph{globular cardinals} in \cite{street2000petit}, \emph{simple $\omega$\nbd categories} in \cite{makkai2001duality}, \emph{cells} in \cite{berger2002cellular}, and \emph{globular pasting schemes} in \cite{weber2004generic}.

The one constant seems to be that their category (usually, the full subcategory of $\omegacat$ on the $\omega$\nbd categories that they present) is denoted by $\thetacat$.
For this reason we propose yet another name based on this denotation, which has at least the advantage of brevity.
\end{comm}

\begin{lem} \label{lem:thetas_basic_properties}
	Let $U$ be a theta.
	Then $U$ is an acyclic molecule.
\end{lem}
\begin{proof}
	By induction on the construction of $U$.
	If $U$ was produced by (\textit{Point}), then it is a 0\nbd dimensional atom and obviously acyclic.
	
	If $U$ was produced by (\textit{Suspend}), then it is equal to $\sus{V}$ for some theta $V$.
	By the inductive hypothesis, $V$ is an acyclic molecule.
	Then $U$ is a molecule by Proposition \ref{prop:suspension_of_molecules} and acyclic by Proposition \ref{prop:suspension_of_acyclic}.

	If $U$ was produced by (\textit{Wedge}), then it is equal to $V \cp{0} W$ for some thetas $V, W$, which are acyclic molecules by the inductive hypothesis.
	Then $U$ is a molecule by definition and acyclic by Lemma \ref{lem:acyclic_closed_under_pasting}.
\end{proof}

\begin{lem} \label{lem:boundaries_of_thetas}
	Let $U$ be a theta, $n \in \mathbb{N}$, $\alpha \in \set{+, -}$, and $x \in U$.
	Then
	\begin{enumerate}
		\item $\bound{n}{\alpha}U$ is a theta,
		\item $\clset{x}$ is a globe.
	\end{enumerate}
\end{lem}
\begin{proof}
	If $n = 0$ or $\dim{x} = 0$, then $\bound{0}{\alpha}U$ and $\clset{x}$ are isomorphic to the point, which is a theta, for all molecules $U$, so we can suppose that $n > 0$ and $\dim{x} > 0$.
	We proceed by induction on the construction of $U$.
	If $U$ was produced by (\textit{Point}), then $U = 1$, $\bound{n}{\alpha}1 = 1$, and $\dim{x} = 0$.
	
	If $U$ was produced by (\textit{Suspend}), then it is equal to $\sus{V}$ for some theta $V$.
	By Corollary \ref{cor:boundary_of_suspension}, $\bound{n}{\alpha}U = \sus{\bound{n-1}{\alpha}V}$, and $\bound{n-1}{\alpha}V$ is a theta by the inductive hypothesis.
	It follows that $\bound{n}{\alpha}U$ is a theta.
	Moreover, $x = \sus{x'}$ for some $x' \in V$, so $\clset{x} = \sus{\clset{x'}}$, which is a globe by the inductive hypothesis.

	If $U$ was produced by (\textit{Wedge}), then it is equal to $V \cp{0} W$ for some thetas $V, W$.
	By Lemma \ref{lem:pasting_higher_boundary}, $\bound{n}{\alpha}U$ is isomorphic to $\bound{n}{\alpha}V \cp{0} \bound{n}{\alpha}W$, and both $\bound{n}{\alpha}V$ and $\bound{n}{\alpha}W$ are thetas by the inductive hypothesis, so $\bound{n}{\alpha}U$ is a theta.
	Moreover, $x$ is either in the image of $V$ or in the image of $W$, so $\clset{x}$ is a globe by the inductive hypothesis.
\end{proof}

\begin{dfn}[Plane tree] \index{plane tree}
The set of \emph{plane trees} is inductively generated by the following clauses.
\begin{enumerate}
	\item (\textit{Root}). The root $\treeroot$ is a plane tree.
	\item (\textit{Branch}). If $m > 0$ and $t_1, \ldots, t_m$ are plane trees, then $\branch{t_1, \ldots, t_m}$ is a plane tree.
\end{enumerate}
\end{dfn}

\begin{comm}
	Plane trees are also known as \emph{level trees} or \emph{ordered trees}.
\end{comm}

\begin{dfn}[Height of a plane tree] \index{plane tree!height} \index{$\height{t}$}
Let $t$ be a plane tree.
The \emph{height} $\height{t}$ of $t$ is defined by structural induction by
\begin{itemize}
	\item $\height{\treeroot} \eqdef 0$,
	\item $\height{\branch{t_1, \ldots, t_m}} \eqdef \max\set{\height{t_1}, \ldots, \height{t_m}} + 1$.
\end{itemize}
\end{dfn}

\begin{dfn}[Theta encoded by a plane tree] \index{theta!encoded by a plane tree} \index{plane tree!theta encoded by} \index{$\thetac{t}$}
	Let $t$ be a plane tree.
	The \emph{theta encoded by $t$} is the theta $\thetac{t}$ defined as follows, by structural induction on $t$.
	\begin{itemize}
		\item If $t = \treeroot$, then $\thetac{t} \eqdef 1$.
		\item If $t = \branch{t_1, \ldots, t_m}$, then $\thetac{t} \eqdef \sus{\thetac{t_1}} \cp{0} \ldots \cp{0} \sus{\thetac{t_m}}$.
	\end{itemize}
\end{dfn}

\begin{lem} \label{lem:plane_tree_height_is_dimension_of_theta}
	Let $t$ be a plane tree.
	Then $\dim{\thetac{t}} = \height{t}$.
\end{lem}
\begin{proof}
	By structural induction on $t$.
	If $t = \treeroot$, then $\dim{\thetac{t}} = \dim{1} = 0 = \height{t}$.
	If $t = \branch{t_1, \ldots, t_m}$, then
	\begin{align*}
		\dim{\thetac{t}} & = \dim{\sus{\thetac{t_1}} \cp{0} \ldots \cp{0} \sus{\thetac{t_m}}} = \max \set{\dim{\sus{\thetac{t_1}}}, \ldots, \dim{\sus{\thetac{t_m}}}} = \\
				 & = \max \set{\dim{\thetac{t_1}} + 1, \ldots, \dim{\thetac{t_m} + 1}} = \\
				 & = \max \set{\height{t_1}, \ldots, \height{t_m}} + 1 = \height{t}
	\end{align*}
	using the inductive hypothesis on $t_1, \ldots, t_m$.
\end{proof}

\begin{lem} \label{lem:points_of_theta_encoded_by_tree}
	Let $t = \branch{t_1, \ldots, t_m}$ be a plane tree. 
	Then $\size{\grade{0}{(\thetac{t})}} = m + 1$.
\end{lem}
\begin{proof}
	By induction on $m \geq 1$.
	If $m = 1$, then $\thetac{t} = \sus{\thetac{t_1}}$, and by construction $\size{\grade{0}{(\thetac{t})}} = 2$.
	If $m > 1$, let $t' \eqdef \branch{t_1, \ldots, t_{m-1}}$.
	Then $\thetac{t} = \thetac{t'} \cp{0} \sus{\thetac{t_m}}$, and
	\begin{align*}
		\size{\grade{0}{(\thetac{t})}} & = \size{\grade{0}{(\thetac{t'})}} + \size{\grade{0}{(\sus{\thetac{t_m}})}} - \size{\grade{0}{(\thetac{t'} \cap \sus{\thetac{t_m}})}} = \\
					       & = m + 2 - 1 = m + 1
	\end{align*}
	since $\thetac{t'} \cap \sus{\thetac{t_m}} = \bound{0}{+}\thetac{t'}$ is a point.
\end{proof}

\begin{lem} \label{lem:uniqueness_of_tree_encoding}
	Let $t, t'$ be plane trees and suppose $\thetac{t}$ is isomorphic to $\thetac{t'}$.
	Then $t = t'$.
\end{lem}
\begin{proof}
	By Lemma \ref{lem:plane_tree_height_is_dimension_of_theta}, we must have $n \eqdef \height{t} = \height{t'}$; we proceed by induction on $n$.
	When $n = 0$, both $t$ and $t'$ must be $\treeroot$.
	Suppose $n > 0$.
	Then $t = \branch{t_1, \ldots, t_m}$ and $t' = \branch{t'_1, \ldots, t'_{m'}}$, and by Lemma \ref{lem:points_of_theta_encoded_by_tree} necessarily $m = m'$.
	Let $\varphi\colon \thetac{t} \to \thetac{t'}$ be an isomorphism, and let
	\[
		\set{x_i} \eqdef \begin{cases}
			\bound{0}{-}\sus{\thetac{t_0}} & \text{if $i = 0$}, \\
			\bound{0}{+}\sus{\thetac{t_i}} & \text{if $i \in \set{1, \ldots, m}$},
		\end{cases} \quad
		\set{x'_i} \eqdef \begin{cases}
			\bound{0}{-}\sus{\thetac{t'_0}} & \text{if $i = 0$}, \\
			\bound{0}{+}\sus{\thetac{t'_i}} & \text{if $i \in \set{1, \ldots, m}$}.
		\end{cases}
	\]
	Then $\grade{0}{(\thetac{t})} = \set{x_0, \ldots, x_m}$ and $\grade{0}{(\thetac{t'})} = \set{x'_0, \ldots, x'_m}$.
	By construction of $\thetac{t}$ and $\thetac{t'}$, if $i \leq j$, then
	\[
		x_i \precflow x_j, \quad \quad x'_i \precflow x'_j,
	\]
	and the reverse inequality only holds if $i = j$, since the flow preorder is a linear order by Proposition \ref{prop:acyclic_flow_is_linear_order} and Lemma \ref{lem:thetas_basic_properties}.
	Since isomorphisms preserve dimensions and the flow preorder, we must have $\varphi(x_i) = x'_i$ for all $i \in \set{0, \ldots, m}$, hence also
	\[
		\varphi(\sus{\thetac{t_i}}) = \sus{\thetac{t'_i}}
	\]
	for all $i \in \set{1, \ldots, m}$.
	Thus $\varphi$ determines isomorphisms between $\sus{\thetac{t_i}}$ and $\sus{\thetac{t'_i}}$,
	which restrict to isomorphisms between $\thetac{t_i}$ and $\thetac{t'_i}$ for all $i \in \set{1, \ldots, m}$.
	By the inductive hypothesis, $t_i = t'_i$ for all $i \in \set{1, \ldots, m}$, hence $t = t'$.
\end{proof}

\begin{prop} \label{prop:thetas_are_encoded_by_trees}
	Let $U$ be a theta.
	Then there exists a unique plane tree $t$ such that $U$ is isomorphic to $\thetac{t}$.
\end{prop}
\begin{proof}
	We prove the existence of $t$ by induction on the construction of $U$; uniqueness will then follow from Lemma \ref{lem:uniqueness_of_tree_encoding}.
	If $U$ was produced by (\textit{Point}), then $U$ is equal to $\thetac{\treeroot}$.

	If $U$ was produced by (\textit{Suspend}), then $U$ is equal to $\sus{V}$ for some theta $V$.
	By the inductive hypothesis, there exists a plane tree $t'$ such that $V$ is isomorphic to $\thetac{t'}$.
	Then $U$ is isomorphic to $\thetac{t}$ for $t \eqdef \branch{t'}$.

	If $U$ was produced by (\textit{Wedge}), then $U$ is equal to $V \cp{0} W$ for some thetas $V, W$ with $0 < \min \set{\dim{V}, \dim{W}}$.
	By the inductive hypothesis, there exist plane trees $t_1$, $t_2$ such that $V$ is isomorphic to $\thetac{t_1}$ and $W$ to $\thetac{t_2}$.
	Because $\dim{V}, \dim{W} > 0$, necessarily $t_1, t_2 \neq \treeroot$, so there exist sequences of plane trees $(\order{i}{t_1})_{i=1}^m$, $(\order{i}{t_2})_{i=1}^p$ such that
	\[
		t_1 = \branch{\order{1}{t_1}, \ldots, \order{m}{t_1}}, \quad t_2 = \branch{\order{1}{t_2}, \ldots, \order{p}{t_2}}.
	\]
	Then $U$ is isomorphic to $\thetac{t}$ for $t \eqdef \branch{\order{1}{t_1}, \ldots, \order{m}{t_1}, \order{1}{t_2}, \ldots, \order{p}{t_2}}$.
\end{proof}

\begin{exm}[The theta encoded by a plane tree] \label{exm:theta_encoding}
	Let $t$ be the plane tree
	\[
		\branch{\branch{\bullet, \bullet}, \bullet, \branch{\bullet}, \bullet}
	\]
	which in graphical form looks like
	\[\begin{tikzcd}
	&& \bullet \\
	& \bullet & \bullet & \bullet & \bullet \\
	\bullet & \bullet && \bullet && {.}
	\arrow[from=1-3, to=2-2]
	\arrow[from=2-2, to=3-1]
	\arrow[from=2-2, to=3-2]
	\arrow[from=1-3, to=2-3]
	\arrow[from=1-3, to=2-4]
	\arrow[from=2-4, to=3-4]
	\arrow[from=1-3, to=2-5]
\end{tikzcd}\]
	Then $\thetac{t}$ is the oriented face poset of the 2\nbd dimensional diagram
\[\begin{tikzcd}
	\bullet & \bullet & \bullet & \bullet & \bullet
	\arrow[""{name=0, anchor=center, inner sep=0}, curve={height=18pt}, from=1-1, to=1-2]
	\arrow[from=1-2, to=1-3]
	\arrow[""{name=1, anchor=center, inner sep=0}, curve={height=12pt}, from=1-3, to=1-4]
	\arrow[from=1-4, to=1-5]
	\arrow[""{name=2, anchor=center, inner sep=0}, curve={height=-12pt}, from=1-3, to=1-4]
	\arrow[""{name=3, anchor=center, inner sep=0}, curve={height=-18pt}, from=1-1, to=1-2]
	\arrow[""{name=4, anchor=center, inner sep=0}, from=1-1, to=1-2]
	\arrow[shorten <=3pt, shorten >=3pt, Rightarrow, from=1, to=2]
	\arrow[shorten <=2pt, shorten >=2pt, Rightarrow, from=0, to=4]
	\arrow[shorten <=2pt, shorten >=2pt, Rightarrow, from=4, to=3]
\end{tikzcd}\]
	which is isomorphic to $\sus{(\sus{1} \cp{0} \sus{1})} \cp{0} \sus{1} \cp{0} \sus{(\sus{1})} \cp{0} \sus{1}$.
\end{exm}

\begin{lem} \label{lem:thetas_closed_under_pasting}
	Let $U, V$ be thetas, $k \in \mathbb{N}$, and suppose $U \cp{k} V$ is defined.
	Then $U \cp{k} V$ is a theta.
\end{lem}
\begin{proof}
	We proceed by induction on $k$.
	If $k = 0$, this holds by definition, so let $k > 0$ and suppose that $k < \min \set{\dim{U}, \dim{V}}$, for otherwise $U \cp{k} V$ is isomorphic to $U$ or $V$.
	By Proposition \ref{prop:thetas_are_encoded_by_trees}, there exist unique plane trees $t, t'$ such that $U$ is isomorphic to $\thetac{t}$ and $V$ is isomorphic to $\thetac{t'}$.
	Since $\dim{U}, \dim{V} > 0$, there exist sequences $(t_i)_{i=1}^m$, $(t'_i)_{i=1}^p$ such that 
	\[
		t = \branch{t_1, \ldots, t_m}, \quad t' = \branch{t'_1, \ldots, t'_p}.
	\]
	Then $U$ and $V$ are isomorphic to 
	\[
		\sus{\thetac{t_1}} \cp{0} \ldots \cp{0} \sus{\thetac{t_m}}, \quad \sus{\thetac{t'_1}} \cp{0} \ldots \cp{0} \sus{\thetac{t'_p}},
	\]
	respectively.
	Since $k > 0$, by Lemma \ref{lem:pasting_higher_boundary} and Corollary \ref{cor:boundary_of_suspension}, $\bound{k}{+}U$ and $\bound{k}{-}V$ are isomorphic to
	\[
		\sus{\bound{k-1}{+}\thetac{t_1}} \cp{0} \ldots \cp{0} \sus{\bound{k-1}{+}\thetac{t_m}}, \quad \sus{\bound{k-1}{-}\thetac{t'_1}} \cp{0} \ldots \cp{0} \sus{\bound{k-1}{-}\thetac{t'_p}},
	\]
	respectively.
	By Lemma \ref{lem:boundaries_of_thetas} and Proposition \ref{prop:thetas_are_encoded_by_trees}, there exist unique sequences $(s_i)_{i=1}^m$, $(s'_i)_{i=1}^p$ of plane trees such that
	\begin{itemize}
		\item $\bound{k-1}{+}\thetac{t_i}$ is isomorphic to $\thetac{s_i}$ for each $i \in \set{1, \ldots, m}$,
		\item $\bound{k-1}{-}\thetac{t'_i}$ is isomorphic to $\thetac{s'_i}$ for each $i \in \set{1, \ldots, p}$.
	\end{itemize}
	Letting $s \eqdef \branch{s_1, \ldots, s_m}$ and $s' \eqdef \branch{s'_1, \ldots, s'_p}$, we deduce that $\bound{k}{+}U$ is isomorphic to $\thetac{s}$ and $\bound{k}{-}V$ is isomorphic to $\thetac{s'}$. 
	Since $U \cp{k} V$ is defined, $\thetac{s}$ is isomorphic to $\thetac{s'}$, and by Lemma \ref{lem:uniqueness_of_tree_encoding} we conclude that $s = s'$, and in particular that $m = p$.

	It follows that $\sus{\thetac{t_i}} \cp{k} \sus{\thetac{t'_i}}$ is defined for all 
	$i \in \set{1, \ldots, m}$, and by Proposition \ref{prop:suspension_of_molecules} it is isomorphic to $\sus{(\thetac{t_i} \cp{k-1} \thetac{t'_i})}$.	
	By the inductive hypothesis, $\order{i}{W} \eqdef \thetac{t_i} \cp{k-1} \thetac{t'_i}$ is a theta for all $i \in \set{1, \ldots, m}$, and
	\[
		\sus{\order{1}{W}} \cp{0} \ldots \cp{0} \sus{\order{m}{W}}
	\]
	is isomorphic to $U \cp{k} V$.
	We conclude that $U \cp{k} V$ is a theta.
\end{proof}

\begin{prop} \label{prop:thetas_are_molecules_of_globes}
	Let $U$ be a molecule.
	The following are equivalent:
	\begin{enumerate}[label=(\alph*)]
		\item $U$ is a theta;
		\item for all $x \in U$, $\clset{x}$ is a globe.
	\end{enumerate}
\end{prop}
\begin{proof}
	One implication is given by Lemma \ref{lem:boundaries_of_thetas}.
	For the converse implication, we proceed by induction on the construction of $U$.
	If $U$ was produced by (\textit{Point}) or (\textit{Atom}), then it is an atom by Lemma
	\ref{lem:atom_greatest_element}, so it has a greatest element $\top$.
	Then by assumption $U = \clset{\top}$ is a globe, hence a theta.

	If $U$ was produced by (\textit{Paste}), then it splits into proper submolecules $V \cup W$ along the $k$\nbd boundary.
	By the inductive hypothesis, $V$ and $W$ are thetas.
	We conclude by Lemma \ref{lem:thetas_closed_under_pasting} that $U$ is a theta.
\end{proof}

\begin{cor} \label{cor:thetas_generated_by_globes_under_pasting}
	The class of thetas is the inductive subclass of oriented graded posets closed under isomorphisms and generated by the clauses
	\begin{enumerate}
		\item (\textit{Globe}). If $U$ is a globe, then $U$ is a theta.
		\item (\textit{Paste}). If $U, V$ are thetas and $k < \min \set{\dim{U}, \dim{V}}$ is such that $U \cp{k} V$ is defined, then $U \cp{k} V$ is a theta.
	\end{enumerate}
\end{cor}

\begin{lem} \label{lem:thetas_closed_under_duals}
	Let $U$ be a theta, $J \subseteq \posnat$.
	Then $\dual{J}{U}$ is a theta.
\end{lem}
\begin{proof}
	Every element of $\dual{J}{U}$ is of the form $\dual{J}{x}$ for some $x \in U$.
	Then $\clset{\dual{J}{x}} = \dual{J}{\clset{x}}$, which is a globe by 
	Lemma \ref{lem:boundaries_of_thetas} combined with Lemma \ref{lem:duals_of_globes}.
	We conclude by Proposition \ref{prop:thetas_are_molecules_of_globes}.
\end{proof}


\section{Oriented simplices} \label{sec:simplices}

\begin{guide}
	Simplices are, by quite a margin, the most well-studied class of shapes for combinatorial topology.
	As a shape category for higher category theory they do not quite have the same strength --- most notably, there is no good directed version of barycentric subdivision --- but they remain a formidable contender, on the back of their simplicity and the wealth of machinery that can be borrowed from simplicial homotopy theory.

	We define oriented simplices as the class of oriented graded posets generated by the point under cones, that is, joins with a point.
	Equivalently, they are the closure of the point under joins (Corollary 
	\ref{cor:simplices_generated_by_joins}).
	We deduce that all oriented simplices are acyclic molecules.
	We classify the elements of an oriented $n$\nbd simplex, as well as their input and output faces (Lemma \ref{lem:elements_of_the_simplex}), and prove that the simplex category can be represented as the full subcategory of $\rdcpxmap$ on the oriented simplices (Proposition 
	\ref{prop:simplex_category_in_rdcpxmap}).
\end{guide}

\begin{dfn}[Oriented simplex] \index{simplex!oriented} \index{$\simplex{n}$}
	The class of \emph{oriented simplices} is the inductive subclass of oriented graded posets closed under isomorphisms and generated by the following clauses.
	\begin{enumerate}
		\item (\textit{Point}). The point is an oriented simplex.
		\item (\textit{Cone}). If $U$ is an oriented simplex, then $1 \join U$ is an oriented simplex.
	\end{enumerate}
	We let
	\[
		\simplex{0} \eqdef 1, \quad \quad \simplex{n} \eqdef 1 \join \simplex{n-1}
	\]
	for $n > 0$, and call $\simplex{n}$ the \emph{oriented $n$\nbd simplex}.
\end{dfn}

\begin{comm}
	Oriented simplices were called \emph{orientals} in \cite{street1987algebra}.
	It can be useful to also define $\simplex{-1} \eqdef \varnothing$.
\end{comm}

\begin{lem} \label{lem:simplices_basic_properties}
	Let $U$ be an oriented simplex, $n \eqdef \dim{U}$.
	Then
	\begin{enumerate}
		\item $U$ is an atom,
		\item $U$ is uniquely isomorphic to $\simplex{n}$,
		\item $U$ is acyclic.
	\end{enumerate}
\end{lem}
\begin{proof}
	By induction on the construction of $U$.
	If $U$ was produced by (\textit{Point}), then $U = \simplex{0} = 1$ is a 0\nbd dimensional atom, and is obviously acyclic.
	If $U$ was produced by (\textit{Cone}), then $U$ is equal to $1 \join V$ for some oriented simplex $V$.
	By the inductive hypothesis, $V$ is an acyclic atom, isomorphic to $\simplex{m}$ for $m \eqdef \dim{V}$.
	Moreover, $\dim{U} = \dim{1} + \dim{V} + 1$, so $m = n-1$.
	By Proposition \ref{prop:join_of_molecules} combined with Proposition \ref{prop:join_of_acyclic}, $U$ is an acyclic atom, isomorphic to $1 \join \simplex{n-1} = \simplex{n}$.
\end{proof}

\begin{lem} \label{lem:simplices_closed_under_join}
	Let $U$, $V$ be oriented simplices.
	Then $U \join V$ is an oriented simplex.
\end{lem}
\begin{proof}
	By Lemma \ref{lem:simplices_basic_properties}, $U$ and $V$ are isomorphic to $\simplex{n}$ and $\simplex{m}$, respectively, for $n \eqdef \dim{U}$ and $m \eqdef \dim{V}$.
	Then, $U \join V$ is isomorphic to
	\[
		\underbrace{(1 \join (1 \join \ldots (1 \join 1)\ldots))}_{\text{$n+1$ times}} 
		\join
		\underbrace{(1 \join (1 \join \ldots (1 \join 1)\ldots))}_{\text{$m+1$ times}},
	\]
	which by associativity of joins is isomorphic to
	\[
		1\join \underbrace{(1 \join \ldots (1 \join 1)\ldots)}_{\text{$n+m+1$ times}} = 1 \join \simplex{n+m} = \simplex{n+m+1},
	\]
	and we conclude.
\end{proof}

\begin{cor} \label{cor:simplices_generated_by_joins}
	The class of oriented simplices is the inductive subclass of oriented graded posets closed under isomorphisms and generated by the clauses
	\begin{enumerate}
		\item (\textit{Point}). The point is an oriented simplex.
		\item (\textit{Join}). If $U$, $V$ are oriented simplices, then $U \join V$ is an oriented simplex.
	\end{enumerate}
\end{cor}

\begin{comm}
	From now on, we will use coherence of the monoidal structure to avoid bracketing iterated joins.
\end{comm}

\begin{lem} \label{lem:faces_of_simplex}
	Let $U$ be an oriented simplex, $z \in U$.
	Then $\clset{z}$ is an oriented simplex.
\end{lem}
\begin{proof}
	By induction on the construction of $U$.
	If $U$ was produced by (\textit{Point}), then $U = 1$ and $z$ is the only element of $U$, so $\clset{z} = U$ is an oriented 0\nbd simplex.
	If $U$ was produced by (\textit{Cone}), then $U$ is equal to $1 \join V$ for some oriented simplex $V$.
	If $z = \inj{x}$ for the unique $x \in 1$, then $\clset{z}$ is isomorphic to 1.
	If $z = \inr{y}$ for some $y \in V$, then $\clset{z}$ is isomorphic to $\clset{y}$, which is an oriented simplex by the inductive hypothesis.
	Finally, if $z = x \join y$ for the unique $x \in 1$ and some $y \in W$, then $\clset{z}$ is isomorphic to $1 \join \clset{y}$, which by the inductive hypothesis is an oriented simplex.
\end{proof}

\begin{lem} \label{lem:odd_dual_of_oriented_simplex}
	Let $U$ be an oriented simplex.
	Then $\opp{U}$ is an oriented simplex isomorphic to $U$.
\end{lem}
\begin{proof}
	We proceed by induction on the construction of $U$ using Corollary \ref{cor:simplices_generated_by_joins}.
	If $U$ was produced by (\textit{Point}), then $U = 1$ and $\opp{1}$ is isomorphic to $1$.
	If $U$ was produced by (\textit{Join}), then $U$ is isomorphic to $V \join W$ for some oriented simplices $V$, $W$.
	By Proposition \ref{prop:odd_duals_and_joins}, $\opp{U}$ is isomorphic to $\opp{W} \join \opp{V}$, which is an oriented simplex by the inductive hypothesis.
	Because duals preserve dimension, we conclude by Lemma 
	\ref{lem:simplices_basic_properties}.
\end{proof}

\begin{prop} \label{prop:simplex_factorisation}
	Let $n, m \in \mathbb{N}$, let $p\colon \simplex{n} \to \simplex{m}$ be a map, and let $k \eqdef \dim{p(\simplex{n})}$.
	Then $p$ factors uniquely as
	\begin{enumerate}
		\item a surjective map $\widehat{p}\colon \simplex{n} \to \simplex{k}$,
		\item followed by an inclusion $\imath\colon \simplex{k} \incl \simplex{m}$.
	\end{enumerate}
\end{prop}
\begin{proof}
	By Lemma \ref{lem:faces_of_simplex}, the image of $p$ is a $k$\nbd simplex.
	The statement then follows from Proposition \ref{prop:em_factorisation_in_rdcpxmap}.
\end{proof}

\begin{dfn}[Elements of oriented simplices] \index{simplex!elements}
	Let $(b_0 \ldots b_n)$ be a string of $(n+1)$ bits, that is, $b_i \in \set{0, 1}$ for all $i \in \set{0, \ldots, n}$, and let $k \eqdef \sum_{i=0}^n b_i - 1$.
	If $k \geq 0$, then there is an inclusion $\simplex{k} \incl \simplex{n}$ given by the composite of
	\begin{enumerate}
		\item the isomorphism between $\simplex{k}$ and $U_0 \join \ldots \join U_n$, where
			\[
				U_i \eqdef \begin{cases}
					\varnothing & \text{if $b_i = 0$}, \\
					1 & \text{if $b_i = 1$},
				\end{cases}
			\]
		\item followed by the inclusion defined by 
			\begin{align*}
				\imath_0 \join \ldots \join \imath_n &\colon U_0 \join \ldots \join U_n \incl \simplex{n}, \\
				\imath_i & \eqdef \begin{cases}
					\varnothing \incl 1 & \text{if $b_i = 0$}, \\
					\idd{1}\colon 1 \to 1 & \text{if $b_i = 1$},
				\end{cases}
			\end{align*}
			for each $i \in \set{0, \ldots, n}$.
	\end{enumerate}
	We denote by $(b_0\ldots b_n)$ the greatest element of the image of this inclusion.
\end{dfn}

\begin{lem} \label{lem:elements_of_the_simplex}
	Let $n \in \mathbb{N}$, $z \in \simplex{n}$, and $k \eqdef \dim{z}$.
	Then 
	\begin{enumerate}
		\item $z = (b_0\ldots b_n)$ for a unique bit string with $k = \sum_{i=0}^n b_i - 1$, 
		\item if $k > 0$, then
		\begin{align*}
			\faces{}{+}z & = \set{(b_0 \ldots b_{j-1} 0 b_{j+1} \ldots b_n) \mid 
				\text{$b_j = 1$ and $\sum_{i=0}^{j-1} b_i$ is even} }, \\
			\faces{}{-}z & = \set{(b_0 \ldots b_{j-1} 0 b_{j+1} \ldots b_n) \mid 
				\text{$b_j = 1$ and $\sum_{i=0}^{j-1} b_i$ is odd} }.
		\end{align*}
	\end{enumerate}
\end{lem}
\begin{proof}
	We proceed by induction on $n$.
	If $n = 0$, then $z$ is the only element of the point, and $z = (1)$.
	Since $\dim{z} = 0$, there is nothing else to prove.
	Suppose $n > 0$.
	Then $\simplex{n} = 1 \join \simplex{n-1}$.
	If $z = \inj{x}$ for some $x \in 1$, then $x$ is the only element of $1$, $\dim{z} = 0$, and $z = (10\ldots0)$.
	If $z = \inr{y}$ for some $y \in \simplex{n-1}$, by the inductive hypothesis $y = (b_0 \ldots b_{n-1})$ for some string $(b_i)_{i=0}^n$ of $n$ bits, hence $z = (0b_0\ldots b_{n-1})$, and by Lemma \ref{lem:join_faces}
	\[
		\faces{}{\alpha}z = \set{\inr{y'} \mid y' \in \faces{}{\alpha}y} = \set{(0b'_0 \ldots b'_{n-1}) \mid (b'_0 \ldots b'_{n-1}) \in \faces{}{\alpha}y},
	\]
	from which the statement promptly follows using the inductive hypothesis.
	Finally, if $z = x \join y$ for the unique element $x \in 1$ and some $y \in \simplex{n-1}$, then $z$ is equal to $(1 b_0 \ldots b_{n-1})$ for some string $(b_i)_{i=0}^n$ of $n$ bits, and
	\begin{align*}
		\faces{}{+}z & = \set{\inr{y}} + \set{x \join y' \mid y' \in \faces{}{-}y} = \\
				   & = \set{(0 b_0 \ldots b_{n-1})} + \set{(1 b'_0 \ldots b'_{n-1}) \mid (b'_0 \ldots b'_{n-1}) \in \faces{}{-}y}, \\
		\faces{}{-}z & = \set{x \join y' \mid y' \in \faces{}{+}y} = \set{(1 b'_0 \ldots b'_{n-1}) \mid (b'_0 \ldots b'_{n-1}) \in \faces{}{+}y},
	\end{align*}
	and, for all $j \in \set{0, \ldots, n}$, $\sum_{i=0}^{j-1} b_i$ has the opposite parity of $1 + \sum_{i=0}^{j-1} b_i$, so we conclude using the inductive hypothesis.
\end{proof}

\begin{cor} \label{cor:number_of_faces_of_simplex}
	Let $n \in \mathbb{N}$, $k \leq n$.
	Then $\size{\grade{k}{(\simplex{n})}} = \displaystyle \binom{n+1}{k+1}$.
\end{cor}
\begin{proof}
	Each element $(b_0 \ldots b_n)$ of dimension $k$ is uniquely identified by a choice of $k + 1$ positions in which $b_i = 1$.
\end{proof}

\begin{dfn}[The simplex category $\simplexcat$] \index{simplex!category} \index{$\simplexcat$}
	The \emph{simplex category} is the full subcategory $\Delta$ of $\poscat$ whose objects are non-empty finite linear orders $[n] \eqdef \set{0 < \ldots < n}$ for all $n \in \mathbb{N}$.
\end{dfn}

\begin{prop} \label{prop:simplex_category_in_rdcpxmap}
	Let $f\colon [n] \to [m]$ be a morphism in the simplex category.
	Then $\simplex{}f\colon \simplex{n} \to \simplex{m}$, defined by
	\[
		(b_0\ldots b_n) \mapsto (b'_0 \ldots b'_m), \quad \quad
		b'_j \eqdef \begin{cases}
			1 & \text{if $b_i = 1$ for some $i \in \invrs{f}(j)$}, \\
			0 & \text{otherwise},
		\end{cases}
	\]
	is a map of regular directed complexes.
	This assignment determines a full and faithful functor $\simplexcat \incl \rdcpxmap$.
\end{prop}
\begin{proof}
	We can factorise $f$ as a surjection $\widehat{f}\colon [n] \to [k]$ followed by an injection $\imath\colon [k] \incl [m]$, where by construction $\sum_{i=0}^m b'_i = k + 1$.
	For each $j > 0$, let $\varepsilon_j\colon \simplex{j-1} \to 1$ be the unique map to the terminal object of $\rdcpxmap$.
	Then $\simplex{}f$ is equal to the composite of
	\[
		\varepsilon_{\size{\invrs{\widehat{f}}(0)}} \join \ldots \join \varepsilon_{\size{\invrs{\widehat{f}}(k)}} \colon \simplex{n} \to \simplex{k}
	\]
	with the inclusion $\simplex{k} \incl \simplex{m}$ whose image is $\clset{(b'_0\ldots b'_m)}$.
	By Proposition \ref{prop:join_monoidal_on_maps}, both of these are maps of regular directed complexes, so $\simplex{}f$ is a map of regular directed complexes.

	Functoriality is straightforward, as is faithfulness, so it suffices to show that every map $p\colon \simplex{n} \to \simplex{m}$ is equal to $\simplex{}f$ for some $f\colon [n] \to [m]$.
	By Proposition \ref{prop:simplex_factorisation}, we can do this separately for inclusions and surjective maps.
	Since $\simplex{n}$ is an atom, if $p$ is an inclusion, then it is isomorphic to the inclusion $\clset{(b_0\ldots b_m)} \incl \simplex{m}$ for some $(b_0\ldots b_m) \in \simplex{m}$.
	Then $p = \simplex{}f$ for the injection $f\colon [n] \to [m]$ defined by
	\[
		f(i) \eqdef \min \set{j \in \set{0, \ldots, m} \mid \sum_{\ell \leq j} b_\ell - 1 = i}.
	\]
	Finally, suppose that $p$ is a surjective map; we proceed by induction on $n$.
	If $n \in \set{0, 1}$, $p$ is either the unique map to the terminal object, or it is a dimension-preserving map $\simplex{1} \to \simplex{1}$, which is necessarily the identity by 
	Proposition \ref{prop:characterisation_of_morphisms_among_maps} and 
	Theorem \ref{thm:morphisms_of_atoms_are_injective}.
	All of these are equal to $\simplex{}f$ for some order-preserving surjection $f\colon [n] \to [m]$.
	Suppose $n > 1$, and let
	\begin{align*}
		v_i \eqdef & (\underbrace{0\ldots 0}_{\text{$i$}}1
		\underbrace{0\ldots 0}_{\text{$n-i$}}) \in \simplex{n}, \\
		w_j \eqdef & (\underbrace{0\ldots 0}_{\text{$j$}}1
		\underbrace{0\ldots 0}_{\text{$m-j$}}) \in \simplex{m}
	\end{align*}
	for each $i \in \set{0, \ldots, n}$ and $j \in \set{0, \ldots, m}$.
	By Lemma \ref{lem:faces_of_simplex}, these are all the 0\nbd dimensional elements of $\simplex{n}$ and $\simplex{m}$.
	We define $f\colon [n] \to [m]$ by \[
		f(i) \eqdef j \quad \text{ if } \quad p(v_i) = w_j
	\]
	for each $i \in \set{0, \ldots, n}$ and $j \in \set{0, \ldots, m}$.

	We claim that $f$ is order-preserving.
	Let $k \leq \ell \in \set{0, \ldots, n}$.
	Because $n > 1$, there exists an order-preserving injection $\imath\colon [n-1] \incl [n]$ such that $k = \imath(k')$ and $\ell = \imath(\ell')$ for some $k', \ell' \in \set{0, \ldots, n-1}$.
	By the inductive hypothesis, the composite $p \after \simplex{}\imath\colon \simplex{n-1} \to \simplex{m}$ is equal to $\simplex{}f'$ for some order-preserving map $f'\colon [n-1] \to [m]$, and necessarily $f \after \imath = f'$, so $f(k) = f'(k') \leq f'(\ell') = f(\ell)$.	
	It follows that $\simplex{}f$ is defined.
	For any order-preserving injection $\imath\colon [n-1] \incl [n]$ we must have $p \after \simplex{}\imath = \simplex{}(f \after \imath)$, so $p$ and $\simplex{}f$ are equal on $\bound{}{}\simplex{n}$.
	Since $p$ is surjective, $f$ is also surjective, so both $p$ and $\simplex{}f$ send the greatest element of $\simplex{n}$ to the greatest element of $\simplex{m}$.
	We conclude that $p = \simplex{}f$.
\end{proof}

\begin{dfn}[Simplicial coface maps] \index{simplex!coface} \index{$\simface{i}$}
	Let $n > 0$.
	The \emph{coface maps of the oriented $n$\nbd simplex} are the inclusions $\simface{i}\colon \simplex{n-1} \incl \simplex{n}$ obtained, for each $i \in \set{0, \ldots, n}$, by composing
	\[
		\idd{} \join \eta \join \idd{}\colon \simplex{i-1} \join \varnothing \join \simplex{n-i-1} \incl \simplex{i-1} \join 1 \join \simplex{n-i-1},
	\]
	where $\eta$ is the unique inclusion of type $\varnothing \incl 1$, with the unique isomorphisms $\simplex{n-1} \iso \simplex{i-1} \join \varnothing \join \simplex{n-i-1}$ and $\simplex{i-1} \join 1 \join \simplex{n-i-1} \iso \simplex{n}$.
\end{dfn}

\begin{dfn}[Simplicial codegeneracy maps] \index{simplex!codegeneracy} \index{$\simdeg{i}$}
	Let $n \in \mathbb{N}$.
	The \emph{codegeneracy maps of the oriented $n$\nbd simplex} are the surjective maps $\simdeg{i}\colon \simplex{n+1} \to \simplex{n}$ obtained, for each $i \in \set{0, \ldots, n}$, by composing
	\[
		\idd{} \join \varepsilon \join \idd{}\colon \simplex{i-1} \join \simplex{1} \join \simplex{n-i-1} \to \simplex{i-1} \join 1 \join \simplex{n-i-1},
	\]
	where $\varepsilon$ is the unique map of type $\simplex{1} \to 1$, with the unique isomorphisms $\simplex{n+1} \iso \simplex{i-1} \join \simplex{1} \join \simplex{n-i-1}$ and $\simplex{i-1} \join 1 \join \simplex{n-i-1} \iso \simplex{n}$.
\end{dfn}


\section{Oriented cubes} \label{sec:cubes}

\begin{guide}
	Cubical combinatorics have historically been less popular in topology than simplicial combinatorics, although they have had a recent resurgence.
	As a shape category for higher category theory, cubes and simplices have somewhat complementary strengths: closure under Gray products makes oriented cubes convenient for describing higher morphisms, but cubes are somewhat ``too symmetrical'' to conveniently express composition, which is an intrinsically asymmetric, many-to-one operation.

	We define oriented cubes as the class of oriented graded posets generated by the point under cylinders, that is, Gray products with the arrow.
	Equivalently, they are the closure of the point and arrow under Gray products
	(Corollary \ref{cor:cubes_generated_by_gray}).
	We deduce that all oriented cubes are acyclic molecules.
	We classify the elements of an oriented $n$\nbd cube, as well as their input and output faces (Lemma \ref{lem:elements_of_the_cube}), and describe a faithful representation of the category of cubes with connections in $\rdcpxmap$, as well as a representation of cubical composition via comaps.
\end{guide}

\begin{dfn}[Oriented cube] \index{cube!oriented} \index{$\cube{n}$}
	The class of \emph{oriented cubes} is the inductive subclass of oriented graded posets closed under isomorphisms and generated by the following clauses.
	\begin{enumerate}
		\item (\textit{Point}). The point is an oriented cube.
		\item (\textit{Cylinder}). If $U$ is an oriented cube, then $\thearrow{} \gray U$ is an oriented cube.
	\end{enumerate}
	We let
	\[
		\cube{0} \eqdef 1, \quad \quad \cube{n} \eqdef \thearrow{} \gray \cube{n-1}
	\]
	for $n > 0$, and call $\cube{n}$ the \emph{oriented $n$\nbd cube}.
\end{dfn}

\begin{lem} \label{lem:cubes_basic_properties}
	Let $U$ be an oriented cube, $n \eqdef \dim{U}$.
	Then
	\begin{enumerate}
		\item $U$ is an atom,
		\item $U$ is uniquely isomorphic to $\cube{n}$,
		\item $U$ is acyclic.
	\end{enumerate}
\end{lem}
\begin{proof}
	By induction on the construction of $U$.
	If $U$ was produced by (\textit{Point}), then $U = \cube{0} = 1$ is a 0\nbd dimensional atom, and is obviously acyclic.
	If $U$ was produced by (\textit{Cylinder}), then $U$ is equal to $\thearrow{} \gray V$ for some oriented cube $V$.
	By the inductive hypothesis, $V$ is an acyclic atom, isomorphic to $\cube{m}$ for $m \eqdef \dim{V}$.
	Moreover, $\dim{U} = \dim{\thearrow{}} + \dim{V}$, so $m = n-1$.
	By Proposition \ref{prop:gray_product_of_molecules} combined with Proposition
	\ref{prop:gray_product_of_acyclic}, $U$ is an acyclic atom, isomorphic to $\thearrow{} \gray \cube{n-1} = \cube{n}$.
\end{proof}

\begin{lem} \label{lem:cubes_closed_under_gray}
	Let $U$, $V$ be oriented cubes.
	Then $U \gray V$ is an oriented cube.
\end{lem}
\begin{proof}
	By Lemma \ref{lem:cubes_basic_properties}, $U$ and $V$ are isomorphic to $\cube{n}$ and $\cube{m}$, respectively, for $n \eqdef \dim{U}$ and $m \eqdef \dim{V}$.
	Then, $U \gray V$ is isomorphic to
	\[
		\underbrace{(\thearrow{} \gray (\thearrow{} \gray \ldots (\thearrow{} \gray \thearrow{})\ldots))}_{\text{$n$ times}} 
		\gray
		\underbrace{(\thearrow{} \gray (\thearrow{} \gray \ldots (\thearrow{} \gray \thearrow{})\ldots))}_{\text{$m$ times}},
	\]
	which by associativity of Gray products is isomorphic to
	\[
		\underbrace{\thearrow{}\gray (\thearrow{} \gray \ldots (\thearrow{} \gray \thearrow{})\ldots)}_{\text{$n+m$ times}} = \cube{n+m},
	\]
	and we conclude.
\end{proof}

\begin{cor} \label{cor:cubes_generated_by_gray}
	The class of oriented cubes is the inductive subclass of oriented graded posets closed under isomorphisms and generated by the clauses
	\begin{enumerate}
		\item (\textit{Point}). The point is an oriented cube.
		\item (\textit{Arrow}). The arrow is an oriented cube.
		\item (\textit{Gray}). If $U$, $V$ are oriented cubes and $0 < \min \set{\dim{U}, \dim{V}}$, then $U \gray V$ is an oriented cube.
	\end{enumerate}
\end{cor}

\begin{comm}
	From now on, we will use coherence of the monoidal structure to avoid bracketing iterated Gray products.
\end{comm}

\begin{lem} \label{lem:faces_of_cube}
	Let $U$ be an oriented cube, $x \in U$.
	Then $\clset{x}$ is an oriented cube.
\end{lem}
\begin{proof}
	By induction on the construction of $U$ using Corollary \ref{cor:cubes_generated_by_gray}.
	If $U$ was produced by (\textit{Point}) or (\textit{Arrow}), this is clear.
	If $U$ was produced by (\textit{Gray}), then $U$ is isomorphic to $V \gray W$ for a pair of lower-dimensional oriented cubes $V$, $W$.
	Then $x = (v, w)$ for a unique pair of elements $v \in V$ and $w \in W$, and by the inductive hypothesis $\clset{v}$ and $\clset{w}$ are oriented cubes.
	Then $\clset{x} = \clset{v} \gray \clset{w}$ is an oriented cube by Lemma \ref{lem:cubes_closed_under_gray}. 
\end{proof}

\begin{lem} \label{lem:cubes_and_duals}
	Let $U$ be an oriented cube.
	Then $\opp{U}$, $\coo{U}$, and $\optot{U}$ are oriented cubes isomorphic to $U$.
\end{lem}
\begin{proof}
	By induction on the construction of $U$ using Corollary \ref{cor:cubes_generated_by_gray}.
	If $U$ was produced by (\textit{Point}) or (\textit{Arrow}), this is easy to check.
	If $U$ was produced by (\textit{Gray}), then $U$ is isomorphic to $V \gray W$ for a pair of lower-dimensional oriented cubes $V$, $W$.
	By Proposition \ref{prop:duals_and_gray_products}, $\opp{U}$ is isomorphic to $\opp{W} \gray \opp{V}$, which is isomorphic to $W \gray V$ by the inductive hypothesis.
	By Lemma \ref{lem:cubes_closed_under_gray}, this is an oriented cube of the same dimension as $U$, so it is isomorphic to $U$ by Lemma \ref{lem:cubes_basic_properties}.
	The cases of $\coo{U}$ and $\optot{U}$ are entirely analogous.
\end{proof}

\begin{dfn}[Elements of oriented cubes] \index{cube!elements}
	Let $n \in \mathbb{N}$, let $(c_1\ldots c_n)$ be a string of $n$ symbols $c_i \in \set{0^-, 0^+, 1}$, and let $k \eqdef \size{\set{i \in \set{1, \ldots, n} \mid c_i = 1}}$.
	There is an inclusion $\cube{k} \incl \cube{n}$ given by the composite of
	\begin{enumerate}
		\item the isomorphism between $\cube{k}$ and $U_1 \gray \ldots \gray U_n$, where
			\[
				U_i \eqdef \begin{cases}
					1 
					& \text{if $c_i \in \set{0^+, 0^-}$}, \\
					\thearrow{}
					& \text{if $c_i = 1$},
				\end{cases}
			\]
		\item followed by the inclusion defined by
			\begin{align*}
				\imath_1 \gray \ldots \gray \imath_n &\colon U_1 \gray \ldots \gray U_n \incl \cube{n}, \\
				\imath_i & \eqdef \begin{cases}
					1 \iso \bound{0}{\alpha}\thearrow{} \incl \thearrow{} & \text{if $c_i = 0^\alpha$, $\alpha \in \set{+, -}$}, \\
					\idd{\thearrow{}}\colon \thearrow{} \to \thearrow{} & \text{if $c_i = 1$},
				\end{cases}
			\end{align*}
			for each $i \in \set{1, \ldots, n}$.
	\end{enumerate}
	We denote by $(c_1\ldots c_n)$ the greatest element of the image of this inclusion.
\end{dfn}

\begin{lem} \label{lem:elements_of_the_cube}
	Let $n \in \mathbb{N}$, $x \in \cube{n}$, and $k \eqdef \dim{x}$.
	Then 
	\begin{enumerate}
		\item $x = (c_1\ldots c_n)$ for a unique string with $k = 
			\size{\set{ i \in \set{1, \ldots, n} \mid c_i = 1 }}$,
		\item for all $\alpha \in \set{+, -}$,
		\[
			\faces{}{\alpha}x = \set{(c_1 \ldots c_{j-1} 0^{\alpha\beta} c_{j+1} \ldots c_n) \mid 			\text{$c_j = 1$, $\beta \eqdef (-)^{\size{ \set{ i < j \mid c_i = 1 } }}$ } }.
		\]
	\end{enumerate}
\end{lem}
\begin{proof}
	We proceed by induction on $n$.
	If $n = 0$, then $x$ is the only element of the point, $x = ()$, and its set of faces is empty, which satisfies the statement.
	Suppose $n > 0$.
	Then $\cube{n} = \thearrow{} \gray \cube{n-1}$ and $x = (y, z)$ for some $y \in \thearrow{}$ and $z \in \cube{n-1}$.
	Identifying $\thearrow{}$ with $\globe{1}$, we have $y \in \set{0^-, 0^+, 1}$, while $z = (c_1 \ldots c_{n-1})$ for a unique string $(c_i)_{i=1}^{n-1}$ by the inductive hypothesis.
	Then $x = (c'_1 \ldots c'_n) \eqdef (y c_1 \ldots c_{n-1})$.
	Given $\alpha \in \set{+, -}$, we have
	\[
		\faces{}{\alpha}x 
		= \faces{}{\alpha}y \times \set{(c_1 \ldots c_{n-1})} + \set{y} \times \faces{}{(-)^{\dim{y}}\alpha}(c_1 \ldots c_{n-1}).
	\]
	If $y = 1$, by the inductive hypothesis, this is equal to
	\begin{align*}
		& \set{ (0^\alpha c_1 \ldots c_{n-1}) } + \\
		& \quad + \set{(1 c_1 \ldots c_{j-1} 0^{-\alpha\beta}
			c_{j+1} \ldots c_{n-1} \mid \text{$c_j = 1$, $\beta = (-)^{\size{ \set{ i < j \mid c_i = 1 } }}$} } = \\
		& \quad = \set{(c'_1 \ldots c'_{j-1} 0^{\alpha\beta} c'_{j+1} \ldots c'_n) \mid 
		\text{$c'_j = 1$, $\beta = (-)^{\size{ \set{ i < j \mid c'_i = 1 } }}$} },
	\end{align*}
	while if $y = 0^\gamma$ for some $\gamma \in \set{+, -}$, this is equal to
	\begin{align*}
		& \set{(0^\gamma c_1 \ldots c_{j-1} 0^{\alpha\beta}
			c_{j+1} \ldots c_{n-1} \mid \text{$c_j = 1$, $\beta = (-)^{\size{ \set{ i < j \mid c_i = 1 } }}$} } = \\
		& \quad = \set{(c'_1 \ldots c'_{j-1} 0^{\alpha\beta} c'_{j+1} \ldots c'_n) \mid 
		\text{$c'_j = 1$, $\beta = (-)^{\size{ \set{ i < j \mid c'_i = 1 } }}$} }. \qedhere
	\end{align*}
\end{proof}

\begin{cor} \label{cor:number_of_elements_of_cube}
	Let $n \in \mathbb{N}$, $k \leq n$.
	Then $\size{\grade{k}{(\cube{n})}} = \displaystyle 2^{n-k} \binom{n}{k}$.
\end{cor}
\begin{proof}
	Each element $(c_1 \ldots c_n)$ of dimension $k$ is uniquely identified by a choice of $k$ positions in which $c_i = 1$ together with a choice of signs in $\set{+, -}$ for the $(n - k)$ positions in which $c_i \in \set{0^+, 0^-}$. 
\end{proof}

\begin{dfn}[Cubical coface maps] \index{cube!coface} \index{$\cubeface{i}{+}, \cubeface{i}{-}$}
	Let $n > 0$.
	The \emph{coface maps of the oriented $n$\nbd cube} are the inclusions $\cubeface{i}{\alpha}\colon \cube{n-1} \incl \cube{n}$ obtained, for each $i \in \set{0, \ldots, n-1}$ and $\alpha \in \set{+, -}$, by composing
	\[
		\idd{} \gray 0^\alpha \gray \idd{}\colon \cube{i} \gray 1 \gray \cube{n-i-1} \incl \cube{i} \gray \thearrow{} \gray \cube{n-i-1},
	\]
	where $0^\alpha$ is the inclusion $1 \iso \bound{0}{\alpha}\thearrow{} \incl \thearrow{}$, with the unique isomorphisms $\cube{n-1} \iso \cube{i} \gray 1 \gray \cube{n-i-1}$ and $\cube{i} \gray \thearrow{} \gray \cube{n-i-1} \iso \cube{n}$.
\end{dfn}

\begin{dfn}[Cubical codegeneracy maps] \index{cube!codegeneracy} \index{$\cubedeg{i}$}
	Let $n \in \mathbb{N}$.
	The \emph{codegeneracy maps of the oriented $n$\nbd cube} are the surjective maps $\cubedeg{i}\colon \cube{n+1} \to \cube{n}$ obtained, for each $i \in \set{0, \ldots, n}$, by composing
	\[
		\idd{} \gray \varepsilon \gray \idd{}\colon \cube{i} \gray \thearrow{} \gray \cube{n-i} \to \cube{i} \gray 1 \gray \cube{n-i},
	\]
	where $\varepsilon$ is the unique map of type $\thearrow{} \to 1$, with the unique isomorphisms $\cube{n+1} \iso \cube{i} \gray \thearrow{} \gray \cube{n-i}$ and $\cube{i} \gray 1 \gray \cube{n-i} \iso \cube{n}$.
\end{dfn}

\begin{dfn}[Cubical coconnection maps] \index{cube!coconnection} \index{$\cubeconn{i}{+}, \cubeconn{i}{-}$}
	Let $n > 0$.
	The \emph{coconnection maps of the oriented $n$\nbd cube} are the surjective maps $\cubeconn{i}{\alpha}\colon \cube{n+1} \to \cube{n}$ obtained, for each $i \in \set{0, \ldots, n-1}$ and $\alpha \in \set{+, -}$, by composing
	\[
		\idd{} \gray \cubeconn{}{\alpha} \gray \idd{}\colon \cube{i} \gray \cube{2} \gray \cube{n-i-1} \to \cube{i} \gray \cube{1} \gray \cube{n-i-1},
	\]
	where $\cubeconn{}{\alpha}$ is defined, for each $c \in \set{0^-, 0^+, 1}$, by
	\[
		\cubeconn{}{\alpha}(0^\alpha c) = \cubeconn{}{\alpha}(c 0^\alpha) \eqdef 0^\alpha, \quad \quad \cubeconn{}{\alpha}(1 1) \eqdef 1, \quad \quad \cubeconn{}{\alpha}(0^{-\alpha} c) = \cubeconn{}{\alpha}(c 0^{-\alpha}) \eqdef c,
	\]
	with isomorphisms $\cube{n+1} \iso \cube{i} \gray \cube{2} \gray \cube{n-i-1}$ and $\cube{i} \gray \cube{1} \gray \cube{n-i-1} \iso \cube{n}$.
\end{dfn}

\begin{comm}
	The coface, codegeneracy, and coconnection maps determine a faithful representation in $\rdcpxmap$ of the category of \emph{cubes with connections}, called the \emph{intermediate cubical site} in \cite{grandis2003cubical}.
	We strongly suspect that an analogue of Proposition 
	\ref{prop:simplex_category_in_rdcpxmap} holds, that is, this representation is also full, but we do not have a proof at this time.
\end{comm}

\begin{dfn}[Cubical pasting] \index{cube!pasting}
	Let $n > 0$ and $i \in \set{0, \ldots, n-1}$.
	The \emph{pasting of two oriented $n$\nbd cubes in the $i$-th direction} is the molecule
	\[
		\cube{n} \cubecp{i} \cube{n} \eqdef \cube{i} \gray (\thearrow{} \cp{0} \thearrow{}) \gray \cube{n-i-1},
	\]
	which by Proposition \ref{prop:gray_product_of_molecules} is a generalised pasting of the form $\cube{n} \gencp{n-1} \cube{n}$.
\end{dfn}

\begin{dfn}[Cubical composition comaps] \index{cube!composition} \index{$\cubecomp{i}$}
	Let $n > 0$.
	The \emph{$n$\nbd dimensional cubical composition comaps} are the comaps $\cubecomp{i}\colon \cube{n} \cubecp{i} \cube{n} \to \cube{n}$ obtained, for each $i \in \set{0, \ldots, n-1}$, by composing
	\[
		\idd{} \gray c \gray \idd{}\colon \cube{i} \gray (\thearrow{} \cp{0} \thearrow{}) \gray \cube{n-i-1} \to \cube{i} \gray \thearrow{} \gray \cube{n-i-1}
	\]
	with the unique isomorphism $\cube{i} \gray \thearrow{} \gray \cube{n-i-1} \iso \cube{n}$, where $c$ is the unique comap of type $\thearrow{} \cp{0} \thearrow{} \to \thearrow{}$ whose existence is granted by Proposition \ref{prop:comap_to_globe} since $\thearrow{}$ is isomorphic to $\globe{1}$.
\end{dfn}

\begin{comm}
	In \cite{al2002multiple}, the authors showed that strict $\omega$\nbd categories are equivalent to \emph{cubical $\omega$\nbd categories}, whose structural operations are cubical faces, degeneracies, connections (dual to cofaces, codegeneracies, coconnections), and compositions.
	This equivalence is witnessed by a \emph{cubical nerve} functor for strict $\omega$\nbd categories.
	Our constructions imply that the operations that determine the cubical nerve can all be modelled by maps and comaps of molecules, represented in $\omegacat$ via the covariant and contravariant $\molecin{-}$ functors.
\end{comm}


\section{Positive opetopes} \label{sec:opetopes}

\begin{guide}
	There are many ways in which opetopes have been defined in the literature, and several articles devoted to comparing different definitions, starting from \cite{cheng2004weak}; we do not attempt such a detailed comparison.

	Instead, we mimick the approach of \cite{thanh2022type} to give an inductive definition of positive opetopes --- or rather \cemph{positive opetope trees}, which is what we call the diagram shapes that can appear as input boundaries of positive opetopes --- as the class of oriented graded posets generated by the point under \emph{``grafting''}, that is, pasting at a submolecule in codimension 1, and the rewrite construction restricted to atoms in the output.
	We deduce that positive opetope trees are acyclic molecules, and characterise them as the round molecules whose atoms are all ``many-to-one'' (Proposition 
	\ref{prop:positive_opetopes_are_many_to_one_molecules}).
	We define the zoom complex associated to a positive opetope tree, and prove that it is a faithful encoding of its shape (Proposition \ref{prop:positive_opetope_reconstruction_zoom_complex}).
	Finally, we prove that the rewritable submolecule problem has a trivial solution for positive opetope trees (Proposition \ref{prop:inclusions_of_positive_opetope_trees_are_submolecule}).
\end{guide}

\begin{dfn}[Positive opetope tree] \index{positive opetope!tree}
	The class of \emph{positive opetope trees} is the inductive subclass of oriented graded posets closed under isomorphisms and generated by the following clauses.
	\begin{enumerate}
		\item (\textit{Point}). The point is a positive opetope tree.
		\item (\textit{Graft}). Let $U$, $V$ be positive opetope trees of the same finite dimension and $x \in \bound{}{-}V$ such that $\bound{}{+}U$ is isomorphic to $\clset{x}$.
			Then $U \cpsub{\clset{x}}{-} V$ is a positive opetope tree.
		\item (\textit{Shift}). Let $U$ be a positive opetope tree.
			Then $U \celto \compos{U}$ is a positive opetope tree.
	\end{enumerate}
\end{dfn}

\begin{dfn}[Positive opetope] \index{positive opetope}
	A \emph{positive opetope} is a positive opetope tree with a greatest element.
\end{dfn}

\begin{comm}
	Positive opetopes were first defined in \cite{zawadowski2007positive}.
	Our presentation is based on \cite{thanh2022type}, omitting the ``introduction of degeneracies'' rule which enables the formation of non-positive (non-regular) opetopes.
\end{comm}

\begin{lem} \label{lem:positive_opetope_basic_properties}
	Let $U$ be a positive opetope tree.
	Then 
	\begin{enumerate}
		\item $U$ is a round molecule,
		\item $U$ is acyclic,
		\item for all $k < \dim{U}$, $\size{\faces{k}{+}U} = 1$,
		\item for all $x \in U$, if $\dim{x} > 0$, then $\size{\faces{}{+}x} = 1$.
	\end{enumerate}
\end{lem}
\begin{proof}
	By induction on the construction of $U$.
	If $U$ was produced by (\textit{Point}), then $U$ is a 0\nbd dimensional atom, and obviously acyclic.
	
	If $U$ was produced by (\textit{Graft}), then $U$ is isomorphic to $V \cpsub{\clset{x}}{-} W$ for some pair of positive opetope trees $V$, $W$ and some $x \in \bound{}{-}W$.
	By the inductive hypothesis, $V$ and $W$ are acyclic round molecules, and $\clset{x} \submol \bound{}{-}W$ by Lemma \ref{lem:downset_is_submolecule}, which also ensures that $V \cpsub{\clset{x}}{-} W$ is well-defined.
	Then $U$ is a molecule by Lemma \ref{lem:pasting_at_submolecule}, is round by 
	Lemma \ref{lem:round_pasting_at_submolecule}, and is acyclic by 
	Lemma \ref{lem:acyclic_closed_under_pasting}.
	Moreover, $\bound{}{+}U$ is isomorphic to $\bound{}{+}V$, and $\dim{U} = \dim{V}$, so by the inductive hypothesis we have $\size{\faces{k}{+}U} = \size{\faces{k}{+}V} = 1$ for all $k < \dim{U}$.
	Finally, any $x \in U$ is either in the image of $V$ or in the image of $W$, so if $\dim{x} > 0$ we have $\size{\faces{}{+}x} = 1$ by the inductive hypothesis.
	
	If $U$ was produced by (\textit{Shift}), then $U$ is isomorphic to $V \celto \compos{V}$ for some positive opetope tree $V$.
	By the inductive hypothesis, $V$ is an acyclic round molecule, which also ensures that $V \celto \compos{V}$ is well-defined.
	Then $U$ is an atom, and by construction $\size{\faces{}{+}U} = \size{\maxel{\compos{V}}} = 1$, while by globularity $\faces{k}{+}U$ is isomorphic to $\faces{k}{+}V$ for all $k < \dim{U} - 1$, so $\size{\faces{k}{+}U} = 1$ by the inductive hypothesis.
	Now, if $\dim{U} = 1$, we can rely on Lemma \ref{lem:only_1_molecules} to conclude.
	Otherwise, let $r$ be the single element of $\faces{k}{+}U$ for $k = \dim{U} - 1$.
	Then any $x \in U$ is either in the image of $V$, and we can use the inductive hypothesis, or is the greatest element $\top_{\compos{V}}$ in the image of $\compos{V}$, in which case $\faces{}{+}x = \set{r}$, or is the greatest element $\top$ of $U$, in which case $\faces{}{+}\top = \set{\top_{\compos{V}}}$.
	Finally, any path in $\hasseo{U}$ either
	\begin{itemize}
		\item only passes through vertices in the image of $V$, or
		\item contains a segment of the path $x \to \top_{\compos{V}} \to r$, where $x$ is in the image of $\faces{}{-}V$, or
		\item contains a segment of the path $x \to \top \to \top_{\compos{V}} \to r$, where $x$ is in the image of $\maxel{V}$.
	\end{itemize}
	In the first case, the path cannot contain any cycles by the inductive hypothesis on $V$.
	In the second and third case, since $t$ is the only element in the image of $\faces{}{+}V$, by
	Lemma \ref{lem:path_in_graph_of_molecule} the path from $x$ to $r$ can be replaced with a path entirely contained in the image of $V$.
	Thus, under the hypothesis that the path contains a cycle, replacing all such segments we obtain a cycle in $\hasseo{V}$, contradicting the inductive hypothesis.
	We conclude that $U$ is acyclic.
\end{proof}

\begin{lem} \label{lem:positive_opetope_last_constructor}
Let $U$ be a positive opetope tree.
The following are equivalent:
\begin{enumerate}[label=(\alph*)]
    \item $U$ is a positive opetope;
    \item the final constructor producing $U$ is (\textit{Point}) or (\textit{Shift}).
\end{enumerate}
\end{lem}
\begin{proof}
	Essentially the same as the proof of Lemma \ref{lem:atom_greatest_element}.
\end{proof}

\begin{dfn}[In-tree] \index{directed graph!in-tree} \index{directed graph!in-tree!root}
	Let $\mathscr{G}$ be a directed graph.
	We say that $\mathscr{G}$ is an \emph{in-tree} if there exists a vertex $r$ of $\mathscr{G}$ with the property that, for all other vertices $x$, there exists a unique path from $x$ to $r$ in $\mathscr{G}$.
	When such an $r$ exists, it is unique with this property, and is called the \emph{root of $\mathscr{G}$}.
\end{dfn}

\begin{rmk}
	The underlying undirected graph of an in-tree is a rooted tree; in particular, it is connected and acyclic.
\end{rmk}

\begin{lem} \label{lem:molecules_whose_graph_is_in_tree}
	Let $U$ be a round molecule, $n \eqdef \dim{U}$, and suppose that, if $n > 0$, then $\size{\faces{}{+}x} = 1$ for all $x \in \grade{n}{U}$.
	Then there exists $r \in \grade{n}{U}$ such that
	\begin{enumerate}
		\item $\flow{n-1}{U}$ is an in-tree whose root is $r$,
		\item $\faces{}{+}U = \faces{}{+}r$,
		\item $\graph{U}$ is an in-tree whose root is the only output face of $r$.
	\end{enumerate}
\end{lem}
\begin{proof}
	If $n = 0$, then $U$ is the point, and all statements are trivially true with $r$ the only element of $U$.
	Suppose $n > 0$ and let $x \in \grade{n}{U}$.
	By Corollary \ref{cor:flow_acyclic_in_codimension_1}, $\flow{n-1}{U}$ is acyclic, and since it is finite, there exists a path of maximal length starting from $x$.
	We claim that this path is unique.
	We may proceed by induction on the length $\ell$ of the maximal path.
	If $\ell = 0$, this is obvious, since the only path of length 0 is the constant path, so suppose $\ell > 0$.
	Then by assumption there exists a unique $y \in \faces{}{+}x$, which by Corollary
	\ref{cor:codimension_1_elements} has a unique input coface $x'$.
	It follows that there is a unique edge $x \to x'$ out of $x$, and the length of the maximal path starting from $x'$ is strictly smaller than $\ell$, so by the inductive hypothesis there is a unique maximal path starting from $x'$, which allows us to conclude.

	It follows that there exists a function $r\colon \grade{n}{U} \to \grade{n}{U}$ sending each $x \in \grade{n}{U}$ to the last vertex visited in the unique maximal path starting from $x$.
	Moreover, if there exists an edge $x \to x'$ in $\flow{n-1}{U}$, then necessarily $r(x) = r(x')$.
	Let $R \eqdef r(\grade{n}{U})$; then $\grade{n}{U} = \sum_{x \in R} \invrs{r}(x)$ is a partition of $\grade{n}{U}$ into sets of vertices that are disconnected in $\flow{n-1}{U}$.
	Since by Proposition \ref{prop:round_molecule_connected_flowgraph} $\flow{n-1}{U}$ is connected, necessarily $\size{R} = 1$, and we let $r$ be the unique element of $R$.
	We conclude that $\flow{n-1}{U}$ is an in-tree whose root is $r$.

	As a consequence, $r$ is the only element of $\grade{n}{U}$ with no edges out of it, hence the only element with $\faces{}{+}r \subseteq \faces{}{+}U$.
	By Lemma \ref{lem:round_is_pure}, $U$ is pure, so $\faces{}{+}U = \faces{}{+}r$.
	Finally, by Lemma \ref{lem:path_in_graph_of_molecule}, for every vertex $x$ of $\graph{U}$, there exists a path from $x$ to the only element of $\faces{}{+}U$.
	Since every wire vertex is the source of at most one edge by Proposition 
	\ref{prop:graph_of_molecule_properties}, and every node vertex is by assumption the source of exactly one edge, this path is necessarily unique.
	We conclude that $\graph{U}$ is an in-tree whose root is the only output face of $r$.
\end{proof}

\begin{prop} \label{prop:positive_opetopes_are_many_to_one_molecules}
	Let $U$ be a molecule.
	The following are equivalent:
	\begin{enumerate}[label=(\alph*)]
		\item $U$ is a positive opetope tree;
		\item $U$ is round and, for all $x \in U$, if $\dim{x} > 0$, then $\size{\faces{}{+}x} = 1$.
	\end{enumerate}
\end{prop}
\begin{proof}
	One implication is part of Lemma \ref{lem:positive_opetope_basic_properties}.
	Conversely, we proceed by induction on $n \eqdef \dim{U}$.
	If $n = 0$, then $U$ is the point, and we are done.
	
	Let $n > 0$, and suppose $U$ is an atom.
	Then $U$ is isomorphic to $V \celto W$ for some round molecules $V$ and $W$ such that $\bound{}{\alpha}V$ is isomorphic to $\bound{}{\alpha}W$ for all $\alpha \in \set{+, -}$.
	Since $V$ also satisfies the assumptions, by the inductive hypothesis $V$ is a positive opetope tree, and since $\size{\maxel{W}} = \size{\faces{}{+}U} = 1$, $W$ is isomorphic to $\compos{V}$.
	It follows that $U$ is isomorphic to $V \celto \compos{V}$, so it is a positive opetope.

	Finally, suppose that $U$ is not an atom, and let $(\order{i}{U})_{i=1}^m$ be an $(n-1)$\nbd layering of $U$, with associated $(n-1)$\nbd ordering $(\order{i}{x})_{i=1}^m$.
	We will identify each layer with its image in $U$.
	For each $i \in \set{0, \ldots, m-1}$, let
	\[
		\order{i}{V} \eqdef \order{m-i}{U} \cp{n-1} \ldots \cp{n-1} \order{m}{U}.
	\]
	We will prove that each $\order{i}{V}$ is a positive opetope tree by recursion on $i$.
	By Lemma \ref{lem:molecules_whose_graph_is_in_tree}, $\flow{n-1}{U}$ is an in-tree, so $\order{m}{x}$ is necessarily its root.
	It follows that $\bound{}{+}U = \bound{}{+}\order{m}{x}$, so $\order{0}{V} = \order{m}{U} = \clset{\order{m}{x}}$, which is a positive opetope.
	Let $i > 0$.
	Then by the dual of Proposition \ref{prop:layering_from_ordering}
	\[
		W \eqdef \bound{}{+}\order{m-i}{x} \submol \bound{}{-}\order{m-i+1}{U} = \bound{}{-}\order{i-1}{V},
	\]
	$\order{i-1}{V}$ is a positive opetope tree by the inductive hypothesis on $i$, and by what we have already proven $\clset{\order{m-i}{x}}$ is a positive opetope. 
	It follows that $\order{i}{V}$ is isomorphic to $\clset{\order{m-i}{x}} \cpsub{W}{-} \order{i-1}{V}$, which is a positive opetope tree.
	Since $U = \order{m-i}{V}$, we conclude.
\end{proof}

\begin{comm}
	Proposition \ref{prop:positive_opetopes_are_many_to_one_molecules} may be read in analogy with Zawadowski's result \cite[Proposition 13.4]{zawadowski2007positive} relating presheaves on positive opetopes to ``positive-to-one'' polygraphs: positive opetope trees are the ``positive-to-one'' round molecules, and positive opetopes are the ``positive-to-one'' atoms.
\end{comm}

\begin{cor} \label{cor:faces_and_boundaries_of_positive_opetope_trees}
	Let $U$ be a positive opetope tree, $n \in \mathbb{N}$, $\alpha \in \set{+, -}$, and $x \in U$.
	Then
	\begin{enumerate}
		\item $\bound{n}{\alpha}U$ is a positive opetope tree,
		\item $\clset{x}$ is a positive opetope.
	\end{enumerate}
\end{cor}
\begin{proof}
	Both $\bound{n}{\alpha}U$ and $\clset{x}$ are round molecules satisfying the conditions of Proposition \ref{prop:positive_opetopes_are_many_to_one_molecules}.
\end{proof}

\begin{comm}
	If we are only interested in positive opetopes, as opposed to positive opetope trees, the results of \cite{leclerc2023poset} imply the following characterisation, strengthening 
	Proposition \ref{prop:positive_opetopes_are_many_to_one_molecules}: an oriented graded poset $U$ is a positive opetope if and only if
	\begin{enumerate}
		\item for all $x \in U$, if $\dim{x} > 0$, then $\size{\faces{}{+}x} = 1$,
		\item $U$ has a greatest element,
		\item $\augm{U}$ is an oriented thin graded poset, and
		\item for all $x \in U$, $\graph{(\bound{}{-}x)}$ is acyclic.
	\end{enumerate}
\end{comm}

\begin{dfn}[Zoom complex] \index{positive opetope!zoom complex}
	A \emph{zoom complex} is a finite sequence $(\grade{i}{\mathscr{T}})_{i=0}^n$ of directed graphs with open edges such that, for all $i \in \set{0, \ldots, n}$, $\grade{i}{\mathscr{T}}$ is an in-tree, together with, for all $i \in \set{1, \ldots, n}$, a bijection
	\[
		\grade{i}{\zeta}\colon \faces{}{-}\grade{i}{\mathscr{T}} \iso N_{\grade{i-1}{\mathscr{T}}}
	\]
	between the wire vertices in the input boundary of $\grade{i}{\mathscr{T}}$ and the node vertices of $\grade{i-1}{\mathscr{T}}$.
\end{dfn}

\begin{dfn}[Zoom complex of a positive opetope tree]
	Let $U$ be a positive opetope tree, $n \eqdef \dim{U}$.
	The \emph{zoom complex of $U$} is the sequence $(\graph{(\bound{i}{-}U)})_{i=0}^n$ of directed graphs with open edges, together with the bijections
	\[
		\grade{i}{\zeta}\colon \faces{}{-}\graph{(\bound{i}{-}U)} \iso \faces{}{-}(\bound{i}{-}U) = \faces{i-1}{-}U = N_{\graph{(\bound{i-1}{-}U)}}
	\]
	determined by Proposition \ref{prop:graph_of_molecule_properties} and globularity of $U$ for each $i \in \set{1, \ldots, n}$.
\end{dfn}

\begin{prop} \label{prop:zoom_complex_of_positive opetope}
	Let $U$ be a positive opetope tree.
	Then the zoom complex of $U$ is well-defined as a zoom complex.	
\end{prop}
\begin{proof}
	It suffices to show that $\graph{(\bound{i}{-}U)}$ is an in-tree for all $i \in \set{0, \ldots, n}$, which follows from Lemma \ref{lem:positive_opetope_basic_properties} and Lemma 
	\ref{lem:molecules_whose_graph_is_in_tree}.
\end{proof}

\begin{lem} \label{lem:faces_of_positive_opetopes}
	Let $U$ be a positive opetope tree, $x \in U$, $n \eqdef \dim{x}$.
	Then either
	\begin{itemize}
		\item $x \in \faces{n}{-}U$, or
		\item there exists $y \in \cofaces{}{+}x \cap \faces{n+1}{-}U$.
	\end{itemize}
\end{lem}
\begin{proof}
	Suppose that $x \not\in \faces{n}{-}U$.
	Then there exists $y \in \cofaces{}{+}x$.
	We construct a sequence $(y_j)_{j \in \mathbb{N}}$ of elements of $\cofaces{}{+}x$ with the property that, for all $j$, there exists a path from $y_{j+1}$ to $y_j$ in $\hasseo{U}$.
	We start by letting $y_0 \eqdef y$.
	For each $j \in \mathbb{N}$, if $y_j \in \faces{n+1}{-}U$, then we let $y_k \eqdef y_j$ for all $k > j$.
	Otherwise, there exists $z \in \cofaces{}{+}y_j$.
	Since $U$ is a molecule, by Proposition \ref{prop:oriented_diamond_rdc} the interval $[x, z]$ is of the form
\[\begin{tikzcd}
	& z \\
	y_j && {y'} \\
	& x
	\arrow["{+}"', from=1-2, to=2-1]
	\arrow["{+}"', from=2-1, to=3-2]
	\arrow["\alpha", from=1-2, to=2-3]
	\arrow["{-\alpha}", from=2-3, to=3-2]
\end{tikzcd}\]
	for a unique $y' \in \grade{n+1}{U}$ and $\alpha \in \set{+, -}$.
	Because $\size{\faces{}{+}z} = 1$, necessarily $\alpha = -$, so there is a path $y' \to z \to y_j$ in $\hasseo{U}$, and we let $y_{j+1} \eqdef y'$.
	Since $\hasseo{U}$ is finite and acyclic by Lemma \ref{lem:positive_opetope_basic_properties}, there exists $j \in \mathbb{N}$ such that $y_j = y_{j+1}$, hence $y_j \in \cofaces{}{+} x \cap \faces{n+1}{-}U$.
\end{proof}

\begin{prop} \label{prop:positive_opetope_reconstruction_zoom_complex}
	Let $U$ be a positive opetope tree.
	Then $U$ can be uniquely reconstructed from the data of its zoom complex.
\end{prop}
\begin{proof}
	Let $x \in U$, $i \eqdef \dim{x}$.
	By Lemma \ref{lem:faces_of_positive_opetopes}, either $x \in \faces{i}{-}U$, in which case $x$ appears as a node vertex in $\graph{(\bound{i}{-}U)}$, or there exists $y \in \cofaces{}{+}x \cap \faces{i+1}{-}U$, in which case $x$ appears as a wire vertex in $\graph{(\bound{i+1}{-}U)}$.
	It follows that we can reconstruct the underlying set of $U$ by taking the amalgamated union of the vertices of all directed graphs in the zoom complex of $U$, identifying those that are related by one of the $\zeta_i$ bijections.
	By Proposition \ref{prop:data_for_ogposets}, it then suffices to show that we can reconstruct the functions $\faces{}{+}$, $\faces{}{-}$.
	
	We proceed by induction on $n \eqdef \dim{U}$.
	If $n \leq 1$, then the data of $\graph{U}$ with the separation into node and wire vertices is equivalent to the data of $\hasseo{U}$ together with the function $\dim$.
	Suppose that $n > 1$.
	Then the truncated zoom complex $(\graph{(\bound{i}{-}U)})_{i=0}^{n-1}$ is the zoom complex of $\bound{}{-}U$, so we may assume, by the inductive hypothesis, that we have reconstructed $\bound{}{-}U$ as an oriented graded poset.

	From $\graph{U}$ we can reconstruct $\flow{n-1}{U}$ and pick an $(n-1)$\nbd ordering $(\order{i}{x})_{i=1}^m$ of $U$.
	Because $U$ is acyclic, it is frame-acyclic by Proposition 
	\ref{prop:acyclicity_implications} and Proposition 
	\ref{prop:dimensionwise_implies_frame_acyclic}, so by Corollary
	\ref{cor:frame_acyclicity_equivalent_conditions} and Proposition 
	\ref{prop:layering_from_ordering} the sequence
	\[
		\order{0}{U} \eqdef \bound{}{-}U, \quad \quad 
		\order{i}{U} \eqdef \bound{n-1}{+}\order{i-1}{U} \cup \clset{{\order{i}{x}}} \quad \text{for $i \in \set{1, \ldots, m}$}
	\]
	determines an $(n-1)$\nbd layering $(\order{i}{U})_{i=1}^m$ of $U$ with $\bound{}{-}\order{i}{x} \submol \bound{}{-}\order{i}{U}$ for all $i \in \set{1, \ldots, m}$.
	We will prove, inductively on $i \in \set{0, \ldots, m}$, that we can reconstruct $\order{i}{U}$ as an oriented graded poset, from which we reconstruct $U$ as $\bigcup_{i=1}^m \order{i}{U}$.
	The case $i = 0$ holds by the inductive hypothesis on dimension, so suppose that $i > 0$.
	By the inductive hypothesis on $i$, we may assume we have reconstructed $\bound{n-1}{+}\order{i-1}{U} = \bound{}{-}\order{i}{U}$.
	Let 
	\[
		x \in \order{i}{U} \setminus \bound{}{-}\order{i}{U} = \clset{\order{i}{x}} \setminus \bound{}{-}\order{i}{x} = \set{\order{i}{x}} \cup \faces{}{+}\order{i}{x}.
	\]
	If $x = \order{i}{x}$, in which case $x$ appears as a node vertex in $\graph{U}$, we have
	\[
		\faces{}{-}x = \set{y \mid y \in s(x)}, \quad \quad \faces{}{+}x = \set{y \mid y \in t(x)}
	\]
	where $s, t$ are the source and target functions of $\graph{U}$.
	Otherwise, $x$ is the only output face of $\order{i}{x}$, and $\bound{}{-}\order{i}{x}$ is known as a closed subset of $\bound{}{-}\order{i}{U}$.
	Then
	\[
		\faces{}{\alpha}x = \faces{}{\alpha}(\bound{}{-}\order{i}{x})
	\]
	for all $\alpha \in \set{+, -}$.
	This completes the reconstruction of $U$.
\end{proof}

\begin{comm}
	Proposition \ref{prop:positive_opetope_reconstruction_zoom_complex} implies that positive opetope trees can be identified with a certain subclass of zoom complexes.
	This is the route taken in \cite{kock2010polynomial}, whose authors characterised a class of zoom complexes corresponding to general, not necessarily positive opetopes.
	It seems likely, but we will not prove, that the class of positive opetopes can be obtained by restricting to those zoom complexes in this class whose node vertices in degrees $> 0$ are targets of at least one edge.
\end{comm}

\begin{exm}[The zoom complex of a positive opetope]
	Let $U$ be the 3\nbd dimensional atom whose input and output boundary are the oriented face posets of 
	\[\begin{tikzcd}[column sep=small]
	{{\scriptstyle 0} \;\bullet} &&& {{\scriptstyle 3} \;\bullet} \\
	& {{\scriptstyle 1} \;\bullet} & {{\scriptstyle 2} \;\bullet}
	\arrow["0"', curve={height=6pt}, from=1-1, to=2-2]
	\arrow[""{name=0, anchor=center, inner sep=0}, "3", curve={height=-12pt}, from=2-2, to=2-3]
	\arrow["2"', curve={height=6pt}, from=2-3, to=1-4]
	\arrow[""{name=1, anchor=center, inner sep=0}, "5", curve={height=-18pt}, from=1-1, to=1-4]
	\arrow[""{name=2, anchor=center, inner sep=0}, "1"', curve={height=6pt}, from=2-2, to=2-3]
	\arrow[""{name=3, anchor=center, inner sep=0}, "4", curve={height=-18pt}, from=1-1, to=2-3]
	\arrow["0", curve={height=-6pt}, shorten >=4pt, Rightarrow, from=2-2, to=3]
	\arrow["2"', curve={height=6pt}, shorten >=6pt, Rightarrow, from=2-3, to=1]
	\arrow["1"', shorten <=2pt, shorten >=2pt, Rightarrow, from=2, to=0]
\end{tikzcd}
\quad \text{and} \quad
\begin{tikzcd}[column sep=small]
	{{\scriptstyle 0} \;\bullet} &&& {{\scriptstyle 3} \;\bullet} \\
	& {{\scriptstyle 1} \;\bullet} & {{\scriptstyle 2} \;\bullet}
	\arrow["0"', curve={height=6pt}, from=1-1, to=2-2]
	\arrow["2"', curve={height=6pt}, from=2-3, to=1-4]
	\arrow[""{name=0, anchor=center, inner sep=0}, "5", curve={height=-18pt}, from=1-1, to=1-4]
	\arrow[""{name=1, anchor=center, inner sep=0}, "1"', curve={height=6pt}, from=2-2, to=2-3]
	\arrow["4"', shorten <=7pt, shorten >=7pt, Rightarrow, from=1, to=0]
\end{tikzcd}\]
	respectively.
	Since all its elements have a single output face, by Proposition 
	\ref{prop:positive_opetopes_are_many_to_one_molecules}, $U$ is a positive opetope.
	The zoom complex of $U$ is the sequence
\[
	\input{img/zoom_complex_0.tex} \quad
	\input{img/zoom_complex_1.tex}
\]
\[
	\input{img/zoom_complex_2.tex} \quad
	\input{img/zoom_complex_3.tex}
\]
	of directed graphs with open edges, all of which are in-trees, pictured as string diagrams.
\end{exm}

\begin{prop} \label{prop:inclusions_of_positive_opetope_trees_are_submolecule}
	Let $\imath\colon V \incl U$ be an inclusion of molecules such that $\dim{U} = \dim{V}$, $U$ is a positive opetope tree, and $V$ is round.
	Then
	\begin{enumerate}
		\item $V$ is a positive opetope tree,
		\item $\imath$ is a submolecule inclusion.
	\end{enumerate}
\end{prop}
\begin{proof}
	First of all, since $V$ is a round molecule embedding into $U$, it satisfies the conditions of Proposition \ref{prop:positive_opetopes_are_many_to_one_molecules}, so it is a positive opetope tree.
	Let $n \eqdef \dim{U} = \dim{V}$.
	By Lemma \ref{lem:molecules_whose_graph_is_in_tree}, both $\flow{n-1}{U}$ and $\flow{n-1}{V}$ are are in-trees.
	Then the induced subgraph of $\flow{n-1}{U}$ on the vertices in the image of $V$ is necessarily path-induced: given two vertices $y, y'$ in $\flow{n-1}{V}$, there exists a unique path from each of them to the root $r$ of $\flow{n-1}{V}$.
	If there existed a path from $\imath(y)$ to $\imath(y')$ in $\flow{n-1}{U}$ which is not the image of a path in $\flow{n-1}{V}$, then there would be two different paths in $\flow{n-1}{U}$ from $\imath(y)$ to $\imath(r)$, a contradiction since $\flow{n-1}{U}$ is an in-tree.

	It follows from Lemma \ref{lem:connected_subgraph_conditions_path_induced} that there is an $(n-1)$\nbd ordering of $U$ in which the vertices in the image of $V$ are consecutive.
	Moreover, both $U$ and $V$ are acyclic, so they are frame-acyclic, and by Corollary 
	\ref{cor:frame_acyclicity_equivalent_conditions} this $(n-1)$\nbd ordering is induced by an $(n-1)$\nbd layering.
	Finally, by Lemma \ref{lem:positive_opetope_basic_properties} $\bound{}{+}V$ is an atom, so by Lemma \ref{lem:downset_is_submolecule} it is a submolecule of \emph{any} molecule that contains it.
	We conclude by the dual of Lemma \ref{lem:round_submolecules_from_layering}.
\end{proof}

\clearpage
\thispagestyle{empty}

%% file: img/zoom_complex_0.tex
\begin{tikzpicture}[xscale=5, yscale=4, baseline={([yshift=-.5ex]current bounding box.center)}]
\path[fill=white] (0, 0) rectangle (1, 1);
\node[circle, fill=black, draw=black, inner sep=1pt] at (0.5, 0.5) {};
\node[text=black, font={\scriptsize \sffamily}, xshift=12pt, yshift=0pt] at (0.5, 0.5) {(0, 0)};
\end{tikzpicture}

%% file: geometric.tex
\chapter{Geometric realisation} \label{chap:geometric}
\thispagestyle{firstpage}

\begin{guide}
	From the very start, we have justified our focus on regular directed complexes with the claim that they are a directed version of regular cell complexes.
	The purpose of this chapter is to prove that those were not empty words: there is a geometric realisation functor which turns each map of regular directed complexes into a cellular map of \cemph{regular CW complexes}.
	The generating $n$\nbd dimensional cells of the geometric realisation of a regular directed complex are indexed by its $n$\nbd dimensional elements.

	This is the content of Theorem \ref{thm:regular_directed_complex_is_regular_cw_complex}, which is almost an immediate consequence of Proposition \ref{prop:order_complex_of_round_molecule}, stating that the \cemph{order complex} of the underlying poset of a round $n$\nbd dimensional molecule is a \cemph{PL $n$\nbd ball}, and the order complex of its boundary is a \cemph{PL $(n-1)$\nbd sphere}.
	We also promptly derive that Gray products, suspensions, and joins of regular directed complexes are mapped to products, suspensions, and joins of their geometric realisations, up to cellular homeomorphism (Proposition \ref{prop:realisation_of_operations}).

	This is not a particularly complicated proof, nor does it use much of what we have developed since Chapter \ref{chap:layerings}.
	However, it does require setting up a lot of notation and terminology from combinatorial and poset topology.
	Much of what we use is very classical, and we have avoided going into too much detail; \cite{bjorner1995topological} and \cite{wachs2007poset} are excellent starting points for the reader who is interested in learning more.
\end{guide}


\section{Elements of poset topology} \label{sec:poset_topology}

\begin{guide}
	In this section, we introduce some basic notions of poset topology, such as the \emph{order complex} of a poset and its \emph{geometric realisation}, as well as a few elementary definitions from piecewise linear (PL) topology.

	We then state one of the foundational results of PL topology (Proposition 
	\ref{prop:gluing_pl_balls_and_spheres}): gluing two PL $n$\nbd balls at their boundaries produces a PL $n$\nbd sphere if and only if their entire boundaries are glued, and a PL $n$\nbd ball if the glued portion is a PL $(n-1)$\nbd ball.
	Surprisingly, these intuitive facts are false in the absence of the PL assumption, or another taming constraint. 
\end{guide}

\begin{dfn}[Simplicial set] \index{simplicial set}
	A \emph{simplicial set} is a presheaf on the simplex category.
\end{dfn}

\begin{dfn}[The category $\sset$] \index{$\sset$}
	We let $\sset$ denote the category of simplicial sets and morphisms of presheaves.
\end{dfn}

\begin{dfn}[Simplicial subset] \index{simplicial set!subset}
	Let $K$ be a simplicial set.
	A \emph{simplicial subset} $L \subseteq K$ is a simplicial set $L$ together with a monomorphism $L \incl K$ whose components are subset inclusions $L[n] \subseteq K[n]$ for all $n \in \mathbb{N}$.
\end{dfn}

\begin{dfn}[Simplex in a simplicial set] \index{simplex!in a simplicial set}
	Let $K$ be a simplicial set, $n \in \mathbb{N}$.
	An \emph{$n$\nbd simplex in $K$} is an element $x \in K[n]$.
	We will identify $[n]$ with its Yoneda embedding into $\sset$, and also write $x\colon [n] \to K$ for an $n$\nbd simplex $x$ in $K$.
\end{dfn}

\begin{dfn}[Non-degenerate simplex in a simplicial set] \index{simplex!in a simplicial set!non-degenerate} \index{$\nondeg{K}$}
	Let $K$ be a simplicial set and $x \colon [n] \to K$.
	We say that $x$ is \emph{non-degenerate} if, given any surjective map $p\colon [n] \to [m]$ and $y\colon [m] \to K$, if $x = y \after p$, then $y = x$ and $p = \idd{[n]}$.
	We write $\nondeg{K}$ for the set of non-degenerate simplices in $K$.
\end{dfn}

\begin{lem} \label{lem:eilenberg_zilber_lemma}
	Let $K$ be a simplicial set, $x\colon [n] \to K$.
	Then there exists a unique pair $(p\colon [n] \to [m], y\colon [m] \to K)$ such that
	\begin{enumerate}
		\item $p$ is surjective,
		\item $y$ is non-degenerate,
		\item $x = y \after p$.
	\end{enumerate}
\end{lem}
\begin{proof}
	See \cite[Section II.3]{gabriel1967calculus}.
\end{proof}

\begin{dfn}[Finite simplicial set] \index{simplicial set!finite}
	A simplicial set is \emph{finite} if $\size{\nondeg{K}}$ is finite.
\end{dfn}

\begin{dfn}[Dimension of a simplex] \index{simplex!in a simplicial set!dimension}
	Let $K$ be a simplicial set, $x$ a simplex in $K$, and let $(p, y\colon [m] \to K)$ be the unique pair of a surjective map and a non-degenerate simplex in $K$ such that $x = y \after p$.
	The \emph{dimension of $x$} is the natural number $\dim{x} \eqdef m$.
\end{dfn}

\begin{dfn}[Order complex of a poset] \index{poset!order complex}
	Let $P$ be a poset.
	The \emph{order complex of $P$} is the simplicial set $\ordcpx{P}$ whose $n$\nbd simplices are chains $c\colon [n] \to P$ and maps in $\simplexcat \incl \poscat$ act by precomposition.
\end{dfn}

\begin{lem} \label{lem:order_complex_is_a_functor}
	Let $f\colon P \to Q$ be an order-preserving map of posets.
	Then
	\begin{align*}
		\ordcpx{f}\colon \ordcpx{P} & \to \ordcpx{Q}, \\
		(c\colon [n] \to P) & \mapsto (f \after c\colon [n] \to Q)
	\end{align*}
	is a morphism of simplicial sets.
	This determines a functor $\ordcpx{-}\colon \poscat \to \sset$.
\end{lem}
\begin{proof}
	Straightforward.
\end{proof}

\begin{lem} \label{lem:order_complex_limits}
	The functor $\ordcpx{-}\colon \poscat \to \sset$ preserves all limits.
\end{lem}
\begin{proof}
	Follows from the fact that $\ordcpx{-}$ is a right adjoint functor, part of a nerve-realisation pair together with the left Kan extension of the inclusion $\simplexcat \incl \poscat$ along the Yoneda embedding $\simplexcat \incl \sset$.	
\end{proof}

\begin{lem} \label{lem:order_complex_preserves_injective}
	Let $\imath\colon P \to Q$ be an injective order-preserving map of posets.
	Then $\ordcpx{\imath}\colon \ordcpx{P} \to \ordcpx{Q}$ is a monomorphism.
\end{lem}
\begin{proof}
	Let $n \in \mathbb{N}$ and let $c, c'\colon [n] \to P$ be chains.
	Since $\imath$ is injective, $\imath \after c = \imath \after c'$ implies $c = c'$.
	It follows that $\ordcpx{\imath}[n]\colon \ordcpx{P}[n] \to \ordcpx{Q}[n]$ is injective, and monomorphisms of presheaves are precisely the componentwise injective morphisms.
\end{proof}

\begin{lem} \label{lem:order_complex_pushouts}
	The functor $\ordcpx{-}\colon \poscat \to \sset$ preserves pushouts of closed embeddings along closed embeddings.
	The image of a pushout square of closed embeddings is both a pushout and a pullback square of monomorphisms.
\end{lem}
\begin{proof}
	Consider a pushout square of closed embeddings
	\[\begin{tikzcd}
	P \cap Q && Q \\
	P && P \cup Q
	\arrow[hook', from=1-1, to=2-1]
	\arrow[hook, from=2-1, to=2-3]
	\arrow[hook, from=1-1, to=1-3]
	\arrow[hook', from=1-3, to=2-3]
	\arrow["\lrcorner"{anchor=center, pos=0.125, rotate=180}, draw=none, from=2-3, to=1-1]
\end{tikzcd}\]
	in $\poscat$.
	By Lemma \ref{lem:pushouts_of_spans_of_closed_embeddings}, this is also a pullback square, so by Lemma \ref{lem:order_complex_limits} its image through $\ordcpx{-}$ is a pullback square, and by Lemma \ref{lem:order_complex_preserves_injective} it is a square of monomorphisms, so it only remains to show that it is a pushout square.
	
	Let $n \in \mathbb{N}$, let $c\colon [n] \to P \cup Q$ be a chain, and let $x \eqdef c(n)$.
	Then the image of $c$ is entirely contained in $\clset{x}$.
	We have $x \in P$ or $x \in Q$; suppose without loss of generality that $x \in P$.
	Then $c$ factors through $P \incl P \cup Q$.
	Moreover, $c$ also factors through $P \incl P \cup Q$ if and only if $x \in P \cap Q$.
	It follows that 	
\[\begin{tikzcd}
	\ordcpx{(P \cap Q)}[n] && \ordcpx{Q}[n] \\
	\ordcpx{P}[n] && \ordcpx{(P \cup Q)}[n]
	\arrow[hook', from=1-1, to=2-1]
	\arrow[hook, from=2-1, to=2-3]
	\arrow[hook, from=1-1, to=1-3]
	\arrow[hook', from=1-3, to=2-3]
	\arrow["\lrcorner"{anchor=center, pos=0.125, rotate=180}, draw=none, from=2-3, to=1-1]
\end{tikzcd}\]
	is a pushout square of sets and functions.
	Because colimits in presheaf categories are computed pointwise, we conclude.
\end{proof}

\begin{dfn}[Ordered simplicial complex] \index{simplicial set!ordered simplicial complex}
	An \emph{ordered simplicial complex} is a simplicial set $K$ with the property that, for all $x, y\colon [n] \to K$, if $x\after \imath = y\after \imath$ for all $\imath\colon [0] \incl [n]$, then $x = y$.
\end{dfn}

\begin{lem} \label{lem:order_complex_is_ordered_simplicial_complex}
	Let $P$ be a poset.
	Then 
	\begin{enumerate}
		\item $\ordcpx{P}$ is an ordered simplicial complex,
		\item the non-degenerate simplices in $\ordcpx{P}$ are the injective chains in $P$,
		\item if $P$ is finite, then $\ordcpx{P}$ is finite.
	\end{enumerate}
\end{lem}
\begin{proof}
	Let $c, c'$ be $n$\nbd simplices in $\ordcpx{P}$, suppose $c \after \imath = c' \after \imath$ for all $\imath\colon [0] \incl [n]$, and let $k \in [n]$.
	There is an injection $[0] \incl [n]$ sending $0$ to $k$, so $c(k) = c'(k)$.
	It follows that $c = c'$ as functions, so $c = c'$ as $n$\nbd simplices.

	Next, let $c\colon [n] \to P$ be a chain.
	Then $c$ factors uniquely up to unique isomorphism as a surjective order-preserving map $\widehat{c}\colon [n] \to c(P)$ followed by an injective order-preserving map $\imath\colon c(P) \incl P$.
	Moreover, the image of $c$ is a finite linear order, so up to isomorphism we can take $c(P) = [m]$ for $m \eqdef \size{c}$.
	If $c$ is already injective, then uniqueness of the factorisation implies that $m = n$, $\widehat{c} = \idd{[n]}$, and $\imath = c$, so $c$ is non-degenerate as a simplex.
	Conversely, if $c$ is not injective, then $c = \imath \after \widehat{c}$ is a non-trivial factorisation, and $c \not\in \nondeg{\ordcpx{P}}$.

	Finally, if $P$ is finite, the number of injective chains in $P$ is bounded by the size of the power set of $P$, which is finite, so $\size{\nondeg{\ordcpx{P}}}$ is finite.
\end{proof}

\begin{comm}
	The order complex is more commonly seen as an \emph{abstract} (unordered) simplicial complex in topological combinatorics.
	When seen as a simplicial set, the order complex is also known as the \emph{nerve} of a poset.
\end{comm}

\begin{dfn}[The category $\ktop$] \index{$\ktop$} \index{space}
	We let $\ktop$ denote the category of compactly generated, weakly Hausdorff topological spaces and continuous maps.
\end{dfn}

\begin{comm}
	We will simply say \emph{space} for an object of $\ktop$.
	We pick $\ktop$ as one of several possible ``convenient'' categories of topological spaces, which have all small limits and colimits, are cartesian closed, contain all CW complexes, and are such that the categorical product of two CW complexes is a CW complex.
	We refer to \cite{strickland2009category} for a survey of its properties.
\end{comm}

\begin{dfn}[Standard geometric simplex] \index{simplex!standard geometric} \index{$\geosim{n}$}
	Let $n \in \mathbb{N}$.
	The \emph{standard geometric $n$\nbd simplex} is the subspace
	\[
		\geosim{n} \eqdef \set{(x_0, x_1, \ldots, x_n) \mid \sum_{i=0}^n x_i = 1} \subseteq \mathbb{R}^{n+1}.
	\]
\end{dfn}

\begin{dfn}[Boundary of the standard geometric simplex] \index{simplex!standard geometric!boundary} \index{$\bound{}{}\geosim{n}$} \index{boundary!of a standard geometric simplex}
	Let $n \in \mathbb{N}$.
	The \emph{boundary of the standard geometric $n$\nbd simplex} is the subspace
	\[
		\bound{}{}\geosim{n} \eqdef \set{(x_0, x_1, \ldots, x_n) \mid \text{$x_i = 0$ for some $i \in \set{0, \ldots, n}$}} \subseteq \geosim{n}.
	\]
\end{dfn}

\begin{rmk}
	For all $n \in \mathbb{N}$, the standard geometric $n$\nbd simplex is a closed $n$\nbd ball, and its boundary is an $(n-1)$\nbd sphere.
\end{rmk}

\begin{dfn}[Geometric realisation of simplicial sets] \index{simplicial set!geometric realisation}
	Let $\realis{-}\colon \simplexcat \to \ktop$ be the functor which sends an order-preserving map $f\colon [n] \to [m]$ to the continuous map $\realis{f}\colon \geosim{n} \to \geosim{m}$ defined by
	\[
		(x_0, \ldots, x_n) \mapsto (x'_0, \ldots, x'_m), \quad \quad x'_j \eqdef \sum_{i \in \invrs{f}(j)} x_i.
	\]
	The \emph{geometric realisation of simplicial sets} is the functor
	\[
		\realis{-}\colon \sset \to \ktop
	\]
	defined as the left Kan extension of $\realis{-}\colon \simplexcat \to \ktop$ along the Yoneda embedding $\simplexcat \incl \sset$.
\end{dfn}

\begin{dfn}[PL map] \index{PL!map}
	Let $K$ be a simplicial set and let $U \subseteq \mathbb{R}^n$ be a subspace of Euclidean space.
	A continuous map $f\colon \realis{K} \to U$ is \emph{PL}, short for \emph{piecewise linear}, if, for all $x\colon [m] \to K$, the composite $f \after \realis{x}\colon \geosim{m} \to U$ extends to a linear map $\mathbb{R}^{m+1} \to \mathbb{R}^n$.
	A \emph{PL homeomorphism} is a PL map which is a homeomorphism.
\end{dfn}

\begin{dfn}[PL ball] \index{PL!ball}
	Let $n \in \mathbb{N}$.
	A finite ordered simplicial complex $K$ is a \emph{PL $n$\nbd ball} if there exists a PL homeomorphism $\realis{K} \iso \geosim{n}$.
\end{dfn}

\begin{dfn}[PL sphere] \index{PL!sphere}
	Let $n \in \mathbb{N}$.
	A finite ordered simplicial complex $K$ is a \emph{PL $n$\nbd sphere} if there exists a PL homeomorphism $\realis{K} \iso \bound{}{}\geosim{n+1}$.
\end{dfn}

\begin{dfn}[Link of a simplex] \index{simplex!in a simplicial set!link} \index{$\lnk{K}{x}$} 
	Let $K$ be a simplicial set and let $x\colon [n] \to K$ be a simplex in $K$.
	The \emph{link of $x$} is the simplicial subset $\lnk{K}{x} \subseteq K$ whose simplices $y\colon [m] \to \lnk{K}{x}$ are those $y\colon [m] \to K$ such that
	\begin{enumerate}
		\item $x$ and $y$ are disjoint, that is,
			\[\begin{tikzcd}
	\varnothing && {[m]} \\
	{[n]} && K
	\arrow[hook', from=1-1, to=2-1]
	\arrow["y", from=1-3, to=2-3]
	\arrow[hook, from=1-1, to=1-3]
	\arrow["x", from=2-1, to=2-3]
	\arrow["\lrcorner"{anchor=center, pos=0.125}, draw=none, from=1-1, to=2-3]
\end{tikzcd}\]
		is a pullback square in $\sset$,
	\item there exists $z\colon [\ell] \to K$ and a cospan $f\colon [n] \to [\ell]$, $g\colon [m] \to [\ell]$ such that $x = z \after f$ and $y = z \after g$.
	\end{enumerate}
\end{dfn}

\begin{lem} \label{lem:links_in_pl_balls_and_spheres}
	Let $n \in \mathbb{N}$, let $K$ be a finite ordered simplicial complex, $x$ a simplex in $K$, and $m \eqdef \dim{x}$.
	Then 
	\begin{enumerate}
		\item if $K$ is a PL $n$\nbd ball, then $\lnk{K}{x}$ is either a PL $(n-m-1)$\nbd ball or a PL $(n-m-1)$\nbd sphere,
		\item if $K$ is a PL $n$\nbd sphere, then $\lnk{K}{x}$ is a PL $(n-m-1)$\nbd sphere.
	\end{enumerate}
\end{lem}
\begin{proof}
	See \cite[Corollary 1.16 and Lemma 1.17]{hudson1969piecewise}.
\end{proof}

\begin{dfn}[Boundary of a PL ball] \index{PL!ball!boundary} \index{boundary!of a PL ball}
	Let $K$ be a PL ball.
	The \emph{boundary of $K$} is the simplicial subset $\bound{}{}K \subseteq K$ whose simplices are those $x$ such that $\lnk{K}{x}$ is a PL ball.
\end{dfn}

\begin{prop} \label{prop:gluing_pl_balls_and_spheres}
	Let $K$, $L$ be finite ordered simplicial complexes, $n \in \mathbb{N}$, and let
\[\begin{tikzcd}
	{K \cap L} && L \\
	K && {K \cup L}
	\arrow[hook', from=1-1, to=2-1]
	\arrow[hook, from=1-1, to=1-3]
	\arrow[hook, from=2-1, to=2-3]
	\arrow[hook', from=1-3, to=2-3]
	\arrow["\lrcorner"{anchor=center, pos=0.125, rotate=180}, draw=none, from=2-3, to=1-1]
\end{tikzcd}\]
	be a pushout diagram of monomorphisms in $\sset$.
	Then
	\begin{enumerate}
		\item if $K$ and $L$ are PL $n$\nbd balls and $K \cap L$ is a PL $(n-1)$\nbd ball such that $K \cap L = \bound{}{}K \cap \bound{}{}L$, then $K \cup L$ is a PL $n$\nbd ball with
			\[
				\bound{}{}(K \cup L) = (\bound{}{}K \setminus L) \cup (\bound{}{}L \setminus K) \cup \bound{}{}(K \cap L),
			\]
		\item if $K$ is a PL $n$\nbd ball and $K \cap L = \bound{}{}K$, then $K \cup L$ is a PL $n$\nbd sphere if and only if $L$ is a PL $n$\nbd ball and $K \cap L = \bound{}{}L$.
	\end{enumerate}
\end{prop}
\begin{proof}
	See \cite[Theorem 2 and Theorem 3]{zeeman1966seminar}.
\end{proof}

\begin{comm}
	Starting from dimension 5, there are classical counterexamples of simplicial complexes whose realisation is homeomorphic, but not PL homeomorphic to a simplex or its boundary.
	Proposition \ref{prop:gluing_pl_balls_and_spheres} does not hold in general when the PL assumptions are dropped.
\end{comm}


\section{Face posets} \label{sec:face_posets}

\begin{guide}
	In this section, we recall the definitions of regular CW complex and face poset, and make precise the sense in which a regular CW complex is determined by its face poset
	(Proposition \ref{prop:ordcpx_right_inverse_to_face_poset} and Proposition 
	\ref{prop:face_poset_has_an_inverse_up_to_isomorphism}).
	
	For CW complexes, we purposely use the same notation and terminology as we used for polygraphs in Section \ref{sec:polygraphs}: polygraphs and CW complexes are conceptually analogous and play the same role in the canonical model structures on strict $\omega$\nbd categories and topological spaces, respectively \cite{lafont2010folk}.
	Regular directed complexes are a bridge between the two.
	We invite the reader to ponder why all regular directed complexes present regular CW complexes, but only those with frame-acyclic molecules present polygraphs.
\end{guide}

\begin{dfn}[Cellular extension of a space] \index{space!cellular extension}
	Let $X$ be a space.
	A \emph{cellular extension of $X$} is a space $X_{\gener{S}}$ together with a pushout diagram
\[
	\begin{tikzcd}
		{\coprod_{e \in \gener{S}} \bound{}{}B_e} && {\coprod_{e \in \gener{S}} B_e} \\
		X && {X_{\gener{S}}}
	\arrow["{(\bound{}{}e)_{e \in \gener{S}}}", from=1-1, to=2-1]
	\arrow["{(e)_{e \in \gener{S}}}", from=1-3, to=2-3]
	\arrow[hook, from=2-1, to=2-3]
	\arrow["{\coprod_{e \in \gener{S}}\imath_e}", hook, from=1-1, to=1-3]
	\arrow["\lrcorner"{anchor=center, pos=0.125, rotate=180}, draw=none, from=2-3, to=1-1]
\end{tikzcd}\]
in $\ktop$, where, for each $e \in \gener{S}$, the map $\imath_e \colon \bound{}{}B_e \incl B_e$ is the embedding of a sphere as the boundary of a closed ball.
\end{dfn}

\begin{dfn}[CW complex] \index{space!CW complex} \index{regular CW complex!generating cell} \index{$\cwcom{X}{S}$}
	A \emph{CW complex} is a space $X$ equipped with a sequential colimit cone of embeddings
\[\begin{tikzcd}
	{\varnothing \equiv\skel{-1}{X}} & {\skel{0}{X}} & {\skel{1}{X}} & \ldots & {\skel{n}{X}} & \ldots \\
	&&&&& X
	\arrow[hook, from=1-1, to=1-2]
	\arrow[hook, from=1-2, to=1-3]
	\arrow[hook, from=1-3, to=1-4]
	\arrow[hook, from=1-4, to=1-5]
	\arrow[hook, from=1-5, to=1-6]
	\arrow[hook, "{\imath_n}", from=1-5, to=2-6]
	\arrow[hook, "{\imath_1}", curve={height=6pt}, from=1-3, to=2-6]
	\arrow[hook, "{\imath_0}", curve={height=12pt}, from=1-2, to=2-6]
	\arrow[hook, "{\imath_{-1}}"', curve={height=18pt}, from=1-1, to=2-6]
\end{tikzcd}\]
and, for each $n \in \mathbb{N}$, a pushout diagram
\[
	\begin{tikzcd}
		{\coprod_{e \in \grade{n}{\gener{S}}} \bound{}{}B_e} && {\coprod_{e \in \grade{n}{\gener{S}}} B_e} \\
		\skel{n-1}{X} && \skel{n}{X}
	\arrow["{(\bound{}{}e)_{e \in \grade{n}{\gener{S}}}}", from=1-1, to=2-1]
	\arrow["{(e)_{e \in \grade{n}{\gener{S}}}}", from=1-3, to=2-3]
	\arrow[hook, from=2-1, to=2-3]
	\arrow["{\coprod_{e \in \grade{n}{\gener{S}}}\imath_e}", hook, from=1-1, to=1-3]
	\arrow["\lrcorner"{anchor=center, pos=0.125, rotate=180}, draw=none, from=2-3, to=1-1]
\end{tikzcd}\]
in $\ktop$, exhibiting $\skel{n}{X}$ as a cellular extension of $\skel{n-1}{X}$, such that $B_e$ is a closed $n$\nbd ball, hence $\bound{}{}B_e$ is an $(n-1)$\nbd sphere, for all $e \in \grade{n}{\gener{S}}$.
The set
\[
	\gener{S} \eqdef \sum_{n \in \mathbb{N}} \set{\imath_n \after e\colon B_e \to X \mid e \in \grade{n}{\gener{S}}}
\]
is called the set of \emph{generating cells} of the CW complex.
We write $\cwcom{X}{S}$ for a CW complex with set $\gener{S}$ of generating cells.
\end{dfn}

\begin{comm}
	For all $n \geq -1$, we will identify $\skel{n}{X}$ with its homeomorphic image in $X$, and given $e \in \grade{n}{\gener{S}}$ write $e\colon B_e \to X$ for $\imath_n \after e$.
\end{comm}

\begin{rmk}
All the data of a CW complex can be reconstructed uniquely up to unique homeomorphism from the pair $\cwcom{X}{S}$.
For example, for all $n \geq -1$, $\skel{n}{X}$ is homeomorphic to the subspace
\[
	\bigcup \set{e(B_e) \mid \text{$e \in \gener{S}$, $B_e$ is a closed $k$\nbd ball, $k \leq n$}}\subseteq X.
\]
\end{rmk}

\begin{dfn}[Interior of a ball]
	Let $B$ be a closed ball with boundary $\bound{}{}B$.
	The \emph{interior of $B$} is the subspace $\inter{B} \eqdef B \setminus \bound{}{}B$.
\end{dfn}

\begin{lem} \label{lem:stratification_of_cw_complex}
	Let $\cwcom{X}{S}$ be a CW complex.
	Then $\set{e(\inter{B_e}) \mid e \in \gener{S}}$ is a partition of $X$ into non-empty, pairwise disjoint open balls.
\end{lem}
\begin{proof}
	See \cite[Chapter II, Lemma 1.2]{lundell1969topology}.
\end{proof}

\begin{dfn}[Regular CW complex] \index{regular CW complex}
	A CW complex $\cwcom{X}{S}$ is \emph{regular} if every $e \in \gener{S}$ is an embedding.
\end{dfn}

\begin{lem} \label{lem:regular_cw_complex_is_normal}
	Let $\cwcom{X}{S}$ be a regular CW complex.
	For all $e, e' \in \gener{S}$, if $e(\inter{B_e}) \cap e'(B_{e'}) \neq \varnothing$, then $e(B_e) \subseteq e'(B_{e'})$.
\end{lem}
\begin{proof}
	See \cite[Chapter III, Theorem 2.1]{lundell1969topology}.
\end{proof}

\begin{dfn}[Face poset of a regular CW complex] \index{regular CW complex!face poset} \index{face poset!of a regular CW complex}
	Let $\cwcom{X}{S}$ be a regular CW complex.
	The \emph{face poset of $\cwcom{X}{S}$} is the poset $\fpos{\cwcom{X}{S}}$ whose
	\begin{itemize}
		\item underlying set is $\gener{S}$,
		\item partial order is defined by $e \leq e'$ if and only if $e(B_e) \subseteq e'(B_{e'})$.
	\end{itemize}
\end{dfn}

\begin{dfn}[Map of regular CW complexes] \index{regular CW complex!map} \index{map!of regular CW complexes}
	Let $\cwcom{X}{S}$, $\cwcom{Y}{T}$ be regular CW complexes.
	A \emph{map} $f\colon \cwcom{X}{S} \to \cwcom{Y}{T}$ is a pair of
	\begin{enumerate}
		\item a continuous map $f\colon X \to Y$,
		\item a function $\fpos{f}\colon \gener{S} \to \gener{T}$,
	\end{enumerate}
	such that, for all $e \in \gener{S}$,
	\[
		f(e(\inter{B_e})) = \fpos{f}(e)(\inter{B_{\fpos{f}(e)}}) \quad \text{and} \quad
		f(e(B_e)) = \fpos{f}(e)(B_{\fpos{f}(e)}).
	\]
\end{dfn}

\begin{comm}
	What we call a map is called a \emph{cellular regular continuous map} in \cite{lundell1969topology}.
\end{comm}

\begin{dfn}[The category $\rcpx$] \index{$\rcpx$}
	We let $\rcpx$ denote the category whose objects are regular CW complexes and morphisms are maps of regular CW complexes.
\end{dfn}

\begin{prop} \label{prop:cw_maps_basic_properties}
	Let $(f, \fpos{f})\colon \cwcom{X}{S} \to \cwcom{Y}{T}$ be a map of regular CW complexes.
	Then
	\begin{enumerate}
		\item $\fpos{f}$ is uniquely determined by $f$,
		\item $\fpos{f}$ determines a closed order-preserving map $\fpos{\cwcom{X}{S}} \to \fpos{\cwcom{Y}{T}}$.
	\end{enumerate}
	This assignment determines a functor $\fpos{}\colon \rcpx \to \posclos$.
\end{prop}
\begin{proof}
	Let $e \in \gener{S}$ and pick $x \in \inter{B_e}$, which is always non-empty.
	By Lemma \ref{lem:stratification_of_cw_complex}, there exists a unique $h \in \gener{T}$ such that $f(e(x)) \in h(\inter{B_h})$.
	Then 
	\[
		f(e(x)) \in h(\inter{B_{h}}) \cap f(e(\inter{B_e})) = h(\inter{B_{h}}) \cap \fpos{f}(e)(\inter{B_{\fpos{f}(e)}}),
	\]
	so necessarily $h = \fpos{f}(e)$.
	It follows that $\fpos{f}$ is uniquely determined by $f$.

	Next, suppose that $e \leq e'$ in $\fpos{\cwcom{X}{S}}$.
	Then
	\[
		\fpos{f}(e)(B_{\fpos{f}(e)}) = f(e(B_e)) \subseteq f(e'(B_{e'})) = \fpos{f}(e')(B_{\fpos{f}(e')}),
	\]
	so $\fpos{f}(e) \leq \fpos{f}(e')$, which proves that $\fpos{f}$ is order-preserving.

	Finally, let $e \in \gener{S}$, and suppose $h \leq \fpos{f}(e)$ in $\fpos{\cwcom{Y}{T}}$.
	Pick $y \in h(\inter{B_{h}})$.
	Because $y \in \fpos{f}(e)(B_{\fpos{f}(e)})$, which is equal to $f(e(B_e))$, there exists $x \in e(B_e)$ such that $f(x) = y$.
	Let $e' \in \gener{S}$ be the unique generating cell such that $x \in e'(\inter{B_{e'}})$.
	Then, as in the first part of the proof, we deduce that $\fpos{f}(e') = h$.
	Moreover, since $x \in e'(\inter{B_{e'}}) \cap e(B_e)$, it follows from Lemma \ref{lem:regular_cw_complex_is_normal} that $e' \leq e$.
	This proves that $\fpos{f}$ is closed.
\end{proof}

\begin{dfn}[CW poset] \index{poset!CW}
	A \emph{CW poset} is a poset $P$ with the property that, for all $x \in P$, $\realis{\ordcpx{(\clset{x} \setminus \set{x})}}$ is homeomorphic to a sphere.
\end{dfn}

\begin{lem} \label{lem:cwposet_basic_properties}
	Let $P$ be a CW poset, $x \in P$.
	Then
	\begin{enumerate}
		\item $P$ is graded,
		\item if $n \eqdef \dim{x}$, then $\realis{\ordcpx{\clset{x}}}$ is homeomorphic to a closed $n$\nbd ball and $\realis{\ordcpx{(\clos{\faces{}{}x})}} = \bound{}{}\realis{\ordcpx{\clset{x}}}$ is homeomorphic to an $(n-1)$\nbd sphere.
	\end{enumerate}
\end{lem}
\begin{proof}
	See \cite[Section 2]{bjorner1984posets}.
\end{proof}

\begin{dfn}[The category $\cwpos$] \index{$\cwpos$}
	We let $\cwpos$ denote the full subcategory of $\posclos$ on the CW posets.
\end{dfn}

\begin{comm}
	The original definition of CW poset in \cite{bjorner1984posets} corresponds to $\augm{P}$ for a non-empty CW poset $P$ according to our definition.
	Morphisms are accordingly required to reflect the least element.
	Of course, the two definitions are equivalent in light of Lemma \ref{lem:augmentation_has_an_inverse}.
\end{comm}

\begin{lem} \label{lem:face_poset_is_cwposet}
	Let $\cwcom{X}{S}$ be a regular CW complex.
	Then $\fpos{\cwcom{X}{S}}$ is a CW poset.
\end{lem}
\begin{proof}
	This is one side of \cite[Proposition 3.1]{bjorner1984posets}.
\end{proof}

\begin{cor} \label{cor:face_poset_functor}
	The functor $\fpos{}\colon \rcpx \to \posclos$ factors through the subcategory inclusion $\cwpos \incl \posclos$.
\end{cor}

\begin{prop} \label{prop:ordcpx_right_inverse_to_face_poset}
	Let $P$ be a CW poset.
	Then $\realis{\ordcpx{P}}$ admits a structure of regular CW complex whose set of generating cells is
	\[
		\set{\realis{\ordcpx{\clset{x}}} \incl \realis{\ordcpx{P}} \mid x \in P}.
	\]
	This assignment extends to a functor $\realis{\ordcpx{-}}\colon \cwpos \to \rcpx$, which is a section of $\fpos{}\colon \rcpx \to \cwpos$ up to natural isomorphism.
\end{prop}
\begin{proof}
	The fact that $\realis{\ordcpx{P}}$ admits a structure of regular CW complex whose face poset is isomorphic to $P$ is the other side of \cite[Proposition 3.1]{bjorner1984posets}.
	For functoriality, let $f\colon P \to Q$ be a closed order-preserving map of CW posets, $x \in P$, and $y \in Q$.
	Then
	\begin{enumerate}
		\item $\realis{\ordcpx{f}}(\realis{\ordcpx{\clset{x}}}) = \realis{\ordcpx{\clset{f(x)}}}$,
		\item $\invrs{\realis{\ordcpx{f}}}(\realis{\ordcpx{\clset{y}}}) = \bigcup_{x \in \invrs{f}(y)} \realis{\ordcpx{\clset{x}}}$
	\end{enumerate}
	follow straightforwardly from the definitions.

	Let $v \in \realis{\ordcpx{\clset{x}}}$.
	Then $v \in \bound{}{}\realis{\ordcpx{\clset{x}}}$ if and only if $v \in \realis{c}$ for some chain $c\colon [n] \to \clset{x}$ whose image is included in $\clos{\faces{}{}x}$, or, equivalently, such that $c(n) < x$.
	Suppose by way of contradiction that 
	\[
		v \in \inter{\realis{\ordcpx{\clset{x}}}} \quad \text{and} \quad \realis{\ordcpx{f}}(v) \in \bound{}{}\realis{\ordcpx{\clset{f(x)}}}.
	\]
	Then $\realis{\ordcpx{f}}(v) \in \realis{c}$ for some chain $c\colon [n] \to \clset{f(x)}$ such that $c(n) < f(x)$.
	It follows that 
	\[
		v \in \inter{\realis{\ordcpx{\clset{x}}}} \cap \bigcup_{x' \in \invrs{f}(c(n))} \realis{\ordcpx{\clset{x'}}},
	\]
	so there exists $x' \in P$ such that $f(x') = c(n)$ and $v \in \inter{\realis{\ordcpx{\clset{x}}}} \cap \realis{\ordcpx{\clset{x'}}}$.
	It follows from Lemma \ref{lem:regular_cw_complex_is_normal} and the first part of the proof that $x \leq x'$, contradicting the fact that $f(x') < f(x)$.
	We conclude that 
	\[
	\realis{\ordcpx{f}}(\inter{\realis{\ordcpx{\clset{x}}}}) = \inter{\realis{\ordcpx{\clset{f(x)}}}},\]
	so $\realis{\ordcpx{f}}$ is a map of regular CW complexes.
	Naturality of the isomorphism between $\fpos{\realis{\ordcpx{-}}}$ and $\bigid{\cwpos}$ is straightforward.
\end{proof}

\begin{prop} \label{prop:face_poset_has_an_inverse_up_to_isomorphism}
	Let $\cwcom{X}{S}$ be a regular CW complex.
	Then $\cwcom{X}{S}$ is isomorphic to $\realis{\ordcpx{\fpos{\cwcom{X}{S}}}}$ in $\rcpx$.
\end{prop}
\begin{proof}
	See \cite[Chapter III, Theorem 1.7]{lundell1969topology}.
\end{proof}

\begin{rmk}
	No claim of naturality is made in Proposition \ref{prop:face_poset_has_an_inverse_up_to_isomorphism}, and for good reason, since $\fpos{}\colon \rcpx \to \cwpos$ is evidently not faithful.
\end{rmk}


\section{Presenting regular CW complexes} \label{sec:regular_cw}

\begin{guide}
	In this section, we prove the main results of this chapter: the order complex of a round $n$\nbd dimensional molecule is a PL $n$\nbd ball, and the order complex of its boundary is a PL $(n-1)$\nbd sphere (Proposition \ref{prop:order_complex_of_round_molecule}), so the geometric realisation of a regular directed complex has a canonical structure of regular CW complex (Theorem 
	\ref{thm:regular_directed_complex_is_regular_cw_complex}).

	One immediate consequence is that we can, in principle, compute the \cemph{cellular homology} of a regular CW complex by giving it a presentation as a regular directed complex, then computing the homology of its augmented chain complex (Corollary \ref{cor:cellular_homology}).
	In general, it seems plausible that working with directed complexes could sometimes be beneficial in ``undirected'' combinatorial topology, thanks to the algebraic grip on pasting given by the input-output subdivision of boundaries, and the ability to construct a valid orientation inductively.
	This potential is still largely unexplored.
\end{guide}

\begin{dfn}[Order complex of an oriented graded poset] \index{oriented graded poset!order complex}
	Let $P$ be an oriented graded poset.
	The \emph{order complex of $P$} is the order complex of its underlying poset.
	We write $\ordcpx{P}$ for $\ordcpx{(\fun{U}P)}$, and $\ordcpx{-}\colon \ogpos \to \sset$ for the composite functor $\ordcpx{(\fun{U}-)}$.
\end{dfn}

\begin{prop} \label{prop:order_complex_of_round_molecule}
	Let $U$ be a round molecule, $n \eqdef \dim{U}$.
	Then
	\begin{enumerate}
		\item $\ordcpx{U}$ is a PL $n$\nbd ball,
		\item $\ordcpx{(\bound{}{}U)} = \bound{}{}(\ordcpx{U})$ is a PL $(n-1)$\nbd sphere.
	\end{enumerate}
\end{prop}
\begin{proof}
	We proceed by induction on $n$.
	If $n = 0$, then $U = 1$ and its order complex is the terminal simplicial set, which is a PL 0\nbd ball.
	
	Let $n > 0$, and suppose that $U$ is an atom.
	By Lemma \ref{lem:atom_merger_of_its_boundary}, $U$ is isomorphic to $\bound{}{-}U \celto \bound{}{+}U$, so there is a pushout diagram
\[\begin{tikzcd}
	{\bound{n-2}{}U} && \bound{}{+}U \\
	\bound{}{-}U && {\bound{}{}U}
	\arrow[hook', from=1-1, to=2-1]
	\arrow[hook, from=1-1, to=1-3]
	\arrow[hook, from=2-1, to=2-3]
	\arrow[hook', from=1-3, to=2-3]
	\arrow["\lrcorner"{anchor=center, pos=0.125, rotate=180}, draw=none, from=2-3, to=1-1]
\end{tikzcd}\]
	of inclusions in $\ogpos$.
	By Lemma \ref{lem:order_complex_pushouts}, $\ordcpx{-}$ sends this to a pushout diagram of monomorphisms in $\sset$.
	By the inductive hypothesis, $\ordcpx{(\bound{}{-}U)}$ and $\ordcpx{(\bound{}{+}U)}$ are PL $(n-1)$\nbd balls, while $\ordcpx{(\bound{n-2}{}U)}$ is a PL $(n-2)$\nbd sphere, equal to $\bound{}{}\ordcpx{(\bound{}{\alpha}U)}$ for each $\alpha \in \set{+, -}$.
	By Proposition \ref{prop:gluing_pl_balls_and_spheres}, $\ordcpx{(\bound{}{}U)}$ is a PL $(n-1)$\nbd sphere.
	Since the underlying poset of $U$ is the posetal join of $\bound{}{}U$ with a single point, $\ordcpx{U}$ is a cone over the PL $(n-1)$\nbd sphere $\ordcpx{(\bound{}{}U)}$, which is a PL $n$\nbd ball with boundary $\bound{}{}(\ordcpx{U})$ by \cite[Lemma 10]{zeeman1966seminar}.

	Next, suppose that $U$ is not an atom, and let $(\order{i}{U})_{i=1}^m$ be an $(n-1)$\nbd layering of $U$, with associated $(n-1)$\nbd ordering $(\order{i}{x})_{i=1}^m$.
	Define, recursively,
	\begin{align*}
		\order{0}{V} & \eqdef \bound{}{-}U \celto \bound{}{-}U, \\
		\order{i}{V} & \eqdef \order{i-1}{V} \cp{n-1} \order{i}{U}.
	\end{align*}
	We will prove that $\ordcpx{\order{i}{V}}$ is a PL $n$\nbd ball by recursion on $i \in \set{0, \ldots, m}$.
	First of all, $\order{0}{V}$ is an $n$\nbd dimensional atom, so its order complex is a PL $n$\nbd ball by what we have already proved.
	Suppose $i > 0$.
	It follows from Proposition \ref{prop:layering_from_ordering} that $\order{i}{V}$ is the pasting of $\clset{\order{i}{x}}$ at the submolecule $\bound{}{-}\order{i}{x} \submol \bound{}{+}\order{i-1}{V}$, so $\order{i}{V}$ can be exhibited by a pushout diagram
\[\begin{tikzcd}
	{\bound{}{-}\order{i}{x}} && \clset{\order{i}{x}} \\
	\order{i-1}{V} && {\order{i}{V}}
	\arrow[hook', from=1-1, to=2-1]
	\arrow[hook, from=1-1, to=1-3]
	\arrow[hook, from=2-1, to=2-3]
	\arrow[hook', from=1-3, to=2-3]
	\arrow["\lrcorner"{anchor=center, pos=0.125, rotate=180}, draw=none, from=2-3, to=1-1]
\end{tikzcd}\]
	in $\ogpos$ with $\bound{}{-}\order{i}{x} = \bound{}{}\order{i-1}{V} \cap \bound{}{}\order{i}{x}$.
	Then $\ordcpx{-}$ sends this to a pushout diagram of monomorphisms in $\sset$.
	By the inductive hypothesis, $\ordcpx{\clset{\order{i}{x}}}$ and $\ordcpx{\order{i-1}{V}}$ are PL $n$\nbd balls, while $\ordcpx{(\bound{}{-}\order{i}{x})}$ is a PL $(n-1)$\nbd ball contained in their boundaries.
	It follows from Proposition \ref{prop:gluing_pl_balls_and_spheres} that $\ordcpx{(\order{i}{V})}$ is a PL $n$\nbd ball with boundary $\ordcpx{(\bound{}{}\order{i}{V})}$.
	Let $V \eqdef \order{m}{V}$; we conclude that $\ordcpx{V}$ is a PL $n$\nbd ball.
	Moreover, $V$ is round and, by construction, isomorphic to $\order{0}{V} \cp{n-1} U$.

	Consider the $(n+1)$\nbd dimensional atom $W \eqdef V \celto \compos{V}$.
	By the same argument proving that $\ordcpx{(\bound{}{}U)}$ is a PL $(n-1)$\nbd sphere, $\ordcpx{(\bound{}{}W)}$ is a PL $n$\nbd sphere.
	Moreover, since $U$ is round, in $\bound{}{}W$ we have
	\begin{enumerate}
		\item $\compos{V} \cap \order{0}{V} = \bound{}{-}V = \bound{}{-}\compos{V} = \bound{}{-}\order{0}{V}$, 
		\item $U \cap (\compos{V} \cup \order{0}{V}) = \bound{}{}U$.
	\end{enumerate}
	Then we have a pushout diagram
	\[\begin{tikzcd}
		{\bound{}{-}V} && {\order{0}{V}} \\
		{\compos{V}} && \compos{V} \cup \bound{}{-}V
		\arrow[hook', from=1-1, to=2-1]
		\arrow[hook, from=1-1, to=1-3]
		\arrow[hook, from=2-1, to=2-3]
		\arrow[hook', from=1-3, to=2-3]
		\arrow["\lrcorner"{anchor=center, pos=0.125, rotate=180}, draw=none, from=2-3, to=1-1]
	\end{tikzcd}\]
	in $\ogpos$, preserved by $\ordcpx{-}$.
	Since $\ordcpx{\compos{V}}$ and $\ordcpx{\order{0}{V}}$ are PL $n$\nbd balls, while $\ordcpx{(\bound{}{-}V)}$ is a PL $(n-1)$\nbd ball contained in their boundaries, it follows that $\ordcpx{(\compos{V} \cup \bound{}{-}V)}$ is a PL $n$\nbd ball.
	Finally, we have a pushout diagram
	\[\begin{tikzcd}
	{\bound{}{}U} && {U} \\
	\compos{V} \cup \order{0}{V} && \bound{}{}W
	\arrow[hook', from=1-1, to=2-1]
	\arrow[hook, from=1-1, to=1-3]
	\arrow[hook, from=2-1, to=2-3]
	\arrow[hook', from=1-3, to=2-3]
	\arrow["\lrcorner"{anchor=center, pos=0.125, rotate=180}, draw=none, from=2-3, to=1-1]
\end{tikzcd}\]
	in $\ogpos$, preserved by $\ordcpx{-}$.
	Since $\ordcpx{(\compos{V} \cup \order{0}{V})}$ is a PL $n$\nbd ball and $\ordcpx{(\bound{}{}W)}$ is a PL $n$\nbd sphere, we conclude by a final application of Proposition \ref{prop:gluing_pl_balls_and_spheres} that $\ordcpx{U}$ is a PL $n$\nbd ball with boundary $\bound{}{}(\ordcpx{U})$.
\end{proof}

\begin{cor} \label{cor:regular_directed_complex_is_cw_poset}
	Let $P$ be a regular directed complex.
	Then $\fun{U}P$ is a CW poset.
\end{cor}
\begin{proof}
	By Proposition \ref{prop:order_complex_of_round_molecule}, for all $x \in P$, $\ordcpx{(\bound{}{}x)} = \ordcpx{(\clset{x} \setminus \set{x})}$ is a PL sphere, so its geometric realisation is homeomorphic to a sphere.
\end{proof}

\begin{cor} \label{cor:forgetful_on_rdcpxmap_factors_through_cwpos}
	The functor $\fun{U}\colon \rdcpxmap \to \posclos$ factors through the subcategory inclusion $\cwpos \incl \posclos$.
\end{cor}

\begin{thm} \label{thm:regular_directed_complex_is_regular_cw_complex}
	Let $P$ be a regular directed complex.
	Then $\realis{\ordcpx{P}}$ admits a structure of regular CW complex whose set of generating cells is 
	\[
		\set{\realis{\ordcpx{\clset{x}}} \incl \realis{\ordcpx{P}} \mid x \in P}.
	\]
	This assignment extends to a functor $\realis{\ordcpx{-}}\colon \rdcpxmap \to \rcpx$, such that the diagram of functors
\[\begin{tikzcd}
	\rdcpxmap && \rcpx \\
	& \cwpos
	\arrow["{\realis{\ordcpx{-}}}", from=1-1, to=1-3]
	\arrow["{\fun{U}}"', from=1-1, to=2-2]
	\arrow["{\fpos{}}", from=1-3, to=2-2]
\end{tikzcd}\]
	commutes up to natural isomorphism.
\end{thm}
\begin{proof}
	A combination of Proposition \ref{prop:ordcpx_right_inverse_to_face_poset} and Corollary \ref{cor:forgetful_on_rdcpxmap_factors_through_cwpos}.
\end{proof}

\begin{cor} \label{cor:cellular_homology}
	Let $P$ be a regular directed complex.
	The cellular homology of the regular CW complex $\realis{\ordcpx{P}}$ is naturally isomorphic to the homology of the chain complex $\freeab{P}$.
\end{cor}

\begin{comm} \index{product!of spaces} \index{suspension!of a space} \index{join!of spaces}
	Given spaces $X$, $Y$, let
	\begin{itemize}
		\item $X \times Y$ denote the product of $X$ and $Y$,
		\item $\sus{X}$ denote the suspension of $X$,
		\item $X \join Y$ denote the join of $X$ and $Y$,
	\end{itemize}
	defined in the standard way, with all limits and colimits taken in $\ktop$.

	When $X$ and $Y$ have a structure of CW complex, then $X \times Y$, $\sus{X}$, and $X \join Y$ also have a canonical structure of CW complex \cite[Chapter I, Section 7]{lundell1969topology}.
	This structure is regular if those on $X$ and $Y$ are.
\end{comm}

\begin{prop} \label{prop:realisation_of_operations}
	Let $P$, $Q$ be regular directed complexes.
	Then
	\begin{enumerate}
		\item $\realis{\ordcpx{(P \gray Q)}}$ is homeomorphic to $\realis{\ordcpx{P}} \times \realis{\ordcpx{Q}}$,
		\item $\realis{\ordcpx{(\sus{P})}}$ is homeomorphic to $\sus{\realis{\ordcpx{P}}}$,
		\item $\realis{\ordcpx{(P \join Q)}}$ is homeomorphic to $\realis{\ordcpx{P}} \join \realis{\ordcpx{Q}}$.
	\end{enumerate}
\end{prop}
\begin{proof}
	By comparison with \cite[Chapter I, Section 7]{lundell1969topology}, the face posets of $\realis{\ordcpx{P}} \times \realis{\ordcpx{Q}}$,  $\sus{\realis{\ordcpx{P}}}$, and $\realis{\ordcpx{P}} \join \realis{\ordcpx{Q}}$ with their canonical structure are isomorphic to the underlying posets of $P \gray Q$, $\sus{P}$, and $P \join Q$, respectively.
	The statement then follows from Proposition \ref{prop:face_poset_has_an_inverse_up_to_isomorphism}.
\end{proof}

\clearpage
\thispagestyle{empty}

%% file: steiner.tex
\chapter{Steiner theory} \label{chap:steiner}
\thispagestyle{firstpage}

\begin{guide}
	We have already seen how to construct a chain complex from a regular directed complex, or more generally from an oriented thin graded poset.
	This construction loses some information (technically, it is not full on isomorphisms): there is no intrinsic way to distinguish a generator $x$ of the free abelian group $\freeab{\grade{n}{P}}$ from its inverse $-x$.
	We can restore this information by ``marking'' the free \emph{commutative monoid} $\freemon{\grade{n}{P}}$ as a distinguished submonoid of $\freeab{\grade{n}{P}}$, at which point we can uniquely reconstruct $P$.

	This structure of an augmented chain complex together with a sequence of distinguished submonoids is called an \cemph{augmented directed chain complex}, and is the fundamental ingredient of \cemph{Steiner theory}, as it has come to be known.
	Even though it starts from generic directed chain complexes, Steiner theory rapidly focusses on very particular subclasses, now known as \cemph{Steiner complexes} and \cemph{strong Steiner complexes}.

	As it turns out, all the objects needed for applications of Steiner theory are in the image of a functor $\dfreeab{-}$ from regular directed complexes: they are the augmented directed chain complexes of thetas, their Gray products, and their joins, all of which are strong Steiner complexes (Proposition \ref{prop:acyclic_rdcpx_gives_strong_steiner}).
	So at the level of objects, there is nothing gained by moving from oriented graded posets to augmented directed chain complexes.

	On the other hand, there is an adjunction between strict $\omega$\nbd categories and augmented directed chain complexes which restricts to a \emph{full and faithful} functor on Steiner complexes 
	(Theorem \ref{thm:steiner_main_theorem}).
	It follows that strong Steiner complexes with their natural notion of morphism determine a \emph{dense subcategory} of $\omegacat$.
	Augmented directed chain complexes inherit a tensor product and join from augmented chain complexes, these restrict to strong Steiner complexes, and by density can then be extended along colimits to all strict $\omega$\nbd categories.
	This is the preferred way, nowadays, of defining the Gray product and join of strict $\omega$\nbd categories.

	Since the natural notions of morphism between regular directed complexes do not cover all functors between their $\omega$\nbd categories of molecules, this is a genuine use case of Steiner theory.
	On the other hand, Steiner complexes have also been used as a formalism for higher-categorical diagrams, a purpose for which they have only disadvantages relative to regular directed complexes.
	The aim of this chapter, besides presenting the main results of Steiner theory, is to clarify the relation between the two formalisms, and in particular when one can safely switch between a regular directed complex and its augmented directed chain complex, and when one cannot.

	The main point of divergence is the following.
	To go from a regular directed complex $P$ to an augmented directed chain complex, it is \emph{always} equivalent to directly apply the functor $\dfreeab{-}$, which endows $\freeab{P}$ with the sequence of free commutative monoids $(\freemon{\grade{n}{P}})_{n \in \mathbb{N}}$, or to apply Steiner's left adjoint functor $\linea{}$ to the strict $\omega$\nbd category $\molecin{P}$ (Theorem
	\ref{thm:two_functors_from_rdcpx_to_dchaug}).

	However, $\molecin{P}$ is \emph{only} equivalent to the application of Steiner's right adjoint functor $\nufun{}$ to $\dfreeab{P}$ when $P$ is \emph{dimension-wise acyclic}, in which case $\dfreeab{P}$ is a Steiner complex (Proposition 
	\ref{prop:dw_acyclic_rdcpx_gives_steiner_complex} and Theorem 
	\ref{thm:two_omegacats_from_dw_acyclic_rdcpx}).
	In general, $\nufun{\dfreeab{P}}$ is a significantly more degenerate quotient of $\molecin{P}$, and the latter construction should be favoured.

	From these results, we deduce that $\molecin{-}$ is compatible with Gray products and joins of \emph{acyclic} regular directed complexes 
	(Proposition \ref{prop:monoidal_functors_from_acyclic}).
	It seems unlikely that this extends to all regular directed complexes: make of it what you will.
\end{guide}


\section{Directed chain complexes} \label{sec:directed_chain}

\begin{guide}
	In this section, we define the category $\dchaug$ of augmented directed chain complexes and homomorphisms.
	We lift the functors $\freeab{-}$ valued in $\chaug$ to functors $\dfreeab{-}$ valued in $\dchaug$.
	We then lift all the constructions of augmented chain complexes that we considered in Chapter \ref{chap:constructions} --- tensor products, suspensions, joins, and duals --- to augmented directed chain complexes, and prove that $\dfreeab{-}$ is compatible with all.

	Finally, we recall Steiner's adjunction consisting of a functor $\linea{}$ from $\omegacat$ to $\dchaug$ and its right adjoint $\nufun{}$, and prove that $\dfreeab{-}$ is naturally isomorphic to $\linea{\molecin{-}}$ on regular directed complexes 
	(Theorem \ref{thm:two_functors_from_rdcpx_to_dchaug}).
\end{guide}

\begin{dfn}[Direction on a chain complex] \index{chain complex!direction} \index{$\grade{n}{\dir{C}}$}
Let $C$ be a chain complex of abelian groups in non-negative degree.
A \emph{direction} on $C$ is a choice of a commutative submonoid $\grade{n}{\dir{C}}$ of $\grade{n}{C}$ for each $n \in \mathbb{N}$.
\end{dfn}

\begin{dfn}[Augmented directed chain complex] \index{directed chain complex} \index{directed chain complex!augmented} \index{chain complex!directed|see {directed chain complex}} \index{augmented directed chain complex|see {directed chain complex}}
	An \emph{augmented directed chain complex} is an augmented chain complex $C$ together with a direction on its underlying chain complex.
\end{dfn}

\begin{dfn}[Homomorphism of augmented directed chain complexes] \index{directed chain complex!homomorphism}
	Let $C, D$ be augmented directed chain complexes.
	A \emph{homomorphism} $f\colon C \to D$ is a homomorphism of the underlying augmented chain complexes such that
	\[
		\grade{n}{f}(\grade{n}{\dir{C}}) \subseteq \grade{n}{\dir{D}}
	\]
	for all $n \in \mathbb{N}$.
\end{dfn}

\begin{dfn}[The category $\dchaug$] \index{$\dchaug$}
	We let $\dchaug$ denote the category whose objects are augmented directed chain complexes and morphisms are homomorphisms of augmented directed chain complexes.
\end{dfn}

\begin{lem} \label{lem:forgetful_dchaug}
	Forgetting the direction determines a faithful functor 
	\[ \fun{U}\colon \dchaug \to \chaug. \]
\end{lem}
\begin{proof}
	By definition.
\end{proof}

\begin{dfn}[Augmented directed chain complex of an oriented thin graded poset] \index{directed chain complex!of an oriented thin graded poset} \index{$\dfreeab{P}$}
	Let $P$ be an oriented graded poset such that $\augm{P}$ is oriented thin.
	The \emph{augmented directed chain complex of $P$} is the augmented directed chain complex $\dfreeab{P}$ whose
	\begin{itemize}
		\item underlying augmented chain complex is $\freeab{P}$,
		\item direction is given by $\grade{n}{\dir{\freeab{P}}} \eqdef \freemon{\grade{n}{P}}$ for each $n \in \mathbb{N}$,
	\end{itemize}
	where $\freemon{\grade{n}{P}}$ is the free commutative monoid on the set $\grade{n}{P}$.
\end{dfn}

\begin{lem} \label{lem:dfreeab_functorial_on_otgpos}
	The assignment $P \mapsto \dfreeab{P}$ extends to a unique functor $\dfreeab{-}\colon \otgpos \to \dchaug$ such that the triangle of functors
	\[\begin{tikzcd}
	\otgpos && \dchaug \\
	&& \chaug
	\arrow["{\dfreeab{-}}", from=1-1, to=1-3]
	\arrow["{\fun{U}}", from=1-3, to=2-3]
	\arrow["{\freeab{}}"', curve={height=6pt}, from=1-1, to=2-3]
\end{tikzcd}\]
	commutes.
\end{lem}
\begin{proof}
	The requirement that $\fun{U} \after \dfreeab{-}$ be strictly equal to $\freeab{}$, together with the specification on objects, determines the functor uniquely, so it suffices to show that, given any morphism $f\colon P \to Q$ of oriented graded posets, $\freeab{f}$ lifts to a homomorphism of augmented directed chain complexes.
	But this is true since $\grade{n}{\freeab{f}}$ maps generators of $\freemon{\grade{n}{P}}$ to generators of $\freemon{\grade{n}{Q}}$.
\end{proof}

\begin{lem} \label{lem:dfreeab_functorial_on_maps_and_comaps}
	There exists a unique pair of functors
	\[
		\dfreeab{-}\colon \rdcpxmap \to \dchaug, \quad \quad \dfreeab{\pb{-}}\colon \opp{\rdcpxcomap} \to \dchaug
	\]
	such that the diagram
\[\begin{tikzcd}
	{\opp{\rdcpxiso}} && \rdcpxiso && \rdcpx \\
	{\opp{\rdcpxcomap}} && \otgpos && \rdcpxmap \\
	&& \dchaug \\
	&& \chaug
	\arrow["{\augm{(-)}}", from=1-3, to=2-3]
	\arrow["{\dfreeab{-}}", from=2-3, to=3-3]
	\arrow["{\invrs{(-)}}", from=1-1, to=1-3]
	\arrow[hook', from=1-1, to=2-1]
	\arrow["{\dfreeab{\pb{-}}}"', from=2-1, to=3-3]
	\arrow[hook, from=1-3, to=1-5]
	\arrow[hook', from=1-5, to=2-5]
	\arrow["{\dfreeab{-}}", from=2-5, to=3-3]
	\arrow["{\fun{U}}", from=3-3, to=4-3]
	\arrow["{\freeab{\pb{-}}}"', curve={height=18pt}, from=2-1, to=4-3]
	\arrow["{\freeab{}}", curve={height=-18pt}, from=2-5, to=4-3]
\end{tikzcd}\]
	commutes.
\end{lem}
\begin{proof}
	Commutativity of the diagram specifies the two functors uniquely both on objects and on morphisms, so it suffices to show that, for all maps $p\colon P \to Q$ and comaps $c\colon P \to Q$ of regular directed complexes, the homomorphisms $\freeab{p}$ and $\freeab{\pb{c}}$ are compatible with the directions on $\freeab{P}$ and $\freeab{Q}$.
	This is evident by inspection of their definitions.
\end{proof}

\begin{dfn}[Tensor product of augmented directed chain complexes] \index{directed chain complex!tensor product} \index{product!tensor}
	Let $C$, $D$ be augmented directed chain complexes.
	The \emph{tensor product of $C$ and $D$} is the tensor product of the underlying augmented chain complexes of $C$ and $D$, with the direction defined by
	\[
		\grade{n}{\dir{(C \otimes D)}} \eqdef \bigoplus_{k=0}^n \grade{k}{\dir{C}} \otimes \grade{n-k}{\dir{D}},
	\]
	where $\oplus$ and $\otimes$ are, respectively, the direct sum and tensor product of commutative monoids.
	The tensor product of augmented directed chain complexes extends to a monoidal structure on $\dchaug$, whose unit is $\mathbb{Z}$ with the direction $\grade{0}{\dir{\mathbb{Z}}} \eqdef \mathbb{N}$, such that the forgetful functor $\fun{U}\colon \dchaug \to \chaug$ lifts to a strict monoidal functor
	\[
		\fun{U}\colon (\dchaug, \otimes, \mathbb{Z}) \to (\chaug, \otimes, \mathbb{Z});
	\]
	see \cite[Example 3.10]{steiner2004omega}.
\end{dfn}

\begin{prop} \label{prop:gray_to_tensor_of_dchaug}
	The functors
	\begin{align*}
		&\dfreeab{-}\colon \otgpos \to \dchaug, \\
		&\dfreeab{-}\colon \rdcpxmap \to \dchaug, \\
		&\dfreeab{\pb{-}}\colon \opp{\rdcpxcomap} \to \dchaug
	\end{align*}
	lift to strong monoidal functors
	\begin{align*}
		&\dfreeab{-} \colon (\otgpos, \augm{\gray}, \augm{1}) \to (\dchaug, \otimes, \mathbb{Z}), \\
		&\dfreeab{-} \colon (\rdcpxmap, \gray, 1) \to (\dchaug, \otimes, \mathbb{Z}), \\
		&\dfreeab{\pb{-}} \colon (\opp{\rdcpxcomap}, \gray, 1) \to (\dchaug, \otimes, \mathbb{Z}).
	\end{align*}
\end{prop}
\begin{proof}
	Immediate after Corollary \ref{cor:monoidal_functors_to_chaug}.
\end{proof}

\begin{dfn}[Suspension of an augmented directed chain complex] \index{directed chain complex!suspension} \index{suspension!of an augmented directed chain complex}
	Let $C$ be an augmented directed chain complex.
	The \emph{suspension of $C$} is the suspension of the underlying augmented chain complex of $C$, with the direction defined by
	\[
		\grade{n}{\dir{\sus{C}}} \eqdef
		\begin{cases}
			\freemon{\set{\bot^+, \bot^-}} & \text{if $n = 0$,} \\
			\grade{n-1}{\dir{C}} & \text{if $n > 0$.}
		\end{cases}
	\]
	The suspension extends to an endofunctor $\fun{S}$ of $\dchaug$, such that the diagram
	\[\begin{tikzcd}
	\dchaug && \dchaug \\
	\chaug && \chaug
	\arrow["{\sus{}}", from=1-1, to=1-3]
	\arrow["{\fun{U}}", from=1-1, to=2-1]
	\arrow["{\sus{}}", from=2-1, to=2-3]
	\arrow["{\fun{U}}", from=1-3, to=2-3]
\end{tikzcd}\]
	commutes; see \cite[\S 1.2]{ozornova2023quillen}.
\end{dfn}

\begin{prop} \label{prop:suspension_otgpos_to_dchaug}
	Let $P$ be an oriented graded poset such that $\augm{P}$ is oriented thin.
	Then $\dfreeab{(\sus{P})}$ is naturally isomorphic to $\sus{(\dfreeab{P})}$.
\end{prop}
\begin{proof}
	Straightforward after Proposition \ref{prop:suspension_of_otgpos_and_chaug}.
\end{proof}

\begin{dfn}[Join of augmented directed chain complexes] \index{directed chain complex!join} \index{join!of augmented directed chain complexes}
	Let $C$, $D$ be augmented directed chain complexes.
	The \emph{join of $C$ and $D$} is the join of the underlying augmented chain complexes of $C$ and $D$, with the direction defined by
	\[
	\grade{n}{\dir{(C \join D)}} \eqdef \grade{n}{\dir{D}} \oplus \left( \bigoplus_{k=0}^{n-1} \grade{k}{\dir{C}} \otimes \grade{n-1-k}{\dir{D}} \right) \oplus \grade{n}{\dir{C}}
	\]
	for each $n \in \mathbb{N}$.
	The join of augmented directed chain complexes extends to a monoidal structure on $\dchaug$, whose unit is $0$ with the trivial direction, such that the forgetful functor $\fun{U}\colon \dchaug \to \chaug$ lifts to a strict monoidal functor
	\[
		\fun{U}\colon (\dchaug, \join, 0) \to (\chaug, \join, 0);
	\]
	see \cite[\S 6.5]{ara2020joint}.
\end{dfn}

\begin{prop} \label{prop:join_to_join_of_dchaug}
	The functors
	\begin{align*}
		&\dfreeab{-}\colon \otgpos \to \dchaug, \\
		&\dfreeab{-}\colon \rdcpxmap \to \dchaug, \\
		&\dfreeab{\pb{-}}\colon \opp{\rdcpxcomap} \to \dchaug
	\end{align*}
	lift to strong monoidal functors
	\begin{align*}
		&\dfreeab{-} \colon (\otgpos, \gray, 1) \to (\dchaug, \join, 0), \\
		&\dfreeab{-} \colon (\rdcpxmap, \join, \varnothing) \to (\dchaug, \join, 0), \\
		&\dfreeab{\pb{-}} \colon (\opp{\rdcpxcomap}, \join, \varnothing) \to (\dchaug, \join, 0).
	\end{align*}
\end{prop}
\begin{proof}
	Immediate after Corollary \ref{cor:join_monoidal_functors_to_chaug}.
\end{proof}

\begin{dfn}[Duals of an augmented directed chain complex] \index{directed chain complex!dual} \index{dual!of an augmented directed chain complex}
	Let $C$ be an augmented directed chain complex, $J \subseteq \posnat$.
	The \emph{$J$-dual of $C$} is the $J$\nbd dual of the underlying augmented chain complex of $C$, with the direction defined by
	\[
		\grade{n}{\dir{\dual{J}{C}}} \eqdef \grade{n}{\dir{C}}.
	\]
	The $J$\nbd dual extends to an endofunctor $\dual{J}{}$ of $\dchaug$, such that the diagram
	\[\begin{tikzcd}
	\dchaug && \dchaug \\
	\chaug && \chaug
	\arrow["{\dual{J}{}}", from=1-1, to=1-3]
	\arrow["{\fun{U}}", from=1-1, to=2-1]
	\arrow["{\dual{J}{}}", from=2-1, to=2-3]
	\arrow["{\fun{U}}", from=1-3, to=2-3]
	\end{tikzcd}\]
	commutes; see \cite[\S 2.18]{ara2020joint}.
\end{dfn}

\begin{prop}
	Let $P$ be an oriented graded poset such that $\augm{P}$ is oriented thin.
	Then $\dfreeab{(\dual{J}{P})}$ is naturally isomorphic to $\dual{J}{(\dfreeab{P})}$.
\end{prop}
\begin{proof}
	Immediate after Proposition \ref{prop:dual_chaug_and_ogtpos}.
\end{proof}

\begin{dfn}[Linearisation of a strict $\omega$-category] \index{strict $\omega$-category!linearisation}
	Let $X$ be a strict $\omega$-category.
	The \emph{linearisation of $X$} is the augmented directed chain complex $\linea{X}$ whose underlying augmented chain complex is defined by
	\[
		\grade{n}{\linea{X}} \eqdef \frac{
		\freeab{(\skel{n}{X})} } {
		\spanset{ t \cp{k} u - t - u \mid \text{$t, u \in \skel{n}{X}$, $k < n$} } }\;
	\]
	for all $n \in \mathbb{N}$, where $\freeab{(\skel{n}{X})}$ denotes the free abelian group on the set of cells of the $n$\nbd skeleton of $X$, with
	\begin{align*}
		\der\colon \grade{n}{\linea{X}} & \to \grade{n-1}{\linea{X}}, \\
			t \in \skel{n}{X} & \mapsto \bound{n-1}{+}t - \bound{n-1}{-}t,
	\end{align*}
	for each $n > 0$, and
	\begin{align*}
		\eau\colon \grade{0}{\linea{X}} & \to \mathbb{Z}, \\
		t \in \skel{0}{X} & \mapsto 1,
	\end{align*}
	together with the direction defined by
	\[
		\grade{n}{\dir{\linea{X}}} \eqdef \Ima \left(
			\freemon{(\skel{n}{X})}
			\incl \freeab{(\skel{n}{X})}
			\to \grade{n}{\linea{X}} \right)
	\]
	for each $n \in \mathbb{N}$, where $\freeab{(\skel{n}{X})} \to \grade{n}{\linea{X}}$ is the canonical quotient homomorphism.
	Linearisation extends to a functor $\omegacat \to \dchaug$, sending a strict functor $f\colon X \to Y$ to the homomorphism defined by
	\begin{align*}
		\grade{n}{\linea{f}}\colon \grade{n}{\linea{X}} & \to \grade{n}{\linea{Y}}, \\
			t \in \skel{n}{X} & \mapsto f(t).
	\end{align*}
\end{dfn}

\begin{dfn}[Globular table in an augmented directed chain complex] \index{directed chain complex!globular table}
	Let $C$ be an augmented directed chain complex.
	A \emph{globular table in $C$} is a double sequence 
	\[
		x \equiv (\gltab{n}{\alpha}{x})_{n \in \mathbb{N}, \, \alpha \in \set{+, -}} 
	\]
	such that
	\begin{enumerate}
		\item $\gltab{n}{\alpha}{x} \in \grade{n}{\dir{C}}$ for all $n \in \mathbb{N}$ and $\alpha \in \set{+, -}$,
		\item $\der \gltab{n}{\alpha}{x} = \gltab{n-1}{+}{x} - \gltab{n-1}{-}{x}$ for all $n > 0$ and $\alpha \in \set{+, -}$,
		\item $\eau \gltab{0}{\alpha}{x} = 1$ for all $\alpha \in \set{+, -}$,
		\item there exists $m \in \mathbb{N}$ such that $\gltab{n}{\alpha}{x} = 0$ for all $n > m$ and $\alpha \in \set{+, -}$.
	\end{enumerate}
\end{dfn}

\begin{dfn}[Strict $\omega$-category of globular tables] \index{directed chain complex!globular table!strict $\omega$-category}
	Let $C$ be an augmented directed chain complex.
	The \emph{strict $\omega$-category of globular tables in $C$} is the strict $\omega$\nbd category $\nufun{C}$ whose set of cells is
	\[
		\set{ x \mid \text{$x$ is a globular table in $C$} },
	\]
	with the boundary operators defined, for each globular table $x$ in $C$, $n \in \mathbb{N}$, and $\alpha \in \set{+, -}$, by
	\[
		\gltab{m}{\beta}{(\bound{n}{\alpha}x)} \eqdef 
		\begin{cases}
			\gltab{m}{\beta}x
			& \text{if $m < n$,} \\
			\gltab{n}{\alpha}x
			& \text{if $m = n$,} \\
			0
			& \text{if $m > n$,}
		\end{cases}
	\]
	and the $k$\nbd composition operations defined, for all $k \in \mathbb{N}$ and $k$\nbd composable pairs $x, y$ of globular tables, by
	\[
		\gltab{n}{\alpha}{(x \cp{k} y)} \eqdef \gltab{n}{\alpha}{x} - \gltab{n}{\alpha}{(\bound{k}{+}x)} + \gltab{n}{\alpha}{y}.
	\]
	This assignment extends to a functor $\nufun{}\colon \dchaug \to \omegacat$, sending a homomorphism $f\colon C \to D$ to the strict functor defined by
	\begin{align*}
		\nufun{f}\colon \nufun{C} & \to \nufun{D}, \\
		x & \mapsto f(x), \quad \gltab{n}{\alpha}{f(x)} \eqdef \grade{n}{f}(\gltab{n}{\alpha}{x}).
	\end{align*}
\end{dfn}

\begin{prop} \label{prop:adjunction_dchaug_omegacat}
	The functor $\linea{}\colon \omegacat \to \dchaug$ is left adjoint to $\nufun{}\colon \dchaug \to \omegacat$. 
\end{prop}
\begin{proof}
	See \cite[Theorem 2.11]{steiner2004omega}.
\end{proof}

\begin{lem} \label{lem:suspension_of_dchaug_and_omegacat}
	Let $C$ be an augmented directed chain complex.
	Then $\nufun{\sus{C}}$ is naturally isomorphic to $\sus{(\nufun{C})}$.
\end{lem}
\begin{proof}
	This is \cite[Lemma 2.12]{ozornova2023quillen}.
\end{proof}

\begin{lem} \label{lem:duals_of_dchaug_and_omegacat}
	Let $C$ be an augmented directed chain complex, $X$ a strict $\omega$\nbd category, and $J \subseteq \posnat$.
	Then
	\begin{enumerate}
		\item $\linea{\dual{J}{X}}$ is naturally isomorphic to $\dual{J}{(\linea{X})}$,
		\item $\nufun{\dual{J}{C}}$ is naturally isomorphic to $\dual{J}{(\nufun{C})}$.
	\end{enumerate}
\end{lem}
\begin{proof}
	See \cite[Proposition 2.19]{ara2020joint}.
\end{proof}

\begin{lem} \label{lem:linearisation_lower_dim_cells}
	Let $t$ be a cell in a strict $\omega$\nbd category $X$, $n > \dim{t}$.
	Then $t = 0$ in $\grade{n}{\linea{X}}$.
\end{lem}
\begin{proof}
	Let $k \eqdef \dim{t}$.
	Then $t = t \cp{k} t$ by 
	Lemma \ref{lem:composition_with_lower_dim}, so $t = t + t$, hence $t = 0$, in $\grade{n}{\linea{X}}$.
\end{proof}

\begin{lem} \label{lem:cells_in_linearisation_of_molecin}
	Let $P$ be an oriented graded poset, $\isocl{f\colon U \to P}$ a molecule over $P$, $n \geq \dim{U}$.
	Then, in $\grade{n}{\left(\linea{\molecin{P}}\right)}$,
	\[
		\isocl{f\colon U \to P} = \sum_{x \in \grade{n}{U}} \isocl{\restr{f}{\clset{x}}\colon \clset{x} \to P}.
	\]
	\end{lem}
\begin{proof}
	We proceed by induction on $k \eqdef \lydim{U}$.
	If $k = -1$, then $U$ is an atom.
	If $\dim{U} = n$, then the equality trivially holds, since $U = \clset{x}$ for the unique $x \in \grade{n}{U}$.
	Otherwise, $\grade{n}{U} = \varnothing$, and $\isocl{f} = 0$ in $\grade{n}{\left(\linea{\molecin{P}}\right)}$ by Lemma \ref{lem:linearisation_lower_dim_cells} combined with Lemma \ref{lem:dim_molecule_as_cell}.

	Suppose that $k \geq 0$.
	Then $U$ admits a $k$\nbd layering $(\order{i}{U})_{i=1}^m$, and
	\[
		\isocl{f} = \isocl{\restr{f}{\order{1}{U}}} \cp{k} \ldots \cp{k} \isocl{\restr{f}{\order{m}{U}}} = \sum_{i=1}^m \isocl{\restr{f}{\order{i}{U}}}
	\]
	in $\grade{n}{\left(\linea{\molecin{P}}\right)}$.
	Since each $\order{i}{U}$ contains at most one $n$\nbd dimensional element, and every element of $U$ is in the image of some $\order{i}{U}$, we conclude by the inductive hypothesis.
\end{proof}

\begin{thm} \label{thm:two_functors_from_rdcpx_to_dchaug}
	Let $P$ be a regular directed complex.
	Then $\linea{\molecin{P}}$ is naturally isomorphic to $\dfreeab{P}$.
\end{thm}
\begin{proof}
	For each $n \in \mathbb{N}$, let
	\begin{align*}
		\grade{n}{\varphi}\colon \grade{n}{\freeab{P}} & \to 
		\grade{n}{\left(\linea{\molecin{P}}\right)}, \\
		x \in \grade{n}{P} & \mapsto \isocl{\clset{x} \incl P} \in \skel{n}{\molecin{P}}.
	\end{align*}
	For each $n > 0$ and $x \in \grade{n}{P}$, we have
	\[
		\der \grade{n}{\varphi}(x) = \isocl{\bound{n-1}{+}x \incl P} - \isocl{\bound{n-1}{-}x \incl P},
	\]
	which by Lemma \ref{lem:maximal_in_boundary} and Lemma \ref{lem:cells_in_linearisation_of_molecin} is equal to
	\[
		\sum_{y \in \faces{}{+}x}\isocl{\clset{y} \incl P} - \sum_{y \in \faces{}{-}x}\isocl{\clset{y} \incl P} = \grade{n-1}{\varphi}(\der x).
	\]
	Meanwhile, for each $x \in \grade{0}{P}$, we have
	\[
		\eau \grade{0}{\varphi}(x) = \eau \isocl{\clset{x} \incl P} = 1 = \eau x.
	\]
	This proves that $(\grade{n}{\varphi})_{n \in \mathbb{N}}$ is well-defined as a homomorphism of augmented chain complexes.
	Moreover, by construction, it is compatible with directions, so it lifts to a homomorphism of augmented directed chain complexes.

	Let $\isocl{f\colon U \to P}$ be an $n$\nbd dimensional cell in $\molecin{P}$.
	By Corollary \ref{cor:morphisms_of_rdcpx_are_local_isomorphisms} combined with
	Lemma \ref{lem:cells_in_linearisation_of_molecin},
	\[
		\isocl{f\colon U \to P} = \sum_{x \in \grade{n}{U}} \isocl{\clset{f(x)} \incl P} = \grade{n}{\varphi}\left( \sum_{x \in \grade{n}{U}} f(x) \right),
	\]
	which proves that $\grade{n}{\varphi}$ is surjective onto $\grade{n}{\left(\linea{\molecin{P}}\right)}$.
	Finally, observe that the assignment
	\[
		\isocl{f\colon U \to P} \mapsto \sum_{x \in \grade{n}{U}} f(x)
	\]
	is compatible with all the defining equations of $\grade{n}{\left(\linea{\molecin{P}}\right)}$, so it determines a section of $\grade{n}{\varphi}$.
	This proves that $\grade{n}{\varphi}$ is also injective.
	We conclude that $(\grade{n}{\varphi})_{n \in \mathbb{N}}$ is an isomorphism of augmented directed chain complexes.
	Finally, naturality over both maps and comaps is a straightforward consequence of Lemma \ref{lem:cells_in_linearisation_of_molecin}.
\end{proof}

\begin{comm}
	In line with our general convention, Theorem 
	\ref{thm:two_functors_from_rdcpx_to_dchaug} states that both the diagrams
	\[\begin{tikzcd}
	\rdcpxmap && \omegacat \\
	&& \dchaug
	\arrow["{\molecin{-}}", from=1-1, to=1-3]
	\arrow["{\linea{}}", from=1-3, to=2-3]
	\arrow["{\dfreeab{-}}"', curve={height=6pt}, from=1-1, to=2-3]
\end{tikzcd}
	\quad \text{and} \quad
\begin{tikzcd}
	\opp{\rdcpxcomap} && \omegacat \\
	&& \dchaug
	\arrow["{\molecin{\pb{-}}}", from=1-1, to=1-3]
	\arrow["{\linea{}}", from=1-3, to=2-3]
	\arrow["{\dfreeab{\pb{-}}}"', curve={height=6pt}, from=1-1, to=2-3]
\end{tikzcd}\]
	commute up to natural isomorphism.	
\end{comm}


\section{Steiner complexes and acyclicity} \label{sec:steiner_acyclic}

\begin{guide}
	In this section, we recall the definitions of \emph{Steiner complex} and \emph{strong Steiner complex}, and state the main theorem of Steiner theory (Theorem \ref{thm:steiner_main_theorem}).
	We prove that $\dfreeab{P}$ is a Steiner complex when $P$ is dimension-wise acyclic (Proposition 
	\ref{prop:dw_acyclic_rdcpx_gives_steiner_complex}), and deduce that in this case $\molecin{P}$ is isomorphic to $\nufun{\dfreeab{P}}$.
	
	We outline the definition of Gray products and joins of strict $\omega$\nbd categories as extensions along colimits of tensor products and joins of strong Steiner complexes.
	Finally, we prove that $\dfreeab{P}$ is a strong Steiner complex when $P$ is acyclic 
	(Proposition \ref{prop:acyclic_rdcpx_gives_strong_steiner}), and deduce that $\molecin{-}$ is compatible with Gray products and joins when restricted to acyclic regular directed complexes (Proposition \ref{prop:monoidal_functors_from_acyclic}).
\end{guide}

\begin{dfn}[Basis of an augmented directed chain complex] \index{directed chain complex!basis} \index{basis!of an augmented directed chain complex}
	Let $C$ be an augmented directed chain complex.
	A \emph{basis for $C$} is a sequence of subsets $(\grade{n}{\gener{B}} \subseteq \grade{n}{C})_{n \in \mathbb{N}}$ such that, for all $n \in \mathbb{N}$,
	\begin{enumerate}
		\item $\grade{n}{C}$ is isomorphic to $\freeab{\grade{n}{\gener{B}}}$,
		\item $\grade{n}{\dir{C}}$ is isomorphic to $\freemon{\grade{n}{\gener{B}}}$.
	\end{enumerate}
\end{dfn}

\begin{rmk}
	Not every augmented directed chain complex admits a basis, but when it does, the basis is unique: the elements of $\grade{n}{\gener{B}}$ can be characterised as the minimal elements in $\grade{n}{C}$ under the partial order defined by $x \leq y$ if and only if $y - x \in \grade{n}{\dir{C}}$.
\end{rmk}

\begin{dfn}[Augmented directed chain complex with basis]
	An \emph{augmented directed chain complex with basis} is an augmented directed chain complex which admits a basis.
\end{dfn}

\begin{lem} \label{lem:dchaug_of_ogtpos_has_a_basis}
	Let $P$ be an oriented graded poset such that $\augm{P}$ is oriented thin.
	Then $(\grade{n}{P} \subseteq \freeab{\grade{n}{P}})_{n \in \mathbb{N}}$ is a basis for $\dfreeab{P}$.
\end{lem}
\begin{proof}
	By construction.
\end{proof}

\begin{dfn}[Support of a chain] \index{directed chain complex!support of a chain}
	Let $C$ be an augmented directed chain complex with basis $(\grade{n}{\gener{B}})_{n \in \mathbb{N}}$, $n \in \mathbb{N}$, and $x \equiv \sum_{b \in \grade{n}{\gener{B}}} x_b b \in \grade{n}{C}$.
	The \emph{support of $x$} is the subset
	\[
		\supp{x} \eqdef \set{b \in \grade{n}{\gener{B}} \mid x_b \neq 0} \subseteq \grade{n}{\gener{B}}.
	\]
\end{dfn}

\begin{lem} \label{lem:positive_part_and_negative_part}
	Let $C$ be an augmented directed chain complex with basis, $n \in \mathbb{N}$, and $x \in \grade{n}{C}$.
	Then there exist unique $x^+, x^- \in \grade{n}{\dir{C}}$ such that
	\begin{enumerate}
		\item $x = x^+ - x^-$,
		\item $\supp{x^+} \cap \supp{x^-} = \varnothing$.
	\end{enumerate}
\end{lem}
\begin{proof}
	See \cite[\S 2.7]{ara2020joint}.
\end{proof}

\begin{dfn}[Positive and negative part of a chain] \index{$x^+, x^-$}
	Let $C$ be an augmented directed chain complex with basis, $n \in \mathbb{N}$, and $x \in \grade{n}{C}$.
	The \emph{positive part} and the \emph{negative part of $x$} are, respectively, the unique $x^+ \in \grade{n}{\dir{C}}$ and the unique $x^- \in \grade{n}{\dir{C}}$ such that $x = x^+ - x^-$ and $\supp{x^+} \cap \supp{x^-} = \varnothing$.
\end{dfn}

\begin{dfn}[Unital basis] \index{directed chain complex!basis!unital} \index{$\gltab{n}{\alpha}{\batom{b}}$}
	Let $C$ be an augmented directed chain complex with basis $(\grade{n}{\gener{B}})_{n \in \mathbb{N}}$.
	For all $n \in \mathbb{N}$ and $b \in \grade{n}{\gener{B}}$, let
	\[
		\gltab{m}{\alpha}{\batom{b}} \eqdef
		\begin{cases}
			0 
			& \text{if $m > n$,} \\
			b
			& \text{if $m = n$,} \\
			(\der \gltab{m+1}{\alpha}{\batom{b}})^\alpha 
			& \text{if $m < n$}
		\end{cases}
	\]
	for each $m \in \mathbb{N}$ and $\alpha \in \set{+, -}$, where the definition is obtained by downward recursion when $m \leq n$.
	We say that the basis $(\grade{n}{\gener{B}})_{n \in \mathbb{N}}$ is \emph{unital} if, for all $n \in \mathbb{N}$ and $b \in \grade{n}{\gener{B}}$, 
	\[
		\batom{b} \equiv (\gltab{m}{\alpha}{\batom{b}})_{m \in \mathbb{N}, \, \alpha \in \set{+, -}}
	\]
	is a globular table, or, equivalently, if $\eau \gltab{0}{+}{\batom{b}} = \eau \gltab{0}{-}{\batom{b}} = 1$.
\end{dfn}

\begin{lem} \label{lem:basis_atoms_in_dchaug}
	Let $P$ be a regular directed complex, $x \in P$, $m \in \mathbb{N}$, and $\alpha \in \set{+, -}$.
	Then, in $\dfreeab{P}$,
	\[
		\gltab{m}{\alpha}{\batom{x}} = \sum_{y \in \faces{m}{\alpha}x} y.
	\]
\end{lem}
\begin{proof}
	Let $n \eqdef \dim{x}$, so $x \in \grade{n}{P}$.
	By definition, for $m > n$, $\gltab{m}{\alpha}{\batom{x}} = 0$, while $\faces{m}{\alpha}x = \varnothing$, and the equality holds.
	For $m \leq n$, we proceed by downward recursion.
	If $m = n$, we have $\gltab{m}{\alpha}{\batom{x}} = x$, while $\faces{m}{\alpha}x = \set{x}$, and the equality holds.
	Let $m < n$.
	Then
	\[
		\der \gltab{m+1}{\alpha}{\batom{x}} = \der \left( \quad 
		\sum_{\mathclap{y \in \faces{m+1}{\alpha} x}} \; y \; \right)
	\]
	by the inductive hypothesis, and $\faces{m+1}{\alpha} x = \grade{m+1}{(\bound{m+1}{\alpha}x)}$ by Lemma \ref{lem:maximal_in_boundary}.
	By Lemma \ref{lem:chain_complex_boundary_of_molecule} and globularity of $\clset{x}$, this is equal to
	\[
		\sum_{\mathclap{y \in \faces{}{+}(\bound{m+1}{\alpha}x)}}\; y \quad - \;\quad
		\sum_{\mathclap{y \in \faces{}{-}(\bound{m+1}{\alpha}x)}}\; y \quad = \quad
		\sum_{\mathclap{y \in \faces{m}{+}x}}\; y \; - \quad
		\sum_{\mathclap{y \in \faces{m}{-}x}}\; y \;,
	\]
	hence by definition 
	\[
		\gltab{m}{+}{\batom{x}} = \; \sum_{\mathclap{y \in \faces{m}{+}x}}\; y\,, \quad \quad
		\gltab{m}{-}{\batom{x}} = \; \sum_{\mathclap{y \in \faces{m}{-}x}}\; y\,.
	\]
	This completes the proof.
\end{proof}

\begin{cor} \label{cor:support_is_faces}
	Let $P$ be a regular directed complex, $x \in P$, $n \in \mathbb{N}$, $\alpha \in \set{+, -}$.
	Then $\supp{\gltab{n}{\alpha}{\batom{x}}} = \faces{n}{\alpha}x$.
\end{cor}

\begin{prop} \label{prop:rdcpx_unital_basis}
	Let $P$ be a regular directed complex.
	Then $\dfreeab{P}$ has a unital basis.
\end{prop}
\begin{proof}
	Let $x \in P$.
	For all $\alpha \in \set{+, -}$, since $\clset{x}$ is an atom, by 
	Lemma \ref{lem:only_0_molecule}, $\faces{0}{\alpha}x = \bound{0}{\alpha}x = \set{x^\alpha}$ for a unique $x^\alpha \in \grade{0}{P}$.
	Then, by Lemma \ref{lem:basis_atoms_in_dchaug},
	\[
		\eau \gltab{0}{\alpha}{\batom{x}} = \eau x^\alpha = 1,
	\]
	so the basis $(\grade{n}{P})_{n \in \mathbb{N}}$ is unital.
\end{proof}

\begin{dfn}[Flow graph of an augmented directed chain complex with basis] \index{directed chain complex!flow graph} \index{flow graph!of an augmented directed chain complex}
	Let $C$ be an augmented directed chain complex with basis $(\grade{n}{\gener{B}})_{n \in \mathbb{N}}$, $k \in \mathbb{N}$.
	The \emph{$k$\nbd flow graph of $C$} is the directed graph $\flow{k}{C}$ whose
	\begin{itemize}
		\item set of vertices is $\bigcup_{i > k} \grade{i}{\gener{B}}$, and
		\item set of edges is
	    \[
		    \set{(b, c) \mid \supp{\gltab{k}{+}{\batom{b}}} \cap \supp{\gltab{k}{-}{\batom{c}}} \neq \varnothing},
		\]
	with $s\colon (b, c) \mapsto b$ and $t\colon (b, c) \mapsto c$.
	\end{itemize}
\end{dfn}

\begin{dfn}[Steiner complex] \index{directed chain complex!Steiner complex|see {Steiner complex}} \index{Steiner complex}
	A \emph{Steiner complex} is an augmented directed chain complex $C$ with a unital basis such that, for all $k \in \mathbb{N}$, $\flow{k}{C}$ is acyclic.
\end{dfn}

\begin{dfn}[The category $\stcpx$] \index{$\stcpx$}
	We let $\stcpx$ denote the full subcategory of $\dchaug$ on the Steiner complexes.
\end{dfn}

\begin{thm} \label{thm:steiner_main_theorem}
	The restriction of $\nufun{}\colon \dchaug \to \omegacat$ to $\stcpx$ is full and faithful.
	Moreover, if $C$ is a Steiner complex with basis $(\grade{n}{\gener{B}})_{n \in \mathbb{N}}$, then $\nufun{C}$ is a polygraph whose set of generating cells is
	\[
		\set{\batom{b} \mid b \in \bigcup_{n \in \mathbb{N}} \grade{n}{\gener{B}}}.
	\]
\end{thm}
\begin{proof}
	See \cite[Theorem 5.6 and Theorem 6.1]{steiner2004omega}.
\end{proof}

\begin{lem} \label{lem:flow_graph_of_dchaug_is_flow_graph_of_rdcpx}
	Let $P$ be a regular directed complex.
	Then $\flow{k}{\dfreeab{P}}$ is isomorphic to $\flow{k}{P}$.
\end{lem}
\begin{proof}
	Immediate from Corollary \ref{cor:support_is_faces}.
\end{proof}

\begin{prop} \label{prop:dw_acyclic_rdcpx_gives_steiner_complex}
	Let $P$ be a dimension-wise acyclic regular directed complex.
	Then $\dfreeab{P}$ is a Steiner complex.
\end{prop}
\begin{proof}
	Follows from Proposition \ref{prop:rdcpx_unital_basis} and Lemma \ref{lem:flow_graph_of_dchaug_is_flow_graph_of_rdcpx}.
\end{proof}

\begin{thm} \label{thm:two_omegacats_from_dw_acyclic_rdcpx}
	Let $P$ be a dimension-wise acyclic regular directed complex.
	Then $\nufun{\dfreeab{P}}$ is naturally isomorphic to $\molecin{P}$.
\end{thm}
\begin{proof}
	Composing the component $\eta\colon \molecin{P} \to \nufun{\linea{\molecin{P}}}$ of the unit of the adjunction between $\linea{}$ and $\nufun{}$ with the natural isomorphism between $\linea{\molecin{P}}$ and $\dfreeab{P}$ from Theorem \ref{thm:two_functors_from_rdcpx_to_dchaug}, we obtain a strict functor
	\[
		\varphi\colon \molecin{P} \to \nufun{\dfreeab{P}}.
	\]
	By Corollary \ref{cor:dw_acyclic_rdcpx_presents_polygraphs}, $\molecin{P}$ is a polygraph whose set of generating cells is $\set{\isocl{\clset{x} \incl P} \mid x \in P}$, while by Theorem \ref{thm:steiner_main_theorem} combined with Proposition \ref{prop:dw_acyclic_rdcpx_gives_steiner_complex}, $\nufun{\dfreeab{P}}$ is a polygraph whose set of generating cells is $\set{\batom{x} \mid x \in P}$.
By sending $\isocl{\clset{x} \incl P}$ to $\batom{x}$, $\varphi$ determines a bijection between the generating cells of $\molecin{P}$ and of $\nufun{\dfreeab{P}}$.
	By \cite[Proposition 16.2.12]{ara2023polygraphs}, we conclude that $\varphi$ is an isomorphism of polygraphs.
\end{proof}

\begin{exm}[A regular directed complex $P$ such that $\molecin{P}$ is not isomorphic to $\nufun{\dfreeab{P}}$] \index[counterex]{A regular directed complex $P$ such that $\molecin{P}$ is not isomorphic to $\nufun{\dfreeab{P}}$}
	We prove that Theorem \ref{thm:two_omegacats_from_dw_acyclic_rdcpx} does not extend beyond dimension-wise acyclic regular directed complexes.
	Let $P$ be the oriented face poset of the 1\nbd dimensional diagram
\begin{equation}\begin{tikzcd} \label{eq:nonacyclic_graph}
	{{\scriptstyle 0} \;\bullet} && {{\scriptstyle 1} \;\bullet}
	\arrow["0", from=1-1, to=1-3]
	\arrow["1"', curve={height=-18pt}, from=1-3, to=1-1]
	\arrow["2"', curve={height=18pt}, from=1-3, to=1-1]
\end{tikzcd}\end{equation}
	which is evidently not dimension-wise acyclic.
	By Corollary \ref{cor:dim3_has_frame_acyclic_molecules}, $\molecin{P}$ is a polygraph, which in this case means that $\molecin{P}$ is isomorphic to the free category on the directed graph (\ref{eq:nonacyclic_graph}).
	However, in $\nufun{\dfreeab{P}}$, let 
	\[
		x \eqdef \batom{(1, 0)} \cp{0} \batom{(1, 1)}, \quad \quad
		y \eqdef \batom{(1, 0)} \cp{0} \batom{(1, 2)},
	\]
	which as globular tables are defined, for all $\alpha \in \set{+, -}$, by
	\[
		\gltab{n}{\alpha}x \eqdef \begin{cases}
			(0, 0) & \text{if $n = 0$,} \\
			(1, 0) + (1, 1) & \text{if $n = 1$}, \\
			0 & \text{if $n > 1$},
		\end{cases} 
		\quad
		\gltab{n}{\alpha}y \eqdef \begin{cases}
			(0, 0) & \text{if $n = 0$,} \\
			(1, 0) + (1, 2) & \text{if $n = 1$}, \\
			0 & \text{if $n > 1$}.
		\end{cases}
	\]
	Then $x \cp{0} y$ and $y \cp{0} x$ are both equal to the globular table $z$ defined, for all $\alpha \in \set{+, -}$, by
	\[
		\gltab{n}{\alpha}z \eqdef  
		\begin{cases}
			(0, 0) & \text{if $n = 0$}, \\
			2(1, 0) + (1, 1) + (1, 2) & \text{if $n = 1$}, \\
			0 & \text{if $n > 1$}.
		\end{cases}
	\]
	We conclude that $\nufun{\dfreeab{P}}$ is not free, so it is not isomorphic to $\molecin{P}$.
\end{exm}

\begin{dfn}[Oriented Hasse diagram of an augmented directed chain complex with basis] \index{directed chain complex!oriented Hasse diagram} \index{oriented Hasse diagram!of an augmented directed chain complex}
	Let $C$ be an augmented directed chain complex with basis $(\grade{n}{\gener{B}})_{n \in \mathbb{N}}$.
	The \emph{oriented Hasse diagram of $C$} is the directed graph $\hasseo{C}$ whose
\begin{itemize}
	\item set of vertices is $\bigcup_{n\in \mathbb{N}} \grade{n}{\gener{B}}$,
	\item set of edges is 
	\[ 
		\set{ (b, c) \mid \text{$b \in \supp{(\der c)^-}$ or $c \in \supp{(\der b)^+}$} }, 
    	\]
	with $s\colon (b, c) \mapsto b$ and $t\colon (b, c) \mapsto c$.
\end{itemize}
\end{dfn}

\begin{dfn}[Strong Steiner complex] \index{Steiner complex!strong}
	A \emph{strong Steiner complex} is an augmented directed chain complex $C$ with a unital basis such that $\hasseo{C}$ is acyclic.
\end{dfn}

\begin{dfn}[The category $\sstcpx$] \index{$\sstcpx$}
	We let $\sstcpx$ denote the full subcategory of $\dchaug$ on the strong Steiner complexes.
\end{dfn}

\begin{prop} \label{prop:strong_steiner_is_steiner}
	Every strong Steiner complex is a Steiner complex.
\end{prop}
\begin{proof}
	See \cite[Proposition 3.7]{steiner2004omega}.
\end{proof}

\begin{comm}
	In \cite{ara2023categorical}, the authors give an internal characterisation of the essential image of $\sstcpx$ in $\omegacat$ through $\nufun{}$.
	At the end of the article, they also sketch an analogous characterisation of the essential image of $\stcpx$.
	This can be used to provide an alternative proof of Theorem 
	\ref{thm:two_omegacats_from_dw_acyclic_rdcpx}, showing directly that $\molecin{P}$ is in the essential image of $\stcpx$ when $P$ is a dimension-wise acyclic regular directed complex, then using the fact that $\linea{}$ is inverse to $\nu{}$ up to natural isomorphism on this subcategory.
\end{comm}

\begin{lem} \label{lem:monoidal_structures_on_strong_steiner}
	The monoidal structures $(\dchaug, \otimes, \mathbb{Z})$ and $(\dchaug, \join, 0)$ restrict to monoidal structures on $\sstcpx$.
\end{lem}
\begin{proof}
	See \cite[Example 3.10]{steiner2004omega} and \cite[Corollary 6.21]{ara2020joint} for the tensor product and join, respectively.
\end{proof}

\begin{prop} \label{prop:gray_product_of_omega_categories}
	There exists an essentially unique monoidal structure $(\omegacat, \gray, 1)$ on $\omegacat$ such that
	\begin{enumerate}
		\item $\nufun{}\colon (\sstcpx, \otimes, \mathbb{Z}) \to (\omegacat, \gray, 1)$ is a strong monoidal functor,
		\item for all strict $\omega$\nbd categories $X$, the functors $X \gray -$ and $- \gray X$ preserve all small colimits.
	\end{enumerate}
\end{prop}
\begin{proof}
	This is \cite[Theorem A.15]{ara2020joint}, filling some gaps in the proof of \cite[Theorem 7.3]{steiner2004omega}.
\end{proof}

\begin{rmk}
	The unit of this monoidal structure is $1 \simeq \nufun{\mathbb{Z}}$, the terminal strict $\omega$\nbd category, making the structure \emph{semicartesian} monoidal.
\end{rmk}

\begin{dfn}[Gray product of strict $\omega$-categories] \index{strict $\omega$-category!Gray product} \index{Gray product!of strict $\omega$-categories}
	Let $X$, $Y$ be strict $\omega$\nbd categories.
	The \emph{Gray product of $X$ and $Y$} is the monoidal product $X \gray Y$ of $X$ and $Y$ in the monoidal structure $(\omegacat, \gray, 1)$.
\end{dfn}

\begin{lem} \label{lem:linearisation_is_monoidal_for_gray}
	The functor $\linea{}\colon (\omegacat, \gray, 1) \to (\dchaug, \otimes, \mathbb{Z})$ is strong monoidal.
\end{lem}
\begin{proof}
	This is \cite[Proposition 2.14]{ozornova2023quillen}.
\end{proof}

\begin{prop} \label{prop:join_of_omega_categories}
	There exists an essentially unique monoidal structure $(\omegacat, \join, \varnothing)$ on $\omegacat$ such that
	\begin{enumerate}
		\item $\nufun{}\colon (\sstcpx, \join, 0) \to (\omegacat, \join, \varnothing)$ is a strong monoidal functor,
		\item for all strict $\omega$\nbd categories $X$, the functors $X \join -$ and $- \join X$ preserve all small connected colimits.
	\end{enumerate}
\end{prop}
\begin{proof}
	This is \cite[Theorem 6.29]{ara2020joint}.
\end{proof}

\begin{rmk}
	The unit of this monoidal structure is $\varnothing \simeq \nufun{0}$, the initial strict $\omega$\nbd category, making the structure \emph{semicocartesian} monoidal.
\end{rmk}

\begin{dfn}[Join of strict $\omega$\nbd categories] \index{strict $\omega$-category!join} \index{join!of strict $\omega$-categories}
	Let $X$, $Y$ be strict $\omega$\nbd categories.
	The \emph{join of $X$ and $Y$} is the monoidal product $X \join Y$ of $X$ and $Y$ in the monoidal structure $(\omegacat, \join, \varnothing)$.
\end{dfn}

\begin{lem} \label{lem:linearisation_is_monoidal_for_join}
	The functor $\linea{}\colon (\omegacat, \join, \varnothing) \to (\dchaug, \join, 0)$ is strong monoidal.
\end{lem}
\begin{proof}
	See \cite[Proposition 6.34]{ara2020joint}.
\end{proof}

\begin{lem} \label{lem:oriented_hasse_of_dchaug_of_rdcpx}
	Let $P$ be an oriented graded poset such that $\augm{P}$ is oriented thin.
	Then $\hasseo{\dfreeab{P}}$ is isomorphic to $\hasseo{P}$.
\end{lem}
\begin{proof}
	By construction, for all $x \in P$,
	\[
		\supp{(\der x)^+} = \faces{}{+}x, \quad \quad \supp{(\der x)^-} = \faces{}{-}x,
	\]
	so the definitions of $\hasseo{\dfreeab{P}}$ and of $\hasseo{P}$ coincide.
\end{proof}

\begin{prop} \label{prop:acyclic_rdcpx_gives_strong_steiner}
	Let $P$ be an acyclic regular directed complex.
	Then $\dfreeab{P}$ is a strong Steiner complex.
\end{prop}
\begin{proof}
	Follows from Proposition \ref{prop:rdcpx_unital_basis} and Lemma 
	\ref{lem:oriented_hasse_of_dchaug_of_rdcpx}.
\end{proof}

\begin{prop} \label{prop:monoidal_functors_from_acyclic}
	The functor $\molecin{-}\colon \rdcpxmapac \to \omegacat$ lifts to strong monoidal functors
	\begin{align*}
		& \molecin{-}\colon (\rdcpxmapac, \gray, 1) \to (\omegacat, \gray, 1), \\
		& \molecin{-}\colon (\rdcpxmapac, \join, \varnothing) \to (\omegacat, \join, \varnothing),
	\end{align*}
	and the functor $\molecin{\pb{-}}\colon \opp{(\rdcpxcomapac)} \to \omegacat$ lifts to strong monoidal functors
	\begin{align*}
		& \molecin{\pb{-}}\colon (\opp{(\rdcpxcomapac)}, \gray, 1) \to (\omegacat, \gray, 1), \\
		& \molecin{\pb{-}}\colon (\opp{(\rdcpxcomapac)}, \join, \varnothing) \to (\omegacat, \join, \varnothing).
	\end{align*}
\end{prop}
\begin{proof}
	Follows from Proposition 
	\ref{prop:acyclic_rdcpx_gives_strong_steiner} and Theorem 
	\ref{thm:two_omegacats_from_dw_acyclic_rdcpx}, which is applicable by Proposition 
	\ref{prop:strong_steiner_is_steiner}, in combination with
	Proposition \ref{prop:gray_to_tensor_of_dchaug} and Proposition
	\ref{prop:gray_product_of_omega_categories} for the Gray product, 
	and in combination with
	Proposition \ref{prop:join_to_join_of_dchaug} and
	Proposition \ref{prop:join_of_omega_categories} for the join.
\end{proof}

\begin{comm}
	We do not have a counterexample, but some evidence makes it seem unlikely that Proposition 
	\ref{prop:monoidal_functors_from_acyclic} extends to all regular directed complexes.
	We have no reason to believe that frame-acyclicity is stable under Gray products, because by Example \ref{exm:dw_acyclic_not_gray_stable} this is not the case for dimension-wise acyclicity, yet polygraphs are closed under Gray products 
	\cite{ara2020folk}.
	Thus a pair $P$, $Q$ of regular directed complexes with frame-acyclic molecules such that $\molecin{P \gray Q}$ is not a polygraph would be proof that $\molecin{-}$ is not, in general, compatible with Gray products.
\end{comm}

\clearpage
\thispagestyle{empty}

%% file: backmatter.tex
\phantomsection
\addcontentsline{toc}{section}{Bibliography}
\bibliographystyle{apalike}
\small \bibliography{main} 
\cleardoublepage
\thispagestyle{empty}

\phantomsection
\addcontentsline{toc}{section}{Index}
\normalsize \printindex

\cleardoublepage
\thispagestyle{empty}
\phantomsection
\addcontentsline{toc}{section}{Index of counterexamples}
\printindex[counterex]